\documentclass[reqno,11pt]{amsart}
\usepackage{amsmath}
\usepackage{amsxtra}
\usepackage{amscd}
\usepackage{amsthm}
\usepackage{amsfonts}
\usepackage{amssymb}
\usepackage{eucal}
\usepackage{scalerel}
\usepackage{verbatim}
\usepackage[matrix,arrow,curve]{xy}

\usepackage{a4wide}

\usepackage[T1]{fontenc}

\usepackage{xr-hyper}
\usepackage[
pdftex,
bookmarks=false,
colorlinks=true,
debug=true,
pdfnewwindow=true]{hyperref}

\theoremstyle{plain}
\newtheorem{Thm}[equation]{Theorem}
\newtheorem{Cor}[equation]{Corollary}
\newtheorem{Lem}[equation]{Lemma}
\newtheorem{Prop}[equation]{Proposition}
\newtheorem{Conj}[equation]{Conjecture}

\theoremstyle{definition}
\newtheorem{Def}[equation]{Definition}

\theoremstyle{remark}

\newtheorem{Rem}[equation]{Remark}

\newtheorem{conj}{Conjecture}

\numberwithin{equation}{section}
\renewcommand{\rm}{\normalshape}

\newif\ifShowLabels
\ShowLabelstrue
\newdimen\theight
\def\TeXref#1{%
    \leavevmode\vadjust{\setbox0=\hbox{{\tt
        \quad\quad  {\small \rm #1}}}%
    \theight=\ht0
    \advance\theight by \lineskip
    \kern -\theight \vbox to
    \theight{\rightline{\rlap{\box0}}%
    \vss}%
    }}%

\ShowLabelsfalse

\renewcommand{\sec}[2]{\section{#2}\label{S:#1}%
    \ifShowLabels \TeXref{{S:#1}} \fi}
\newcommand{\ssec}[2]{\subsection{#2}\label{SS:#1}%
    \ifShowLabels \TeXref{{SS:#1}} \fi}

\newcommand{\refs}[1]{Section ~\ref{S:#1}}
\newcommand{\refss}[1]{Section ~\ref{SS:#1}}

\newcommand{\reft}[1]{Theorem ~\ref{T:#1}}
\newcommand{\refl}[1]{Lemma ~\ref{L:#1}}
\newcommand{\refp}[1]{Proposition ~\ref{P:#1}}
\newcommand{\refc}[1]{Corollary ~\ref{C:#1}}
\newcommand{\refd}[1]{Definition ~\ref{D:#1}}
\newcommand{\refr}[1]{Remark ~\ref{R:#1}}
\newcommand{\refe}[1]{\eqref{E:#1}}
\newcommand{\refco}[1]{Conjecture ~\ref{Co:#1}}

\newenvironment{thm}[1]%
    { \begin{Thm} \label{T:#1}  \ifShowLabels \TeXref{T:#1} \fi }%
    { \end{Thm} }

\renewcommand{\th}[1]{\begin{thm}{#1} \sl }
\renewcommand{\eth}{\end{thm} }

\newenvironment{lemma}[1]%
    { \begin{Lem} \label{L:#1}  \ifShowLabels \TeXref{L:#1} \fi }%
    { \end{Lem} }
\newcommand{\lem}[1]{\begin{lemma}{#1} \sl}
\newcommand{\elem}{\end{lemma}}

\newenvironment{propos}[1]%
    { \begin{Prop} \label{P:#1}  \ifShowLabels \TeXref{P:#1} \fi }%
    { \end{Prop} }
\newcommand{\prop}[1]{\begin{propos}{#1}\sl }
\newcommand{\eprop}{\end{propos}}

\newenvironment{corol}[1]%
    { \begin{Cor} \label{C:#1}  \ifShowLabels \TeXref{C:#1} \fi }%
    { \end{Cor} }
\newcommand{\cor}[1]{\begin{corol}{#1} \sl }
\newcommand{\ecor}{\end{corol}}

\newenvironment{defeni}[1]%
    { \begin{Def} \label{D:#1}  \ifShowLabels \TeXref{D:#1} \fi }%
    { \end{Def} }
\newcommand{\defe}[1]{\begin{defeni}{#1} \sl }
\newcommand{\edefe}{\end{defeni}}

\newenvironment{remark}[1]%
    { \begin{Rem} \label{R:#1}  \ifShowLabels \TeXref{R:#1} \fi }%
    { \end{Rem} }
\newcommand{\rem}[1]{\begin{remark}{#1}}
\newcommand{\erem}{\end{remark}}

\newenvironment{conjec}[1]%
    { \begin{Conj} \label{Co:#1}  \ifShowLabels \TeXref{Co:#1} \fi }%
    { \end{Conj} }
\renewcommand{\conj}[1]{\begin{conjec}{#1} \sl }
\newcommand{\econj}{\end{conjec}}

\newcommand{\eq}[1]%
    { \ifShowLabels \TeXref{E:#1} \fi
       \begin{equation} \label{E:#1} }
\newcommand{\eeq}{ \end{equation} }

\newcommand{\prf}{ \begin{proof} }
\newcommand{\epr}{ \end{proof} }

\usepackage[dvipsnames]{xcolor}

\newcommand\nc{\newcommand}
\nc{\unl}{\underline}
\nc{\ol}{\overline}
\nc{\on}{\operatorname}

\nc{\BA}{{\mathbb{A}}}
\nc{\BC}{{\mathbb{C}}}
\nc{\BD}{{\mathbb{D}}}
\nc{\BF}{{\mathbb{F}}}
\nc{\BG}{{\mathbb{G}}}
\nc{\BM}{{\mathbb{M}}}
\nc{\BN}{{\mathbb{N}}}
\nc{\BO}{{\mathbb{O}}}
\nc{\BQ}{{\mathbb{Q}}}
\nc{\BP}{{\mathbb{P}}}
\nc{\BR}{{\mathbb{R}}}
\nc{\BZ}{{\mathbb{Z}}}
\nc{\BS}{{\mathbb{S}}}

\nc{\CA}{{\mathcal{A}}} \nc{\CB}{{\mathcal{B}}} \nc{\CalC}{{\mathcal
C}} \nc{\CalD}{{\mathcal D}} \nc{\CE}{{\mathcal{E}}}
\nc{\CF}{{\mathcal{F}}} \nc{\CG}{{\mathcal{G}}}
\nc{\CH}{{\mathcal{H}}} \nc{\CI}{{\mathcal{I}}}
\nc{\CK}{{\mathcal{K}}} \nc{\CL}{{\mathcal{L}}}
\nc{\CM}{{\mathcal{M}}} \nc{\CN}{{\mathcal{N}}}
\nc{\CO}{{\mathcal{O}}} \nc{\CP}{{\mathcal{P}}}
\nc{\CQ}{{\mathcal{Q}}} \nc{\CR}{{\mathcal{R}}}
\nc{\CS}{{\mathcal{S}}} \nc{\CT}{{\mathcal{T}}}
\nc{\CU}{{\mathcal{U}}} \nc{\CV}{{\mathcal{V}}}
\nc{\CW}{{\mathcal{W}}} \nc{\CX}{{\mathcal{X}}}
\nc{\CY}{{\mathcal{Y}}} \nc{\CZ}{{\mathcal{Z}}}

\nc{\fa}{{\mathfrak{a}}}
\nc{\fb}{{\mathfrak{b}}}
\nc{\fg}{{\mathfrak{g}}}
\nc{\fgl}{{\mathfrak{gl}}}
\nc{\fh}{{\mathfrak{h}}}
\nc{\fj}{{\mathfrak{j}}}
\nc{\fl}{{\mathfrak{l}}}
\nc{\fm}{{\mathfrak{m}}}
\nc{\fn}{{\mathfrak{n}}}
\nc{\fu}{{\mathfrak{u}}}
\nc{\fp}{{\mathfrak{p}}}
\nc{\frr}{{\mathfrak{r}}}
\nc{\fs}{{\mathfrak{s}}}
\nc{\ft}{{\mathfrak{t}}}
\nc{\fw}{{\mathfrak{w}}}
\nc{\fz}{{\mathfrak{z}}}

\nc{\fA}{{\mathfrak{A}}}
\nc{\fB}{{\mathfrak{B}}}
\nc{\fD}{{\mathfrak{D}}}
\nc{\fE}{{\mathfrak{E}}}
\nc{\fF}{{\mathfrak{F}}}
\nc{\fG}{{\mathfrak{G}}}
\nc{\fI}{{\mathfrak{I}}}
\nc{\fJ}{{\mathfrak{J}}}
\nc{\fK}{{\mathfrak{K}}}
\nc{\fL}{{\mathfrak{L}}}
\nc{\fM}{{\mathfrak{M}}}
\nc{\fN}{{\mathfrak{N}}}
\nc{\frP}{{\mathfrak{P}}}
\nc{\fQ}{{\mathfrak Q}}
\nc{\fR}{{\mathfrak R}}
\nc{\fS}{{\mathfrak S}}
\nc{\fT}{{\mathfrak{T}}}
\nc{\fU}{{\mathfrak{U}}}
\nc{\fW}{{\mathfrak{W}}}
\nc{\fY}{{\mathfrak{Y}}}
\nc{\fZ}{{\mathfrak{Z}}}

\nc{\ba}{{\mathbf{a}}}
\nc{\bb}{{\mathbf{b}}}
\nc{\bc}{{\mathbf{c}}}
\nc{\bd}{{\mathbf{d}}}
\nc{\be}{{\mathbf{e}}}
\nc{\bi}{{\mathbf{i}}}
\nc{\bj}{{\mathbf{j}}}
\nc{\bn}{{\mathbf{n}}}
\nc{\bp}{{\mathbf{p}}}
\nc{\bq}{{\mathbf{q}}}
\nc{\bu}{{\mathbf{u}}}
\nc{\bv}{{\mathbf{v}}}
\nc{\bw}{{\mathbf{w}}}
\nc{\bx}{{\mathbf{x}}}
\nc{\by}{{\mathbf{y}}}
\nc{\bz}{{\mathbf{z}}}

\nc{\bA}{{\mathbf{A}}}
\nc{\bB}{{\mathbf{B}}}
\nc{\bC}{{\mathbf{C}}}
\nc{\bD}{{\mathbf{D}}}
\nc{\bE}{{\mathbf{E}}}
\nc{\bI}{{\mathbf{I}}}
\nc{\bK}{{\mathbf{K}}}
\nc{\bH}{{\mathbf{H}}}
\nc{\bM}{{\mathbf{M}}}
\nc{\bN}{{\mathbf{N}}}
\nc{\bO}{{\mathbf{O}}}
\nc{\bQ}{{\mathbf Q}}
\nc{\bS}{{\mathbf{S}}}
\nc{\bT}{{\mathbf{T}}}
\nc{\bV}{{\mathbf{V}}}
\nc{\bW}{{\mathbf{W}}}
\nc{\bX}{{\mathbf{X}}}
\nc{\bP}{{\mathbf{P}}}
\nc{\bY}{{\mathbf{Y}}}
\nc{\bZ}{{\mathbf{Z}}}

\nc{\sA}{{\mathsf{A}}}
\nc{\sB}{{\mathsf{B}}}
\nc{\sC}{{\mathsf{C}}}
\nc{\sD}{{\mathsf{D}}}
\nc{\sF}{{\mathsf{F}}}
\nc{\sK}{{\mathsf{K}}}
\nc{\sM}{{\mathsf{M}}}
\nc{\sO}{{\mathsf{O}}}
\nc{\sQ}{{\mathsf{Q}}}
\nc{\sP}{{\mathsf{P}}}
\nc{\sT}{{\mathsf{T}}}
\nc{\sV}{{\mathsf{V}}}
\nc{\sW}{{\mathsf{W}}}
\nc{\sX}{{\mathsf{X}}}
\nc{\sZ}{{\mathsf{Z}}}

\nc{\sfb}{{\mathsf{b}}}
\nc{\sfc}{{\mathsf{c}}}
\nc{\sd}{{\mathsf{d}}}
\nc{\sg}{{\mathsf{g}}}
\nc{\sk}{{\mathsf{k}}}
\nc{\sfl}{{\mathsf{l}}}
\nc{\sfp}{{\mathsf{p}}}
\nc{\sr}{{\mathsf{r}}}
\nc{\st}{{\mathsf{t}}}
\nc{\sfu}{{\mathsf{u}}}
\nc{\sw}{{\mathsf{w}}}
\nc{\sz}{{\mathsf{z}}}
\nc{\sx}{{\mathsf{x}}}

\nc{\bLambda}{{\boldsymbol{\Lambda}}}
\nc{\vv}{{\boldsymbol{v}}}
\nc{\Fl}{{{\mathcal F}\ell}}
\nc{\Gr}{{\on{Gr}}}
\nc{\CHH}{{\CH\!\!\CH}}
\nc{\lambdavee}{{\lambda^{\!\scriptscriptstyle\vee}}}
\nc{\alphavee}{\alpha^{\!\scriptscriptstyle\vee}}
\nc{\rhovee}{{\rho^{\!\scriptscriptstyle\vee}}}
\newcommand\iso{\,\vphantom{j^{X^2}}\smash{\overset{\sim}{\vphantom{\rule{0pt}{0.20em}}\smash{\longrightarrow}}}\,}
\nc{\oQM}{\vphantom{j^{X^2}}\smash{\overset{\circ}{\vphantom{\vstretch{0.7}{A}}\smash{\QM}}}}
\nc{\oZ}{{}^\dagger\!\vphantom{j^{X^2}}\smash{\overset{\circ}{\vphantom{\vstretch{0.7}{A}}\smash{Z}}}}
\nc{\odZ}{{}^\dagger\!\vphantom{j^{X^2}}\smash{\overset{\circ}{\vphantom{\vstretch{0.7}{A}}\smash{\mathfrak Z}}}^{c',c}}
\nc{\bdZ}{{}^\dagger\!\vphantom{j^{X^2}}\smash{\overset{\bullet}{\vphantom{\vstretch{0.7}{A}}\smash{\mathfrak Z}}}^{c',c}}
\nc{\oS}{\vphantom{j^{X^2}}\smash{\overset{\circ}{\vphantom{\vstretch{0.7}{A}}\smash{S}}}}
\nc{\buM}{\vphantom{j^{X^2}}\smash{\overset{\bullet}{\vphantom{\vstretch{0.7}{A}}\smash{M}}}}
\nc{\dW}{{}^\dagger\ol\CW{}}
\nc{\hW}{{}^\dagger\hat\CW{}}
\nc{\wW}{{}^\dagger\widetilde\CW{}}
\nc{\dZ}{{}^\dagger\!\fZ^{c',c}}
\nc{\dZc}{{}^\dagger\!\fZ^{c,c}}
\nc{\tZ}{{}^\dagger\!\tilde{Z}{}}
\nc{\hZ}{{}^\dagger\!\hat{Z}{}}

\nc{\ssl}{\mathfrak{sl}} \nc{\gl}{\mathfrak{gl}}
\nc{\wt}{\widetilde} \nc{\Sym}{\mathrm{Sym}} \nc{\Res}{\mathrm{Res}}
\nc{\sE}{{\mathsf{E}}} \nc{\bs}{{\mathbf{s}}}
\nc{\trig}{\mathrm{trig}} \nc{\rat}{\mathrm{rat}}
\nc{\sign}{\mathrm{sign}} \nc{\sL}{{\mathsf{L}}}
\nc{\fv}{{\mathfrak{v}}} \nc{\ad}{\mathrm{ad}}
\nc{\spsi}{{\mathsf{\psi}}} \nc{\sh}{{\mathsf{h}}}
\nc{\rtt}{\mathrm{rtt}} \nc{\qdet}{\mathrm{qdet}}
\nc{\M}{\mathrm{M}} \nc{\Ker}{\mathrm{Ker}} \nc{\ssc}{\mathrm{sc}}
\nc{\loc}{\mathrm{loc}} \nc{\fra}{\mathrm{frac}}
\nc{\ddj}{\mathrm{DJ}}

\begin{document}
\title[Multiplicative slices, relativistic Toda and shifted quantum algebras]
{Multiplicative slices, relativistic Toda and shifted quantum affine algebras}

\author{Michael Finkelberg}
 \address{M.F.:
  National Research University Higher School of Economics, Russian Federation,
  Department of Mathematics, 6 Usacheva st., Moscow 119048;
  Skolkovo Institute of Science and Technology;
  Institute for Information Transmission Problems}
 \email{fnklberg@gmail.com}

\author[Alexander Tsymbaliuk]{Alexander Tsymbaliuk}
 \address{A.T.:  Yale University, Department of Mathematics, New Haven, CT 06511, USA}
 \email{sashikts@gmail.com}

\dedicatory{To Tony Joseph on his 75th birthday, with admiration}

\begin{abstract}
We introduce the shifted quantum affine algebras. They map homomorphically into
the quantized $K$-theoretic Coulomb branches of $3d\ {\mathcal N}=4$ SUSY quiver
gauge theories. In type $A$, they are endowed with a coproduct,
and they act on the equivariant
$K$-theory of parabolic Laumon spaces. In type $A_1$, they are closely related
to the type $A$ open relativistic quantum Toda system.
\end{abstract}
\maketitle
\tableofcontents


\sec{int}{Introduction}

\ssec{11}{Summary}
\

The goal of this paper is to initiate the study of
{\em shifted quantum affine algebras}\footnote{They were introduced by
B.~Feigin in 2010.} and {\em shifted $\vv$-Yangians}.
They arise as a tool to write down
via generators and relations the quantized $K$-theoretic Coulomb branches of
$3d\ \CN=4$ SUSY quiver gauge theories (see~\cite[Remark~3.9(2)]{bfn}),
similarly to the appearance of shifted Yangians in the
study of the quantized Coulomb branches of $3d\ \CN=4$ SUSY quiver gauge
theories~\cite{bfn}.\footnote{We must admit right away that we were not able to
prove the desired presentation of the quantized Coulomb branch for a single quiver.}
Similarly to~\cite{fkp}, the shifted quantum affine algebras carry
a coproduct, see~Section~\ref{Copro} for partial results in this direction.
The multiplicative analogue of the construction~\cite{bdghk}
equips the equivariant $K$-theory of parabolic Laumon spaces with an action
of the quantized $K$-theoretic Coulomb branch for a type $A$ quiver, and hence
with an action of a shifted quantum affine algebra of type $A$.
Similarly to~\cite{fkp}, the unframed case of type $A_1$ quiver is closely
related to the open relativistic quantum Toda system of type~$A$.


\ssec{12}{Outline of the paper}
\

\noindent
$\bullet$
In~\refs{RLT}, we give a construction of the completed phase space of the
(quasiclassical) relativistic open Toda system for arbitrary simply-connected
semisimple algebraic group $G$ via quasihamiltonian and Poisson
reductions. It is a direct multiplicative analogue of the Kazhdan-Kostant
construction of the (nonrelativistic) open Toda integrable system. We want
to stress right away that it depends on a choice of a pair of Coxeter elements
in the Weyl group $W$ of $G$, via a choice of Steinberg's
cross-section.\footnote{The appearance of Coxeter elements in the construction of
relativistic Toda lattice goes back at least to~\cite{stss}.} In the case when the
two Coxeter elements coincide, the resulting
completed phase space is isomorphic to the universal centralizer $\fZ^G_G$,
see~\refss{unicentr}. In the case $G=SL(n)$, the universal centralizer is
isomorphic to a natural $n$-fold cover of the moduli space of centered
periodic $SU(2)$-monopoles of charge $n$, see~\refc{consol}.

\medskip
\noindent
$\bullet$
The conjectural quantization of the above construction of the completed phase
space of the relativistic open Toda is described in~\refss{quantum}.
We conjecture that it is isomorphic to the corresponding spherical symmetric nil-DAHA which
is realized as an equivariant $K$-theory of a twisted affine Grassmannian,
i.e.\ as a sort of twisted quantized Coulomb branch (the twist is necessary
in the case of non-simply-laced $G$). The bulk of~\refs{Nildaha} is occupied
by the review of Cherednik's definition of symmetric nil-DAHA, its residue
construction, and its realization as the equivariant $K$-theory of a twisted
affine flag variety. In the simply-laced case no twist is required, and the
spherical nil-DAHA in question is isomorphic to the convolution algebra
$K^{G(\CO)\rtimes\BC^\times}(\Gr_G)$ up to some finite extension. This convolution
algebra is defined for arbitrary reductive $G$. In case $G=GL(n)$, this
convolution algebra is likely to have a presentation via generators and
relations (as a truncated shifted quantum affine algebra of type $A_1$),
see~\refs{sl2}. From this presentation and~Proposition~\ref{descent for sl2}
we obtain a homomorphism
$K^{G(\CO)\rtimes\BC^\times}(\Gr_G)\to K^{L(\CO)\rtimes\BC^\times}(\Gr_L)$
for any Levi subgroup $L\subset G=GL(n)$. We conjecture an existence of such
a homomorphism for arbitrary Levi subgroup $L$ in arbitrary reductive group
$G$, but we have no clue as to a geometric construction of such a
homomorphism. It would be important for a study of equivariant quantum
$K$-theory of the flag variety $\CB$ of $G$. Its analogue for the
equivariant Borel-Moore homology convolution algebra
$H_\bullet^{G(\CO)\rtimes\BC^\times}(\Gr_G)\to H_\bullet^{L(\CO)\rtimes\BC^\times}(\Gr_L)$
is constructed in~\cite{fkp}. However, the construction is not geometric;
it uses an isomorphism with the quantum open (nonrelativistic) Toda lattice.

\medskip
\noindent
$\bullet$
Recall that for an arbitrary $3d\ \CN=4$ SUSY quiver gauge theory of type $ADE$,
the non-quantized $K$-theoretic Coulomb branch is identified with a
multiplicative generalized
slice in the corresponding affine Grassmannian~\cite[Remarks~3.9(2),~3.17]{bfn}.
These multiplicative slices are studied in detail in~\refs{slices}
(in the unframed case, they were studied in detail in~\cite{fkr}).
In particular, they embed into the loop group $G(z)$, and it is likely that
the image coincides with the space of scattering matrices of singular periodic
monopoles~\cite{ck}. The multiplication in the loop group gives rise to
the multiplication of slices, which is conjecturally
quantized by the coproduct of the corresponding shifted quantum affine algebras.

\medskip
\noindent
$\bullet$
In Section~\ref{sqaa}, we introduce the \emph{shifted quantum affine
algebras} $\CU^\ssc_{\mu^+,\mu^-}$ and $\CU^\ad_{\mu^+,\mu^-}$
(\emph{simply-connected} and \emph{adjoint} versions, respectively)
for any simple Lie algebra $\fg$ and its two coweights
$\mu^+,\mu^-$ (these algebras depend only on $\mu=\mu^+ + \mu^-$
up to an isomorphism). For $\mu^+=\mu^-=0$, they are central
extensions of the standard quantum loop algebra $U_\vv(L\fg)$ and
its adjoint version $U^\ad_\vv(L\fg)$. These algebras can be viewed
as trigonometric versions of the shifted Yangians $\bY_\mu$,
see~\cite{kwwy,bfn,fkp}.

An alternative (but equivalent) definition of $\CU^\ssc_{\mu^+,\mu^-}$
was suggested to us by B.~Feigin in Spring 2010 in an attempt to generalize
the results of~\cite{bffr} to the $K$-theoretic setting (which is the subject of
Section~\ref{Laumon} of the present paper). In this approach, we
consider an algebra with the same generators and defining relations
as $U_\vv(L\fg)$ in the new Drinfeld realization with just one
modification: the relation
  $[e_i(z),f_j(w)]=\frac{\delta_{ij}\delta(z/w)}{\vv_i-\vv_i^{-1}}\left(\psi^+_i(z)-\psi^-_i(z)\right)$
is replaced by
  $p_i(z)[e_i(z),f_j(w)]=\frac{\delta_{ij}\delta(z/w)}{\vv_i-\vv_i^{-1}}\left(\psi^+_i(z)-\psi^-_i(z)\right)$
for any collection of rational functions $\{p_i(z)\}_{i\in I}$ (here
$I$ parametrizes the set of vertices of the Dynkin diagram of $\fg$).
For $\fg=\ssl_2$ and $\mu^+=\mu^-\in-\BN$, the algebra $\CU^\ssc_{\mu^+,\mu^-}$ appeared
in~\cite[5.2]{dfk}.

We also provide an alternative presentation of the antidominantly shifted quantum affine algebras
with a finite number of generators and defining relations, see Theorem~\ref{Levendorskii-type} and
Appendix~\ref{Proof of Levendorskii-type} for its proof.
We note that this result (and its proof) also holds for any affine Lie algebra, except for type
$A^{(1)}_1$. In the unshifted case, more precisely for $U_\vv(L\fg)$, it can be viewed
as a $\vv$-version of the famous Levendorskii presentation of the
Yangian $Y(\fg)$, see~\cite{le}. Motivated by~\cite{gnw}, we also
provide a slight modification of this presentation in
Theorem~\ref{Levendorskii-type update}.

\medskip
\noindent
$\bullet$
In Section~\ref{ABCD generators}, we introduce other generators of
$\CU^\ad_{\mu^+,\mu^-}$, which can be encoded by the generating
series $\{A^\pm_i(z),B^\pm_i(z),C^\pm_i(z),D^\pm_i(z)\}_{i\in I}$.
We provide a complete list of the defining relations between these
generators for antidominant $\mu^+,\mu^-\in \Lambda^-$ (we use $\Lambda^-$ to denote
the submonoid of the coweight lattice $\Lambda$ spanned by antidominant coweights),
see Theorem~\ref{Relations bw ABCD} and Appendix~\ref{Proof of Relations
bw ABCD} for its proof. This should be viewed as a shifted
$\vv$-version of the corresponding construction for Yangians
of~\cite{gklo1}. We note that while some of the relations were
established (without a proof) in \emph{loc.~cit.}, the authors did
not aim at providing a complete list of the defining relations.
However, a rational analogue of Theorem~\ref{Relations bw ABCD}
provides such a list.

We would like to point out that this is one of the few places where it is
essential to work with the \emph{adjoint} version. In the simplest case,
that is of $U^\ad_\vv(L\ssl_2)$, these generating series coincide with the
entries of the matrices $T^\pm(z)$ from the RTT realization of
$U^\ad_\vv(L\ssl_2)$, see~\cite{df} and our discussion in
Section~\ref{section RTT algebras}.

\medskip
\noindent
$\bullet$
In Section~\ref{Section with Homomorphism}, we construct homomorphisms
  $$\wt{\Phi}^{\unl\lambda}_\mu\colon \CU^\ad_{0,\mu}[\sz^{\pm 1}_1,\ldots,\sz^{\pm 1}_N]
    \longrightarrow \wt{\CA}^\vv_\fra[\sz^{\pm 1}_1,\ldots,\sz^{\pm 1}_N]$$
from the adjoint version of shifted quantum affine algebras to the
$\BC(\vv)[\sz^{\pm 1}_1,\ldots,\sz^{\pm 1}_N]$-algebras
$\wt{\CA}^\vv_\fra[\sz^{\pm 1}_1,\ldots,\sz^{\pm 1}_N]$ of
difference operators on multidimensional tori, see
Theorem~\ref{Homomorphism} and Appendix~\ref{Proof Homomorphism} for
its proof. Here $\unl{\lambda}=(\omega_{i_1},\ldots,\omega_{i_N})$
is a sequence of fundamental coweights, such that $\lambda-\mu$ is a
sum of simple coroots with coefficients in $\BN$, where
$\lambda:=\sum_{s=1}^N\omega_{i_s}$. This result can be viewed as a
$\vv$-version of the corresponding construction for shifted Yangians
of~\cite[Theorem B.15]{bfn}, while the \emph{unshifted case} of it,
more precisely the case of $U_\vv(L\fg)$, appeared (without a proof) in~\cite{gklo}.
For $\fg=\ssl_2,\ N=0$ and antidominant shift, the above homomorphism made its
first appearance in~\cite[Section~6]{dfk}.

\medskip
\noindent
$\bullet$
In Section~\ref{Coulomb}, we consider the quantized $K$-theoretic Coulomb branch
$\CA^\vv$ in the particular case of quiver gauge theories of ADE
type (a straightforward generalization of the constructions
of~\cite{bfna,bfn}, with the equivariant Borel-Moore homology replaced
by the equivariant $K$-theory). There is a natural embedding $\bz^*(\iota_*)^{-1}
\colon \CA^\vv\hookrightarrow\wt\CA^\vv[\sz^{\pm 1}_1,\ldots,\sz^{\pm 1}_N]$.
In~Theorem~\ref{Factored Homomorphism}, we show that our homomorphism
$\wt\Phi{}^{\unl\lambda}_\mu$ of Section~\ref{Section with Homomorphism}
factors through the above embedding (with $\BC[\vv^{\pm 1}]$ extended to $\BC(\vv)$),
giving rise to a homomorphism
  $$\ol\Phi{}^{\unl\lambda}_\mu\colon  \CU^\ad_{0,\mu}[\sz^{\pm 1}_1,\ldots,\sz^{\pm 1}_N]
    \longrightarrow \CA^\vv_\fra.$$
This is a $\vv$-version of the corresponding result for shifted
Yangians of~\cite[Theorem B.18]{bfn}.

In~Section~\ref{trun}, we add certain truncation relations to the relations
defining $\CU^\ad_{0,\mu}[\sz^{\pm 1}_1,\ldots,\sz^{\pm 1}_N]$ to obtain the
{\em truncated shifted quantum affine algebras}
$\CU^{\unl\lambda}_\mu$ such that the homomorphism $\ol\Phi{}^{\unl\lambda}_\mu$
factors through the projection and the same named homomorphism
  $\CU^\ad_{0,\mu}[\sz^{\pm 1}_1,\ldots,\sz^{\pm 1}_N]\twoheadrightarrow
   \CU^{\unl\lambda}_\mu\stackrel{\ol{\Phi}^{\unl\lambda}_\mu}{\longrightarrow}
   \CA^\vv_\fra$.
We expect that $\ol\Phi{}^{\unl\lambda}_\mu\colon\CU^{\unl\lambda}_\mu\to\CA^\vv_\fra$
is an isomorphism, see~\refco{Phi iso}.

In~\refss{yang}, we define the {\em shifted $\vv$-Yangians}
  $_\bi\!\CY^\vv_\mu[\sz^{\pm 1}_1,\ldots,\sz^{\pm 1}_N]\subset
   \CU^\ad_{0,\mu}[\sz^{\pm 1}_1,\ldots,\sz^{\pm 1}_N]$
and their truncated quotients
$_\bi\!\CY^{\unl{\lambda}}_\mu\subset\CU^{\unl\lambda}_\mu$. We conjecture that
$\ol\Phi{}^{\unl\lambda}_\mu\colon_\bi\!\CY^{\unl{\lambda}}_\mu\to\CA^\vv_\fra$
is an isomorphism, see~\refco{Phi iso Yang}.

One of our biggest failures is the failure to define the integral forms
$_\bi\!\fY^{\unl\lambda}_\mu\subset{}_\bi\!\CY^{\unl\lambda}_\mu$ and
$\fU^{\unl\lambda}_\mu\subset\CU^{\unl\lambda}_\mu$ over
$\BC[\vv^{\pm1}]\subset\BC(\vv)$ that would (at least conjecturally)
map isomorphically onto $\CA^\vv\subset\CA^\vv_\fra$.
Only in the case of $\fg={\mathfrak{sl}}_2$,
making use of the $ABCD$-generators of~Section~\ref{ABCD generators},
we are able to introduce the desired integral form
in~\refss{integral form} (see also~\cite{fits} for the integral forms for $\fg=\ssl_n$).
It is worth noting that for arbitrary simply-laced $\fg$ and any $i\in I$,
the images under $\ol\Phi{}^{\unl\lambda}_\mu$ of the generators
$B_{i,r}^+$ and $e_{i,r}$ (resp.\ $C_{i,r}^+$ and $f_{i,r}$) are the
classes of dual exceptional collections of vector bundles on the
corresponding minuscule Schubert varieties in the affine
Grassmannian, see~\refr{exceptional}.

The desired integral forms $_\bi\!\fY^{\unl{\lambda}}_\mu$ and $\fU^{\unl\lambda}_\mu$
are expected to be quantizations of a certain
cover $\hW^{\unl\lambda^*}_{\mu^*}$ of a multiplicative slice introduced
in~\refss{cover}, see~\refco{master conj}. Here $*$ stands for the involution
$\mu\mapsto-w_0\mu$ of the coweight lattice $\Lambda$.

\medskip
\noindent
$\bullet$
In~\refs{sl2}, we prove the surjectivity of the homomorphism
$\ol{\Phi}^0_{-n\alpha}$ in the simplest case of $\fg=\ssl_2$ and
antidominant shifts, see Theorem~\ref{surj}. This identifies the
slightly localized and extended quantized $K$-theoretic Coulomb branch
$K^{\wt{GL}(n,\CO)\rtimes\wt\BC^\times}_{\on{loc}}(\Gr_{GL(n)})$ with a quotient
of the localized version of the truncated shifted quantum affine algebra
$\fU^0_{-n\alpha,\loc}$ (where $\wt{GL}(n)$ and $\wt\BC^\times$ stand for the
two-fold covers of $GL(n),\BC^\times$; while the localization is obtained by
inverting $1-\vv^{2m},\ 1\leq m\leq n$). We reduce
the proof of the isomorphism $\fU^0_{-n\alpha,\loc}\iso
K^{\wt{GL}(n,\CO)\rtimes\wt\BC^\times}_{\on{loc}}(\Gr_{GL(n)})$ to a verification
of an identity with quantum resultants in $\CU^0_{-n\alpha}$,
see~\refr{quantum result},~\refr{identity}. It would be interesting to
describe explicitly a basis of $\fU^0_{-n\alpha,\loc}$ projecting to the
`canonical' basis of $K^{\wt{GL}(n,\CO)\rtimes\wt\BC^\times}_{\on{loc}}(\Gr_{GL(n)})$
formed by the classes of irreducible equivariant perverse coherent sheaves~\cite{bfm}.

\medskip
\noindent
$\bullet$
In Section~\ref{Copro}, we discuss generalizations of the
classical coproducts on $U_\vv(L\fg)$ to the shifted setting.
We start by considering the simplest case $\fg=\ssl_2$.
We will denote $\CU^\ssc_{0,b\alpha/2}$ simply by $\CU^\ssc_{0,b}$
(here $b\in \BZ$ and $\alpha$ is the simple positive coroot). We
construct homomorphisms
  $$\Delta_{b_1,b_2}\colon  \CU^\ssc_{0,b}\longrightarrow
    \CU^\ssc_{0,b_1}\otimes \CU^\ssc_{0,b_2}$$
for any $b_1,b_2\in \BZ$, which recover the classical
Drinfeld-Jimbo coproduct for $b_1=b_2=0$.
Our construction is parallel to the one for shifted Yangians
of~\cite{fkp} and proceeds in two steps.
First, we define such homomorphisms in the \emph{antidominant}
case $b_1,b_2\in \BZ_{\leq 0}$, see Theorem~\ref{Coproduct for sl2} and
Appendix~\ref{Proof Coproduct} for its proof.
The proof is crucially based on the aforementioned
alternative presentation of the antidominantly shifted quantum affine algebras
with a finite number of generators and defining relations of Theorem~\ref{Levendorskii-type}.
Second, we use the algebra embeddings
  $\iota_{n,m_1,m_2}\colon  \CU^\ssc_{0,n}\hookrightarrow \CU^\ssc_{0,n+m_1+m_2}$
(here $m_1,m_2\leq 0$) to reduce the general case to the
antidominant one, see Theorem~\ref{Coproduct for sl2 in general} and
Appendix~\ref{Proof Coproduct General} for its proof.
We note that our proof of injectivity of the {\em shift homomorphisms} $\iota_{n,m_1,m_2}$
is based on the PBW property of the shifted quantum affine algebras of $\ssl_2$, see
Lemma~\ref{Properties of iota} and Theorem~\ref{PBW 2} of Appendix~\ref{Proof of PBW}.

In Sections~\ref{Coproduct on sl_n for antidominant} and~\ref{General shifted coproduct sl_n},
we generalize the aforementioned case of $\ssl_2$ to the case of $\ssl_n\ (n\geq 2)$.
The idea is again to treat first the case of antidominant shifts and then deduce the general case.
To achieve the former goal, it is essential to have explicit formulas for the action of the
Drinfeld-Jimbo coproduct on the generators $\{e_{i,-1}, f_{i,1}, h_{i,\pm 1}\}_{i\in I}$ of
$U_\vv(L\ssl_n)$. This is the key technical result,
stated in Theorem~\ref{Coproduct on quantum loop sln} and proved in
Appendix~\ref{Proof of coproduct on loop sln}. Once this is established,
it is easy to guess the formulas for the homomorphism
  $\Delta_{\mu_1,\mu_2}\colon \CU^\ssc_{0,\mu_1+\mu_2}\to
   \CU^\ssc_{0,\mu_1}\otimes \CU^\ssc_{0,\mu_2}$
in the case $\mu_1,\mu_2\in \Lambda^-$ (antidominant), see
Theorem~\ref{coproduct on antidominant sln} and its proof
in Appendix~\ref{Proof of coproduct antidominant sl_n}.
In~Theorem~\ref{Coproduct for sln in dominant case}
we derive the construction of $\Delta_{\mu_1,\mu_2}$ for general $\mu_1,\mu_2\in \Lambda$
by utilizing the algebra embeddings
  $\iota_{\mu,\nu_1,\nu_2}\colon  \CU^\ssc_{0,\mu}\hookrightarrow \CU^\ssc_{0,\mu+\nu_1+\nu_2}$
for $\mu\in \Lambda, \nu_1,\nu_2\in \Lambda^-$, see~Theorem~\ref{1026} and its proof
in Appendix~\ref{Proof of injectivity} (the latter is based on the shuffle realization
of $U_\vv(L\ssl_n)$ of~\cite{n,t3}).

Motivated by~\cite{fkp}, we expect that our construction of homomorphisms
$\Delta_{\mu_1,\mu_2}$ can be generalized to any simply-laced $\fg$ and its two
coweights $\mu_1,\mu_2\in \Lambda$. However, we failed to achieve this due to a lack of
explicit formulas for the Drinfeld-Jimbo coproduct of the
generators $\{e_{i,-1},f_{i,1},h_{i,\pm 1}\}_{i\in I}$ of $U_\vv(L\fg)$
(even for $\fg=\ssl_n$, the formulas of~Theorem~\ref{Coproduct on quantum loop sln} seem to be new,
to our surprise).

Moreover, we expect that this coproduct extends to
  $$\Delta^{\on{ad}}_{\mu_1,\mu_2}\colon
    \CU^{\on{ad}}_{0,\mu_+\mu_2}[\sz_1^{\pm1},\ldots,\sz_{N_1+N_2}^{\pm1}]\longrightarrow
    \CU^{\on{ad}}_{0,\mu_1}[\sz_1^{\pm1},\ldots,\sz_{N_1}^{\pm1}]\otimes
    \CU^{\on{ad}}_{0,\mu_2}[\sz_{N_1+1}^{\pm1},\ldots,\sz_{N_1+N_2}^{\pm1}],$$
which descends to the same named homomorphism
  $\Delta^{\on{ad}}_{\mu_1,\mu_2}\colon \CU^{\unl{\lambda}}_{\mu_1+\mu_2}\to
   \CU^{\unl{\lambda}\!^{(1)}}_{\mu_1}\otimes\CU^{\unl{\lambda}\!^{(2)}}_{\mu_2}$
between truncated algebras, see~\refco{truncopr}. We check a particular case of
this conjecture for $\fg=\ssl_2$ in~Proposition~\ref{descent for sl2}, using
the RTT realization of $\CU^{\on{ad}}_{0,2b}$ of~Theorem~\ref{Ding-Frenkel 2}.

\medskip
\noindent
$\bullet$
In Section~\ref{Ubiquity}, we discuss relativistic/trigonometric
Lax matrices, the shifted RTT algebras of $\ssl_2$ and their
relation to the shifted quantum affine algebras of $\ssl_2$.  This
yields a link between two seemingly different appearances of the
RTT relations (both trigonometric and rational).

In Section~\ref{relativistic}, we recall the Kuznetsov-Tsyganov~\cite{kt} local
relativistic Lax matrix $L^{\vv,0}_i(z)$ satisfying the trigonometric
RTT-relation. The complete monodromy matrix
$T^{\vv,0}_n(z)=L^{\vv,0}_n(z)\cdots L^{\vv,0}_1(z)$ also satisfies the same relation,
and its matrix coefficient $T^{\vv,0}_n(z)_{11}$ encodes all the
hamiltonians of the $q$-difference quantum open Toda lattice for
$GL(n)$~\cite{e,se}.

We introduce two more local Lax matrices $L^{\vv,\pm1}_i(z)$ satisfying the same
trigonometric RTT-relation. They give rise to the plethora of $3^n$ complete
monodromy matrices $T^\vv_{\vec{k}}(z), \vec{k}\in \{-1,0,1\}^n,$ given by the length $n$ products of the
three local Lax matrices in arbitrary order. The matrix coefficient
$T^\vv_{\vec{k}}(z)_{11}$ encodes the hamiltonians of the corresponding
modified quantum difference Toda lattice; the quadratic hamiltonians are
given by the formula~\refe{mixedH2}. At the quasiclassical level, these
integrable systems go back to~\cite{fg}. We show that among these $3^n$
integrable systems there are no more than $3^{n-2}$ nonequivalent, see Lemma~\ref{conjugation}.
It is shown in~\cite{fgt} that they are all obtained by the construction of~\cite{se}
using arbitrary pairs of orientations of the $A_{n-1}$ Dynkin diagram,
see~\refr{Et-Se}.



In Section~\ref{section RTT algebras}, we introduce the
\emph{shifted RTT algebras of $\ssl_2$}, denoted by
$\CU^\rtt_{0,-2n}$, and construct isomorphisms
$\Upsilon_{0,-2n}\colon \CU^\ad_{0,-2n}\iso \CU^\rtt_{0,-2n}$
for any $n\in \BN$, see Theorem~\ref{Ding-Frenkel}
and~Theorem~\ref{Ding-Frenkel 2}. For $n=0$, this recovers the isomorphism of
the \emph{new Drinfeld} and the \emph{RTT} realizations of the quantum
loop algebra $U^\ad_\vv(L\ssl_2)$, due to~\cite{df}. We also
identify the $ABCD$ generators of $\CU^\ad_{0,-2n}$
of~Section~\ref{ABCD generators} with the generators of
$\CU^\rtt_{0,-2n}$, see~Corollary~\ref{ABCD in rtt}.

Viewing the Lax matrix $L^{\vv,-1}_1(z)$ as a homomorphism from
$\CU^\rtt_{0,-2}$ to the algebra of difference operators on
$\BC^\times$ and composing it with $\Upsilon_{0,-2}$, we recover the
homomorphism $\wt\Phi^0_{-2}$ of Section~\ref{Section with Homomorphism}.
More generally, among all pairwise isomorphic shifted algebras $\{\CU^\ad_{b,-2-b}|b\in \BZ\}$
only those with $b,-2-b\leq 0$ admit an RTT realization, i.e., there are analogous isomorphisms
$\Upsilon_{b,-2-b}\colon \CU^\ad_{b,-2-b}\iso \CU^\rtt_{b,-2-b}$. Moreover,
recasting the homomorphisms $\wt{\Phi}_{b,-2-b}$ (generalizations of $\wt\Phi^0_{-2}$ for $b=0$)
as the homomorphisms $\CU^{\rtt}_{b,-2-b}\to \hat{\CA}^\vv_1$, we recover the other two Lax matrices
$L^{\vv,0}_1(z)$ (for $b=-1$) and $L^{\vv,1}_1(z)$ (for $b=-2$).

Finally, we use the RTT presentation of $U^\ad_\vv(L\ssl_2)$ to derive
explicit formulas for the action of the Drinfeld-Jimbo coproduct on the
Drinfeld \emph{half-currents}, see Proposition~\ref{nonshifted
coproduct via currents} and Appendix~\ref{Proof of coproduct via
Drinfeld currents} for its proof. We also show that the same
formulas hold in the antidominantly shifted setting for the
homomorphisms $\Delta_{b_1,b_2}$, see Proposition~\ref{shifted
coproduct via currents}. As a consequence of the latter, the homomorphism
$\Delta_{2b_1,2b_2}^\ad$ is intertwined with the RTT coproduct $\Delta^\rtt_{2b_1,2b_2}$,
see Corollary~\ref{compatibility of coproducts}, which is used to prove the
aforementioned Proposition~\ref{descent for sl2} on the descent of
$\Delta_{2b_1,2b_2}^\ad$ to the truncated versions.

\medskip
\noindent
$\bullet$
In Section~\ref{Laumon}, we provide yet another geometric realization of
the shifted quantum affine algebras (resp.\ shifted Yangians) of
$\ssl_n$ via the parabolic Laumon spaces. Roughly speaking, this
arises by combining our homomorphism $\ol{\Phi}^{\unl\lambda}_\mu$
of Section~\ref{Coulomb} (resp.\ $\ol{\Phi}^{\unl\lambda}_\mu$
of~\cite[Theorem~B.18]{bfn}) with an action of the quantized $K$-theoretic (resp.\
cohomological) Coulomb branch $\CA^\vv_\fra$ on the localized
equivariant $K$-theory (resp.\ cohomology) of parabolic Laumon spaces,
constructed in~\cite{bdghk}, see Remark~\ref{Relation to bdghk}(c).

For any $\pi=(p_1,\ldots,p_n)\in \BZ_{>0}^n$, we construct an action of
$\CU^\ssc_{0,\mu}$, the simply-connected shifted quantum affine algebra
of $\ssl_n$ with the shift $\mu=\sum_{j=1}^{n-1}(p_{j+1}-p_j)\omega_j$,
on $M(\pi)$: the direct sum of localized equivariant $K$-theory of $\fQ_{\unl{d}}$,
see Theorem~\ref{Action on K-theory of Laumon}. Here $\fQ_{\unl{d}}$ is
the type $\pi$ Laumon based parabolic quasiflags' space, which we
recall in Section~\ref{Parabolic Laumon spaces}.
In Theorem~\ref{action of gl_n on K-theory}, we slightly generalize this
by constructing an action of the shifted quantum affine algebra of
$\gl_n$ (defined in Section~\ref{shifted gl_n}) on $M(\pi)$.
In Theorem~\ref{Tensor products}, we establish an isomorphism
$M(\pi')\otimes M(\pi'')\iso M(\pi)$ (here $\pi=\pi'+\pi''$) of
$\CU^\ssc_{0,\mu}$-modules, where the action on the source arises from
the formal coproduct
  $\wt{\Delta}\colon  \CU^\ssc_{0,\mu}\to
   \CU^\ssc_{0,\mu'}\widehat{\otimes} \CU^\ssc_{0,\mu''},$
constructed in Section~\ref{Drinfeld coproduct} (an analogue of the
Drinfeld formal coproduct on $U_\vv(L\fg)$).

The \emph{rational} counterpart of these results is established in
Theorem~\ref{Action on cohomology of Laumon}, where we construct an
action of $\CY^\hbar_\mu$ (the shifted Yangian of $\ssl_n$ with scalars
extended to $\BC(\hbar)$) on $V(\pi)$: the sum of localized equivariant cohomology of
$\fQ_{\unl{d}}$. The dominant case $(p_1\leq\ldots\leq p_n)$ of this result
was treated in~\cite{bffr}, where the proof was deduced from the
Gelfand-Tsetlin formulas of~\cite{fmo}. In contrast, our straightforward proof
is valid for any $\pi$ and, thus,  gives an alternative proof of the above
Gelfand-Tsetlin formulas. We also propose a $\vv$-analogue of the
Gelfand-Tsetlin formulas of~\cite{fmo}, see~Proposition~\ref{qGT}.

Our construction can be also naturally generalized to provide the
actions of the shifted quantum toroidal (resp.\ affine Yangian) algebras of
$\ssl_n$ on the sum of localized equivariant $K$-theory (resp.\ cohomology) of
the parabolic affine Laumon spaces, see Section~\ref{affine Laumon}.

In Section~\ref{section on whittaker}, we introduce the \emph{Whittaker} vectors
in the completions of $M(\pi)$ and $V(\pi)$:
  $$\fm:=\sum_{\unl{d}} [\CO_{\fQ_{\unl{d}}}]\in M(\pi)^\wedge
    \ \mathrm{and}\
    \fv:=\sum_{\unl{d}} [\fQ_{\unl{d}}]\in V(\pi)^\wedge.$$
This name is motivated by their \emph{eigenvector properties} of
Proposition~\ref{Whittaker}, Remark~\ref{remark on whittaker}(c).

Motivated by the work of Brundan-Kleshchev, see~\cite{bk}, we expect that
the truncated shifted quantum affine algebras $\CU^{N\omega_{n-1}}_\mu$ of $\ssl_n$
should be $\vv$-analogues of the finite W-algebras $W(\ssl_N,e_\pi)$, see~\cite{sev},
where $N:=\sum p_i$ and $e_\pi\in\ssl_N$ is a nilpotent element of Jordan type $\pi$.


\ssec{acknowledge}{Acknowledgments}
\

We are deeply grateful to M.~Bershtein, R.~Bezrukavnikov,
A.~Braverman, A.~Brochier, S.~Cautis, S.~Cherkis, T.~Dimofte,
P.~Etingof, B.~Feigin, S.~Gautam, M.~Gekhtman, V.~Ginzburg,
D.~Hernandez, D.~Jordan, J.~Kamnitzer, S.~Khoroshkin, A.~Marshakov,
A.~Molev, H.~Nakajima, V.~Pestun, L.~Rybnikov, A.~Shapiro, A.~Weekes,
M.~Yakimov, whose generous help and advice was crucial in the
process of our work. We also thank the anonymous referee.

M.F.\ was partially supported by the Russian Academic Excellence Project `5-100'.

A.T.\ gratefully acknowledges support from the
Simons Center for Geometry and Physics, Stony Brook University
and is extremely grateful to IHES for invitation and wonderful working
conditions in the summer 2017, where most of the research for this paper
was performed. A.T.\ also thanks Yale University and Max Planck Institute for
Mathematics in Bonn for the hospitality and support in the summer 2018, where
the final version of this paper was completed.
A.T.\ was partially supported by the NSF Grants DMS--1502497, DMS--1821185.


\sec{RLT}{Relativistic open Toda lattice}


\ssec{quasiham}{Quasihamiltonian reduction}
\

Let $G\supset B\supset T$ be a reductive group with a Borel and Cartan subgroups.
Let $T\subset B_-\subset G$ be the opposite Borel subgroup; let $U$
(resp.~$U_-$) be the unipotent radical of $B$ (resp.~$B_-$).
We consider the {\em double} $D(G)=G\times G$ (see e.g.~\cite[3.2]{amm})
equipped with an action of
$G\times G\colon (u_1,u_2)\cdot(g_1,g_2)=(u_1g_1u_2^{-1},u_2g_2u_2^{-1})$,
and with a moment map
  $\mu=(\mu_1,\mu_2)\colon D(G)\to G\times G,\
   \mu(g_1,g_2)=(g_1g_2g_1^{-1},g_2^{-1})$
(see~\cite[Remark~3.2]{amm}). The double $D(G)$ carries a (non-closed) 2-form
  $\omega_D=\frac{1}{2}(\on{Ad}_{g_2}g_1^*\theta,g_1^*\theta)+
   \frac{1}{2}(g_1^*\theta,g_2^*\theta+g_2^*\ol{\theta})$ where $(\cdot,\cdot)$
is a nondegenerate invariant symmetric bilinear form on $\fg$, and $\theta$
(resp.\ $\ol\theta$) is the left- (resp.\ right-) invariant Maurer-Cartan
form on $G$.

We choose a pair of Coxeter elements $c,c'\in W=N_G(T)/T$, and their
representatives $\dot{c},\dot{c}'\in N_G(T)$. Steinberg's cross-section
$\Sigma^{\dot{c}}_G\subset G$ is defined
as $Z^0(G)\cdot (U_-\dot{c}\cap\dot{c}U)$. If $G$ is semisimple simply-connected,
then the composed morphism
  $\Sigma^{\dot{c}}_G\hookrightarrow G\to G/\!\!/\!\on{Ad}_G=T/W$
is an isomorphism~\cite[Theorem~1.4]{s}. For arbitrary $G$, the composed morphism
$\varrho\colon \Sigma^{\dot{c}}_G\to T/W$ is a ramified Galois cover with Galois group
$\pi_1(G/Z^0(G))$. Furthermore, we consider
$\Xi^{\dot{c}}_G:=Z^0(G)\cdot U_-\dot{c}U_-\supset\Sigma^{\dot{c}}_G$.
According to~\cite[8.9]{s} (for a proof, see e.g.~\cite{hl}),
$\Sigma^{\dot{c}}_G$ meets any $U_-$-orbit (with respect to the conjugation action)
on $\Xi^{\dot{c}}_G$ in exactly one point, and the conjugation action of $U_-$ on
$\Xi^{\dot{c}}_G$ is free, so that $\Xi^{\dot{c}}_G/\!\on{Ad}_{U_-}\simeq\Sigma^{\dot{c}}_G$.

For example, according to~\cite[7.4b)]{s}, for an appropriate choice of $\dot{c}$,
the Steinberg cross-section $\Sigma^{\dot{c}}_{SL(n)}$ consists of the matrices with 1's
just above the main diagonal, $(-1)^{n-1}$ in the bottom left corner, arbitrary
entries elsewhere in the first column, and zeros everywhere else
(in our conventions, $B$ (resp.\ $B_-$) is the subgroup of upper triangular
(resp.\ lower triangular) matrices in $SL(n)$). Hence
$\Xi^{\dot{c}}_{SL(n)}$ consists of matrices with 1's just above the main
diagonal, and zeros everywhere above that.

Following~\cite{fm}, we define the phase space of the open relativistic Toda
lattice as the quasihamiltonian reduction
  $\dZ(G):=\mu^{-1}(\Xi^{\dot{c}'}_G\times\on{inv}(\Xi^{\dot{c}}_G))/U_-\times U_-$
where $\on{inv}\colon G\to G$ is the inversion $g\mapsto g^{-1}$.
The composed projection
  $$\mu^{-1}(\Xi^{\dot{c}'}_G\times\on{inv}(\Xi^{\dot{c}}_G))
    \to\on{inv}(\Xi^{\dot{c}}_G)
    \hookrightarrow G\twoheadrightarrow G/\!\!/\!\on{Ad}_G=T/W$$
gives rise to an integrable system $\varpi\colon \dZ(G)\to T/W$ which factors through
$\dZ(G)\stackrel{\tilde\varpi}{\longrightarrow}\Sigma^{\dot{c}}_G
\stackrel{\varrho}{\longrightarrow}T/W$.

\lem{symp}
If $G$ is semisimple simply-connected, then $\dZ(G)$ is smooth, and $\omega_D$
gives rise to a symplectic form on $\dZ(G)$.
\elem

\prf
The morphism $\Xi^{\dot{c}}_G\to\Sigma^{\dot{c}}_G=T/W$ is smooth by~\cite[Theorem~1.5]{s},
so the fibered product
$\Xi^{\dot{c}'}_G\times_{T/W}\Xi^{\dot{c}}_G\subset\Xi^{\dot{c}'}_G\times\Xi^{\dot{c}}_G$
is smooth. But
  $$\mu\colon D(G)\supset\mu^{-1}(\Xi^{\dot{c}'}_G\times\on{inv}(\Xi^{\dot{c}}_G))
    \to\Xi^{\dot{c}'}_G\times\on{inv}(\Xi^{\dot{c}}_G)\simeq
    \Xi^{\dot{c}'}_G\times\Xi^{\dot{c}}_G$$
is a submersion onto $\Xi^{\dot{c}'}_G\times_{T/W}\Xi^{\dot{c}}_G$, hence
$M:=\mu^{-1}(\Xi^{\dot{c}'}_G\times\on{inv}(\Xi^{\dot{c}}_G))$ is smooth, and its
quotient modulo the free action of $U_-\times U_-$ is smooth as well.

The restriction of $\omega_D$ to $M$ is
$U_-\times U_-$-invariant, so it descends to a 2-form $\omega$ on $\dZ(G)$.
This 2-form is closed since the differential $d\omega_D=-\mu^*(\chi_1+\chi_2)$
(see~\cite[Definition~2.2(B1)]{amm}) where
$\chi=\frac{1}{12}(\theta,[\theta,\theta])$ is the canonical closed biinvariant
3-form on $G$, and $\chi_1$ (resp.\ $\chi_2$) is its pull-back from the first
(resp.\ second) copy of $G$. But the restriction $\chi|_{\Xi^{\dot{c}}_G}$ vanishes
identically since $(\fb_-,[\fb_-,\fb_-])=0$.

It remains to check the nondegeneracy of $\omega$, that is given
$(g_1,g_2)\in M$ to check that $\on{Ker}\omega_D|_M(g_1,g_2)$ is contained in
the span $v(\fn_-\oplus\fn_-)$ of tangent vectors at $(g_1,g_2)$ arising from
the action of $U_-\times U_-$.
The argument in the proof of~\cite[Theorem~5.1]{amm} shows that
$\on{Ker}\omega_D|_M(g_1,g_2)\subset v(\fg\oplus\fg)$. However, it is clear
that $T_{(g_1,g_2)}M\cap v(\fg\oplus\fg)=v(\fn_-\oplus\fn_-)$.

The lemma is proved.
\epr


\ssec{ham}{Poisson reduction}
\

Note that
  $T\cdot\Xi^{\dot{c}}_G=\Xi^{\dot{c}}_G\cdot T=
   \on{Ad}_T(\Xi^{\dot{c}}_G)= B_-\cdot\dot{c}\cdot B_-=:C_c$
(a Coxeter Bruhat cell). One can check that the natural morphism
  $$\dZ(G)=\mu^{-1}(\Xi^{\dot{c}'}_G\times\on{inv}(\Xi^{\dot{c}}_G))/U_-\times U_-\to
    \mu^{-1}(C_{c'}\times\on{inv}(C_c))/B_-\times B_-$$
is an isomorphism. Moreover, the action of $B_-\times B_-$ on
$\mu^{-1}(C_{c'}\times\on{inv}(C_c))$ factors through the free action of
$(B_-\times B_-)/\Delta_{Z(G)}$: the quotient modulo the diagonal copy
of the center of $G$.

The double $D(G)=G\times G$ carries the Semenov-Tian-Shansky Poisson
structure~\cite[Section~2]{sts}. Following {\em loc.~cit.}, $G\times G$ with
this Poisson structure is denoted by $(D_+(G), \{,\!\}_+)$, the
{\em Heisenberg double}. Another Poisson structure on $G\times G$ denoted $\{,\!\}_-$ in
{\em loc.\ cit}.\ is the {\em Drinfeld double} $D_-(G)$. The diagonal embedding
$G\hookrightarrow D_-(G)$ is Poisson with respect to the standard Poisson structure
on $G$ denoted $\pi_G$ in~\cite[2.1]{el}. The dual (Semenov-Tian-Shansky) Poisson
structure on $G$ is denoted $\pi$ in~\cite[2.2]{el}.

The Heisenberg double $D_+(G)$ is equipped with two commuting (left and right) dressing
Poisson actions of the Drinfeld double $D_-(G)$. Restricting to the diagonal
$G\hookrightarrow D_-(G)$ we obtain two commuting Poisson actions of $(G,\pi_G)$ on
$D_+(G)$. The multiplicative moment map of this action is nothing but
$\mu\colon D_+(G)\to (G,\pi)\times(G,\pi)$ of~\refss{quasiham} (a Poisson morphism).
Now $C_c\subset G$ is a coisotropic subvariety~\cite[6.2]{el} of $(G,\pi)$, and
$\mu^{-1}(C_{c'}\times\on{inv}(C_c))\hookrightarrow D(G)$ is a coisotropic
subvariety of ($D_+(G),\ \{,\}_+$).
The action of $G\times G$ on $(D_+(G), \{,\!\}_+)$ is Poisson if $G\times G$
is equipped with the direct product of the standard Poisson-Lie structures
denoted $\pi_G$ in~\cite[2.1]{el}. Note that $B_-\times B_-\subset G\times G$
is a Poisson-Lie subgroup; its Poisson structure will be denoted
$\pi_{B_-}\times\pi_{B_-}$.

The characteristic distribution~\cite[6.2]{el} of the coisotropic subvariety
$\mu^{-1}(C_{c'}\times\on{inv}(C_c))\subset(D_+(G), \{,\!\}_+)$ coincides with
the distribution defined by the tangent spaces to the $B_-\times B_-$-orbits
in $\mu^{-1}(C_{c'}\times\on{inv}(C_c))$. By~\cite[Proposition~6.7]{el} we obtain
a Poisson structure on
$\mu^{-1}(C_{c'}\times\on{inv}(C_c))/(B_-\times B_-)\simeq\dZ(G)$.
This Poisson structure coincides with the one arising from the symplectic form
$\omega$ on $\dZ(G)$.


\ssec{unicentr}{The universal centralizer}
\

Recall that the universal centralizer~\cite[Section 8]{lus1}
$\fZ^G_G\subset G\times\Sigma^{\dot{c}}_G$ is defined as
$\fZ^G_G=\{(g,x) : gxg^{-1}=x\}$. In case $c=c'$ and
$\dot{c}=\dot{c}'$, we have an evident embedding
  $\fZ^G_G\hookrightarrow\mu^{-1}(\Xi^{\dot{c}}_G\times\on{inv}(\Xi^{\dot{c}}_G))$,
and the composed morphism
  $\eta\colon\fZ^G_G\hookrightarrow
  \mu^{-1}(\Xi^{\dot{c}}_G\times\on{inv}(\Xi^{\dot{c}}_G))\twoheadrightarrow\dZc(G)$.
Clearly, the following diagram commutes:
  $$\begin{CD}
    \fZ^G_G @>\eta>> \dZc(G)\\
    @VV{\on{pr}}V @V{\tilde\varpi}VV\\
    \Sigma^{\dot{c}}_G @= \Sigma^{\dot{c}}_G
    \end{CD}$$

\prop{bezr}
For semisimple simply-connected $G$, the morphism
$\eta\colon \fZ^G_G\to\dZc(G)$ is an isomorphism.
\eprop

\prf
First we prove the surjectivity of $\eta$.
We use the equality $U_-\times U_-=(U_-\times\{e\})\times\Delta_{U_-}$. Given
$(g_1,g_2)\in\mu^{-1}(\Xi^{\dot{c}}_G\times\on{inv}(\Xi^{\dot{c}}_G))$ we first act by
$(u_2,u_2)\in\Delta_{U_-}\colon (g_1,g_2)\mapsto(u_2g_1u_2^{-1},u_2g_2u_2^{-1})$.
We can find a unique $u_2$ such that $u_2g_2u_2^{-1}\in\Sigma^{\dot{c}}_G$.
Let us denote the resulting $(u_2g_1u_2^{-1},u_2g_2u_2^{-1})$ by $(h_1,h_2)$
for brevity. Now we act by
the left shift $h_1\mapsto u_1h_1$ which takes $h_1h_2h_1^{-1}$ to
$u_1h_1h_2h_1^{-1}u_1^{-1}$. We can find a unique $u_1$ such that
$u_1h_1h_2h_1^{-1}u_1^{-1}\in\Sigma^{\dot{c}}_G$. Now both $h_2=u_2g_2u_2^{-1}$ and
$u_1h_1h_2h_1^{-1}u_1^{-1}$ are in $\Sigma^{\dot{c}}_G$. Being conjugate they must
coincide, hence $(u_1h_1,h_2)\in\fZ^G_G$.

Now if $\eta(g,x)=\eta(g',x')$, then there is $u_2\in U_-$ such that
$u_2xu_2^{-1}=x'$, hence $x=x'$ and $u_2=e$. Then $g'=u_1g$ for some
$u_1\in U_-$, and both $g$ and $g'$ commute with $x$, hence $u_1xu_1^{-1}=x$,
hence $u_1=e$, so that $g=g'$.

So $\eta$ is bijective at the level of $\BC$-points. But
$\dZc(G)$ is smooth, hence $\eta$ is an isomorphism.
\epr

\rem{uvy}
For arbitrary reductive $G$ the morphism $\eta$ is an affine embedding, but it
fails to be surjective already for $G=PGL(2)$ where the class of
$(g_1,g_2)$ such that $g_2=\left(\begin{array}{cc} a&-1\\ 1&0\end{array}\right)$
and $g_1g_2g_1^{-1}=\left(\begin{array}{cc} -a&-1\\ 1&0\end{array}\right)$
does not lie in the image of $\eta$ when $a\ne0$. Similarly, for $G=GL(2)$,
the class of $(g_1,g_2)$ such that
$g_2=\left(\begin{array}{cc} a&-1\\ 1&0\end{array}\right)$ and
$g_1g_2g_1^{-1}=\left(\begin{array}{cc} a&1\\ -1&0\end{array}\right)$
does not lie in the image of $\eta$. It also follows that the natural
projection $\dZc(SL(2))\to\dZc(PGL(2))$ is {\em not} surjective.
\erem

\rem{compar}
For $G$ semisimple simply-connected, the reduction
$({\mathbf D}(G),\omega_{{\mathbf D}(G)})/\!\!/\!\on{diag}(G)$
\cite[Example~6.1,~Remark~6.2]{amm}
inherits a symplectic structure on its
nonsingular locus. We have a natural morphism
$\fZ^G_G\to({\mathbf D}(G),\omega_{{\mathbf D}(G)})/\!\!/\!\on{diag}(G)$
which is a birational isomorphism
(but not an isomorphism: e.g.\ it contracts the centralizer of a regular
unipotent element). Thus an open subvariety of $\fZ^G_G$ is equipped with a
symplectic form pulled back from
$({\mathbf D}(G),\omega_{{\mathbf D}(G)})/\!\!/\!\on{diag}(G)$.
This form extends
to a symplectic form on the entire $\fZ^G_G$~\cite[2.4]{bfm}.
The isomorphism $\eta\colon\fZ^G_G\iso\dZc(G)$ is a symplectomorphism.
\erem


\ssec{cox}{Comparison with the Coxeter-Toda lattice}
\

We compare $\dZ(G)$ with the construction of~\cite{hkn}. Throughout this section
we assume $G$ to be semisimple simply-connected. The left action of the
center $Z(G)$ on $D(G),\ \xi\cdot(g_1,g_2)=(\xi g_1,g_2)$ gives rise to the
action of $Z(G)$ on $\dZ(G)=M/U_-\times U_-$ where
$M=\mu^{-1}(\Xi^{\dot{c}'}_G\times\on{inv}(\Xi^{\dot{c}}_G))\subset D(G)=G\times G$.
We consider an open subset
$M\supset\buM:=(U_-\cdot T\cdot\dot{w}_0\cdot U_-\times G)\cap M$ given by
the condition that $g_1$ lies in the big Bruhat cell $C_{w_0}\subset G$.
Clearly, $\buM\subset M$ is $U_-\times U_-$-invariant, and we define
$\bdZ(G):=\buM/U_-\times U_-$, an open subvariety of $\dZ(G)$.
Let $S\subset\buM$ be given by the condition
$g_1\in T\cdot\dot{w}_0$. Then the composed projection
$S\hookrightarrow\buM\twoheadrightarrow\bdZ(G)$ is an isomorphism.
Moreover, the projection $\on{pr}_2\colon S\to G$ is a $Z(G)$-torsor over its image
  $\Xi_G^{\dot{c}}\cap\on{Ad}_T(\dot{w}_0\Xi_G^{\dot{c}'}\dot{w}_0^{-1})=
   \Xi_G^{\dot{c}}\cap\on{Ad}_T(U\dot{w}_0\dot{c}'\dot{w}_0^{-1}U)$.
Finally, note that the composed projection
  $$\Xi_G^{\dot{c}}\cap\on{Ad}_T(U\dot{w}_0\dot{c}'\dot{w}_0^{-1}U)\hookrightarrow
    T\cdot U_-\cdot\dot{c}\cdot U_-\cdot T\cap T\cdot
    U\cdot\dot{w}_0\dot{c}'\dot{w}_0^{-1}\cdot U\cdot T\twoheadrightarrow$$
  $$\twoheadrightarrow(T\cdot U_-\cdot\dot{c}\cdot U_-\cdot T\cap
    T\cdot U\cdot\dot{w}_0\dot{c}'\dot{w}_0^{-1}\cdot U\cdot T)/\!\on{Ad}_T
    =:G^{\dot{c},\dot{w}_0\dot{c}'\dot{w}_0^{-1}}\!/\!\on{Ad}_T$$
is an isomorphism. But according to~\cite{hkn} (see also~\cite{gsv}),
$G^{\dot{c},\dot{w}_0\dot{c}'\dot{w}_0^{-1}}\!/\!\on{Ad}_T$ is the phase space of the
Coxeter-Toda lattice.
All in all, we obtain an isomorphism (respecting the symplectic structures)
$Z(G)\backslash\bdZ\iso G^{\dot{c},\dot{w}_0\dot{c}'\dot{w}_0^{-1}}/\!\on{Ad}_T$.

For example, for an appropriate choice of $\dot{c},\dot{c}'\in SL(n)$, the
slice $S$ is formed by all the tridiagonal matrices of determinant 1 with 1's
just above the main diagonal, and with the invertible entries just below the
main diagonal (see~\cite[Introduction]{gsv}).

We also define an open subset
  $S\supset\oS:=\{(g_1,g_2)\in M : g_1\in T\cdot\dot{w}_0,\ g_2\in U_-\cdot T\cdot U\}$.
It is equipped with a projection $\on{pr}_1\colon \oS\to T\cdot\dot{w}_0\iso T$, and
with another projection $\on{pr}_2\colon \oS\to U_-\cdot T\cdot U\twoheadrightarrow T$.
One can check that $(\on{pr}_1,\on{pr}_2)\colon \oS\iso T\times T$.
We define an open subvariety $\dZ(G)\supset\bdZ(G)\supset\odZ(G)$ as the
isomorphic image of $\oS$. Thus $\odZ(G)\simeq T\times T$.


\ssec{SL2}{Trigonometric zastava for $SL(2)$}
\

Recall the degree $n$ trigonometric open zastava $\oZ^n$ for the group $SL(2)$
(see~\cite{fkr}). This is the moduli space of pairs of relatively
prime polynomials
  $(Q=z^n+q_1z^{n-1}+\ldots+q_n,\ R=r_1z^{n-1}+r_2z^{n-2}+\ldots+r_n)$
such that $q_n\ne0$.
We have a morphism $\zeta\colon \fZ^{GL(n)}_{GL(n)}\to{}\oZ^n$
taking a pair $(g,x)\in\fZ^{GL(n)}_{GL(n)}$ to $(Q,R)$ where $Q$ is
the characteristic polynomial of $x$, and $R$ is a unique polynomial
of degree less than $n$ such that $R(x)=g$.
We denote by $\on{pr}\colon \oZ^n\to(\BC^\times)^{(n)}$ the morphism taking
$(Q,R)$ to the set of roots of $Q$.

Recall that $\Sigma^{\dot{c}}_{GL(n)}=Z^0(GL(n))\cdot\Sigma^{\dot{c}}_{SL(n)}\simeq
Z(GL(n))\times\Sigma^{\dot{c}}_{SL(n)}=\BC^\times\times\Sigma^{\dot{c}}_{SL(n)}$.
We denote by $\sfp\colon\fZ^{GL(n)}_{GL(n)}\to\BC^\times$ the composed projection
$\fZ^{GL(n)}_{GL(n)}\to\Sigma^{\dot{c}}_{GL(n)}\to\BC^\times$.

\prop{ryb}
The following square is Cartesian:
  $$\begin{CD}
    \fZ^{GL(n)}_{GL(n)} @>\zeta>> \oZ^n\\
    @VV{\sfp}V @V{q_n}VV\\
    \BC^\times @>n>> \BC^\times
    \end{CD}$$
Thus $\fZ^{GL(n)}_{GL(n)}$ is an unramified $\BZ/n\BZ$-cover of $\oZ^n$.
\eprop

\prf
Clear from the above discussion.
\epr

Following~\cite[end of chapter~2]{ah}, we consider the subvariety
$\tZ_1^n\hookrightarrow\oZ^n$ formed by the pairs $(Q,R)$ such that
$q_n=1$ and the resultant of $Q$ and $R$, denoted $\mathrm{Result(Q,R)}$, equals 1. Note that we have an
evident embedding $\fZ^{SL(n)}_{SL(n)}\hookrightarrow\fZ^{GL(n)}_{GL(n)}$.

\cor{consol}
The restriction of the morphism $\zeta$ to $\fZ^{SL(n)}_{SL(n)}\subset\fZ^{GL(n)}_{GL(n)}$
gives rise to an isomorphism $\zeta\colon \fZ^{SL(n)}_{SL(n)}\iso\tZ^n_1$.
\ecor

\begin{proof}
For $(g,x)\in \fZ^{GL(n)}_{GL(n)}$, the inclusion $x\in SL(n)$ is equivalent to $q_n=1$, while we claim that the inclusion
$g\in SL(n)$ is equivalent to $\mathrm{Result}(Q,R)=1$. The latter follows by combining the equalities $gx=xg$ and $g=R(x)$
with the standard equality $\mathrm{Result}(Q,R)=\prod_{i=1}^n R(\xi_i)$, where $\{\xi_i\}_{i=1}^n$ are the roots of $Q$.
Since $\{\xi_i\}_{i=1}^n$ are the generalized eigenvalues (taken with corresponding multiplicities) of $x$, it is easy to see that
$\{R(\xi_i)\}_{i=1}^n$ are the generalized eigenvalues of $g$, hence, $\mathrm{det}(g)=\prod_{i=1}^n R(\xi_i)$.
\end{proof}

For a future use we define an unramified $\BZ/2\BZ$-cover
$\hZ^n\to\oZ^n$ where $\hZ^n$ is the moduli space of pairs of relatively prime
polynomials $(Q=q_0z^n+q_1z^{n-1}+\ldots+q_n,\ R=r_1z^{n-1}+r_2z^{n-2}+\ldots+r_n)$
such that $q_n\cdot q_0=(-1)^n$. The projection $\hZ^n\to\oZ^n$ takes
$(Q,R)$ to $(q_0^{-1}Q,R)$.

Finally, there are important embeddings $\Psi\colon \oZ^n, \hZ^n\hookrightarrow SL(2,\BC[z])$
taking $(Q,R)$ to a unique matrix
  $\left(\begin{array}{cc} Q&\tilde{R}\\ R&\tilde{Q}\end{array}\right)$
such that $\deg\tilde{R}\leq n>\deg\tilde{Q}$, and $\tilde{R}(0)=0$, that is
$\tilde{R}=\tilde{r}_0z^n+\tilde{r}_1z^{n-1}+\ldots+\tilde{r}_{n-1}z$.
Identifying $\hZ^n$ and $\oZ^n$ with their images inside
$SL(2,\BC[z])$, the matrix multiplication gives rise to the multiplication
morphisms $\hZ^k\times\hZ^l\to\hZ^{k+l}, \oZ^k\times\oZ^l\to\oZ^{k+l}$.


\sec{Nildaha}{Quantum relativistic open Toda and nil-DAHA}
\label{Nildaha}

Throughout this section (with the exception of~\refss{gln} dealing with $G=GL(n)$)
$G$ is an almost simple simply-connected complex algebraic group.


\ssec{fold}{Root systems and foldings}
\

Let $G^\vee$ be the Langlands dual (adjoint) group with a Cartan torus $T^\vee$.
We choose a Borel subgroup $B^\vee\supset T^\vee$. It defines the set of simple
positive roots $\{\alpha_i,\ i\in I\}$. Let $\fg^\vee$ be the Lie algebra of $G^\vee$.
We realize $\fg^\vee$ as a {\em folding} of a simple simply-laced Lie algebra
$\fg'{}^\vee$, i.e.\ as invariants of an outer automorphism $\sigma$ of
$\fg'{}^\vee$ preserving a Cartan subalgebra $\ft'{}^\vee\subset\fg'{}^\vee$
and acting on the root system of $(\fg'{}^\vee,\ft'{}^\vee)$.
In particular, $\sigma$ gives rise to the same named automorphism of the Langlands dual
Lie algebras $\fg'\supset\ft'$ (note that say, if $\fg$ is of type
$B_n$, then $\fg'$ is of type $A_{2n-1}$, while if $\fg$ is of type $C_n$,
then $\fg'$ is of type $D_{n+1}$; in particular, $\fg\not\subset\fg'$).
We choose a $\sigma$-invariant Borel
subalgebra $\ft'\subset\fb'\subset\fg'$ such that $\fb=(\fb')^\sigma$.
The corresponding set of simple roots is denoted by $I'$.
We denote by $\Xi$ the finite cyclic group generated by $\sigma$, and
$d:=|\Xi|$. Let $G'\supset T'$ denote the corresponding simply-connected
Lie group and its Cartan torus.
The {\em coinvariants} $X_*(T')_\sigma$ of $\sigma$ on the coroot
lattice $X_*(T')$ of $(\fg',\ft')$ coincide with the root lattice
of $\fg^\vee$. We have an injective map $a:\ X_*(T')_\sigma\to X_*(T')^\sigma$
from coinvariants to invariants defined as follows: given a coinvariant
$\alpha$ with a representative $\tilde\alpha\in X_*(T')$ we set
$a(\alpha):=\sum_{\xi\in\Xi}\xi(\tilde\alpha)$.

To compare with the notations of~\cite[4.4, Remark~4.5]{h}, we are in the symmetric case
with $Q'_0=Y:=X^*(T^\vee)=X_*(T)=X_*(T')_\sigma$, and
$Q_0\subset X:=X^*(T')_\sigma$ generated by the classes of
simple roots of $T'\subset B'\subset G'$. Note that $Q'_0$ is generated
by the classes of simple coroots of $T'\subset B'\subset G'$, and we have a
canonical identification $Q_0=Q'_0$ sending a coroot $\tilde\alpha$
to the corresponding root $\tilde\alpha{}^{\!\scriptscriptstyle\vee}$.
The Weyl group $W$ of $G\supset T$ coincides with the invariants
$(W')^\sigma$ of $\sigma$ on the Weyl group $W'$ of $G'\supset T'$
(our $W$ is denoted $W_0$ in~\cite{h}).
The $W$-invariant pairing $X\times Y\to\BQ$ defined in~\cite[4.4]{h} is actually
integer valued: $X\times Y\to\BZ$, so that $m=1$ (notations of {\em loc.~cit.}).
To compare with notations of~\cite[Section~1]{cf}, $P:=X,\ Q:=Q_0$,
and the natural pairing $P\times P\to\BQ$ gives rise to the embedding
$Q=Y\hookrightarrow P$. We will also need an extended lattice
$Y_{\on{ad}}:=X_*(T_{\on{ad}})=X_*(T'_{\on{ad}})_\sigma\supset Y$. Note that
$\Pi:=Y_{\on{ad}}/Y=(X_*(T'_{\on{ad}})/X_*(T'))_\sigma$.
Also note that the above $W$-invariant identification $Q_0=Q'_0$ extends
to the $W$-invariant identification $Q_0\subset X=Y_{\on{ad}}\supset Q'_0$.
The extended pairing $X\times Y_{\on{ad}}\to\BQ$ is no more integer valued in general,
and we denote by $m_{\on{ad}}$ the maximal denominator appearing in the values of
this pairing. Finally, $R\subset X$ stands for the set of roots.


\ssec{flag}{Affine flags}
\

We fix a primitive root of unity $\zeta$ of order $d=\on{ord}(\sigma)$.
We set $\CK=\BC((\st))\supset\CO=\BC[[\st]]$.
The group ind-scheme $G'(\CK)$ is equipped with an automorphism $\varsigma$ defined as the
composition of two automorphisms: a) $\sigma$ on $G'$; b) $\st\mapsto\zeta\st$.
This automorphism preserves the Iwahori subgroup $\bI'\subset G'(\CK)$.
We denote by $\Fl$ the twisted affine flag space
$G'(\CK)^\varsigma/(\bI')^\varsigma$: an ind-proper ind-scheme of ind-finite type,
see~\cite{pr}. We denote by $\fu\subset\on{Lie}(\bI')^\varsigma$
its pronilpotent radical.
The trivial (Tate) bundle $\unl{\fg'(\CK)^\varsigma}$ with the fiber
$\fg'(\CK)^\varsigma$ over $\Fl$ has a structure of an ind-scheme. It contains
a profinite dimensional vector subbundle $\unl\fu$ whose fiber over a point
$b\in\Fl$ represented by a compact subalgebra in $\fg'(\CK)^\varsigma$ is the
pronilpotent radical of this subalgebra. The trivial vector bundle
$\unl{\fg'(\CK)^\varsigma}$ also contains a trivial vector subbundle
$\fu\times\Fl$. We will call $\unl\fu$ the cotangent bundle of $\Fl$,
and we will call the intersection $\bLambda:=\unl\fu\cap(\fu\times\Fl)$ the
affine Steinberg variety.

To simplify the notations we will write $\bI$ for $(\bI')^\varsigma$,
and $\bK$ for $G'(\CO)^\varsigma$. The convolution product on the
complexified equivariant coherent $K$-theory
$K^{\BC^\times\times\bI\rtimes\BC^\times}(\bLambda)$ is defined
as in~\cite[Remark~3.9(3)]{bfna} (cf.~\cite[7.1]{bfm} and~\cite[2.2,~2.3]{vv}).
Here the first copy of $\BC^\times$ acts by dilations in
fibers of $\unl\fu$, while the second one acts by loop rotations,
and $K_{\BC^\times\times\BC^\times}(\on{pt})=\BC[t^{\pm1},q^{\pm1}]$.


\ssec{daha}{DAHA, symmetric case}
\

Following~\cite[4.6]{h}, we set
$\widetilde{X}:=X\oplus\BZ\delta=X^*(T')_\sigma\oplus\BZ\delta$. This is the
group of characters of $\bI\rtimes\BC^\times$. Note that the Picard group
$\on{Pic}(\Fl)$ is canonically isomorphic to $X\oplus\BZ\omega_0$.
The $\bI$-orbits on $\Fl$ are parametrized by the affine
Weyl group $W_a\simeq Y\rtimes W=X_*(T')_\sigma\rtimes W$.
We denote by $\bLambda_e\simeq\fu$ the closed subscheme of $\bLambda$:
the preimage of the one-point $\bI$-orbit $\Fl_e\subset\Fl$.
For $\tilde\lambda=(\check\lambda,k)\in\widetilde{X}$ we denote by
  $\CO_{\bLambda_e}\langle\tilde\lambda\rangle\in K^{\BC^\times\times\bI\rtimes\BC^\times}(\bLambda)$
the (class of the) direct image of the structure sheaf of $\bLambda_e$
twisted by the character $\tilde\lambda$ of $\bI\rtimes\BC^\times$.
Let $\widetilde{I}\subset W_a$ be the set of one-dimensional $\bI$-orbits
on $\Fl$. For $i\in\widetilde{I}$ we denote by $\Fl_i$ the corresponding
orbit, and by $\ol\Fl_i$ its closure, isomorphic to a projective line.
We denote by $\bLambda_i\subset\bLambda$ the closed subscheme of $\bLambda$:
the closure of the preimage of $\Fl_i$. We denote by $\omega_{\bLambda_i}$
the (class of the) direct image (wrt the closed embedding
$\bLambda_i\hookrightarrow\bLambda$) of the inverse image (wrt the smooth
projection $\bLambda_i\to\ol\Fl_i$) of the canonical line bundle on
$\ol\Fl_i\simeq\BP^1$ equipped with the natural
$\BC^\times\times\bI\rtimes\BC^\times$-equivariant structure. Finally, we set
$\sT_i:=-1-t\omega_{\bLambda_i}\in K^{\BC^\times\times\bI\rtimes\BC^\times}(\bLambda)$.

\defe{dacha}
(Cf.~\cite[Definition 5.6]{h}.)
The {\em double affine Hecke algebra (DAHA)} $\CH(W_a,\widetilde{X})$ is
the $\BC[q^{\pm1},t^{\pm1}]$-algebra generated by
$\{\sX_{\tilde\lambda},T_w|\tilde\lambda\in\widetilde{X}, w\in W_a\}$
with the following defining relations:

(a) $T_w$'s satisfy the braid relations of $W_a$;

(b) $\sX_{\tilde\lambda}\sX_{\tilde\mu}=\sX_{\tilde\lambda+\tilde\mu}$, and $\sX_\delta=q$;

(c) $(T_i-t)(T_i+1)=0$ for $i\in\widetilde{I}$, where we set $T_i=T_{s_i}$;

(d) $\sX_{\tilde\lambda}T_i-T_i\sX_{\tilde\lambda-r\check\alpha_i}=(t-1)\sX_{\tilde\lambda}
(1+\sX_{-\check\alpha_i}+\ldots+\sX_{-\check\alpha_i}^{r-1})$ where
$\langle\tilde\lambda,\alpha_i\rangle=r\geq0$.
\edefe

\th{vas}
There is a unique isomorphism
  $\Phi\colon \CH(W_a,\widetilde{X})\iso K^{\BC^\times\times\bI\rtimes\BC^\times}(\bLambda)$
such that
$\Phi(\sX_{\tilde\lambda})=\CO_{\bLambda_e}\langle\tilde\lambda\rangle$, and
$\Phi(T_i)=\sT_i$, for any $i\in\widetilde{I}$.
\eth

\prf Same as the one of~\cite[Theorem~2.5.6]{vv}. \epr


\ssec{nildaha}{Nil-DAHA, symmetric case}
\

The complexified equivariant $K$-theory $K^{\bI\rtimes\BC^\times}(\Fl)$ forms a
$\BC[q^{\pm1}]$-algebra with respect to the convolution. We denote
by $\CO_{\Fl_e}\langle\tilde\lambda\rangle$ the (class of the)
structure sheaf of the point orbit $\Fl_e\in\Fl$ twisted by a
character $\tilde\lambda\in\widetilde{X}$. We denote by
$\omega_{\ol\Fl_i}$ the (class of the) direct image (wrt the closed
embedding $\ol\Fl_i\hookrightarrow\Fl$) of the canonical line bundle
on $\ol\Fl_i$ equipped with the natural
$\bI\rtimes\BC^\times$-equivariant structure. We set
$\bT_i:=-1-\omega_{\ol\Fl_i}\in K^{\bI\rtimes\BC^\times}(\Fl)$.

\defe{nildacha}
(Cf.~\cite[1.1]{cf}.)
The {\em nil-DAHA} $\CHH(W_a,\widetilde{X})$ is the $\BC[q^{\pm1}]$-algebra
generated by $\{\sX_{\tilde\lambda},\CT_w|\tilde\lambda\in\widetilde{X},w\in W_a\}$
with the following defining relations:

(a) $\CT_w$'s satisfy the braid relations of $W_a$;

(b) $\sX_{\tilde\lambda}\sX_{\tilde\mu}=\sX_{\tilde\lambda+\tilde\mu}$, and $\sX_\delta=q$;

(c) $\CT_i(\CT_i+1)=0$ for $i\in\widetilde{I}$, where we set $\CT_i=\CT_{s_i}$;

(d) $\sX_{\tilde\lambda}\CT_i-\CT_i\sX_{\tilde\lambda-r\check\alpha_i}=-\sX_{\tilde\lambda}
(1+\sX_{-\check\alpha_i}+\ldots+\sX_{-\check\alpha_i}^{r-1})$ where
$\langle\tilde\lambda,\alpha_i\rangle=r\geq0$.
\edefe

\th{vass}
There is a unique isomorphism
  $\Phi\colon \CH\!\!\CH(W_a,\widetilde{X})\iso K^{\bI\rtimes\BC^\times}(\Fl)$
such that\\ $\Phi(\sX_{\tilde\lambda})=\CO_{\Fl_e}\langle\tilde\lambda\rangle$, and
$\Phi(\CT_i)=\bT_i$, for any $i\in\widetilde{I}$.
\eth

\prf Same as the one of~\cite[Theorem~2.5.6]{vv}. \epr


\ssec{exten}{Extended nil-DAHA}
\

We consider the $2m_{\on{ad}}$-fold cover $\wt\BC{}^\times\to\BC^\times$ of the loop rotation
group (see the end of~\refss{fold}), and set $\widehat\bI:=\bI\rtimes\wt\BC{}^\times$.
The group of characters of $T\times\wt\BC{}^\times$ is
$\widehat{X}:=X\oplus\BZ\frac{1}{2m_{\on{ad}}}\delta$.
The extended affine Weyl group is $W_e=Y_{\on{ad}}\rtimes W=W_a\rtimes\Pi$.
The extended nil-DAHA $\CHH(W_e,\widehat{X})$ is the (extended) semidirect product
$(\CHH(W_a,\widehat{X})\rtimes\Pi)\otimes_{\BC[q^{\pm 1}]}\BC[q^{\frac{\pm1}{2m_{\on{ad}}}}]$.
That is, it has generators
$\sX_{\hat\lambda},\ \hat\lambda\in\widehat{X}$, and $\CT_i,\ i\in\widetilde{I}$,
and $\pi\in\Pi$; with additional relations $\pi\CT_i\pi^{-1}=\CT_{\pi(i)}$,
and $\pi\sX_{\hat\lambda}\pi^{-1}=\sX_{\pi(\hat\lambda)}$.

\rem{cher}
The definition of~\cite[1.1]{cf} is equivalent to our~\refss{exten}:
the generators $T_i$ of {\em loc.~cit.} correspond to $-\CT_i-1$; geometrically,
$T_i=[\omega_{\ol\Fl_i}]$.
\erem


\ssec{residue}{Residue construction}
\

Let $\CA:=\BC[q^{\frac{\pm1}{2m_{\on{ad}}}}]$, and $\CQ:=\BC(q^{\frac{1}{2m_{\on{ad}}}})$.
Let $\CO_q(T\times T)$ be an $\CA$-algebra
with generators $[\lambda,\mu],\ \lambda,\mu\in X$, and relations
  $[\lambda,\mu]\cdot[\lambda',\mu']=
   q^{\frac{(\mu,\lambda')-(\mu',\lambda)}{2}}[\lambda+\lambda',\mu+\mu']$.
This is the subalgebra of endomorphisms of $\CA[T]$ generated by multiplications
by $\sX_\lambda,\ \lambda\in X$, and $q$-shift operators
$D^\mu_qf(t):=f(q^\mu t)$ where we view $q^\mu$ as a homomorphism $\wt\BC{}^\times\to T$.
In other words, $D^\mu_q\sX_\lambda=q^{(\mu,\lambda)}\sX_\lambda$.
We may and will view $\CO_q(T\times T)$ as a subalgebra of endomorphisms of
the field of rational functions $\CQ(T)$ as well. It embeds into the subalgebra
$\BC_q(T\times T)\subset\on{End}(\CQ(T))$ generated by $D^\mu_q,\ \mu\in X$,
and multiplications by $f\in\CQ(T)$. We consider the semidirect
product $\BC_q(T\times T)\rtimes\BC[W]$ with respect to the diagonal action
of $W$ on $T\times T$. Inside we consider the linear subspace
$\CHH_{\on{res}}(W_e,\widehat{X})$ formed by the finite sums
$\sum_{w\in W}^{\mu\in X}h_{w,\mu}D^\mu_q\cdot[w],\ h_{w,\mu}\in\CQ(T)$, satisfying the
following conditions:

\textup{(a)} $h_{w,\mu}$ is regular except at the divisors
$T_{\alpha,q^k}:=\{t : \alpha(t)=q^k\},\ \alpha\in R,\ k\in\BZ$, where they
are allowed to have only first order poles.

\textup{(b)} $\on{Res}_{T_{\alpha,q^{-k}}}(h_{w,\mu})+
\on{Res}_{T_{\alpha,q^{-k}}}(h_{s_\alpha w,k\alpha+s_\alpha\mu})=0$ for any $\alpha\in R$.

The algebra of regular functions $\BC[T\times\wt\BC{}^\times]$ is embedded into
$\CHH_{\on{res}}(W_e,\widehat{X})$ via the assignment $f\mapsto f\cdot[1]$.
Furthermore, for $i\in I\subset\widetilde{I}$, we consider the
{\em Demazure operator}~\cite[1.3]{cf}
  $\tau_i:=\frac{1}{1-\sX_{\alpha_i}}\cdot ([s_i]-[1])\in\CHH_{\on{res}}(W_e,\widehat{X})$,
and for $i_0\in\widetilde{I}\setminus I$ we consider the
{\em Demazure operator}~\cite[1.3]{cf}
  $\tau_{i_0}:=\frac{1}{1-q\sX_\theta^{-1}}\cdot([s_\theta]\cdot D^\theta_q-[1])
   \in\CHH_{\on{res}}(W_e,\widehat{X})$,
where $\theta\in R$ is the dominant short root, $(\theta,\theta)=2$.

\th{gin}
\textup{(a)} $\CHH_{\on{res}}(W_e,\widehat{X})$ is a subalgebra of
$\BC_q(T\times T)\rtimes\BC[W]$.

\noindent
\textup{(b)} The assignment $f\mapsto f\cdot[1];\ \CT_i\mapsto\tau_i,\
i\in\widetilde{I};\ \Pi\ni\pi\mapsto$ the corresponding automorphism
of $\CQ(T)=\BQ(\widehat{X}\otimes\BC^\times)$ (arising from the automorphism of
the extended Dynkin diagram), defines an isomorphism
$\varphi\colon \CHH(W_e,\widehat{X})\iso\CHH_{\on{res}}(W_e,\widehat{X})$.
\eth

\prf
Same as the one of~\cite[Theorem~7.2]{beg}.
\epr

\rem{ginb}
Nil-DAHA $\CHH(W_e,\widehat{X})$ is {\em not} isomorphic to
the degeneration $\ddot{\mathsf H}|_{v=0}$ of~\cite[Section~6]{beg}.
\erem


\ssec{dicson}{$K$-theory of disconnected flags}
\

We define $\bI_{\on{ad}}$ as the image of $\bI$ in $G'_{\on{ad}}(\CK)^\varsigma$,
and we consider the adjoint version of the affine flags
$\Fl_{\on{ad}}:=G'_{\on{ad}}(\CK)^\varsigma/\bI_{\on{ad}}$.
This is an ind-scheme having $|\Pi|$ connected components, each one isomorphic to $\Fl$.
The isomorphism of~\reft{vass} extends to the same named isomorphism
$\CHH(W_e,\widehat{X})\iso K^{\widehat\bI}(\Fl_{\on{ad}})$. Let us explain why the
RHS forms an algebra. We consider an algebra
  $K(\widehat\bI\backslash G'_{\on{ad}}(\CK)^\varsigma/\widehat\bI)=
   K^{\widehat\bI}(\Fl_{\on{ad}}/\Pi)$. Here we view $\Pi=Z(G'{}^\sigma)$
as the center of the simply-connected group $G'{}^\sigma$ acting trivially on $\Fl_{\on{ad}}$.
Now $K^{\widehat\bI}(\Fl_{\on{ad}}/\Pi)$ contains a subalgebra
$K^{\widehat\bI}(\Fl_{\on{ad}}/\Pi)_{\on{diag}}$ formed by the classes of
bi-equivariant sheaves on $\Fl_{\on{ad}}$ such that the $\Pi$-equivariance
coincides with the $Z(G'{}^\sigma)$-equivariance obtained by the restriction
of $\widehat\bI$-equivariance. Finally,
  $K^{\widehat\bI}(\Fl_{\on{ad}}/\Pi)_{\on{diag}}\simeq K^{\widehat\bI}(\Fl_{\on{ad}})$.


\ssec{snd}{Spherical nil-DAHA}
\

We define the new generators $\hat\CT_i:=-\CT_i-1,\ i\in\widetilde{I}$
(they correspond to the generators $T_i$ of~\cite[Definition~1.1]{cf}).
Geometrically, $\hat\CT_i=[\omega_{\ol\Fl_i}]$. They still satisfy
the braid relations of $W_a$. So for any $w\in W_a$ we have a well-defined
element (product of the generators) $\hat\CT_w$.
We also define $\hat\CT{}'_i:=\hat\CT_i+1=-\CT_i,\ i\in\widetilde{I}$. Geometrically,
for $i\in I\subset\widetilde{I}$, we have
$\hat\CT{}'_i=\sX_\rhovee[\CO_{\ol\Fl_i}]\sX_\rhovee^{-1}$.
These generators also satisfy
the braid relations of $W_a$, so for any $w\in W_a$ we have a well-defined
element (product of the generators) $\hat\CT{}'_w$.

Given a reduced decomposition $w=s_{i_1}\cdots s_{i_l}$ we have for the class
of the structure sheaf of the Schubert variety
$[\CO_{\ol\Fl_w}]=[\CO_{\ol\Fl_{i_1}}]\cdots[\CO_{\ol\Fl_{i_l}}]$ since
$\ol\Fl_w$ has rational singularities. Hence, for $w\in W\subset W_a$, we have
$[\CO_{\ol\Fl_w}]=\sX_\rhovee^{-1}\hat\CT{}'_w\sX_\rhovee$.
In particular, for the longest element $w_0\in W$ we set
$\be:=[\CO_{\ol\Fl_{w_0}}]=\sX_\rhovee^{-1}\hat\CT{}'_{w_0}\sX_\rhovee$, an
idempotent in $\CHH(W_e,\widehat{X})$. Indeed, calculating
$[\CO_{\ol\Fl_{w_0}}][\CO_{\ol\Fl_{w_0}}]$ as the pushforward of the
structure sheaf from the convolution diagram
$\ol\Fl_{w_0}\widetilde\times\ol\Fl_{w_0}\to\ol\Fl_{w_0}$ we get
$\CO_{\ol\Fl_{w_0}}$ since $R\Gamma(\ol\Fl_{w_0},\CO_{\ol\Fl_{w_0}})=\BC$.

We define the spherical nil-DAHA $\CHH^{\on{sph}}(W_a,\widetilde{X}):=
\be\CHH(W_a,\widetilde{X})\be$, and the spherical extended nil-DAHA
$\CHH^{\on{sph}}(W_e,\widehat{X}):=\be\CHH(W_e,\widehat{X})\be$.


\ssec{affgr}{Equivariant $K$-theory of the affine Grassmannian}
\

We denote by $\Gr_{\on{ad}}$ the twisted affine Grassmannian
$G'_{\on{ad}}(\CK)^\varsigma/G'_{\on{ad}}(\CO)^\varsigma$:
an ind-proper ind-scheme of ind-finite type,
see~\cite{pr}. The complexified equivariant coherent $K$-theory
  $K^{\bK\rtimes\wt\BC{}^\times}(\Gr_{\on{ad}})=
   K^{G'(\CO)^\varsigma\rtimes\wt\BC{}^\times}(\Gr_{\on{ad}})$
forms a $\BC[q^{\pm\frac{1}{2m_{\on{ad}}}}]$-algebra
with respect to the convolution (see~\refss{dicson}).
We have the smooth projection
$p\colon \Fl_{\on{ad}}\to\Gr_{\on{ad}}$, and the natural embedding
  $K^{\bK\rtimes\wt\BC{}^\times}(\Gr_{\on{ad}})\hookrightarrow
   K^{\bI\rtimes\wt\BC{}^\times}(\Gr_{\on{ad}})
   \stackrel{p^*}{\hookrightarrow}K^{\bI\rtimes\wt\BC{}^\times}(\Fl_{\on{ad}})$.

\cor{vasse}
The isomorphism $\Phi$ of~\refss{dicson} takes the spherical subalgebra\\
$\CHH^{\on{sph}}(W_e,\widehat{X})\subset\CHH(W_e,\widehat{X})$
isomorphically onto
  $K^{\bK\rtimes\wt\BC{}^\times}(\Gr_{\on{ad}})\subset
   K^{\bI\rtimes\wt\BC{}^\times}(\Fl_{\on{ad}})$.
The right ideal $\be\CHH(W_e,\widehat{X})$ corresponds to
  $K^{\bK\rtimes\BC^\times}(\Fl_{\on{ad}})=
   (K^{\bI\rtimes\wt\BC{}^\times}(\Fl_{\on{ad}}))^W\subset
   K^{\bI\rtimes\wt\BC{}^\times}(\Fl_{\on{ad}})$.
\qed
\ecor


\ssec{class}{Classical limit}
\

The following theorem is proved similarly to~\cite[Theorem~2.15]{bfm}:

\th{twist}
\textup{(a)} The algebra $K^\bK(\Gr_{\on{ad}})$ is commutative.

\noindent
\textup{(b)} Its spectrum together with the projection onto $T/W$ is
naturally isomorphic to $\fZ^G_G\stackrel{\on{pr}}{\longrightarrow}T/W$.

\noindent
\textup{(c)} The Poisson structure on $K^\bK(\Gr_{\on{ad}})$ arising from
the deformation $K^{\bK\rtimes\BC^\times}(\Gr_{\on{ad}})$ corresponds under the
above identification to the Poisson (symplectic) structure of~\refr{compar}
on $\fZ^G_G$.
\qed
\eth

\cor{shout}
\textup{(a)} The algebra $\CHH^{\on{sph}}(W_e,\widehat{X})|_{q=1}$ is commutative.

\noindent
\textup{(b)} This algebra with the subalgebra $\BC[X]^W$ is
naturally isomorphic to $\BC[\fZ^G_G]\supset\BC[T/W]$.

\noindent
\textup{(c)} The Poisson structure on $\CHH^{\on{sph}}(W_e,\widehat{X})|_{q=1}$
arising from the deformation $\CHH^{\on{sph}}(W_e,\widehat{X})$
corresponds under the above
identification to the Poisson (symplectic) structure of~\refr{compar}
on $\fZ^G_G$.
\qed
\ecor


\ssec{gln}{Nil-DAHA, general linear group}
\

In case $G=GL(n)\simeq G^\vee$, the general definition of $\CHH(W_e,\widehat{X})$
takes a particularly explicit form.

\defe{GLn}
The nil-DAHA $\CHH(GL(n))$ is the $\BC[q^{\pm1}]$-algebra with generators\\
$\CT_0,\ldots,\CT_{n-1},\ \sX_1^{\pm1},\ldots,\sX_n^{\pm1},\ \pi^{\pm1}$, and
the following relations:

\textup{(a)} $\CT_i$'s for $i\in\BZ/n\BZ$ satisfy the braid relations of the
affine braid group of type $\tilde{A}_{n-1}$;

\textup{(b)} $\sX_i^{\pm1},\ i=1,\ldots,n$, all commute;

\textup{(c)} $\CT_i(\CT_i+1)=0$ for $i\in\BZ/n\BZ$;

\textup{(d)} $\pi\sX_i\pi^{-1}=\sX_{i+1}$ for $i=1,\ldots,n-1$, and
$\pi\sX_n\pi^{-1}=q\sX_1$;

\textup{(e)} $\pi\CT_i\pi^{-1}=\CT_{i+1}$ for $i\in\BZ/n\BZ$;

\textup{(f)} $\sX_{i+1}\CT_i-\CT_i\sX_i=\sX_i$, and
$\sX_i^{-1}\CT_i-\CT_i\sX_{i+1}^{-1}=\sX_{i+1}^{-1}$ for $i=1,\ldots,n-1$;

\textup{(h)} $q\sX_1\CT_0-\CT_0\sX_n=\sX_n$, and
$q\sX_n^{-1}\CT_0-\CT_0\sX_1^{-1}=\sX_1^{-1}$;

\textup{(fh)} $\sX_i^{\pm1}$ and $\CT_j$ commute for all the pairs $i,j$
not listed in \textup{(f,h)} above.
\edefe

Note that $\sX:=\sX_1\cdots\sX_n$ commutes with all the $\CT_i$'s, while
$\pi\sX\pi^{-1}=q\sX$. For a future use we give the following

\defe{Gln'}
The extended nil-DAHA $\CHH_e(GL(n))$ is the $\BC[\vv^{\pm1}]$-algebra, $q=\vv^2$,
with generators
  $\CT_0,\ldots,\CT_{n-1},\ \sX_1^{\pm1},\ldots,\sX_n^{\pm1},\ \pi^{\pm1},\ \sqrt\sX^{\pm1}$,
and relations \textup{(a--fh)} of~\refd{GLn} plus

\textup{(i)} $(\sqrt\sX{}^{\pm1})^2=\sX^{\pm1}:=\sX_1^{\pm1}\cdots\sX_n^{\pm1}$;

\textup{(j)} $\sqrt\sX{}^{\pm1}$ commutes with all the $\sX_i^{\pm1}$ and all
the $\CT_i$;

\textup{(k)} $\pi\sqrt\sX\pi^{-1}=\vv\sqrt\sX$.
\edefe

We interpret $\sX_i,\ i=1,\ldots,n$, as the $i$-th
diagonal matrix entry character of the
diagonal torus $T\subset GL(n)$. It gives rise to the same named character
of the Iwahori subgroup $\bI\subset GL(n,\CK)$. We denote by
$\CO_{\Fl_e}\langle\sX_i\rangle$ the (class of the) structure sheaf of the point
orbit $\Fl_e\subset\Fl=\Fl_{GL(n)}$ (the affine flag variety of $GL(n)$) twisted
by the character $\sX_i$. We denote by $\omega_{\ol\Fl_i},\ i=0,\ldots,n-1$, the
(class of the) direct image (wrt the closed embedding
$\ol\Fl_i\hookrightarrow\Fl_{SL(n)}\hookrightarrow\Fl_{GL(n)}$)
of the canonical line bundle on $\ol\Fl_i$
equipped with the natural $\bI\rtimes\BC^\times$-equivariant structure. We set
$\bT_i:=-1-\omega_{\ol\Fl_i}\in K^{\bI\rtimes\BC^\times}(\Fl)$ as in~\refss{nildaha}.
Finally, note that the fixed point set $\Fl^T$ is naturally identified with
the extended affine Weyl group of $GL(n)$, that is the group of $n$-periodic
permutations of $\BZ\colon \sigma(k+n)=\sigma(k)+n$, and the fixed point
$\varpi$ corresponding to the shift permutation $\sigma(k)=k+1$ is a point
$\bI\rtimes\BC^\times$-orbit $\Fl_\varpi$. We denote by
$\varpi\in K^{\bI\rtimes\BC^\times}(\Fl)$ the class of the structure sheaf
$\CO_{\Fl_\varpi}$.

\th{vasser}
There is a unique isomorphism
  $\Phi\colon \CH\!\!\CH(GL(n))\iso K^{\bI\rtimes\BC^\times}(\Fl)$
such that $\Phi(\sX_i)=\CO_{\Fl_e}\langle\sX_i\rangle,\ i=1,\ldots,n$, and
$\Phi(\CT_i)=\bT_i,\ i=0,\ldots,n-1$, and $\Phi(\pi)=\varpi$.
\eth

\prf Same as the one of~\cite[Theorem~2.5.6]{vv}. \epr

As in~\refss{snd}, we have an idempotent
  $\be=[\CO_{\ol\Fl_{w_0}}]\in K^{\bI\rtimes\BC^\times}(\Fl_{SL(n)})\subset
   K^{\bI\rtimes\BC^\times}(\Fl)\simeq\CH\!\!\CH(GL(n))$,
and we define the spherical subalgebras
  $\CH\!\!\CH^{\on{sph}}(GL(n)):=\be\CH\!\!\CH(GL(n))\be$,
and
  $\CH\!\!\CH_e^{\on{sph}}(GL(n)):=\be\CH\!\!\CH_e(GL(n))\be$.
We also define a two-fold cover
  $\tilde{G}:=\{(g\in GL(n),\ y\in\BC^\times) :
   \det(g)=y^2\}\twoheadrightarrow G,\ \bK:=GL(n,\CO),\ \tilde\bK:=\tilde{G}(\CO)$,
and finally $\wt\BC{}^\times$ as the two-fold cover (with coordinate $\vv$) of
$\BC^\times$ (with coordinate $q$).

\cor{vassero}
The isomorphism $\Phi$ of~\reft{vasser} takes the spherical subalgebra\\
$\CHH^{\on{sph}}(GL(n))\subset\CHH(GL(n))$ isomorphically onto
$K^{\bK\rtimes\BC^\times}(\Gr_{GL(n)})\subset K^{\bI\rtimes\BC^\times}(\Fl_{GL(n)})$.
This isomorphism extends uniquely to
  $\CHH^{\on{sph}}_e(GL(n))\iso K^{\tilde\bK\rtimes\wt\BC{}^\times}(\Gr_{GL(n)})$
where the right-hand side is equipped with the algebra structure as in~\refss{dicson}.
\qed
\ecor

The following theorem is proved similarly to~\cite[Theorem~3.1,~Proposition~3.18]{bfn}:

\th{zast}
\textup{(a)} The algebras $K^\bK(\Gr_{GL(n)}),\ K^{\tilde\bK}(\Gr_{GL(n)})$ are
commutative.

\noindent
\textup{(b)} The spectrum of $K^\bK(\Gr_{GL(n)})$ together with the projection
onto $(\BC^\times)^{(n)}=\on{Spec}(K_{GL(n)}(\on{pt}))$ is
naturally isomorphic to $\oZ^n\stackrel{\on{pr}}{\longrightarrow}(\BC^\times)^{(n)}$
(see~\refss{SL2}).

\noindent
\textup{(c)} The spectrum of $K^{\tilde\bK}(\Gr_{GL(n)})$ together with the projection
onto $\on{Spec}(K^\bK(\Gr_{GL(n)}))$ is naturally isomorphic to
$\hZ^n\to\oZ^n$ (see~\refss{SL2}).

\noindent
\textup{(d)} The Poisson structure on $K^\bK(\Gr_{GL(n)})$ arising from
the deformation $K^{\bK\rtimes\BC^\times}(\Gr_{GL(n)})$ corresponds under the above
identification to the negative of the Poisson (symplectic) structure
of~\cite{gsv},~\cite{fkr} on $\oZ^n$.
The Poisson (symplectic) structure on $K^{\tilde\bK}(\Gr_{GL(n)})$ arising
from the deformation $K^{\tilde\bK\rtimes\wt\BC{}^\times}(\Gr_{GL(n)})$ is the negative of the
pull-back of the symplectic structure on $\oZ^n$.
\qed
\eth

\cor{ava}
\textup{(a)} The algebras $\CHH^{\on{sph}}(GL(n))|_{q=1},\
\CHH_e^{\on{sph}}(GL(n))|_{\vv=1}$ are commutative.

\noindent
\textup{(b)} The algebra $\CHH^{\on{sph}}(GL(n))|_{q=1}$
with the subalgebra $\BC[\sX_1^{\pm1},\ldots,\sX_n^{\pm1}]^{\fS_n}$ is
naturally isomorphic to $\BC[\oZ^n]\supset\BC[(\BC^\times)^{(n)}]$.

\noindent
\textup{(c)} The Poisson structures on $\CHH^{\on{sph}}(GL(n))|_{q=1},\
\CHH_e^{\on{sph}}(GL(n))|_{\vv=1}$ arising from the deformations
$\CHH^{\on{sph}}(GL(n)),\ \CHH_e^{\on{sph}}(GL(n))$
correspond under the above identification to the negative of the Poisson
(symplectic) structures of~\cite{gsv},~\cite{fkr} on $\oZ^n,\ \hZ^n$.
\qed
\ecor


\ssec{quantum}{Quantum Poisson reduction}
\

Now again $G$ is an almost simple simply-connected algebraic group.
We consider Lusztig's integral form $U_q(\fg)$ of the quantized universal
enveloping algebra over $\BC[q^{\pm1}]$ with Cartan elements $K_\lambda,\
\lambda\in X$. It is denoted $\dot{\mathbb U}_\CA$ in~\cite[2.2]{v}.
We extend the scalars to $\BC[q^{\frac{\pm1}{2m_{\on{ad}}}}]$ and
consider the integrable representations of $U_q(\fg)$ with weights
in $X$. We consider the reflection equation algebra $\CO_q(G)$
spanned by the matrix coefficients of integrable $U_q(\fg)$-modules
(with weights in $X$); it is denoted $\BF_\CA$ in~\cite[2.2]{v}. The
corresponding integral form $\CalD_q(G)$ of the Heisenberg
double~\cite[Section~3]{sts} (quantum differential operators) is
denoted ${\mathbb D}_\CA$ in~\cite[2.2]{v}. The quasiclassical limit
of $\CalD_q(G)$ is $D_+(G)$ with the Poisson structure $\{,\!\}_+$
considered in~\refss{ham}. The moment map
$\mu\colon (D_+(G),\{,\!\}_+)\to(G,\pi)\times(G,\pi)$ is the quasiclassical
limit of $\mu_q\colon U_q(\fg)\otimes U_q(\fg)\to\CalD_q(G)$ (see
e.g.~\cite{lu}). The Poisson action of $(G,\pi_G)\times(G,\pi_G)$ on
$D_+(G)$ is the quasiclassical limit of the comodule structure of
$\CalD_q(G)$ over $\CO_q(G)\otimes\CO_q(G)$.

Recall the subalgebra $U^c_q(\fn)\subset U_q(\fg)$~\cite[2.2]{se} associated
to a Coxeter element $c$ (we shall omit its dependence on $\{n_{ij}\}_{i,j\in I}$
satisfying~\cite[2.2.2]{se}). The $U_q(\fg)$-module
$U_q(\fg)/(U_q(\fg)\cdot[U^c_q(\fn),U^c_q(\fn)])$ is the quantization of
the coisotropic subvariety $C_c\subset(G,\pi)$ of~\refss{ham}.
Given a pair of Coxeter elements $c,c'$, we consider the left ideal $\CI_{c',c}$
of $\CalD_q(G)$ generated by
$\mu_q([U^{c'}_q(\fn),U^{c'}_q(\fn)]\otimes S[U^c_q(\fn),U^c_q(\fn)])$
where $S$ stands for the antipode. The invariants of $\CalD_q(G)/\CI_{c',c}$ with
respect to the coaction of $\CO_q(B_-)\otimes\CO_q(B_-)$ form an algebra
denoted $\CO_q(\dZ(G))$.

\conj{december}
There is an isomorphism $\CHH^{\on{sph}}(W_e,\widehat{X})\iso\CO_q(\dZc(G))$
equal to $\on{id}_{\dZc(G)}$ at $q=1$.
\econj


\sec{slices}{Multiplicative slices}
\label{slices}


\ssec{ass}{Asymmetric definition}
\

We closely follow the exposition in~\cite[Section~2]{bfn}.
Let $G$ be an adjoint simple complex algebraic group. We fix a Borel and a
Cartan subgroup $G\supset B\supset T$. Let $\Lambda$ be the coweight lattice,
and let $\Lambda_+\subset\Lambda$ be the submonoid spanned by the simple
coroots $\alpha_i,\ i\in I$. The involution $\alpha\mapsto-w_0\alpha$ of
$\Lambda$ restricts to an involution of $\Lambda_+$ and induces an involution
$\alpha_i\mapsto\alpha_{i^*}$ of the set of simple coroots. We will sometimes
write $\alpha^*:=-w_0\alpha$ for short. Let $\lambda$ be a dominant coweight
of $G$, and $\mu\leq\lambda$ an arbitrary coweight of $G$, not necessarily
dominant, such that $\alpha:=\lambda-\mu=\sum_{i\in I}a_i\alpha_i,\ a_i\in\BN$.
We will define the {\em multiplicative (trigonometric)} analogues
$\dW^\lambda_\mu$ of the generalized slices $\ol\CW{}^\lambda_\mu$
of~\cite[2(ii)]{bfn}.

Namely, $\dW^\lambda_\mu$ is the moduli space of the following data:
\par
(a) a $G$-bundle $\CP$ on $\BP^1$;
\par
(b) a trivialization
  $\sigma\colon \CP_{\on{triv}}|_{\BP^1\setminus\{1\}}\iso\CP|_{\BP^1\setminus\{1\}}$
having a pole of degree $\leq\lambda$ at $1\in\BP^1$. This means that for an irreducible
$G$-module $V^{\lambdavee}$ and the associated vector bundle
$\CV^{\lambdavee}_\CP$ on $\BP^1$ we have
  $V^{\lambdavee}\otimes\CO_{\BP^1}(-\langle\lambda,\lambdavee\rangle\cdot1)
   \subset\CV^{\lambdavee}_\CP\subset
   V^{\lambdavee}\otimes\CO_{\BP^1}(-\langle w_0\lambda,\lambdavee\rangle\cdot1)$;
\par
(c) a reduction $\phi$ of $\CP$ to a $B$-bundle ($B$-structure $\phi$ on $\CP$)
such that the induced $T$-bundle $\phi^T$ has degree $w_0\mu$, and the fiber
of $\phi$ at $\infty\in\BP^1$ is $B_-\subset G$
(with respect to the trivialization $\sigma$ of $\CP$ at $\infty\in\BP^1$).
This means in particular that for an
irreducible $G$-module $V^{\lambdavee}$ and the associated vector bundle
$\CV^{\lambdavee}_\CP$ on $\BP^1$ we are given an invertible subsheaf
$\CL_{\lambdavee}\subset\CV^{\lambdavee}_\CP$ of degree
$-\langle w_0\mu,\lambdavee\rangle$.
\par
We require $\phi$ to be transversal at $0\in\BP^1$ to the trivial $B$-structure
$B$ in $\CP_{\on{triv}}$.


\ssec{BD}{Multiplicative BD slices}
\

Let $\unl{\lambda}=(\omega_{i_1},\ldots,\omega_{i_N})$ be a sequence of
fundamental coweights of $G$ such that $\sum_{s=1}^N\omega_{i_s}=\lambda$.
We define $\dW^{\unl\lambda}_\mu$ as the moduli space of the following data:
\par
(a) a collection of points $(\sz_1,\ldots,\sz_N)\in(\BC^\times)^N$;
\par
(b) a $G$-bundle $\CP$ on $\BP^1$;
\par
(c) a trivialization (a section) $\sigma$ of $\CP$ on
$\BP^1\setminus\{\sz_1,\ldots,\sz_N\}$ with a pole of degree
$\leq\sum_{s=1}^N\omega_{i_s}\cdot \sz_s$ on the complement;
\par
(d) a reduction $\phi$ of $\CP$ to a $B$-bundle ($B$-structure $\phi$ on $\CP$)
such that the induced $T$-bundle $\phi^T$ has degree $w_0\mu$, and the fiber
of $\phi$ at $\infty\in\BP^1$ is $B_-\subset G$ and transversal to $B$ at $0\in\BP^1$
(with respect to the trivialization $\sigma$).

\rem{2x}
The definition of multiplicative BD slices differs from the definition
of BD slices in~\cite[2(x)]{bfn} only by the open condition of transversality
at $0\in\BP^1$. Thus $\dW^{\unl\lambda}_\mu$ is an open subvariety in
$\ol\CW{}^{\unl\lambda}_\mu$ (and similarly, $\dW^\lambda_\mu$ is an open
subvariety in $\ol\CW{}^\lambda_\mu)$. Hence, the favorable properties of
the slices of~\cite{bfn} (like e.g.\ the Cohen-Macaulay property) are
inherited by the multiplicative slices.
\erem


\ssec{sym}{A symmetric definition}
\

Given arbitrary coweights $\mu_-,\mu_+$ such that $\mu_-+\mu_+=\mu$,
we consider the moduli space $\dW^{\unl\lambda}_{\mu_-,\mu_+}$ of the following data:
\par
(a) a collection of points $(\sz_1,\ldots,\sz_N)\in(\BC^\times)^N$;
\par
(b) $G$-bundles $\CP_-,\CP_+$ on $\BP^1$;
\par
(c) an isomorphism
  $\sigma\colon \CP_-|_{\BP^1\setminus\{\sz_1,\ldots,\sz_N\}}\iso
   \CP_+|_{\BP^1\setminus\{\sz_1,\ldots,\sz_N\}}$
with a pole of degree $\leq\sum_{s=1}^N\omega_{i_s}\cdot \sz_s$ on the complement;
\par
(d) a trivialization of $\CP_-=\CP_+$ at $\infty\in\BP^1$;
\par
(e) a reduction $\phi_-$ of $\CP_-$ to a $B_-$-bundle
(a $B_-$-structure on $\CP_-$) such that the induced $T$-bundle
$\phi_-^T$ has degree $-w_0\mu_-$, and the fiber of $\phi_-$ at
$\infty\in\BP^1$ is $B\subset G$;
\par
(f) a reduction $\phi_+$ of $\CP_+$ to a $B$-bundle (a $B$-structure on $\CP_+$)
such that the induced $T$-bundle $\phi_+^T$ has degree $w_0\mu_+$,
and the fiber of $\phi_+$ at $\infty\in\BP^1$ is $B_-\subset G$.
\par
We require $\phi_-$ and $\phi_+$ to be transversal at $0\in\BP^1$
(with respect to the isomorphism $\sigma$).

Note that the trivial $G$-bundle on $\BP^1$ has a unique $B_-$-reduction
of degree 0 with fiber $B$ at $\infty$. Conversely, a $G$-bundle $\CP_-$
with a $B_-$-structure of degree 0 is necessarily trivial, and its
trivialization at $\infty$ uniquely extends to the whole of $\BP^1$.
Hence $\dW^{\unl\lambda}_{0,\mu}=\ \dW^{\unl\lambda}_\mu$.

For arbitrary $\dW^{\unl\lambda}_{\mu_-,\mu_+}$, the $G$-bundles
$\CP_-,\CP_+$ are identified via $\sigma$ on
$\BP^1\setminus\{\sz_1,\ldots,\sz_N\}$, so they are both equipped with
$B$ and $B_-$-structures transversal around $0,\infty\in\BP^1$, that is they
are both equipped with a reduction to a $T$-bundle around $0,\infty\in\BP^1$.
So $\CP_\pm=\CP_\pm^T\times^TG$
for certain $T$-bundles $\CP_\pm^T$ around $0,\infty\in\BP^1$, trivialized at
$\infty\in\BP^1$. The modified $T$-bundles
$'\CP_\pm^T:=\CP_\pm^T(w_0\mu_-\cdot\infty)$ are canonically isomorphic to
$\CP_\pm^T$ off $\infty\in\BP^1$ and trivialized at $\infty\in\BP^1$.
We define $'\CP_\pm$ as the
result of gluing $\CP_\pm$ and $'\CP_\pm^T\times^TG$ in the punctured
neighbourhood of $\infty\in\BP^1$. Then the isomorphism
  $\sigma\colon '\CP_-|_{\BP^1\setminus\{\infty,\sz_1,\ldots,\sz_N\}}\iso\
   '\CP_+|_{\BP^1\setminus\{\infty,\sz_1,\ldots,\sz_N\}}$
extends to $\BP^1\setminus\{\sz_1,\ldots,\sz_N\}$,
and $\phi_\pm$ also extends from $\BP^1\setminus\{\infty\}$ to a
$B$-structure $'\phi_+$ in $'\CP_+$ of degree $w_0\mu$
(resp.\ a $B_-$-structure $'\phi_-$ on $'\CP_-$ of degree 0).

This defines an isomorphism
$\dW^{\unl\lambda}_{\mu_-,\mu_+}\simeq\dW^{\unl\lambda}_\mu$.
Similarly, for the nondeformed slices we have an isomorphism
$\dW^\lambda_{\mu_-,\mu_+}\simeq\dW^\lambda_\mu$.



\ssec{mult}{Multiplication of slices}
\

Given $\lambda_1\geq\mu_1$ and $\lambda_2\geq\mu_2$ with $\lambda_1,\lambda_2$
dominant, we think of $\dW^{\lambda_1}_{\mu_1}$ (resp.\ $\dW^{\lambda_2}_{\mu_2}$)
in the incarnation $\dW^{\lambda_1}_{\mu_1,0}$ (resp.\ $\dW^{\lambda_2}_{0,\mu_2}$).
Note that $\CP^2_-$ is canonically trivialized as in~\refss{sym},
and $\CP^1_+$ is canonically trivialized for the same reason.
Given $(\CP^1_\pm,\sigma_1,\phi^1_\pm)\in\dW^{\lambda_1}_{\mu_1,0}$,
we change the trivialization of $\CP^1_+$ by a (uniquely determined)
element of $U_-$ (the unipotent radical of $B_-$) so that the value
$\phi^1_-(0)$ becomes $B$ (while $\phi^1_+(0)$ remains equal to $B_-$).
Now the value $\phi^1_-(\infty)$ is not $B$ anymore; it is only transversal
to $B_-$. In order to distinguish the data obtained by the composition with
the above trivialization change, we denote them by
$({}'\CP^1_\pm,{}'\!\sigma_1,{}'\!\phi^1_\pm)$.
Given $(\CP^2_\pm,\sigma_2,\phi^2_\pm)\in\ \dW^{\lambda_2}_{0,\mu_2}$, we consider
$({}'\CP^1_-,\CP^2_+,\sigma_2\circ{}'\!\sigma_1,{}'\!\phi^1_-,\phi^2_+)$
(recall that $'\CP^1_+=\CP_{\on{triv}}=\CP^2_-$). These data {\em do not} lie in
$\dW^{\lambda_1+\lambda_2}_{\mu_1,\mu_2}$ since the value $'\!\phi^1_-(\infty)$ is not
necessarily equal to $B$, it is only transversal to $B_-$. However, we change
the trivialization of $'\CP^1_-(\infty)=\CP^2_+(\infty)$ by a (uniquely
determined) element of $U_-$, so that the value of $'\!\phi^1_-(\infty)$
becomes $B$, and we end up in
$\dW^{\lambda_1+\lambda_2}_{\mu_1,\mu_2}=\dW^{\lambda_1+\lambda_2}_{\mu_1+\mu_2}$.

This defines a multiplication morphism
  $\dW^{\lambda_1}_{\mu_1}\times{}
   \dW^{\lambda_2}_{\mu_2}\to\dW^{\lambda_1+\lambda_2}_{\mu_1+\mu_2}$.

In particular, taking $\mu_2=\lambda_2$ so that $\dW^{\lambda_2}_{\lambda_2}$ is a point
and $\dW^{\lambda_1}_{\mu_1}\times\dW^{\lambda_2}_{\lambda_2}={} \dW^{\lambda_1}_{\mu_1}$,
we get a stabilization morphism $\dW^{\lambda_1}_{\mu_1}\to\dW^{\lambda_1+\lambda_2}_{\mu_1+\lambda_2}$.

\rem{2vi}
The multiplication of slices in~\cite[2(vi)]{bfn} {\em does not} preserve the
multiplicative slices viewed as open subvarieties according to~\refr{2x}
(in particular, it {\em does not} induce the above multiplication on multiplicative slices).
\erem


\ssec{scat}{Scattering matrix}
\

Given a collection $(\sz_1,\ldots,\sz_N)\in(\BC^\times)^N$, we define
$P_{\unl\sz}(z):=\prod_{s=1}^N(z-\sz_s)\in\BC[z]$.
We also define a closed subvariety
$\dW^{\unl{\lambda},\unl\sz}_\mu\subset\dW^{\unl{\lambda}}_\mu$ as the fiber of the
latter over $\unl{\sz}=(\sz_1,\ldots,\sz_N)$. We construct a locally closed embedding
$\Psi\colon \dW^{\unl{\lambda},\unl\sz}_\mu\hookrightarrow G[z,P^{-1}]$
into an ind-affine scheme as follows. According to~\refss{sym}, we have
an isomorphism
  $\zeta\colon \dW^{\unl{\lambda},\unl\sz}_\mu={}
   \dW^{\unl{\lambda},\unl\sz}_{0,\mu}\iso\dW^{\unl{\lambda},\unl\sz}_{\mu,0}$.
We denote $\zeta(\CP_\pm,\sigma,\phi_\pm)$ by
$(\CP'_\pm,\sigma',\phi'_\pm)$. Note that $\CP_-$ and $\CP'_+$
are trivialized, and $\CP'_+$ is obtained from $\CP_+$ by an application of
a certain Hecke transformation at $\infty\in\BP^1$. In particular, we obtain
an isomorphism $\CP_+|_{\BA^1}\iso\CP'_+|_{\BA^1}=\CP_{\on{triv}}|_{\BA^1}$.
As in~\refss{mult}, we change the trivialization of $\CP'_+$ by a uniquely
defined element of $U_-$ so that the value of $\phi'_-(0)$ becomes $B$.
Now we compose this change of trivialization with the above isomorphism
$\CP_+|_{\BA^1}\iso\CP'_+|_{\BA^1}=\CP_{\on{triv}}|_{\BA^1}$ and with
  $\sigma\colon \CP_{\on{triv}}|_{\BA^1\setminus\unl\sz}=\CP_-|_{\BA^1\setminus\unl\sz}\iso
   \CP_+|_{\BA^1\setminus\unl\sz}$
to obtain an isomorphism
$\CP_{\on{triv}}|_{\BA^1\setminus\unl\sz}\iso\CP_{\on{triv}}|_{\BA^1\setminus\unl\sz}$,
i.e.\ an element of $G[z,P^{-1}]$.

Here is an equivalent construction of the above embedding. Given
  $(\CP_{\pm},\sigma,\phi_\pm)\in\dW^{\unl{\lambda},\unl\sz}_{\mu_-,\mu_+}={}
   \dW^{\unl{\lambda},\unl\sz}_\mu$,
we choose a trivialization of the $B$-bundle
$\phi_+|_{\BA^1}$ (resp.\ of the $B_-$-bundle $\phi_-|_{\BA^1}$).
This trivialization gives rise to a trivialization of the $G$-bundle
$\CP_+|_{\BA^1}$ (resp.\ of $\CP_-|_{\BA^1}$), so that $\sigma$ becomes an
element of $G(z)$ regular at $0\in\BP^1$; moreover, the value of $\sigma(0)$
lies in the big Bruhat cell $B\cdot B_-\subset G$. We require that
$\sigma(0)\in B\subset G$. Then $\sigma$ is
well-defined up to the left multiplication by an element of $B[z]$ and the
right multiplication by an element of $B_{-,1}[z]$ (the kernel of evaluation at
$0\in\BP^1\colon B_-[z]\to B_-$), i.e.\ $\sigma$ is a well-defined element of
$B[z]\backslash G(z)/B_{-,1}[z]$. Clearly, this element of $G(z)$ lies in
the closure of the double coset $\ol{G[z]z^{\unl{\lambda},\unl\sz}G[z]}$ where
$z^{\unl{\lambda},\unl\sz}:=\prod_{s=1}^N(z-\sz_s)^{\omega_{i_s}}$.
Moreover, it lies in $\ol{G[z]z^{\unl{\lambda},\unl\sz}G[z]}\cap\on{ev}_0^{-1}(B)$.
Thus we have constructed an embedding
  $\Psi'\colon \dW^{\unl{\lambda},\unl\sz}_\mu\hookrightarrow
   B[z]\backslash(\ol{G[z]z^{\unl{\lambda},\unl\sz}G[z]}\cap
   \on{ev}_0^{-1}(B))/B_{-,1}[z]$.
If we compose with an embedding $G(z)\hookrightarrow G((z^{-1}))$, then the image of
$\Psi'$ lies in $B[z]\backslash U_1[[z^{-1}]]T_1[[z^{-1}]]z^\mu U_-[[z^{-1}]]/U_{-,1}[z]$
where $U_1[[z^{-1}]]\subset U[[z^{-1}]]$ (resp.\ $T_1[[z^{-1}]]\subset T[[z^{-1}]]$)
stands for the kernel of evaluation at $\infty\in\BP^1$.
However, the projection
\begin{equation*}
  U_1[[z^{-1}]]T_1[[z^{-1}]]z^\mu U_-[[z^{-1}]]\to
  B[z]\backslash U_1[[z^{-1}]]T_1[[z^{-1}]]z^\mu U_-[[z^{-1}]]/U_{-,1}[z]
\end{equation*}
is clearly one-to-one.
Summing up, we obtain an embedding
  $$\Psi\colon
    \dW^{\unl{\lambda},\unl\sz}_\mu\to U_1[[z^{-1}]]T_1[[z^{-1}]]z^\mu
    U_-[[z^{-1}]]\cap\ol{G[z]z^{\unl{\lambda},\unl\sz}G[z]}\cap\on{ev}_0^{-1}(B).$$
We claim that $\Psi$ is an isomorphism. To see it, we construct the inverse
map to $\dW^{\unl{\lambda},\unl\sz}_{0,\mu}$: given
  $g(z)\in U_1[[z^{-1}]]T_1[[z^{-1}]]z^\mu
   U_-[[z^{-1}]]\cap\ol{G[z]z^{\unl{\lambda},\unl\sz}G[z]}\cap\on{ev}_0^{-1}(B)$
we use it to glue $\CP_+$ together with a rational isomorphism
$\sigma\colon \CP_{\on{triv}}=\CP_-\to\CP_+$, and define $\phi_+$ as the image of
the standard trivial $B$-structure in $\CP_{\on{triv}}$ under $\sigma$.

\rem{2xi}
The embedding $\ol\CW{}^{\unl{\lambda},\unl\sz}_\mu\hookrightarrow G(z)$
of~\cite[2(xi)]{bfn} restricted to the open subvariety
$\dW^{\unl{\lambda},\unl\sz}_\mu\subset\ol\CW{}^{\unl{\lambda},\unl\sz}_\mu$
{\em does not} give the above embedding
$\dW^{\unl{\lambda},\unl\sz}_\mu\hookrightarrow G(z)$.
\erem


\ssec{cover}{A cover of a slice}
\

We define a $T$-torsor $\wW^{\unl\lambda}_\mu\to\dW^{\unl\lambda}_\mu$ as the moduli
space of data (a-d) as in~\refss{BD} plus

\textup{(e)} a collection of nowhere vanishing sections
$u_\lambdavee\in\Gamma(\BP^1\setminus\{\infty\},\CL_\lambdavee)$ satisfying
Pl\"ucker relations (cf.~\refss{ass}(c)).

The construction of~\refss{scat} defines an isomorphism
  $$\widetilde{\Psi}\colon \wW^{\unl{\lambda},\unl\sz}_\mu\iso
    U_1[[z^{-1}]]T[[z^{-1}]]z^\mu
    U_-[[z^{-1}]]\cap\ol{G[z]z^{\unl{\lambda},\unl\sz}G[z]}\cap\on{ev}_0^{-1}(B).$$

Let $T_{[2]}\subset T$ be the subgroup of 2-torsion.
For a future use we define a $T_{[2]}$-torsor
  $\wW^{\unl{\lambda},\unl\sz}_\mu\supset\hW^{\unl{\lambda},\unl\sz}_\mu
   \to\dW^{\unl{\lambda},\unl\sz}_\mu$
as follows.
The evaluation at $0\in\BP^1$ gives rise to a projection
$\on{pr}_0\colon \ol{G[z]z^{\unl{\lambda},\unl\sz}G[z]}\cap\on{ev}_0^{-1}(B)\to B\to T$.
The leading coefficient (at $z^\mu$) gives rise to a projection
$\on{pr}_\infty\colon U_1[[z^{-1}]]T[[z^{-1}]]z^\mu U_-[[z^{-1}]]\to T$, and
$\hW^{\unl{\lambda},\unl\sz}_\mu$ is cut out by the equation
$\on{pr}_0\cdot\on{pr}_\infty=(-1)^{\lambda-\mu}\in T_{[2]}$, where
$\lambda=\sum_{s=1}^N\omega_{i_s}$, see~\refss{BD}.
As $\sz$ varies, we obtain a $T_{[2]}$-torsor
$\wW^{\unl{\lambda}}_\mu\supset\hW^{\unl{\lambda}}_\mu\to\dW^{\unl{\lambda}}_\mu$.


\ssec{example}{An example}
\

This section is parallel to~\cite[2(xii)]{bfn}, but our present conventions
are slightly different. Let $G=GL(2)=GL(V)$ with $V=\BC e_1\oplus\BC e_2$. Let
$B$ be the stabilizer of $\BC e_2$ (the lower triangular matrices),
and let $B_-$ be the stabilizer of $\BC e_1$ (the upper triangular
matrices). Let $N,m\in\BN;\ \unl{\lambda}$ be an $N$-tuple of
fundamental coweights $(0,1)$, and $\mu=(m,N-m)$, so that
$w_0\mu=(N-m,m)$. Let $\CO:=\CO_{\BP^1}$. We fix a collection
$(\sz_1,\ldots,\sz_N)\in(\BC^\times)^N$ and define
$P_{\unl\sz}(z):=\prod_{s=1}^N(z-\sz_s)\in\BC[z]$. Then
$\dW^{\unl{\lambda},\unl\sz}_\mu$ is the moduli space of flags
$(\CO\otimes V\supset\CV\supset\CL)$, where

\textup{(a)} $\CV$ is a 2-dimensional
locally free subsheaf in $\CO\otimes V$ coinciding with $\CO\otimes V$
around $0,\infty\in\BP^1$ and such that on $\BA^1\subset\BP^1$ the global
sections of $\det\CV$ coincide with $P_{\unl\sz}\BC[z]e_1\wedge e_2$
as a $\BC[z]$-submodule of
$\Gamma(\BA^1,\det(\CO_{\BA^1}\otimes V))=\BC[z]e_1\wedge e_2$.

\textup{(b)} $\CL$ is a line subbundle in $\CV$ of degree
$-m$, assuming the value $\BC e_1$ at $\infty\in\BP^1$, and such that
the value of $\CL$ at $0\in\BP^1$ is transversal to $\BC e_2$.
In particular, $\deg\CV/\CL=m-N$.

On the other hand, let us introduce a closed subvariety
$\hW^{\unl{\lambda},\unl\sz}_\mu$ in $\on{Mat}_2[z]$
formed by all the matrices $\sM=\begin{pmatrix}A&B\\ C&D\end{pmatrix}$
such that $A(z)=a_mz^m+\ldots+a_0$, and $a_m\cdot a_0=(-1)^m$, while
$\deg C(z)<m\geq\deg B(z)$, and $B(0)=0$; furthermore, $\det\sM=P_{\unl\sz}(z)$.

Then we have a two-fold cover $\mho\colon
\hW^{\unl{\lambda},\unl\sz}_\mu\to\dW^{\unl{\lambda},\unl\sz}_\mu$:
given $\sM\in\hW^{\unl{\lambda},\unl\sz}_\mu$ we view it as a transition matrix
in a punctured neighbourhood of $\infty\in\BP^1$ to glue a vector bundle
$\CV$ which embeds, by construction, as a locally free subsheaf into
$\CO\otimes V$. The morphism
$\sM\CO_{\BA^1}e_1\hookrightarrow\CO_{\BA^1}\otimes V$ naturally extends to
$\infty\in\BP^1$ with a pole of degree $m$, hence it extends to an embedding
of $\CO(-m\cdot\infty)$ into $\CV\subset\CO\otimes V$. The image of this
embedding is the desired line subbundle $\CL\subset\CV$.


\ssec{thick}{Thick slices}
\

We define {\em thick} multiplicative
(trigonometric) slices $^\dagger\CW_\mu$ as the moduli space of the
following data:

\textup{(a)} a $G$-bundle $\CP$ on $\BP^1$;

\textup{(b)} a trivialization
$\sigma\colon \CP_{\on{triv}}|_{\widehat\BP{}^1_\infty}\iso\CP|_{\widehat\BP{}^1_\infty}$
in the formal neighbourhood of $\infty\in\BP^1$;

\textup{(c)} a reduction $\phi$ of $\CP$ to a $B$-bundle
($B$-structure $\phi$ on $\CP$)
such that the induced $T$-bundle $\phi^T$ has degree $w_0\mu$, and the fiber
of $\phi$ at $\infty\in\BP^1$ is
transversal to $B$ (with respect to the trivialization
$\sigma$ of $\CP$ at $\infty\in\BP^1$);

\textup{(d)} a collection of nowhere vanishing sections
$u_\lambdavee\in\Gamma(\BP^1\setminus\{\infty\},\CL_\lambdavee)$ satisfying
Pl\"ucker relations (cf.~\refss{ass}(c)).

The construction of~\refss{cover} identifies $^\dagger\CW_\mu$ with the
infinite type scheme (cf.~\cite[5.9]{fkp})
\eq{53}
  ^\dagger\CW_\mu\simeq U_1[[z^{-1}]]T[[z^{-1}]]z^\mu U_-[[z^{-1}]]\subset G((z^{-1})).
\end{equation}
As the inclusion $U_1[[z^{-1}]]\hookrightarrow U((z^{-1}))$ gives rise
to an isomorphism $U_1[[z^{-1}]]\simeq U[z]\backslash U((z^{-1}))$,
we can identify $^\dagger\CW_\mu$ with the quotient
$U[z]\backslash U((z^{-1}))T[[z^{-1}]]z^\mu U_-((z^{-1}))/U_{-,1}[z],$ and we
write $\pi$ for this isomorphism.
The construction of~\refss{scat} (resp.\ of~\refss{cover}) defines a closed
embedding $\dW^\lambda_\mu\hookrightarrow{}^\dagger\CW_\mu$
(resp.~$\hW^\lambda_\mu\hookrightarrow{}^\dagger\CW_\mu$).
We define the multiplication morphism
  $m_{\mu_1,\mu_2}\colon
   ^\dagger\CW_{\mu_1}\times{}^\dagger\CW_{\mu_2}\to{}^\dagger\CW_{\mu_1+\mu_2}$
by the formula $m_{\mu_1,\mu_2}(g_1,g_2)=\pi(g_1g_2)$. Then the multiplication morphism
  $m^{\lambda_1,\lambda_2}_{\mu_1,\mu_2}\colon \dW^{\lambda_1}_{\mu_1}\times
   \dW^{\lambda_2}_{\mu_2}\to\dW^{\lambda_1+\lambda_2}_{\mu_1+\mu_2}$
of~\refss{mult} is the restriction of $m_{\mu_1,\mu_2}$. Similarly, $m_{\mu_1,\mu_2}$
restricts to a multiplication $\hW^{\lambda_1}_{\mu_1}\times
\hW^{\lambda_2}_{\mu_2}\to\hW^{\lambda_1+\lambda_2}_{\mu_1+\mu_2}$.

For $\nu_1,\nu_2$ antidominant, we define the {\em shift maps}
$\iota_{\mu,\nu_1,\nu_2}\colon
^\dagger\CW_{\mu+\nu_1+\nu_2}\to{}^\dagger\CW_\mu$ by $g\mapsto\pi(z^{-\nu_1}gz^{-\nu_2})$.


\section{Shifted quantum affine algebras}
\label{sqaa}

Let $\fg$ be a simple Lie algebra, $\fh\subset \fg$ be a Cartan
subalgebra of $\fg$, and $(\cdot,\cdot)$ be a non-degenerate
invariant bilinear symmetric form on $\fg$ (with a square length of
the shortest root equal to $2$). Let
$\{\alphavee_i\}_{i\in I}\subset \fh^*$
be the simple positive roots of $\fg$ relative to $\fh$, and
$c_{ij}=2\frac{(\alphavee_i,\alphavee_j)}{(\alphavee_i,\alphavee_i)}$--the
entries of the corresponding Cartan matrix. Set
$d_i:=\frac{(\alphavee_i,\alphavee_i)}{2}\in \BZ_{>0}$
so that $d_ic_{ij}=d_jc_{ji}$ for any $i,j\in I$. Let
$\nu\colon \fh\iso\fh^*$ be the isomorphism determined by the symmetric
form $(\cdot,\cdot)$ so that $\alpha_i=h_i=\nu^{-1}(\alphavee_i)/d_i$
are the simple coroots of $\fg$.


\subsection{Algebras $\CU^\ssc_{\mu_1,\mu_2}$ and $\CU^\ad_{\mu_1,\mu_2}$}
\

Given coweights $\mu^+,\mu^-\in \Lambda$, set
$\unl{b}^\pm=\{b^\pm_i\}_{i\in I}\in \BZ^I$ with $b^\pm_i:=\alphavee_i(\mu^\pm)$.
Define the \emph{simply-connected version of shifted quantum affine algebra},
denoted by $\CU^\ssc_{\mu^+,\mu^-}$ or $\CU^\ssc_{\unl{b}^+,\unl{b}^-}$,
to be the associative $\BC(\vv)$-algebra generated by
  $\{e_{i,r},f_{i,r},\psi^\pm_{i,\pm s^\pm_i}, (\psi^\pm_{i,\mp b^\pm_i})^{-1}\}_
    {i\in I}^{r\in \BZ, s^\pm_i\geq -b^\pm_i}$
with the following defining relations
(for all $i,j\in I$ and $\epsilon,\epsilon'\in \{\pm\}$):
\begin{equation}\tag{U1} \label{U1}
  [\psi_i^\epsilon(z),\psi_j^{\epsilon'}(w)]=0,\
  \psi^\pm_{i,\mp b^\pm_i}\cdot (\psi^\pm_{i,\mp b^\pm_i})^{-1}=
  (\psi^\pm_{i,\mp b^\pm_i})^{-1}\cdot \psi^\pm_{i,\mp b^\pm_i}=1,
\end{equation}
\begin{equation}\tag{U2} \label{U2}
  (z-\vv_i^{c_{ij}}w)e_i(z)e_j(w)=(\vv_i^{c_{ij}}z-w)e_j(w)e_i(z),
\end{equation}
\begin{equation}\tag{U3} \label{U3}
  (\vv_i^{c_{ij}}z-w)f_i(z)f_j(w)=(z-\vv_i^{c_{ij}}w)f_j(w)f_i(z),
\end{equation}
\begin{equation}\tag{U4} \label{U4}
  (z-\vv_i^{c_{ij}}w)\psi^\epsilon_i(z)e_j(w)=(\vv_i^{c_{ij}}z-w)e_j(w)\psi^\epsilon_i(z),
\end{equation}
\begin{equation}\tag{U5} \label{U5}
  (\vv_i^{c_{ij}}z-w)\psi^\epsilon_i(z)f_j(w)=(z-\vv_i^{c_{ij}}w)f_j(w)\psi^\epsilon_i(z),
\end{equation}
\begin{equation}\tag{U6} \label{U6}
  [e_i(z),f_j(w)]=\frac{\delta_{ij}}{\vv_i-\vv_i^{-1}}\delta\left(\frac{z}{w}\right)\left(\psi^+_i(z)-\psi^-_i(z)\right),
\end{equation}
\begin{equation}\tag{U7}\label{U7}
 \underset{z_1,\ldots,z_{1-c_{ij}}} \Sym\ \sum_{r=0}^{1-c_{ij}}(-1)^r{1-c_{ij}\brack r}_{\vv_i}
  e_i(z_1)\cdots e_i(z_r) e_j(w) e_i(z_{r+1})\cdots e_i(z_{1-c_{ij}})=0,
\end{equation}
\begin{equation}\tag{U8}\label{U8}
 \underset{z_1,\ldots,z_{1-c_{ij}}} \Sym\ \sum_{r=0}^{1-c_{ij}}(-1)^r{1-c_{ij}\brack r}_{\vv_i}
  f_i(z_1)\cdots f_i(z_r) f_j(w) f_i(z_{r+1})\cdots f_i(z_{1-c_{ij}})=0,
\end{equation}
where
  $\vv_i:=\vv^{d_i}, [a,b]_x:=ab-x\cdot ba,
   [m]_\vv:=\frac{\vv^m-\vv^{-m}}{\vv-\vv^{-1}},
   {a\brack b}_{\vv}:=\frac{[a-b+1]_\vv\cdots [a]_\vv}{[1]_\vv\cdots [b]_\vv},$
$\underset{z_1,\ldots,z_{s}} \Sym$ stands for the symmetrization in $z_1,\ldots,z_s$,
and the generating series are defined as follows:
\begin{equation*}
    e_i(z):=\sum_{r\in \BZ}{e_{i,r}z^{-r}},\
    f_i(z):=\sum_{r\in \BZ}{f_{i,r}z^{-r}},\
    \psi_i^{\pm}(z):=\sum_{r\geq -b^\pm_i}{\psi^\pm_{i,\pm r}z^{\mp r}},\
    \delta(z):=\sum_{r\in \BZ}{z^r}.
\end{equation*}
Let us introduce another set of Cartan generators $\{h_{i,\pm r}\}_{i\in I}^{r>0}$
instead of $\{\psi^\pm_{i,\pm s^\pm_i}\}_{i\in I}^{s^\pm_i>-b^\pm_i}$ via
\begin{equation*}
   (\psi^\pm_{i,\mp b^\pm_i}z^{\pm b^\pm_i})^{-1}\psi^\pm_i(z)=
   \exp\left(\pm (\vv_i-\vv_i^{-1}) \sum_{r>0} h_{i,\pm r}z^{\mp r}\right).
\end{equation*}
Then, relations (U4, U5) are equivalent to the following:
\begin{equation}\tag{U4$'$}\label{U4'}
  \psi^\pm_{i,\mp b^\pm_i}e_{j,s}=\vv_i^{\pm c_{ij}}e_{j,s}\psi^\pm_{i,\mp b^\pm_i},\
  [h_{i,r}, e_{j,s}]=\frac{[rc_{ij}]_{\vv_i}}{r}\cdot e_{j,s+r}\ \on{for}\ r\ne 0,
\end{equation}
\begin{equation}\tag{U5$'$}\label{U5'}
  \psi^\pm_{i,\mp b^\pm_i}f_{j,s}=\vv_i^{\mp c_{ij}}f_{j,s}\psi^\pm_{i,\mp b^\pm_i},\
  [h_{i,r}, f_{j,s}]=-\frac{[rc_{ij}]_{\vv_i}}{r}\cdot f_{j,s+r}\ \on{for}\ r\ne 0.
\end{equation}

Let
  $\CU^{\ssc,<}_{\mu^+,\mu^-},\ \CU^{\ssc,>}_{\mu^+,\mu^-},$ and $\CU^{\ssc,0}_{\mu^+,\mu^-}$
be the $\BC(\vv)$-subalgebras of $\CU^\ssc_{\mu^+,\mu^-}$ generated by
  $\{f_{i,r}\}_{i\in I}^{r\in \BZ}$, $\{e_{i,r}\}_{i\in I}^{r\in \BZ}$,
and
  $\{\psi^\pm_{i,\pm s^\pm_i}, (\psi^\pm_{i,\mp b^\pm_i})^{-1}\}_{i\in I}^{s^\pm_i\geq -b^\pm_i},$
respectively. The following is proved completely analogously
to~\cite[Theorem 2]{he}:

\begin{Prop}\label{Triangular decomposition}
(a) (Triangular decomposition of $\CU^\ssc_{\mu^+,\mu^-}$)
The multiplication map
  $$m\colon \CU^{\ssc,<}_{\mu^+,\mu^-}\otimes \CU^{\ssc,0}_{\mu^+,\mu^-}\otimes \CU^{\ssc,>}_{\mu^+,\mu^-}
    \longrightarrow \CU^\ssc_{\mu^+,\mu^-}$$
is an isomorphism of $\BC(\vv)$-vector spaces.

\noindent
(b) The algebra $\CU^{\ssc,0}_{\mu^+,\mu^-}$ (resp.\ $\CU^{\ssc,<}_{\mu^+,\mu^-}$
and $\CU^{\ssc,>}_{\mu^+,\mu^-}$) is isomorphic to the $\BC(\vv)$-algebra generated by
  $\{\psi^\pm_{i,\pm s^\pm_i}, (\psi^\pm_{i,\mp b^\pm_i})^{-1}\}_{i\in I}^{s^\pm_i\geq -b^\pm_i}$
(resp.\ $\{f_{i,r}\}_{i\in I}^{r\in \BZ}$ and $\{e_{i,r}\}_{i\in I}^{r\in \BZ}$)
with the defining relations~(\ref{U1}) (resp.~(\ref{U3}, \ref{U8}) and~(\ref{U2}, \ref{U7})).
In particular, $\CU^{\ssc,<}_{\mu^+,\mu^-}$ and $\CU^{\ssc,>}_{\mu^+,\mu^-}$
are independent of $\mu^\pm$.
\end{Prop}

Following the terminology of~\cite{lus}, we also define the
\emph{adjoint version of shifted quantum affine algebra}, denoted by
$\CU^\ad_{\mu^+,\mu^-}$ or $\CU^\ad_{\unl{b}^+,\unl{b}^-}$, by
adding extra generators $\{(\phi^\pm_i)^{\pm 1}\}_{i\in I}$ to
$\CU^\ssc_{\mu^+,\mu^-}$, which satisfy the following extra relations:
\begin{equation}\tag{U9} \label{U9}
  (\psi^\epsilon_{i,-\epsilon b^\epsilon_i})^{\pm 1}=
  (\phi^\epsilon_i)^{\pm 2}\cdot \prod_{j - i}(\phi^\epsilon_j)^{\pm c_{ji}},\
  (\phi^\epsilon_i)^{\pm 1}\cdot (\phi^\epsilon_i)^{\mp 1}=1,\
  [\phi^\epsilon_i,\phi^{\epsilon'}_j]=0,
\end{equation}
\begin{equation}\tag{U10}\label{U10}
  \phi^\epsilon_i\psi^{\epsilon'}_j(z)=\psi^{\epsilon'}_j(z)\phi^\epsilon_i,\
  \phi^\epsilon_ie_j(z)=\vv_i^{\epsilon\delta_{ij}}e_j(z)\phi^\epsilon_i,\
  \phi^\epsilon_if_j(z)=\vv_i^{-\epsilon\delta_{ij}}f_j(z)\phi^\epsilon_i,
\end{equation}
for any $i,j\in I$ and $\epsilon,\epsilon'\in \{\pm\}$.

Both algebras $\CU^\ssc_{\mu^+,\mu^-}$ and $\CU^\ad_{\mu^+,\mu^-}$ depend
only on $\mu:=\mu^+ + \mu^-$ up to an isomorphism\footnote{For example, there is
an isomorphism $\CU^\ssc_{\mu^+,\mu^-}\iso \CU^\ssc_{0,\mu^++\mu^-}$ such that
$f_i(z)\mapsto f_i(z), e_i(z)\mapsto z^{b^+_i}e_i(z), \psi^\pm_i(z)\mapsto z^{b^+_i}\psi^\pm_i(z)$.}.
Let $\Lambda^\pm \subset \Lambda$ be the submonoids spanned by
$\{\pm \omega_i\}_{i\in I}$, that
is, $\Lambda^+$ (resp.\ $\Lambda^-$) consists of dominant (resp.
antidominant) coweights of $\Lambda$. We will say that the
algebras $\CU^\ssc_{\mu^+,\mu^-}, \CU^\ad_{\mu^+,\mu^-}$ are
dominantly (resp.\ antidominantly) shifted if $\mu\in \Lambda^+$
(resp.\ $\mu\in \Lambda^-$). We note that
  $\mu\in \Lambda^+\Leftrightarrow b^+_i+b^-_i=\alphavee_i(\mu)\geq 0,\
   \mu\in \Lambda^-\Leftrightarrow b^+_i+b^-_i=\alphavee_i(\mu)\leq 0$
for all $i\in I$.

\begin{Rem}
One of the key reasons to consider $\CU^\ad_{\mu^+,\mu^-}$, not only
$\CU^\ssc_{\mu^+,\mu^-}$, is to construct quantizations of the
thick slices $^\dagger\CW_{\mu^*}$ of~\refss{thick} and the
multiplicative slice covers $\hW^{\unl\lambda^*}_{\mu^*}$ of~\refss{cover},
see our~\refco{master conj}. On the technical side, we also need
an alternative set of Cartan generators, whose generating series
$A^\pm_i(z)$ are defined via~(\ref{A}) of Section~\ref{ABCD
generators} and whose definition requires to work with
$\CU^\ad_{\mu^+,\mu^-}$ (see also Remark~\ref{Remark on ABCD}(b)).
\end{Rem}

\begin{Rem}\label{Basic Properties}
(a) The elements $\{\psi^+_{i,-b^+_i}\psi^-_{i,b^-_i}\}_{i\in I}$
(resp.\ $\{\phi^+_i\phi^-_i\}_{i\in I}$) and their inverses are central
elements of $\CU^\ssc_{\mu^+,\mu^-}$ (resp.\ $\CU^\ad_{\mu^+,\mu^-}$).

\noindent
(b) We have
  $\CU_{0,0}^\ssc/(\psi^+_{i,0}\psi^-_{i,0}-1)\simeq U_\vv(L\fg)$--the
standard quantum loop algebra of $\fg$, while
 $\CU_{0,0}^\ad/(\phi^+_i\phi^-_i-1)\simeq U^\ad_\vv(L\fg)$--the
adjoint version of $U_\vv(L\fg)$.

\noindent
(c) We note that defining relations~(\ref{U1}--\ref{U8},~\ref{U10})
are independent of $\mu^+,\mu^-$.

\noindent
(d) An equivalent definition of $\CU^\ssc_{\mu_1,\mu_2}$ was suggested
to us by Boris Feigin in Spring 2010. In this definition, we take the same
generators as for $U_\vv(L\fg)$ and just modify relation~(\ref{U6}) by requesting
  $p_i(z)[e_i(z),f_j(w)]=
   \frac{\delta_{ij}\delta(z/w)}{\vv_i-\vv_i^{-1}}\left(\psi^+_i(z)-\psi^-_i(z)\right)$
for any collection $\{p_i(z)\}_{i\in I}$ of rational functions.
\end{Rem}

\subsection{Levendorskii type presentation of $\CU^\ssc_{0,\mu}$ for $\mu\in \Lambda^-$}
\label{Levendorskii} $\ $

In Section~\ref{Copro}, we will crucially need a presentation of the shifted
quantum affine algebras via a finite number of generators and defining relations.
This is the purpose of this subsection.

Fix antidominant coweights $\mu_1,\mu_2\in \Lambda^-$ and set $\mu:=\mu_1+\mu_2$.
Define
  $b_{1,i}:=\alphavee_i(\mu_1), b_{2,i}:=\alphavee_i(\mu_2), b_i:=b_{1,i}+b_{2,i}$.
Denote by $\hat{\CU}_{\mu_1,\mu_2}$ the associative $\BC(\vv)$-algebra generated by
  $$\{e_{i,r},f_{i,s},(\psi^+_{i,0})^{\pm 1}, (\psi^-_{i,b_i})^{\pm 1}, h_{i,\pm 1}|
    i\in I, b_{2,i}-1\leq r\leq 0, b_{1,i}\leq s\leq 1\}$$
and with the following defining relations:
\begin{equation}\tag{\^{U}1} \label{u1}
\begin{split}
  \{(\psi^+_{i,0})^{\pm 1}, (\psi^-_{i,b_i})^{\pm 1}, h_{i,\pm 1}\}_{i\in I}
  \ \mathrm{pairwise\ commute},\\
  (\psi^+_{i,0})^{\pm 1}\cdot (\psi^+_{i,0})^{\mp 1}=
  (\psi^-_{i,b_i})^{\pm 1}\cdot (\psi^-_{i,b_i})^{\mp 1}=1,
\end{split}
\end{equation}
\begin{equation}\tag{\^{U}2} \label{u2}
  e_{i,r+1}e_{j,s}-\vv_i^{c_{ij}}e_{i,r}e_{j,s+1}=
  \vv_i^{c_{ij}}e_{j,s}e_{i,r+1}-e_{j,s+1}e_{i,r},
\end{equation}
\begin{equation}\tag{\^{U}3} \label{u3}
  \vv_i^{c_{ij}}f_{i,r+1}f_{j,s}-f_{i,r}f_{j,s+1}=
  f_{j,s}f_{i,r+1}-\vv_i^{c_{ij}}f_{j,s+1}f_{i,r},
\end{equation}
\begin{equation}\tag{\^{U}4} \label{u4}
  \psi^+_{i,0}e_{j,r}=\vv_i^{c_{ij}}e_{j,r}\psi^+_{i,0},\
  \psi^-_{i,b_i}e_{j,r}=\vv_i^{-c_{ij}}e_{j,r}\psi^-_{i,b_i},\
  [h_{i,\pm 1}, e_{j,r}]=[c_{ij}]_{\vv_i}\cdot e_{j,r\pm 1},
\end{equation}
\begin{equation}\tag{\^{U}5} \label{u5}
  \psi^+_{i,0}f_{j,s}=\vv_i^{-c_{ij}}f_{j,s}\psi^+_{i,0},\
  \psi^-_{i,b_i}f_{j,s}=\vv_i^{c_{ij}}f_{j,s}\psi^-_{i,b_i},\
  [h_{i,\pm 1}, f_{j,s}]=-[c_{ij}]_{\vv_i}\cdot f_{j,s\pm 1},
\end{equation}
\begin{equation}\tag{\^{U}6} \label{u6}
  [e_{i,r},f_{j,s}]=0\ \on{if}\ i\ne j \ \ \mathrm{and}\ \ [e_{i,r},f_{i,s}]=
  \begin{cases}
       \psi^+_{i,0}h_{i,1} & \text{if}\ \ r+s=1, \\
       \psi^-_{i,b_i}h_{i,-1} & \text{if}\ \ r+s=b_i-1,\\
\frac{\psi^+_{i,0}-\delta_{b_i,0}\psi^-_{i,b_i}}{\vv_i-\vv_i^{-1}}& \text{if}\ \ r+s=0, \\
\frac{-\psi^-_{i,b_i}+\delta_{b_i,0}\psi^+_{i,0}}{\vv_i-\vv_i^{-1}} & \text{if}\ \ r+s=b_i, \\
       0 & \text{if}\ \ b_i<r+s<0,
  \end{cases}
\end{equation}
\begin{equation}\tag{\^{U}7} \label{u7}
  [e_{i,0},[e_{i,0},\cdots,[e_{i,0},e_{j,0}]_{\vv^{c_{ij}}_i}\cdots]_{\vv^{-c_{ij}-2}_i}]_{\vv^{-c_{ij}}_i}=0\
  \mathrm{for}\ i\ne j,
\end{equation}
\begin{equation}\tag{\^{U}8} \label{u8}
  [f_{i,0},[f_{i,0},\cdots,[f_{i,0},f_{j,0}]_{\vv^{c_{ij}}_i}\cdots]_{\vv^{-c_{ij}-2}_i}]_{\vv^{-c_{ij}}_i}=0\
  \mathrm{for}\ i\ne j,
\end{equation}
\begin{equation}\tag{\^{U}9} \label{u9}
  [h_{i,1},[f_{i,1},[h_{i,1},e_{i,0}]]]=0,\
  [h_{i,-1},[e_{i,b_{2,i}-1},[h_{i,-1},f_{i,b_{1,i}}]]]=0,
\end{equation}
for any $i,j\in I$ and $r,s$ such that above relations make
sense.

\begin{Rem}
\label{56}
One can rewrite relations~(\ref{u7},~\ref{u8}) in the form similar
to~(\ref{U7},~\ref{U8}) as
$$\sum^{1-c_{ij}}_{r=0}(-1)^r{1-c_{ij}\brack r}_{\vv_i}e_{i,0}^re_{j,0}e_{i,0}^{1-c_{ij}-r}=0,\
\sum^{1-c_{ij}}_{r=0}(-1)^r{1-c_{ij}\brack r}_{\vv_i}f_{i,0}^rf_{j,0}f_{i,0}^{1-c_{ij}-r}=0.$$
\end{Rem}

Define inductively
\begin{equation*}
    e_{i,r}:=[2]^{-1}_{\vv_i}\cdot
    \begin{cases}
       [h_{i,1},e_{i,r-1}] & \text{if}\ \ r>0,\\
       [h_{i,-1},e_{i,r+1}] & \text{if}\ \ r<b_{2,i}-1,
    \end{cases}
\end{equation*}
\begin{equation*}
    f_{i,r}:=-[2]^{-1}_{\vv_i}\cdot
    \begin{cases}
       [h_{i,1},f_{i,r-1}] & \text{if}\ \ r>1,\\
       [h_{i,-1},f_{i,r+1}] & \text{if}\ \ r<b_{1,i},
    \end{cases}
\end{equation*}
\begin{equation*}
    \psi^+_{i,r}:=(\vv_i-\vv_i^{-1})\cdot [e_{i,r-1},f_{i,1}]\ \on{for}\ r>0,
\end{equation*}
\begin{equation*}
    \psi^-_{i,r}:=(\vv_i^{-1}-\vv_i)\cdot [e_{i,r-b_{1,i}},f_{i,b_{1,i}}]\
    \on{for}\ r<b_i.
\end{equation*}

\begin{Thm}\label{Levendorskii-type}
There is a unique $\BC(\vv)$-algebra isomorphism
$\hat{\CU}_{\mu_1,\mu_2}\iso \CU^\ssc_{0,\mu}$, such that
\begin{equation*}
    e_{i,r}\mapsto e_{i,r},\
    f_{i,r}\mapsto f_{i,r},\
    \psi^\pm_{i,\pm s^\pm_i}\mapsto \psi^\pm_{i,\pm s^\pm_i}\ \
    \mathrm{for}\ i\in I,r\in \BZ, s^+_i\geq 0, s^-_i\geq -b_i.
\end{equation*}
\end{Thm}

This provides a new presentation of $\CU^\ssc_{0,\mu}$ via a finite
number of generators and relations. The proof of this result is
presented in Appendix~\ref{Proof of Levendorskii-type}. Motivated
by~\cite{gnw}, we also provide a slight modification of this
presentation of $\CU^\ssc_{0,\mu}$ in Theorem~\ref{Levendorskii-type update}.

\begin{Rem}
Theorem~\ref{Levendorskii-type} can be viewed as a $\vv$-version of the
corresponding result for the shifted Yangians of~\cite[Theorem~4.3]{fkp}.
In the particular case $\mu_1=\mu_2=0$, the latter is the standard
Levendorskii presentation of the Yangian, see~\cite{le}. However, we are
not aware of the reference for Theorem~\ref{Levendorskii-type} even
in the {\em unshifted case} $\mu_1=\mu_2=0$.
\end{Rem}


\section{$ABCD$ generators of $\CU^\ad_{\mu^+,\mu^-}$}
\label{ABCD generators}

In this section, we introduce an alternative set of generators of
$\CU^\ad_{\mu^+,\mu^-}$, which will be used later in the paper (they
are also of independent interest), and deduce the defining relations
among them. While the definition works for any two coweights
$\mu^+,\mu^-\in \Lambda$, the relations hold only for antidominant
$\mu^+,\mu^-\in \Lambda^-$, which we assume from now on.

First, we define the Cartan generators
$\{A^\pm_{i,\pm r}\}_{i\in I}^{r\geq 0}$ via
\begin{equation}\label{A}
   z^{\mp b^\pm_i}\psi^\pm_i(z)=
   \frac{\prod_{j-i}\prod_{p=1}^{-c_{ji}}A^\pm_j(\vv_j^{-c_{ji}-2p}z)}{A^\pm_i(z)A^\pm_i(\vv_i^{-2}z)}
   \ \mathrm{with}\  A^\pm_{i,0}:=(\phi^\pm_i)^{-1},
\end{equation}
where we set $A^\pm_i(z)=\sum_{r\geq 0}A^\pm_{i,\pm r}z^{\mp r}$.
Using non-degeneracy of the $\vv$-version of the Cartan matrix $(c_{ij})$
and arguing by induction in $r>0$, one can easily see that
relations~(\ref{A}) for all $i\in I$ determine uniquely all
$A^\pm_{i,\pm r}$, see Remark~\ref{uniqueness of A} (cf.~\cite[Lemma~2.1]{gklo1}).
An explicit formula for $A^\pm_i(z)$ is given by~(\ref{explicit A}) in
Appendix~\ref{Proof of Relations bw ABCD}.

Next, we introduce the generating series
$B^\pm_i(z), C^\pm_i(z), D^\pm_i(z)$ via
\begin{equation}\label{B}
  B^\pm_i(z):=(\vv_i-\vv_i^{-1})A^\pm_i(z)e^\pm_i(z),
\end{equation}
\begin{equation}\label{C}
  C^\pm_i(z):=(\vv_i-\vv_i^{-1})f^\pm_i(z)A^\pm_i(z),
\end{equation}
\begin{equation}\label{D}
  D^\pm_i(z):=A^\pm_i(z)\psi^\pm_i(z)+(\vv_i-\vv_i^{-1})^2 f^\pm_i(z)A^\pm_i(z)e^\pm_i(z),
\end{equation}
where the \emph{Drinfeld half-currents} are defined as follows:
\begin{equation}\label{Drinfeld currents}
\begin{split}
  & e^+_i(z):=\sum_{r\geq 0} e_{i,r}z^{-r},\ e^-_i(z):=-\sum_{r<0} e_{i,r}z^{-r},\\
  & f^+_i(z):=\sum_{r>0} f_{i,r}z^{-r},\ f^-_i(z):=-\sum_{r\leq 0} f_{i,r}z^{-r}.
\end{split}
\end{equation}

It is clear that coefficients of the generating series
  $\{A^\pm_i(z), B^\pm_i(z), C^\pm_i(z), D^\pm_i(z)\}_{i\in I}$
together with $\{\phi^\pm_i\}_{i\in I}$ generate (over $\BC(\vv)$)
the shifted quantum affine algebra $\CU^\ad_{\mu^+,\mu^-}$. The
following is the key result of this section.

\begin{Thm}\label{Relations bw ABCD}
Assume $\mu^+,\mu^-\in \Lambda^-$ and define $\{b^\pm_i\}_{i\in I}$ via
$b^\pm_i:=\alphavee_i(\mu^\pm)$ as before.

\noindent
(a) The generating series $A^\pm_i(z), B^\pm_i(z), C^\pm_i(z), D^\pm_i(z)$
satisfy the following relations:
\begin{equation}\label{ABCD0}
\begin{split}
  \phi^{\epsilon}_i A^{\epsilon'}_j(w)=A^{\epsilon'}_j(w)\phi^{\epsilon}_i,\
  \phi^{\epsilon}_i D^{\epsilon'}_j(w)=D^{\epsilon'}_j(w)\phi^{\epsilon}_i,\\
  \phi^{\epsilon}_i B^{\epsilon'}_j(w)=\vv_i^{\epsilon\delta_{ij}}B^{\epsilon'}_j(w)\phi^{\epsilon}_i,\
  \phi^{\epsilon}_i C^{\epsilon'}_j(w)=\vv_i^{-\epsilon\delta_{ij}}C^{\epsilon'}_j(w)\phi^{\epsilon}_i,\
\end{split}
\end{equation}
\begin{equation}\label{ABCD1}
  [A^{\epsilon}_i(z), A^{\epsilon'}_j(w)]=0,
\end{equation}
\begin{equation}\label{ABCD2}
  [A^{\epsilon}_i(z),B^{\epsilon'}_j(w)]=
  [A^{\epsilon}_i(z),C^{\epsilon'}_j(w)]=
  [B^{\epsilon}_i(z),C^{\epsilon'}_j(w)]=0\ \mathrm{for}\ i\ne j,
\end{equation}
\begin{equation}\label{ABCD5}
  [B^{\epsilon}_i(z),B^{\epsilon'}_i(w)]=[C^{\epsilon}_i(z),C^{\epsilon'}_i(w)]=
  [D^{\epsilon}_i(z),D^{\epsilon'}_i(w)]=0,
\end{equation}
\begin{equation}\label{ABCD6}
  (z-w)[B^{\epsilon'}_i(w),A^{\epsilon}_i(z)]_{\vv_i^{-1}}=
  (\vv_i-\vv_i^{-1})\left(zA^{\epsilon}_i(z)B^{\epsilon'}_i(w)-wA^{\epsilon'}_i(w)B^{\epsilon}_i(z)\right),
\end{equation}
\begin{equation}\label{ABCD7}
  (z-w)[A^{\epsilon}_i(z),C^{\epsilon'}_i(w)]_{\vv_i}=
  (\vv_i-\vv_i^{-1})\left(wC^{\epsilon'}_i(w)A^{\epsilon}_i(z)-zC^{\epsilon}_i(z)A^{\epsilon'}_i(w)\right),
\end{equation}
\begin{equation}\label{ABCD8}
  (z-w)[B^{\epsilon}_i(z),C^{\epsilon'}_i(w)]=
  (\vv_i-\vv_i^{-1})z\left(D^{\epsilon'}_i(w)A^{\epsilon}_i(z)-D^{\epsilon}_i(z)A^{\epsilon'}_i(w)\right),
\end{equation}
\begin{equation}\label{ABCD9}
  (z-w)[B^{\epsilon}_i(z),D^{\epsilon'}_i(w)]_{\vv_i}=
  (\vv_i-\vv_i^{-1})\left(wD^{\epsilon'}_i(w)B^{\epsilon}_i(z)-zD^{\epsilon}_i(z)B^{\epsilon'}_i(w)\right),
\end{equation}
\begin{equation}\label{ABCD10}
  (z-w)[D^{\epsilon'}_i(w),C^{\epsilon}_i(z)]_{\vv_i^{-1}}=
  (\vv_i-\vv_i^{-1})\left(zC^{\epsilon}_i(z)D^{\epsilon'}_i(w)-wC^{\epsilon'}_i(w)D^{\epsilon}_i(z)\right),
\end{equation}
\begin{equation}\label{ABCD11}
  (z-w)[A^{\epsilon}_i(z),D^{\epsilon'}_i(w)]=
  (\vv_i-\vv_i^{-1})\left(wC^{\epsilon'}_i(w)B^{\epsilon}_i(z)-zC^{\epsilon}_i(z)B^{\epsilon'}_i(w)\right),
\end{equation}
\begin{equation}\label{ABCD12}
  A^\epsilon_i(z)D^\epsilon_i(\vv_i^{-2}z)-\vv_i^{-1}B^\epsilon_i(z)C^\epsilon_i(\vv_i^{-2}z)=
  z^{\epsilon b^\epsilon_i}\cdot \prod_{j-i}\prod_{p=1}^{-c_{ji}} A^\epsilon_j(\vv_j^{-c_{ji}-2p}z),
\end{equation}
\begin{equation}\label{ABCD13}
\begin{split}
  & (z-\vv_i^{c_{ij}}w)B^\epsilon_i(z)B^{\epsilon'}_j(w)-(\vv_i^{c_{ij}}z-w)B^{\epsilon'}_j(w)B^\epsilon_i(z)=\\
  & zA^\epsilon_i(z)[\phi^+_iB^+_{i,0},B^{\epsilon'}_j(w)]_{\vv_i^{c_{ij}}}+
  wA^{\epsilon'}_j(w)[\phi^+_jB^+_{j,0}, B^\epsilon_i(z)]_{\vv_i^{c_{ij}}}
  \ \mathrm{for}\ i\ne j,
\end{split}
\end{equation}
\begin{equation}\label{ABCD14}
\begin{split}
  & (\vv_i^{c_{ij}}z-w)C^\epsilon_i(z)C^{\epsilon'}_j(w)-(z-\vv_i^{c_{ij}}w)C^{\epsilon'}_j(w)C^\epsilon_i(z)=\\
  & -[C^\epsilon_i(z),C^+_{j,1}\phi^+_j]_{\vv_i^{c_{ij}}}A^{\epsilon'}_j(w)-
  [C^{\epsilon'}_j(w),C^+_{i,1}\phi^+_i]_{\vv_i^{c_{ij}}}A^\epsilon_i(z)
  \ \mathrm{for}\ i\ne j,
\end{split}
\end{equation}
\begin{equation}\label{ABCD15}
\begin{split}
  & \underset{z_1,\ldots,z_{1-c_{ij}}} \Sym \left\{
  \prod_{a<b} (\vv_i z_a- \vv_i^{-1} z_b)(z_a-z_b)\cdot \right.\\
  & \left.\sum_{r=0}^{1-c_{ij}}(-1)^{r}{1-c_{ij}\brack r}_{\vv_i}
  B^{\epsilon_1}_i(z_1)\cdots B^{\epsilon_{r}}_i(z_{r})  B^{\epsilon'}_j(w)
  B^{\epsilon_{r+1}}_i(z_{r+1})\cdots B^{\epsilon_{1-c_{ij}}}_i(z_{1-c_{ij}})\right\}=0,
\end{split}
\end{equation}
\begin{equation}\label{ABCD16}
\begin{split}
  & \underset{z_1,\ldots,z_{1-c_{ij}}} \Sym \left\{
  \prod_{a<b} (\vv_i z_b- \vv_i^{-1} z_a)(z_b-z_a)\cdot \right.\\
  & \left.\sum_{r=0}^{1-c_{ij}}(-1)^{r}{1-c_{ij}\brack r}_{\vv_i}
  C^{\epsilon_1}_i(z_1)\cdots C^{\epsilon_r}_i(z_{r})  C^{\epsilon'}_j(w)
  C^{\epsilon_{r+1}}_i(z_{r+1})\cdots  C^{\epsilon_{1-c_{ij}}}_i(z_{1-c_{ij}})\right\}=0,
\end{split}
\end{equation}
for any $i,j\in I$ and $\epsilon,\epsilon',\epsilon_1,\ldots,\epsilon_{1-c_{ij}}\in \{\pm\}$.

\noindent
(b) Relations~(\ref{ABCD0}--\ref{ABCD16}) are the defining
relations. In other words, the associative $\BC(v)$-algebra
generated by
  $\{\phi^\pm_i, A^\pm_{i,\pm r},B^+_{i,r},B^-_{i,-r-1},C^+_{i,r+1},
   C^-_{i,-r},D^\pm_{i,\pm r\pm b^\pm_i}\}_{i\in I}^{r\in \BN}$
with the defining relations~(\ref{ABCD0}--\ref{ABCD16}) is
isomorphic to $\CU^\ad_{\mu^+,\mu^-}$.
\end{Thm}

We sketch the proof in Appendix~\ref{Proof of Relations bw ABCD}.
In the unshifted case, more precisely for $U^\ad_\vv(L\fg)$,
the above construction should be viewed as a $\vv$-version of
that of~\cite{gklo1}. In \emph{loc.cit.}, the authors
introduced analogous generating series
$\{A_i(u), B_i(u),C_i(u),D_i(u)\}_{i\in I}$ with coefficients in the
Yangian $Y(\fg)$ and stated (without a proof) the relations between them,
similar to~(\ref{ABCD1}--\ref{ABCD12}).\footnote{We note that the
relation $[D_i(u),D_i(v)]=0$ was missing in their list.} Meanwhile,
we note that adding rational analogues of~(\ref{ABCD13}--\ref{ABCD16})
to their list of relations, we get a complete list of the defining relations
among these generating series.

\begin{Rem}\label{Remark on ABCD}
\noindent
(a) For $\fg=\ssl_2$, relations~(\ref{ABCD1},~\ref{ABCD5}--\ref{ABCD11})
are equivalent to the RTT-relations (with the trigonometric
$R$-matrix of~(\ref{R-matrix})), see our proof
of~Theorem~\ref{Ding-Frenkel 2} below.

\noindent
(b) This construction can be adapted to the setting of
$\CU^\ssc_{\mu^+,\mu^-}$. First, we redefine the generating series
$A^\pm_i(z)=1+\sum_{r>0}A^\pm_{i,\pm r}z^{\mp r}$ which have to satisfy
\begin{equation}\label{A'}
   z^{\mp b^\pm_i}(\psi^\pm_{i,\mp b^\pm_i})^{-1}\psi^\pm_i(z)=
   \frac{\prod_{j-i}\prod_{p=1}^{-c_{ji}}A^\pm_j(\vv_j^{-c_{ji}-2p}z)}{A^\pm_i(z)A^\pm_i(\vv_i^{-2}z)}.
\end{equation}
Next, we define $B^\pm_i(z),C^\pm_i(z)$ via
formulas~(\ref{B}, \ref{C}). Finally, we define $D^\pm_i(z)$ via
\begin{equation}\label{D'}
  D^\pm_i(z):=A^\pm_i(z)\psi^\pm_i(z)+\vv^{\mp 1}_i(\vv_i-\vv_i^{-1})^2 f^\pm_i(z)A^\pm_i(z)e^\pm_i(z).
\end{equation}
The coefficients of these generating series together with
$\{(\psi^\epsilon_{i,-\epsilon b^\epsilon_i})^{\pm 1}\}_{i\in I}^{\epsilon=\pm}$
generate $\CU^\ssc_{\mu^+,\mu^-}$. For $\mu^+,\mu^-\in \Lambda^-$
one can write a complete list of the defining relations among these
generators, which look similar to~(\ref{ABCD1}--\ref{ABCD16}).
\end{Rem}


\section{Homomorphism to difference operators}
\label{Section with Homomorphism}

In this section, we construct homomorphisms from the shifted quantum
affine algebras to the algebras of difference operators.


\subsection{Homomorphism $\wt{\Phi}^{\lambda}_\mu$}
$\ $

Let $\mathrm{Dyn}(\fg)$ be the graph obtained from the Dynkin diagram of $\fg$ by replacing
all multiple edges by simple ones. We fix an orientation of $\mathrm{Dyn}(\fg)$ and we fix a
dominant coweight $\lambda\in \Lambda^+$ and a coweight $\mu\in \Lambda$,
such that $\lambda-\mu=\sum_{i\in I} a_i\alpha_i\ \mathrm{with}\ a_i\in \BN$.
We also fix a sequence
$\unl{\lambda}=(\omega_{i_1},\ldots,\omega_{i_N})$
of fundamental coweights, such that $\sum_{s=1}^N\omega_{i_s}=\lambda$.

Consider the associative $\BC[\vv^{\pm 1}]$-algebra $\hat{\CA}^\vv$ generated by
  $\{D_{i,r}^{\pm 1}, \sw_{i,r}^{\pm 1/2}\}_{i\in I}^{1\leq r\leq a_i}$
with the defining relations (for all $i,j\in I,\ 1\leq r\leq a_i,\ 1\leq s\leq a_j)$:
\begin{equation*}
    [D_{i,r},D_{j,s}]=[\sw^{1/2}_{i,r},\sw^{1/2}_{j,s}]=0,\
    D_{i,r}^{\pm 1}D_{i,r}^{\mp 1}=\sw_{i,r}^{\pm 1/2}\sw_{i,r}^{\mp 1/2}=1,\
    D_{i,r}\sw^{1/2}_{j,s}=\vv_i^{\delta_{ij}\delta_{rs}}\sw^{1/2}_{j,s}D_{i,r}.
\end{equation*}
Let $\wt{\CA}^\vv$ be the localization of $\hat{\CA}^\vv$ by the
multiplicative set generated by
  $\{\sw_{i,r}-\vv_i^m\sw_{i,s}\}_{i\in I, m\in \BZ}^{1\leq r\ne s\leq a_i}\\
\cup\{1-\vv^m\}_{m\in\BZ\setminus\{0\}}$
(which obviously satisfies Ore conditions). We also define their
$\BC(\vv)$-counterparts
  $\hat{\CA}^\vv_\fra:=\hat{\CA}^\vv\otimes_{\BC[\vv^{\pm 1}]} \BC(\vv)$
and
  $\wt{\CA}^\vv_\fra:=\wt{\CA}^\vv\otimes_{\BC[\vv^{\pm 1}]} \BC(\vv)$.

In what follows, we will work with the larger algebra
$\CU^\ad_{0,\mu}[\sz^{\pm 1}_1,\ldots,\sz^{\pm 1}_N]$,
which is obtained from
  $\CU^\ssc_{0,\mu}[\sz^{\pm 1}_1,\ldots,\sz^{\pm 1}_N]:=
   \CU^\ssc_{0,\mu}\otimes_{\BC(\vv)} \BC(\vv)[\sz^{\pm 1}_1,\ldots,\sz^{\pm 1}_N]$
by adding extra generators
  $\{(\phi^\epsilon_i)^{\pm 1}\}_{i\in I}^{\epsilon=\pm}$
satisfying relations~(\ref{U9}, \ref{U10}) with the only change:
\begin{equation*}
    \prod_{s:i_s=i} (-\vv_i\sz_s)^{\mp 1}\cdot (\psi^-_{i,\alphavee_i(\mu)})^{\pm 1}=
   (\phi^-_i)^{\pm 2}\cdot \prod_{j - i}(\phi^-_j)^{\pm c_{ji}}.
\end{equation*}
We will also work with the larger algebras
  $\wt{\CA}^\vv[\sz^{\pm 1}_1,\ldots,\sz^{\pm 1}_N]:=
   \wt{\CA}^\vv\otimes_{\BC[\vv^{\pm 1}]} \BC[\vv^{\pm 1}][\sz^{\pm 1}_1,\ldots,\sz^{\pm 1}_N]$
and
  $\wt{\CA}^\vv_\fra[\sz^{\pm 1}_1,\ldots,\sz^{\pm 1}_N]:=
   \wt{\CA}^\vv_\fra\otimes_{\BC(\vv)} \BC(\vv)[\sz^{\pm 1}_1,\ldots,\sz^{\pm 1}_N]$.

Define
\begin{equation*}
\begin{split}
    &  Z_i(z):=\prod_{1\leq s\leq N}^{i_s=i} \left(1-\frac{\vv_i \sz_s}{z}\right),\
    W_i(z):=\prod_{r=1}^{a_i} \left(1-\frac{\sw_{i,r}}{z}\right),\
    W_{i,r}(z):=\prod_{1\leq s\leq a_i}^{s\ne r} \left(1-\frac{\sw_{i,s}}{z}\right),\\
    & \hat{Z}_i(z):=\prod_{1\leq s\leq N}^{i_s=i}\left(1-\frac{z}{\vv_i \sz_s}\right),\
    \hat{W}_i(z):=\prod_{r=1}^{a_i}\left(1-\frac{z}{\sw_{i,r}}\right),\
    \hat{W}_{i,r}(z):=\prod_{1\leq s\leq a_i}^{s\ne r} \left(1-\frac{z}{\sw_{i,s}}\right).
\end{split}
\end{equation*}

The following is the key result of this section.

\begin{Thm}\label{Homomorphism}
There exists a unique $\BC(\vv)[\sz^{\pm 1}_1,\ldots,\sz^{\pm 1}_N]$-algebra homomorphism
\begin{equation*}
    \wt{\Phi}^{\unl\lambda}_\mu\colon
    \CU^\ad_{0,\mu}[\sz^{\pm 1}_1,\ldots,\sz^{\pm 1}_N]\longrightarrow
    \wt{\CA}^\vv_\fra[\sz^{\pm 1}_1,\ldots,\sz^{\pm 1}_N],
\end{equation*}
such that
 $$e_i(z)\mapsto \frac{-\vv_i}{1-\vv_i^2}
   \prod_{t=1}^{a_i}\sw_{i,t} \prod_{j\to i} \prod_{t=1}^{a_j} \sw_{j,t}^{c_{ji}/2}\cdot
   \sum_{r=1}^{a_i} \delta\left(\frac{\sw_{i,r}}{z}\right)\frac{Z_i(\sw_{i,r})}{W_{i,r}(\sw_{i,r})}
   \prod_{j\to i} \prod_{p=1}^{-c_{ji}} W_j(\vv_j^{-c_{ji}-2p}z)D_{i,r}^{-1},$$
 $$f_i(z)\mapsto \frac{1}{1-\vv_i^2}\prod_{j\leftarrow i}\prod_{t=1}^{a_j} \sw_{j,t}^{c_{ji}/2}\cdot
   \sum_{r=1}^{a_i} \delta\left(\frac{\vv_i^2\sw_{i,r}}{z}\right)\frac{1}{W_{i,r}(\sw_{i,r})}
   \prod_{j\leftarrow i}\prod_{p=1}^{-c_{ji}} W_j(\vv_j^{-c_{ji}-2p}z)D_{i,r},$$
 $$\psi^\pm_i(z)\mapsto \prod_{t=1}^{a_i}\sw_{i,t} \prod_{j - i} \prod_{t=1}^{a_j} \sw_{j,t}^{c_{ji}/2}\cdot
   \left(\frac{Z_i(z)}{W_i(z)W_i(\vv_i^{-2}z)}
   \prod_{j - i} \prod_{p=1}^{-c_{ji}} W_j(\vv_j^{-c_{ji}-2p}z)\right)^\pm,$$
 $$(\phi^+_i)^{\pm 1}\mapsto \prod_{t=1}^{a_i} \sw_{i,t}^{\pm 1/2},\
   (\phi^-_i)^{\pm 1}\mapsto (-\vv_i)^{\mp a_i}\prod_{t=1}^{a_i} \sw_{i,t}^{\mp 1/2}.$$
We write $\gamma(z)^\pm$ for the expansion of a rational function $\gamma(z)$ in $z^{\mp 1}$, respectively.
\end{Thm}

In the {\em unshifted case}, more precisely for $U_\vv(L\fg)$, this result
was stated (without a proof) in~\cite{gklo}. The above formulas simplify
for simply-laced $\fg$, in which case this result can be
viewed as a $\vv$-version of~\cite[Corollary B.17]{bfn}. We present
the proof in Appendix~\ref{Proof Homomorphism}.


\subsection{Homomorphism $\wt{\Phi}^{\lambda}_\mu$ in $ABC$ generators}
$\ $

Generalizing the construction of Section~\ref{ABCD generators}, we
define new Cartan generators $\{A^\pm_{i,\pm r}\}_{i\in I}^{r\geq 0}$
of $\CU^\ad_{0,\mu}[\sz^{\pm 1}_1,\ldots,\sz^{\pm 1}_N]$ via
  $$A^\pm_{i,0}:=(\phi^\pm_i)^{-1},$$
  $$\frac{\psi^+_i(z)}{Z_i(z)}=
    \frac{\prod_{j-i}\prod_{p=1}^{-c_{ji}}A^+_j(\vv_j^{-c_{ji}-2p}z)}{A^+_i(z)A^+_i(\vv_i^{-2}z)},$$
  $$\frac{z^{\alphavee_i(\mu)}\psi^-_i(z)}{\prod_{s:i_s=i} (-\vv_i\sz_s)\cdot \hat{Z}_i(z)}=
    \frac{\prod_{j-i}\prod_{p=1}^{-c_{ji}}A^-_j(\vv_j^{-c_{ji}-2p}z)}{A^-_i(z)A^-_i(\vv_i^{-2}z)},$$
where we set
  $A^\pm_i(z):=\sum_{r\geq 0} A^\pm_{i,\pm r}z^{\mp r}.$
We also define the generating series $B^\pm_i(z)$, $C^\pm_i(z)$, and
$D^\pm_i(z)$ via formulas~(\ref{B}),~(\ref{C}), and~(\ref{D}), respectively.

\begin{Lem}\label{Relations bw ABCD generalized}
For antidominant $\mu\in\Lambda^-$, the generating series
$A^\pm_i(z),B^\pm_i(z),C^\pm_i(z),D^\pm_i(z)$ satisfy
relations~(\ref{ABCD1}--\ref{ABCD11}).
\end{Lem}

\begin{proof}
Let $c$ be the determinant of the Cartan matrix of $\fg$. Choose unique
$\lambda_i^+(z)\in 1+z^{-1}\BC(\vv)[\sz_1^{\pm1},\ldots,\sz_N^{\pm1}][[z^{-1}]]$,
such that
  $Z_i(z)=\frac{\lambda^+_i(z)\lambda^+_i(\vv_i^{-2}z)}
   {\prod_{j-i}\prod_{p=1}^{-c_{ji}}\lambda^+_j(\vv_j^{-c_{ji}-2p}z)}.$
Also choose $\lambda^-_i(z)\in \BC(\vv^{1/c})[\sz_1^{\pm1/c},\ldots,\sz_N^{\pm1/c}][[z]]$,
such that
  $\hat{Z}_i(z)\cdot\prod_{s:i_s=i}(-\vv_i\sz_s)=
   \frac{\lambda^-_i(z)\lambda^-_i(\vv_i^{-2}z)}
   {\prod_{j-i}\prod_{p=1}^{-c_{ji}}\lambda^-_j(\vv_j^{-c_{ji}-2p}z)}.$

Then, the series $\lambda_i^\pm(z)^{-1}X^\pm_i(z)$ for $X=A,B,C,D$ are
those of Section~\ref{ABCD generators}. The result follows
from Theorem~\ref{Relations bw ABCD}(a) (compare with the proof
of~\cite[Proposition 5.5]{ktwwy}).
\end{proof}

\begin{Cor}\label{BC via A}
The following equalities hold in $\CU^\ad_{0,\mu}[\sz_1^{\pm1},\ldots,\sz_N^{\pm1}]$:
\begin{equation*}
\begin{split}
  & B^+_i(z)=[e_{i,0},A^+_i(z)]_{\vv_i^{-1}},\
    C^+_i(z)=[z^{-1}A^+_i(z),f_{i,1}]_{\vv_i^{-1}},\\
  & B^-_i(z)=[e_{i,-1},zA^-_i(z)]_{\vv_i},\
    C^-_i(z)=[A^-_i(z),f_{i,0}]_{\vv_i}.
\end{split}
\end{equation*}
\end{Cor}

\begin{proof}
The above formula for $B^+_i(z)$ (resp.\ $C^+_i(z)$) follows by evaluating
the terms of degree $1$ (resp.\ $0$) in $w$ in the
equality~(\ref{ABCD6}) (resp.~(\ref{ABCD7})) with $\epsilon=\epsilon'=+$.

The formulas for $B^-_i(z), C^-_i(z)$ are proved analogously.
\end{proof}

The following result is straightforward.

\begin{Prop}\label{Homom ABC}
The homomorphism $\wt{\Phi}^{\unl\lambda}_\mu$ maps the $ABC$ currents as follows:
  $$A^+_i(z)\mapsto \prod_{t=1}^{a_i} \sw_{i,t}^{-1/2}\cdot W_i(z),\
    A^-_i(z)\mapsto (-\vv_i)^{a_i}\prod_{t=1}^{a_i} \sw_{i,t}^{1/2}\cdot \hat{W}_i(z),$$
  $$B^+_i(z)\mapsto
    \prod_{t=1}^{a_i}\sw_{i,t}^{1/2}\prod_{j\to i} \prod_{t=1}^{a_j} \sw_{j,t}^{c_{ji}/2}\cdot
    \sum_{r=1}^{a_i} \frac{W_{i,r}(z)Z_i(\sw_{i,r})}{W_{i,r}(\sw_{i,r})}
    \prod_{j\to i} \prod_{p=1}^{-c_{ji}} W_j(\vv_j^{-c_{ji}-2p}\sw_{i,r})D_{i,r}^{-1},$$
  $$B^-_i(z)\mapsto
    -(-\vv_i)^{a_i}\prod_{t=1}^{a_i}\sw_{i,t}^{3/2}\prod_{j\to i} \prod_{t=1}^{a_j} \sw_{j,t}^{c_{ji}/2}\cdot
    \sum_{r=1}^{a_i} \frac{z\hat{W}_{i,r}(z)Z_i(\sw_{i,r})}{\sw_{i,r}W_{i,r}(\sw_{i,r})}
    \prod_{j\to i} \prod_{p=1}^{-c_{ji}} W_j(\vv_j^{-c_{ji}-2p}\sw_{i,r})D_{i,r}^{-1},$$
  $$C^+_i(z)\mapsto -\prod_{t=1}^{a_i} \sw_{i,t}^{-1/2}\prod_{j\leftarrow i}\prod_{t=1}^{a_j} \sw_{j,t}^{c_{ji}/2}\cdot
    \sum_{r=1}^{a_i} \frac{\sw_{i,r}W_{i,r}(z)}{zW_{i,r}(\sw_{i,r})}
    \prod_{j\leftarrow i}\prod_{p=1}^{-c_{ji}} W_j(\vv_j^{-c_{ji}-2p}\vv_i^2\sw_{i,r})D_{i,r},$$
  $$C^-_i(z)\mapsto (-\vv_i)^{a_i} \prod_{t=1}^{a_i} \sw_{i,t}^{1/2}
    \prod_{j\leftarrow i}\prod_{t=1}^{a_j} \sw_{j,t}^{c_{ji}/2}\cdot
    \sum_{r=1}^{a_i} \frac{\hat{W}_{i,r}(z)}{W_{i,r}(\sw_{i,r})}
    \prod_{j\leftarrow i}\prod_{p=1}^{-c_{ji}} W_j(\vv_j^{-c_{ji}-2p}\vv_i^2\sw_{i,r})D_{i,r}.$$
In particular, all these images belong to
  $\wt{\CA}^\vv[\sz^{\pm 1}_1,\ldots,\sz^{\pm 1}_N]\subset
   \wt{\CA}^\vv_\fra[\sz^{\pm 1}_1,\ldots,\sz^{\pm 1}_N]$.
\end{Prop}

\section{$K$-theoretic Coulomb branch}
\label{Coulomb}


\ssec{quiv}{Quiver gauge theories}
$\ $

We follow the notations and setup of~\cite[Appendix~A]{bfn}, so that
$(GL(V),\bN)$ is a quiver gauge theory. As in~Section~\ref{Section with Homomorphism},
we fix a sequence $(\omega_{i_1},\ldots,\omega_{i_N})$ of fundamental coweights
of $G$ which is assumed to be simply-laced for the current discussion.
We choose a basis $w_1,\ldots,w_N$ in $W=\bigoplus_{i\in I}W_i$ such that
$w_s\in W_{i_s}$. This defines a maximal torus $T_W\subset\prod_i GL(W_i)$,
and $K_{T_W}(\on{pt})=\BC[\sz_1^{\pm1},\ldots,\sz_N^{\pm1}]$.
We consider the (quantized)
$K$-theoretic Coulomb branch with flavor deformation
  $\CA^q=K^{(GL(V)\times T_W)_\CO\rtimes\BC^\times}(\CR_{GL(V),\bN})$
equipped with the convolution algebra structure as
in~\cite[Remark~3.9(3)]{bfna}. It is a $K_{\BC^\times\times T_W}(\on{pt})$-algebra;
we denote $K_{\BC^\times}(\on{pt})=\BC[q^{\pm1}]$. We will also need
$\vv=q^{1/2}$, the generator of the equivariant $K$-theory of a point with
respect to the two-fold cover $\wt\BC{}^\times\to\BC^\times$.
Recall that $GL(V)=\prod_{i\in I}GL(V_i)$. We will need its $2^I$-cover
$\wt{GL}(V)=\prod_{i\in I}\wt{GL}(V_i)$ where
  $\wt{GL}(V_i):=\{(g\in GL(V_i),\ y\in\BC^\times) : \det(g)=y^2\}$.
We consider the extended Coulomb branch
  $\CA^\vv:=K^{(\wt{GL}(V)\times T_W)_\CO\rtimes\wt\BC^\times}(\CR_{GL(V),\bN})=
   \CA^q\otimes_{K_{GL(V)\times\BC^\times}(\on{pt})}K_{\wt{GL}(V)\times\wt\BC^\times}(\on{pt})$.
It is equipped with an algebra structure as in~\refss{dicson}.

Recall from~\cite{bfn} that $w_{i,r}^*$ is the cocharacter of the Lie algebra of
$GL(V)=\prod GL(V_i)$, which is equal to 0 except at the vertex $i$,
and is $(0,\ldots,0,1,0,\ldots,0)$ at $i$. Here 1 is at the $r$-th
entry ($r=1,\ldots,a_i=\dim V_i$). We denote the corresponding
coordinates of $T_V$ and $T_V^\vee$ by $\sw_{i,r}$ and
$D_{i,r}\ (i\in I,\ 1\leq r\leq a_i)$. The roots are
$\sw_{i,r}\sw_{i,s}^{-1}\ (r\ne s)$. Furthermore,
$K^{(T_V\times T_W)_\CO\rtimes\BC^\times}(\CR_{T_V,0})$
with scalars extended by $\vv,\sw_{i,r}^{\pm1/2}$ is nothing but the algebra
  $\hat\CA^\vv[\sz_1^{\pm1},\ldots,\sz_N^{\pm1}]:=
   \hat\CA^\vv\otimes_{\BC[\vv^{\pm 1}]} \BC[\vv^{\pm 1}][\sz^{\pm 1}_1,\ldots,\sz^{\pm 1}_N]$,
where $\hat\CA^\vv$ was defined in~Section~\ref{Section with Homomorphism}.
We thus have an algebra embedding
\begin{equation*}
  \bz^*(\iota_*)^{-1}\colon
  \CA^\vv\hookrightarrow\wt\CA{}^\vv[\sz_1^{\pm1},\ldots,\sz_N^{\pm1}].
\end{equation*}

Let $\varpi_{i,n}$ be the $n$-th fundamental coweight of the factor
$GL(V_i)$, i.e.,\
  $w_{i,1}^*+\ldots+w_{i,n}^*=(1,\ldots,1,0,\ldots,0)$
where 1 appears $n$ times $(1\leq n\leq a_i)$. Then
$\Gr_{GL(V)}^{\varpi_{i,n}}$ is closed and isomorphic to the
Grassmannian $\Gr(V_i,n)$ of $n$-dimensional quotients of $V_i$. Let
$\CQ_i$ be the tautological rank $n$ quotient bundle on
$\Gr_{GL(V)}^{\varpi_{i,n}}$. Its pull-back to $\CR_{\varpi_{i,n}}$ is
also denoted by $\CQ_i$ for brevity. Let $\Lambda^p(\CQ_i)$ denote
the class of its $p$-th external power in $\CA^\vv$. More generally,
we can consider a class $f(\CQ_i)$ for a symmetric function $f$ in
$n$ variables so that $\Lambda^p(\CQ_i)$ corresponds to the $p$-th
elementary symmetric polynomial $e_p$.

Similarly, we consider $\varpi_{i,n}^*=-w_0\varpi_{i,n}$, where the corresponding
orbit $\Gr_{GL(V)}^{\varpi_{i,n}^*}$ is closed and isomorphic to the
Grassmannian $\Gr(n,V_i)$ of $n$-dimensional subspaces in $V_i$. Let
$\CS_i$ be the tautological rank $n$ subbundle on
$\Gr_{GL(V)}^{\varpi_{i,n}^*}$. Its pull-back to $\CR_{\varpi_{i,n}^*}$
is also denoted by $\CS_i$. Now similarly to~\cite[(A.3),(A.5)]{bfn},
cf.~\cite[Remark A.8]{bfn}, we obtain

\begin{equation}
    \label{A3}
        \bz^* (\iota_*)^{-1} \left(f(\CQ_i) \otimes \CO_{\CR_{\varpi_{i,n}}}\right)
    = 
    \sum_{\substack{J\subset \{1,\dots, a_i\} \\ \# J = n}}
    f(\sw_{i,J})
    \frac{\displaystyle\prod_{\substack{j\leftarrow i\\ r\in J}}
      \prod_{\substack{s=1 \\ (j,s)\neq (i,r)}}^{a_j}
      (1-\vv\sw_{i,r}\sw_{j,s}^{-1})}
    {\displaystyle\prod_{r\in J, s\notin J} (1-\sw_{i,s}\sw_{i,r}^{-1})}
    \prod_{r\in J} D_{i,r}
\end{equation}
(the appearance of $\vv$ is due to the convention before~\cite[Remark~2.1]{bfna});

\begin{equation}
    \label{A4}
\begin{split}
&        \bz^* (\iota_*)^{-1} \left(f(\CS_i)\otimes \CO_{\CR_{\varpi_{i,n}^*}}\right)= \\
&    \sum_{\substack{J\subset \{1,\dots, a_i\} \\ \# J = n}}
    f(\vv^{-2}\sw_{i,J})
    \prod_{\substack{r\in J\\ t:i_t=i}} (1-\vv\sz_t\sw_{i,r}^{-1})
    \frac{\displaystyle\prod_{\substack{j\rightarrow i\\ r\in J}}
      \prod_{\substack{s=1 \\ (j,s)\neq (i,r)}}^{a_j}
      (1-\vv\sw_{j,s}\sw_{i,r}^{-1})}
    {\displaystyle\prod_{r\in J, s\notin J} (1-\sw_{i,r}\sw_{i,s}^{-1})}
    \prod_{r\in J} D_{i,r}^{-1},
\end{split}
\end{equation}
where $f(\vv^{-2}\sw_{i,J})$ means that we substitute
$\{\vv^{-2}\sw_{i,r}\}_{r\in J}$ to $f$.

Also, for the vector bundles $\Omega_{\varpi_{i,1}}^p,\Omega_{\varpi_{i,1}^*}^p$ of
$p$-forms on $\Gr_{GL(V)}^{\varpi_{i,1}},\Gr_{GL(V)}^{\varpi_{i,1}^*}$ we obtain

\begin{equation}
    \label{A3'}
\begin{split}
&        \bz^* (\iota_*)^{-1}
    \left(\Omega_{\varpi_{i,1}}^p\otimes\CQ_i^{\otimes p'} \otimes \CO_{\CR_{\varpi_{i,1}}}\right)= \\
&    \sum_{1\leq r\leq a_i}
    \sw_{i,r}^{p'-p}\left(\sum_{\substack{J\subset\{1,\ldots,a_i\}\setminus\{r\}\\
    \# J=p}}\prod_{s\in J}\sw_{i,s}\right)
    \frac{\displaystyle\prod_{j\leftarrow i}
      \prod_{\substack{s=1 \\ (j,s)\neq (i,r)}}^{a_j}
      (1-\vv\sw_{i,r}\sw_{j,s}^{-1})}
    {\displaystyle\prod_{s\ne r} (1-\sw_{i,s}\sw_{i,r}^{-1})}
    D_{i,r},
\end{split}
\end{equation}

\begin{equation}
    \label{A4'}
\begin{split}
&        \bz^* (\iota_*)^{-1}
    \left(\Omega_{\varpi_{i,1}^*}^p\otimes\CS_i^{\otimes p'}\otimes \CO_{\CR_{\varpi_{i,1}^*}}\right)= \\
&    \sum_{1\leq r\leq a_i}
    \vv^{-2p'}\sw_{i,r}^{p'+p}\prod_{t:i_t=i} (1-\vv\sz_t\sw_{i,r}^{-1})
    \left(\sum_{\substack{J\subset\{1,\ldots,a_i\}\setminus\{r\}\\
    \# J=p}}\prod_{s\in J}\sw_{i,s}^{-1}\right)
    \frac{\displaystyle\prod_{j\rightarrow i}
      \prod_{\substack{s=1 \\ (j,s)\neq (i,r)}}^{a_j}
      (1-\vv\sw_{j,s}\sw_{i,r}^{-1})}
    {\displaystyle\prod_{s\ne r} (1-\sw_{i,r}\sw_{i,s}^{-1})}
    D_{i,r}^{-1}.
\end{split}
\end{equation}


\subsection{Homomorphism $\ol{\Phi}{}^{\lambda}_\mu$}
$\ $

We set $\CA^\vv_\fra:=\CA^\vv\otimes_{\BC[\vv^{\pm 1}]} \BC(\vv)$.
The key result of this section asserts that the homomorphism
$\wt{\Phi}^{\unl\lambda}_\mu$ of Theorem~\ref{Homomorphism}
factors through the above embedding
  $\bz^*(\iota_*)^{-1}\colon \CA^\vv_\fra
   \hookrightarrow\wt\CA^\vv_\fra[\sz^{\pm 1}_1,\ldots,\sz^{\pm 1}_N]$, similarly
to~\cite[Theorem~B.18]{bfn}.

\begin{Thm}\label{Factored Homomorphism}
There exists a unique $\BC(\vv)[\sz^{\pm 1}_1,\ldots,\sz^{\pm 1}_N]$-algebra homomorphism
\begin{equation*}
    \ol{\Phi}{}^{\unl\lambda}_\mu\colon
    \CU^\ad_{0,\mu}[\sz^{\pm 1}_1,\ldots,\sz^{\pm 1}_N]\longrightarrow
    \CA^\vv_\fra,
\end{equation*}
such that the following diagram commutes:

\setlength{\unitlength}{0.9cm}
\begin{picture}(4,3.2)
 \put(4.5,0.5){$\CA^\vv_\fra$}
 \put(3,2.5){$\CU^\ad_{0,\mu}[\sz^{\pm 1}_1,\ldots,\sz^{\pm 1}_N]$}
 \put(9.5,0.5){$\wt\CA^\vv_\fra[\sz^{\pm 1}_1,\ldots,\sz^{\pm 1}_N]$}

 \put(6,2.2){\vector (3,-1){3.5}}
 \put(5,2.2){\vector (0,-1){1.2}}
 \put(5.5,0.6){\vector (1,0){3.8}}

 \put(7.9,1.7){$\wt{\Phi}^{\unl\lambda}_\mu$}
 \put(4.2,1.5){$\ol{\Phi}{}^{\unl\lambda}_\mu$}
 \put(7,0.2){$\bz^*(\iota_*)^{-1}$}
\end{picture}

\noindent
Explicitly, $\ol{\Phi}{}^{\unl\lambda}_\mu$ maps the generators as follows:
  $$e_{i,r}\mapsto
    \frac{(-1)^{a_i}\vv}{1-\vv^2}\prod_{j\to i}\prod_{t=1}^{a_j}\sw_{j,t}^{-1/2}\cdot
    (\vv^2\CS_i)^{\otimes (r+a_i)}\otimes \CO_{\CR_{\varpi^*_{i,1}}},$$
  $$f_{i,r}\mapsto
    \frac{(-\vv)^{-\sum_{j\leftarrow i}a_j}}{1-\vv^2}\prod_{j\leftarrow i}\prod_{t=1}^{a_j}\sw_{j,t}^{1/2}\cdot
    \CQ_i^{\otimes (-\sum_{j\leftarrow i}a_j)}\otimes (\vv^2\CQ_i)^{\otimes r}\otimes \CO_{\CR_{\varpi_{i,1}}},$$
  $$A^+_{i,r}\mapsto (-1)^r \prod_{t=1}^{a_i} \sw_{i,t}^{-1/2}\cdot e_r(\{\sw_{i,t}\}_{t=1}^{a_i}),$$
  $$A^-_{i,-r}\mapsto (-1)^r(-\vv)^{a_i} \prod_{t=1}^{a_i} \sw_{i,t}^{1/2}\cdot e_r(\{\sw^{-1}_{i,t}\}_{t=1}^{a_i}),$$
  $$\phi^+_i\mapsto \prod_{t=1}^{a_i}\sw_{i,t}^{1/2},\
    \phi^-_i\mapsto (-\vv)^{-a_i}\prod_{t=1}^{a_i}\sw_{i,t}^{-1/2}.$$
\end{Thm}

\begin{proof}
For
  $X\in \{e_{i,r}, f_{i,r}, A^\pm_{i,\pm s}, \phi^\pm_i|i\in I, r\in \BZ, s\in \BN\}$
consider the assignment $X\mapsto \ol{\Phi}{}^{\unl\lambda}_\mu(X)$
with the right-hand side defined as above. Since
  $\bz^*(\iota_*)^{-1}\colon \CA^\vv_\fra
   \hookrightarrow\wt\CA^\vv_\fra[\sz^{\pm 1}_1,\ldots,\sz^{\pm 1}_N]$
is injective and
  $\wt{\Phi}^{\unl\lambda}_\mu\colon
   \CU^\ad_{0,\mu}[\sz^{\pm 1}_1,\ldots,\sz^{\pm 1}_N]
   \to \wt{\CA}^\vv_\fra[\sz^{\pm 1}_1,\ldots,\sz^{\pm 1}_N]$
is an algebra homomorphism, it suffices to check that
  $\bz^*(\iota_*)^{-1}\left(\ol{\Phi}{}^{\unl\lambda}_\mu(X)\right)=\wt{\Phi}^{\unl\lambda}_\mu(X)$
for $X$ as above. This is a straightforward verification based on
formulas~(\ref{A3}) and~(\ref{A4}).
\end{proof}

Combining Proposition~\ref{Homom ABC} with formulas~(\ref{A3'})
and~(\ref{A4'}), we immediately find the images of the generators
$\{B^+_{i,r},C^+_{i,r+1}\}_{i\in I}^{r\geq 0}$ under $\ol{\Phi}{}^{\unl\lambda}_\mu$.

\cor{Phi integral}
For $r\in \BN$, we have
  $$\ol{\Phi}{}^{\unl\lambda}_\mu(B^+_{i,r})=
    (-1)^{r+a_i+1}\vv^{2r}\prod_{t=1}^{a_i}\sw_{i,t}^{1/2}
    \prod_{j\rightarrow i}\prod_{t=1}^{a_j}\sw_{j,t}^{-1/2}\cdot
    \left(\Omega_{\varpi^*_{i,1}}^{a_i-1-r}\otimes\CS_i^{\otimes r}
    \otimes \CO_{\CR_{\varpi^*_{i,1}}}\right),$$
  $$\ol{\Phi}{}^{\unl\lambda}_\mu(C^+_{i,r+1})=
    (-1)^{r+1}(-\vv)^{-\sum_{j\leftarrow i} a_j}\prod_{t=1}^{a_i}\sw_{i,t}^{-1/2}
    \prod_{j\leftarrow i}\prod_{t=1}^{a_j}\sw_{j,t}^{1/2}\cdot
    \left(\Omega_{\varpi_{i,1}}^r\otimes\CQ_i^{\otimes (r+1-\sum_{j\leftarrow i} a_j)}
    \otimes \CO_{\CR_{\varpi_{i,1}}}\right).$$
In particular, the images of
  $\{A^+_{i,r}, B^+_{i,r}, C^+_{i,r+1},\phi^+_i\}_{i\in I}^{r\in\BN}$
under $\ol{\Phi}{}^{\unl\lambda}_\mu$ belong to $\CA^\vv\subset\CA^\vv_\fra$.
In fact, the images of
  $\{A^-_{i,-r}, B^-_{i,-r-1}, C^-_{i,-r},\phi^-_i\}_{i\in I}^{r\in\BN}$
under $\ol{\Phi}{}^{\unl\lambda}_\mu$ also belong to $\CA^\vv$.
\ecor

%

\begin{Rem}\label{BC in cohomological Coulomb}(A.~Weekes)
In the case of shifted Yangians, the images of the generating series
$B_i(z), C_i(z)$~\cite[Section 5.3]{ktwwy} in the quantized
(cohomological) Coulomb branch $\CA_\hbar$ under the homomorphism
$\ol{\Phi}{}^{\unl\lambda}_\mu$ of~\cite[Theorem~B.18]{bfn} are equal to
  $$\ol{\Phi}{}^{\unl\lambda}_\mu(B_i(z))=(-1)^{a_i}z^{-1}\cdot
    c(\wt{\CQ}_i,-z^{-1})\cap[\CR_{\varpi^*_{i,1}}],$$
  $$\ol{\Phi}{}^{\unl\lambda}_\mu(C_i(z))=(-1)^{\sum_{j\leftarrow i} a_j}z^{-1}\cdot
    c(\CS_i,-z^{-1})\cap[\CR_{\varpi_{i,1}}],$$
where $c(\mathcal{F},z)$ denotes the Chern polynomial of a vector
bundle $\mathcal{F}$. Here we view $\CQ_i,\CS_i$ as rank $n-1$
vector bundles on $\CR_{\varpi^*_{i,1}},\CR_{\varpi_{i,1}}$,
respectively, while $\wt{\CQ}_i$ denotes the vector bundle $\CQ_i$
with the equivariance structure twisted by $\hbar$.
\end{Rem}

\rem{exceptional}
Note that
  $\Gr_{GL(V)}^{\varpi_{i,1}}\simeq\BP^{a_i-1}\simeq\Gr_{GL(V)}^{\varpi_{i,1}^*}$,
and if we forget the equivariance, then up to sign,
$\ol{\Phi}{}^{\unl\lambda}_\mu(f_{i,r}),\ 1\leq r\leq a_i$, is the
collection of classes of pull-backs of the line bundles
  $\CO_{\BP^{a_i-1}}(1-\sum_{j\leftarrow i}a_j),\ldots,
   \CO_{\BP^{a_i-1}}(a_i-\sum_{j\leftarrow i}a_j)$,
while $\ol{\Phi}{}^{\unl\lambda}_\mu(C_{i,r}^+),\ 1\leq r\leq a_i$,
is the collection of classes of pull-backs of the vector bundles
$\Omega_{\BP^{a_i-1}}^{r-1}(r-\sum_{j\leftarrow i}a_j)$. These two
collections are the dual exceptional collections of vector bundles
on $\BP^{a_i-1}$ (more precisely, the former collection is left dual
to the latter one). In fact, this is the historically first example
of dual exceptional collections,~\cite{be}. Similarly, up to sign
and forgetting equivariance,
$\ol{\Phi}{}^{\unl\lambda}_\mu(e_{i,r}),\ 0\leq r<a_i$,
are the classes of the exceptional collection of line bundles right dual
to the exceptional collection of vector bundles whose classes are
$\ol{\Phi}{}^{\unl\lambda}_\mu(B_{i,r}^+),\ 0\leq r<a_i$.
\erem

\begin{Rem}\label{Relation to Dimofte et al}
An action of the quantized $K$-theoretic Coulomb branch
$\CA^\vv_\fra$ of the type $A$ quiver gauge theory on the localized
equivariant $K$-theory of parabolic Laumon spaces was constructed
in~\cite{bdghk}. Combining this construction with
Theorem~\ref{Factored Homomorphism}, we see that there should be a
natural action of
$\CU^\ad_{0,\mu}[\sz_1^{\pm1},\ldots,\sz_N^{\pm1}]$ (with $\fg=\ssl_n$)
on the aforementioned $K$-theory. We construct explicitly such an action of
$\CU^\ssc_{0,\mu}$ in Theorem~\ref{Action on K-theory of Laumon} by adapting
the arguments of~\cite{t} to the current setting (the adjoint version is achieved
by considering equivariant $K$-theory with respect to a larger torus).
\end{Rem}


\ssec{trun}{Truncated shifted quantum affine algebras}\label{trun}
\

We consider a $2$-sided ideal $\CI^{\unl\lambda}_\mu$ of
$\CU^\ad_{0,\mu}[\sz^{\pm 1}_1,\ldots,\sz^{\pm 1}_N]$ generated over
$\BC(\vv)[\sz^{\pm 1}_1,\ldots,\sz^{\pm 1}_N]$ by the following elements:
\begin{equation}\label{kernel 1}
  A^\pm_{i,\pm s}\ (s>a_i),\
  A^+_{i,0}A^+_{i,a_i}-(-1)^{a_i},\
  A^-_{i,0}A^-_{i,-a_i}-(-1)^{a_i}\vv_i^{2a_i},
\end{equation}
\begin{equation}\label{kernel 2}
  A^-_{i,-r}-\vv_i^{a_i}A^+_{i,a_i-r}\ (0\leq r\leq a_i).
\end{equation}

\defe{trunc}
$\CU^{\unl\lambda}_\mu:=\CU^\ad_{0,\mu}[\sz^{\pm 1}_1,\ldots,\sz^{\pm 1}_N]/\CI^{\unl\lambda}_\mu$
is called the truncated shifted quantum affine algebra.
\edefe


Note that the homomorphism
  $\wt\Phi{}^{\unl\lambda}_\mu\colon \CU^\ad_{0,\mu}[\sz^{\pm 1}_1,\ldots,\sz^{\pm 1}_N]
   \to\wt\CA_\fra^\vv[\sz^{\pm 1}_1,\ldots,\sz^{\pm 1}_N]$
factors through the same named homomorphism
  $\wt\Phi{}^{\unl\lambda}_\mu\colon
   \CU^{\unl\lambda}_\mu\to\wt\CA_\fra^\vv[\sz^{\pm 1}_1,\ldots,\sz^{\pm 1}_N]$,
due to Proposition~\ref{Homom ABC}.
Similarly to~\cite[Remark~B.21]{bfn}, we expect this homomorphism to be injective:

\begin{Conj}\label{description of truncated}
  $\wt\Phi{}^{\unl\lambda}_\mu\colon\CU^{\unl\lambda}_\mu
   \hookrightarrow\wt\CA_\fra^\vv[\sz^{\pm 1}_1,\ldots,\sz^{\pm 1}_N]$.
\end{Conj}

\begin{Rem}\label{first indication of conjecture}
As a first indication of the validity of this conjecture,
we note that the elements
  $\{B^+_{i,r},C^+_{i,r+1},B^-_{i,-r-1},C^-_{i,-r}\}_{i\in I}^{r\geq a_i}$
which belong to $\Ker(\wt{\Phi}^{\unl\lambda}_\mu)$ (due to
Proposition~\ref{Homom ABC}) also belong to $\CI^{\unl\lambda}_\mu$,
due to Corollary~\ref{BC via A} and relation~(\ref{U10}).
\end{Rem}

Moreover, we expect the following result:

\conj{Phi iso}
  $\ol{\Phi}{}^{\unl\lambda}_\mu\colon \CU^{\unl\lambda}_\mu\iso\CA_\fra^\vv$.
\econj


\ssec{yang}{Truncated shifted $\vv$-Yangians}
\

Recall that $\fg$ is assumed to be simply-laced.
Recall an explicit identification of the Drinfeld-Jimbo
and the new Drinfeld realizations of the standard quantum loop algebra $U_\vv(L\fg)$.
To this end, choose a decomposition of the highest root $\theta$ of $\fg$ into a sum
of simple roots $\theta=\alphavee_{i_1}+\alphavee_{i_2}+\ldots+\alphavee_{i_{h-1}}$
such that
  $\epsilon_k:=\langle\alpha_{i_{k+1}},\alphavee_{i_1}+\ldots+\alphavee_{i_k}\rangle
   \in \BZ_{<0}$ for any $1\leq k\leq h-2$
(here $h$ is the Coxeter number of $\fg$). We encode a choice of such a decomposition
by a sequence $\bi=(i_1,\ldots,i_{h-1})$.
Let $U^{\ddj}_\vv(L\fg)$ denote the Drinfeld-Jimbo quantum group of $\widehat{\fg}$
(affinization of $\fg$) with a trivial central charge, generated by
$\{E_i, F_i, K^{\pm 1}_i\}_{i\in \widetilde{I}}$
(here $\widetilde{I}=I\cup\{i_0\}$ is the vertex set of the extended Dynkin diagram),
see~\cite{lus}. The following result is due to~\cite{d} (proved in~\cite{jin}).

\begin{Thm}\label{identification of DJ and Dr}
There is a $\BC(\vv)$-algebra isomorphism $U^\ddj_\vv(L\fg)\iso U_\vv(L\fg)$, such that
\begin{equation*}
  E_i\mapsto e_{i,0},\ F_i\mapsto f_{i,0},\ K^{\pm 1}_i\mapsto \psi^\pm_{i,0}\ \mathrm{for}\ i\in I,
\end{equation*}
\begin{equation*}
  E_{i_0}\mapsto [f_{i_{h-1},0},[f_{i_{h-2},0},\cdots,[f_{i_2,0}, f_{i_1,1}]_{\vv^{\epsilon_1}}\cdots]_{\vv^{\epsilon_{h-3}}}]_{\vv^{\epsilon_{h-2}}}\cdot
  \psi^-_\theta,
\end{equation*}
\begin{equation*}
  F_{i_0}\mapsto (-\vv)^{-\epsilon}\psi^+_\theta\cdot [e_{i_{h-1},0},[e_{i_{h-2},0},\cdots,[e_{i_2,0}, e_{i_1,-1}]_{\vv^{\epsilon_1}}\cdots]_{\vv^{\epsilon_{h-3}}}]_{\vv^{\epsilon_{h-2}}},
\end{equation*}
\begin{equation*}
  K^{\pm}_{i_0}\mapsto \psi^\mp_\theta,
\end{equation*}
where
  $\psi^\pm_\theta:=\psi^\pm_{i_1,0}\cdots \psi^\pm_{i_{h-1},0},\ \epsilon:=\epsilon_1+\ldots+\epsilon_{h-2}$.
\end{Thm}

In particular, the image of the negative Drinfeld-Jimbo Borel subalgebra of $U^\ddj_\vv(L\fg)$
generated by $\{F_i,K^{\pm 1}_i\}_{i\in \wt{I}}$ under the above isomorphism is
the subalgebra $U^-_\vv$ of $U_\vv(L\fg)$, generated by
$\{f_{i,0},(\psi^-_{i,0})^{\pm 1}, F\}_{i\in I}$ with
  $F:=[e_{i_{h-1},0},[e_{i_{h-2},0},\cdots,[e_{i_2,0}, e_{i_1,-1}]_{\vv^{\epsilon_1}}\cdots]_{\vv^{\epsilon_{h-3}}}]_{\vv^{\epsilon_{h-2}}}$.
Motivated by this observation, we introduce the following definition.

\defe{trunc yang}
(a) Fix $\bi=(i_1,\ldots,i_{h-1})$ as above. The shifted $\vv$-Yangian
  $_\bi\!\CY^\vv_\mu[\sz^{\pm 1}_1,\ldots,\sz^{\pm 1}_N]$
is the $\BC(\vv)[\sz^{\pm 1}_1,\ldots,\sz^{\pm 1}_N]$-subalgebra of
$\CU^\ad_{0,\mu}[\sz^{\pm 1}_1,\ldots,\sz^{\pm 1}_N]$ generated by
$\{f_{i,0},(\psi^-_{i,b_i})^{\pm 1},\hat{F}\}_{i\in I}$, where
  $\hat{F}:=[e_{i_{h-1},b_{i_{h-1}}},[e_{i_{h-2},b_{i_{h-2}}},\cdots,[e_{i_2,b_{i_2}}, e_{i_1,b_{i_1}-1}]_{\vv^{\epsilon_1}}\cdots]_{\vv^{\epsilon_{h-3}}}]_{\vv^{\epsilon_{h-2}}}$
and $b_i:=\alphavee_i(\mu)$.

\noindent
(b) The truncated shifted $\vv$-Yangian $_\bi\!\CY^{\unl{\lambda}}_\mu$ is the
quotient of $_\bi\!\CY^\vv_\mu[\sz^{\pm 1}_1,\ldots,\sz^{\pm 1}_N]$ by the $2$-sided ideal
  $_\bi\CI^{\unl\lambda,+}_\mu:=\CI^{\unl\lambda}_\mu\cap{}_\bi\!\CY^\vv_\mu[\sz^{\pm 1}_1,\ldots,\sz^{\pm 1}_N]$.
%
%
\edefe

\begin{Rem}
For $\fg=\gl_n$ and $\mu=0$, our definition of the $\vv$-Yangian is consistent with that
of the \emph{quantum Yangian} $Y_q(\gl_n)$ of~\cite{nt} (in particular, independent of the choice of $\bi$).
The latter is defined via
the RTT presentation, see our discussion in Appendix~\ref{Proof of coproduct on loop sln},
and corresponds to the subalgebra generated by the coefficients of the matrix $T^-(z)$.
\end{Rem}

\conj{Phi iso Yang}
  $\ol{\Phi}{}^{\unl\lambda}_\mu\colon {}_\bi\!\CY^{\unl\lambda}_\mu\iso\CA_\fra^\vv$.
\econj


\ssec{integral forms}{Integral forms}
\

If we believe~\refco{Phi iso} and~\refco{Phi iso Yang}, we can
transfer the integral forms $\CA^\vv\subset\CA_\fra^\vv$ to the
truncated shifted quantum affine algebras and the truncated shifted
$\vv$-Yangians to obtain the $\BC[\vv^{\pm1}]$-subalgebras
$'\fU^{\unl\lambda}_\mu\subset\CU^{\unl\lambda}_\mu$ and
$'_\bi\!\fY^{\unl\lambda}_\mu\subset{}_\bi\!\CY^{\unl\lambda}_\mu$.
Finally, we define the integral form
$'\fU_{0,\mu}^{\on{ad}}\subset\CU_{0,\mu}^{\on{ad}}$ as an intersection of all the
preimages of ${'\fU^{\unl\lambda}_\mu}|_{{\sz_1=\ldots=\sz_N=1}}$ under projections
  $\CU_{0,\mu}^{\on{ad}}[\sz_1^{\pm1},\ldots,\sz_N^{\pm1}]\twoheadrightarrow
   \CU^{\unl\lambda}_\mu$ as $\unl{\lambda}$ varies, and $'_\bi\!\fY^\vv_\mu:={}
   '\fU_{0,\mu}^{\on{ad}}\cap{}_\bi\!\CY^\vv_\mu[\sz^{\pm 1}_1,\ldots,\sz^{\pm 1}_N]|_{{\sz_1=\ldots=\sz_N=1}}$.
Unfortunately, we can not define these integral forms by generators and relations
in general. In the case of $\ssl_2$ see~\refss{integral form}.

Recall that $*$ stands for the involution
$\mu\mapsto-w_0\mu$ of the coweight lattice $\Lambda$. Similarly
to~\cite[Remark~3.17]{bfn}, one can construct an isomorphism from the
non-quantized extended $K$-theoretic Coulomb branch
$\on{Spec}K^{(\wt{GL}(V)\times T_W)_\CO}(\CR_{GL(V),\bN})$ of~\refss{quiv} to the multiplicative
slice cover $\hW^{\unl\lambda^*}_{\mu^*}$ of~\refss{cover}. Its quantization is the subject
of the following

\conj{master conj}
(a) The shifted $\vv$-Yangian $'_\bi\!\fY^\vv_\mu$ is a quantization of
the thick multiplicative slice $^\dagger\CW_{\mu^*}$ of~\refss{thick}, that is
  ${'_\bi\!\fY^\vv_\mu}|_{{\vv=1}}\simeq\BC[{}^\dagger\CW_{\mu^*}]$.

\noindent
(b) The truncated shifted $\vv$-Yangian $'_\bi\!\fY^{\unl\lambda}_\mu$
and the truncated shifted quantum affine algebra
$'\fU^{\unl\lambda}_\mu$ are quantizations of the multiplicative
slice cover $\hW^{\unl\lambda^*}_{\mu^*}$ of~\refss{cover}, that is
  ${'_\bi\!\fY^{\unl\lambda}_\mu}|_{{\vv=1}}\simeq{}
   {'\fU^{\unl\lambda}_\mu}|_{{\vv=1}}\simeq\BC[\hW^{\unl\lambda^*}_{\mu^*}]$.
\econj


\ssec{braverman}{An example}
\

Let $\fg=\ssl_n,\ \mu=0,\ \unl{\lambda}=(\omega_1,\ldots,\omega_1)$
(the first fundamental coweight taken $n$ times). Note that the symmetric
group $\fS_n$ acts naturally on $'\fU^{\unl\lambda}_\mu$, permuting the parameters
$\sz_1,\ldots,\sz_n$. This action induces the one on the quotient algebra
$'\ol\fU{}^{\unl\lambda}_\mu$ by the relation $\sz_1\cdots\sz_n=1$.
Then we expect that the evaluation homomorphism
$U_\vv(L{\mathfrak{sl}}_n)\twoheadrightarrow U_\vv({\mathfrak{sl}}_n)$~\cite{j}
gives rise to an isomorphism
$({}'\ol\fU{}^{\unl\lambda}_\mu)^{\fS_n}\iso{}_\CA\!{\mathbf O}_{\on{loc}}$,
where $_\CA\!{\mathbf O}$ is the integral form of the quantum coordinate algebra of
$SL(N)$ introduced in~\cite[29.5.2]{lus}, and $_\CA\!{\mathbf O}_{\on{loc}}$
stands for its localization by inverting the quantum minors
$\{c_\nu\}_{\nu\in\Lambda^+}$, see~\cite[9.1.10]{jo}.





\sec{sl2}{Shifted quantum affine ${\mathfrak{sl}}_2$ and nil-DAHA for $GL(n)$}


\ssec{integral form}{Integral form}
\

In this section $\fg={\mathfrak{sl}}_2$, whence we
denote $A^\pm_{i,r}, B^\pm_{i,r}, C^\pm_{i,r},\phi^\pm_i$ simply by
$A^\pm_r,B^\pm_r,C^\pm_r,\phi^\pm$. The shift $\mu\in\Lambda=\BZ$ is an integer.
Furthermore, $\underline{\lambda}=(\omega_1,\ldots,\omega_1)$ (a collection of $N$ copies of
the fundamental coweight). The corresponding shifted quantum affine
algebra is $\CU_{0,\mu}^\ad[\sz_1^{\pm1},\ldots,\sz_N^{\pm1}]$. We define a
$\BC[\vv^{\pm1}]$-subalgebra
  $\fU^\ad_{0,\mu}[\sz_1^{\pm1},\ldots,\sz_N^{\pm1}]\subset
   \CU^\ad_{0,\mu}[\sz_1^{\pm1},\ldots,\sz_N^{\pm1}]$
generated by
  $\{A^\pm_{\pm r}, B^+_{r}, B^-_{-r-1}, C^+_{r+1}, C^-_{-r}, \phi^\pm\}_{r\in \BN}$
and its quotient algebra (an integral version of the truncated shifted quantum affine algebra)
  $\fU^{\unl{\lambda}}_\mu:=\fU^\ad_{0,\mu}[\sz_1^{\pm1},\ldots,\sz_N^{\pm1}]/
   (\CI^{\unl{\lambda}}_\mu\cap\fU^\ad_{0,\mu}[\sz_1^{\pm1},\ldots,\sz_N^{\pm1}])$.
Let $V=\BC^n,\ W=\BC^N$. According to~\refc{Phi integral}, the homomorphism
\begin{equation*}
    \ol{\Phi}{}^N_{N-2n}\colon \CU_{0,N-2n}^\ad[\sz_1^{\pm1},\ldots,\sz_N^{\pm1}]
    \longrightarrow\CA_\fra^\vv=K^{(\wt{GL}(V)\times T_W)_\CO\rtimes\wt\BC^\times}
    (\CR_{GL(V),\on{Hom}(W,V)})\otimes_{\BC[\vv^{\pm1}]}\BC(\vv)
\end{equation*}
takes
  $\fU^\ad_{0,N-2n}[\sz_1^{\pm1},\ldots,\sz_N^{\pm1}]\subset
   \CU^\ad_{0,N-2n}[\sz_1^{\pm1},\ldots,\sz_N^{\pm1}]$
to $\CA^\vv\subset\CA_\fra^\vv$. In particular, we have
  $\fU^\ad_{0,N-2n}[\sz_1^{\pm1},\ldots,\sz_N^{\pm1}]\subset{}
   '\fU^\ad_{0,N-2n}[\sz_1^{\pm1},\ldots,\sz_N^{\pm1}]$ (cf.~\refss{integral forms}).
We also define a $\BC[\vv^{\pm1}]$-subalgebra
  $\fY^\vv_{N-2n}[\sz_1^{\pm1},\ldots,\sz_N^{\pm1}]
   \subset\CY^\vv_{N-2n}[\sz_1^{\pm1},\ldots,\sz_N^{\pm1}]$
generated by
  $\{A^-_{-r}, B^-_{-r-1}, C^-_{-r}, \phi^-\}_{r\in \BN}$.
Furthermore, we define the shifted Borel $\vv$-Yangian
$\fY^\vv_{N-2n,-}[\sz_1^{\pm1},\ldots,\sz_N^{\pm1}]$ as the
$\BC[\vv^{\pm1}]$-subalgebra of $\fY^\vv_{N-2n}[\sz_1^{\pm1},\ldots,\sz_N^{\pm1}]$
generated by $\{A^-_{-r}, C^-_{-r}, \phi^-\}_{r\in \BN}$.
Finally, we have their truncated quotients
$\fY^{\unl{\lambda}}_{N-2n},\ \fY^{\unl{\lambda}}_{N-2n,-}$.
We expect that
\begin{equation*}
\begin{split}
  & \fU^\ad_{0,N-2n}[\sz_1^{\pm1},\ldots,\sz_N^{\pm1}]={}
    '\fU^\ad_{0,N-2n}[\sz_1^{\pm1},\ldots,\sz_N^{\pm1}],\
    \fY^\vv_{N-2n}[\sz_1^{\pm1},\ldots,\sz_N^{\pm1}]={}
    '\fY^\vv_{N-2n}[\sz_1^{\pm1},\ldots,\sz_N^{\pm1}],\\
  & \fU^{\unl{\lambda}}_{N-2n}={} '\fU^{\unl{\lambda}}_{N-2n},\
   \fY^{\unl{\lambda}}_{N-2n}={} '\fY^{\unl{\lambda}}_{N-2n}.
\end{split}
\end{equation*}

\begin{Conj}\label{Conj}
The natural homomorphisms induce isomorphisms
\begin{equation*}
  \fY^{\unl{\lambda}}_{N-2n}\iso \fU^{\unl{\lambda}}_{N-2n}\iso \CA^\vv.
\end{equation*}
\end{Conj}

From now on, we specialize to the case $N=0, \mu=-2n$. According
to~\refc{vassero}, the corresponding Coulomb branch
$\CA^\vv=K^{\wt{GL}(n,\CO)\rtimes\wt\BC^\times}(\Gr_{GL(n)})$
is nothing but the spherical extended nil-DAHA $\CHH_e^{\on{sph}}(GL(n))$.
We define $\BC[\vv^{\pm1}]_{\loc}$ inverting $(1-\vv^{2m}),\ m=1,2,\ldots,n$.
We extend the scalars to $\BC[\vv^{\pm1}]_{\loc}$ to obtain
\begin{equation*}
    \ol{\Phi}{}^0_{-2n,\on{loc}}\colon \fU_{0,-2n,\on{loc}}^\ad\longrightarrow
    K^{\wt{GL}(n,\CO)\rtimes\wt\BC^\times}_{\on{loc}}(\Gr_{GL(n)}).
\end{equation*}
The following theorem and Proposition~\ref{quarter} is a supportive
evidence in favor of Conjecture~\ref{Conj}.

\th{surj}\label{surj}
  $\ol{\Phi}{}^0_{-2n,\on{loc}}\colon \fU_{0,-2n,\on{loc}}^\ad\to
   K^{\wt{GL}(n,\CO)\rtimes\wt\BC^\times}_{\on{loc}}(\Gr_{GL(n)})$
   is surjective.\footnote{A stronger version of the theorem
     (over $\BZ[\vv^{\pm1}]$ as opposed to over $\BC[\vv^{\pm1}]_{\on{loc}}$)
     is proved independently
     in~\cite[Corollary~2.21,~Remark~2.22]{cw}.}
\eth

\begin{proof}
We must prove that
$K^{GL(n,\CO)\rtimes\BC^\times}_{\on{loc}}(\Gr_{GL(n)})$ is generated by
$K_{GL(n)}(\on{pt})=R(GL(n))$, and $\CO(a)_{\varpi_1},\CO(a)_{\varpi_1^*},\ a\in\BZ$.
Here $\varpi_1=(1,0,\ldots,0)$ denotes the first fundamental
coweight of $GL(n)$, and $\Gr^{\varpi_1}\simeq\BP^{n-1}$ is the
corresponding minuscule orbit, so that
$\Gr^{\varpi^*_1}\simeq\check\BP^{n-1}$. Finally, $\CQ$ is the
tautological quotient bundle on $\Gr^{\varpi_1}$, isomorphic to the
ample line bundle $\CO(1)$ on $\BP^{n-1}$, and $\CO(a)_{\varpi_1}$
stands for $\CQ^{\otimes a}$. Similarly, $\CS$ is the tautological
line subbundle on $\Gr^{\varpi^*_1}$ isomorphic to $\CO(-1)$ on
$\check\BP^{n-1}$, and $\CO(a)_{\varpi^*_1}$ stands for
$\CS^{\otimes-a}$. Note that $\CO(1)_{\varpi_1},\CO(1)_{\varpi_1^*}$ are isomorphic
to the restrictions of the determinant line bundle on $\Gr_{GL(n)}$.

Given an arbitrary sequence $\nu_1,\ldots,\nu_N$ with
  $\nu_i\in\{\varpi_1,\ldots,\varpi_{n},\varpi_1^*,\ldots,\varpi_{n}^*\}$,
the equivariant $K$-theory of the iterated convolution diagram
  $K^{GL(n,\CO)\rtimes\BC^\times}(\Gr^{\nu_1}\widetilde{\times}\ldots\widetilde{\times} \Gr^{\nu_N})$
is isomorphic to
  $K^{GL(n,\CO)\rtimes\BC^\times}(\Gr^{\nu_1})\otimes_{K_{GL(n,\CO)\rtimes\BC^\times}(\on{pt})}\cdots
   \otimes_{K_{GL(n,\CO)\rtimes\BC^\times}(\on{pt})}K^{GL(n,\CO)\rtimes\BC^\times}(\Gr^{\nu_N}).$
By the projection formula and rationality of singularities of
$\ol\Gr{}^{\nu_1+\ldots+\nu_N}$, the convolution pushforward morphism
\begin{equation*}
   m_*\colon
   K^{GL(n,\CO)\rtimes\BC^\times}(\Gr^{\nu_1}\widetilde{\times}\ldots\widetilde{\times} \Gr^{\nu_N})\longrightarrow
   K^{GL(n,\CO)\rtimes\BC^\times}(\ol\Gr{}^{\nu_1+\ldots+\nu_N})
\end{equation*}
is surjective. Hence in order to prove the surjectivity statement of
the theorem, it suffices to express
  $K^{GL(n,\CO)\rtimes\BC^\times}_{\on{loc}}(\Gr^\nu),\ \nu\in\{\varpi_1,
   \ldots, \varpi_{n}, \varpi^*_1, \ldots,\varpi^*_{n}\}$,
in terms of $\CO(a)_{\varpi_1},\CO(a)_{\varpi_1^*},\ a\in\BZ$, and $K_{GL(n)}(\on{pt})$.
We will consider $\nu=\varpi_m,\ 1\leq m\leq n$, the case of $\varpi^*_m$ being
similar. Note that $\CO_{\varpi_n}$ is the
structure sheaf of a point $GL(n,\CO)$-orbit corresponding to the
coweight $(1,\ldots,1)$. We argue by induction in $m$.

For $\nu$ as above, the Picard group of $\Gr^\nu$ is $\BZ$, and we denote the
ample generator by $\CO(1)_\nu$. It is isomorphic to the restriction of
the determinant line bundle on $\Gr_{GL(n)}$. We start with an
explicit expression for $\CO_{\varpi_m}:=\CO_{\Gr^{\varpi_m}},\ 1\leq m\leq n$,
in terms of $\CO(a)_{\varpi_1},\ a\in\BZ$. Recall that
  $\ol{\Phi}{}^0_{-2n}(f_r)=\frac{\vv^{2r}}{1-\vv^2}\CO(r)_{\varpi_1}$
and
  $\ol{\Phi}{}^0_{-2n}(e_r)=\frac{(-1)^n\vv^{2r+2n+1}}{1-\vv^2}\CO(-r-n)_{\varpi^*_1}$.
We denote $\on{ad}_x^{\vv^{r}}y:=[x,y]_{\vv^r}=xy-\vv^ryx$.

\prop{11.3} \label{11.3}
For any $1\leq m\leq n$, we have
\eq{8.2}
  \CO_{\varpi_m}=(-1)^{\frac{m(m-1)}{2}}(1-\vv^2)\ol{\Phi}{}^0_{-2n}(\on{ad}_{f_{1-m}}^{\vv^{2m}}
  \on{ad}_{f_{3-m}}^{\vv^{2(m-1)}}\cdots\on{ad}_{f_{m-3}}^{\vv^4}f_{m-1}),
\end{equation}

\eq{8.3}
  \CO_{\varpi^*_m}=(-1)^{nm+\frac{m(m+1)}{2}+1}\vv^{m^2-2}(1-\vv^2)\ol{\Phi}{}^0_{-2n}
  (\on{ad}_{e_{-n+1-m}}^{\vv^{-2m}}
  \on{ad}_{e_{-n+3-m}}^{\vv^{-2(m-1)}}\cdots\on{ad}_{e_{-n+m-3}}^{\vv^{-4}}e_{-n+m-1}).
\end{equation}
\eprop

\prf
We prove~\refe{8.2}; the proof of~\refe{8.3} is similar. We will
compare the images of the LHS and the RHS in
$\wt\CA^\vv_\fra$. According to~(\ref{A3}), the image of the LHS equals
\eq{lhs}
  \sum_{\# J=m}\prod_{r\in J}^{s\not\in J}(1-\sw_s\sw_r^{-1})^{-1}\prod_{r\in J}D_r.
\end{equation}
Here $J\subset\{1,\ldots,n\}$ is a subset of cardinality $m$. Let us
denote the iterated $\vv$-commutator
  $\on{ad}_{f_{1-m}}^{\vv^{2m}}\on{ad}_{f_{3-m}}^{\vv^{2(m-1)}}\cdots\on{ad}_{f_{m-3}}^{\vv^4}f_{m-1}$
by $F_m$. We want to prove
\eq{rhs}
  \wt\Phi^0_{-2n}(F_m)=(-1)^{\frac{m(m-1)}{2}}(1-\vv^2)^{-1}\cdot
  \sum_{\# J=m}\prod_{r\in J}^{s\not\in J}(1-\sw_s\sw_r^{-1})^{-1}\prod_{r\in J}D_r.
\end{equation}
The proof proceeds by induction in $m$. So we assume~\refe{rhs}
known for an integer $k<n$, and want to deduce~\refe{rhs} for
$m=k+1$. We introduce a ``shifted'' $\vv$-commutator
  $F'_k:=\on{ad}^{\vv^{2k}}_{f_{2-k}}\on{ad}^{\vv^{2(k-1)}}_{f_{4-k}}\cdots\on{ad}^{\vv^4}_{f_{k-2}}f_k$.
Then
\begin{equation*}
    \wt\Phi^0_{-2n}(F'_k)=(-1)^{\frac{k(k-1)}{2}}(1-\vv^2)^{-1}\vv^{2k}\cdot\sum_{\# J=k}
    \prod_{r\in J}\sw_r\prod_{r\in J}^{s\not\in J}\left(1-\frac{\sw_s}{\sw_r}\right)^{-1}\prod_{r\in J}D_r.
\end{equation*}
Now
\begin{equation*}
\begin{split}
  & \wt\Phi^0_{-2n}(F_{k+1})=\wt\Phi^0_{-2n}([f_{-k},F'_k]_{\vv^{2(k+1)}})=
    [\wt\Phi^0_{-2n}(f_{-k}),\wt\Phi^0_{-2n}(F'_k)]_{\vv^{2(k+1)}}=\\
  & (-1)^{\frac{k(k-1)}{2}}(1-\vv^2)^{-2}\vv^{2k}\cdot
    \left[\sum_{p=1}^n\frac{(\vv^2\sw_p)^{-k}}{\prod\limits_{t\ne p}
    \left(1-\frac{\sw_t}{\sw_p}\right)}D_p, \sum_{\# J=k}\prod_{r\in J}\sw_r\prod_{r\in J}^{s\not\in J}
    \left(1-\frac{\sw_s}{\sw_r}\right)^{-1}\prod_{r\in J}D_r\right]_{\vv^{2(k+1)}}.
\end{split}
\end{equation*}
First we check that the summands corresponding to $p\in J$ vanish.
Due to the symmetry reasons, we may assume $p=1,\ J=\{1,2,\ldots,k\}$. Then
\begin{equation*}
\begin{split}
  & \left[\frac{(\vv^2\sw_1)^{-k}}{\prod\limits_{t>1}\left(1-\frac{\sw_t}{\sw_1}\right)}D_1,
    \prod\limits_{r=1}^k\sw_r\prod\limits_{r\leq k}^{s>k}
    \left(1-\frac{\sw_s}{\sw_r}\right)^{-1}D_1\cdots D_k\right]_{\vv^{2(k+1)}}=\\
  & \left[\frac{(\vv^2\sw_1)^{-k}}{\prod\limits_{t>k}\left(1-\frac{\sw_t}{\sw_1}\right)
    \prod\limits_{1<r\leq k}\left(1-\frac{\sw_r}{\sw_1}\right)}D_1,
    \frac{\sw_1\cdots\sw_k}{\prod\limits_{s>k}\left(1-\frac{\sw_s}{\sw_1}\right)
    \prod\limits_{1<r\leq k}^{s>k}\left(1-\frac{\sw_s}{\sw_r}\right)}
    D_1\cdots D_k\right]_{\vv^{2(k+1)}}=\\
  & \left(\frac{(\vv^2\sw_1)^{-k}\vv^2\sw_1\cdots\sw_k}
    {\prod\limits_{t>k}\left(1-\frac{\sw_t}{\sw_1}\right)
    \prod\limits_{1<r\leq k}\left(1-\frac{\sw_r}{\sw_1}\right)
    \prod\limits_{s>k}\left(1-\vv^{-2}\frac{\sw_s}{\sw_1}\right)
    \prod\limits_{1<r\leq k}^{s>k}\left(1-\frac{\sw_s}{\sw_r}\right)}-\right.\\
  & \left.-\frac{\vv^{2(k+1)}\sw_1\cdots\sw_k(\vv^2\sw_1)^{-k}\vv^{-2k}}
    {\prod\limits_{s>k}\left(1-\frac{\sw_s}{\sw_1}\right)
    \prod\limits_{1<r\leq k}^{s>k}\left(1-\frac{\sw_s}{\sw_r}\right)
    \prod\limits_{t>k}\left(1-\vv^{-2}\frac{\sw_t}{\sw_1}\right)
    \prod\limits_{1<r\leq k}\left(1-\frac{\sw_r}{\sw_1}\right)}\right)D_1^2D_2\cdots D_k=0.
\end{split}
\end{equation*}
Therefore,
\begin{equation*}
\begin{split}
  & (-1)^{\frac{k(k-1)}{2}}(1-\vv^2)^2\wt\Phi^0_{-2n}(F_{k+1})=\\
  & \sum_{\# J=k}^{p\not\in J}\left[\frac{(\vv^2\sw_p)^{-k}}{\prod_{t\ne p}
    \left(1-\frac{\sw_t}{\sw_p}\right)}D_p, \vv^{2k}\prod_{r\in J}\sw_r\prod_{r\in J}^{s\not\in J}
    \left(1-\frac{\sw_s}{\sw_r}\right)^{-1}\prod_{r\in J}D_r\right]_{\vv^{2(k+1)}}.
\end{split}
\end{equation*}
We expand this combination of $\vv^{2(k+1)}$-commutators as a sum
$\sum\limits_{\# J=k+1}\phi_J(\sw_1,\ldots,\sw_n)\prod\limits_{r\in J}D_r.$ For the
symmetry reasons, it suffices to calculate the rational function
$\phi_J$ for a single $J=\{1,\ldots,k+1\}$. We have
  $$\phi_J(\sw_1,\ldots,\sw_n)D_1\cdots D_{k+1}=$$
  $$\sum_{r=1}^{k+1}\left[\frac{(\vv^2\sw_r)^{-k}}{\prod\limits_{t\ne r}
    \left(1-\frac{\sw_t}{\sw_r}\right)}D_r,
    \frac{\vv^{2k}\sw_1\cdots\widehat\sw_r\cdots\sw_{k+1}}{\prod\limits_{r\ne p\leq k+1}^{t>k+1}
    \left(1-\frac{\sw_t}{\sw_p}\right) \prod\limits_{r\ne p\leq k+1}\left(1-\frac{\sw_r}{\sw_p}\right)}
    D_1\cdots\widehat{D_r}\cdots D_{k+1}\right]_{\vv^{2(k+1)}}=$$
  $$\sum_{r=1}^{k+1}\left(\frac{\sw_r^{-k}\sw_1\cdots\widehat\sw_r\cdots\sw_{k+1}}
    {\prod\limits_{t>k+1}\left(1-\frac{\sw_t}{\sw_r}\right)
    \prod\limits_{r\ne p\leq k+1}\left(1-\frac{\sw_p}{\sw_r}\right)
    \prod\limits_{r\ne p\leq k+1}^{t>k+1}\left(1-\frac{\sw_t}{\sw_p}\right) \prod\limits_{r\ne p\leq k+1}
    \left(1-\frac{\vv^2\sw_r}{\sw_p}\right)}-\right.$$
  $$\left.-\frac{\vv^{2(k+1)}\sw_r^{-k}\sw_1\cdots\widehat\sw_r\cdots\sw_{k+1}}
    {\prod\limits_{r\ne p\leq k+1}^{t>k+1}\left(1-\frac{\sw_t}{\sw_p}\right) \prod\limits_{r\ne p\leq k+1}
    \left(1-\frac{\sw_r}{\sw_p}\right) \prod\limits_{t>k+1}\left(1-\frac{\sw_t}{\sw_r}\right) \prod\limits_{r\ne p\leq k+1}
    \left(1-\frac{\vv^2\sw_p}{\sw_r}\right)}\right) D_1\cdots D_{k+1}=$$
  $$-\vv^{2(k+1)}\sw_1\cdots\sw_{k+1}\prod_{r\leq k+1}^{t>k+1}\left(1-\frac{\sw_t}{\sw_r}\right)^{-1}\times$$
  $$\sum_{r=1}^{k+1}\left(\frac{\sw_r^{-k-1}}{\prod\limits_{r\ne p\leq k+1}
    \left(1-\frac{\sw_r}{\sw_p}\right)\left(1-\frac{\vv^2\sw_p}{\sw_r}\right)}-\frac{\vv^{-2(k+1)}\sw_r^{-k-1}}{\prod\limits_{r\ne p\leq k+1}
    \left(1-\frac{\sw_p}{\sw_r}\right)\left(1-\frac{\vv^2\sw_r}{\sw_p}\right)}\right)
    D_1\cdots D_{k+1}.$$
This is equal to the following expression, by~\refl{jpg} below:
\begin{equation*}
\begin{split}
  & -\vv^{2(k+1)}\sw_1\cdots\sw_{k+1}
    \prod_{r\leq k+1}^{t>k+1}\left(1-\frac{\sw_t}{\sw_r}\right)^{-1}
    \frac{(-1)^k(\vv^2-1)}{\vv^{2(k+1)}\prod\limits_{r\leq k+1}\sw_r}D_1\cdots D_{k+1}=\\
  & (-1)^k(1-\vv^2)\prod_{r\leq k+1}^{t>k+1}\left(1-\frac{\sw_t}{\sw_r}\right)^{-1}D_1\cdots D_{k+1}.
\end{split}
\end{equation*}
We conclude that
\begin{equation*}
    \wt\Phi^0_{-2n}(F_{k+1})=(-1)^{\frac{k(k+1)}{2}}(1-\vv^2)^{-1}\cdot
    \sum_{\# J=k+1}\prod_{r\in J}^{s\not\in J}(1-\sw_s\sw_r^{-1})^{-1}\prod_{r\in J}D_r,
\end{equation*}
and~\refe{rhs} is proved. It remains to check
\lem{jpg}\label{alg.identity}
We have
\begin{equation*}
    \sum^{k+1}_{r=1}\left(\frac{\sw_r^{-k-1}}{\prod\limits_{s\ne r}(1-\sw_r/\sw_s)
    (1-\vv^2\sw_s/\sw_r)}-\frac{\vv^{-2(k+1)}\sw_r^{-k-1}}{\prod\limits_{s\ne r}(1-\sw_s/\sw_r)
    (1-\vv^2\sw_r/\sw_s)}\right)=\frac{(-1)^k(\vv^2-1)}{\vv^{2(k+1)}\prod\limits^{k+1}_{r=1}\sw_r}.
\end{equation*}
\elem

\prf
The LHS is a degree $-k-1$ rational function of
$\sw_1,\ldots,\sw_{k+1}$ with poles at the hyperplanes given by equations
  $\sw_r-\sw_s,\ \sw_r-\vv^2\sw_s,\ \sw_r\ (1\leq r\ne s\leq k+1)$.
One can check
  $\on{Res}_{\sw_r-\sw_s}\on{LHS}=\on{Res}_{\sw_r-\vv^2\sw_s}\on{LHS}=0$,
so that $\on{LHS}=f\cdot\prod_{1\leq r\leq k+1}\sw_r^{-1}$ for a
rational function $f\in\BC(\vv)$. To compute $f$, we specialize
$\sw_1\mapsto0$ in the equality
\begin{equation*}
    f=\prod_{t=1}^{k+1}\sw_t\cdot\sum_{r=1}^{k+1}\left(\frac{\prod\limits_{s\ne r}\sw_s}
    {\prod\limits_{s\ne r}(\sw_s-\sw_r)(\sw_r-\vv^2\sw_s)}\cdot\frac{1}{\sw_r}-\frac{\vv^{-2(k+1)}\prod\limits_{s\ne r}\sw_s}
    {\prod\limits_{s\ne r}(\sw_r-\sw_s)(\sw_s-\vv^2\sw_r)}\cdot\frac{1}{\sw_r}\right).
\end{equation*}
The only summands surviving under this specialization correspond to
$r=1$, and so we get
\begin{equation*}
    f=\prod_{t=2}^{k+1}\sw_t\cdot\left(\frac{\prod\limits_{s=2}^{k+1}\sw_s}{(-\vv^2)^k\cdot
    \prod\limits_{s=2}^{k+1}\sw_s^2}-\frac{\vv^{-2(k+1)}\cdot\prod\limits_{s=2}^{k+1}\sw_s}{(-1)^k\cdot
    \prod\limits_{s=2}^{k+1}\sw_s^2}\right)=
    (-1)^k(\vv^{-2k}-\vv^{-2(k+1)}).
\end{equation*}
 The lemma is proved.
\epr
 The proposition is proved.
\epr

Returning to the proof of~\reft{surj}, we need to prove that
$K_{\on{loc}}^{GL(n,\CO)\rtimes\BC^\times}(\Gr^{\varpi_m})$
lies in the image $\ol{\Phi}{}_{-2n,\on{loc}}^0(\fU^\ad_{0,-2n,\on{loc}})$ for $1\leq m\leq n$.
We know that the class of the structure sheaf
$\CO_{\varpi_m}\in K_{\on{loc}}^{GL(n,\CO)\rtimes\BC^\times}(\Gr^{\varpi_m})$
lies in $\ol{\Phi}{}_{-2n,\on{loc}}^0(\fU^\ad_{0,-2n,\on{loc}})$. It is also known that
$K^{GL(n,\CO)\rtimes\BC^\times}(\Gr^{\varpi_m})$ as a left
$K_{GL(n,\CO)\rtimes\BC^\times}(\on{pt})$-module is generated by the classes
$\Sigma^\lambda(\CQ)$ where $\CQ$ is the tautological quotient bundle on
$\Gr^{\varpi_m}\simeq\Gr(m,n)$, and $\Sigma^\lambda$ is the polynomial Schur
functor corresponding to a Young diagram $\lambda$ with $\leq m$ rows
(in fact, it is enough to consider $\lambda$'s with $\leq n-m$ columns).
Given such $\lambda$, it suffices to check that
  $\on{Sym}\left(\sw_1^{\lambda_1}\cdots\sw_m^{\lambda_m}\prod\limits_{r\leq m}^{s>m}
   \left(1-\frac{\sw_s}{\sw_r}\right)^{-1}D_1\cdots D_m\right)$
lies in $\wt\Phi{}_{-2n,\on{loc}}^0(\fU^\ad_{0,-2n,\on{loc}})$
(here $\on{Sym}$ stands for the
symmetrization with respect to the symmetric group $\fS_n$).
More generally, for a Young diagram $\mu$ with $\leq n$ rows we will show that
  $\on{Sym}\left(\sw_1^{\mu_1}\cdots\sw_n^{\mu_n}\cdot\prod\limits_{r\leq m}^{s>m}
   \left(1-\frac{\sw_s}{\sw_r}\right)^{-1}D_1\cdots D_m\right)$
lies in $\wt\Phi{}_{-2n,\on{loc}}^0(\fU^\ad_{0,-2n,\on{loc}})$.
To this end, we use the {\em right}
multiplication by $K_{GL(n,\CO)\rtimes\BC^\times}(\on{pt})$. It suffices to check
that the ${K_{GL(n,\CO)\rtimes\BC^\times}(\on{pt})}_{\on{loc}}$-bimodule generated by
$X_{1,m}:=\on{Sym}\left(\prod\limits_{r\leq m}^{s>m}
\left(1-\frac{\sw_s}{\sw_r}\right)^{-1}D_1\cdots D_m\right)$
contains elements
  $X_{F,m}:=\on{Sym}\left(F\prod\limits_{r\leq m}^{s>m}
   \left(1-\frac{\sw_s}{\sw_r}\right)^{-1} D_1\cdots D_m\right)$
for any polynomial $F\in \BC[\sw_1,\ldots,\sw_n]$. We can assume that
$F\in \BC[\sw_1,\ldots,\sw_n]^{\fS_m\times \fS_{n-m}},$
where the symmetric groups act by permuting $\{\sw_r,\ 1\leq r\leq m\}$ and
$\{\sw_s,\ m+1\leq s\leq n\}$. Note that
  $\BC[\sw_1,\ldots,\sw_n]^{\fS_m\times \fS_{n-m}}$
is generated by $\BC[\sw_1,\ldots,\sw_m]^{\fS_m}$ as a left
  $\BC[\sw_1,\ldots,\sw_n]^{\fS_n}$-module.
Hence, it suffices to treat the case
$F\in \BC[\sw_1,\ldots,\sw_m]^{\fS_m}=\BC[p_1,\ldots,p_m]$,
where $p_k:=\sum_{r=1}^m \sw_r^k$. The latter case follows from the equality
\begin{equation*}
  \left[\sum_{r=1}^n\sw_r^k,X_{F,m}\right]=(1-\vv^{2k})X_{Fp_k,m}
\end{equation*}
for $F\in\BC[\sw_1,\ldots,\sw_m]^{\fS_m}$.

 The theorem is proved.
\end{proof}


\rem{agree}
The end of our proof of~\reft{surj} is a variation of
the following argument we learned from P.~Etingof. We define
$\BC[\vv^{\pm1}]_{\on{Loc}}$ inverting $(1-\vv^m),\ m\in \BZ$.
We consider a $\BC[\vv^{\pm1}]_{\on{Loc}}$-algebra $\fA$ of finite
difference operators with generators
$\{\sw_i^{\pm1},D_i^{\pm1}\}_{i=1}^n$ and defining relations
$D_i\sw_j=\vv^{2\delta_{ij}}\sw_jD_i,\ [D_i,D_j]=[\sw_i,\sw_j]=0$.
Then the algebra
of $\fS_n$-invariants $\fA^{\fS_n}$ is generated by its subalgebras
$\BC[\vv^{\pm1}]_{\on{Loc}}[D_1^{\pm1},\ldots,D_n^{\pm1}]^{\fS_n}$
and
$\BC[\vv^{\pm1}]_{\on{Loc}}[\sw_1^{\pm1},\ldots,\sw_n^{\pm1}]^{\fS_n}$.

Indeed, let $\fB$ be the $\BC[\vv^{\pm1}]_{\on{Loc}}$-algebra
generated by $\sw^{\pm1},D^{\pm1}$ subject to $D\sw=\vv^2\sw D$.
Then $\fA=\fB^{\otimes n}$ (tensor product over $\BC[\vv^{\pm1}]_{\on{Loc}}$),
and $\fA^{\fS_n}=\on{Sym}^n\fB$ (symmetric power over $\BC[\vv^{\pm1}]_{\on{Loc}}$).
Now $\on{Sym}^n\fB$ is spanned by the
elements $\{b^{\otimes n}\}_{b\in\fB}$, and hence $\on{Sym}^n\fB$ is
generated by the elements $\{b_{(1)}+\ldots+b_{(n)}\}_{b\in\fB}$, where
  $b_{(r)}=1\otimes\cdots\otimes1\otimes b \otimes1\otimes\cdots\otimes1\ (b$ at the $r$-th entry).
Indeed, it suffices to verify the generation claim for an algebra
$\BC[\vv^{\pm1}]_{\on{Loc}}[b]$ where it is nothing but the
fundamental theorem on symmetric functions.

We conclude that $\on{Sym}^n\fB$ is generated by the elements
$\{p_{m,k}=\sum_{r=1}^n\sw_r^mD_r^k\}_{m,k\in\BZ}$. However,
$p_{m,k}=(\vv^{2mk}-1)^{-1}[\sum_{r=1}^nD_r^k,\sum_{s=1}^n\sw_s^m]$
for $m\ne0\ne k$.
\erem


\rem{quantum result}
Motivated by~\cite[Remark~3.5]{bfn} we call
$\CO_{\varpi_n}\in K^{GL(n,\CO)\rtimes\BC^\times}(\Gr_{GL(n)})$
{\em the quantum resultant}. In fact, it is a quantization of the boundary
equation for the trigonometric zastava $\oZ^n_{SL(2)}$ which is
nothing but the resultant of two polynomials. Note that, up to
multiplication by an element of $\BC[\vv^{\pm1}]$, the quantum
resultant is uniquely characterized by the property
\eq{qcom}
  \CO_{\varpi_n}\ol{\Phi}{}^0_{-2n}(A^\pm_{\pm r})=
  \vv^{\pm (2r-n)}\ol{\Phi}{}^0_{-2n}(A^\pm_{\pm r})\CO_{\varpi_n},\
  \CO_{\varpi_n}\ol{\Phi}{}^0_{-2n}(f_p)=\vv^{2p}\ol{\Phi}{}^0_{-2n}(f_p)\CO_{\varpi_n}.
\end{equation}
\erem


\rem{perverse}
Here is a geometric explanation of the equality
\eq{calL}
  \CO(-k-1)_{\varpi_1}*\CO_{\varpi_k}-\vv^{2(k+1)}\CO_{\varpi_k}*\CO(-k-1)_{\varpi_1}=
  (-1)^k(1-\vv^2)\vv^{-2(k+1)}\CO(-1)_{\varpi_{k+1}},
\end{equation}
established as an induction step during our proof of~\refp{11.3}.
We have the convolution morphisms
  $$\Gr^{\varpi_1}\widetilde{\times}\Gr^{\varpi_k}\stackrel{m}{\longrightarrow}
    \ol\Gr{}^{\varpi_1+\varpi_k}
    \stackrel{m'}{\longleftarrow}\Gr^{\varpi_k}\widetilde{\times}\Gr^{\varpi_1},$$
and
$\ol\Gr{}^{\varpi_1+\varpi_k}=\Gr^{\varpi_1+\varpi_k}\sqcup\Gr^{\varpi_{k+1}}$.
Let us consider the transversal slice
  $\ol\CW{}^{\varpi_1+\varpi_k}_{\varpi_{k+1}}\subset\ol\Gr{}^{\varpi_1+\varpi_k}$
through the point $\varpi_{k+1}=(1,\ldots,1,0,\ldots,0)\ (k+1$ 1's).
It suffices to check that
\begin{multline*}
  m_*\left({\CO(-k-1)_{\varpi_1}\widetilde\boxtimes
  \CO_{\varpi_k}}|_{{m^{-1}\ol\CW{}^{\varpi_1+\varpi_k}_{\varpi_{k+1}}}}\right)-\vv^{2(k+1)}m'_*\left({\CO_{\varpi_k}
  \widetilde\boxtimes
  \CO(-k-1)_{\varpi_1}}|_{{m'{}^{-1}\ol\CW{}^{\varpi_1+\varpi_k}_{\varpi_{k+1}}}}\right)=\\
  (-1)^k(1-\vv^2)\vv^{-2(k+1)}\sw_1^{-1}\cdots\sw_{k+1}^{-1},
\end{multline*}
where we view $\vv^{-2(k+1)}\sw_1^{-1}\cdots\sw_{k+1}^{-1}$ as a
character of $T\times\BC^\times\ (T\subset GL(n)$ is the diagonal Cartan
torus). According to~\cite[Corollary~3.4]{mov},
$\ol\CW{}^{\varpi_1+\varpi_k}_{\varpi_{k+1}}$ is naturally
isomorphic to the slice
$\ol\CW{}^\theta_0\subset\Gr_{GL(k+1)\times(\BC^\times)^{n-k-1}}$ where
$\theta=(1,0,\ldots,0,-1)$ is the highest coroot of $GL(k+1)$.
Moreover, the preimages of
$\ol\CW{}^{\varpi_1+\varpi_k}_{\varpi_{k+1}}$ in the two convolution
diagrams are isomorphic to the cotangent bundles $T^*\BP^k$ and
$T^*\check\BP^k$, respectively. We will keep the following notation
for the convolution morphisms restricted to the slice:
  $$T^*\BP^k\stackrel{m}{\longrightarrow}\ol\CW{}^\theta_0
    \stackrel{m'}{\longleftarrow}T^*\check\BP^k.$$
Note also that $\ol\CW{}^\theta_0$ is isomorphic to the minimal
nilpotent orbit closure
$\ol\BO_{\on{min}}\subset{\mathfrak{sl}}_{k+1}$. Finally,
  ${\CO(-k-1)_{\varpi_1}\widetilde\boxtimes
   \CO_{\varpi_k}}|_{{m^{-1}\ol\CW{}^{\varpi_1+\varpi_k}_{\varpi_{k+1}}}}$
and
  ${\CO_{\varpi_k}\widetilde\boxtimes
   \CO(-k-1)_{\varpi_1}}|_{{m'{}^{-1}\ol\CW{}^{\varpi_1+\varpi_k}_{\varpi_{k+1}}}}$
are isomorphic to the pull-backs of $\CO_{\BP^k}(-k-1)$ and
$\CO_{\check\BP^k}(-k-1)$ respectively, but with {\em nontrivial}
$\BC^\times$-equivariant structures.

Let us explain our choice of the line bundles. According
to~\cite[Proposition~8.2]{bfm}, the convolutions in question are
$GL(k+1)\times\BC^\times$-equivariant perverse coherent sheaves on
$\ol\BO_{\on{min}}\subset{\mathfrak{sl}}_{k+1}$. Since
$\dim H^k(T^*\BP^k,\CO_{T^*\BP^k}(-k-1))=1$, while
$H^k(T^*\BP^k,\CO_{T^*\BP^k}(k+1))=0$, we have an exact sequence of
perverse coherent sheaves\footnote{We are grateful to
R.~Bezrukavnikov for his explanations about perverse coherent
sheaves.} on $\ol\BO_{\on{min}}\subset{\mathfrak{sl}}_{k+1}$:
  $$0\to j_{!*}\CO_{\BO_{\on{min}}}(-k-1)[k]\to m_*\CO_{T^*\BP^k}(-k-1)[k]\to\delta_0\to0,$$
where $j\colon \BO_{\on{min}}\hookrightarrow\ol\BO_{\on{min}}$ is
the open embedding, and $\delta_0$ is an irreducible skyscraper
sheaf at $0\in\ol\BO_{\on{min}}$ with certain $\BC^\times$-equivariant
structure. The same exact sequence holds for
$m'_*\CO_{T^*\check\BP^k}(-k-1)[k]$, but the quotient $\delta_0$ has
a {\em different} $\BC^\times$-equivariant structure.
\erem


\prop{quarter}\label{quarter}
The restriction of $\wt\Phi{}^0_{-2n}$ to $\fY^0_{-2n,-}$ is injective.
\eprop

\prf
Consider an ordering
  $A^-_0\prec A^-_{-1}\prec\ldots\prec A^-_{-n+1}\prec C^-_{0}\prec\ldots\prec C^-_{-n+1}$.
We set $(A^-_0)^{-k}:=((-\vv^2)^{-n}A^-_{-n})^k$ for $k>0$.
For $\vec{r}=(r_1,\ldots,r_{2n})\in \BZ\times \BN^{2n-1}$, we define the ordered monomial
  $m_{\vec{r}}:=(A^-_0)^{r_1}(A^-_{-1})^{r_2}\cdots (A^-_{-n+1})^{r_n}(C^-_0)^{r_{n+1}}\cdots (C^-_{-n+1})^{r_{2n}}$.

\begin{Lem}\label{PBW span}
The ordered monomials $\{m_{\vec{r}}\}$ span $\fY^0_{-2n,-}$.
\end{Lem}

\begin{proof}
According to relations~(\ref{ABCD1}, \ref{ABCD5}), we have
$[A^-_t,A^-_s]=[C^-_t,C^-_s]=0$ for $s,t\leq 0$. Due to
Remark~\ref{first indication of conjecture}, we also have $C^-_s=0$
for $s\leq -n$. It remains to prove that all $A^-_t$ can be taken to the
left of all $C^-_s$. This is implied by the fact that $C^-_sA^-_t$ can be written as
a linear combination of normally ordered monomials $A^-_{t'}C^-_{s'}$.
The latter claim follows from relation~(\ref{ABCD7}) by induction in $\min\{-t,-s\}$.
The lemma is proved.
\end{proof}

The following result will be proved in Section~\ref{Proof of PBW basis}:

\lem{PBW basis}
(a) The ordered monomials $\{m_{\vec{r}}\}$ form a $K_{\BC^\times}(\on{pt})$-basis of $\fY^0_{-2n,-}$.

\noindent
(b) $\{\ol\Phi{}^0_{-2n}(m_{\vec{r}})\}$ form a $K_{\BC^\times}(\on{pt})$-basis of $\ol\Phi{}^0_{-2n}(\fY^0_{-2n,-})$.
\elem

The proposition is proved.
\epr


\subsection{Positive Grassmannian}
\label{Proof of PBW basis}
$\ $

Recall the positive part of the affine Grassmannian
$\Gr_{GL(n)}^+\subset\Gr_{GL(n)}$~\cite[3(ii)]{bfn} parametrizing
the {\em sublattices} in the standard one.
Recall also that
 $K^{GL(n,\CO)\rtimes\BC^\times}_{\on{loc}}(\Gr^{\varpi_1})=
  K^{GL(n,\CO)\rtimes\BC^\times}_{\on{loc}}(\BP^{n-1})$
is generated over $K_{GL(n)}(\on{pt})$ by the classes of
$\CO(a)_{\varpi_1}, -n+1\leq a\leq 0$.
The proof of~\reft{surj}
shows that
  $\ol{\Phi}{}^0_{-2n,\on{loc}}\colon \fU_{0,-2n,\on{loc}}^\ad
   \twoheadrightarrow K^{\wt{GL}(n,\CO)\rtimes\wt\BC^\times}_{\on{loc}}(\Gr_{GL(n)})$
restricts to a surjective homomorphism
  $\ol{\Phi}{}^0_{-2n,\on{loc}}\colon\
   \fY_{-2n,-,\on{loc}}^0\twoheadrightarrow
   K^{\wt{GL}(n,\CO)\rtimes\wt\BC^\times}_{\on{loc}}(\Gr_{GL(n)}^+)$.

\prop{isom}\label{isom}
  $\ol{\Phi}{}^0_{-2n,\on{loc}}\colon \fY_{-2n,-,\on{loc}}^0\iso
   K^{\wt{GL}(n,\CO)\rtimes\wt\BC^\times}_{\on{loc}}(\Gr_{GL(n)}^+)$.
\eprop

\prf
We have to check that
  $\ol{\Phi}{}^0_{-2n,\on{loc}}\colon \fY_{-2n,-,\on{loc}}^0\to
   K^{\wt{GL}(n,\CO)\rtimes\wt\BC^\times}_{\on{loc}}(\Gr_{GL(n)}^+)$
is injective. To this end, note that $\Gr_{GL(n)}^+$ is a union of
connected components numbered by nonnegative integers:
$\Gr_{GL(n)}^+=\bigsqcup_{r\in\BN}\Gr_{GL(n)}^{+,r},$
where $\Gr_{GL(n)}^{+,r}$ parametrizes the sublattices of codimension $r$ in
the standard one. The direct sum decomposition
  $K^{\wt{GL}(n,\CO)\rtimes\wt\BC^\times}_{\on{loc}}(\Gr_{GL(n)}^+)=
   \bigoplus_{r\in\BN}K^{\wt{GL}(n,\CO)\rtimes\wt\BC^\times}_{\on{loc}}(\Gr_{GL(n)}^{+,r})$
is a grading of the convolution algebra. For any connected component,
$K^{\wt{GL}(n,\CO)\rtimes\wt\BC^\times}_{\on{loc}}(\Gr_{GL(n)}^{+,r})$ is a free
${K_{\wt{GL}(n,\CO)\rtimes\wt\BC^\times}(\on{pt})}_{\on{loc}}$-module of rank $d_r$, where
$d_r$ is the number of $T$-fixed points in $\Gr_{GL(n)}^{+,r}$, that is the number
of weights of the irreducible $GL(n)$-module with the highest weight
$(r,0,\ldots,0)$, isomorphic to $\on{Sym}^r(\BC^n)$. Note that all the weights
of $\on{Sym}^r(\BC^n)$ have multiplicity one; in other words,
$d_r=\dim\on{Sym}^r(\BC^n)$.

According to~Lemma~\ref{PBW span}, we can introduce a grading
$\fY_{-2n,-,\on{loc}}^0=\bigoplus_{r\in\BN}\fY_{-2n,-,\on{loc}}^{0,r}$: a monomial
$m_{\vec{r}}$ has degree $r$ if $r_{n+1}+\ldots+r_{2n}=r$. It is immediate from
the relations between $A^-_\bullet,C^-_\bullet$-generators that this grading is
well-defined. Also, it is clear that
  $\ol{\Phi}{}^0_{-2n,\on{loc}}(\fY_{-2n,-,\on{loc}}^{0,r})\subset
   K^{\wt{GL}(n,\CO)\rtimes\wt\BC^\times}_{\on{loc}}(\Gr_{GL(n)}^{+,r})$.
Meanwhile, we know from~\reft{surj} that
  $\ol{\Phi}{}^0_{-2n,\on{loc}}(\fY_{-2n,-,\on{loc}}^{0,r})=
   K^{\wt{GL}(n,\CO)\rtimes\wt\BC^\times}_{\on{loc}}(\Gr_{GL(n)}^{+,r})$.
On the other hand, we know from~Lemma~\ref{PBW span} that $\fY_{-2n,-,\on{loc}}^{0,r}$
as a left ${K_{\wt{GL}(n,\CO)\rtimes\wt\BC^\times}(\on{pt})}_{\on{loc}}$-module
has no more than $d'_r$ generators, where $d'_r$ is the number of compositions of $r$
into $n$ (ordered) summands. Since $d_r=d'_r$, we conclude that
  $\ol{\Phi}{}^0_{-2n,\on{loc}}\colon \fY_{-2n,-,\on{loc}}^{0,r}\to
   K^{\wt{GL}(n,\CO)\rtimes\wt\BC^\times}_{\on{loc}}(\Gr_{GL(n)}^{+,r})$
must be an isomorphism, and $\fY_{-2n,-,\on{loc}}^{0,r}$ is a free left
${K_{\wt{GL}(n,\CO)\rtimes\wt\BC^\times}(\on{pt})}_{\on{loc}}$-module of rank $d_r=d'_r$.
This completes the proof of~\refp{isom},~\refl{PBW basis} (and~\refp{quarter}).
\epr

\rem{identity}
One can check that the natural morphism
$K^{\wt{GL}(n,\CO)\rtimes\wt\BC^\times}(\Gr_{GL(n)}^+)[\CO_{\varpi_n}^{-1}]\to
K^{\wt{GL}(n,\CO)\rtimes\wt\BC^\times}(\Gr_{GL(n)})$ is an isomorphism.
Now it follows from the proof of~\refp{isom} and~\reft{surj} that in order
to check~Conjecture~\ref{description of truncated},~\refco{Phi iso} and~\refco{Phi iso Yang}
in our case: $\on{Ker}(\wt\Phi{}^0_{-2n,\loc})=\CI^0_{-2n,\loc}$, it suffices to
check the following equality in $\CU^\ad_{0,-2n}/\CI^0_{-2n}$:
\begin{equation*}
  -\vv^{n^2-2}(1-\vv^2)^2\cdot
  (\on{ad}_{f_{1-n}}^{\vv^{2n}}\on{ad}_{f_{3-n}}^{\vv^{2(n-1)}}\cdots\on{ad}_{f_{n-3}}^{\vv^4}f_{n-1})
  (\on{ad}_{e_{1-2n}}^{\vv^{-2n}}
  \on{ad}_{e_{3-2n}}^{\vv^{-2(n-1)}}\cdots\on{ad}_{e_{-3}}^{\vv^{-4}}e_{-1})=1.
\end{equation*}
\erem

\rem{kedem}
Consider a subalgebra $\fU^<_{0,-2n}\subset\fU^{\on{ad}}_{0,-2n}$ generated by
$\{(\vv-\vv^{-1})f_s\}_{s\in\BZ}$. Note that it is independent of $n$,
cf.~Proposition~\ref{Triangular decomposition}. The image
$\ol{\Phi}{}^0_{-2n}(\fU^<_{0,-2n})$ in
$K^{\wt{GL}(n,\CO)\rtimes\wt\BC^\times}(\Gr_{GL(n)})$ is isomorphic to the {\em $M$-system
  algebra} $\CU'_{n-1}$ of~\cite{dfk}. In particular, the generators
$M_{m,s}\in\CU'_{n-1}$ of~\cite[2.1]{dfk} correspond to scalar multiples of the classes
$\CO(-s)_{\varpi_m}\in K^{\wt{GL}(n,\CO)\rtimes\wt\BC^\times}(\Gr_{GL(n)})$,
cf.~\refe{8.2} and~\cite[(2.23)]{dfk}.
\erem

\section{Coproducts on shifted quantum affine algebras}
\label{Copro}

Throughout this section, we work mainly with simply-connected shifted
quantum affine algebras. However, all the results can be obviously
generalized to the adjoint versions.


\subsection{Drinfeld formal coproduct}
$\ $

The standard quantum loop algebra $U_\vv(L\fg)$ admits the Drinfeld
formal coproduct
\begin{equation*}
  \wt{\Delta}\colon U_\vv(L\fg)\longrightarrow U_\vv(L\fg)\widehat{\otimes} U_\vv(L\fg),
\end{equation*}
defined in the new Drinfeld realization of $U_\vv(L\fg)$ via
\begin{equation}\label{Drinfeld coproduct}
\begin{split}
  & \wt{\Delta}(e_i(z)):=e_i(z)\otimes 1 + \psi^-_i(z)\otimes e_i(z),\\
  & \wt{\Delta}(f_i(z)):=f_i(z)\otimes \psi^+_i(z) + 1\otimes f_i(z),\\
  & \wt{\Delta}(\psi^\pm_i(z)):=\psi^\pm_i(z)\otimes \psi^\pm_i(z).
\end{split}
\end{equation}

\begin{Rem}
Composing $\wt{\Delta}$ with the $\BC^\times$-action on the first factor,
D.~Hernandez obtained a deformed coproduct
$\Delta_\zeta\colon  U_\vv(L\fg)\to U_\vv(L\fg)\otimes U_\vv(L\fg)((\zeta))$,
where $\zeta$ is a formal variable, see~\cite[Section 6]{he}.
\end{Rem}

This can be obviously generalized to the shifted setting.

\begin{Lem}\label{new Drinfeld coproduct}
For any coweights $\mu^\pm_1, \mu^\pm_2\in \Lambda$, there is a $\BC(\vv)$-algebra homomorphism
\begin{equation*}
    \wt{\Delta}\colon \CU^\ssc_{\mu^+_1+\mu^+_2, \mu^-_1+\mu^-_2}\longrightarrow
    \CU^\ssc_{\mu^+_1,\mu^-_1}\widehat{\otimes}\CU^\ssc_{\mu^+_2,\mu^-_2},
\end{equation*}
defined via~(\ref{Drinfeld coproduct}).
\end{Lem}

We call this homomorphism a formal coproduct for shifted quantum
affine algebras. Given two representations $V_1, V_2$ of
$\CU^\ssc_{\mu^+_1,\mu^-_1}, \CU^\ssc_{\mu^+_2,\mu^-_2}$,
respectively, we will use $V_1\wt{\otimes} V_2$ to denote the
representation of $\CU^\ssc_{\mu^+_1+\mu^+_2, \mu^-_1+\mu^-_2}$ on
the vector space $V_1\otimes V_2$ induced by $\wt{\Delta}$, whenever
the action of the infinite sums representing
$\wt{\Delta}(e_{i,r}), \wt{\Delta}(f_{i,r})$
are well-defined. We will discuss a particular example of this
construction in Section~\ref{Factorizing Laumons}.


\subsection{Drinfeld-Jimbo coproduct}
\label{Drinfeld-Jimbo coproduct}
$\ $

\medskip
The standard quantum loop algebra $U_\vv(L\fg)$ also admits the Drinfeld-Jimbo coproduct
\begin{equation*}
  \Delta\colon U_\vv(L\fg)\longrightarrow U_\vv(L\fg)\otimes U_\vv(L\fg),
\end{equation*}
defined in the Drinfeld-Jimbo realization of $U_\vv(L\fg)$ via
\begin{equation*}
    \Delta\colon
    E_i\mapsto E_i\otimes K_i+1\otimes E_i,\
    F_i\mapsto F_i\otimes 1 + K_i^{-1}\otimes F_i,\
    K^{\pm 1}_i\mapsto K^{\pm 1}_i\otimes K^{\pm 1}_i,\ i\in \wt{I}.
\end{equation*}
Recall that $\wt{I}=I\cup\{i_0\}$ is the vertex set of the extended Dynkin
diagram and $\{E_i, F_i, K^{\pm 1}_i\}_{i\in \wt{I}}$ are the
standard Drinfeld-Jimbo generators of $U^\ddj_\vv(L\fg)\simeq U_\vv(L\fg)$.

We also denote the Drinfeld-Jimbo coproduct on $U^\ad_\vv(L\fg)$ by $\Delta^\ad$:
the natural inclusion $U_\vv(L\fg)\hookrightarrow U^\ad_\vv(L\fg)$ intertwines
$\Delta$ and $\Delta^\ad$, while $\Delta^\ad(\phi^\pm_i)=\phi^\pm_i\otimes \phi^\pm_i$.

\medskip
The goal of this section is to generalize these coproducts to the shifted
setting. In other words, given $\fg$ and coweights $\mu_1,\mu_2\in \Lambda$,
we would like to construct homomorphisms
\begin{equation*}
    \Delta_{\mu_1,\mu_2}\colon  \CU^\ssc_{0,\mu_1+\mu_2}\longrightarrow
    \CU^\ssc_{0,\mu_1}\otimes \CU^\ssc_{0,\mu_2},
\end{equation*}
which coincide with $\Delta$ in the particular case $\mu_1=\mu_2=0$.
We provide such a construction for the simplest case $\fg=\ssl_2$ in
Section~\ref{coproduct for sl_2 antidominant} ($\mu_1,\mu_2\in \Lambda^-$)
and Section~\ref{Coproduct generalization 1} (general $\mu_1,\mu_2$).
Using the RTT presentation of $U_\vv(L\ssl_n)$, we generalize this to obtain
$\Delta_{\mu_1,\mu_2}$ for $\fg=\ssl_n$ in
Section~\ref{Coproduct on sl_n for antidominant} ($\mu_1,\mu_2\in \Lambda^-$)
and Section~\ref{General shifted coproduct sl_n} (general $\mu_1,\mu_2$).

\begin{Rem}\label{comment}
(a) This result is nontrivial due to an absence of the Drinfeld-Jimbo type
presentation of shifted quantum affine algebras.

\noindent
(b) A similar coproduct for the shifted Yangians has been constructed
in~\cite{fkp} for arbitrary simply-laced $\fg$.

\noindent
(c) Once $\Delta_{\mu_1,\mu_2}$ is constructed, one should
be able to immediately extend it to the homomorphism
  $\Delta^\ad_{\mu_1,\mu_2}\colon  \CU^\ad_{0,\mu_1+\mu_2}\to \CU^\ad_{0,\mu_1}\otimes \CU^\ad_{0,\mu_2}$
by setting $\Delta^\ad_{\mu_1,\mu_2}(\phi^\pm_i)=\phi^\pm_i \otimes \phi^\pm_i$.
\end{Rem}


\subsection{Homomorphisms $\Delta_{b_1,b_2}$ for $b_1,b_2\in \BZ_{\leq 0},\ \fg=\ssl_2$}
\label{coproduct for sl_2 antidominant}
$\ $

We start this subsection by explicitly computing the
Drinfeld-Jimbo coproduct of the Drinfeld generators
$e_0, e_{-1}, f_0, f_1, \psi^\pm_0$ of $U_\vv(L\ssl_2)$ and
$h_{\pm 1}=\pm \frac{\psi^\mp_0\psi^\pm_{\pm 1}}{\vv-\vv^{-1}}$, which
generate the quantum loop algebra $U_\vv(L\ssl_2)$.

\begin{Lem}\label{Explicit sl2}
We have
  $$\Delta(e_0)=e_0\otimes \psi^+_0+1\otimes e_0,\
    \Delta(e_{-1})=e_{-1}\otimes \psi^-_0+1\otimes e_{-1},$$
  $$\Delta(f_0)=f_0\otimes 1+\psi^-_0\otimes f_0,\
    \Delta(f_1)=f_1\otimes 1+\psi^+_0\otimes f_1,\
    \Delta(\psi^\pm_0)=\psi^\pm_0\otimes\psi^\pm_0,$$
  $$\Delta(h_1)=h_1\otimes 1+1\otimes h_1-(\vv^2-\vv^{-2})e_0\otimes f_1,\
    \Delta(h_{-1})=h_{-1}\otimes 1 +1\otimes h_{-1}+(\vv^2-\vv^{-2})e_{-1}\otimes f_0.$$
\end{Lem}

\begin{proof}
This is a straightforward computation based on the explicit identification
between the Drinfeld-Jimbo and the new Drinfeld realizations
of the quantum loop algebra $U_\vv(L\ssl_2)$
of~Theorem~\ref{identification of DJ and Dr}:
  $e_0=E_{i_1}, f_0=F_{i_1}, \psi^\pm_0=K^{\pm 1}_{i_1}, e_{-1}=K^{-1}_{i_1}F_{i_0}, f_1=E_{i_0}K_{i_1}.$
\end{proof}

The key result of this subsection provides analogues of $\Delta$
for antidominantly shifted quantum affine algebras of $\ssl_2$.
For $\mu_1,\mu_2\in \Lambda^-$, we construct homomorphisms
  $\Delta_{b_1,b_2}\colon \CU^\ssc_{0,b_1+b_2}\to \CU^\ssc_{0,b_1}\otimes \CU^\ssc_{0,b_2},$
where $b_1:=\alphavee(\mu_1), b_2:=\alphavee(\mu_2)$ (so that $b_1,b_2\in \BZ_{\leq 0}$).

\begin{Thm}\label{Coproduct for sl2}
For any $b_1,b_2\in \BZ_{\leq 0}$, there is a unique $\BC(\vv)$-algebra homomorphism
\begin{equation*}
    \Delta_{b_1,b_2}\colon
    \CU^\ssc_{0,b_1+b_2}\longrightarrow \CU^\ssc_{0,b_1}\otimes \CU^\ssc_{0,b_2}
\end{equation*}
(we will denote $\Delta=\Delta_{b_1,b_2}$ when the algebras involved are clear),
such that
  $$\Delta(e_r)=1\otimes e_r,\ \Delta(f_s)=f_s\otimes 1
    \ \mathrm{for}\ b_2\leq r<0, b_1<s\leq 0,$$
  $$\Delta(e_0)=e_0\otimes \psi^+_0+1\otimes e_0,\
    \Delta(e_{b_2-1})=e_{-1}\otimes \psi^-_{b_2}+1\otimes e_{b_2-1},$$
  $$\Delta(f_1)=f_1\otimes 1+\psi^+_0\otimes f_1,\
    \Delta(f_{b_1})=f_{b_1}\otimes 1+ \psi^-_{b_1}\otimes f_0,$$
  $$\Delta((\psi^+_0)^{\pm 1})=(\psi^+_0)^{\pm 1}\otimes (\psi^+_0)^{\pm 1},\
    \Delta((\psi^-_{b_1+b_2})^{\pm 1})=(\psi^-_{b_1})^{\pm 1}\otimes (\psi^-_{b_2})^{\pm 1},$$
  $$\Delta(h_1)=h_1\otimes 1+1\otimes h_1-(\vv^2-\vv^{-2})e_0\otimes f_1,\
    \Delta(h_{-1})=h_{-1}\otimes 1 +1\otimes h_{-1}+(\vv^2-\vv^{-2})e_{-1}\otimes f_0.$$
\end{Thm}

These homomorphisms generalize the Drinfeld-Jimbo coproduct, since
we recover the formulas of Lemma~\ref{Explicit sl2} for $b_1=b_2=0$.
The proof of~Theorem~\ref{Coproduct for sl2} is presented
in~Appendix~\ref{Proof Coproduct}
and is crucially based on~Theorem~\ref{Levendorskii-type} which
provides a presentation of the shifted quantum affine algebras via a finite
number of generators and relations.

\begin{Rem}\label{Maps jmath}
The similarity between the formulas for $\Delta_{b_1,b_2}$ of
Theorem~\ref{Coproduct for sl2} and Drinfeld-Jimbo coproduct $\Delta$ of
Lemma~\ref{Explicit sl2} can be explained as follows.
Let $U^-_\vv$ (resp.\ $\CU^{\ssc,-}_{0,b_1,b_2}$) be the
subalgebra of $U_\vv(L\ssl_2)$ (resp.\ $\CU^\ssc_{0,b_1+b_2}$)
generated by $\{e_{-1},f_0,(\psi^-_0)^{\pm 1}\}$, or equivalently,
by $\{e_{-r-1},f_{-r},(\psi^-_0)^{\pm 1},\psi^-_{-r-1}\}_{r\in \BN}$
(resp.\ by $\{e_{b_2-1},f_{b_1},(\psi^-_{b_1+b_2})^{\pm 1}\}$,
or equivalently, by
  $\{e_{b_2-r-1},f_{b_1-r},(\psi^-_{b_1+b_2})^{\pm1},
   \psi^-_{b_1+b_2-r-1}\}_{r\in \BN}$).
Analogously, let $U^+_\vv$ (resp.\ $\CU^{\ssc,+}_{0,b_1,b_2}$) be
the subalgebra of $U_\vv(L\ssl_2)$ (resp.\ $\CU^\ssc_{0,b_1+b_2}$)
generated by $\{e_0,f_1,(\psi^+_0)^{\pm 1}\}$ in both cases, or equivalently, by
$\{e_{r},f_{r+1},(\psi^+_0)^{\pm 1},\psi^+_{r+1}\}_{r\in \BN}$.
Then, there are unique $\BC(\vv)$-algebra homomorphisms
$\jmath^\pm_{b_1,b_2}\colon  U^\pm_\vv\to \CU^{\ssc,\pm}_{0,b_1,b_2}$,
such that
\begin{equation*}
  \jmath^+_{b_1,b_2}\colon  e_0\mapsto e_0,\ f_1\mapsto f_1,\ (\psi^+_0)^{\pm 1}\mapsto (\psi^+_0)^{\pm 1},
\end{equation*}
\begin{equation*}
  \jmath^-_{b_1,b_2}\colon  e_{-1}\mapsto e_{b_2-1},\ f_0\mapsto f_{b_1},\ (\psi^-_0)^{\pm 1}\mapsto (\psi^-_{b_1+b_2})^{\pm 1}.
\end{equation*}
Moreover, the following diagram is commutative:

\setlength{\unitlength}{1cm}
\begin{picture}(4,3.2)
 \put(3.3,2.5){$U^\pm_\vv$}
 \put(2.9,0.5){$\CU^{\ssc,\pm}_{0,b_1,b_2}$}
 \put(9.3,2.5){$U^\pm_\vv\otimes U^\pm_\vv$}
 \put(8.8,0.5){$\CU^{\ssc,\pm}_{0,b_1,0}\otimes \CU^{\ssc,\pm}_{0,0,b_2}$}

 \put(4.0,2.6){\vector (1,0){5.2}}
 \put(4.2,0.6){\vector (1,0){4.3}}
 \put(3.5,2.3){\vector (0,-1){1.3}}
 \put(10.1,2.3){\vector (0,-1){1.4}}

 \put(6.4,2.75){$\Delta$}
 \put(6.1,0.25){$\Delta_{b_1,b_2}$}
 \put(2.6,1.5){$\jmath^\pm_{b_1,b_2}$}
 \put(10.2,1.5){$\jmath^\pm_{b_1,0}\otimes \jmath^\pm_{0,b_2}$}
\end{picture}
\end{Rem}

\begin{Rem}\label{ad-coproduct on sl2}
The aforementioned homomorphism $\Delta_{b_1,b_2}$ can be naturally extended to the homomorphism
$\Delta^\ad_{b_1,b_2}\colon  \CU^\ad_{0,b_1+b_2}\to \CU^\ad_{0,b_1}\otimes \CU^\ad_{0,b_2}$
by setting $\Delta^\ad_{b_1,b_2}(\phi^\pm)=\phi^\pm\otimes \phi^\pm$.
\end{Rem}


\subsection{Homomorphisms $\Delta_{b_1,b_2}$ for arbitrary $b_1,b_2\in \BZ,\ \fg=\ssl_2$}
\label{Coproduct generalization 1}
$\ $

In this subsection, we generalize the construction of
$\Delta_{b_1,b_2}$ of Theorem~\ref{Coproduct for sl2}
($b_1,b_2\in \BZ_{\leq 0}$) to the general case $b_1,b_2\in \BZ$.
We follow the corresponding construction for the shifted
Yangians of~\cite[Theorem 4.12]{fkp}.

The key ingredient of our approach are the \emph{shift homomorphisms}
$\iota_{n,m_1,m_2}$ (the trigonometric analogues of the shift
homomorphisms of~\cite{fkp}).

\begin{Prop}\label{Maps iota}
For any $n\in \BZ$ and $m_1,m_2\in \BZ_{\leq 0}$, there is a unique
$\BC(\vv)$-algebra homomorphism
$\iota_{n,m_1,m_2}\colon  \CU^\ssc_{0,n}\to \CU^\ssc_{0,n+m_1+m_2}$,
which maps the currents as follows
\begin{equation*}
    e(z)\mapsto (1-z^{-1})^{-m_1}e(z),\
    f(z)\mapsto (1-z^{-1})^{-m_2}f(z),\
    \psi^\pm(z)\mapsto (1-z^{-1})^{-m_1-m_2}\psi^\pm(z).
\end{equation*}
\end{Prop}

\begin{proof}
The above assignment is obviously compatible with defining
relations~(\ref{U1}--\ref{U8}). Moreover, we have
  $\iota_{n,m_1,m_2}\colon  \psi^+_0\mapsto \psi^+_0,
   \psi^-_{n}\mapsto (-1)^{m_1+m_2}\psi^-_{n+m_1+m_2}$.
\end{proof}

These homomorphisms satisfy two important properties:

\begin{Lem}\label{Properties of iota}
(a) We have
  $\iota_{n+m_1+m_2,m'_1,m'_2}\circ \iota_{n,m_1,m_2}=\iota_{n,m_1+m'_1,m_2+m'_2}$
for any $n\in \BZ$ and $m_1,m_2,m'_1,m'_2\in \BZ_{\leq 0}$.

\noindent
(b) The homomorphism $\iota_{n,m_1,m_2}$ is injective
for any $n\in \BZ$ and $m_1,m_2\in \BZ_{\leq 0}$.
\end{Lem}

Part (a) is obvious, while part (b) is proved in
Appendix~\ref{Proof of PBW} and follows from the PBW property for
$\CU^\ssc_{0,n}$ (cf.~Theorem~\ref{1026}). The following is the key result of this subsection.

\begin{Thm}\label{Coproduct for sl2 in general}
For any $b_1,b_2\in \BZ$ and $b:=b_1+b_2$, there is a unique $\BC(\vv)$-algebra
homomorphism
\begin{equation*}
  \Delta_{b_1,b_2}\colon  \CU^\ssc_{0,b}\longrightarrow \CU^\ssc_{0,b_1}\otimes \CU^\ssc_{0,b_2},
\end{equation*}
such that for any $m_1,m_2\in \BZ_{\leq 0}$ the following diagram is commutative:

\setlength{\unitlength}{1cm}
\begin{picture}(4,3.2)
 \put(3.2,2.5){$\CU^\ssc_{0,b}$}
 \put(2.6,0.5){$\CU^\ssc_{0,b+m_1+m_2}$}
 \put(9.0,2.5){$\CU^\ssc_{0,b_1}\otimes \CU^\ssc_{0,b_2}$}
 \put(8.4,0.5){$\CU^\ssc_{0,b_1+m_1}\otimes \CU^\ssc_{0,b_2+m_2}$}

 \put(4.1,2.6){\vector (1,0){4.6}}
 \put(4.6,0.6){\vector (1,0){3.7}}
 \put(3.5,2.3){\vector (0,-1){1.4}}
 \put(10.05,2.3){\vector (0,-1){1.4}}

 \put(6.0,2.75){$\Delta_{b_1,b_2}$}
 \put(5.4,0.25){$\Delta_{b_1+m_1,b_2+m_2}$}
 \put(2.2,1.5){$\iota_{b,m_2,m_1}$}
 \put(10.1,1.5){$\iota_{b_1,0,m_1}\otimes \iota_{b_2,m_2,0}$}
\end{picture}
\end{Thm}

The proof of this theorem is presented in~Appendix~\ref{Proof Coproduct General}
and is similar to the proof of~\cite[Theorem 4.12]{fkp}.

\begin{Cor}\label{uniform coproduct h_1}
For any $b_1,b_2\in \BZ$, we have
\begin{equation*}
  \Delta_{b_1,b_2}(h_1)=h_1\otimes 1+1\otimes h_1-(\vv^2-\vv^{-2})e_0\otimes f_1,
\end{equation*}
\begin{equation*}
  \Delta_{b_1,b_2}(h_{-1})=h_{-1}\otimes 1 +1\otimes h_{-1}+(\vv^2-\vv^{-2})e_{-1}\otimes f_0.
\end{equation*}
\end{Cor}

\begin{proof}
In the antidominant case $b_1,b_2\in \BZ_{\leq 0}$, both equalities
are due to our definition of $\Delta_{b_1,b_2}$ of Theorem~\ref{Coproduct for sl2}.
For general $b_1,b_2$, choose $m_1,m_2\in \BZ_{\leq 0}$ such that
$b_1+m_1,b_2+m_2\in \BZ_{\leq 0}$. By the definition of $\iota_{b,m_2,m_1}$, we have
  $\iota_{b,m_2,m_1}(h_{\pm 1})=h_{\pm 1}\pm \frac{m_1+m_2}{\vv-\vv^{-1}}$.
Meanwhile, we also have
\begin{equation*}
   \iota_{b_1,0,m_1}\otimes \iota_{b_2,m_2,0}(h_{\pm 1}\otimes 1+1\otimes h_{\pm 1})=
   h_{\pm 1}\otimes 1 + 1\otimes h_{\pm 1}\pm \frac{m_1+m_2}{\vv-\vv^{-1}},
\end{equation*}
while $\iota_{b_1,0,m_1}(e_r)=e_r,\ \iota_{b_2,m_2,0}(f_s)=f_s$ for any $r,s\in \BZ$.
The result follows by combining the formula for $\Delta_{b_1+m_1,b_2+m_2}(h_{\pm 1})$
with the commutativity of the diagram of~Theorem~\ref{Coproduct for sl2 in general}
(we also use injectivity of the vertical arrows, due to~Lemma~\ref{Properties of iota}(b)).
\end{proof}

The following result is analogous to~\cite[Proposition 4.14]{fkp} and we
leave its proof to the interested reader.

\begin{Lem}\label{Comultiplication property}
For $b=b_1+b_2+b_3$ with $b_1,b_3\in \BZ, b_2\in \BZ_{\leq 0}$,
the following diagram is commutative:

\setlength{\unitlength}{1cm}
\begin{picture}(4,3.2)
 \put(3.2,2.5){$\CU^\ssc_{0,b}$}
 \put(1.95,0.5){$\CU^\ssc_{0,b_1+b_2}\otimes \CU^\ssc_{0,b_3}$}
 \put(9.0,2.5){$\CU^\ssc_{0,b_1}\otimes \CU^\ssc_{0,b_2+b_3}$}
 \put(8.6,0.5){$\CU^\ssc_{0,b_1}\otimes \CU^\ssc_{0,b_2}\otimes \CU^\ssc_{0,b_3}$}

 \put(4.0,2.6){\vector (1,0){4.8}}
 \put(4.6,0.6){\vector (1,0){3.8}}
 \put(3.5,2.3){\vector (0,-1){1.4}}
 \put(10.05,2.3){\vector (0,-1){1.4}}

 \put(6.0,2.75){$\Delta_{b_1,b_2+b_3}$}
 \put(5.8,0.25){$\Delta_{b_1,b_2}\otimes\on{Id}$}
 \put(2.0,1.5){$\Delta_{b_1+b_2,b_3}$}
 \put(10.2,1.5){$\on{Id}\otimes \Delta_{b_2,b_3}$}
\end{picture}
\end{Lem}

\subsection{Drinfeld-Jimbo coproduct on $U_\vv(L\ssl_n)$ via Drinfeld generators}
$\ $

According to Theorem~\ref{Levendorskii-type}, the quantum loop algebra $U_\vv(L\ssl_n)$
is generated by the elements
$\{e_{i,0}, f_{i,0}, e_{i,-1}, f_{i,1},\psi^{\pm}_{i,0}, h_{i,\pm 1}\}_{i=1}^{n-1}$.
The key result of this subsection provides explicit formulas for the action
of the Drinfeld-Jimbo coproduct $\Delta$ on these generators of $U_\vv(L\ssl_n)$.
Since $e_{i,0}=E_i, f_{i,0}=F_i, \psi^\pm_{i,0}=K_i^{\pm 1}$ (for $i\in I=\{1,2,\cdots,n-1\}$),
we obviously have
\begin{equation*}
  \Delta(e_{i,0})=1\otimes e_{i,0}+e_{i,0}\otimes \psi^+_{i,0},\
  \Delta(f_{i,0})=f_{i,0}\otimes 1 + \psi^-_{i,0}\otimes f_{i,0},\
  \Delta(\psi^\pm_{i,0})=\psi^\pm_{i,0}\otimes \psi^\pm_{i,0}.
\end{equation*}
It remains to compute the coproduct of the remaining generators above.

\begin{Thm}\label{Coproduct on quantum loop sln}
Let $\Delta$ be the Drinfeld-Jimbo coproduct on $U_\vv(L\ssl_n)$.
Then, we have
\begin{equation}\label{h_i1}
\begin{split}
  & \Delta(h_{i,1})=\\
  & h_{i,1}\otimes 1+1\otimes h_{i,1}-
  (\vv^2-\vv^{-2})E_{i,i+1}^{(0)}\otimes F^{(1)}_{i+1,i}+
  (\vv-\vv^{-1})\sum_{l>i+1}E^{(0)}_{i+1,l}\otimes F^{(1)}_{l,i+1}+\\
  & (\vv-\vv^{-1})\sum_{k<i}\vv^{k+1-i}\wt{E}^{(0)}_{ki}\otimes F^{(1)}_{ik}+
  \vv^{-2}(\vv-\vv^{-1})\sum_{l>i+1}[E^{(0)}_{i,i+1},E^{(0)}_{i+1,l}]_{\vv^3}\otimes F^{(1)}_{li}-\\
  & (\vv-\vv^{-1})\sum_{k<i} \vv^{k-i-1} [E^{(0)}_{i,i+1},\wt{E}^{(0)}_{ki}]_{\vv^3}\otimes F^{(1)}_{i+1,k}+\\
  & (\vv-\vv^{-1})^2\sum_{l>i+1}^{k<i}\vv^{k-i}\left(E^{(0)}_{il}\wt{E}^{(0)}_{ki}-E^{(0)}_{i+1,l}\wt{E}^{(0)}_{k,i+1}\right)\otimes F^{(1)}_{lk},
\end{split}
\end{equation}
\begin{equation}\label{h_i-1}
\begin{split}
  & \Delta(h_{i,-1})=\\
  & h_{i,-1}\otimes 1+1\otimes h_{i,-1}+
  (\vv^2-\vv^{-2})E_{i,i+1}^{(-1)}\otimes F^{(0)}_{i+1,i}-
  (\vv-\vv^{-1})\sum_{l>i+1}E^{(-1)}_{i+1,l}\otimes F^{(0)}_{l,i+1}-\\
  & (\vv-\vv^{-1})\sum_{k<i}\vv^{i-k-1}E^{(-1)}_{ki}\otimes \wt{F}^{(0)}_{ik}-
  \vv^{2}(\vv-\vv^{-1})\sum_{l>i+1}E^{(-1)}_{il}\otimes [F^{(0)}_{l,i+1},F^{(0)}_{i+1,i}]_{\vv^{-3}}+\\
  & (\vv-\vv^{-1})\sum_{k<i} \vv^{i+1-k} E^{(-1)}_{k,i+1}\otimes [\wt{F}^{(0)}_{ik},F^{(0)}_{i+1,i}]_{\vv^{-3}}-\\
  & (\vv-\vv^{-1})^2\sum_{l>i+1}^{k<i}\vv^{i-k}E^{(-1)}_{kl}\otimes \left(\wt{F}^{(0)}_{i+1,k}F^{(0)}_{l,i+1}-\wt{F}^{(0)}_{ik}F^{(0)}_{li}\right),
\end{split}
\end{equation}
\begin{equation}\label{e_i-1}
\begin{split}
  & \Delta(e_{i,-1})=1\otimes e_{i,-1}+e_{i,-1}\otimes \psi^-_{i,0}-
  (\vv-\vv^{-1})\sum_{l>i+1}E^{(-1)}_{il}\otimes F^{(0)}_{l,i+1}\psi^-_{i,0}+\\
  & (\vv-\vv^{-1})\sum_{k<i}\vv^{i-k-1}E^{(-1)}_{k,i+1}\otimes \wt{F}^{(0)}_{ik}\psi^-_{i,0}-
  (\vv-\vv^{-1})^2\sum_{l>i+1}^{k<i} \vv^{i-k-1} E^{(-1)}_{kl}\otimes \wt{F}^{(0)}_{ik}F^{(0)}_{l,i+1}\psi^-_{i,0},
\end{split}
\end{equation}
\begin{equation}\label{f_i1}
\begin{split}
  & \Delta(f_{i,1})=f_{i,1}\otimes 1 + \psi^+_{i,0}\otimes f_{i,1} +
  \vv^{-1}(\vv-\vv^{-1})\sum_{l>i+1} E^{(0)}_{i+1,l}\psi^+_{i,0}\otimes F^{(1)}_{li}-\\
  & (\vv-\vv^{-1})\sum_{k<i} \vv^{k-i} \wt{E}^{(0)}_{ki}\psi^+_{i,0}\otimes F^{(1)}_{i+1,k}-
  (\vv-\vv^{-1})^2\sum_{l>i+1}^{k<i}\vv^{k-i-1}E^{(0)}_{i+1,l}\wt{E}^{(0)}_{ki}\psi^+_{i,0}\otimes F^{(1)}_{lk},
\end{split}
\end{equation}
where for $1\leq j<i\leq n$ we set
\begin{equation}\label{long commutators}
\begin{split}
  & E^{(0)}_{ji}:=[e_{i-1,0},\cdots,[e_{j+1,0},e_{j,0}]_{\vv^{-1}}\cdots]_{\vv^{-1}}=
                 [\cdots[e_{i-1,0},e_{i-2,0}]_{\vv^{-1}},\cdots,e_{j,0}]_{\vv^{-1}},\\
  & F^{(0)}_{ij}:=[f_{j,0},\cdots,[f_{i-2,0},f_{i-1,0}]_{\vv}\cdots]_{\vv}=
                 [\cdots[f_{j,0},f_{j+1,0}]_{\vv},\cdots,f_{i-1,0}]_{\vv},\\
  & E^{(-1)}_{ji}:=[e_{i-1,0},\cdots,[e_{j+1,0},e_{j,-1}]_{\vv^{-1}}\cdots]_{\vv^{-1}}=
                  [[\cdots[e_{i-1,0},e_{i-2,0}]_{\vv^{-1}},\cdots,e_{j+1,0}]_{\vv^{-1}}, e_{j,-1}]_{\vv^{-1}},\\
  & F^{(1)}_{ij}:=[f_{j,1},[f_{j+1,0},\cdots,[f_{i-2,0},f_{i-1,0}]_{\vv}\cdots]_{\vv}]_{\vv}=
                 [\cdots[f_{j,1},f_{j+1,0}]_{\vv},\cdots,f_{i-1,0}]_{\vv},\\
  & \wt{E}^{(0)}_{ji}:=[e_{i-1,0},\cdots,[e_{j+1,0},e_{j,0}]_{\vv}\cdots]_{\vv}=
                      [\cdots[e_{i-1,0},e_{i-2,0}]_{\vv},\cdots,e_{j,0}]_{\vv},\\
  & \wt{F}^{(0)}_{ij}:=[f_{j,0},\cdots,[f_{i-2,0},f_{i-1,0}]_{\vv^{-1}}\cdots]_{\vv^{-1}}=
                      [\cdots[f_{j,0},f_{j+1,0}]_{\vv^{-1}},\cdots,f_{i-1,0}]_{\vv^{-1}}.
\end{split}
\end{equation}
\end{Thm}

The proof of this result is based on the RTT realization of $U_\vv(L\ssl_n)$ and
is presented in Appendix~\ref{Proof of coproduct on loop sln}.

\begin{Rem}
The right equalities in each of the lines of~(\ref{long commutators}) are not
obvious and are established during our proof of Theorem~\ref{Coproduct on quantum loop sln}.
They play an important role in the proof of Theorem~\ref{coproduct on antidominant sln} below.
\end{Rem}

Let $U^>_\vv(L\fg)$ and $U^\geq_\vv(L\fg)$ (resp.\ $U^<_\vv(L\fg)$ and
$U^\leq_\vv(L\fg)$) be the $\BC(\vv)$-subalgebras
of $U_\vv(L\fg)$ generated by $\{e_{i,r}\}_{i\in I}^{r\in\BZ}$ and
$\{e_{i,r},\psi^\pm_{i,\pm s}\}_{i\in I}^{r\in \BZ,s\in \BN}$ (resp.\
$\{f_{i,r}\}_{i\in I}^{r\in\BZ}$ and
$\{f_{i,r},\psi^\pm_{i,\pm s}\}_{i\in I}^{r\in \BZ,s\in \BN}$).

\begin{Cor}\label{useful}
For any $1\leq i<n$ and $r\in\BZ$, we have
\begin{equation*}
   \Delta(h_{i,\pm 1})-h_{i,\pm 1}\otimes 1-1\otimes h_{i,\pm 1}\in U^>_\vv(L\ssl_n)\otimes U^<_\vv(L\ssl_n),
\end{equation*}
\begin{equation*}
   \Delta(e_{i,r})-1\otimes e_{i,r}\in U^>_\vv(L\ssl_n)\otimes U^\leq_\vv(L\ssl_n),
\end{equation*}
\begin{equation*}
   \Delta(f_{i,r})-f_{i,r}\otimes 1\in U^\geq_\vv(L\ssl_n)\otimes U^<_\vv(L\ssl_n).
\end{equation*}
\end{Cor}

\begin{proof}
The claim is clear for $\Delta(h_{i,\pm 1}), \Delta(e_{i,-1}), \Delta(f_{i,1})$,
due to~(\ref{h_i1}--\ref{f_i1}). Applying iteratively
  $[\Delta(h_{i,\pm 1}),\Delta(e_{i,r})]=[2]_\vv\cdot \Delta(e_{i,r\pm 1}),
   [\Delta(h_{i,\pm 1}),\Delta(f_{i,r})]=-[2]_\vv\cdot \Delta(f_{i,r\pm 1})$,
we deduce the claim for $\Delta(e_{i,r})$ and $\Delta(f_{i,r})$.
\end{proof}


\subsection{Homomorphisms $\Delta_{\mu_1,\mu_2}$ for $\mu_1,\mu_2\in \Lambda^-,\ \fg=\ssl_n$}
\label{Coproduct on sl_n for antidominant}
$\ $

In this subsection, we construct homomorphisms
  $\Delta_{\mu_1,\mu_2}\colon \CU^\ssc_{0,\mu_1+\mu_2}\to \CU^\ssc_{0,\mu_1}\otimes \CU^\ssc_{0,\mu_2}$
for $\mu_1,\mu_2\in \Lambda^-$, which coincide with the Drinfeld-Jimbo coproduct
on $U_\vv(L\ssl_n)$ for $\mu_1=\mu_2=0$.
Set $b_{1,i}:=\alphavee_i(\mu_1)$ and $b_{2,i}:=\alphavee_i(\mu_2)$
(so that $b_{1,i}, b_{2,i}\in \BZ_{\leq 0}$).

\begin{Thm}\label{coproduct on antidominant sln}
For any $\mu_1,\mu_2\in \Lambda^-$, there is a unique $\BC(\vv)$-algebra homomorphism
\begin{equation*}
  \Delta_{\mu_1,\mu_2}\colon \CU^\ssc_{0,\mu_1+\mu_2}\longrightarrow \CU^\ssc_{0,\mu_1}\otimes \CU^\ssc_{0,\mu_2}
\end{equation*}
(we will denote $\Delta=\Delta_{\mu_1,\mu_2}$ when the algebras involved are clear), such that
\begin{equation*}
    \Delta(e_{i,r})=1\otimes e_{i,r},\
    \Delta(f_{i,s})=f_{i,s}\otimes 1 \
    \mathrm{for}\ b_{2,i}\leq r<0, b_{1,i}<s\leq 0,
\end{equation*}
\begin{equation*}
    \Delta(e_{i,0})=1\otimes e_{i,0}+e_{i,0}\otimes \psi^+_{i,0},\
    \Delta(f_{i,b_{1,i}})=f_{i,b_{1,i}}\otimes 1 + \psi^-_{i,b_{1,i}}\otimes f_{i,0},
\end{equation*}
\begin{equation*}
\begin{split}
    & \Delta(e_{i,b_{2,i}-1})=1\otimes e_{i,b_{2,i}-1}+e_{i,-1}\otimes \psi^-_{i,b_{2,i}}-
    (\vv-\vv^{-1})\sum_{l>i+1}E^{(-1)}_{il}\otimes F^{(0)}_{l,i+1}\psi^-_{i,b_{2,i}}+\\
    & (\vv-\vv^{-1})\sum_{k<i}\vv^{i-k-1}E^{(-1)}_{k,i+1}\otimes \wt{F}^{(0)}_{ik}\psi^-_{i,b_{2,i}}-
    (\vv-\vv^{-1})^2\sum_{l>i+1}^{k<i} \vv^{i-k-1} E^{(-1)}_{kl}\otimes \wt{F}^{(0)}_{ik}F^{(0)}_{l,i+1}\psi^-_{i,b_{2,i}},
\end{split}
\end{equation*}
\begin{equation*}
\begin{split}
    & \Delta(f_{i,1})=f_{i,1}\otimes 1 + \psi^+_{i,0}\otimes f_{i,1} +
    \vv^{-1}(\vv-\vv^{-1})\sum_{l>i+1} E^{(0)}_{i+1,l}\psi^+_{i,0}\otimes F^{(1)}_{li}-\\
    & (\vv-\vv^{-1})\sum_{k<i} \vv^{k-i} \wt{E}^{(0)}_{ki}\psi^+_{i,0}\otimes F^{(1)}_{i+1,k}-
    (\vv-\vv^{-1})^2\sum_{l>i+1}^{k<i} \vv^{k-i-1}E^{(0)}_{i+1,l}\wt{E}^{(0)}_{ki}\psi^+_{i,0}\otimes F^{(1)}_{lk},
\end{split}
\end{equation*}
\begin{equation*}
    \Delta((\psi^+_{i,0})^{\pm 1})=(\psi^+_{i,0})^{\pm 1}\otimes (\psi^+_{i,0})^{\pm 1},\
    \Delta((\psi^-_{i,b_{1,i}+b_{2,i}})^{\pm 1})=(\psi^-_{i,b_{1,i}})^{\pm 1}\otimes (\psi^-_{i,b_{2,i}})^{\pm 1},
\end{equation*}
\begin{equation*}
\begin{split}
    & \Delta(h_{i,1})=h_{i,1}\otimes 1+1\otimes h_{i,1}-
    (\vv^2-\vv^{-2})E_{i,i+1}^{(0)}\otimes F^{(1)}_{i+1,i}+
    (\vv-\vv^{-1})\sum_{l>i+1}E^{(0)}_{i+1,l}\otimes F^{(1)}_{l,i+1}+\\
    & (\vv-\vv^{-1})\sum_{k<i}\vv^{k+1-i}\wt{E}^{(0)}_{ki}\otimes F^{(1)}_{ik}+
    \vv^{-2}(\vv-\vv^{-1})\sum_{l>i+1}[E^{(0)}_{i,i+1},E^{(0)}_{i+1,l}]_{\vv^3}\otimes F^{(1)}_{li}-\\
    & (\vv-\vv^{-1})\sum_{k<i} \vv^{k-i-1} [E^{(0)}_{i,i+1},\wt{E}^{(0)}_{ki}]_{\vv^3}\otimes F^{(1)}_{i+1,k}+\\
    & (\vv-\vv^{-1})^2\sum_{l>i+1}^{k<i}\vv^{k-i}\left(E^{(0)}_{il}\wt{E}^{(0)}_{ki}-E^{(0)}_{i+1,l}\wt{E}^{(0)}_{k,i+1}\right)\otimes F^{(1)}_{lk},
\end{split}
\end{equation*}
\begin{equation*}
\begin{split}
    & \Delta(h_{i,-1})=h_{i,-1}\otimes 1+1\otimes h_{i,-1}+
    (\vv^2-\vv^{-2})E_{i,i+1}^{(-1)}\otimes F^{(0)}_{i+1,i}-
    (\vv-\vv^{-1})\sum_{l>i+1}E^{(-1)}_{i+1,l}\otimes F^{(0)}_{l,i+1}-\\
    & (\vv-\vv^{-1})\sum_{k<i}\vv^{i-k-1}E^{(-1)}_{ki}\otimes \wt{F}^{(0)}_{ik}-
    \vv^{2}(\vv-\vv^{-1})\sum_{l>i+1}E^{(-1)}_{il}\otimes [F^{(0)}_{l,i+1},F^{(0)}_{i+1,i}]_{\vv^{-3}}+\\
    & (\vv-\vv^{-1})\sum_{k<i} \vv^{i+1-k} E^{(-1)}_{k,i+1}\otimes [\wt{F}^{(0)}_{ik},F^{(0)}_{i+1,i}]_{\vv^{-3}}-\\
    & (\vv-\vv^{-1})^2\sum_{l>i+1}^{k<i}\vv^{i-k}E^{(-1)}_{kl}\otimes \left(\wt{F}^{(0)}_{i+1,k}F^{(0)}_{l,i+1}-\wt{F}^{(0)}_{ik}F^{(0)}_{li}\right),
\end{split}
\end{equation*}
where $E^{(0)}_{ji}, \wt{E}^{(0)}_{ji}, E^{(-1)}_{ji}, F^{(0)}_{ij}, \wt{F}^{(0)}_{ij}, F^{(1)}_{ij}$
are defined as in~(\ref{long commutators}).
\end{Thm}

The proof of this result is similar to our proof of Theorem~\ref{Coproduct for sl2}, but is much more tedious;
we sketch it in Appendix~\ref{Proof of coproduct antidominant sl_n}.

\begin{Rem}\label{Maps jmath for sln}
The similarity between the formulas for $\Delta_{\mu_1,\mu_2}$ of
Theorem~\ref{coproduct on antidominant sln} and $\Delta$ of
Theorem~\ref{Coproduct on quantum loop sln} can be explained via an
analogue of~Remark~\ref{Maps jmath}. To be more precise, let
$U^\pm_\vv$ be the positive/negative Borel subalgebras
in the Drinfeld-Jimbo presentation of $U_\vv(L\ssl_n)$,
while their analogues $\CU^{\ssc,\pm}_{0,\mu_1,\mu_2}$ (subalgebras of $\CU^\ssc_{0,\mu_1+\mu_2}$)
will be introduced in Appendix~\ref{Proof of coproduct antidominant sl_n}.
There are natural $\BC(\vv)$-algebra homomorphisms
  $\jmath^\pm_{\mu_1,\mu_2}\colon  U^\pm_\vv\to \CU^{\ssc,\pm}_{0,\mu_1,\mu_2}$,
see Proposition~\ref{jmath homomorphisms}.
According to Proposition~\ref{commutativity for jmath},
the following diagram is commutative:

\setlength{\unitlength}{1cm}
\begin{picture}(4,3.2)
 \put(3.3,2.5){$U^\pm_\vv$}
 \put(2.9,0.5){$\CU^{\ssc,\pm}_{0,\mu_1,\mu_2}$}
 \put(9.3,2.5){$U^\pm_\vv\otimes U^\pm_\vv$}
 \put(8.8,0.5){$\CU^{\ssc,\pm}_{0,\mu_1,0}\otimes \CU^{\ssc,\pm}_{0,0,\mu_2}$}

 \put(4.0,2.6){\vector (1,0){5.2}}
 \put(4.2,0.6){\vector (1,0){4.3}}
 \put(3.5,2.3){\vector (0,-1){1.3}}
 \put(10.1,2.3){\vector (0,-1){1.4}}

 \put(6.4,2.75){$\Delta$}
 \put(6.1,0.25){$\Delta_{\mu_1,\mu_2}$}
 \put(2.6,1.5){$\jmath^\pm_{\mu_1,\mu_2}$}
 \put(10.2,1.5){$\jmath^\pm_{\mu_1,0}\otimes \jmath^\pm_{0,\mu_2}$}
\end{picture}
\end{Rem}


\subsection{Homomorphisms $\Delta_{\mu_1,\mu_2}$ for arbitrary
$\mu_1,\mu_2\in \Lambda,\ \fg=\ssl_n$}
\label{General shifted coproduct sl_n}
$\ $

Let us first generalize the {\em shift homomorphisms} of~Proposition~\ref{Maps iota}.

\begin{Lem}\label{Maps iota for sln}
For any $\mu\in \Lambda$ and $\nu_1,\nu_2\in \Lambda^-$, there is a
unique $\BC(\vv)$-algebra homomorphism
$\iota_{\mu,\nu_1,\nu_2}\colon  \CU^\ssc_{0,\mu}\to \CU^\ssc_{0,\mu+\nu_1+\nu_2}$,
which maps the currents as follows
\begin{multline*}
  \iota_{\mu,\nu_1,\nu_2}\colon
  e_i(z)\mapsto (1-z^{-1})^{-\alphavee_i(\nu_1)}e_i(z),\
  f_i(z)\mapsto (1-z^{-1})^{-\alphavee_i(\nu_2)}f_i(z),\\
    \psi^\pm_i(z)\mapsto (1-z^{-1})^{-\alphavee_i(\nu_1+\nu_2)}\psi^\pm_i(z).
\end{multline*}
\end{Lem}

\begin{proof}
The proof is analogous to that of Proposition~\ref{Maps iota}.
\end{proof}

The proof of the following technical result is presented in Appendix~\ref{Proof of injectivity}
and is based on the shuffle realization of the quantum loop algebra $U_\vv(L\ssl_n)$, see~\cite{n} (cf.~\cite{t3}).

\begin{Thm}\label{1026}
The homomorphism $\iota_{\mu,\nu_1,\nu_2}$ is injective for any
$\mu\in \Lambda,\ \nu_1,\nu_2\in \Lambda^-$.
\end{Thm}

Combining this theorem with Corollary~\ref{useful} and our arguments from
the proof of Theorem~\ref{Coproduct for sl2 in general},
we get the key result of this section.

\begin{Thm}\label{Coproduct for sln in dominant case}
For any $\mu_1,\mu_2\in \Lambda$ and $\mu:=\mu_1+\mu_2$, there is a unique
$\BC(\vv)$-algebra homomorphism
\begin{equation*}
  \Delta_{\mu_1,\mu_2}\colon  \CU^\ssc_{0,\mu}\longrightarrow \CU^\ssc_{0,\mu_1}\otimes \CU^\ssc_{0,\mu_2},
\end{equation*}
such that for any $\nu_1,\nu_2\in \Lambda^-$ the following diagram is commutative:

\setlength{\unitlength}{1cm}
\begin{picture}(4,3.2)
 \put(3.2,2.5){$\CU^\ssc_{0,\mu}$}
 \put(2.7,0.5){$\CU^\ssc_{0,\mu+\nu_1+\nu_2}$}
 \put(9.0,2.5){$\CU^\ssc_{0,\mu_1}\otimes \CU^\ssc_{0,\mu_2}$}
 \put(8.5,0.5){$\CU^\ssc_{0,\mu_1+\nu_1}\otimes \CU^\ssc_{0,\mu_2+\nu_2}$}

 \put(4.1,2.6){\vector (1,0){4.6}}
 \put(4.5,0.6){\vector (1,0){3.6}}
 \put(3.5,2.3){\vector (0,-1){1.4}}
 \put(10.1,2.3){\vector (0,-1){1.4}}

 \put(6.0,2.75){$\Delta_{\mu_1,\mu_2}$}
 \put(5.6,0.25){$\Delta_{\mu_1+\nu_1,\mu_2+\nu_2}$}
 \put(2.3,1.5){$\iota_{\mu,\nu_2,\nu_1}$}
 \put(10.15,1.5){$\iota_{\mu_1,0,\nu_1}\otimes \iota_{\mu_2,\nu_2,0}$}
\end{picture}
\end{Thm}

The following is proved analogously to Corollary~\ref{uniform coproduct h_1}:

\begin{Prop}
For arbitrary $\mu_1,\mu_2\in \Lambda$, the images $\Delta_{\mu_1,\mu_2}(h_{i,\pm 1})$
are given by formulas~(\ref{h_i1}) and~(\ref{h_i-1}).
\end{Prop}


\subsection{Open Problems}
\label{open problems}
$\ $

Following~\cite{fkp}, we expect that homomorphisms
  $\Delta_{\mu_1,\mu_2}\colon \CU^\ssc_{0,\mu_1+\mu_2}\to \CU^\ssc_{0,\mu_1}\otimes \CU^\ssc_{0,\mu_2}$
(specializing to the Drinfeld-Jimbo coproduct for $\mu_1=\mu_2=0$)
exist for any simply-laced Lie algebra $\fg$ and its two coweights $\mu_1,\mu_2\in \Lambda$.
Moreover, their construction should proceed in the same way as for the
aforementioned case $\fg=\ssl_n$. To be more precise, for antidominant $\mu_1,\mu_2\in \Lambda^-$,
we expect that the homomorphism $\Delta_{\mu_1,\mu_2}$
is characterized by the following two properties:

(a) $\Delta_{\mu_1,\mu_2}(e_{i,r})=1\otimes e_{i,r},\ \Delta_{\mu_1,\mu_2}(f_{i,s})=f_{i,s}\otimes 1$
for $\alphavee_i(\mu_2)\leq r<0, \alphavee_i(\mu_1)<s\leq 0$;

(b) an analogue of the commutative diagram of Remark~\ref{Maps jmath for sln} holds.

\noindent
For general $\mu_1,\mu_2$, we expect that the construction of $\Delta_{\mu_1,\mu_2}$
should be easily deduced from the antidominant case with the help of
{\em shift homomorphisms}
$\iota_{\mu,\nu_1,\nu_2}\ (\mu\in \Lambda, \nu_1,\nu_2\in \Lambda^-)$ as
in~Theorems~\ref{Coproduct for sl2 in general},~\ref{Coproduct for sln in dominant case}.

The outlined construction of $\Delta_{\mu_1,\mu_2}$ for a general $\fg$ lacks
explicit formulas for the Drinfeld-Jimbo coproduct of
$\{e_{i,0}, e_{i,-1}, f_{i,0}, f_{i,1}, \psi^\pm_{i,0}, h_{i,\pm 1}\}_{i\in I}$--the generators of $U_\vv(L\fg)$,
similar to those of Lemma~\ref{Explicit sl2} and Theorem~\ref{Coproduct on quantum loop sln}.

\section{Ubiquity of RTT relations}
\label{Ubiquity}


\subsection{Rational Lax matrix}
$\ $

Before we proceed to the \emph{trigonometric} setting, let us recall
the classical relation between rational Lax matrices and type $A$ quantum
open Toda systems, which goes back to~\cite{ft}.

Let $R_\rat(z)\in \mathrm{End}(\BC^2\otimes \BC^2)$ be the standard
rational $R$-matrix:
\begin{equation*}
    R_\rat(z)=\mathrm{Id}+\frac{\hbar}{z}P,\
    \mathrm{where}\  P\in \mathrm{End}(\BC^2\otimes \BC^2)\ \mathrm{is\ the\ permutation\ map}.
\end{equation*}
Let $\hat{\CA}^\hbar_n$ be the associative $\BC[\hbar]$-algebra
generated by $\{\sfu^{\pm 1}_i,w_i\}_{i=1}^n$
with the defining relations
  $[\sfu_i,\sfu_j]=[w_i,w_j]=0, \sfu^{\pm 1}_i\sfu^{\mp 1}_i=1,
   [\sfu_i,w_j]=\delta_{ij}\hbar\sfu_i.$
Define the \emph{(local) rational Lax matrix}
\eq{Classical Lax}
  L^\hbar_i(z)=\left(%
  \begin{array}{cc}
    z-w_i & \sfu_i^{-1} \\
    -\sfu_i & 0 \\
  \end{array}%
  \right)\in \mathrm{Mat}(2,\hat{\CA}^\hbar_n[z])
\end{equation}
and introduce the \emph{complete monodromy matrix}
$T^\hbar_n(z):=L^\hbar_n(z)\cdots L^\hbar_1(z)$.
Then, the monodromy matrix $T^\hbar_n(z)$ satisfies the rational
RTT-relation:
\begin{equation*}
    R_\rat(z-w)(T^\hbar_n(z)\otimes 1)(1\otimes T^\hbar_n(w))=
    (1\otimes T^\hbar_n(w))(T^\hbar_n(z)\otimes 1)R_\rat(z-w).
\end{equation*}
Due to this relation, the coefficients (in $z$) of the matrix
element $T^\hbar_n(z)_{11}$ generate a commutative subalgebra of
$\hat{\CA}^\hbar_n$, known as the quantum open Toda system of $\gl_n$. The
coefficient of $z^{n-2}$ equals
\eq{sH2}
  {\mathsf H}_2^{\on{rat}}=\frac{1}{2}\left(\sum_{i=1}^nw_i\right)^2-
  \frac{1}{2}\sum_{i=1}^nw^2_i-\sum_{i=1}^{n-1}\sfu_i\sfu_{i+1}^{-1}.
\end{equation}
We recover the standard quantum open Toda hamiltonian of $\ssl_n$ once we set
$w_1+\ldots+w_n=0$.


\subsection{Trigonometric/Relativistic Lax matrices}\label{relativistic}
$\ $

Let $R_\trig(z)\in \mathrm{End}(\BC^2\otimes \BC^2)$
be the standard trigonometric $R$-matrix (see~\cite[(3.7)]{df}):
\begin{equation}\label{R-matrix}
  R_\trig(z)=\left(%
    \begin{array}{cccc}
      1 & 0 & 0 & 0 \\
      0 & \frac{z-1}{\vv z-\vv^{-1}} & \frac{z(\vv-\vv^{-1})}{\vv z-\vv^{-1}} & 0 \\
      0 & \frac{\vv-\vv^{-1}}{\vv z-\vv^{-1}} & \frac{z-1}{\vv z-\vv^{-1}} & 0 \\
      0 & 0 & 0 & 1 \\
    \end{array}%
  \right).
\end{equation}
Let $\hat{\CA}^\vv_n$ be the associative $\BC(\vv)$-algebra
generated by $\{\wt\sw^{\pm 1}_i, D^{\pm 1}_i\}_{i=1}^n$ with the
defining relations
  $[\wt\sw_i,\wt\sw_j]=[D_i,D_j]=0,\
   \wt\sw^{\pm 1}_i\wt\sw^{\mp 1}_i=D^{\pm 1}_iD^{\mp 1}_i=1,\
   D_i\wt\sw_j=\vv^{\delta_{ij}}\wt\sw_jD_i$.
If we set $\sw_i^{\pm1}=\wt\sw_i^{\pm2}$, we see that $\hat{\CA}^\vv_n$ is
a particular example of the algebras $\hat{\CA}^\vv_\fra$ of
Section~\ref{Section with Homomorphism}.
Define the \emph{(local) relativistic Lax matrix}
\eq{hatL}
  L^{\vv,0}_i(z)=
  \left(%
  \begin{array}{cc}
    \wt\sw_i^{-1}z^{1/2}-\wt\sw_i z^{-1/2} & D^{-1}_iz^{1/2} \\
    -D_iz^{-1/2} & 0 \\
  \end{array}%
  \right)\in{\on{Mat}}(2,z^{-1/2}\hat{\CA}^\vv_n[z])
\end{equation}
and introduce the \emph{complete monodromy matrix}
$T^{\vv,0}_n(z):=L^{\vv,0}_n(z)\cdots L^{\vv,0}_1(z)$.

\begin{Lem}\label{RTT_monodromy}
The monodromy matrix $T^{\vv,0}_n(z)$ satisfies the trigonometric RTT-relation:
\begin{equation*}
   R_\trig(z/w)(T^{\vv,0}_n(z)\otimes 1)(1\otimes T^{\vv,0}_n(w))=
   (1\otimes T^{\vv,0}_n(w))(T^{\vv,0}_n(z)\otimes 1)R_\trig(z/w).
\end{equation*}
\end{Lem}

\begin{proof}
It suffices to check the above relation for $n=1$. The proof in the
latter case is straightforward.
\end{proof}

\begin{Cor}\label{relativistic hamiltonian}
The coefficients (in $z$) of the matrix element
$z^{n/2}T^{\vv,0}_n(z)_{11}$ generate a commutative subalgebra of
$\hat{\CA}^\vv_n$. The coefficient of $z$ equals
\eq{hatH2}
  {\mathsf H}^{0}_2=(-1)^{n-1}\wt\sw_1\cdots\wt\sw_n\cdot
  \left(\sum_{i=1}^n\wt\sw_i^{-2}+\sum_{i=1}^{n-1}\wt\sw_i^{-1}\wt\sw_{i+1}^{-1}D_iD^{-1}_{i+1}\right).
\end{equation}
This hamiltonian is equivalent to the quadratic hamiltonian of the
$q$-difference quantum Toda lattice of~\cite[(5.7)]{e} (see
also~\cite{se}) once we set $\wt\sw_1\cdots \wt\sw_n=1$.
\end{Cor}

\begin{Rem}\label{Kuznetsov-Tsyganov}
The notion of a \emph{relativistic Lax matrix} goes back
to~\cite{kt}. In particular, our choice of $L^{\vv,0}_i(z)$ is a
slight variation of their construction, which is adapted to a
different choice of the trigonometric $R$-matrix.
\end{Rem}

Now let us consider two \emph{(local) trigonometric Lax matrices}
\eq{relativistic Lax}
  L^{\vv,-1}_i(z)=\left(%
  \begin{array}{cc}
    \wt\sw_i^{-1}-\wt\sw_iz^{-1} & \wt\sw_iD_i^{-1} \\
    -\wt\sw_iD_iz^{-1} & \wt\sw_i \\
  \end{array}%
  \right)\in \mathrm{Mat}(2, z^{-1}\hat{\CA}^\vv_n[z]),
\end{equation}
\eq{new Lax matrix}
  L^{\vv,1}_i(z)=\left(%
  \begin{array}{cc}
    \wt\sw_i^{-1}z-\wt\sw_i & \wt\sw_i^{-1}D_i^{-1}z \\
    -\wt\sw_i^{-1}D_i & -\wt\sw_i^{-1} \\
  \end{array}%
  \right)\in \mathrm{Mat}(2, \hat{\CA}^\vv_n[z]).
\end{equation}

\begin{Lem}\label{RTT_monodromy_second}
The Lax matrices $L^{\vv,\pm 1}_i(z)$ satisfy the trigonometric RTT-relation:
\begin{equation*}
    R_\trig(z/w)(L^{\vv,\pm 1}_i(z)\otimes 1)(1\otimes L^{\vv,\pm 1}_i(w))=
    (1\otimes L^{\vv,\pm 1}_i(w))(L^{\vv,\pm 1}_i(z)\otimes 1)R_\trig(z/w).
\end{equation*}
\end{Lem}

\begin{proof}
The proof is straightforward.
\end{proof}

%
%


\subsection{Mixed Toda hamiltonians}\label{section mixed Toda}
$\ $

Now we construct $3^n$ Hamiltonians generalizing ${\mathsf H}^0_2$
in spirit of~\cite[(90)]{fg}, cf.\ also~\cite[(1.1) and~Section~2]{bg}.
For any $\vec{k}=(k_n,\ldots,k_1)\in \{-1,0,1\}^n$, define
the {\em mixed complete monodromy matrix}
\begin{equation*}
  T^\vv_{\vec{k}}(z):=L^{\vv,k_n}_n(z)\cdots L^{\vv,k_1}_{1}(z).
\end{equation*}
In particular, $T^\vv_{\vec{0}}(z)=T^{\vv,0}_n(z)$.
Since all three matrices $L^{\vv,-1}_i(z), L^{\vv,0}_i(z), L^{\vv,1}_i(z)$
satisfy the RTT-relation with the same $R$-matrix
$R_\trig(z)$, the same is true for $T^\vv_{\vec{k}}(z)$. Hence, the
coefficients (in $z$) of the matrix element
$T^\vv_{\vec{k}}(z)_{11}$ generate a commutative subalgebra of
$\hat{\CA}^\vv_n$. We have
\begin{equation*}
   T^\vv_{\vec{k}}(z)_{11}=
   {\mathsf H}^{\vec{k}}_1z^s+{\mathsf H}^{\vec{k}}_2z^{s+1}+\mathrm{higher\ powers\ of}\ z,
\end{equation*}
where $s=\sum_{i=1}^n \frac{k_i-1}{2}$.
Here
  ${\mathsf H}^{\vec{k}}_1=(-1)^n\wt\sw_1\cdots\wt\sw_n$,
while the hamiltonian ${\mathsf H}^{\vec{k}}_2$ equals
\eq{mixedH2}
  {\mathsf H}^{\vec{k}}_2=
  (-1)^{n-1}\wt\sw_1\cdots\wt\sw_n\cdot \left(\sum_{i=1}^n \wt\sw^{-2}_i+
  \sum_{i=1}^{n-1} \sigma_{i,i+1}D_iD^{-1}_{i+1}+
  \sum_{1\leq i<j-1< n}^{k_{i+1}=\ldots=k_{j-1}=1} \sigma_{i,j}D_iD^{-1}_j\right),
\end{equation}
where $\sigma_{i,j}:=\wt\sw_i^{-k_i-1}\wt\sw_{i+1}^{-k_{i+1}-1}\cdots\wt\sw_j^{-k_j-1}$.

\begin{Rem}
At the classical level, the birational B\"acklund-Darboux transformations
interchanging various hamiltonians ${\mathsf H}^{\vec{k}}_2$
are given in~\cite[Theorem~6.1]{gsv}.
\end{Rem}

\begin{Lem}\label{conjugation}
For any $\vec{k}$, set $\vec{k}'=(0,k_{n-1},\ldots,k_2,0)$.
Then, ${\mathsf H}^{\vec{k}}_2$ is equivalent to ${\mathsf H}^{\vec{k}'}_2$.
\end{Lem}

\begin{proof}
It is straightforward to see that
  ${\mathsf H}^{\vec{k}'}_2=\mathrm{Ad}(F(\wt{\sw}_1,\ldots,\wt{\sw}_n)){\mathsf H}^{\vec{k}}_2,$
where
  $F(\wt{\sw}_1,\ldots,\wt{\sw}_n)=\exp(k_1f_-(\log(\wt{\sw}_1))+k_nf_+(\log(\wt{\sw}_n)))$
with
  $f_\pm(t)=\pm\frac{t^2}{2\log(\vv)}+\frac{t}{2}$.
\end{proof}

\rem{Et-Se}
It follows that among the aforementioned $3^n$ \emph{mixed Toda hamiltonians}
${\mathsf H}^{\vec{k}}_2$, parameterized by $\vec{k}\in \{-1,0,1\}^n$,
there are no more than $3^{n-2}$ different up to equivalence.
In~\cite{fgt} these hamiltonians are identified
with the modified versions of the $q$-Toda hamiltonian
in~\cite{e,se}, which now depend on a choice of two orientations of the
Dynkin diagram of type $A_{n-1}$ (equivalently, a choice of a pair of Coxeter elements).
There are $4^{n-2}$ such choices, but some of them are equivalent leading to exactly $3^{n-2}$
inequivalent hamiltonians, which turn out to be equivalent to the aforementioned
${\mathsf H}^{\vec{k}}_2$. All the $q$-Toda hamiltonians
of~\cite{e,se} correspond to the pairs of coinciding orientations, i.e.\ to
$\vec{k}=(0,\ldots,0)$, and they share the same eigenfunction
$J$~\cite[Section 3]{fjm},
while our mixed Toda hamiltonians do not admit the common eigenfunctions.
We are grateful to P.~Etingof for his suggestion to study the construction
of~\cite{se} for pairs of different orientations.
\erem

%
%


\subsection{Shifted RTT algebras of $\ssl_2$}
\label{section RTT algebras}
$\ $

Fix $n\in \BN$. Following~\cite{df} (cf. also Remark~\ref{4vs3 RTT}), we introduce the
\emph{(trigonometric) shifted RTT algebras of $\ssl_2$}, denoted by
$\CU^{\rtt}_{0,-2n}$. These are associative $\BC(\vv)$-algebras generated by
\begin{equation*}
\begin{split}
  & \{t^+_{11}[r], t^+_{12}[r], t^+_{21}[r+1], t^+_{22}[r],
    t^-_{11}[-m], t^-_{12}[-m-1], t^-_{21}[-m], t^-_{22}[-m-1+\delta_{n,0}]\}_{r\geq 0}^{m\geq -n} \cup\\
  & \{(t^+_{11}[0])^{-1}, (t^-_{11}[n])^{-1}\}
\end{split}
\end{equation*}
subject to the following defining relations:
\begin{equation}\tag{R1}\label{RTT1}
  (t^+_{11}[0])^{\pm 1}(t^+_{11}[0])^{\mp 1}=1,\
  (t^-_{11}[n])^{\pm 1}(t^-_{11}[n])^{\mp 1}=1,
\end{equation}
\begin{equation}\tag{R2}\label{RTT2}
  R_\trig(z/w)(T^\epsilon(z)\otimes 1)(1\otimes T^{\epsilon'}(w))=
  (1\otimes T^{\epsilon'}(w))(T^\epsilon(z)\otimes 1)R_\trig(z/w),
\end{equation}
\begin{equation}\tag{R3}\label{RTT3}
  \qdet\ T^\pm(z)=1
\end{equation}
for all $\epsilon,\epsilon'\in \{\pm\}$, where the two-by-two
matrices $T^\pm(z)$ are given by
\begin{equation*}
   T^\pm(z)=
   \left(%
   \begin{array}{cc}
     T^\pm_{11}(z) & T^\pm_{12}(z) \\
     T^\pm_{21}(z) & T^\pm_{22}(z) \\
   \end{array}%
   \right)\ \mathrm{with}\
    T^\pm_{ij}(z):=\sum_{r}t^\pm_{ij}[r]z^{-r},
\end{equation*}
and the quantum determinant $\qdet$ is defined in a standard way
as\footnote{It is instructive to point out the difference with~\cite{m}, where
the author uses a different trigonometric $R$-matrix given by
$R^{\M}_\trig(z/w)=(R_\trig(z/w)^t)^{-1}$ as well as
${T^{\M,\pm}(z)}=T^\pm(z)^t$. For this reason, the
quantum determinant $\qdet^\M$ of~\cite[Exercise~1.6]{m} is consistent
with our definition of $\qdet$, that is,
  $\qdet^{\M}\ T^{\M,\pm}(z):=T^{\M,\pm}_{11}(z)T^{\M,\pm}_{22}(\vv^{-2}z)-
   \vv^{-1}T^{\M,\pm}_{21}(z)T^{\M,\pm}_{12}(\vv^{-2}z)=\qdet\ T^{\pm}(z)$.}
\begin{equation*}
  \qdet\ T^\pm(z):=T^\pm_{11}(z)T^\pm_{22}(\vv^{-2}z)-\vv^{-1}T^\pm_{12}(z)T^\pm_{21}(\vv^{-2}z).
\end{equation*}

Note that $T^\pm(z)$ admits the following unique \emph{Gauss} decomposition:
\begin{equation*}
    T^\pm(z)=
    \left(%
    \begin{array}{cc}
      1 & 0 \\
      \tilde{f}^\pm(z) & 1 \\
    \end{array}%
    \right)
    \left(%
    \begin{array}{cc}
      \tilde{g}^\pm_1(z) & 0 \\
      0 & \tilde{g}^\pm_2(z) \\
    \end{array}%
    \right)
    \left(%
    \begin{array}{cc}
      1 & \tilde{e}^\pm(z) \\
      0 & 1 \\
    \end{array}%
    \right),
\end{equation*}
where coefficients of the \emph{half-currents}
$\tilde{e}^\pm(z), \tilde{f}^\pm(z), \tilde{g}^\pm_1(z), \tilde{g}^\pm_2(z)$
are elements of $\CU^\rtt_{0,-2n}$.

To establish the relation between $\CU^\rtt_{0,-2n}$ and
$\CU^\ad_{0,-2n}$ (adjoint version of the shifted quantum affine
algebra of $\ssl_2$), recall \emph{Drinfeld half-currents}
$e^\pm(z), f^\pm(z)$ of~(\ref{Drinfeld currents}).


\begin{Thm}\label{Ding-Frenkel}
(a) The currents $\tilde{g}^\pm_1(z), \tilde{g}^\pm_2(z)$ pairwise commute and satisfy
\begin{equation*}
  \tilde{g}^\pm_2(z)\tilde{g}^\pm_1(\vv^{-2}z)=1.
\end{equation*}

\noindent
(b) There exists a unique $\BC(\vv)$-algebra homomorphism
  $\Upsilon_{0,-2n}\colon  \CU^\ad_{0,-2n}\to \CU^\rtt_{0,-2n}$,
defined by
\begin{equation*}
\begin{split}
  & e^\pm(z)\mapsto \tilde{e}^\pm(z)/(\vv-\vv^{-1}),\
  f^\pm(z)\mapsto \tilde{f}^\pm(z)/(\vv-\vv^{-1}),\\
  & \psi^\pm(z)\mapsto \tilde{g}^\pm_2(z)\tilde{g}^\pm_1(z)^{-1},\
  (\phi^+)^{\pm 1}\mapsto (t^+_{11}[0])^{\mp 1},\
  (\phi^-)^{\pm 1}\mapsto \vv^{\mp n}(t^-_{11}[n])^{\mp 1}.
\end{split}
\end{equation*}

\noindent
(c) For any $b_1,b_2\in \BZ_{\leq 0}$, there exists a unique $\BC(\vv)$-algebra
homomorphism
\begin{equation*}
  \Delta^\rtt_{2b_1,2b_2}\colon \CU^\rtt_{0,2b_1+2b_2}\longrightarrow \CU^\rtt_{0,2b_1}\otimes \CU^\rtt_{0,2b_2},
\end{equation*}
defined by $T^\pm(z)\mapsto T^\pm(z)\otimes T^\pm(z)$.
\end{Thm}

\begin{Rem}
The $n=0$ case of this theorem was proved in~\cite{df}, cf.~Remark~\ref{4vs3 RTT}.
\end{Rem}

\begin{proof}
The verification of part (b) is analogous to the one for $n=0$,
dealt with in~\cite{df}. Once (b) is established, it is easy to see
that $\qdet\ T^\pm(z)=\tilde{g}^\pm_2(z)\tilde{g}^\pm_1(\vv^{-2}z)$,
hence (a). It is clear that $\Delta^\rtt_{2b_1,2b_2}$ is well-defined on
the generators. The compatibility of $\Delta^\rtt_{2b_1,2b_2}$ with
the defining relations~(\ref{RTT1}--\ref{RTT3}) is checked
analogously to the case $n=0$.
\end{proof}

Recall the generating series $A^\pm(z), B^\pm(z), C^\pm(z),
D^\pm(z)$ with coefficients in $\CU^\ad_{0,-2n}$, introduced in
Section~\ref{ABCD generators}.

\begin{Cor}\label{ABCD in rtt}
The homomorphism $\Upsilon_{0,-2n}$ maps these generating series as follows:
\begin{equation*}
    A^+(z)\mapsto T^+_{11}(z), B^+(z)\mapsto T^+_{12}(z),
    C^+(z)\mapsto T^+_{21}(z), D^+(z)\mapsto T^+_{22}(z),
\end{equation*}
\begin{equation*}
    A^-(z)\mapsto (\vv z)^nT^-_{11}(z), B^-(z)\mapsto (\vv z)^nT^-_{12}(z),
    C^-(z)\mapsto (\vv z)^nT^-_{21}(z), D^-(z)\mapsto (\vv z)^nT^-_{22}(z).
\end{equation*}
\end{Cor}

\begin{proof}
Due to Theorem~\ref{Ding-Frenkel}(a, b), we have
\begin{equation*}
    \Upsilon_{0,-2n}(\psi^\pm(z))=1/\tilde{g}^\pm_1(z)\tilde{g}^\pm_1(\vv^{-2}z),\
    \Upsilon_{0,-2n}((\phi^+)^{-1})=t^+_{11}[0],\ \Upsilon_{0,-2n}((\phi^-)^{-1})=\vv^nt^-_{11}[n].
\end{equation*}
Combining this with
  $\psi^+(z)=\frac{1}{A^+(z)A^+(\vv^{-2}z)},
   \psi^-(z)=\frac{z^{2n}}{A^-(z)A^-(\vv^{-2}z)},
   \ \mathrm{and}\ A^\pm_0=(\phi^\pm)^{-1}$, we get
  $\Upsilon_{0,-2n}(A^+(z))=\tilde{g}^+_1(z)=T^+_{11}(z),\
   \Upsilon_{0,-2n}(A^-(z))=(\vv z)^n\tilde{g}^-_1(z)=(\vv z)^nT^-_{11}(z)$.
The computation of the images of the remaining generating series is
straightforward, e.g.\\
  $\Upsilon_{0,-2n}(B^-(z))=(\vv-\vv^{-1})\Upsilon_{0,-2n}(A^-(z))\Upsilon_{0,-2n}(e^-(z))=
   (\vv z)^n\tilde{g}^-_1(z)\tilde{e}^-(z)=(\vv z)^nT^-_{12}(z)$.
\end{proof}

The following is the key result of this subsection.

\begin{Thm}\label{Ding-Frenkel 2}
For $n\in\BN$, $\Upsilon_{0,-2n}\colon \CU^\ad_{0,-2n}\to \CU^\rtt_{0,-2n}$
is an isomorphism of $\BC(\vv)$-algebras.
\end{Thm}

\begin{proof}
Due to Theorem~\ref{Ding-Frenkel} and Corollary~\ref{ABCD in rtt},
it suffices to prove that there exists a $\BC(\vv)$-algebra
homomorphism $\CU^\rtt_{0,-2n}\to\CU^\ad_{0,-2n}$, such that
\begin{equation}\label{inverse}
\begin{split}
  & (t^+_{11}[0])^{-1}\mapsto \phi^+,\ (t^-_{11}[n])^{-1}\mapsto \vv^n\phi^-,\\
  &T^+_{11}(z)\mapsto A^+(z),\ T^+_{12}(z)\mapsto B^+(z),\
   T^+_{21}(z)\mapsto C^+(z),\ T^+_{22}(z)\mapsto D^+(z), \\
  & T^-_{11}(z)\mapsto (\vv z)^{-n}A^-(z),\ T^-_{12}(z)\mapsto (\vv z)^{-n}B^-(z),\\
  & T^-_{21}(z)\mapsto (\vv z)^{-n}C^-(z),\ T^-_{22}(z)\mapsto (\vv z)^{-n}D^-(z).
\end{split}
\end{equation}
This amounts to verifying that the assignment~(\ref{inverse})
preserves defining relations~(\ref{RTT1}--\ref{RTT3}). Relation~(\ref{RTT1})
is preserved, due to
$A^\pm_0\phi^\pm=\phi^\pm A^\pm_0=1$, while~(\ref{RTT3}) is
preserved, due to relation~(\ref{ABCD12}). Finally,~(\ref{RTT2}) is an
equality in $\mathrm{End}(\BC^2\otimes \BC^2)\otimes\CU^\ad_{0,-2n}$ and thus
can be viewed as a collection of $16$ relations in $\CU^\ad_{0,-2n}$ for each
choice of $\epsilon,\epsilon'\in \{\pm\}$. It is straightforward to see that
$6$ of these relations follow from the rest, while the remaining
$10$ relations exactly match the $10$ relations of~(\ref{ABCD1},~\ref{ABCD5}--\ref{ABCD11})
under the assignment~(\ref{inverse}).
\end{proof}

\begin{Rem}\label{RTT general shifts}
The results of this subsection admit natural generalizations to the case of arbitrary
$b_1,b_2\in \BZ_{\leq 0}$ such that $b_1+b_2$ is even. In other words, one can define an
analogous shifted RTT algebra of $\ssl_2$, denoted $\CU^{\rtt}_{b_1,b_2}$, and construct
a $\BC(\vv)$-algebra isomorphism
$\Upsilon_{b_1,b_2}\colon \CU^\ad_{b_1,b_2}\iso \CU^\rtt_{b_1,b_2}$.
This observation is used in Remark~\ref{explanation of other Lax} below, where we provide an
alternative interpretation of the Lax matrices
$L^{\vv,-1}_1(z), L^{\vv,0}_1(z), L^{\vv,1}_1(z)$ from Section~\ref{relativistic}.
\end{Rem}


\subsection{Relation between two different appearances of RTT}
\

Recall the local trigonometric Lax matrix $L^{\vv,-1}_1(z)$ of~\refe{relativistic Lax}.
Combining the equality $\qdet\ L^{\vv,-1}_1(z)=1$ with
Lemma~\ref{RTT_monodromy_second}, we see that $L^{\vv,-1}_1(z)$ gives rise to
an algebra homomorphism
  $\Phi^\rtt_{0,-2}\colon  \CU^\rtt_{0,-2}\to
   \hat{\CA}^\vv_1$ defined by $T^\pm(z)\mapsto L^{\vv,-1}_1(z)$.
Recall the homomorphism
  $\wt{\Phi}^0_{-2}\colon \CU^{\ad}_{0,-2}\to \hat{\CA}^\vv_1$
of Theorem~\ref{Homomorphism} (where $\sw_1^{1/2}=\wt{\sw}_1$). The following is straightforward.

\begin{Lem}\label{Lax via homomorphism 1}
The composition $\Phi^\rtt_{0,-2}\circ \Upsilon_{0,-2}$ coincides with $\wt{\Phi}^0_{-2}$.
\end{Lem}

\begin{Rem}\label{explanation of other Lax}
Let us provide a similar interpretation of the other two Lax matrices
$L^{\vv,0}_1(z)$ and $L^{\vv,1}_1(z)$. Recall that the algebras $\CU^{\ad}_{0,-2}$
and $\CU^{\ad}_{b,-2-b}$ are isomorphic for any $b\in \BZ$.
In particular, one can pull-back the homomorphism $\wt{\Phi}^0_{-2}$ to obtain a
homomorphism $\wt{\Phi}_{b,-2-b}\colon \CU^{\ad}_{b,-2-b}\to \hat{\CA}^\vv_1$, explicitly given by
\begin{equation*}
\begin{split}
    & e(z)\mapsto \frac{\wt\sw^{2+b}}{\vv-\vv^{-1}} \delta\left(\frac{\wt\sw^2}{z}\right)D^{-1},\
    f(z)\mapsto \frac{\wt\sw^{b}}{1-\vv^2}\delta\left(\frac{\vv^2\wt\sw^2}{z}\right)D,\\
    & \psi^\pm(z)\mapsto \left(\frac{\vv^{-b}\wt\sw^2 z^b}{(1-\wt\sw^2/z)(1-\vv^2\wt\sw^2/z)}\right)^\pm,\
    (\phi^+)^{\pm 1}\mapsto \vv^{\mp b/2}\wt{\sw}^{\pm 1},\
    (\phi^-)^{\pm 1}\mapsto -\vv^{\mp (b/2+1)}\wt{\sw}^{\mp 1}.
\end{split}
\end{equation*}
Due to Remark~\ref{RTT general shifts}, the algebra $\CU^{\ad}_{b,-2-b}$ admits an RTT realization,
that is there is an isomorphism $\Upsilon_{b,-2-b}\colon  \CU^{\ad}_{b,-2-b}\iso \CU^{\rtt}_{b,-2-b}$, only for $b=0,-1,-2$.
Analogously to Lemma~\ref{Lax via homomorphism 1}, recasting the homomorphisms $\wt{\Phi}_{b,-2-b}$
as the homomorphisms $\CU^{\rtt}_{b,-2-b}\to \hat{\CA}^\vv_1$, we recover the Lax matrix
$L^{\vv,0}_1(z)$ (for $b=-1$) and $L^{\vv,1}_1(z)$ (for $b=-2$). Moreover, this also explains
why we had exactly three Lax matrices in Section~\ref{relativistic}.
\end{Rem}

Fix $n\geq 1$ and consider the \emph{complete monodromy matrix}
$T^{\vv,-1}_n(z)=L^{\vv,-1}_n(z)\cdots L^{\vv,-1}_1(z)$.
Applying iteratively $\Delta^\rtt_{\bullet,\bullet}$ of
Theorem~\ref{Ding-Frenkel}(c), we get
  $\Delta^\rtt_{n}\colon \CU^\rtt_{0,-2n}\to (\CU^\rtt_{0,-2})^{\otimes n}$.
Composing it with the homomorphism
  $(\Phi^\rtt_{0,-2})^{\otimes n}\colon  (\CU^\rtt_{0,-2})^{\otimes n}\to (\hat{\CA}^\vv_1)^{\otimes n}\simeq  \hat{\CA}^\vv_n$,
we obtain the homomorphism
  $\Phi^\rtt_{0,-2n}\colon {} \CU^\rtt_{0,-2n}\to \hat{\CA}^\vv_n$.
The following in straightforward.

\begin{Lem}
We have $\Phi^\rtt_{0,-2n}(T^\pm(z))=T^{\vv,-1}_n(z)$.
\end{Lem}

\begin{Rem}
For $n>1$, the composition $\Phi^\rtt_{0,-2n}\circ \Upsilon_{0,-2n}$
{\em does not} coincide with the homomorphism $\wt{\Phi}_{-2n}^0$ of
Theorem~\ref{Homomorphism}.
\end{Rem}

\begin{Rem}
\label{classical Toda}
The result of~Lemma~\ref{Lax via homomorphism 1} admits a natural
{\em rational} counterpart. Let $\bY_{-2}$ be the shifted Yangian
of $\ssl_2$ with the shift $-\alpha$. Recall the homomorphism
  $\Phi^0_{-2}\colon \bY_{-2}\to \hat{\CA}^\hbar_1$
of~\cite[Corollary~B.17]{bfn}. Consider a slight modification of it
\begin{equation*}
    \hat{\Phi}_{-2}\colon
    E(z)\mapsto (z-w)^{-1}\sfu^{-1},\
    F(z)\mapsto -(z-w-\hbar)^{-1}\sfu,\
    H(z)\mapsto (z-w)^{-1}(z-w-\hbar)^{-1}.
\end{equation*}
One can also define a \emph{(rational) shifted RTT algebra of
$\ssl_2$}, denoted by $\CY^\rtt_{-2}$. This is an associative
$\BC[\hbar]$-algebra generated by
  $\{t_{11}[r-1],t_{12}[r],t_{21}[r],t_{22}[r+1], (t_{11}[-1])^{-1}\}_{r\geq0}$
and with the defining relations
  $(t_{11}[-1])^{\pm 1} (t_{11}[-1])^{\mp 1}=1,\
   T_{11}(z)T_{22}(z-\hbar)-T_{12}(z)T_{21}(z-\hbar)=1,\
   R_\rat(z-w)(T(z)\otimes 1)(1\otimes T(w))=
   (1\otimes T(w))(T(z)\otimes 1)R_\rat(z-w)$,
where
  $T(z)=(T_{ij}(z))_{i,j=1}^2$ with $T_{ij}(z):=\sum_r t_{ij}[r]z^{-r}$.
Consider the Gauss decomposition of $T(z)$:
\begin{equation*}
  T(z)=
  \left(%
  \begin{array}{cc}
    1 & 0 \\
    \tilde{f}(z) & 1 \\
  \end{array}%
  \right)
  \left(%
  \begin{array}{cc}
    \tilde{g}_1(z) & 0 \\
    0 & \tilde{g}_2(z) \\
  \end{array}%
  \right)
  \left(%
  \begin{array}{cc}
    1 & \tilde{e}(z) \\
    0 & 1 \\
  \end{array}%
  \right).
\end{equation*}
Analogously to Theorem~\ref{Ding-Frenkel}(b), there is a $\BC[\hbar]$-algebra
homomorphism
  $\Upsilon^\rat_{-2}\colon  \bY_{-2}\to \CY^\rtt_{-2}$,
defined by
  $E(z)\mapsto \tilde{e}(z), F(z)\mapsto \tilde{f}(z), H(z)\mapsto \tilde{g}_2(z)\tilde{g}_1(z)^{-1}$.
Composing $\Upsilon^\rat_{-2}$ with the homomorphism
  $\CY^\rtt_{-2}\to \hat{\CA}^\hbar_1$
given by $T(z)\mapsto L^\hbar_1(z)$, we recover $\hat{\Phi}_{-2}$
from above.
\end{Rem}


\subsection{Homomorphism $\Delta_{b_1,b_2}\ (b_1,b_2\in \BZ_{\leq 0})$ via Drinfeld half-currents, $\fg=\ssl_2$}
$\ $

Recall the currents $e^\pm(z), f^\pm(z), \psi^\pm (z)$ of~(\ref{Drinfeld currents}).

\begin{Prop}\label{nonshifted coproduct via currents}
Let $\Delta$ be the Drinfeld-Jimbo coproduct on $U_\vv(L\ssl_2)$. Then, we have
\begin{equation}\label{coproduct 1}
  \Delta(e^\pm(z))= 1\otimes e^\pm(z)+
  \sum_{r=0}^\infty (-\vv)^r(\vv-\vv^{-1})^{2r}\cdot e^\pm(z)^{r+1}\otimes f^\pm (\vv^2z)^r\psi^\pm(z),
\end{equation}
\begin{equation}\label{coproduct 2}
  \Delta(f^\pm(z))= f^\pm(z)\otimes 1+
  \sum_{r=0}^\infty (-\vv)^{-r}(\vv-\vv^{-1})^{2r}\cdot \psi^\pm(z)e^\pm(\vv^2z)^r\otimes f^\pm (z)^{r+1},
\end{equation}
\begin{equation}\label{coproduct 3}
  \Delta(\psi^\pm(z))=\sum_{r=0}^\infty
  (-1)^r[r+1]_\vv(\vv-\vv^{-1})^{2r}\cdot \psi^\pm(z)e^\pm(\vv^2z)^r\otimes f^\pm(\vv^2z)^r\psi^\pm(z).
\end{equation}
\end{Prop}

These formulas are analogous to those for the Yangian
$Y_\hbar(\ssl_2)$ of~\cite[Exercise 3.2]{m}. The proof of this
result is based on the RTT realization of $U_\vv(L\ssl_2)$ and
is presented in Appendix~\ref{Proof of coproduct via Drinfeld currents}.

\begin{Prop}\label{shifted coproduct via currents}
Let $b_1, b_2\in \BZ_{\leq 0}$ and $b=b_1+b_2$. Then, the homomorphism
  $\Delta_{b_1,b_2}\colon  \CU^\ssc_{0,b}\to \CU^\ssc_{0,b_1}\otimes \CU^\ssc_{0,b_2}$
from Theorem~\ref{Coproduct for sl2}
also satisfies the formulas~(\ref{coproduct 1}--\ref{coproduct
3}), where by abuse of notation $e^\pm (z), f^\pm (z), \psi^\pm(z)$
denote the generating series for each respective algebra.
\end{Prop}

\begin{proof}
Our proof is based on the commutative diagram of Remark~\ref{Maps jmath}:

\setlength{\unitlength}{1cm}
\begin{picture}(4,3.2)
 \put(3.3,2.5){$U^\pm_\vv$}
 \put(2.9,0.5){$\CU^{\ssc,\pm}_{0,b_1,b_2}$}
 \put(9.5,2.5){$U^\pm_\vv\otimes U^\pm_\vv$}
 \put(9.1,0.5){$\CU^{\ssc,\pm}_{0,b_1,0}\otimes \CU^{\ssc,\pm}_{0,0,b_2}$}

 \put(4.0,2.6){\vector (1,0){5.2}}
 \put(4.2,0.6){\vector (1,0){4.8}}
 \put(3.5,2.3){\vector (0,-1){1.3}}
 \put(10.3,2.3){\vector (0,-1){1.4}}

 \put(6.4,2.75){$\Delta$}
 \put(6.1,0.25){$\Delta_{b_1,b_2}$}
 \put(2.5,1.5){$\jmath^\pm_{b_1,b_2}$}
 \put(10.4,1.5){$\jmath^\pm_{b_1,0}\otimes \jmath^\pm_{0,b_2}$}
\end{picture}

Since $\jmath^+_{\bullet,\bullet}\colon  e^+(z)\mapsto e^+(z),
f^+(z)\mapsto f^+(z), \psi^+(z)\mapsto \psi^+(z)$, we immediately
get the validity of~(\ref{coproduct 1}--\ref{coproduct 3}) for the
currents $e^+(z),f^+(z),\psi^+(z)$ and the homomorphism $\Delta_{b_1,b_2}$.

Let us now treat the case of $e^-(z), f^-(z), \psi^-(z)$.
Combining the commutativity of the above diagram (in the ``$-$'' case)
with equality~(\ref{coproduct 1}) yields
\begin{equation*}
    \Delta_{b_1,b_2,}(\unl{e}^-(z))=1\otimes \unl{e}^-(z)+
    \sum_{r=0}^\infty (-\vv)^r(\vv-\vv^{-1})^{2r}\cdot e^-(z)^{r+1}\otimes f^- (\vv^2z)^r\psi^-(z),
\end{equation*}
where $\unl{e}^-(z):=e^-(z)+\sum_{r=b_2}^{-1}e_rz^{-r}$. Meanwhile,
$\Delta_{b_1,b_2}(e_r)=1\otimes e_r$ for $b_2\leq r\leq -1$. Hence,
$\Delta_{b_1,b_2}(e^-(z))$ is given by the right-hand side of~(\ref{coproduct 1}).
Likewise, we get the validity of~(\ref{coproduct 2}),~(\ref{coproduct 3}) for
the currents $f^-(z),\psi^-(z)$ and the homomorphism $\Delta_{b_1,b_2}$.
\end{proof}

Since our proof of~(\ref{coproduct 1}--\ref{coproduct 3}) in
Appendix~\ref{Proof of coproduct via Drinfeld currents} is based on
the RTT-type coproduct $\Delta^\rtt_{0,0}$, we immediately get

\begin{Cor}\label{compatibility of coproducts}
Let $b_1,b_2\in \BZ_{\leq 0}$ and $b=b_1+b_2$. The following diagram is commutative:
  $$\begin{CD}
    \CU^\ad_{0,2b} @>{\Delta^\ad_{2b_1,2b_2}}>> \CU^\ad_{0,2b_1}\otimes \CU^\ad_{0,2b_2}\\
    @V{\Upsilon_{0,2b}}V{\wr}V   @V{\wr}V{\Upsilon_{0,2b_1}\otimes \Upsilon_{0,2b_2}}V\\
    \CU^{\rtt}_{0,2b} @>>{\Delta^{\rtt}_{2b_1,2b_2}}> \CU^{\rtt}_{0,2b_1}\otimes \CU^{\rtt}_{0,2b_2}
    \end{CD}$$



\end{Cor}


\subsection{Coproduct for truncated shifted algebras, $\fg=\ssl_2$}
\label{descent coproduct}
\

For $b_1,b_2\in \BZ_{\leq 0}$ and $b=b_1+b_2$, recall the homomorphism
$\Delta^\ad_{2b_1,2b_2}\colon  \CU^\ad_{0,2b}\to \CU^\ad_{0,2b_1}\otimes \CU^\ad_{0,2b_2}$
of Remark~\ref{ad-coproduct on sl2}.
Consider the truncated versions of the algebras involved
$\CU^0_{2b},\CU^0_{2b_1},\CU^0_{2b_2}$, see~\refd{trunc}.
The goal of this subsection is to prove the following result.

\begin{Prop}\label{descent for sl2}
For $b_1,b_2\leq 0$, the homomorphism $\Delta^\ad_{2b_1,2b_2}$ descends to the same named
homomorphism $\CU^0_{2b}\to \CU^0_{2b_1}\otimes \CU^0_{2b_2}$.
\end{Prop}

\begin{proof}
Define a $2$-sided ideal $\CI\subset \CU^\ad_{0,2b_1}\otimes \CU^\ad_{0,2b_2}$ via
$\CI:=\CI^0_{2b_1}\otimes \CU^\ad_{0,2b_2}+\CU^\ad_{0,2b_1}\otimes \CI^0_{2b_2}$.
It suffices to show that $\Delta^\ad_{2b_1,2b_2}(X)\in \CI$ for every generator $X$
of the ideal $\CI^0_{2b}$ of~(\ref{kernel 1}--\ref{kernel 2}).
To achieve this, recall the commutative diagram of Corollary~\ref{compatibility of coproducts}.

\medskip
\noindent
$\circ$ \emph{Case $X=A^+_s (s>-b)$.}
Applying the aforementioned commutative diagram to the equality
  $\Delta^\rtt_{2b_1,2b_2}(t^+_{11}[s])=
   \sum_{s_1,s_2\geq 0}^{s_1+s_2=s} t^+_{11}[s_1]\otimes t^+_{11}[s_2]+
   \sum_{s_1,s_2\geq 0}^{s_1+s_2=s} t^+_{12}[s_1]\otimes t^+_{21}[s_2]$,
we get
  $\Delta^\ad_{2b_1,2b_2}(A^+_s)=
   \sum_{s_1,s_2\geq 0}^{s_1+s_2=s} A^+_{s_1}\otimes A^+_{s_2}+
   \sum_{s_1,s_2\geq 0}^{s_1+s_2=s} B^+_{s_1}\otimes C^+_{s_2}$.
For $s_1+s_2=s>-b$, either $s_1>-b_1$ or $s_2>-b_2$.
Hence, each summand in the right-hand side belongs to $\CI$, due to
Remark~\ref{first indication of conjecture}.

\medskip
\noindent
$\circ$ \emph{Case $X=A^+_0A^+_{-b}-(-1)^b$.}
As above $\Delta^\ad_{2b_1,2b_2}(A^+_{-b})\equiv A^+_{-b_1}\otimes A^+_{-b_2}$, where
the notation $x\equiv y$ is used to denote $x-y\in \CI$.
We also have $\Delta^\ad_{2b_1,2b_2}(A^+_0)=A^+_0\otimes A^+_0$. Thus
  $\Delta^\ad_{2b_1,2b_2}(A^+_0A^+_{-b}-(-1)^b)\equiv A^+_0A^+_{-b_1}\otimes A^+_0A^+_{-b_2}-(-1)^b=
   (A^+_0A^+_{-b_1}-(-1)^{b_1})\otimes A^+_0A^+_{-b_2}+(-1)^{b_1}\otimes (A^+_0A^+_{-b_2}-(-1)^{b_2})\equiv 0$.
Hence, $\Delta^\ad_{2b_1,2b_2}(A^+_0A^+_{-b}-(-1)^b)\in \CI$.

%

\medskip
\noindent
$\circ$ \emph{Case $X=A^-_{-r}-\vv^{-b}A^+_{-b-r}\ (0\leq r\leq -b)$.}
Analogously to the first case considered above, we have
  $\Delta^\ad_{2b_1,2b_2}(A^+_{-b-r})\equiv
   \sum_{\substack{0\leq r_1\leq -b_1\\ 0\leq r_2\leq -b_2}}^{r_1+r_2=r} A^+_{-b_1-r_1}\otimes A^+_{-b_2-r_2}+
   \sum_{\substack{1\leq r_1\leq -b_1\\ 0\leq r_2\leq -b_2-1}}^{r_1+r_2=r} B^+_{-b_1-r_1}\otimes C^+_{-b_2-r_2}$,
where the lower bounds on $r_1,r_2$ are due to Remark~\ref{first indication of conjecture}.
Completely analogously, we obtain
  $\Delta^\ad_{2b_1,2b_2}(A^-_{-r})
   \equiv
   \sum_{\substack{0\leq r_1\leq -b_1\\ 0\leq r_2\leq -b_2}}^{r_1+r_2=r} A^-_{-r_1}\otimes A^-_{-r_2}+
   \sum_{\substack{1\leq r_1\leq -b_1\\ 0\leq r_2\leq -b_2-1}}^{r_1+r_2=r} B^-_{-r_1}\otimes C^-_{-r_2}$.
Hence,
\begin{equation}\label{truncated key}
\begin{split}
  \Delta^\ad_{2b_1,2b_2}(A^-_{-r}-\vv^{-b}A^+_{-b-r})\equiv
  \sum_{\substack{0\leq r_1\leq -b_1\\ 0\leq r_2\leq -b_2}}^{r_1+r_2=r}(A^-_{-r_1}\otimes A^-_{-r_2}-\vv^{-b}A^+_{-b_1-r_1}\otimes A^+_{-b_2-r_2})+\\
  \sum_{\substack{1\leq r_1\leq -b_1\\ 0\leq r_2\leq -b_2-1}}^{r_1+r_2=r}(B^-_{-r_1}\otimes C^-_{-r_2}-\vv^{-b}B^+_{-b_1-r_1}\otimes C^+_{-b_2-r_2}).
\end{split}
\end{equation}
The first sum of~(\ref{truncated key}) belongs to $\CI$ as
  $A^-_{-r_1}\otimes A^-_{-r_2}-\vv^{-b}A^+_{-b_1-r_1}\otimes A^+_{-b_2-r_2}=
   (A^-_{-r_1}-\vv^{-b_1}A^+_{-b_1-r_1})\otimes A^-_{-r_2}+\vv^{-b_1}A^+_{-b_1-r_1}\otimes (A^-_{-r_2}-\vv^{-b_2}A^+_{-b_2-r_2})\in \CI$.
Completely analogously,
  $B^-_{-r_1}\otimes C^-_{-r_2}-\vv^{-b}B^+_{-b_1-r_1}\otimes C^+_{-b_2-r_2}=
   (B^-_{-r_1}-\vv^{-b_1}B^+_{-b_1-r_1})\otimes C^-_{-r_2}+\vv^{-b_1}B^+_{-b_1-r_1}\otimes (C^-_{-r_2}-\vv^{-b_2}C^+_{-b_2-r_2})$.
To complete the proof, it suffices to show
\begin{equation}\label{inclusion}
\begin{split}
  & B^-_{-r_1}-\vv^{-b_1}B^+_{-b_1-r_1}\in \CI^0_{2b_1}\ \mathrm{for}\ 1\leq r_1\leq -b_1,\\
  & C^-_{-r_2}-\vv^{-b_2}C^+_{-b_2-r_2}\in \CI^0_{2b_2}\ \mathrm{for}\ 0\leq r_2\leq -b_2-1.
\end{split}
\end{equation}
To prove the first inclusion of~(\ref{inclusion}), recall that
$B^+(z)=[e_0,A^+(z)]_{\vv^{-1}}$, due to Corollary~\ref{BC via A}.
Likewise (comparing the terms of degree $1$ in $w$ in the equality~(\ref{ABCD6})
with $\epsilon=-,\epsilon'=+$), we obtain $B^-(z)=[e_0,A^-(z)]_{\vv^{-1}}$. Therefore,
\begin{equation*}
  B^-_{-r_1}-\vv^{-b_1}B^+_{-b_1-r_1}=[e_0,A^-_{-r_1}-\vv^{-b_1}A^+_{-b_1-r_1}]_{\vv^{-1}}\in \CI^0_{2b_1}.
\end{equation*}
Similarly, applying the equalities $zC^\pm(z)=[A^\pm(z),f_1]_{\vv^{-1}}$, we obtain
\begin{equation*}
  C^-_{-r_2}-\vv^{-b_2}C^+_{-b_2-r_2}=[A^-_{-r_2-1}-\vv^{-b_2}A^+_{-b_2-r_2-1},f_1]_{\vv^{-1}}\in \CI^0_{2b_2},
\end{equation*}
which implies the second inclusion of~(\ref{inclusion}). Thus,
$\Delta^\ad_{2b_1,2b_2}(A^-_{-r}-\vv^{-b}A^+_{-b-r})\in \CI$.

\medskip
\noindent
The cases when $X$ is one of $A^-_{-s} (s>-b),\ A^-_0A^-_b-(-\vv^2)^{-b}$
are treated analogously to the above first two cases.
This completes our proof.
\end{proof}


\ssec{truncop}{Coproduct for truncated shifted algebras, general $\fg$}
\

Recall the homomorphism
  $\Delta_{\mu_1,\mu_2}\colon \CU^{\on{sc}}_{0,\mu}\to
   \CU^{\on{sc}}_{0,\mu_1}\otimes\CU^{\on{sc}}_{0,\mu_2}$
of~Theorem~\ref{Coproduct for sln in dominant case} ($\mu=\mu_1+\mu_2,\ \fg=\ssl_n$).
Given $N=N_1+N_2$, this coproduct extends to
\begin{equation*}
    \Delta^{\on{ad}}_{\mu_1,\mu_2}\colon
    \CU^{\on{ad}}_{0,\mu}[\sz_1^{\pm1},\ldots,\sz_N^{\pm1}]\longrightarrow
    \CU^{\on{ad}}_{0,\mu_1}[\sz_1^{\pm1},\ldots,\sz_{N_1}^{\pm1}]\otimes
    \CU^{\on{ad}}_{0,\mu_2}[\sz_{N_1+1}^{\pm1},\ldots,\sz_N^{\pm1}]
\end{equation*}
as in~Remark~\ref{comment}(c). Given two sequences
$\unl{\lambda}\!^{(1)}=(\omega_{i_1},\ldots,\omega_{i_{N_1}}),
\unl{\lambda}\!^{(2)}=(\omega_{i_{N_1+1}},\ldots,\omega_{i_N})$, we concatenate them to
$\unl{\lambda}=(\omega_{i_1},\ldots,\omega_{i_N})$ and consider the corresponding
truncated shifted algebras
$\CU^{\unl{\lambda}}_\mu,\CU^{\unl{\lambda}\!^{(1)}}_{\mu_1},\CU^{\unl{\lambda}\!^{(2)}}_{\mu_2}$
as in~\refd{trunc}.

\conj{truncopr}
The aforementioned homomorphism $\Delta^{\on{ad}}_{\mu_1,\mu_2}$
descends to the same named homomorphism $\Delta^{\on{ad}}_{\mu_1,\mu_2}\colon \CU^{\unl{\lambda}}_\mu
\to\CU^{\unl{\lambda}\!^{(1)}}_{\mu_1}\otimes\CU^{\unl{\lambda}\!^{(2)}}_{\mu_2}$.
\econj

We hope that the comultiplication $\Delta^{\on{ad}}_{\mu_1,\mu_2}$ can be defined
for arbitrary simply-laced $\fg$ (see~Section~\ref{open problems}) and descends
to the truncated shifted algebras.


\section{$K$-theory of parabolic Laumon spaces}
\label{Laumon}


\subsection{Parabolic Laumon spaces}
\label{Parabolic Laumon spaces}
$\ $

We recall the setup of~\cite{bffr}. Let $\bC$ be a smooth projective curve of
genus zero. We fix a coordinate $z$ on $\bC$, and consider the action of
$\BC^\times$ on $\bC$ such that $\vv(z)=\vv^{-2}z$.
We have $\bC^{\BC^\times}=\{0,\infty\}$.

We consider an $N$-dimensional vector space $W$ with a basis
$w_1,\ldots,w_N$. This defines a Cartan torus $T\subset
G=GL(N)=GL(W)$. We also consider its $2^N$-fold cover, the
bigger torus $\widetilde{T}$, acting on $W$ as follows: for
$\widetilde{T}\ni\unl{t}=(t_1,\ldots,t_N)$ we have
$\unl{t}(w_i)=t_i^2w_i$.

We fix an $n$-tuple of positive integers $\pi=(p_1,\ldots,p_n)\in \BZ_{>0}^n$
such that $p_1+\ldots+p_n=N$. Let $P\subset G$ be a parabolic subgroup
preserving the flag
  $0\subset W_1:=\langle w_1,\ldots,w_{p_1}\rangle\subset
   W_2:=\langle w_1,\ldots,w_{p_1+p_2}\rangle\subset \cdots\subset
   W_{n-1}:=\langle w_1,\ldots,w_{p_1+\ldots+p_{n-1}}\rangle\subset W_n:=W.$
Let $\CB:=G/P$ be the corresponding partial flag variety.

Given an $(n-1)$-tuple of nonnegative integers
$\unl{d}=(d_1,\ldots,d_{n-1})\in \BN^{n-1}$, we consider the Laumon parabolic
quasiflags' space $\CQ_{\unl{d}}$, see~\cite[4.2]{l}. It is the
moduli space of flags of locally free subsheaves
\begin{equation*}
  0\subset\CW_1\subset\cdots\subset\CW_{n-1}\subset\CW=W\otimes\CO_\bC
\end{equation*}
such that $\on{rank}(\CW_i)=p_1+\ldots+p_i$ and $\deg(\CW_i)=-d_i$.
It is known to be a smooth connected projective variety of dimension
$\dim\CB+\sum_{i=1}^{n-1} d_i(p_i+p_{i+1})$, see ~\cite[2.10]{l}.

We consider the following locally closed subvariety
$\fQ_{\unl{d}}\subset\CQ_{\unl{d}}$ (parabolic quasiflags based at
$\infty\in\bC$) formed by the flags
\begin{equation*}
  0\subset\CW_1\subset\cdots\subset\CW_{n-1}\subset\CW=W\otimes\CO_\bC
\end{equation*}
such that $\CW_i\subset\CW$ is a vector subbundle in a neighborhood
of $\infty\in\bC$, and the fiber of $\CW_i$ at $\infty$ equals the
span $\langle w_1,\ldots,w_{p_1+\ldots+p_i}\rangle\subset W$. It is
known to be a smooth connected quasiprojective variety of dimension
$\sum_{i=1}^{n-1} d_i(p_i+p_{i+1})$.


\subsection{Fixed points}\label{fixed points}
$\ $

The group $G\times\BC^\times$ acts naturally on $\CQ_{\unl{d}}$, and the
group $\widetilde{T}\times\BC^\times$ acts naturally on $\fQ_{\unl{d}}$.
The set of fixed points of $\widetilde{T}\times\BC^\times$ on
$\fQ_{\unl{d}}$ is finite; its description is given in~\cite[4.4]{bffr}.

Let $\vec{\unl{d}}$ be a collection of nonnegative integral vectors
  $\vec{d}_{ij}=(d^{(1)}_{ij},\ldots, d^{(p_j)}_{ij}),\ n-1\geq i\geq j\geq 1$,
such that
  $d_i=\sum_{j=1}^i |d_{ij}|=\sum_{j=1}^i \sum_{a=1}^{p_j}d^{(a)}_{ij}$,
and for $i\geq k\geq j$ we have $\vec{d}_{kj}\geq \vec{d}_{ij}$,
i.e., $d^{(a)}_{kj}\geq d^{(a)}_{ij}$ for any $1\leq a\leq p_j$.
Abusing notation, we denote by $\vec{\unl{d}}$ the corresponding
$\widetilde{T}\times\BC^\times$-fixed point in $\fQ_{\unl{d}}$:

\noindent
 $\CW_1=\CO_\bC(-d^{(1)}_{11}\cdot0)w_1\oplus \cdots \oplus \CO_\bC(-d^{(p_1)}_{11}\cdot0)w_{p_1},$

\noindent
 $\CW_2=\CO_\bC(-d^{(1)}_{21}\cdot0)w_1\oplus \cdots \oplus \CO_\bC(-d^{(p_1)}_{21}\cdot0)w_{p_1}
  \oplus\CO_\bC(-d^{(1)}_{22}\cdot0)w_{p_1+1}\oplus \cdots \oplus \CO_\bC(-d^{(p_2)}_{22}\cdot0)w_{p_1+p_2},$

 $\vdots$

\noindent
 $\CW_{n-1}=\CO_\bC(-d^{(1)}_{n-1,1}\cdot0)w_1\oplus \cdots \oplus \CO_\bC(-d^{(p_1)}_{n-1,1}\cdot0)w_{p_1}
  \oplus\cdots\\ \cdots \oplus \CO_\bC(-d^{(1)}_{n-1,n-1}\cdot0)w_{p_1+\ldots+p_{n-2}+1}
  \oplus\cdots\oplus\CO_\bC(-d^{(p_{n-1})}_{n-1,n-1}\cdot0)w_{p_1+\ldots+p_{n-1}}.$

\medskip
\noindent
 $\unl{N\!otation}$: Given a collection $\vec{\unl{d}}$ as above, we will
denote by $\vec{\unl{d}}\pm \delta^{(p)}_{ij}$ the collection
$\vec{\unl{d}}{}'$, such that
 ${d{}}'^{(p)}_{ij}=d^{(p)}_{ij}\pm 1$,
while
 ${d{}}'^{(a)}_{kl}=d^{(a)}_{kl}$ for $(a,k,l)\ne (p,i,j)$.


\subsection{Correspondences}\label{classic}
$\ $

For $i\in\{1,\ldots,n-1\}$ and $\unl{d}=(d_1,\ldots,d_{n-1})$, we
set $\unl{d}+i:=(d_1,\ldots,d_i+1,\ldots,d_{n-1})$. We have a correspondence
  $\CE_{\unl{d},i}\subset\CQ_{\unl{d}}\times \CQ_{\unl{d}+i}$
formed by the pairs $(\CW_\bullet,\CW'_\bullet)$ such that
$\CW'_i\subset\CW_i$ and we have $\CW_j=\CW'_j$ for $j\ne i$, see
~\cite[4.5]{bffr}. In other words, $\CE_{\unl{d},i}$ is the moduli
space of flags of locally free sheaves
\begin{equation*}
    0\subset\CW_1\subset\cdots\subset\CW_{i-1}\subset\CW'_i\subset\CW_i\subset
    \CW_{i+1}\subset\cdots\subset\CW_{n-1}\subset\CW
\end{equation*}
such that $\on{rank}(\CW_j)=p_1+\ldots+p_j$ and $\deg(\CW_j)=-d_j$,
while $\on{rank}(\CW'_i)=p_1+\ldots+p_i$ and $\deg(\CW'_i)=-d_i-1$.
According to~\cite[2.10]{l}, $\CE_{\unl{d},i}$ is a smooth
projective algebraic variety of dimension
  $\dim \CB+\sum_{i=1}^{n-1} d_i(p_i+p_{i+1})+p_i$.

We denote by $\bp$ (resp.\ $\bq$) the natural projection
$\CE_{\unl{d},i}\to\CQ_{\unl{d}}$ (resp.
$\CE_{\unl{d},i}\to\CQ_{\unl{d}+i}$). We also have a map
$\bs\colon\CE_{\unl{d},i}\to\bC,$
\begin{equation*}
    (0\subset\CW_1\subset\cdots\subset\CW_{i-1}\subset\CW'_i\subset\CW_i\subset
    \CW_{i+1}\subset\cdots\subset\CW_{n-1}\subset\CW)\mapsto\on{supp}(\CW_i/\CW'_i).
\end{equation*}
The correspondence $\CE_{\unl{d},i}$ comes equipped with a natural
line bundle $\CL_i$ whose fiber at a point
\begin{equation*}
    (0\subset\CW_1\subset\cdots\subset\CW_{i-1}\subset\CW'_i\subset\CW_i\subset
    \CW_{i+1}\subset\cdots\subset\CW_{n-1}\subset\CW)
\end{equation*}
equals $\Gamma(\bC,\CW_i/\CW'_i)$. Finally, we have a transposed correspondence
  $^\sT\CE_{\unl{d},i}\subset \CQ_{\unl{d}+i}\times\CQ_{\unl{d}}$.

Restricting to $\fQ_{\unl{d}}\subset\CQ_{\unl{d}}$, we obtain the correspondence
$\sE_{\unl{d},i}\subset\fQ_{\unl{d}}\times\fQ_{\unl{d}+i}$ together
with the line bundle $\sL_i$ and the natural maps
  $\bp\colon\sE_{\unl{d},i}\to\fQ_{\unl{d}},\
   \bq\colon\sE_{\unl{d},i}\to\fQ_{\unl{d}+i},\
   \bs\colon\sE_{\unl{d},i}\to\bC\backslash\{\infty\}$.
We also have a transposed correspondence
  $^\sT\sE_{\unl{d},i}\subset \fQ_{\unl{d}+i}\times\fQ_{\unl{d}}$.
It is a smooth quasiprojective variety of dimension
$\sum_{i=1}^{n-1} d_i(p_i+p_{i+1})+p_i$.


\subsection{Equivariant $K$-groups}
\label{Equivariant K-groups}
$\ $

We denote by ${}'M(\pi)$ the direct sum of equivariant (complexified) $K$-groups:
\begin{equation*}
  {}'M(\pi)=\bigoplus_{\unl{d}}K^{\widetilde{T}\times\BC^\times}(\fQ_{\unl{d}}).
\end{equation*}
It is a module over
  $K_{\widetilde{T}\times\BC^\times}(\on{pt})=\BC[\widetilde{T}\times\BC^\times]=
   \BC[t^{\pm 1}_1,\ldots,t^{\pm 1}_N,\vv^{\pm 1}]$.
We define
\begin{equation*}
    M(\pi):=\ {}'M(\pi)\otimes_{K_{\widetilde{T}\times\BC^\times}(\on{pt})}
    \on{Frac}(K_{\widetilde{T}\times\BC^\times}(\on{pt})).
\end{equation*}
It is naturally graded
\begin{equation*}
    M(\pi)=\oplus_{\unl{d}}M(\pi)_{\unl{d}},\ \mathrm{where}\
    M(\pi)_{\unl{d}}=K^{\widetilde{T}\times\BC^\times}(\fQ_{\unl{d}})
    \otimes_{K_{\widetilde{T}\times\BC^\times}(\on{pt})}\on{Frac}(K_{\widetilde{T}\times\BC^\times}(\on{pt})).
\end{equation*}

According to the Thomason localization theorem, restriction to the
$\widetilde{T}\times\BC^\times$-fixed point set induces an isomorphism
\begin{equation*}
    K^{\widetilde{T}\times\BC^\times}(\fQ_{\unl{d}})
    \otimes_{K_{\widetilde{T}\times\BC^\times}(\on{pt})} \on{Frac}(K_{\widetilde{T}\times\BC^\times}(\on{pt}))\iso
    K^{\widetilde{T}\times\BC^\times}(\fQ_{\unl{d}}^{\widetilde{T}\times\BC^\times})
    \otimes_{K_{\widetilde{T}\times\BC^\times}(\on{pt})}\on{Frac}(K_{\widetilde{T}\times\BC^\times}(\on{pt})).
\end{equation*}

The classes of the structure sheaves $[\vec{\unl{d}}]$ of the
$\widetilde{T}\times\BC^\times$-fixed points $\vec{\unl{d}}$
(see~Section~\ref{fixed points}) form a basis in
  $\bigoplus_{\unl{d}}K^{\widetilde{T}\times\BC^\times}(\fQ_{\unl{d}}^{\widetilde{T}\times\BC^\times})
   \otimes_{K_{\widetilde{T}\times\BC^\times}(\on{pt})}\on{Frac} (K_{\widetilde{T}\times\BC^\times}(\on{pt})).$
The embedding of a point $\vec{\unl{d}}$ into $\fQ_{\unl{d}}$ is a proper morphism,
so the direct image in the equivariant $K$-theory is well-defined,
and we will denote by $[\vec{\unl{d}}]\in M(\pi)_{\unl{d}}$ the direct image
of the structure sheaf of the point $\vec{\unl{d}}$.
The set $\{[\vec{\unl{d}}]\}$ forms a basis of $M(\pi)$.


\subsection{Action of $\CU^\vv_\pi$ on $M(\pi)$}
\label{Section on action}
$\ $

From now on, we will denote by $\CU^\vv_\pi$ the shifted quantum
affine algebra $\CU^\ssc_{0,\mu}$ for $\fg=\ssl_n$ and
$\mu=\sum_{j=1}^{n-1} (p_{j+1}-p_j)\omega_j$. We will also need the characters
$T_i$ of $\widetilde{T}\times\BC^\times$ defined via $T_i:=\prod_{j=p_1+\ldots+p_{i-1}+1}^{p_1+\ldots+p_i} t_j$.
Let $\vv$ stand for the character of $\widetilde{T}\times\BC^\times:\ (\unl{t},\vv)\mapsto \vv$.

For any $0\leq i\leq n$, we will denote by $\unl{\CW}{}_i$ the
tautological $(p_1+\ldots+p_i)$-dimensional vector bundle on
$\fQ_{\unl{d}}\times\bC$. Let
  $\varpi\colon \fQ_{\unl{d}}\times(\bC\backslash\{\infty\})\rightarrow \fQ_{\unl{d}}$
denote the standard projection. We define the generating series $\bb_i(z)$
with coefficients in the equivariant $K$-theory of $\fQ_{\unl{d}}$ as follows:
\begin{equation*}
    \bb_i(z):=
    \Lambda^{\bullet}_{-1/z}(\varpi_*({\unl{\CW}{}_i}|_{{\bC\backslash\{\infty\}}}))=
    1+\sum_{r\geq 1}{\Lambda^r(\varpi_*({\unl{\CW}{}_i}|_{{\bC\backslash\{\infty\}}}))(-z^{-1})^r}.
\end{equation*}
We also define the operators
\begin{equation}\label{ddvas}
  e_{i,r}:=T_{i+1}^{-1}\vv^{d_{i+1}-d_i+2-i}\bp_*((\vv^i\sL_i)^{\otimes r}\otimes \bq^*)
  \colon M(\pi)_{\unl{d}}\to M(\pi)_{\unl{d}-i},
\end{equation}
\begin{equation}\label{ttris}
  f_{i,r}:=T_i^{-1}\vv^{d_i-d_{i-1}+i}\bq_*((-\sL_i)^{\otimes p_i}\otimes
  (\vv^i\sL_i)^{\otimes r}\otimes \bp^*)
  \colon M(\pi)_{\unl{d}}\to M(\pi)_{\unl{d}+i},
\end{equation}
and consider the following generating series of operators on $M(\pi)$:
\begin{equation}\label{dvas}
  e_i(z)=\sum_{r=-\infty}^\infty e_{i,r}z^{-r}
  \colon  M(\pi)_{\unl{d}}\to M(\pi)_{\unl{d}-i}[[z,z^{-1}]],
\end{equation}
\begin{equation}\label{tris}
  f_i(z)=\sum_{r=-\infty}^\infty f_{i,r}z^{-r}
  \colon M(\pi)_{\unl{d}}\to M(\pi)_{\unl{d}+i}[[z,z^{-1}]].
\end{equation}
We define
  $\psi^+_i(z)\colon M(\pi)_{\unl{d}}\to M(\pi)_{\unl{d}}[[z^{-1}]]$ and
  $\psi^-_i(z)\colon M(\pi)_{\unl{d}}\to z^{p_i-p_{i+1}}M(\pi)_{\unl{d}}[[z]]$
via
\begin{equation}\label{raz1}
  \psi^\pm_i(z):=T_{i+1}^{-1}T_i\vv^{d_{i+1}-2d_i+d_{i-1}}\cdot
  \left(\frac{\bb_{i+1}(z\vv^{-i-2})\bb_{i-1}(z\vv^{-i})}{\bb_i(z\vv^{-i-2})\bb_i(z\vv^{-i})}\right)^{\pm},
\end{equation}
where as before $\gamma(z)^\pm$ denotes the expansion of a rational function $\gamma(z)$ in $z^{\mp 1}$, respectively.

\noindent
$\unl{N\!otation}$:
To each $\vec{\unl{d}}$, we assign a collection of
$\widetilde{T}\times\BC^\times$-weights
$s^{(a)}_{ij}:=t_{p_1+\ldots+p_{j-1}+a}^2\vv^{-2d^{(a)}_{ij}}$.

\begin{Prop}\label{matrix_elements}
(a) The matrix coefficients of the operators $f_{i,r}, e_{i,r}$ in
the fixed point basis $\{[\vec{\unl{d}}]\}$ of $M(\pi)$ are as follows:
\begin{equation*}
    f_{{i,r}[\vec{\unl{d}},\vec{\unl{d}}{}']}=
    T_i^{-1}\vv^{d_i-d_{i-1}+i}(1-\vv^2)^{-1}(-s^{(a)}_{ij})^{p_i}(s^{(a)}_{ij}\vv^i)^r
    \frac{\prod_{j'\leq i-1}^{a'\leq p_{j'}}(1-s^{(a)}_{ij}/s^{(a')}_{i-1,j'})}
         {\prod_{j'\leq i, a'\leq p_{j'}}^{(j',a')\ne (j,a)}(1-s^{(a)}_{ij}/s^{(a')}_{ij'})}
\end{equation*}
if $\vec{\unl{d}}{}'=\vec{\unl{d}}+\delta^{(a)}_{ij}$ for certain
$j\leq i,1\leq a\leq p_j$;
\begin{equation*}
    e_{{i,r}[\vec{\unl{d}},\vec{\unl{d}}{}']}=
    T_{i+1}^{-1}\vv^{d_{i+1}-d_i+2-i}(1-\vv^2)^{-1}(s^{(a)}_{ij}\vv^{i+2})^r
    \frac{\prod_{j'\leq i+1}^{a'\leq p_{j'}}(1-s^{(a')}_{i+1,j'}/s^{(a)}_{ij})}
         {\prod_{j'\leq i, a'\leq p_{j'}}^{(j',a')\ne (j,a)}(1-s^{(a')}_{ij'}/s^{(a)}_{ij})}
\end{equation*}
if $\vec{\unl{d}}{}'=\vec{\unl{d}}-\delta^{(a)}_{ij}$ for certain
$j\leq i,1\leq a\leq p_j$.

All the other matrix coefficients of $e_{i,r},f_{i,r}$ vanish.

\noindent
(b) The eigenvalue $\psi^\pm_i(z)|_{{\vec{\unl{d}}}}$ of $\psi^\pm_i(z)$ on $[\vec{\unl{d}}]$ equals
\begin{equation*}
    T_{i+1}^{-1}T_i\vv^{d_{i+1}-2d_{i}+d_{i-1}}
    \left(\frac{\prod_{j\leq i+1}^{a\leq p_j}(1-z^{-1}\vv^{i+2}s^{(a)}_{i+1,j})
                \prod_{j\leq i-1}^{a\leq p_j}(1-z^{-1}\vv^is^{(a)}_{i-1,j})}
               {\prod_{j\leq i}^{a\leq p_j}(1-z^{-1}\vv^{i+2}s^{(a)}_{ij})
                \prod_{j\leq i}^{a\leq p_j}(1-z^{-1}\vv^is^{(a)}_{ij})}\right)^\pm.
\end{equation*}
\end{Prop}

The proof is straightforward and is analogous to that of~\cite[Proposition 2.15]{t}.

\medskip
The following is the key result of this section.

\begin{Thm}\label{Action on K-theory of Laumon}
The generating series of operators $\{\psi^\pm_i(z), e_i(z),f_i(z)\}_{i=1}^{n-1}$
of~(\ref{dvas}--\ref{raz1}) acting on $M(\pi)$ satisfy the relations in $\CU^\vv_\pi$,
i.e., they give rise to the action of $\CU^\vv_\pi$ on $M(\pi)$.
\end{Thm}

In the particular case $\pi=1^n$, we recover~\cite[Theorem 2.12]{t}.

\begin{proof}
First, note that $\psi^+_i(z)$ contains only nonpositive powers
of $z$, while $\psi^-_i(z)$ contains only powers of $z$ bigger or
equal to $p_i-p_{i+1}$ (this follows from
Proposition~\ref{matrix_elements}(b)). Moreover, the coefficients of
$z^0$ in $\psi^+_i(z)$ and of $z^{p_i-p_{i+1}}$ in $\psi^-_i(z)$ are invertible operators.

Applying Proposition~\ref{matrix_elements}, the verification of all
the defining relations of $\CU^\vv_\pi$, except for (U6), boils down
to routine straightforward computations in the fixed point basis
(compare to the proof of~\cite[Theorem 2.12]{t}). The same arguments can
be used to show that $[e_i(z),f_j(w)]=0$ for $i\ne j$. It remains to
prove $(\vv-\vv^{-1})[e_i(z),f_i(w)]=\delta\left(\frac{z}{w}\right)\left(\psi^+_i(z)-\psi^-_i(z)\right)$.
Applying Proposition~\ref{matrix_elements}(a), we see that the
left-hand side is diagonal in the fixed point basis and its
eigenvalue on $[\vec{\unl{d}}]$ equals
  $$T^{-1}_{i+1}T^{-1}_i\vv^{d_{i+1}-d_{i-1}}(1-\vv^2)^{-1}\cdot\delta\left(\frac{z}{w}\right)\times$$
  $$\sum_{j\leq i}^{a\leq p_j}(-s^{(a)}_{ij})^{p_i}
    \left\{
           \vv^{2p_i}\frac{\prod_{j'\leq i+1}^{a'\leq p_{j'}}(1-s^{(a')}_{i+1,j'}/s^{(a)}_{ij})
                   \prod_{j'\leq i-1}^{a'\leq p_{j'}}(1-\vv^2s^{(a)}_{ij}/s^{(a')}_{i-1,j'})}
                  {\prod_{j'\leq i,a'\leq p_{j'}}^{(j',a')\ne (j,a)}
                   (1-s^{(a')}_{ij'}/s^{(a)}_{ij})(1-\vv^2s^{(a)}_{ij}/s^{(a')}_{ij'})}
          \delta\left(\frac{z}{\vv^{i+2}s^{(a)}_{ij}}\right)-\right.$$
  $$\left.
             \frac{\prod_{j'\leq i+1}^{a'\leq p_{j'}}(1-\vv^2s^{(a')}_{i+1,j'}/s^{(a)}_{ij})
                   \prod_{j'\leq i-1}^{a'\leq p_{j'}}(1-s^{(a)}_{ij}/s^{(a')}_{i-1,j'})}
                  {\prod_{j'\leq i,a'\leq p_{j'}}^{(j',a')\ne (j,a)}
                   (1-\vv^2s^{(a')}_{ij'}/s^{(a)}_{ij})(1-s^{(a)}_{ij}/s^{(a')}_{ij'})}
             \delta\left(\frac{z}{\vv^is^{(a)}_{ij}}\right)\right\}.$$
To compare this expression with the eigenvalue of $\psi^+_i(z)-\psi^-_i(z)$
on $[\vec{\unl{d}}]$, it suffices
to apply Lemma~\ref{key equality} below to the particular case of
$\gamma(z)$ chosen to be the rational function of
Proposition~\ref{matrix_elements}(b).

 The theorem is proved.
\end{proof}

\begin{Rem}\label{Relation to bdghk}
(a) The above verification of (U6) by applying
Lemma~\ref{key equality} significantly simplifies our original
indirect proof of this relation in~\cite{t}.

\noindent
(b) For $\pi=p^n$, this produces the action of the quantum loop algebra
$U_\vv(L\ssl_n)$ on $M(\pi)$.

\noindent
(c) According to~\cite{bdghk}, there is an action of $\CA^\vv_\fra$
on $M(\pi)$. Its pull-back along the homomorphism
$\ol{\Phi}{}^{\unl\lambda}_\mu$ ($\unl\lambda=(\omega_{n-1},\ldots,\omega_{n-1})$
taken $N$ times) yields essentially the action of $\CU^\vv_\pi$
on $M(\pi)$ established above. In particular, the kernel
$\Ker(\ol{\Phi}{}^{\unl\lambda}_\mu)=\Ker(\wt{\Phi}^{\unl\lambda}_\mu)$ acts
trivially on $M(\pi)$. The first instance of that is the fact that the generators
$\{A^\pm_{i,\pm r}: r>p_1+\ldots+p_i\}$ of $\CU^\vv_\pi$
(see~Remark~\ref{Remark on ABCD}(b)) act trivially on $M(\pi)$,
due to the observation that the eigenvalue of $A^\pm_i(z)$ on
$[\vec{\unl{d}}]$ equals $\prod_{j\leq i}^{a\leq p_j}
(1-(z^{-1}\vv^is^{(a)}_{ij})^{\pm 1})$.
%
\end{Rem}


\subsection{Tensor products}
\label{Factorizing Laumons}
$\ $

Fix two $n$-tuples
  $\pi'=(p'_1,\ldots,p'_n), \pi''=(p''_1,\ldots,p''_n)\in \BZ_{>0}^n$
and define $\pi=(p_1,\ldots,p_n)$ via $p_i:=p'_i+p''_i\in \BZ_{>0}$.
Let $\CU^\vv_{\pi'}, \CU^\vv_{\pi''}, \CU^\vv_\pi$ be the
corresponding shifted quantum affine algebras of $\ssl_n$ as
defined in Section~\ref{Section on action}. According to
Theorem~\ref{Action on K-theory of Laumon}, we have natural actions
of $\CU^\vv_\pi$ on $M(\pi)$, of $\CU^\vv_{\pi'}$ on $M(\pi')$, and
of $\CU^\vv_{\pi''}$ on $M(\pi'')$. The vector spaces $M(\pi)$ and
$M(\pi')\otimes M(\pi'')$ have natural fixed point bases
$\{[\vec{\unl{d}}]\}$ and $\{[\vec{\unl{d}}{}']\otimes
[\vec{\unl{d}}{}'']\}$, parameterized by $\vec{\unl{d}}$ and pairs
$(\vec{\unl{d}}{}', \vec{\unl{d}}{}'')$ with
$\vec{\unl{d}},\vec{\unl{d}}{}',\vec{\unl{d}}{}''$ satisfying the
conditions of Section~\ref{fixed points}.
The assignment $(\vec{\unl{d}}{}', \vec{\unl{d}}{}'')\mapsto \vec{\unl{d}}{}'\cup\vec{\unl{d}}{}''$
defined via
  $(d'\cup d'')^{(a)}_{ij}={d{}}'^{(a)}_{ij},
   (d'\cup d'')^{(p'_j+b)}_{ij}={d{}}''^{(b)}_{ij}$
for $1\leq a\leq p'_j, 1\leq b\leq p''_j$ provides a bijection
between such pairs $(\vec{\unl{d}}{}',\vec{\unl{d}}{}'')$ and $\vec{\unl{d}}$.
We also identify $\wt{T}'\times \wt{T}''\iso \wt{T}$ via
  $t_{p_1+\ldots+p_{j-1}+a}=t'_{p'_1+\ldots+p'_{j-1}+a},
   t_{p_1+\ldots+p_{j-1}+p'_j+b}=t''_{p''_1+\ldots+p''_{j-1}+b}$
for $a,b$ as above. Finally, we use $\vec{\unl{0}}$ to denote the
collection of zero vectors.

Recall the Drinfeld formal coproduct
  $\wt{\Delta}\colon \CU^\vv_\pi\to \CU^\vv_{\pi'}\widehat{\otimes} \CU^\vv_{\pi''}$
of Lemma~\ref{new Drinfeld coproduct}.

\begin{Thm}\label{Tensor products}
There is a unique collection of
$c_{\vec{\unl{d}}{}',\vec{\unl{d}}{}''}\in\on{Frac}(K_{\wt{T}\times\BC^\times}(\on{pt}))$
with $c_{\vec{\unl{0}},\vec{\unl{0}}}=1$, such that the map
  $[\vec{\unl{d}}{}']\otimes [\vec{\unl{d}}{}'']\mapsto
   c_{\vec{\unl{d}}{}',\vec{\unl{d}}{}''}\cdot [\vec{\unl{d}}{}'\cup\vec{\unl{d}}{}'']$
induces an isomorphism $M(\pi')\wt{\otimes} M(\pi'')\iso M(\pi)$ of
$\CU^\vv_\pi$-representations.
\end{Thm}

First let us make sense of the $\CU^\vv_\pi$-module $M(\pi')\wt{\otimes} M(\pi'')$. The
action of $e_i(z)$ in the fixed point basis
$\{[\vec{\unl{d}}{}'']\}$ of $M(\pi'')$ can be written as
  $e_i(z)[\vec{\unl{d}}{}'']=
   \sum_{j\leq i}^{a\leq p_j} a_{\vec{\unl{d}}{}'',\delta^{(a)}_{ij}}
   \delta(s^{(a)}_{ij}\vv^{i+2}/z)[\vec{\unl{d}}{}''-\delta^{(a)}_{ij}]$
for certain
  $a_{\vec{\unl{d}}{}'',\delta^{(a)}_{ij}}\in \on{Frac}(K_{\wt{T}''\times\BC^\times}(\on{pt}))$.
According to the comultiplication formula~(\ref{Drinfeld coproduct}), we have
  $\wt{\Delta}(e_i(z))([\vec{\unl{d}}{}'] \otimes [\vec{\unl{d}}{}''])=
   e_i(z)([\vec{\unl{d}}{}'])\otimes [\vec{\unl{d}}{}'']+
   \psi^-_i(z)([\vec{\unl{d}}{}'])\otimes e_i(z)([\vec{\unl{d}}{}'']).$
The first summand is well-defined. To make sense of the second
summand, we just need to apply the formula
$\gamma(z)\delta(a/z)=\gamma(a)\delta(a/z)$ to the rational function $\gamma(z)$
chosen to be the eigenvalue of $\psi^-_i(z)$ on
$[\vec{\unl{d}}{}']$. The action of $f_i(z)$ on
$M(\pi')\wt{\otimes} M(\pi'')$ is defined analogously.
Finally, the formula
 $\wt{\Delta}(\psi^\pm_i(z))=\psi^\pm_i(z)\otimes \psi^\pm_i(z)$
provides a well-defined action of $\psi^\pm_i(z)$.
These formulas endow $M(\pi')\otimes M(\pi'')$ with a well-defined action of $\CU^\vv_\pi$.

\begin{proof}
According to Proposition~\ref{matrix_elements}(b), the eigenvalue of
$\wt\Delta(\psi^\pm_i(z))=\psi^\pm_i(z)\otimes\psi^\pm_i(z)$ on
  $[\vec{\unl{d}}{}'] \otimes [\vec{\unl{d}}{}'']\in M(\pi')\otimes M(\pi)''$
equals the eigenvalue of $\psi^\pm_i(z)$ on
$[\vec{\unl{d}}{}'\cup\vec{\unl{d}}{}'']\in M(\pi)$. Hence, the map
  $[\vec{\unl{d}}{}']\otimes [\vec{\unl{d}}{}'']\mapsto
   c_{\vec{\unl{d}}{}',\vec{\unl{d}}{}''}\cdot [\vec{\unl{d}}{}'\cup\vec{\unl{d}}{}'']$
intertwines actions of $\psi^\pm_i(z)$ for any
$c_{\vec{\unl{d}}{}',\vec{\unl{d}}{}''}\in\on{Frac}(K_{\wt{T}\times\BC^\times}(\on{pt}))$.

Consider
$c_{\vec{\unl{d}}{}',\vec{\unl{d}}{}''}\in\on{Frac}(K_{\wt{T}\times\BC^\times}(\on{pt}))$
such that $c_{\vec{\unl{0}},\vec{\unl{0}}}=1$ and
\begin{equation}\label{compatibility}
\begin{split}
  & \frac{c_{\vec{\unl{d}}{}'-\delta^{(a)}_{ij},\vec{\unl{d}}{}''}}{c_{\vec{\unl{d}}{}',\vec{\unl{d}}{}''}}=
  (T''_{i+1})^{-1}\vv^{d''_{i+1}-d''_i}\cdot
  \frac{\prod_{j'\leq i+1}^{a'\leq p''_{j'}}(1-{s{}}''^{(a')}_{i+1,j'}/{s{}}'^{(a)}_{ij})}
       {\prod_{j'\leq i}^{a'\leq p''_{j'}}(1-{s{}}''^{(a')}_{ij'}/{s{}}'^{(a)}_{ij})},\\
  & \frac{c_{\vec{\unl{d}}{}',\vec{\unl{d}}{}''-\delta^{(a)}_{ij}}}{c_{\vec{\unl{d}}{}',\vec{\unl{d}}{}''}}=
  (T'_i)^{-1}\vv^{d'_i-d'_{i-1}}\cdot
  \frac{\prod_{j'\leq i}^{a'\leq p'_{j'}}(1-\vv^{-2}{s{}}'^{(a')}_{ij'}/{s{}}''^{(a)}_{ij})}
       {\prod_{j'\leq i-1}^{a'\leq p'_{j'}}(1-\vv^{-2}{s{}}'^{(a')}_{i-1,j'}/{s{}}''^{(a)}_{ij})}.
\end{split}
\end{equation}
The existence of $c_{\vec{\unl{d}}{}',\vec{\unl{d}}{}''}$ satisfying
these relations as well as a verification that
  $[\vec{\unl{d}}{}']\otimes [\vec{\unl{d}}{}'']\mapsto
   c_{\vec{\unl{d}}{}',\vec{\unl{d}}{}''}\cdot [\vec{\unl{d}}{}'\cup \vec{\unl{d}}{}'']$
intertwines actions of $e_{i,r}$ and $f_{i,r}$ are left to the interested reader.
\end{proof}

\begin{Rem}
In the particular case $p_1=\ldots=p_n=p$, this implies the isomorphism
$M(p^n)\simeq M(1^n)^{\wt{\otimes} p}$ of $U_\vv(L\ssl_n)$-representations.
This isomorphism is reminiscent of the isomorphism between the action of
the quantum toroidal algebra of $\gl_1$ on the equivariant $K$-theory of
the Gieseker moduli spaces $M(r,n)$ and the $r$-fold tensor product of such
representation for $r=1$, see~\cite[Theorem 4.6]{t2}.
\end{Rem}


\subsection{Shifted quantum affine algebras of $\gl_n$}\label{shifted gl_n}
$\ $

Let $U_\vv(\widehat{\gl_n})$ be the quantum affine algebra of $\gl_n$
as defined in~\cite[Definition 3.1]{df}, and let $U_\vv(L\gl_n)$ be the quantum
loop algebra of $\gl_n$, that is,
  $U_\vv(L\gl_n):=U_\vv(\widehat{\gl_n})/(\vv^{\pm c/2}-1)$.
This is an associative $\BC(\vv)$-algebra generated by
\begin{equation*}
  \{X^\pm_{i,r}, k^\pm_{j,\mp s^\pm_j}|i=1,\ldots,n-1,\ j=1,\ldots,n,\ r\in \BZ,\ s^\pm_j \in \BN\}
\end{equation*}
and with the defining relations as in~\cite[(3.3, 3.4)]{df}.
There is a natural injective $\BC(\vv)$-algebra homomorphism
$U_\vv(L\ssl_n)\hookrightarrow U_\vv(L\gl_n)$, defined by
\begin{equation}\label{relating sln and gln}
  e_i(z)\mapsto \frac{X^-_i(\vv^iz)}{\vv-\vv^{-1}},\
  f_i(z)\mapsto \frac{X^+_i(\vv^iz)}{\vv-\vv^{-1}},\
  \psi^\pm_i(z)\mapsto (k^\mp_i(\vv^iz))^{-1}k^\mp_{i+1}(\vv^iz).
\end{equation}

For $\pi=(p_1,\ldots,p_n)\in \BZ^n_{>0}$, define the shifted quantum affine
algebra $\CU^\vv_\pi(\gl_n)$ in the same way as
$U_\vv(L\gl_n)$ except that now $s^+_j\geq-p_j$ and we
formally add inverse elements $\{(k^-_{j,0})^{-1}, (k^+_{j,p_j})^{-1}\}_{j=1}^n$
(as we no longer require $k^-_{j,0}k^+_{j,p_j}=1$).
Note that the assignment~(\ref{relating sln and gln}) still gives  rise to
an injective\footnote{One can prove the
injectivity of $\varrho$ by using Proposition~\ref{Triangular
decomposition} for both algebras. Indeed, the homomorphism $\varrho$
is `glued' from three homomorphisms:
  $\varrho^>\colon \CU^{\vv,>}_\pi\to \CU^{\vv,>}_\pi(\gl_n),
   \varrho^<\colon \CU^{\vv,<}_\pi\to \CU^{\vv,<}_\pi(\gl_n),
   \varrho^0\colon \CU^{\vv,0}_\pi\to \CU^{\vv,0}_\pi(\gl_n)$.
The homomorphisms $\varrho^>,\varrho^<$ are isomorphisms due to
Proposition~\ref{Triangular decomposition}(b), while the injectivity
of $\varrho^0$ is clear.}
homomorphism $\varrho\colon \CU^\vv_\pi\hookrightarrow\CU^\vv_\pi(\gl_n)$.

Consider the following generating series of operators on $M(\pi)$:
\begin{equation*}
  X^+_i(z):=(\vv-\vv^{-1})f_i(\vv^{-i}z)\colon M(\pi)_{\unl{d}}\to M(\pi)_{\unl{d}+i}[[z,z^{-1}]],
\end{equation*}
\begin{equation*}
  X^-_i(z):=(\vv-\vv^{-1})e_i(\vv^{-i}z)\colon M(\pi)_{\unl{d}}\to M(\pi)_{\unl{d}-i}[[z,z^{-1}]],
\end{equation*}
\begin{equation*}
  k^-_j(z):=T_j^{-1}\vv^{d_j-d_{j-1}}\cdot (\bb_j(z\vv^{-2j})/\bb_{j-1}(z\vv^{-2j}))^+
    \colon M(\pi)_{\unl{d}}\to M(\pi)_{\unl{d}}[[z^{-1}]],
\end{equation*}
\begin{equation*}
  k^+_j(z):=T_j^{-1}\vv^{d_j-d_{j-1}}\cdot (\bb_j(z\vv^{-2j})/\bb_{j-1}(z\vv^{-2j}))^-
    \colon M(\pi)_{\unl{d}}\to z^{-p_j}M(\pi)_{\unl{d}}[[z]]
\end{equation*}
with $e_i(z), f_i(z), \bb_j(z)$ defined in Section~\ref{Section on action}.

The following is a simple generalization of Theorem~\ref{Action on
K-theory of Laumon}.

\begin{Thm}\label{action of gl_n on K-theory}
The generating series of operators $X^\pm_i(z),k^\pm_j(z)$ acting
on $M(\pi)$ satisfy the relations of $\CU^\vv_\pi(\gl_n)$, i.e., they
give rise to the action of $\CU^\vv_\pi(\gl_n)$ on $M(\pi)$.
\end{Thm}

The restriction of this action to the subalgebra $\CU_\pi^\vv$
(embedded into $\CU^\vv_\pi(\gl_n)$ via $\varrho$)
recovers the action of $\CU^\vv_\pi$
on $M(\pi)$ of Theorem~\ref{Action on K-theory of Laumon}.


\subsection{The cohomology case revisited}
$\ $

The above results can be immediately
generalized to the cohomological setting. Let $V(\pi)$ be the direct
sum of localized $T\times \BC^\times$-equivariant cohomology of type
$\pi$ Laumon parabolic based quasiflags' spaces:
\begin{equation*}
    V(\pi):=
    \bigoplus_{\unl{d}}H^\bullet_{T\times\BC^\times}(\fQ_{\unl{d}})
    \otimes_{H^\bullet_{T\times\BC^\times}(\on{pt})} \on{Frac}(H^\bullet_{T\times\BC^\times}(\on{pt})).
\end{equation*}
It is a module over $\on{Frac}(H^\bullet_{T\times\BC^\times}(\on{pt}))$, where
  $H^\bullet_{T\times\BC^\times}(\on{pt})=\BC[\mathrm{Lie}(T\times\BC^\times)]=
   \BC[x_1,\ldots,x_N,\hbar]$.

Let $\CY^\hbar_\pi=\bY_\pi\otimes_{\BC[\hbar]}\BC(\hbar)$, where $\bY_\pi$ is the
shifted Yangian of $\ssl_n$ in the sense of~\cite[Appendix B(i)]{bfn}.
It is the associative $\BC(\hbar)$-algebra generated by
  $\{E^{(r+1)}_i,F^{(r+1)}_i,H^{(r+1+p_i-p_{i+1})}_i\}_{1\leq i<n}^{r\in \BN}$
with the same defining relations as in the standard
Yangian $Y_\hbar(\ssl_n)$.

We define the generating series
  $\ba_i(z)$
with coefficients in the equivariant cohomology of $\fQ_{\unl{d}}$ as follows:
  $$\ba_i(z):=z^{p_1+\ldots+p_i}\cdot c(\varpi_*({\unl{\CW}{}_i}|_{{\bC\backslash\{\infty\}}}),(-z\hbar)^{-1}),$$
where $c(\mathcal{V},x)$ denotes the Chern polynomial (in $x$) of $\mathcal{V}$.
We also define the operators
\begin{equation}\label{Yddvas}
  E^{(r+1)}_i:=\bp_*((c_1(\sL_i)+i\hbar/2)^r\cdot \bq^*)\colon V(\pi)_{\unl{d}}\to V(\pi)_{\unl{d}-i},
\end{equation}
\begin{equation}\label{Yttris}
  F^{(r+1)}_i:=(-1)^{p_i}\bq_*((c_1(\sL_i)+i\hbar/2)^r\cdot \bp^*)\colon V(\pi)_{\unl{d}}\to V(\pi)_{\unl{d}+i}.
\end{equation}
We define
  $H_i(z)=z^{p_{i+1}-p_i}+\sum_{r>p_i-p_{i+1}}H^{(r)}_i\hbar^{-r+p_i-p_{i+1}+1}z^{-r}$
via
\begin{equation}\label{Yraz}
  H_i(z):=\left(\frac{\ba_{i+1}(z-\frac{i+2}{2})\ba_{i-1}(z-\frac{i}{2})}
                    {\ba_i(z-\frac{i+2}{2})\ba_i(z-\frac{i}{2})}\right)^+
  \colon V(\pi)_{\unl{d}}\to z^{p_{i+1}-p_i}V(\pi)_{\unl{d}}[[z^{-1}]].
\end{equation}

The following result is completely analogous to Theorem~\ref{Action
on K-theory of Laumon}.

\begin{Thm}\label{Action on cohomology of Laumon}
The operators
  $\{E^{(r+1)}_i, F^{(r+1)}_i, H^{(r+1+p_i-p_{i+1})}_i\}_{1\leq i<n}^{r\in\BN}$
of~(\ref{Yddvas}--\ref{Yraz}) acting on $V(\pi)$ satisfy the defining
relations of $\CY^\hbar_\pi$, i.e., they give rise to the action of
$\CY^\hbar_\pi$ on $V(\pi)$.
\end{Thm}

A slight refinement of this theorem in the {\em dominant}
case $p_1\leq\ldots\leq p_n$ constituted the key result
of~\cite{bffr}. In \emph{loc.~cit.}, the authors constructed the
action of the shifted Yangian of $\gl_n$, denoted by
$\CY^\hbar_\pi(\gl_n)$, on $V(\pi)$. There is a natural (injective)
homomorphism $\CY^\hbar_\pi\to \CY^\hbar_\pi(\gl_n)$, such that
  $F^{(r+1)}_i\mapsto \sum_{s=0}^r \binom{r}{s}\left(\frac{2-i}{2}\hbar\right)^{r-s}\mathsf{f}^{(s+1)}_i,\
   E^{(r+1)}_i\mapsto \sum_{s=0}^r \binom{r}{s}\left(\frac{2-i}{2}\hbar\right)^{r-s}\mathsf{e}^{(s+1+p_{i+1}-p_i)}_i$.
The pull-back of the action of~\cite{bffr} along this homomorphism
recovers the action $\CY^\hbar_\pi$ on $V(\pi)$ of Theorem~\ref{Action on cohomology of
Laumon}.

The proof of~\cite{bffr} was based on an explicit identification of the
geometric action in the fixed point basis with the formulas
of~\cite{fmo} for the action of $\CY^\hbar_\pi(\gl_n)$ in the
Gelfand-Tsetlin basis.
The benefits of our straightforward proof of Theorem~\ref{Action on
cohomology of Laumon} are two-fold:

(1) we eliminate the crucial assumption $p_1\leq \ldots\leq p_n$ of~\cite{bffr},

(2) we obtain an alternative proof of the formulas of~\cite{fmo}
(cf. Proposition~\ref{qGT} below).

Moreover, we can derive $\vv$-analogues of the Gelfand-Tsetlin
formulas of~\cite{fmo} via a certain specialization of the parameters
in Proposition~\ref{matrix_elements} as explained below. We set
$t_l=\vv^{\beta_l}$ for $1\leq l\leq N$. To a collection
  $\vec{\unl{d}}=(d^{(a)}_{ij})_{1\leq j\leq i\leq n-1}^{1\leq a\leq p_j}$,
we associate a Gelfand-Tsetlin pattern
  $\Lambda=\Lambda(\vec{\unl{d}})=(\lambda^{(a)}_{ij})_{1\leq j\leq i\leq n}^{1\leq a\leq p_j}$
as follows:
  $\lambda^{(a)}_{nj}=\beta_{p_1+\ldots+p_{j-1}+a}+j-1,\
   \lambda^{(a)}_{ij}=\beta_{p_1+\ldots+p_{j-1}+a}+j-1-d^{(a)}_{ij}.$
Set $\lambda^{(a)}_j:=\lambda^{(a)}_{nj}$, which is independent of
$\vec{\unl{d}}$. Note that the vector space $M(\pi)$ has a basis
$\{[\Lambda]\}$ 
parametrized by $\Lambda=(\lambda^{(a)}_{ij})^{1\leq a\leq p_j}_{1\leq j\leq i\leq n}$
with $\lambda^{(a)}_{nj}=\lambda^{(a)}_j$ and
$\lambda^{(a)}_{i+1,j}-\lambda^{(a)}_{ij}\in \BN$.
Consider a specialization of $\{\beta_l\}_{1\leq l\leq N}$ such that
$\lambda^{(a)}_j-\lambda^{(a)}_{j+1}\in \BN$, while
$\lambda^{(a)}_i-\lambda^{(b)}_j\notin \BZ$ if $a\ne b$.
Let $S$ be the subset of those $\Lambda$ from above such that
$\lambda^{(a)}_{ij}-\lambda^{(a)}_{i+1,j+1}\in \BN$ (note that $S$
is finite), while $\bar{S}$ will denote the set of the remaining
Gelfand-Tsetlin patterns $\Lambda$.

As before, we define
\begin{equation*}
  A_i^\pm(z):=k_1^\mp(\vv^{2-i}z)k_2^\mp(\vv^{4-i}z)\cdots k_i^\mp(\vv^iz),
\end{equation*}
\begin{equation*}
  B_i^\pm(z):=(\vv-\vv^{-1})A_i^\pm(z)e_i^\pm(z),
\end{equation*}
\begin{equation*}
  C_i^\pm(z):=(\vv-\vv^{-1})f_i^\pm(z)A_i^\pm(z).
\end{equation*}
We set $\lambda_{ij}(z):=\prod_{a=1}^{p_j}(\vv^{-\lambda_{ij}^{(a)}}-\vv^{\lambda_{ij}^{(a)}}z^{-1})$.
The next result follows from Proposition~\ref{matrix_elements}.

\begin{Prop}\label{qGT}
(a) The vector subspace of $M(\pi)$ spanned by
$\{[\Lambda]\}_{\Lambda\in \bar{S}}$ is $\CU^\vv_\pi(\gl_n)$-invariant.
We denote by $L(\pi)$ the corresponding quotient of $M(\pi)$.

\noindent
(b) Let $\{\xi_\Lambda\}_{\Lambda\in S}$  be the basis of $L(\pi)$
inherited from $\{[\Lambda]\}_{\Lambda\in S}$. Then, we have:
\begin{equation*}
    A^\pm_i(\vv^iz)\xi_\Lambda=\vv^{m_i}\lambda_{i1}(z)\lambda_{i2}(\vv^2z)\cdots
    \lambda_{ii}(\vv^{2(i-1)}z)\xi_\Lambda,
\end{equation*}
\begin{equation*}
    B^\pm_i(\vv^i\cdot\vv^{2l^{(a)}_{ij}})\xi_\Lambda=-\vv^{m_{i+1}-i}\cdot
    \lambda_{i+1,1}(\vv^{2l^{(a)}_{ij}})\lambda_{i+1,2}(\vv^{2(l^{(a)}_{ij}+1)})\cdots
    \lambda_{i+1,i+1}(\vv^{2(l^{(a)}_{ij}+i)})\xi_{\Lambda+\delta^{(a)}_{ij}},
\end{equation*}
\begin{equation*}
    C^\pm_i(\vv^i\cdot\vv^{2l^{(a)}_{ij}})\xi_\Lambda=\vv^{m_{i-1}+i-1}\cdot
    \lambda_{i-1,1}(\vv^{2l^{(a)}_{ij}})\lambda_{i-1,2}(\vv^{2(l^{(a)}_{ij}+1)})\cdots
    \lambda_{i-1,i-1}(\vv^{2(l^{(a)}_{ij}+i-2)})\xi_{\Lambda-\delta^{(a)}_{ij}},
\end{equation*}
where $m_j:=\sum_{j'=1}^{j}(j'-1)p_{j'}$ and $l_{ij}^{(a)}:=\lambda_{ij}^{(a)}-j+1$.
\end{Prop}

\begin{Rem}
(a) In the simplest case $\pi=1^n$, the above homomorphism
$\CY^\hbar_\pi\to \CY^\hbar_\pi(\gl_n)$ is the classical embedding
of the Yangian of $\ssl_n$ into the Yangian of $\gl_n$.

\noindent
(b) The injectivity of the above homomorphism
 $\CY^\hbar_\pi\to \CY^\hbar_\pi(\gl_n)$
follows from the PBW property for $\CY^\hbar_\pi$
(see~\cite[Corollary 3.15]{fkp}) and its analogue for
$\CY^\hbar_\pi(\gl_n)$.

\noindent
(c) We take this opportunity to correct the sign
in~\cite[(4.2)]{bffr}, where the `$-$' sign should be replaced by
$(-1)^{p_k}$, that is,
$\mathsf{f}_k^{(r+1)}:=(-1)^{p_k}\bq_*(c_1(\CL'_k)^r\cdot \bp^*)$.

\noindent (d)
 We take this opportunity to correct the typos
in~\cite{ffnr}. First, the formulas for the eigenvalues of
$\mathbf{h}_i(u)$ and $\ba_{mi}(u)$ of Theorem 3.20 and its proof
should be corrected by replacing $p_{i'j'}\rightsquigarrow\hbar^{-1}p_{i'j'}$.
Second, the formulas defining $\ba_m(u)$ (Section 2.11),
$\ba_{mi}(u)$ (Section 2.13), $\ba_{mi}(u)$ (Section 3.17) should be
modified by ignoring $\bp_*,\bq^*$.
%
\end{Rem}

\begin{Rem}
Let $e_\pi\in \gl_N$ be a nilpotent element of Jordan type $\pi$. For
$p_1\leq\ldots\leq p_n$, Brundan-Kleshchev proved that the finite
W-algebra $W(\gl_N,e_\pi)$ is the quotient of $\CY^\hbar_\pi(\gl_n)$
by the $2$-sided ideal generated by $\{d_1^{(r)}\}_{r>p_1}$, see~\cite{bk}.
Together with Theorem~\ref{Action on cohomology of Laumon} this yields a natural
action of $W(\gl_N,e_\pi)$ on $V(\pi)$, referred to as a finite
analogue of the AGT relation in~\cite{bffr}. We expect that the
truncated version of  $\CU^\vv_\pi(\gl_n)$ with $\lambda=N\omega_{n-1}$ should
be isomorphic to the $\vv$-version of the $W$-algebra $W(\gl_N,e_\pi)$
as defined by Sevostyanov in~\cite{sev}.
\end{Rem}


\subsection{Shifted quantum toroidal $\ssl_n$ and parabolic affine Laumon spaces}
\label{affine Laumon}
$\ $

The second main result of~\cite{t} provides the action of the
quantum toroidal algebra $U_{\vv,u}(\widehat{\ssl}_n)$ (denoted
$\ddot{\mathrm{U}}_\vv(\widehat{\ssl}_n)$ in \emph{loc.~cit.}) on
the direct sum of localized equivariant $K$-groups of the affine Laumon
spaces $\CP_{\unl{d}}$. The cohomological counterpart of this was
established in~\cite{ffnr}, where the action
of the affine Yangian $Y_{\hbar,\hbar'}(\widehat{\ssl}_n)$ (denoted
$\widehat{Y}$ in \emph{loc.~cit.}) on the direct sum of localized
equivariant cohomology of $\CP_{\unl{d}}$ was constructed.

Likewise, the results of Theorem~\ref{Action on K-theory of Laumon}
and Theorem~\ref{Action on cohomology of Laumon} can be naturally
generalized to provide the actions of the shifted quantum toroidal
algebra $\CU^{\vv,u}_\pi$ (resp.\ shifted affine Yangian
$\CY^{\hbar,\hbar'}_\pi$) on the direct sum of localized equivariant
$K$-groups (resp.\ cohomology) of parabolic affine Laumon spaces.
Here $\CU^{\vv,u}_\pi$ is the associative $\BC(\vv,u)$-algebra generated by
  $\{e_{i,r}, f_{i,r}, \psi^\pm_{i,\pm s^\pm_i}|1\leq i\leq n, r\in \BZ, s^+_i\geq 0, s^-_i\geq p_i-p_{i+1}\}$
and with the same defining relations as for $U_{\vv,u}(\widehat{\ssl}_n)$,
while $\CY^{\hbar,\hbar'}_\pi$ is the associative $\BC(\hbar,\hbar')$-algebra
generated by
 $\{E^{(r+1)}_i, F^{(r+1)}_i, H^{(r+1+p_i-p_{i+1})}_i|1\leq i\leq n, r\in \BN\}$
and with the same defining relations as for
$Y_{\hbar,\hbar'}(\widehat{\ssl}_n)$ (here we set $p_{n+1}:=p_1$).
On the geometric side, the parabolic affine Laumon spaces of type $\pi$
are defined similarly to the case $\pi=1^n$. We leave details to
the interested reader.



\subsection{Whittaker vector}\label{section on whittaker}
$\ $

Consider the \emph{Whittaker vector}
\begin{equation*}
  \fm:=\sum_{\unl{d}} [\CO_{\fQ_{\unl{d}}}]\in M(\pi)^\wedge,
\end{equation*}
where
  $M(\pi)^\wedge:=\prod_{\unl{d}} M(\pi)_{\unl{d}}$.
We also define the operators
\begin{equation*}
   e'_{i,r}:=\bp_*((\vv^i\sL_i)^{\otimes r}\otimes \bq^*)=
   \vv^{i-1}(k^-_{i+1,0})^{-1}e_{i,r}\colon  M(\pi)_{\unl{d}}\to M(\pi)_{\unl{d}-i}.
\end{equation*}

\begin{Prop}\label{Whittaker}
For $1\leq i\leq n-1$, we have
  $$e'_{i,0}(\fm)=(1-\vv^2)^{-1}\fm \ \mathrm{and}\
    e'_{i,1}(\fm)=\ldots=e'_{i,p_i-1}(\fm)=0.$$
\end{Prop}

\begin{proof}
According to the Bott-Lefschetz formula, we have:

\noindent
(1) $\fm=\sum_{\vec{\unl{d}}}a_{\vec{\unl{d}}}[\vec{\unl{d}}]$,
where
  $a_{\vec{\unl{d}}}=\prod_{w\in T_{\vec{\unl{d}}}\fQ_{\unl{d}}}(1-w)^{-1}$;

\noindent
(2) $\frac{a_{\vec{\unl{d}}{}'}}{a_{\vec{\unl{d}}}}
      \bp_*((\vv^i\sL_i)^{\otimes r}\otimes \bq^*)_{[\vec{\unl{d}}{}',\vec{\unl{d}}]}=
      \bq_*((\vv^i\sL_i)^{\otimes r}\otimes \bp^*)_{[\vec{\unl{d}},\vec{\unl{d}}{}']}$.

Set $C_{i,0}:=(1-\vv^2)^{-1}$ and $C_{i,r}:=0\ \mathrm{for}\ 0<r<p_i$.
It suffices to prove the equality
  $C_{i,r}=\sum_{j\leq i}^{a\leq p_j}
   \bq_*((\vv^i\sL_i)^{\otimes r}\otimes \bp^*)_{[\vec{\unl{d}},\vec{\unl{d}}+\delta^{(a)}_{ij}]}$
for any $\vec{\unl{d}}$ and any $1\leq i\leq n-1,0\leq r\leq p_i-1$.
According to Proposition~\ref{matrix_elements}(a), we have
\begin{equation*}
\begin{split}
  & \bq_*((\vv^i\sL_i)^{\otimes r}\otimes \bp^*)_{[\vec{\unl{d}},\vec{\unl{d}}+\delta^{(a)}_{ij}]}=
    (1-\vv^2)^{-1}(s^{(a)}_{ij}\vv^i)^r
    \frac{\prod_{j'\leq i-1}^{a'\leq p_{j'}}(1-s^{(a)}_{ij}/s^{(a')}_{i-1,j'})}
         {\prod_{j'\leq i, a'\leq p_{j'}}^{(j',a')\ne (j,a)}(1-s^{(a)}_{ij}/s^{(a')}_{ij'})}=\\
  & \frac{\vv^i}{1-\vv^2}
    \frac{\prod_{j'\leq i}^{a'\leq p_{j'}} s^{(a')}_{ij'}}
         {\prod_{j'\leq i-1}^{a'\leq p_{j'}} s^{(a')}_{i-1,j'}}\cdot
    (s^{(a)}_{ij}\vv^i)^{r-1}
    \frac{\prod_{j'\leq i-1}^{a'\leq p_{j'}} (s^{(a')}_{i-1,j'}-s^{(a)}_{ij})}
         {\prod_{j'\leq i,a'\leq p_{j'}}^{(j',a')\ne (j,a)} (s^{(a')}_{ij'}-s^{(a)}_{ij})}.
\end{split}
\end{equation*}
For $1\leq r\leq p_i-1$, the sum
\begin{equation*}
    \sum_{j\leq i}^{a\leq p_j} (s^{(a)}_{ij}\vv^i)^{r-1}
    \frac{\prod_{j'\leq i-1}^{a'\leq p_{j'}} (s^{(a')}_{i-1,j'}-s^{(a)}_{ij})}
         {\prod_{j'\leq i,a'\leq p_{j'}}^{(j',a')\ne (j,a)} (s^{(a')}_{ij'}-s^{(a)}_{ij})}
\end{equation*}
is a rational function in
  $\{s^{(a')}_{ij'}\}_{j'\leq i}^{a'\leq p_{j'}}$
of degree $r-p_i<0$ and without poles. Hence, it is zero.
For $r=0$, the same arguments imply
  $$\sum_{j\leq i}^{a\leq p_j} (s^{(a)}_{ij}\vv^i)^{-1}
    \frac{\prod_{j'\leq i-1}^{a'\leq p_{j'}} (s^{(a')}_{i-1,j'}-s^{(a)}_{ij})}
         {\prod_{j'\leq i,a'\leq p_{j'}}^{(j',a')\ne (j,a)} (s^{(a')}_{ij'}-s^{(a)}_{ij})}=
    \sum_{j\leq i}^{a\leq p_j} (s^{(a)}_{ij}\vv^i)^{-1}
    \frac{\prod_{j'\leq i-1}^{a'\leq p_{j'}} s^{(a')}_{i-1,j'}}
         {\prod_{j'\leq i,a'\leq p_{j'}}^{(j',a')\ne (j,a)} (s^{(a')}_{ij'}-s^{(a)}_{ij})}.$$
It remains to compute
  $\sum_{j\leq i}^{a\leq p_j} \prod_{j'\leq i, a'\leq p_{j'}}^{(j',a')\ne (j,a)}
   \frac{s^{(a')}_{ij'}}{s^{(a')}_{ij'}-s^{(a)}_{ij}}$,
which is a rational function in
  $\{s^{(a')}_{ij'}\}_{j'\leq i}^{a'\leq p_{j'}}$
of degree $0$ and without poles, hence, a constant. Specializing
$s^{(1)}_{i1}\mapsto 0$, we see that this constant is equal to 1
(note that only one summand is nonzero under this specialization).

The proposition is proved.
\end{proof}

\begin{Rem}\label{remark on whittaker}
(a) For $\pi=1^n$, this result was proved in~\cite[Proposition 2.31]{bf}.

\noindent
(b) By the same arguments, we also find
  $e''_{i,p_i}(\fm)=\frac{(-1)^{p_i-1}\vv^{ip_i}}{1-\vv^2}\fm$,
where $e''_{i,r}:=(k^-_{i,0})^{2}e'_{i,r}$.

\noindent
(c) Likewise, one can prove that
 $E_i^{(1)}(\fv)=\ldots=E_i^{(p_i-1)}(\fv)=0, E_i^{(p_i)}(\fv)=\hbar^{-1}\fv,$
where $\fv:=\sum_{\unl{d}}[\fQ_{\unl{d}}]\in V(\pi)^\wedge$. This
result was established in~\cite[Proposition 5.1]{bffr}.
\end{Rem}

\appendix


\section{Proof of Theorem~\ref{Levendorskii-type} and its modification}
\label{Proof of Levendorskii-type}

To prove Theorem~\ref{Levendorskii-type}, let us first note that
relations~(\ref{u1}--\ref{u9}) hold in $\CU^\ssc_{0,\mu}$.
Hence, there exists an algebra homomorphism
$\varepsilon\colon \hat{\CU}_{\mu_1,\mu_2}\to \CU^\ssc_{0,\mu}$
such that
  $e_{i,r}\mapsto e_{i,r}, f_{i,s}\mapsto f_{i,s},
   (\psi^+_{i,0})^{\pm 1}\mapsto (\psi^+_{i,0})^{\pm 1},
   (\psi^-_{i,b_i})^{\pm 1}\mapsto (\psi^-_{i,b_i})^{\pm 1},
   h_{i,\pm 1}\mapsto h_{i,\pm 1}$
for $i\in I, b_{2,i}-1\leq r\leq 0, b_{1,i}\leq s\leq 1$.
Moreover, the way we defined
  $e_{i,r}, f_{i,r}, \psi^\pm_{i,r}\in \hat{\CU}_{\mu_1,\mu_2}$
right before Theorem~\ref{Levendorskii-type}, it is clear that
  $\varepsilon\colon  e_{i,r}\mapsto e_{i,r}, f_{i,r}\mapsto f_{i,r},
   \psi^\pm_{i,\pm s^{\pm}_i}\mapsto \psi^\pm_{i,\pm s^{\pm}_i}$
for $i\in I, r\in \BZ, s^+_i\geq 0, s^-_i\geq -b_i$. In particular,
$\varepsilon$ is surjective. Injectivity of $\varepsilon$ is
equivalent to showing that relations~(\ref{U1}--\ref{U8}) hold
in $\hat{\CU}_{\mu_1,\mu_2}$. This occupies the rest of this Appendix
until~\ref{Lev-update},
where we consider a slight modification of this presentation,
see Theorem~\ref{Levendorskii-type update} and its proof.


\subsection{Derivation of some useful relations in $\hat{\CU}_{\mu_1,\mu_2}$}
\

First, we note that~(\ref{u1}, \ref{u4}, \ref{u5}) together
with our definition of $e_{i,r}, f_{i,r}, \psi^+_{i,r}$ imply:
\begin{equation}\tag{v1} \label{v1}
  \psi^+_{i,0}e_{j,r}=\vv_i^{c_{ij}}e_{j,r}\psi^+_{i,0},\
  \psi^-_{i,b_i}e_{j,r}=\vv_i^{-c_{ij}}e_{j,r}\psi^-_{i,b_i},\
  [h_{i,\pm 1}, e_{j,r}]=[c_{ij}]_{\vv_i}\cdot e_{j,r\pm 1},
\end{equation}
\begin{equation}\tag{v2} \label{v2}
  \psi^+_{i,0}f_{j,r}=\vv_i^{-c_{ij}}f_{j,r}\psi^+_{i,0},\
  \psi^-_{i,b_i}f_{j,r}=\vv_i^{c_{ij}}f_{j,r}\psi^-_{i,b_i},\
  [h_{i,\pm 1}, f_{j,r}]=-[c_{ij}]_{\vv_i}\cdot f_{j,r\pm 1},
\end{equation}
\begin{equation}\tag{v3}\label{v3}
  [\psi^+_{i,0},\psi^\pm_{j,\pm s^\pm_j}]=0,\
  [\psi^-_{i,b_i},\psi^\pm_{j,\pm s^\pm_j}]=0
\end{equation}
for any $i,j\in I, r\in \BZ, s^+_j\geq 0, s^-_j\geq -b_j$.

Second, combining relations~(\ref{u1}, \ref{u4}, \ref{u5}, \ref{u6}), we get
\begin{equation}\tag{v4} \label{v4}
\begin{split}
  & [e_{i,1},f_{i,0}]=[e_{i,0},f_{i,1}]=\psi^+_{i,1}/(\vv_i-\vv_i^{-1}),\\
  & [e_{i,b_{2,i}},f_{i,b_{1,i}-1}]=[e_{i,b_{2,i}-1},f_{i,b_{1,i}}]=\psi^-_{i,b_i-1}/(\vv_i^{-1}-\vv_i).
\end{split}
\end{equation}
Note that
  $\psi^+_{i,1}=(\vv_i-\vv_i^{-1})[e_{i,0},f_{i,1}]=(\vv_i-\vv_i^{-1})\psi^+_{i,0}h_{i,1}$.
Hence, $[h_{i,1},\psi^+_{i,1}]=0$. Combining this further
with~(\ref{v1}, \ref{v2}, \ref{v4}) and our definition of $\psi^+_{i,2}$, we obtain
\begin{equation}\tag{v5}\label{v5}
  [e_{i,2},f_{i,0}]=[e_{i,1},f_{i,1}]=[e_{i,0},f_{i,2}]=\psi^+_{i,2}/(\vv_i-\vv_i^{-1}).
\end{equation}
Likewise, we also get
\begin{equation}\tag{v6}\label{v6}
  [e_{i,b_{2,i}},f_{i,b_{1,i}-2}]=[e_{i,b_{2,i}-1},f_{i,b_{1,i}-1}]=
  [e_{i,b_{2,i}-2},f_{i,b_{1,i}}]=\psi^-_{i,b_i-2}/(\vv_i^{-1}-\vv_i).
\end{equation}

Third, let us point out that relation~(\ref{u9}) is equivalent to
\begin{equation}\tag{v7} \label{v7}
  [h_{i,1}, \psi^+_{i,2}]=0,\ [h_{i,-1},\psi^-_{i,b_i-2}]=0.
\end{equation}
According to the above relations, for any $i,j\in I$ we also have
\begin{equation}\tag{v8} \label{v8}
  [h_{j,-1},\psi^+_{i,2}]=0,\ [h_{j,1},\psi^-_{i,b_i-2}]=0.
\end{equation}

Finally, we define elements $h_{i,\pm 2}\in \hat{\CU}_{\mu_1,\mu_2}$
as follows:
\begin{equation} \label{Def-n of h_2}
\begin{split}
  & h_{i,2}:=(\psi^+_{i,0})^{-1}\psi^+_{i,2}/(\vv_i-\vv_i^{-1})-(\vv_i-\vv_i^{-1})h^2_{i,1}/2,\\
  & h_{i,-2}:=(\psi^-_{i,b_i})^{-1}\psi^-_{i,b_i-2}/(\vv_i^{-1}-\vv_i)-(\vv_i^{-1}-\vv_i)h^2_{i,-1}/2.
\end{split}
\end{equation}
 Due to relations~(\ref{u1}, \ref{v7}, \ref{v8}), for every $i,j\in I$ we have
\begin{equation}\tag{v9}\label{v9}
  [h_{i,\pm 1},h_{i,\pm 2}]=0,\ [h_{j,\mp 1},h_{i,\pm 2}]=0.
\end{equation}

\begin{Lem}\label{Shifts by 2}
For any $i\in I, r\in \BZ$, we have
  $$[h_{i,\pm 2},e_{i,r}]=\frac{[4]_{\vv_i}}{2}\cdot e_{i,r\pm 2},\
    [h_{i,\pm 2},f_{i,r}]=-\frac{[4]_{\vv_i}}{2}\cdot f_{i,r\pm 2}.$$
\end{Lem}

\begin{proof}
Due to~(\ref{u2}), we have $[e_{i,0},e_{i,-1}]_{\vv_i^2}=0$.
Commuting this with $h_{i,1}$ and applying relation~(\ref{u4}), we obtain
  $e_{i,1}e_{i,-1}-\vv_i^2e_{i,0}^2=\vv_i^2e_{i,-1}e_{i,1}-e^2_{i,0}.$
Commuting this further with $f_{i,1}$ and applying
relation~(\ref{u6}), we obtain
\begin{equation*}
\begin{split}
  & \psi^+_{i,2}e_{i,-1}-\vv_i^2\psi^+_{i,1}e_{i,0}+e_{i,1}\psi^+_{i,0}-
  \vv_i^2e_{i,0}\psi^+_{i,1}-\delta_{b_i,0}e_{i,1}\psi^-_{i,b_i}=\\
  & \vv_i^2e_{i,-1}\psi^+_{i,2}-e_{i,0}\psi^+_{i,1}+\vv_i^2\psi^+_{i,0}e_{i,1}-
  \psi^+_{i,1}e_{i,0}-\vv_i^2\delta_{b_i,0}\psi^-_{i,b_i}e_{i,1}.
\end{split}
\end{equation*}
First, note that
  $e_{i,1}\psi^-_{i,b_i}=\vv_i^2\psi^-_{i,b_i}e_{i,1}$,
due to~(\ref{u4}). Second, we have
\begin{equation}\tag{v10} \label{v10}
  e_{i,1}\psi^+_{i,0}-\vv_i^2e_{i,0}\psi^+_{i,1}=\vv_i^2\psi^+_{i,0}e_{i,1}-\psi^+_{i,1}e_{i,0}.
\end{equation}
Indeed, due to the equality $\psi^+_{i,1}=(\vv_i-\vv_i^{-1})\psi^+_{i,0}h_{i,1}$ and
relations~(\ref{u1}, \ref{v1}), we have
\begin{equation*}
  \psi^+_{i,1}e_{i,0}-\vv_i^2e_{i,0}\psi^+_{i,1}=\vv_i^2(\vv_i-\vv_i^{-1})[2]_{\vv_i}\cdot e_{i,1}\psi^+_{i,0}=
  (\vv_i^4-1)e_{i,1}\psi^+_{i,0}=\vv_i^2\psi^+_{i,0}e_{i,1}-e_{i,1}\psi^+_{i,0}.
\end{equation*}
\noindent
 Therefore, we get
\begin{equation}\tag{v11} \label{v11}
  \psi^+_{i,2}e_{i,-1}-\vv_i^2\psi^+_{i,1}e_{i,0}=
  \vv_i^2e_{i,-1}\psi^+_{i,2}-e_{i,0}\psi^+_{i,1}.
\end{equation}
Combining the formulas
  $\psi^+_{i,1}=(\vv_i-\vv_i^{-1})\psi^+_{i,0}h_{i,1},\
   \psi^+_{i,2}=(\vv_i-\vv_i^{-1})\psi^+_{i,0}(h_{i,2}+\frac{\vv_i-\vv_i^{-1}}{2}h^2_{i,1})$
with relations~(\ref{u1}, \ref{v1}, \ref{v11}), we finally get
$[h_{i,2},e_{i,-1}]=\frac{[4]_{\vv_i}}{2}e_{i,1}$. Commuting this
relation with $h_{i,\pm 1}$ and using~(\ref{v1}, \ref{v9}), we obtain
$[h_{i,2},e_{i,r}]=\frac{[4]_{\vv_i}}{2}e_{i,r+2}$ for any $r\in \BZ$.

Likewise, starting from the relation
$[e_{i,b_{2,i}},e_{i,b_{2,i}-1}]_{\vv_i^2}=0$ and commuting it
first with $h_{i,-1}$ and then with $f_{i,b_{1,i}}$, we recover
$[h_{i,-2},e_{i,b_{2,i}}]=\frac{[4]_{\vv_i}}{2}e_{i,b_{2,i}-2}$.
Commuting this further with $h_{i,\pm 1}$, we get
 $[h_{i,-2},e_{i,r}]=\frac{[4]_{\vv_i}}{2}e_{i,r-2}$ for any $r\in \BZ$.

 The proof of $[h_{i,\pm 2},f_{i,r}]=-\frac{[4]_{\vv_i}}{2}\cdot f_{i,r\pm 2}$ is completely analogous.
\end{proof}


\subsection{Verification of relations~(\ref{U1}--\ref{U6}) with $i=j$ for $\hat{\CU}_{\mu_1,\mu_2}$}
\


\subsubsection{Verification of~(\ref{U2})}
\

We need to prove $X^+(i;r,s)=0$ for any $r,s\in \BZ$, where
\begin{equation*}
  X^+(i;r,s):=[e_{i,r+1},e_{i,s}]_{\vv_i^2}+[e_{i,s+1},e_{i,r}]_{\vv_i^2}.
\end{equation*}
Note that $X^+(i;r,s)=X^+(i;s,r)$, and $X^+(i;-1,-1)=0$ due to relation~(\ref{u2}).

For $a\in \{\pm 1, \pm 2\}$, we define
  $L_{i,a}:=a/[2a]_{\vv_i}\cdot \ad(h_{i,a})\in \mathrm{End}(\hat{\CU}_{\mu_1,\mu_2})$.
Then, we have $L_{i,a}(X^+(i;r,s))=X^+(i;r+a,s)+X^+(i;r,s+a)$. Set
$L^\pm_i:=\frac{1}{2}(L^2_{i,\pm 1}-L_{i,\pm 2})$. Then
$L^\pm_i(X^+(i;r,s))=X^+(i;r\pm 1,s\pm 1)$.
Applying iteratively $L^+_i$ to the equality $X^+(i;-1,-1)=0$,
we get $X^+(i;r,r)=0$ for any $r\geq -1$. Since
$2X^+(i;-1,0)=L_{i,1}(X^+(i;-1,-1))=0$, we analogously get
$X^+(i;r,r+1)=0$ for $r\geq -1$.
Fix $s\in \BZ_{>0}$ and assume by induction that $X^+(i;r,r+N)=0$
for any $r\geq -1, 0\leq N\leq s$.
Then $X^+(i;-1,s)=L_{i,1}(X^+(i;-1,s-1))-X^+(i;0,s-1)=0$, due to the
above assumption. Applying $(L^+_i)^{r+1}$ to the latter equality,
we get $X^+(i;r,r+s+1)=0$ for $r\geq -1$. An induction in $s$
completes the proof of $X^+(i;r,s)=0$ for any $r,s\geq -1$. Finally,
applying iteratively $L^-_i$, we obtain $X^+(i;r,s)=0$ for any
$r,s\in \BZ$.


\subsubsection{Verification of~(\ref{U3})}
\

 This relation is verified completely analogously to~(\ref{U2}).


\subsubsection{Verification of~(\ref{U4})}
\

We consider the case $\epsilon=+$ (the case $\epsilon=-$ is completely analogous).
We need to prove $Y^+(i;r,s)=0$ for any $r\in \BN,s\in \BZ$, where
\begin{equation*}
  Y^+(i;r,s):=[\psi^+_{i,r+1},e_{i,s}]_{\vv_i^2}+[e_{i,s+1},\psi^+_{i,r}]_{\vv_i^2}.
\end{equation*}
The $r=s=0$ case is due to~(\ref{v10}) from our proof of
Lemma~\ref{Shifts by 2}. Moreover, the same argument also yields
$Y^+(i;0,s)=0$ for any $s\in \BZ$.

Note that
  $Y^+(i;r,s-1)+Y^+(i;s,r-1)=(\vv_i-\vv_i^{-1})[X^+(i;r-1,s-1),f_{i,1}]=0$
for $r,s\geq 0$. The first equality is due to~(\ref{v1}) and our
definition of $\psi^+_{i,r}$, while the second equality follows from
$X^+(i;r-1,s-1)=0$ proved above. In particular,
$Y^+(i;r,-1)+Y^+(i;0,r-1)=0$ for $r\in \BN$.

Combining the above two observations, we find
\begin{equation}\tag{v12} \label{v12}
  Y^+(i;r,-1)=0\ \mathrm{for\ any}\ r\in \BN.
\end{equation}
 Commuting iteratively the equality $Y^+(i;1,-1)=0$ with $h_{i,\pm 1}$,
we get $Y^+(i;1,s)=0$ for any $s\in \BZ$, due to~(\ref{u1}, \ref{v1}, \ref{v9}).


\medskip
 Next, we prove the following five statements by induction in $N\in \BZ_+$:

\noindent
($A_N$) $[h_{i,1},\psi^+_{i,r}]=0$ for $0\leq r\leq N+1$;

\noindent
($B_N$) $[h_{i,-1},\psi^+_{i,r}]=0$ for $0\leq r\leq N+1$;

\noindent
($C_N$) $[e_{i,r},f_{i,s}]=\psi^+_{i,r+s}/(\vv_i-\vv_i^{-1})$
for any $r,s\in \BN$ with $1\leq r+s\leq N+2$;

\noindent
($D_N$) $Y^+(i;r,s)=0$ for any $0\leq r\leq N, s\in \BZ$;

\noindent
($E_N$) $[\psi^+_{i,r},\psi^+_{i,s}]=0$
for any $r,s\geq 0$ with $r+s\leq N+2$.

\medskip
\noindent
\emph{\underline{Base of Induction ($N=1$)}}:
The assertions ($A_1$, $B_1$, $D_1$, $E_1$) have been already proved
above, while ($C_1$) follows immediately from
$[h_{i,1},\psi^+_{i,2}]=0$ (cf.~(\ref{v7})) and~(\ref{v1}, \ref{v2}, \ref{v4}, \ref{v5}).

\medskip
\noindent
\emph{\underline{Induction Step}}:
Assuming ($A_N$--$E_N$) for a given $N\in \BZ_{>0}$, we prove ($A_{N+1}$--$E_{N+1}$).

\begin{proof}[Proof of the induction step]
\

Consider a polynomial algebra $B:=\BC(\vv)[\{x_r\}_{r=1}^\infty]$,
which is $\BN$-graded via $\deg(x_r)=r$. Define elements
$\{\sh_r\}_{r=1}^\infty$ of $B$ via
  $\exp\left((\vv_i-\vv_i^{-1})\sum_{r=1}^\infty \sh_rz^{-r}\right)=1+\sum_{r=1}^\infty x_r z^{-r}.$
Then, $\sh_r=\frac{x_r}{\vv_i-\vv_i^{-1}}+p_r(x_1,\ldots,x_{r-1})$
with polynomials $p_r$ satisfying $\deg(p_r(x_1,\ldots,x_{r-1}))=r$.

Using the above polynomials $p_r$, we define
$h_{i,1},\ldots,h_{i,N+1}\in \hat{\CU}_{\mu_1,\mu_2}$ via
\begin{equation}\label{Def-n of h_k}
  h_{i,r}:=\frac{(\psi^+_{i,0})^{-1}\psi^+_{i,r}}{\vv_i-\vv_i^{-1}}+
           p_r((\psi^+_{i,0})^{-1}\psi^+_{i,1},\ldots, (\psi^+_{i,0})^{-1}\psi^+_{i,r-1})
  \ \mathrm{for}\ 1\leq r\leq N+1.
\end{equation}
These $h_{i,r}$ are well-defined and are independent of the choice of $N>r-1$,
due to the assumption ($E_N$) and the
aforementioned degree condition on $p_r$. The following is
straightforward:\footnote{If we knew that
$[\psi^+_{i,a},\psi^+_{i,b}]=0$ for any $0\leq a,b\leq N+1$, then~(\ref{Shifts by k})
would immediately follow from ($D_N$) by the standard arguments.
However, every monomial appearing in $p_r$ involves only pairwise
commuting $\psi^+_{i,a}$'s, due to the degree condition on $p_r$ and
the assumption ($E_N$). Hence, the equality~(\ref{Shifts by k})
follows formally from its validity in the aforementioned simpler
case ($[\psi^+_{i,a},\psi^+_{i,b}]=0$ for any $0\leq a,b\leq N+1$).}
\begin{equation}\tag{v13}\label{Shifts by k}
  [h_{i,r},e_{i,s}]=\frac{[2r]_{\vv_i}}{r}\cdot e_{i,s+r}
  \ \mathrm{for}\ 1\leq r\leq N+1, s\in \BZ.
\end{equation}


\medskip
\noindent
$\circ$ \emph{Validity of ($A_{N+1}$).}
We need to prove $[h_{i,1},\psi^+_{i,N+2}]=0$. According to ($C_N$), we have
$\psi^+_{i,N+2}=(\vv_i-\vv_i^{-1})[e_{i,N+2-r},f_{i,r}]$ for $0\leq r\leq N+2$.
Hence,
\begin{equation}\tag{v14}\label{v14}
  [h_{i,1},\psi^+_{i,N+2}]/(\vv_i^2-\vv_i^{-2})=[e_{i,N+3-r},f_{i,r}]-[e_{i,N+2-r},f_{i,r+1}]
  \ \mathrm{for}\ 0\leq r\leq N+2.
\end{equation}
Adding up these equalities for $r=0,1$ and using~Lemma~\ref{Shifts by 2}
together with the assumption $(C_N)$, we get
\begin{equation*}
    \frac{2[h_{i,1},\psi^+_{i,N+2}]}{\vv_i^2-\vv_i^{-2}}=[e_{i,N+3},f_{i,0}]-[e_{i,N+1},f_{i,2}]=
    \frac{2}{[4]_{\vv_i}}\cdot [h_{i,2},[e_{i,N+1},f_{i,0}]]=\frac{2[h_{i,2},\psi^+_{i,N+1}]}{\vv_i^4-\vv_i^{-4}}.
\end{equation*}
Likewise, adding up the equality~(\ref{v14}) for $r=0,1,\ldots,N$
and using~(\ref{Shifts by k}), we obtain
\begin{equation*}
    \frac{N+1}{\vv_i^2-\vv_i^{-2}}[h_{i,1},\psi^+_{i,N+2}]=
    \frac{N+1}{[2(N+1)]_{\vv_i}}\cdot [h_{i,N+1},[e_{i,2},f_{i,0}]]=
    \frac{(N+1)[h_{i,N+1},\psi^+_{i,2}]}{\vv_i^{2(N+1)}-\vv_i^{-2(N+1)}}.
\end{equation*}
Comparing the above two equalities, we find
\begin{equation}\tag{v15}\label{v15}
  [h_{i,1},\psi^+_{i,N+2}]=
  \frac{\vv_i^2-\vv_i^{-2}}{\vv_i^4-\vv_i^{-4}}[h_{i,2},\psi^+_{i,N+1}]=
  \frac{\vv_i^2-\vv_i^{-2}}{\vv_i^{2(N+1)}-\vv_i^{-2(N+1)}}[h_{i,N+1},\psi^+_{i,2}].
\end{equation}

On the other hand, combining~(\ref{Def-n of h_k}) with the assumption ($E_N$), we get
\begin{equation*}
    [h_{i,s},\psi^+_{i,N+3-s}]=(\psi^+_{i,0})^{-1}[\psi^+_{i,s},\psi^+_{i,N+3-s}]/(\vv_i-\vv_i^{-1})
    \ \mathrm{for}\ 1\leq s\leq N+1.
\end{equation*}
Hence,
\begin{equation}\tag{v16}\label{v16}
   [h_{i,1},\psi^+_{i,N+2}]=
   \frac{(\psi^+_{i,0})^{-1}[\psi^+_{i,2},\psi^+_{i,N+1}]}{(\vv_i-\vv_i^{-1})[2]_{\vv_i^2}}=
   \frac{(\psi^+_{i,0})^{-1}[\psi^+_{i,2},\psi^+_{i,N+1}]}{(\vv_i-\vv_i^{-1})[-N-1]_{\vv_i^2}}.
\end{equation}
Since $[2]_{\vv_i^2}\ne [-N-1]_{\vv_i^2}$, the second equality
of~(\ref{v16}) implies $[\psi^+_{i,2},\psi^+_{i,N+1}]=0$. Hence,
$[h_{i,1},\psi^+_{i,N+2}]=0$, and ($A_{N+1}$) follows.


\medskip
\noindent
$\circ$ \emph{Validity of ($B_{N+1}$).}
We need to prove $[h_{i,-1}, \psi^+_{i,N+2}]=0$. This follows from
 $[h_{i,-1}, \psi^+_{i,N+2}]=(\vv_i-\vv_i^{-1})[2]_{\vv_i}\cdot ([e_{i,N},f_{i,1}]-[e_{i,N+1},f_{i,0}])=0,$
where we used~(\ref{v1}, \ref{v2}) in the first equality and ($C_N$) in
the second one. Hence, ($B_{N+1}$) holds.


\medskip
\noindent
$\circ$ \emph{Validity of ($C_{N+1}$).}
According to ($C_N$), we have
$\psi^+_{i,N+2}=(\vv_i-\vv_i^{-1})[e_{i,r},f_{i,N+2-r}]$
for any $0\leq r\leq N+2$. Therefore,
  $[h_{i,1},\psi^+_{i,N+2}]=(\vv_i^2-\vv_i^{-2})([e_{i,r+1},f_{i,N+2-r}]-[e_{i,r},f_{i,N+3-r}])$ due to~(\ref{v1}, \ref{v2}).
The left-hand side is zero due to ($A_{N+1}$) established above, hence
\begin{equation*}
  [e_{i,N+3},f_{i,0}]=[e_{i,N+2},f_{i,1}]=\ldots=[e_{i,1},f_{i,N+2}]=[e_{i,0},f_{i,N+3}].
\end{equation*}
Combining this with our definition
  $\psi^+_{i,N+3}=(\vv_i-\vv_i^{-1})[e_{i,N+2},f_{i,1}]$
yields ($C_{N+1}$).


\medskip
\noindent
$\circ$ \emph{Validity of ($D_{N+1}$).}
Due to ($A_{N+1}$) and ($B_{N+1}$) established above, we have
$[h_{i,\pm 1},Y^+(i;N+1,s)]=[2]_{\vv_i}\cdot Y^+(i;N+1,s\pm 1)$.
Combining this with~(\ref{v12}), we see that $Y^+(i;N+1,s)=0$ for
any $s\in \BZ$. Hence, ($D_{N+1}$) holds.


\medskip
\noindent
$\circ$ \emph{Validity of ($E_{N+1}$).}
We need to prove  $[\psi^+_{i,r},\psi^+_{i,N+3-r}]=0$ for any $1\leq r\leq N+1$.
Equivalently, it suffices to prove $[h_{i,r},\psi^+_{i,N+3-r}]=0$
for $1\leq r\leq N+1$. According to ($C_N$), we have
$\psi^+_{i,N+3-r}=(\vv_i-\vv_i^{-1})[e_{i,N+3-r},f_{i,0}]$.
Therefore,
  $[h_{i,r}, \psi^+_{i,N+3-r}]=\frac{\vv_i^{2r}-\vv_i^{-2r}}{r}\cdot
   ([e_{i,N+3},f_{i,0}]-[e_{N+3-r},f_{i,r}])=0,$
due to~(\ref{Shifts by k}) and the assertion ($C_{N+1}$) proved above.
\end{proof}

The induction step is accomplished.
In particular, ($D_N$) completes our verification of~(\ref{U4}) with $i=j$.


\subsubsection{Verification of~(\ref{U5})}
\

This relation is verified completely analogously to~(\ref{U4}).


\subsubsection{Verification of~(\ref{U6})}
\

We need to prove
  $$[e_{i,r}, f_{i,N-r}]=\frac{1}{\vv_i-\vv_i^{-1}}\cdot
   \begin{cases}
      \psi^+_{i,N}-\delta_{N,0}\delta_{b_i,0}\psi^-_{i,b_i} & \text{if}\ \ N\geq 0, \\
      -\psi^-_{i,N}+\delta_{N,0}\delta_{b_i,0}\psi^+_{i,0}& \text{if}\ \ N\leq b_i, \\
      0 & \text{if}\ \ b_i<N<0.
   \end{cases}$$
Note that given any value of $N\in \BZ$, we know this equality for a
certain value of $r\in \BZ$.

\medskip
\noindent
$\circ$ \emph{Case $N>0$}.
If $0\leq r\leq N$, then
$[e_{i,r},f_{i,N-r}]=\psi^+_{i,N}/(\vv_i-\vv_i^{-1})$, due to ($C_N$).
For $r<0$, we proceed by induction in $|r|$.
Due to~(\ref{v1}, \ref{v2}), we have
  $[e_{i,r},f_{i,N-r}]=[2]_{\vv_i}^{-1}\cdot [[h_{i,-1},e_{i,r+1}],f_{i,N-r}]=
   [2]_{\vv_i}^{-1}\cdot [h_{i,-1},[e_{i,r+1},f_{i,N-r}]]+[e_{i,r+1},f_{i,N-r-1}]=\psi^+_{i,N},$
where in the last equality we used the induction assumption and the
equality $[h_{i,-1}, \psi^+_{i,N+1}]=0$, due to ($B_N$). The case
$l:=N-r<0$ is treated in the same way.

\medskip
\noindent
$\circ$ \emph{Case $N\leq 0$}.
We proceed by induction in $|N|$. For any $r\in\BZ$, we have
\begin{equation*}
    [e_{i,r},f_{i,N-r}]=[2]_{\vv_i}^{-1}\cdot [h_{i,-1},[e_{i,r+1},f_{i,N-r}]]+[e_{i,r+1},f_{i,N-r-1}]=
    [e_{i,r+1},f_{i,N-r-1}],
\end{equation*}
where we used the induction assumption together
with~(\ref{u1}, \ref{v1}, \ref{v2}) and $[h_{i,-1},\psi^-_i(z)]=0$ (the latter is
proved completely analogously to ($A_N$)).
Hence, the expression $[e_{i,r},f_{i,N-r}]$ is independent of $r\in \BZ$.
The result follows since we know the equality holds for a certain value of $r$.


\subsubsection{Verification of~(\ref{U1})}
\

We consider the case $\epsilon=+$ (the case $\epsilon=-$ is
completely analogous). We need to prove
  $[\psi^+_{i,r},\psi^+_{i,s_i^+}]=[\psi^+_{i,r},\psi^-_{i,-s_i^-}]=0$
for any $r,s_i^+\geq 0, s_i^-\geq -b_i$. This is clear for $r=0$ or $s_i^+=0$,
or $s_i^-=-b_i$, due to~(\ref{v3}). Therefore, it remains to
prove $[h_{i,r},\psi^+_{i,s_i^+}]=0$ and
$[h_{i,r},\psi^-_{i,-s_i^-}]=0$ for $r>0, s_i^+>0, s_i^->-b_i$.

For $s_i^+>0$, we have
  $\psi^+_{i,s_i^+}=(\vv_i-\vv_i^{-1})[e_{i,s_i^+-1},f_{i,1}]$,
so that
\begin{equation*}
    [h_{i,r},\psi^+_{i,s_i^+}]=
    \frac{[2r]_{\vv_i}}{r}(\vv_i-\vv_i^{-1})\cdot ([e_{i,s_i^+ + r-1},f_{i,1}]-[e_{i,s_i^+ -1},f_{i,r+1}])=0,
\end{equation*}
where the first equality is due to~(\ref{Shifts by k}), while the
second equality is due to relation~(\ref{U6}) with $i=j$ proved above.

For $s_i^->-b_i$, we have
  $\psi^-_{i,-s_i^-}=(\vv_i^{-1}-\vv_i)[e_{i,-b_{1,i}-s_i^-},f_{i,b_{1,i}}]$,
so that
\begin{equation*}
    [h_{i,r},\psi^-_{i,-s_i^-}]=\frac{[2r]_{\vv_i}}{r}(\vv_i^{-1}-\vv_i)\cdot
    ([e_{i,r-b_{1,i}-s_i^-},f_{i,b_{1,i}}]-[e_{i,-b_{1,i}-s_i^-},f_{i,r+b_{1,i}}])=0,
\end{equation*}
where the first equality is due to~(\ref{Shifts by k}), while the
second equality is due to relation~(\ref{U6}) with $i=j$ proved above.

\medskip
This completes our verification of relations~(\ref{U1}--\ref{U6}) with $i=j$
for $\hat\CU_{\mu_1,\mu_2}$.


\subsection{Verification of relations~(\ref{U1}--\ref{U8}) with $i\ne j$ for $\hat{\CU}_{\mu_1,\mu_2}$}
\


\subsubsection{Verification of~(\ref{U2})}
\

We need to prove $X^+(i,j;r,s)=0$ for any $r,s\in \BZ$, where
\begin{equation*}
  X^+(i,j;r,s):=[e_{i,r+1},e_{j,s}]_{\vv_i^{c_{ij}}}+[e_{j,s+1},e_{i,r}]_{\vv_i^{c_{ij}}}.
\end{equation*}
First, the equality $X^+(i,j;-1,-1)=0$ follows from~(\ref{u2}).
Second, due to~(\ref{v1}) we have
\begin{equation*}
\begin{split}
  & [h_{i,1},X^+(i,j;r,s)]=[c_{ii}]_{\vv_i}\cdot X^+(i,j;r+1,s)+[c_{ij}]_{\vv_i}\cdot X^+(i,j;r,s+1),\\
  & [h_{j,1},X^+(i,j;r,s)]=[c_{ji}]_{\vv_j}\cdot X^+(i,j;r+1,s)+[c_{jj}]_{\vv_j}\cdot X^+(i,j;r,s+1).
\end{split}
\end{equation*}
Combining these equalities with nondegeneracy of the  matrix
  $A_{ij}:=
   \begin{bmatrix}
     [c_{ii}]_{\vv_i} & [c_{ij}]_{\vv_i} \\
     [c_{ji}]_{\vv_j} & [c_{jj}]_{\vv_j}
   \end{bmatrix}$,
we see that
  $X^+(i,j;r,s)=0\Rightarrow X^+(i,j;r+1,s)=0, X^+(i,j;r,s+1)=0$.
Since $X^+(i,j;-1,-1)=0$, we get $X^+(i,j;r,s)=0$ for $r,s\geq -1$
by induction in $r,s$.

A similar reasoning with $h_{i,-1}, h_{j,-1}$ used instead of
$h_{i,1},h_{j,1}$ yields the implication
\begin{equation*}
  X^+(i,j;r,s)=0\Longrightarrow X^+(i,j;r-1,s)=0,\ X^+(i,j;r,s-1)=0.
\end{equation*}
Hence, an induction argument completes the proof of $X^+(i,j;r,s)=0$
for any $r,s\in \BZ$.


\subsubsection{Verification of~(\ref{U3})}
\

We need to prove $X^-(i,j;r,s)=0$ for any $r,s\in \BZ$, where
\begin{equation*}
  X^-(i,j;r,s):=[f_{i,r+1},f_{j,s}]_{\vv_i^{-c_{ij}}}+[f_{j,s+1},f_{i,r}]_{\vv_i^{-c_{ij}}}.
\end{equation*}
The $r=s=0$ case follows from~(\ref{u3}). The general
case follows from
\begin{equation*}
  X^-(i,j;r,s)=0\Longrightarrow X^-(i,j;r\pm 1,s)=0,\ X^-(i,j;r,s\pm 1)=0
\end{equation*}
applied iteratively to $X^-(i,j;0,0)=0$, in the same vein as in the above
verification of~(\ref{U2}).


\subsubsection{Verification of~(\ref{U6})}
\

We need to prove $X(i,j;r,s)=0$ for any $r,s\in \BZ$, where
\begin{equation*}
  X(i,j;r,s):=[e_{i,r},f_{j,s}].
\end{equation*}
First, the equality $X(i,j;0,0)=0$ follows from~(\ref{u6}). Second,
due to~(\ref{v1}, \ref{v2}) we have
\begin{equation*}
\begin{split}
  & [h_{i,\pm 1},X(i,j;r,s)]=[c_{ii}]_{\vv_i}\cdot X(i,j;r\pm 1,s)-[c_{ij}]_{\vv_i}\cdot X(i,j;r,s\pm 1),\\
  & [h_{j,\pm 1},X(i,j;r,s)]=[c_{ji}]_{\vv_j}\cdot X(i,j;r\pm 1,s)-[c_{jj}]_{\vv_j}\cdot X(i,j;r,s\pm 1).
\end{split}
\end{equation*}
Combining these equalities with nondegeneracy of the  matrix
  $B_{ij}:=
   \begin{bmatrix}
     [c_{ii}]_{\vv_i} & -[c_{ij}]_{\vv_i} \\
     [c_{ji}]_{\vv_j} & -[c_{jj}]_{\vv_j}
   \end{bmatrix}$,
we see that
  $X(i,j;r,s)=0\Rightarrow X(i,j;r\pm 1,s)=0, X(i,j;r,s\pm 1)=0$.
Hence, the equality $X(i,j;r,s)=0$ for any $r,s\in\BZ$
follows from the $r=s=0$ case considered above.


\subsubsection{Verification of~(\ref{U4})}
\

We consider the case $\epsilon=+$ (the case $\epsilon=-$ is
completely analogous). We need to prove $Y^+(i,j;r,s)=0$ for any
$r\in \BN, s\in \BZ$, where
\begin{equation*}
  Y^+(i,j;r,s):=[\psi^+_{i,r+1},e_{j,s}]_{\vv_i^{c_{ij}}}+[e_{j,s+1},\psi^+_{i,r}]_{\vv_i^{c_{ij}}}.
\end{equation*}
Due to relation~(\ref{U6}) (established already both for $i=j$ and $i\ne j$),
we have
\begin{equation*}
  (\vv_i-\vv_i^{-1})[[e_{i,r+1},e_{j,s}]_{\vv_i^{c_{ij}}},f_{i,0}]=[\psi^+_{i,r+1},e_{j,s}]_{\vv_i^{c_{ij}}},
\end{equation*}
\begin{equation*}
    (\vv_i-\vv_i^{-1})[[e_{j,s+1},e_{i,r}]_{\vv_i^{c_{ij}}},f_{i,0}]=
    [e_{j,s+1}, \psi^+_{i,r}-\delta_{r,0}\delta_{b_i,0}\psi^-_{i,-b_i}]_{\vv_i^{c_{ij}}}=
    [e_{j,s+1}, \psi^+_{i,r}]_{\vv_i^{c_{ij}}}.
\end{equation*}
Therefore,
  $Y^+(i,j;r,s)=(\vv_i-\vv_i^{-1})[X^+(i,j;r,s),f_{i,0}]=0$,
where the last equality follows from $X^+(i,j;r,s)=0$ proved above.


\subsubsection{Verification of~(\ref{U5})}
\

We consider the case $\epsilon=+$ (the case $\epsilon=-$ is
completely analogous). We need to prove $Y^-(i,j;r,s)=0$ for any
$r\in \BN, s\in \BZ$, where
\begin{equation*}
  Y^-(i,j;r,s):=[\psi^+_{i,r+1},f_{j,s}]_{\vv_i^{-c_{ij}}}+[f_{j,s+1},\psi^+_{i,r}]_{\vv_i^{-c_{ij}}}.
\end{equation*}
Analogously to our verification of~(\ref{U4}), we have
  $Y^-(i,j;r,s)=(\vv_i-\vv_i^{-1})[e_{i,0},X^-(i,j;r,s)]$.
Thus, the equality $Y^-(i,j;r,s)=0$ follows from $X^-(i,j;r,s)=0$
proved above.


\subsubsection{Verification of~(\ref{U1})}
\

We consider the case $\epsilon=\epsilon'=+$ (other cases are
completely analogous). Due to relation~(\ref{v3}), it suffices
to prove $[h_{i,r},\psi^+_{j,s}]=0$ for $r,s\in \BZ_{>0}$, where the
elements $\{h_{i,r}\}_{r=1}^\infty$ were defined in~(\ref{Def-n of h_k}).

Analogously to~(\ref{Shifts by k}), relations~(\ref{U4}, \ref{U5}) imply
\begin{equation*}
    [h_{i,r},e_{j,s}]=\frac{[rc_{ij}]_{\vv_i}}{r}\cdot e_{j,s+r},\
    [h_{i,r},f_{j,s}]=-\frac{[rc_{ij}]_{\vv_i}}{r}\cdot f_{j,s+r}\ \mathrm{for\ any}\ r\in \BZ_{>0}, s\in \BZ.
\end{equation*}
Hence, we have
\begin{equation*}
    [h_{i,r}, \psi^+_{j,s}]=(\vv_i-\vv_i^{-1})[h_{i,r},[e_{j,s},f_{j,0}]]=
    (\vv_i-\vv_i^{-1})\frac{[rc_{ij}]_{\vv_i}}{r}\cdot ([e_{j,s+r},f_{j,0}]-[e_{j,s},f_{j,r}])=0,
\end{equation*}
where the first and the last equalities follow from~(\ref{U6}) with
$i=j$ established above.


\subsubsection{Verification of~(\ref{U7})}
\

In the simplest case $c_{ij}=0$, we need to prove $[e_{i,r},e_{j,s}]=0$ for any
$r,s\in \BZ$. The equality $[e_{i,0},e_{j,0}]=0$ is due to~(\ref{u7}), while
commuting it iteratively with $h_{i,\pm1},h_{j,\pm1}$, we get $[e_{i,r},e_{j,s}]=0$,
due to~(\ref{v1}, \ref{v2}).

In general, we set $m:=1-c_{ij}$. For any
$\vec{r}=(r_1,\ldots,r_m)\in \BZ^m$ and $s\in \BZ$, define
\begin{equation*}
    Z^+(i,j;\vec{r},s):=\sum_{\pi\in\fS_m}\sum_{t=0}^m(-1)^t{m\brack t}_{\vv_i}
    e_{i,r_{\pi(1)}}\cdots e_{i,r_{\pi(t)}}e_{j,s}e_{i,r_{\pi(t+1)}}\cdots e_{i,r_{\pi(m)}}.
\end{equation*}
To check~(\ref{U7}), we need to prove $Z^+(i,j;\vec{r},s)=0$ for any
$\vec{r}\in \BZ^m,s\in \BZ$.

Let $\vec{0}=(0,\ldots,0)\in \BZ^m$. The equality
$Z^+(i,j;\vec{0},0)=0$ follows from~(\ref{u7}) (cf.~Remark~\ref{56}). Commuting
$Z^+(i,j;\vec{0},s)$ with $h_{i,\pm 1}, h_{j,\pm 1}$, and using
nondegeneracy of the matrix $A_{ij}$, we get
$Z^+(i,j;\vec{0},s)=0\Rightarrow Z^+(i,j;\vec{0},s\pm 1)=0$.
Therefore, $Z^+(i,j;\vec{0},s)=0$ for any $s\in \BZ$.

Next, we prove that $Z^+(i,j;\vec{r},s)=0$ for any
$\vec{r}=(r_1,\ldots,r_k,0,\ldots,0)\in \BZ^m,\ s\in \BZ$ by induction in $0\leq k\leq m$.
The base case $k=0$ was just treated above.
For the induction step, note that the commutator
$[h_{i,r'},Z^+(i,j;\vec{r},s)]$ equals
$\frac{(m-k)\cdot[2r']_{\vv_i}}{r'}Z^+(i,j;(r_1,\ldots,r_k,r',0,\ldots,0),s)$
plus some other terms which are zero by the induction assumption.
Hence, $Z^+(i,j;\vec{r},s)=0$ for any $\vec{r}\in \BZ^m, s\in \BZ$.

%


\subsubsection{Verification of~(\ref{U8})}
\

Set $m:=1-c_{ij}$. For any $\vec{r}\in \BZ^m,s\in \BZ$, define
\begin{equation*}
    Z^-(i,j;\vec{r},s):=\sum_{\pi\in\fS_m}\sum_{t=0}^m(-1)^t{m\brack t}_{\vv_i}
    f_{i,r_{\pi(1)}}\cdots f_{i,r_{\pi(t)}}f_{j,s}f_{i,r_{\pi(t+1)}}\cdots f_{i,r_{\pi(m)}}.
\end{equation*}
Then, we need to show $Z^-(i,j;\vec{r},s)=0$.
This is proved completely analogously to~(\ref{U7}).

\medskip
 This completes our proof of Theorem~\ref{Levendorskii-type}.

\begin{Rem}
(a) Specializing $\vv\mapsto v\in \BC^\times$ from the beginning
and viewing all algebras as $\BC$-algebras, the statement of
Theorem~\ref{Levendorskii-type} still holds as long as $v$ is not a root of unity.

\noindent
(b) A slightly different proof can be obtained by following  the
arguments in~\cite{le}.

\noindent
(c) We note that both Theorem~\ref{Levendorskii-type} and its proof
are valid also for all affine Lie algebras, except for the type $A_1^{(1)}$.
\end{Rem}


\subsection{An alternative presentation of $\CU^\ssc_{0,\mu}$ for $\mu\in\Lambda^-$}
\label{Lev-update}
$\ $

Inspired by the recent result~\cite[Theorem 2.13]{gnw}, we provide
another realization of $\CU^\ssc_{0,\mu}$ (with $\mu\in\Lambda^-$) without the
defining relation~(\ref{u9}). Following the notations of
Section~\ref{Levendorskii}, denote by $\wt{\CU}_{\mu_1,\mu_2}$ the
associative $\BC(\vv)$-algebra generated by
  $$\{e_{i,r},f_{i,s},(\psi^+_{i,0})^{\pm 1}, (\psi^-_{i,b_i})^{\pm 1}, h_{i,\pm 1}|
    i\in I, b_{2,i}-1\leq r\leq 1, b_{1,i}-1\leq s\leq 1\}$$
with the defining relations~(\ref{u1}--\ref{u8}).
Define inductively $e_{i,r}, f_{i,r}, \psi^\pm_{i,r}$ as it was done for
$\hat{\CU}_{\mu_1,\mu_2}$ right before Theorem~\ref{Levendorskii-type}.

\begin{Thm}\label{Levendorskii-type update}
There is a unique $\BC(\vv)$-algebra isomorphism
$\wt{\CU}_{\mu_1,\mu_2}\iso \CU^\ssc_{0,\mu}$, such that
  $$e_{i,r}\mapsto e_{i,r},\
    f_{i,r}\mapsto f_{i,r},\
    \psi^\pm_{i,\pm s^\pm_i}\mapsto \psi^\pm_{i,\pm s^\pm_i}\ \
    \mathrm{for}\ i\in I,r\in \BZ, s^+_i\geq 0, s^-_i\geq -b_i.$$
\end{Thm}

\begin{proof}
Due to Theorem~\ref{Levendorskii-type},
it suffices to show that~(\ref{u9}) can be derived from~(\ref{u1}--\ref{u8}).
We will treat only the first relation of~(\ref{u9})
(the second is completely analogous).

First, we note that relations~(\ref{v1}--\ref{v5}) and~(\ref{U2}, \ref{U3}, \ref{U6}) with $i\ne j$
hold in $\wt{\CU}_{\mu_1,\mu_2}$, since their proofs for the algebra
$\hat{\CU}_{\mu_1,\mu_2}$ were solely based on relations~(\ref{u1}--\ref{u6}).
Likewise, the equalities $Y^\pm(i,j;r,s)=0$ from our verifications of~(\ref{U4}, \ref{U5}) for $i\ne j$
still hold for $r\in \{0,1\}, s\in \BZ$.

Second, we have
\begin{equation}\tag{v18}\label{v18}
  [\psi^+_{i,2}, e_{i,0}]_{\vv_i^2}+[e_{i,1},\psi^+_{i,1}]_{\vv_i^2}=0,\
  [\psi^+_{i,2}, f_{i,0}]_{\vv_i^{-2}}+[f_{i,1},\psi^+_{i,1}]_{\vv_i^{-2}}=0.
\end{equation}
These equalities are proved completely analogously to~(\ref{v11})
from our proof of Lemma~\ref{Shifts by 2}, but now we start from the
equality $[e_{i,1},e_{i,0}]_{\vv_i^2}=0$ rather than
$[e_{i,0},e_{i,-1}]_{\vv_i^2}=0$ (commuting it first with $h_{i,1}$
and then further with $f_{i,0}$).

Recall $h_{i,2}$ of~(\ref{Def-n of h_2}).
Analogously to Lemma~\ref{Shifts by 2}, we see that~(\ref{v18})
implies\footnote{Note that we can not deduce the statement of
Lemma~\ref{Shifts by 2} due to the absence of~(\ref{u9}).}
\begin{equation}\tag{v19}\label{v19}
  [h_{i,2}, e_{i,0}]=\frac{[4]_{\vv_i}}{2}\cdot e_{i,2},\
  [h_{i,2}, f_{i,0}]=-\frac{[4]_{\vv_i}}{2}\cdot f_{i,2}.
\end{equation}
Likewise, the aforementioned equalities $Y^\pm(i,j;1,s)=0$ for $i\ne j, s\in \BZ$, also imply
\begin{equation}\tag{v20}\label{v20}
  [h_{i,2}, e_{j,s}]=\frac{[2c_{ij}]_{\vv_i}}{2}\cdot e_{j,s+2},\
  [h_{i,2}, f_{j,s}]=-\frac{[2c_{ij}]_{\vv_i}}{2}\cdot f_{j,s+2}\
  \mathrm{for}\ i\ne j, s\in \BZ.
\end{equation}
Finally, due to~(\ref{u7}, \ref{u8}, \ref{v1}, \ref{v2}, \ref{v19}, \ref{v20}), we also get
$[e_{i,r},e_{j,s}]=[f_{i,r},f_{j,s}]=0$ if $c_{ij}=0$ and
$Z^\pm(i,j;r',0,s)=Z^\pm(i,j;1,1,s)=0$ if $c_{ij}=-1$ for $r,s\in \BZ, r'\in\{0,1,2\}$.

\medskip
In the simply-laced case, the rest of the proof follows from the next result.

\begin{Lem}\label{Derivation 1}
Let $i,j\in I$ be such that $c_{ij}=-1$.
Then $[\psi^+_{i,1},\psi^+_{i,2}]=0$.
\end{Lem}

\begin{proof}
As just proved, we have
  $[f_{i,1},[f_{i,1},f_{j,0}]_{\vv_i^{-1}}]_{\vv_i}=0$.
Commuting this equality with $e_{j,1}$ and applying~(\ref{v4})
together with~(\ref{U6}) for $i\ne j$, we get
  $[f_{i,1},[f_{i,1},\psi^+_{j,1}]_{\vv_i^{-1}}]_{\vv_i}=0.$
Combining the latter equality with
 $\psi^+_{j,1}=(\vv_j-\vv_j^{-1})\psi^+_{j,0}h_{j,1}=(\vv_i-\vv_i^{-1})\psi^+_{j,0}h_{j,1}$
and using~(\ref{v2}), we find
\begin{equation*}
    [f_{i,1},[f_{i,1},h_{j,1}]]_{\vv_i^2}=0\Longrightarrow
    [f_{i,1},f_{i,2}]_{\vv_i^2}=0\Longrightarrow
    [f_{i,2},f_{i,1}]_{\vv_i^{-2}}=0.
\end{equation*}
Commuting this further with $e_{i,0}$, we obtain
\begin{equation*}
  [\psi^+_{i,2},f_{i,1}]_{\vv_i^{-2}}+[f_{i,2},\psi^+_{i,1}]_{\vv_i^{-2}}=0.
\end{equation*}
Finally, we apply $[e_{i,0},-]_{\vv_i^{-2}}$ to the latter equality.
In the left-hand side we get two summands computed below.

\textup{(1)} We have
  $[e_{i,0},[f_{i,2},\psi^+_{i,1}]_{\vv_i^{-2}}]_{\vv_i^{-2}}=
   [[e_{i,0},f_{i,2}],\psi^+_{i,1}]_{\vv_i^{-4}}+[f_{i,2}, [e_{i,0},\psi^+_{i,1}]_{\vv_i^{-2}}]_{\vv_i^{-2}}.$
Due to~(\ref{u4}),
  $[e_{i,0},\psi^+_{i,1}]_{\vv_i^{-2}}=(\vv_i^{-2}-\vv_i^2)e_{i,1}\psi^+_{i,0}\Rightarrow
   [f_{i,2}, [e_{i,0},\psi^+_{i,1}]_{\vv_i^{-2}}]_{\vv_i^{-2}}=
   (\vv_i^{-2}-\vv_i^2)[f_{i,2},e_{i,1}]_{\vv_i^{-4}}\psi^+_{i,0}$.
Combining this with~(\ref{v5}), we thus get
\begin{equation}\tag{v21}\label{v21}
  [e_{i,0},[f_{i,2},\psi^+_{i,1}]_{\vv_i^{-2}}]_{\vv_i^{-2}}=
  [\psi^+_{i,2},\psi^+_{i,1}]_{\vv_i^{-4}}/(\vv_i-\vv_i^{-1})+
  (\vv_i^{-2}-\vv_i^2)[f_{i,2},e_{i,1}]_{\vv_i^{-4}}\psi^+_{i,0}.
\end{equation}

\textup{(2)} We have
  $[e_{i,0}, [\psi^+_{i,2},f_{i,1}]_{\vv_i^{-2}}]_{\vv_i^{-2}}=
   [[e_{i,0},\psi^+_{i,2}]_{\vv_i^{-2}},f_{i,1}]_{\vv_i^{-2}}+\vv_i^{-2}[\psi^+_{i,2},[e_{i,0},f_{i,1}]].$
By (\ref{v18}):
  $[e_{i,0},\psi^+_{i,2}]_{\vv_i^{-2}}=-\vv_i^{-2}[\psi^+_{i,2},e_{i,0}]_{\vv_i^2}=
    \vv_i^{-2}[e_{i,1},\psi^+_{i,1}]_{\vv_i^2}=\vv_i^{-2}(\vv_i-\vv_i^{-1})[e_{i,1},h_{i,1}]_{\vv_i^4}\psi^+_{i,0}$.
Hence,
\begin{equation*}
\begin{split}
  & [[e_{i,0},\psi^+_{i,2}]_{\vv_i^{-2}},f_{i,1}]_{\vv_i^{-2}}=
    \vv_i^{-4}(\vv_i-\vv_i^{-1})[e_{i,1}h_{i,1}-\vv_i^4h_{i,1}e_{i,1},f_{i,1}]\psi^+_{i,0}=\\
  & \vv_i^{-4}(\vv_i-\vv_i^{-1})([\psi^+_{i,2},h_{i,1}]_{\vv_i^4}/(\vv_i-\vv_i^{-1})-(\vv_i+\vv_i^{-1})[e_{i,1},f_{i,2}]_{\vv_i^4})\psi^+_{i,0}.
\end{split}
\end{equation*}
Therefore,
\begin{equation}\tag{v22}\label{v22}
  [e_{i,0}, [\psi^+_{i,2},f_{i,1}]_{\vv_i^{-2}}]_{\vv_i^{-2}}=
  \frac{[\psi^+_{i,2},\psi^+_{i,1}]}{\vv_i^2(\vv_i-\vv_i^{-1})}+
  (\vv_i^2-\vv_i^{-2})[f_{i,2},e_{i,1}]_{\vv_i^{-4}}\psi^+_{i,0}+
   \frac{[\psi^+_{i,2},\psi^+_{i,1}]_{\vv_i^4}}{\vv_i^4(\vv_i-\vv_i^{-1})}.
\end{equation}

Substituting~(\ref{v22}) and~(\ref{v21}) into
  $[e_{i,0},[\psi^+_{i,2},f_{i,1}]_{\vv_i^{-2}}+[f_{i,2},\psi^+_{i,1}]_{\vv_i^{-2}}]_{\vv_i^{-2}}=0$,
we find
\begin{equation*}
    [\psi^+_{i,2},\psi^+_{i,1}]_{\vv_i^{-4}}+
    \vv_i^{-2}[\psi^+_{i,2},\psi^+_{i,1}]+
    \vv_i^{-4}[\psi^+_{i,2},\psi^+_{i,1}]_{\vv_i^4}=0.
\end{equation*}
The left-hand side of this equality equals
  $\frac{1-\vv_i^{-6}}{1-\vv_i^{-2}}\cdot [\psi^+_{i,2},\psi^+_{i,1}]$.
Hence, $[\psi^+_{i,1},\psi^+_{i,2}]=0$.
\end{proof}

Our next result completes the proof for non-simply-laced $\fg$.

\begin{Lem}\label{Derivation 2}
If $c_{ij}\ne 0$ and $[\psi^+_{i,1},\psi^+_{i,2}]=0$, then $[\psi^+_{j,1},\psi^+_{j,2}]=0$.
\end{Lem}

\begin{proof}
Due to~(\ref{v1}, \ref{v2}):
  $[h_{i,1},e_{i,r}]=\frac{[2]_{\vv_i}}{[c_{ji}]_{\vv_j}}\cdot [h_{j,1},e_{i,r}],\
   [h_{i,1},f_{i,r}]=\frac{[2]_{\vv_i}}{[c_{ji}]_{\vv_j}}\cdot [h_{j,1},f_{i,r}].$
Hence
  $[h_{i,1},\psi^+_{i,2}]=
   (\vv_i-\vv_i^{-1})([[h_{i,1},e_{i,1}],f_{i,1}]+[e_{i,1},[h_{i,1},f_{i,1}]])=
   [2]_{\vv_i}/[c_{ji}]_{\vv_j}\cdot [h_{j,1},\psi^+_{i,2}].$
Therefore,
  $[\psi^+_{i,1},\psi^+_{i,2}]=0\Rightarrow
   [h_{j,1},\psi^+_{i,2}]=0\Rightarrow
   [h_{j,1},h_{i,2}]=0$
with the second implication due to~(\ref{u1}).
Commuting the latter equality with $f_{j,0}$, we get
\begin{equation*}
  0=[f_{j,0},[h_{j,1},h_{i,2}]]=
    [c_{jj}]_{\vv_j}\cdot [f_{j,1},h_{i,2}]+\frac{[2c_{ij}]_{\vv_i}}{2}\cdot [h_{j,1},f_{j,2}].
\end{equation*}
Commuting this further with $e_{j,0}$, we obtain
\begin{equation}\tag{v23}\label{v23}
  [c_{jj}]_{\vv_j}\cdot [e_{j,0},[f_{j,1},h_{i,2}]]+
  \frac{[2c_{ij}]_{\vv_i}}{2}\cdot [e_{j,0},[h_{j,1},f_{j,2}]]=0.
\end{equation}
Note that
\begin{equation*}
\begin{split}
  & [e_{j,0},[f_{j,1},h_{i,2}]]=[\psi^+_{j,1},h_{i,2}]/(\vv_j-\vv_j^{-1})-
    \frac{[2c_{ij}]_{\vv_i}}{2}\cdot [f_{j,1}, e_{j,2}]=
    -\frac{[2c_{ij}]_{\vv_i}}{2}\cdot [f_{j,1}, e_{j,2}],\\
  & [e_{j,0},[h_{j,1},f_{j,2}]]=-[c_{jj}]_{\vv_j}\cdot [e_{j,1},f_{j,2}]+[h_{j,1},\psi^+_{j,2}]/(\vv_j-\vv_j^{-1}),\\
  & [e_{j,2},f_{j,1}]-[e_{j,1},f_{j,2}]=[c_{jj}]_{\vv_j}^{-1}\cdot [h_{j,1},[e_{j,1},f_{j,1}]]=
    [c_{jj}]_{\vv_j}^{-1}\cdot [h_{j,1},\psi^+_{j,2}]/(\vv_j-\vv_j^{-1}).
\end{split}
\end{equation*}
Substituting the last three equalities into~(\ref{v23}), we get
  $\frac{[2c_{ij}]_{\vv_i}}{\vv_j-\vv_j^{-1}}\cdot [h_{j,1},\psi^+_{j,2}]=0$.
Thus,
  $[h_{j,1},\psi^+_{j,2}]=0\Rightarrow [\psi^+_{j,1},\psi^+_{j,2}]=0$.
\end{proof}
 This completes our proof of Theorem~\ref{Levendorskii-type update}.
\end{proof}


\section{Proof of Theorem~\ref{Relations bw ABCD}}\label{Proof of Relations bw ABCD}

The proof of part (a) proceeds in two steps.
First, we consider the simplest case $\fg=\ssl_2$.
Then, we show how a general case can be easily reduced
to the case of $\ssl_2$.


\subsection{Proof of Theorem~\ref{Relations bw ABCD}(a) for $\fg=\ssl_2$}
\

First, let us derive an explicit formula for $A^\pm(z)$. Recall the
elements $\{h_{\pm r}\}_{r=1}^\infty$ of~Section~\ref{sqaa}, such that
  $z^{\mp b^\pm}(\psi^\pm_{\mp b^\pm})^{-1}\psi^\pm(z)=
   \exp\left(\pm(\vv-\vv^{-1})\sum_{r>0}h_{\pm r}z^{\mp r}\right).$
For $r\ne 0$, define $t_r:=-h_r/(1+\vv^{2r})$, and set
\begin{equation}\label{A for sl2}
  A^\pm(z):=(\phi^\pm)^{-1}\cdot \exp\left(\pm (\vv-\vv^{-1})\sum_{r>0}t_{\pm r}z^{\mp r}\right).
\end{equation}
Then,
  $z^{\mp b^\pm}\psi^\pm(z)=\frac{1}{A^\pm(z)A^\pm(\vv^{-2}z)}$
and $A^\pm(z)$ is the unique solution with $A^\pm_0:=(\phi^\pm)^{-1}$.

Relations~(\ref{ABCD0}) and~(\ref{ABCD1}) follow immediately
from~(\ref{U10}) and (\ref{U1}), respectively, while the verification
of~(\ref{ABCD5}--\ref{ABCD12}) is based on the following result.

\begin{Lem}\label{Properties of half-currents}
For any $\epsilon,\epsilon'\in \{\pm\}$, we have:

\noindent
(a1) $(\vv z-\vv^{-1}w)A^{\epsilon}(z)e(w)=(z-w)e(w)A^{\epsilon}(z).$

\noindent
(a2) $(\vv z-\vv^{-1}w)A^{\epsilon}(z)e^{\epsilon'}(w)-(z-w)e^{\epsilon'}(w)A^{\epsilon}(z)=
      (\vv-\vv^{-1})wA^{\epsilon}(z)e^{\epsilon}(z).$

\noindent
(a3) $(\vv z-\vv^{-1}w)A^{\epsilon}(z)e^{\epsilon'}(w)-(z-w)e^{\epsilon'}(w)A^{\epsilon}(z)=
      (1-\vv^{-2})we^{\epsilon}(\vv^2z)A^{\epsilon}(z).$

\noindent
(b1) $(z-w)A^{\epsilon}(z)f(w)=(\vv z-\vv^{-1}w)f(w)A^{\epsilon}(z).$

\noindent
(b2) $(z-w)A^{\epsilon}(z)f^{\epsilon'}(w)-(\vv z-\vv^{-1}w)f^{\epsilon'}(w)A^{\epsilon}(z)=
      (\vv^{-1}-\vv)zf^{\epsilon}(z)A^{\epsilon}(z).$

\noindent
(b3) $(z-w)A^{\epsilon}(z)f^{\epsilon'}(w)-(\vv z-\vv^{-1}w)f^{\epsilon'}(w)A^{\epsilon}(z)=
      (1-\vv^2)zA^{\epsilon}(z)f^{\epsilon}(\vv^2z).$

\noindent
(c) $(z-w)[e^{\epsilon}(z),f^{\epsilon'}(w)]=z(\psi^{\epsilon'}(w)-\psi^{\epsilon}(z))/(\vv-\vv^{-1}).$

\noindent
(d1) $(z-\vv^2w)e^{\epsilon}(z)e^{\epsilon'}(w)-(\vv^2z-w)e^{\epsilon'}(w)e^{\epsilon}(z)=
      z[e_0,e^{\epsilon'}(w)]_{\vv^2}+w[e_0,e^{\epsilon}(z)]_{\vv^2}.$

\noindent
(d2) $(z-\vv^2w)e^{\epsilon}(z)e^{\epsilon'}(w)-(\vv^2z-w)e^{\epsilon'}(w)e^{\epsilon}(z)=
      (1-\vv^2)(we^{\epsilon}(z)^2+ze^{\epsilon'}(w)^2).$

\noindent
(e1) $(\vv^2z-w)f^{\epsilon}(z)f^{\epsilon'}(w)-(z-\vv^2w)f^{\epsilon'}(w)f^{\epsilon}(z)=
      \vv^2[f_1,f^{\epsilon'}(w)]_{\vv^{-2}}+\vv^2[f_1,f^{\epsilon}(z)]_{\vv^{-2}}.$

\noindent
(e2) $(\vv^2z-w)f^{\epsilon}(z)f^{\epsilon'}(w)-(z-\vv^2w)f^{\epsilon'}(w)f^{\epsilon}(z)=
      (\vv^2-1)(zf^{\epsilon}(z)^2+wf^{\epsilon'}(w)^2).$

\noindent
(f1) $(z-\vv^2w)\psi^{\epsilon}(z)e^{\epsilon'}(w)-(\vv^2z-w)e^{\epsilon'}(w)\psi^{\epsilon}(z)=
      (\vv^{-2}-\vv^2)w\psi^{\epsilon}(z)e^{\epsilon}(\vv^2z).$

\noindent
(f2) $(z-\vv^2w)\psi^{\epsilon}(z)e^{\epsilon'}(w)-(\vv^2z-w)e^{\epsilon'}(w)\psi^{\epsilon}(z)=
      (1-\vv^4)we^{\epsilon}(\vv^{-2}z)\psi^{\epsilon}(z).$

\noindent
(g1) $(\vv^2z-w)\psi^{\epsilon}(z)f^{\epsilon'}(w)-(z-\vv^2w)f^{\epsilon'}(w)\psi^{\epsilon}(z)=
      (\vv^2-\vv^{-2})z\psi^{\epsilon}(z)f^{\epsilon}(\vv^{-2}z).$

\noindent
(g2) $(\vv^2z-w)\psi^{\epsilon}(z)f^{\epsilon'}(w)-(z-\vv^2w)f^{\epsilon'}(w)\psi^{\epsilon}(z)=
      (\vv^4-1)zf^{\epsilon}(\vv^2z)\psi^{\epsilon}(z).$
\end{Lem}

\begin{proof}
(a1) According to~(\ref{U4'}), we have
$[t_r,e_s]=\frac{\vv^{-2r}-1}{r(\vv-\vv^{-1})}e_{s+r}$
for $r\ne 0, s\in \BZ$. Combining this with~(\ref{A for sl2}), we find
  $A^\pm(z)e(w)=e(w)A^\pm(z)\vv^{\mp 1}\exp\left(\sum_{r>0}\frac{\vv^{\mp 2r}-1}{r}(w/z)^{\pm r}\right).$
The latter exponent equals $\frac{z-w}{z-\vv^{-2}w}$ (in the ``$+$''
case) or $\frac{z-w}{\vv^2z-w}$ (in the ``$-$'' case), hence, (a1).

(a2, a3) First, we consider the case $\epsilon=\epsilon'=+$.
Due to (a1), we have
  $\vv A^+_{r+1}e_s-\vv^{-1}A^+_re_{s+1}=e_sA^+_{r+1}-e_{s+1}A^+_r$
for any $r\in \BN, s\in \BZ$. Multiplying this equality by
$z^{-r}w^{-s-1}$ and summing over all $r,s\in \BN$, we find
  $w^{-1}((\vv z-\vv^{-1}w)A^+(z)e^+(w)-(z-w)e^+(w)A^+(z))=[e_0,A^+(z)]_{\vv^{-1}}$.
Note that the right-hand side is independent of $w$. Substituting
either $w=z$ or $w=\vv^2z$ into the left-hand side, we get the
equalities (a2) and (a3) for $\epsilon=\epsilon'=+$, respectively.

Next, we consider the case $\epsilon=\epsilon'=-$. Due to (a1), we have
  $\vv A^-_{-r+1}e_{-s}-\vv^{-1}A^-_re_{-s+1}=e_{-s}A^-_{-r+1}-e_{-s+1}A^-_{-r}$
for any $r\in \BN, s\in \BZ$, where we set $A^-_1:=0$. Multiplying this equality
by $-z^rw^{s-1}$ and summing over all $r\in \BN,s\in\BZ_{>0}$, we find
  $w^{-1}((\vv z-\vv^{-1}w)A^-(z)e^-(w)-(z-w)e^-(w)A^-(z))=[e_0,A^-(z)]_{\vv^{-1}}$.
Note that the right-hand side is independent of $w$. Substituting
either $w=z$ or $w=\vv^2z$ into the left-hand side, we get the
equalities (a2) and (a3) for $\epsilon=\epsilon'=-$, respectively.

The case $\epsilon'\ne\epsilon$ follows by combining the formula
$e^{\epsilon'}(w)=e^\epsilon(w)+\epsilon'e(w)$ with part (a1) and
the cases $\epsilon=\epsilon'$ of parts (a2, a3), established above.

(b1--b3) Parts (b1, b2, b3) are proved completely analogously to (a1, a2, a3), respectively.

(c) First, we consider the case $\epsilon=\epsilon'$. According
to~(\ref{U6}), we have $[e_r,f_s]=\frac{\psi^+_{r+s}}{\vv-\vv^{-1}}$ for
$r\geq 0,s>0$. For $N>0$, we have
$(z-w)\sum_{s=1}^Nw^{-s}z^{s-N}=z(w^{-N}-z^{-N})$. Hence,
  $(z-w)[e^+(z),f^+(w)]=\sum_{N>0}z(w^{-N}-z^{-N})\frac{\psi^+_N}{\vv-\vv^{-1}}=
   z\frac{\psi^+(w)-\psi^+(z)}{\vv-\vv^{-1}}$.
Likewise, we have $[e_{-r},f_{-s}]=-\frac{\psi^-_{-r-s}}{\vv-\vv^{-1}}$
for $r>0,s\geq 0$. For $N>0$, we have
$(z-w)\sum_{s=1}^Nz^sw^{N-s}=z(z^N-w^N)$. Hence,
  $(z-w)[e^-(z),f^-(w)]=-\sum_{N>0}z(z^N-w^N)\frac{\psi^-_{-N}}{\vv-\vv^{-1}}=z\frac{\psi^-(w)-\psi^-(z)}{\vv-\vv^{-1}}$.

Next, we consider the case $\epsilon\ne \epsilon'$. According to~(\ref{U6}), we have
  $[e(z),f(w)]=\frac{\delta(z/w)}{\vv-\vv^{-1}}(\psi^+(z)-\psi^-(z))=
   \frac{\delta(z/w)}{\vv-\vv^{-1}}(\psi^+(w)-\psi^-(w))$.
Taking the terms with negative powers of $w$, we find
  $[e(z),f^+(w)]=\frac{z/w}{1-z/w}\frac{\psi^+(z)-\psi^-(z)}{\vv-\vv^{-1}}\Rightarrow
   (z-w)[e(z),f^+(w)]=z\frac{\psi^-(z)-\psi^+(z)}{\vv-\vv^{-1}}$,
while taking the terms with nonpositive powers of $z$, we find
  $[e^+(z),f(w)]=\frac{1}{1-w/z}\frac{\psi^+(w)-\psi^-(w)}{\vv-\vv^{-1}}\Rightarrow
   (z-w)[e^+(z),f(w)]=z\frac{\psi^+(w)-\psi^-(w)}{\vv-\vv^{-1}}$.
Combining these equalities with
  $(z-w)[e^+(z),f^+(w)]=z\frac{\psi^+(w)-\psi^+(z)}{\vv-\vv^{-1}}$
from above and $e^-(z)=e^+(z)-e(z), f^-(z)=f^+(z)-f(z)$, we obtain
the $\epsilon\ne \epsilon'$ cases of part (c).

(d1) Comparing the coefficients of $z^{-\epsilon r}w^{-\epsilon's}$
in both sides of relation~(\ref{U2}), we find
  $e_{\epsilon r+1}e_{\epsilon' s}-\vv^2e_{\epsilon r}e_{\epsilon' s+1}=
   \vv^2e_{\epsilon' s}e_{\epsilon r+1}-e_{\epsilon' s+1}e_{\epsilon r}$
for any $r,s\in \BZ$. Multiplying this equality by
$\epsilon\epsilon'\cdot z^{-\epsilon r}w^{-\epsilon' s}$ and summing over
$r\geq \delta_{\epsilon,-},\ s\geq \delta_{\epsilon',-}$, we get (d1).

(d2) Substituting $w=z$ into the $\epsilon=\epsilon'$ case of (d1),
we find
  $[e_0,e^\pm(z)]_{\vv^2}=(1-\vv^2)e^\pm(z)^2$.
Replacing accordingly the right-hand side of (d1), we obtain (d2).

(e1, e2) Parts (e1, e2) are proved completely analogously to (d1, d2), respectively.

(f1, f2) Parts (f1, f2) are deduced from relation~(\ref{U4}) in
the same way as we deduced parts (a2, a3) from (a1).

(g1, g2) Parts (g1, g2) are proved completely analogously to (f1, f2), respectively.
\end{proof}

Now let us verify relations~(\ref{ABCD5}--\ref{ABCD12})
using~Lemma~\ref{Properties of half-currents}. The idea is first to
use parts (a3, b2) of~Lemma~\ref{Properties of half-currents}
(resp.\ parts (a2, b3)) to move
all the series $A^\bullet(\cdot)$ to the right (resp.\ to the left), and then to
use Lemma~\ref{Properties of half-currents}(c--g2) to simplify the
remaining part. Since $\fg=\ssl_2$ we will drop the
index $i$ from our notation.


\subsubsection{Verification of the first relation in~(\ref{ABCD5})}
\

We need to prove $[B^\epsilon(z),B^{\epsilon'}(w)]=0$, or
equivalently, $(z-w)[B^\epsilon(z),B^{\epsilon'}(w)]=0$. By
definition,
  $B^\epsilon(z)B^{\epsilon'}(w)=
   (\vv-\vv^{-1})^2A^\epsilon(z)e^\epsilon(z)A^{\epsilon'}(w)e^{\epsilon'}(w)$.
Applying Lemma~\ref{Properties of half-currents}(a2), we see that
\begin{equation*}
    (z-w)B^\epsilon(z)B^{\epsilon'}(w)=
    (\vv-\vv^{-1})^2 A^\epsilon(z)A^{\epsilon'}(w)
    ((\vv^{-1}z-\vv w)e^\epsilon(z)e^{\epsilon'}(w)+(\vv-\vv^{-1})ze^{\epsilon'}(w)^2).
\end{equation*}
Hence, the equality $(z-w)[B^\epsilon(z),B^{\epsilon'}(w)]=0$ boils
down to the vanishing of
\begin{equation*}
    (\vv^{-1}z-\vv w)e^\epsilon(z)e^{\epsilon'}(w)+(\vv-\vv^{-1})ze^{\epsilon'}(w)^2+
    (\vv^{-1}w-\vv z)e^{\epsilon'}(w)e^\epsilon(z)+(\vv-\vv^{-1})we^\epsilon(z)^2,
\end{equation*}
which is exactly the statement of Lemma~\ref{Properties of
half-currents}(d2).


\subsubsection{Verification of the second relation in~(\ref{ABCD5})}
\

We need to prove $[C^\epsilon(z),C^{\epsilon'}(w)]=0$, or
equivalently, $(z-w)[C^\epsilon(z),C^{\epsilon'}(w)]=0$.
By definition,
  $C^\epsilon(z)C^{\epsilon'}(w)=
   (\vv-\vv^{-1})^2f^\epsilon(z)A^\epsilon(z)f^{\epsilon'}(w)A^{\epsilon'}(w)$.
Applying Lemma~\ref{Properties of half-currents}(b2), we see that
\begin{equation*}
    (z-w)C^\epsilon(z)C^{\epsilon'}(w)=(\vv-\vv^{-1})^2
    ((\vv z-\vv^{-1} w)f^\epsilon(z)f^{\epsilon'}(w)+(\vv^{-1}-\vv)zf^{\epsilon}(z)^2)A^\epsilon(z)A^{\epsilon'}(w).
\end{equation*}
Hence, the equality $(z-w)[C^\epsilon(z),C^{\epsilon'}(w)]=0$ boils
down to the vanishing of
\begin{equation*}
    (\vv z-\vv^{-1} w)f^\epsilon(z)f^{\epsilon'}(w)+(\vv^{-1}-\vv)zf^{\epsilon}(z)^2+
    (\vv w-\vv^{-1} z)f^{\epsilon'}(w)f^\epsilon(z)+(\vv^{-1}-\vv)wf^{\epsilon'}(w)^2,
\end{equation*}
which is exactly the statement of Lemma~\ref{Properties of
half-currents}(e2).


\subsubsection{Verification of the third relation in~(\ref{ABCD5})}
\

The verification of the equality
  $[D^\epsilon(z),D^{\epsilon'}(w)]=0$
is much more cumbersome and is left to the interested reader.


\subsubsection{Verification of~(\ref{ABCD6})}
\

We need to prove
  $(z-w)[B^{\epsilon'}(w),A^\epsilon(z)]_{\vv^{-1}}=
   (\vv-\vv^{-1})(zA^\epsilon(z)B^{\epsilon'}(w)-wA^{\epsilon'}(w)B^{\epsilon}(z))$.
By definition and~(\ref{ABCD1}), the RHS equals
  $(\vv-\vv^{-1})^2A^\epsilon(z)A^{\epsilon'}(w)(ze^{\epsilon'}(w)-we^\epsilon(z))$.
Meanwhile, the LHS equals
  $(\vv-\vv^{-1})(z-w)
   (A^{\epsilon'}(w)e^{\epsilon'}(w)A^\epsilon(z)-\vv^{-1}A^\epsilon(z)A^{\epsilon'}(w)e^{\epsilon'}(w)).$
We use Lemma~\ref{Properties of half-currents}(a2) to replace the
first term, so that the LHS equals
\begin{equation*}
    (\vv-\vv^{-1})A^\epsilon(z)A^{\epsilon'}(w)
    \left((\vv z-\vv^{-1}w)e^{\epsilon'}(w)-(\vv-\vv^{-1})we^\epsilon(z)-\vv^{-1}(z-w)e^{\epsilon'}(w)\right),
\end{equation*}
which exactly coincides with the above formula for the RHS.


\subsubsection{Verification of~(\ref{ABCD7})}
\

We need to prove
  $(z-w)[A^\epsilon(z),C^{\epsilon'}(w)]_\vv=
   (\vv-\vv^{-1})(wC^{\epsilon'}(w)A^\epsilon(z)-zC^{\epsilon}(z)A^{\epsilon'}(w))$.
By definition and~(\ref{ABCD1}), the RHS equals
  $(\vv-\vv^{-1})^2(wf^{\epsilon'}(w)-zf^\epsilon(z))A^\epsilon(z)A^{\epsilon'}(w)$.
Meanwhile, the LHS equals
  $(\vv-\vv^{-1})(z-w)
   (A^\epsilon(z)f^{\epsilon'}(w)A^{\epsilon'}(w)-\vv f^{\epsilon'}(w)A^{\epsilon'}(w)A^\epsilon(z))$.
We use Lemma~\ref{Properties of half-currents}(b2) to replace the
first term, so that the LHS equals
\begin{equation*}
    (\vv-\vv^{-1})\left((\vv z-\vv^{-1}w)f^{\epsilon'}(w)+(\vv^{-1}-\vv)zf^{\epsilon}(z)-\vv(z-w)f^{\epsilon'}(w)\right)
    A^\epsilon(z)A^{\epsilon'}(w),
\end{equation*}
which exactly coincides with the above formula for the RHS.


\subsubsection{Verification of~(\ref{ABCD8})}
\

We need to prove
  $(z-w)[B^\epsilon(z),C^{\epsilon'}(w)]=
   (\vv-\vv^{-1})z(D^{\epsilon'}(w)A^\epsilon(z)-D^\epsilon(z)A^{\epsilon'}(w))$.
Applying the equality
  $A^\epsilon(z)e^\epsilon(z)=\vv^{-1}e^\epsilon(\vv^2z)A^\epsilon(z)$,
which follows from Lemma~\ref{Properties of half-currents}(a2), we
see that the LHS equals
\begin{equation*}
    \vv^{-1}(\vv-\vv^{-1})^2(z-w)
    \left(e^\epsilon(\vv^2z)A^\epsilon(z)f^{\epsilon'}(w)A^{\epsilon'}(w)-
          f^{\epsilon'}(w)A^{\epsilon'}(w)e^\epsilon(\vv^2z)A^\epsilon(z)\right).
\end{equation*}
Applying~Lemma~\ref{Properties of half-currents}(a3, b2) to move both
$A^\epsilon(z),A^{\epsilon'}(w)$ to the right and simplifying the resulting
expression, we find that the LHS equals
\begin{equation*}
\begin{split}
  & \vv^{-1}(\vv-\vv^{-1})^2\left((\vv z-\vv^{-1}w)[e^\epsilon(\vv^2z),f^{\epsilon'}(w)]+\right.\\
  & \left.(\vv-\vv^{-1})(zf^{\epsilon'}(w)e^{\epsilon'}(\vv^2 w)-(\vv-\vv^{-1})ze^\epsilon(\vv^2z)f^\epsilon(z))\right)
    A^\epsilon(z)A^{\epsilon'}(w).
\end{split}
\end{equation*}
Meanwhile,
  $D^\epsilon(z)=\psi^\epsilon(z)A^\epsilon(z)+
   \vv^{-1}(\vv-\vv^{-1})^2f^\epsilon(z)e^\epsilon(\vv^2z)A^\epsilon(z)$,
so that the RHS equals
\begin{equation*}
    (\vv-\vv^{-1})\left(z(\psi^{\epsilon'}(w)-\psi^\epsilon(z))+
     \vv^{-1}(\vv-\vv^{-1})^2z(f^{\epsilon'}(w)e^{\epsilon'}(\vv^2w)-
     f^\epsilon(z)e^\epsilon(\vv^2z))\right)A^\epsilon(z)A^{\epsilon'}(w).
\end{equation*}
Thus, the equality LHS$=$RHS boils down to proving
\begin{equation*}
    \vv^{-2}(\vv^2z-w)[e^{\epsilon}(\vv^2z),f^{\epsilon'}(w)]-(1-\vv^{-2})z[e^\epsilon(\vv^2z),f^\epsilon(z)]=
    \frac{z}{\vv-\vv^{-1}}(\psi^{\epsilon'}(w)-\psi^\epsilon(z)),
\end{equation*}
which immediately follows by applying Lemma~\ref{Properties of
half-currents}(c) to both terms on the left.


\subsubsection{Verification of~(\ref{ABCD9})}
\

We need to prove
  $(z-w)[B^\epsilon(z),D^{\epsilon'}(w)]_\vv=
   (\vv-\vv^{-1})(wD^{\epsilon'}(w)B^\epsilon(z)-zD^\epsilon(z)B^{\epsilon'}(w))$.
Combining the aforementioned equality
$A^\epsilon(z)e^\epsilon(z)=\vv^{-1}e^\epsilon(\vv^2z)A^\epsilon(z)$
with~Lemma~\ref{Properties of half-currents}(a3), we find that
$(w-z)\cdot \mathrm{RHS}$ equals
\begin{equation*}
\begin{split}
  & \left(\vv^{-1}(\vv-\vv^{-1})^2
    (w(\vv^{-1}w-\vv z)\psi^{\epsilon'}(w)e^{\epsilon}(\vv^2z)+
     z(\vv^{-1}z-\vv w)\psi^\epsilon(z)e^{\epsilon'}(\vv^2w))+\right.\\
  & \vv^{-1}(\vv-\vv^{-1})^3zw (\psi^{\epsilon'}(w)e^{\epsilon'}(\vv^2 w)+
    \psi^\epsilon(z)e^\epsilon(\vv^2z))+\\
  & \vv^{-2}(\vv-\vv^{-1})^4(w(\vv^{-1}w-\vv z)f^{\epsilon'}(w)e^{\epsilon'}(\vv^2w)e^\epsilon(\vv^2z)+
    z(\vv^{-1}z-\vv w)f^\epsilon(z)e^\epsilon(\vv^2z)e^{\epsilon'}(\vv^2w))+\\
  & \left.\vv^{-2}(\vv-\vv^{-1})^5zw (f^{\epsilon'}(w)e^{\epsilon'}(\vv^2 w)^2+f^\epsilon(z)e^\epsilon(\vv^2z)^2)\right)
    A^\epsilon(z)A^{\epsilon'}(w).
\end{split}
\end{equation*}
Meanwhile, using Lemma~\ref{Properties of half-currents}(a3, b2) to
move $A^\epsilon(z)$ to the right of
  $f^{\epsilon'}(w)e^{\epsilon'}(\vv^2 w)$,
we find that $(w-z)\cdot \mathrm{LHS}$ equals
\begin{equation*}
\begin{split}
  & \vv^{-1}(\vv-\vv^{-1})(w-z)\cdot
    \left((w-\vv^2z)\psi^{\epsilon'}(w)e^\epsilon(\vv^2z)+(\vv^2-1)z\psi^{\epsilon'}(w)e^{\epsilon'}(\vv^2w)+\right.\\
  & \left.(z-w)e^\epsilon(\vv^2z)\psi^{\epsilon'}(w))\right)A^\epsilon(z)A^{\epsilon'}(w)+\\
  & \vv^{-2}(\vv-\vv^{-1})^3\cdot \left(
    (w-\vv^2z)(w-z)f^{\epsilon'}(w)e^{\epsilon'}(\vv^2w)e^\epsilon(\vv^2z)+(\vv^2-1)z(w-z)f^{\epsilon'}(w)e^{\epsilon'}(\vv^2w)^2-\right.\\
  & (\vv z-\vv^{-1}w)(\vv^{-1}z-\vv w)e^\epsilon(\vv^2z)f^{\epsilon'}(w)e^{\epsilon'}(\vv^2w)-
    (\vv^{-1}-\vv)z(\vv^{-1}z-\vv w)e^\epsilon(\vv^2z)f^\epsilon(z)e^{\epsilon'}(\vv^2w)-\\
  & \left.(\vv z- \vv^{-1}w)(\vv-\vv^{-1})w e^\epsilon(\vv^2 z)f^{\epsilon'}(w)e^\epsilon(\vv^2z)-
    (\vv^{-1}-\vv)(\vv-\vv^{-1})zwe^{\epsilon}(\vv^2z)f^\epsilon(z)e^\epsilon(\vv^2z)\right)\times\\
  & A^\epsilon(z)A^{\epsilon'}(w).
\end{split}
\end{equation*}
To check that the above two big expressions coincide, we first
reorder some of the terms. We use Lemma~\ref{Properties of
half-currents}(f1) to move $\psi^{\epsilon'}(w)$ to the left of
$e^\epsilon(\vv^2z)$ via
\begin{equation*}
    (w-z)e^\epsilon(\vv^2z)\psi^{\epsilon'}(w)=
    \psi^{\epsilon'}(w)\left((\vv^{-2}w-\vv^2z)e^\epsilon(\vv^2z)-(\vv^{-2}-\vv^2)ze^{\epsilon'}(\vv^2w)\right).
\end{equation*}
We also use Lemma~\ref{Properties of half-currents}(c) to move
$f^\bullet(\cdot)$ to the left of $e^\bullet(\cdot)$. After obvious
cancelations, everything boils down to proving
\begin{equation*}
    (\vv^{-1}z-\vv w)e^\epsilon(\vv^2z)e^{\epsilon'}(\vv^2w)-
    (\vv z-\vv^{-1}w)e^{\epsilon'}(\vv^2w)e^\epsilon(\vv^2z)=
    (\vv^{-1}-\vv)(ze^{\epsilon'}(\vv^2w)^2+we^\epsilon(\vv^2z)^2),
\end{equation*}
which is exactly the statement of Lemma~\ref{Properties of
half-currents}(d2).


\subsubsection{Verification of~(\ref{ABCD10})}
\

This verification is completely analogous to the above verification
of~(\ref{ABCD9}) and is left to the interested reader.


\subsubsection{Verification of~(\ref{ABCD11})}
\

We need to prove
  $(z-w)[A^\epsilon(z),D^{\epsilon'}(w)]=
   (\vv-\vv^{-1})(wC^{\epsilon'}(w)B^\epsilon(z)-zC^\epsilon(z)B^{\epsilon'}(w))$.
The LHS equals
  $(\vv-\vv^{-1})^2(z-w)
   (A^\epsilon(z)f^{\epsilon'}(w)A^{\epsilon'}(w)e^{\epsilon'}(w)-
    f^{\epsilon'}(w)A^{\epsilon'}(w)e^{\epsilon'}(w)A^\epsilon(z))$.
Applying Lemma~\ref{Properties of half-currents}(b2) to the first
summand and Lemma~\ref{Properties of half-currents}(a2) to the
second summand, we see that the LHS equals
\begin{equation*}
\begin{split}
  & (\vv-\vv^{-1})^2((\vv z-\vv^{-1}w)f^{\epsilon'}(w)+(\vv^{-1}-v)zf^{\epsilon}(z))
    A^\epsilon(z)A^{\epsilon'}(w)e^{\epsilon'}(w)-\\
  & (\vv-\vv^{-1})^2 f^{\epsilon'}(w)A^\epsilon(z)A^{\epsilon'}(w)
    ((\vv z-\vv^{-1}w)e^{\epsilon'}(w)-(\vv-\vv^{-1})we^{\epsilon}(z))=\\
  & (\vv-\vv^{-1})^3(wf^{\epsilon'}(w)A^{\epsilon'}(w)A^\epsilon(z)e^\epsilon(z)-
    zf^\epsilon(z)A^\epsilon(z)A^{\epsilon'}(w)e^{\epsilon'}(w)),
\end{split}
\end{equation*}
which obviously coincides with the RHS.


\subsubsection{Verification of~(\ref{ABCD12})}
\

We need to prove
  $A^\epsilon(z)D^\epsilon(\vv^{-2}z)-\vv^{-1}B^\epsilon(z)C^\epsilon(\vv^{-2}z)=
   z^{\epsilon b^\epsilon}$.
Due to~Lemma~\ref{Properties of half-currents}(b3), we have
  $f^\epsilon(\vv^{-2}z)A^\epsilon(\vv^{-2}z)=\vv A^\epsilon(\vv^{-2}z)f^\epsilon(z)$.
Thus,
\begin{equation*}
    A^\epsilon(z)D^\epsilon(\vv^{-2}z)=A^\epsilon(z)A^\epsilon(\vv^{-2}z)
    (\psi^\epsilon(\vv^{-2}z)+\vv(\vv-\vv^{-1})^2f^\epsilon(z)e^\epsilon(\vv^{-2}z)),
\end{equation*}
\begin{equation*}
    B^\epsilon(z)C^\epsilon(\vv^{-2}z)=
    \vv(\vv-\vv^{-1})^2A^\epsilon(z)e^\epsilon(z)A^\epsilon(\vv^{-2}z)f^\epsilon(z).
\end{equation*}
According to Lemma~\ref{Properties of half-currents}(a2), we have
  $e^\epsilon(z)A^\epsilon(\vv^{-2}z)=\vv A^\epsilon(\vv^{-2}z)e^\epsilon(\vv^{-2}z)$.
Hence,
\begin{equation*}
    B^\epsilon(z)C^\epsilon(\vv^{-2}z)=
    \vv^2(\vv-\vv^{-1})^2A^\epsilon(z)A^\epsilon(\vv^{-2}z)e^\epsilon(\vv^{-2}z)f^\epsilon(z).
\end{equation*}
Due to Lemma~\ref{Properties of half-currents}(c), we have
  $-\vv(\vv-\vv^{-1})^2[e^\epsilon(\vv^{-2}z),f^\epsilon(z)]=
   \psi^\epsilon(z)-\psi^\epsilon(\vv^{-2}z)$.
Therefore, we finally get
\begin{equation*}
\begin{split}
  & A^\epsilon(z)D^\epsilon(\vv^{-2}z)-\vv^{-1}B^\epsilon(z)C^\epsilon(\vv^{-2}z)=
  A^\epsilon(z)A^\epsilon(\vv^{-2}z)\psi^\epsilon(\vv^{-2}z)-\\
  & \vv(\vv-\vv^{-1})^2A^\epsilon(z)A^\epsilon(\vv^{-2}z)[e^\epsilon(\vv^{-2}z),f^\epsilon(z)]=
  A^\epsilon(z)A^\epsilon(\vv^{-2}z)\psi^\epsilon(z)=z^{\epsilon b^\epsilon},
\end{split}
\end{equation*}
which completes our verification of~(\ref{ABCD12}).


\subsection{Proof of Theorem~\ref{Relations bw ABCD}(a) for a general $\fg$}
\

First, let us derive an explicit formula for  $A^\pm_i(z)$. Recall the
elements $\{h_{i,\pm r}\}_{i\in I}^{r>0}$ of~Section~\ref{sqaa}, such that
  $z^{\mp b^\pm_i}(\psi^\pm_{i,\mp b^\pm_i})^{-1}\psi^\pm_i(z)=
    \exp\left(\pm(\vv_i-\vv_i^{-1})\sum_{r>0}h_{i,\pm r}z^{\mp r}\right).$
For $r\ne 0$, consider the following $I\times I$ matrix $C_\vv(r)$:
  $$C_\vv(r)_{ij}=
    \begin{cases}
       0 & \text{if}\ \ c_{ij}=0, \\
       -1-\vv_i^{2r} & \text{if}\ \ j=i, \\
       \frac{\vv_j-\vv_j^{-1}}{\vv_i-\vv_i^{-1}}\sum_{p=1}^{-c_{ji}}\vv_j^{r(c_{ji}+2p)} & \text{if}\ \  j - i.
    \end{cases}$$
Set
  $t_{i,r}:=\sum_{j\in I}(C_\vv(r)^{-1})_{ij}h_{j,r}$
(matrix $C_\vv(r)$ is invertible, due to~Lemma~\ref{Matrices} below).
Define
\begin{equation}\label{explicit A}
  A^\pm_i(z):=(\phi^\pm_i)^{-1}\cdot \exp\left(\pm (\vv_i-\vv_i^{-1})\sum_{r>0}t_{i,\pm r}z^{\mp r}\right).
\end{equation}
These $A^\pm_i(z)$ satisfy
  $z^{\mp b^\pm_i}\psi^\pm_i(z)=
    \frac{\prod_{j-i}\prod_{p=1}^{-c_{ji}}A^\pm_j(\vv_j^{-c_{ji}-2p}z)}{A^\pm_i(z)A^\pm_i(\vv_i^{-2}z)}$
as well as $A^\pm_{i,0}=(\phi^\pm_i)^{-1}$.
This provides an explicit formula for $A^\pm_i(z)$, which we
referred to in Section~\ref{ABCD generators}.

\begin{Rem}\label{uniqueness of A}
Comparing the coefficients of $z^{\mp r}\ (r>0)$ in the system of
equations~(\ref{A}) for all $i$, we see that $A_{i,\pm r}$ are recovered uniquely
modulo the values of $A_{i, \pm s} (0\leq s<r)$, due to invertibility of $C_\vv(r)$.
Therefore, an induction in $r$ implies that the system of equations~(\ref{A})
has a unique solution $\{A^\pm_i(z)\}_{i\in I}$, hence, given by~(\ref{explicit A}).
\end{Rem}

Define auxiliary $I\times I$ matrices $B_\vv(r), D_\vv(r)$ via
  $B_\vv(r)_{ij}=\frac{[rc_{ij}]_{\vv_i}}{r},
   D_\vv(r)_{ij}=\delta_{ij}\frac{\vv_j^{-2r}-1}{r(\vv_j-\vv_j^{-1})}.$
The matrix $B_\vv(r)$ is a $\vv$-version of the Cartan matrix of $\fg$ and it
is known to be invertible for any $r\ne 0$. The following is straightforward.

\begin{Lem}\label{Matrices}
For $r\ne 0$, we have $B_\vv(r)=C_\vv(r)D_\vv(r)$.
In particular, $C_\vv(r)$ is invertible.
\end{Lem}

The following result is an immediate corollary of Lemma~\ref{Matrices} and
relations~(\ref{U4'}, \ref{U5'}).

\begin{Lem}\label{splitting A+}
For $\epsilon\in \{\pm\}$, we have:

\noindent
(a)
  $(\vv_iz-\vv_i^{-1}w)A^\epsilon_i(z)e_i(w)=(z-w)e_i(w)A^\epsilon_i(z)$, while
  $A^\epsilon_i(z)e_j(w)=e_j(w)A^\epsilon_i(z)$ for $j\ne i$.

\noindent
(b)
  $(z-w)A^\epsilon_i(z)f_i(w)=(\vv_iz-\vv_i^{-1}w)f_i(w)A^\epsilon_i(z)$, while
  $A^\epsilon_i(z)f_j(w)=f_j(w)A^\epsilon_i(z)$ for $j\ne i.$
\end{Lem}

Now we are ready to sketch the proof of Theorem~\ref{Relations bw
 ABCD}(a) for a general $\fg$.


\subsubsection{Verification of~(\ref{ABCD1}) and (\ref{ABCD2})}
\

Relations~(\ref{ABCD1}, \ref{ABCD2}) follow from~Lemma~\ref{splitting A+}
and relations~(\ref{U1}, \ref{U6}).


\subsubsection{Verification of~(\ref{ABCD5}--\ref{ABCD12})}
\

Let us introduce the series $\bar{A}^\pm_i(z)$ via
  $z^{\mp b^\pm_i}\psi^\pm_i(z)=\frac{1}{\bar{A}^\pm_i(z)\bar{A}^\pm_i(\vv_i^{-2}z)}$,
and define the generating series
$\bar{B}^\pm_i(z),\bar{C}^\pm_i(z),\bar{D}^\pm_i(z)$ by using
formulas~(\ref{B}--\ref{D}) but with $\bar{A}^\pm_i(z)$ instead
of $A^\pm_i(z)$. For a fixed $i$, these series satisfy the
corresponding relations~(\ref{ABCD5}--\ref{ABCD12}) of the
$\ssl_2$ case. However, $A^\pm_i(z)\bar{A}^\pm_i(z)^{-1}$ is
expressed through $\{A^\pm_j(z)\}_{j\ne i}$, hence, commutes with
$e^\epsilon_i(z),f^\epsilon_i(z), A^\epsilon_i(z)$, due to
Lemma~\ref{splitting A+}. Relations~(\ref{ABCD5}--\ref{ABCD12}) follow
(this also explains the RHS of~(\ref{ABCD12})).


\subsubsection{Verification of~(\ref{ABCD13})}
\

Analogously to Lemma~\ref{Properties of half-currents}(d1),
relation~(\ref{U2}) implies the following equality:
\begin{equation*}
    (z-\vv_i^{c_{ij}}w)e^\epsilon_i(z)e^{\epsilon'}_j(w)-(\vv_i^{c_{ij}}z-w)e^{\epsilon'}_j(w)e^\epsilon_i(z)=
    z[e_{i,0},e^{\epsilon'}_j(w)]_{\vv_i^{c_{ij}}}+w[e_{j,0},e^\epsilon_i(z)]_{\vv_i^{c_{ij}}}
\end{equation*}
for any $\epsilon,\epsilon'\in \{\pm\}$ (wee also note that these
equalities for all possible $\epsilon,\epsilon'$ imply~(\ref{U2})).
Multiplying the above equality by
$(\vv_i-\vv_i^{-1})(\vv_j-\vv_j^{-1})A^\epsilon_i(z)A^{\epsilon'}_j(w)$
on the left and using Lemma~\ref{splitting A+}(a),
relation~(\ref{ABCD1}), and an equality $(\vv_i-\vv_i^{-1})e_{i,0}=\phi^+_iB^+_{i,0}$,
we obtain~(\ref{ABCD13}).


\subsubsection{Verification of~(\ref{ABCD14})}
\

Analogously to Lemma~\ref{Properties of half-currents}(e1),
relation~(\ref{U3}) implies the following equality:
\begin{equation*}
    (\vv_i^{c_{ij}}z-w)f^\epsilon_i(z)f^{\epsilon'}_j(w)-(z-\vv_i^{c_{ij}}w)f^{\epsilon'}_j(w)f^\epsilon_i(z)=
    -[f^{\epsilon'}_j(w),f_{i,1}]_{\vv_i^{c_{ij}}}-[f^\epsilon_i(z),f_{j,1}]_{\vv_i^{c_{ij}}}
\end{equation*}
for any $\epsilon,\epsilon'\in \{\pm\}$ (wee also note that these
equalities for all possible $\epsilon,\epsilon'$ imply~(\ref{U3})).
Multiplying the above equality by
  $(\vv_i-\vv_i^{-1})(\vv_j-\vv_j^{-1})A^\epsilon_i(z)A^{\epsilon'}_j(w)$
on the right and using Lemma~\ref{splitting A+}(b),
relation~(\ref{ABCD1}), and an equality
$(\vv_i-\vv_i^{-1})f_{i,1}=C^+_{i,1}\phi^+_i$, we obtain~(\ref{ABCD14}).


\subsubsection{Verification of~(\ref{ABCD15})}
\

\medskip
\noindent
$\circ$ \emph{Case $c_{ij}=0$.}
The equality $[B^\epsilon_i(z),B^{\epsilon'}_j(w)]=0$
follows immediately from Lemma~\ref{splitting A+}(a) and
$[e^\epsilon_i(z),e^{\epsilon'}_j(w)]=0$, which is a consequence of
the corresponding Serre relation~(\ref{U7}).

\medskip
\noindent
$\circ$ \emph{Case $c_{ij}=-1$.}
The corresponding Serre relation~(\ref{U7}) is equivalent to
\begin{equation*}
    \{e^{\epsilon_1}_i(z_1)e^{\epsilon_2}_i(z_2)e^{\epsilon'}_j(w)-
      (\vv_i+\vv_i^{-1})e^{\epsilon_1}_i(z_1)e^{\epsilon'}_j(w)e^{\epsilon_2}_i(z_2)+
      e^{\epsilon'}_j(w)e^{\epsilon_1}_i(z_1)e^{\epsilon_2}_i(z_2)\}+
      \{z_1\leftrightarrow z_2\}=0
\end{equation*}
for any $\epsilon_1,\epsilon_2,\epsilon'\in \{\pm\}$. Let us denote the
first $\{\cdots\}$ in the LHS by $J^{\epsilon_1,\epsilon_2,\epsilon'}(z_1,z_2,w)$.
Set
\begin{equation*}
    M:=(\vv_i-\vv_i^{-1})^2(\vv_j-\vv_j^{-1})(\vv_iz_1-\vv_i^{-1}z_2)(\vv_iz_2-\vv_i^{-1}z_1)
    A^{\epsilon_1}_i(z_1)A^{\epsilon_2}_i(z_2)A^{\epsilon'}_j(w).
\end{equation*}
Combining the equality
\begin{equation*}
    (\vv_i z_2-\vv_i^{-1} z_1)A^{\epsilon_2}_i(z_2)e^{\epsilon_1}_i(z_1)=
    (z_2-z_1)e^{\epsilon_1}_i(z_1)A^{\epsilon_2}_i(z_2)+
    (\vv_i-\vv_i^{-1})z_1A^{\epsilon_2}_i(z_2)e^{\epsilon_2}_i(z_2)
\end{equation*}
(see Lemma~\ref{Properties of half-currents}(a2)) with
Lemma~\ref{splitting A+}(a), we find
\begin{equation*}
\begin{split}
  & M\cdot J^{\epsilon_1,\epsilon_2,\epsilon'}(z_1,z_2,w)=
    \frac{(\vv_i-\vv_i^{-1})z_1}{\vv_iz_2-\vv_i^{-1}z_1}M\cdot J^{\epsilon_2,\epsilon_2,\epsilon'}(z_2,z_2,w)+
    (z_2-z_1)(\vv_iz_1-\vv_i^{-1}z_2)\times\\
  & \{B^{\epsilon_1}_i(z_1)B^{\epsilon_2}_i(z_2)B^{\epsilon'}_j(w)-
    (\vv_i+\vv_i^{-1})B^{\epsilon_1}_i(z_1)B^{\epsilon'}_j(w)B^{\epsilon_2}_i(z_2)+
    B^{\epsilon'}_j(w)B^{\epsilon_1}_i(z_1)B^{\epsilon_2}_i(z_2)\}.
\end{split}
\end{equation*}
The first summand in the RHS is zero as
  $J^{\epsilon_2,\epsilon_2,\epsilon'}(z_2,z_2,w)=0$.
Therefore, multiplying
  $J^{\epsilon_1,\epsilon_2,\epsilon'}(z_1,z_2,w)+J^{\epsilon_2,\epsilon_1,\epsilon'}(z_2,z_1,w)=0$
by $M$ on the left, we obtain~(\ref{ABCD15}).

\medskip
\noindent
$\circ$ \emph{Case $c_{ij}=-2,-3$.}
These cases are treated similarly to $c_{ij}=-1$, but the corresponding computations become
more cumbersome. We verified these cases using \textbf{MATLAB}.


\subsubsection{Verification of~(\ref{ABCD16})}
\

This verification is analogous to that of~(\ref{ABCD15})
and is left to the interested reader.

\subsection{Proof of Theorem~\ref{Relations bw ABCD}(b)}
\

Part (b) of Theorem~\ref{Relations bw ABCD} can be obtained by reversing the above arguments.
In other words, starting from the algebra generated by $(A^\pm_{i,0})^{-1}$
and the coefficients of the currents $A^\pm_i(z), B^\pm_i(z), C^\pm_i(z), D^\pm_i(z)$
with the defining relations~(\ref{ABCD0}--\ref{ABCD16}), we need to show
that the elements $\phi^\pm_i$ and currents $e_i(z), f_i(z), \psi^\pm_i(z)$, defined
via~(\ref{A}--\ref{D}), satisfy relations~(\ref{U1}--\ref{U10}).

\medskip
This completes our proof of Theorem~\ref{Relations bw ABCD}.


\section{Proof of Theorem~\ref{Homomorphism}}\label{Proof Homomorphism}

We denote the images of $e_i(z), f_i(z), \psi^\pm_i(z)$ under
$\wt{\Phi}^{\unl\lambda}_\mu$ by $E_i(z), F_i(z), \Psi_i(z)^\pm$. It
suffices to prove that they satisfy relations~(\ref{U1}--\ref{U8}), since
relations~(\ref{U9},~\ref{U10}) are obviously preserved by
$\wt{\Phi}^{\unl\lambda}_\mu$. While checking these relations, we
will use LHS and RHS when referring to their left-hand and
right-hand sides. Set
  $\rho^+_i:=\frac{-\vv_i}{1-\vv_i^2},
   \rho^-_i:=\frac{1}{1-\vv_i^2},
   W_{i,rs}(z):=\prod\limits_{1\leq t\leq a_i}^{r\ne t\ne s}(1-\frac{\sw_{i,t}}{z})$.



\subsection{Compatibility with~(\ref{U1})}
\

First, we check that the range of powers of $z$ in $\psi^\pm_i(z)$
and $\Psi_i(z)^\pm$ agree. Note that
\begin{equation*}
    (1-\nu/z)^+=1-\nu\cdot z^{-1}\in \BC[[z^{-1}]],\
    (1/(1-\nu/z))^+=1+\nu z^{-1}+\nu^2 z^{-2}+\ldots \in \BC[[z^{-1}]],
\end{equation*}
\begin{equation*}
    (1-\nu/z)^-=-\nu\cdot z^{-1}(1-z/\nu)\in z^{-1}\BC[[z]],\
    (1/(1-\nu/z))^-=-z/\nu-z^2/\nu^2-\ldots \in z\BC[[z]].
\end{equation*}
Therefore, $\Psi_i(z)^+$ contains only nonpositive powers of $z$,
while $\Psi_i(z)^-$ contains only powers of $z$ bigger of equal to
\begin{equation*}
    -\#\{s:i_s=i\} + 2a_i -\sum_{j-i} a_j(-c_{ji})=
    -\alphavee_i(\lambda)+ \alphavee_i(\lambda-\mu)=
    -\alphavee_i(\mu)=-\alphavee_i(\mu^-)=-b^-_i.
\end{equation*}
Moreover, the coefficients of $z^0$ in $\Psi_i(z)^+$ and of
$z^{-b^-_i}$ in $\Psi_i(z)^-$ are invertible.

The equality $[\Psi_i(z)^\epsilon,\Psi_j(w)^{\epsilon'}]=0$ follows
from the commutativity of $\{\sw^{\pm 1/2}_{i,r}\}_{i\in I}^{1\leq r\leq a_i}$.

%


\subsection{Compatibility with~(\ref{U2})}
\

\medskip
\noindent
$\circ$ \emph{Case $c_{ij}=0$}.
The equality $[E_i(z),E_j(w)]=0$ is obvious in this case, since $D^{-1}_{i,r}$
commute with $\sw^{\pm 1/2}_{k,s}$ for $k=j$ or $k\to j$, while
$D^{-1}_{j,s}$ commute with $\sw^{\pm 1/2}_{k,r}$ for $k=i$ or $k\to i$.

\medskip
\noindent
$\circ$ \emph{Case $c_{ij}=2$.}
We may assume $\fg=\ssl_2$ and we will drop the index $i$ from our
notation. We need to prove
  $(z-\vv^2w)E(z)E(w)/(\rho^+)^2=-(w-\vv^2z)E(w)E(z)/(\rho^+)^2$.
The LHS equals
\begin{equation*}
\begin{split}
  & \vv^{-2}\prod_{t=1}^a \sw_t^2\cdot (z-\vv^2w)\cdot \sum_{r=1}^a
    \delta\left(\frac{\sw_r}{z}\right)\delta\left(\frac{\vv^{-2}\sw_r}{w}\right)
    \frac{Z(\sw_r)Z(\vv^{-2}\sw_r)}{W_r(\sw_r)W_r(\vv^{-2}\sw_r)}D_r^{-2}+\\
  & \vv^{-2}\prod_{t=1}^a \sw_t^2\cdot (z-\vv^2w)\cdot \sum_{1\leq r\ne s\leq a}
    \delta\left(\frac{\sw_r}{z}\right)\delta\left(\frac{\sw_s}{w}\right)
    \frac{Z(\sw_r)Z(\sw_s)}{W_r(\sw_r)W_{rs}(\sw_s)(1-\vv^{-2}\sw_r/\sw_s)}D_r^{-1}D_s^{-1}.
\end{split}
\end{equation*}
Using the equality
\begin{equation}\label{delta property 1}
   G(z,w)\delta\left(\frac{\nu_1}{z}\right)\delta\left(\frac{\nu_2}{w}\right)=
   G(\nu_1,\nu_2)\delta\left(\frac{\nu_1}{z}\right)\delta\left(\frac{\nu_2}{w}\right),
\end{equation}
we see that the first sum is zero, while the second sum equals
\begin{equation*}
\begin{split}
  & \prod_{t=1}^a \sw_t^2\cdot \sum_{1\leq r\ne s\leq a}
    \delta\left(\frac{\sw_r}{z}\right)\delta\left(\frac{\sw_s}{w}\right)
    \frac{Z(\sw_r)Z(\sw_s)}{W_{rs}(\sw_r)W_{rs}(\sw_s)}
    \frac{\vv^{-2}(\sw_r-\vv^2\sw_s)}{(1-\sw_s/\sw_r)(1-\vv^{-2}\sw_r/\sw_s)}D_r^{-1}D_s^{-1}=\\
  & \prod_{t=1}^a \sw_t^2\cdot \sum_{1\leq r\ne s\leq a}
    \delta\left(\frac{\sw_r}{z}\right)\delta\left(\frac{\sw_s}{w}\right)
    \frac{Z(\sw_r)Z(\sw_s)}{W_{rs}(\sw_r)W_{rs}(\sw_s)}\frac{\sw_r\sw_s}{\sw_s-\sw_r}D_r^{-1}D_s^{-1}.
\end{split}
\end{equation*}
Swapping $z$ and $w$, we see that $-(w-\vv^2z)E(w)E(z)/(\rho^+)^2$ equals
  $$-\prod_{t=1}^a \sw_t^2\cdot \sum_{1\leq r\ne s\leq a}
    \delta\left(\frac{\sw_r}{w}\right)\delta\left(\frac{\sw_s}{z}\right)
    \frac{Z(\sw_r)Z(\sw_s)}{W_{rs}(\sw_r)W_{rs}(\sw_s)}\frac{\sw_r\sw_s}{\sw_s-\sw_r}D_r^{-1}D_s^{-1}.$$
Swapping $r$ and $s$ in the latter sum, we get exactly the same expression as for the LHS.

\medskip
\noindent
$\circ$ \emph{Case $c_{ij}<0$.}
In this case, we can assume $I=\{i,j\}$ and $i\to j$. We need to prove
  $(z-\vv_i^{c_{ij}}w)E_i(z)E_j(w)/(\rho^+_i\rho^+_j)=
   (\vv_i^{c_{ij}}z-w)E_j(w)E_i(z)/(\rho^+_i\rho^+_j)$.
The LHS equals
\begin{equation*}
\begin{split}
   & \vv_i^{-c_{ij}}\prod_{t=1}^{a_i} \sw^{1+c_{ij}/2}_{i,t}
   \prod_{t=1}^{a_j} \sw_{j,t}\cdot (z-\vv_i^{c_{ij}}w)\times\\
   & \sum_{1\leq r\leq a_i}^{1\leq s\leq a_j}
   \delta\left(\frac{\sw_{i,r}}{z}\right)\delta\left(\frac{\sw_{j,s}}{w}\right)
   \frac{Z_i(\sw_{i,r})}{W_{i,r}(\sw_{i,r})}D_{i,r}^{-1}\frac{Z_j(\sw_{j,s})}{W_{j,s}(\sw_{j,s})}
   \prod_{p=1}^{-c_{ij}}W_i(\vv_i^{-c_{ij}-2p}w)D_{j,s}^{-1}=\\
   & \prod_{t=1}^{a_i} \sw^{1+c_{ij}/2}_{i,t}\prod_{t=1}^{a_j} \sw_{j,t}\cdot A(z,w)\times\\
   & \sum_{1\leq r\leq a_i}^{1\leq s\leq a_j}
   \delta\left(\frac{\sw_{i,r}}{z}\right)\delta\left(\frac{\sw_{j,s}}{w}\right)
   \frac{Z_i(\sw_{i,r})Z_j(\sw_{j,s})\prod_{p=1}^{-c_{ij}}W_{i,r}(\vv_i^{-c_{ij}-2p}\sw_{j,s})}
        {W_{i,r}(\sw_{i,r})W_{j,s}(\sw_{j,s})}
   D_{i,r}^{-1}D_{j,s}^{-1},
\end{split}
\end{equation*}
where
  $A(z,w)=\vv_i^{-c_{ij}}(z-\vv_i^{c_{ij}}w)\prod_{p=1}^{-c_{ij}}\left(1-\frac{\vv_i^{-2}z}{\vv_i^{-c_{ij}-2p}w}\right)$,
due to~(\ref{delta property 1}).
Likewise, the RHS equals
\begin{equation*}
\begin{split}
  & \prod_{t=1}^{a_i} \sw^{1+c_{ij}/2}_{i,t}\prod_{t=1}^{a_j} \sw_{j,t}\cdot (\vv_i^{c_{ij}}z-w)\times\\
  & \sum_{1\leq r\leq a_i}^{1\leq s\leq a_j}
  \delta\left(\frac{\sw_{j,s}}{w}\right)\delta\left(\frac{\sw_{i,r}}{z}\right)
  \frac{Z_j(\sw_{j,s})}{W_{j,s}(\sw_{j,s})}\prod_{p=1}^{-c_{ij}}W_i(\vv_i^{-c_{ij}-2p}w)
  D_{j,s}^{-1}\frac{Z_i(\sw_{i,r})}{W_{i,r}(\sw_{i,r})}D_{i,r}^{-1}=\\
  & \prod_{t=1}^{a_i} \sw^{1+c_{ij}/2}_{i,t}\prod_{t=1}^{a_j} \sw_{j,t}\cdot B(z,w)\times\\
  & \sum_{1\leq r\leq a_i}^{1\leq s\leq a_j}
  \delta\left(\frac{\sw_{i,r}}{z}\right)\delta\left(\frac{\sw_{j,s}}{w}\right)
  \frac{Z_i(\sw_{i,r})Z_j(\sw_{j,s})\prod_{p=1}^{-c_{ij}}W_{i,r}(\vv_i^{-c_{ij}-2p}\sw_{j,s})}
       {W_{i,r}(\sw_{i,r})W_{j,s}(\sw_{j,s})}
  D_{i,r}^{-1}D_{j,s}^{-1},
\end{split}
\end{equation*}
where
  $B(z,w)=(\vv_i^{c_{ij}}z-w)\prod_{p=1}^{-c_{ij}}\left(1-\frac{z}{\vv_i^{-c_{ij}-2p}w}\right)$,
due to~(\ref{delta property 1}).

The equality LHS$=$RHS follows from $A(z,w)=B(z,w)$.

\subsection{Compatibility with~(\ref{U3})}
\

\medskip
\noindent
$\circ$ \emph{Case $c_{ij}=0$}.
The equality $[F_i(z),F_j(w)]=0$ is obvious in this case, since $D_{i,r}$ commute
with $\sw^{\pm 1/2}_{k,s}$ for $k=j$ or $k\leftarrow j$, while
$D_{j,s}$ commute with $\sw^{\pm 1/2}_{k,r}$ for $k=i$ or $k\leftarrow i$.

\medskip
\noindent
$\circ$ \emph{Case $c_{ij}=2$.}
We may assume $\fg=\ssl_2$ and we will drop the index $i$ from our
notation. We need to prove
  $(\vv^2z-w)F(z)F(w)/(\rho^-)^2=-(\vv^2w-z)F(w)F(z)/(\rho^-)^2$.
The LHS equals
\begin{equation*}
\begin{split}
  & (\vv^2z-w)\cdot \sum_{r=1}^a
    \delta\left(\frac{\vv^2\sw_r}{z}\right)\delta\left(\frac{\vv^4\sw_r}{w}\right)
    \frac{1}{W_r(\sw_r)W_r(\vv^2\sw_r)}D_r^2+\\
  & (\vv^2z-w)\cdot \sum_{1\leq r\ne s\leq a}
    \delta\left(\frac{\vv^2\sw_r}{z}\right)\delta\left(\frac{\vv^2\sw_s}{w}\right)
    \frac{1}{W_r(\sw_r)W_{rs}(\sw_s)(1-\vv^2\sw_r/\sw_s)}D_rD_s.
\end{split}
\end{equation*}
Using equality~(\ref{delta property 1}), we see that the first
sum is zero, while the second sum equals
\begin{equation*}
\begin{split}
  & \sum_{1\leq r\ne s\leq a} \delta\left(\frac{\vv^2\sw_r}{z}\right)\delta\left(\frac{\vv^2\sw_s}{w}\right)
    \frac{1}{W_{rs}(\sw_r)W_{rs}(\sw_s)}\frac{\vv^4\sw_r-\vv^2\sw_s}{(1-\sw_s/\sw_r)(1-\vv^2\sw_r/\sw_s)}D_rD_s=\\
  & \sum_{1\leq r\ne s\leq a} \delta\left(\frac{\vv^2\sw_r}{z}\right)\delta\left(\frac{\vv^2\sw_s}{w}\right)
    \frac{1}{W_{rs}(\sw_r)W_{rs}(\sw_s)}\frac{\vv^2\sw_r\sw_s}{\sw_s-\sw_r}D_rD_s.
\end{split}
\end{equation*}
Swapping $z$ and $w$, we see that $-(\vv^2w-z)F(w)F(z)/(\rho^-)^2$ equals
  $$-\sum_{1\leq r\ne s\leq a} \delta\left(\frac{\vv^2\sw_r}{w}\right)\delta\left(\frac{\vv^2\sw_s}{z}\right)
    \frac{1}{W_{rs}(\sw_r)W_{rs}(\sw_s)}\frac{\vv^2\sw_r\sw_s}{\sw_s-\sw_r}D_rD_s.$$
Swapping $r$ and $s$ in this sum, we get exactly the same expression as for the LHS.

\medskip
\noindent
$\circ$ \emph{Case $c_{ij}<0$.}
In this case, we can assume $I=\{i,j\}$ and $i\to j$. Recall that
$\vv_i^{c_{ij}}=\vv_j^{c_{ji}}$. We need to prove
  $(\vv_j^{c_{ji}}z-w)F_i(z)F_j(w)/(\rho^-_i\rho^-_j)=
   (z-\vv_j^{c_{ji}}w)F_j(w)F_i(z)/(\rho^-_i\rho^-_j)$.
The LHS equals
\begin{equation*}
\begin{split}
  & \prod_{t=1}^{a_j} \sw^{c_{ji}/2}_{j,t}\cdot (\vv_j^{c_{ji}}z-w)\times\\
  & \sum_{1\leq r\leq a_i}^{1\leq s\leq a_j}
  \delta\left(\frac{\vv_i^2\sw_{i,r}}{z}\right)\delta\left(\frac{\vv_j^2\sw_{j,s}}{w}\right)
  \frac{1}{W_{i,r}(\sw_{i,r})}\prod_{p=1}^{-c_{ji}}W_j(\vv_j^{-c_{ji}-2p}z)D_{i,r}\frac{1}{W_{j,s}(\sw_{j,s})}D_{j,s}=\\
  & \prod_{t=1}^{a_j} \sw^{c_{ji}/2}_{j,t}\cdot A(z,w)\cdot \sum_{1\leq r\leq a_i}^{1\leq s\leq a_j}
  \delta\left(\frac{\vv_i^2\sw_{i,r}}{z}\right)\delta\left(\frac{\vv_j^2\sw_{j,s}}{w}\right)
  \frac{\prod_{p=1}^{-c_{ji}}W_{j,s}(\vv_j^{-c_{ji}-2p}z)}{W_{i,r}(\sw_{i,r})W_{j,s}(\sw_{j,s})}
  D_{i,r}D_{j,s},
\end{split}
\end{equation*}
where
  $A(z,w)=(\vv_j^{c_{ji}}z-w)\prod_{p=1}^{-c_{ji}}\left(1-\frac{\vv_j^{-2}w}{\vv_j^{-c_{ji}-2p} z}\right)$,
due to~(\ref{delta property 1}).
Likewise, the RHS equals
\begin{equation*}
\begin{split}
  & \vv_j^{c_{ji}}\prod_{t=1}^{a_j} \sw^{c_{ji}/2}_{j,t}\cdot (z-\vv_j^{c_{ji}}w)\times\\
  & \sum_{1\leq r\leq a_i}^{1\leq s\leq a_j}
  \delta\left(\frac{\vv_j^2\sw_{j,s}}{w}\right)\delta\left(\frac{\vv_i^2\sw_{i,r}}{z}\right)
  \frac{1}{W_{j,s}(\sw_{j,s})}D_{j,s}\frac{1}{W_{i,r}(\sw_{i,r})}\prod_{p=1}^{-c_{ji}}W_j(\vv_j^{-c_{ji}-2p}z)D_{i,r}=\\
  & \prod_{t=1}^{a_j} \sw^{c_{ji}/2}_{j,t}\cdot B(z,w)\cdot \sum_{1\leq r\leq a_i}^{1\leq s\leq a_j}
  \delta\left(\frac{\vv_i^2\sw_{i,r}}{z}\right)\delta\left(\frac{\vv_j^2\sw_{j,s}}{w}\right)
  \frac{\prod_{p=1}^{-c_{ji}}W_{j,s}(\vv_j^{-c_{ji}-2p}z)}{W_{i,r}(\sw_{i,r})W_{j,s}(\sw_{j,s})}
  D_{i,r}D_{j,s},
\end{split}
\end{equation*}
where
  $B(z,w)=\vv_j^{c_{ji}}(z-\vv_j^{c_{ji}}w)\prod_{p=1}^{-c_{ji}}\left(1-\frac{w}{\vv_j^{-c_{ji}-2p}z}\right)$,
due to~(\ref{delta property 1}).

The equality LHS$=$RHS follows from $A(z,w)=B(z,w)$.

\subsection{Compatibility with~(\ref{U4})}
\

\medskip
\noindent
$\circ$ \emph{Case $c_{ij}=0$}.
The equality $[\Psi_i(z),E_j(w)]=0$ is obvious in this case, since $D^{-1}_{j,s}$
commute with $\sw^{\pm 1/2}_{k,r}$ for $k=i$ or $k - i$.

\medskip
\noindent
$\circ$ \emph{Case $c_{ij}=2$.}
We may assume $\fg=\ssl_2$ and we will drop the index $i$ from our
notation. We need to prove
  $(z-\vv^2w)\Psi(z)E(w)/\rho^+=(\vv^2z-w)E(w)\Psi(z)/\rho^+$.
The LHS equals
\begin{equation*}
\begin{split}
  & \prod_{t=1}^a \sw_t^2\cdot (z-\vv^2w)\cdot
  \frac{Z(z)}{W(z)W(\vv^{-2}z)} \sum_{r=1}^a
  \delta\left(\frac{\sw_r}{w}\right)\frac{Z(\sw_r)}{W_r(\sw_r)}D_r^{-1}=\\
  & \prod_{t=1}^a \sw_t^2\cdot \sum_{r=1}^a \delta\left(\frac{\sw_r}{w}\right)
  \frac{Z(z)Z(\sw_r)}{W_r(\sw_r)W_r(z)W_r(\vv^{-2}z)}
  \frac{z-\vv^2 w}{(1-w/z)(1-w/\vv^{-2}z)}D_r^{-1},
\end{split}
\end{equation*}
due to~(\ref{delta property 1}). Likewise, the RHS equals
\begin{equation*}
\begin{split}
  & \vv^{-2}\prod_{t=1}^a \sw_t^2\cdot (\vv^2z-w)\cdot \sum_{r=1}^a
  \delta\left(\frac{\sw_r}{w}\right)
  \frac{Z(\sw_r)}{W_r(\sw_r)}D_r^{-1}\frac{Z(z)}{W(z)W(\vv^{-2}z)}=\\
  & \prod_{t=1}^a \sw_t^2\cdot \sum_{r=1}^a \delta\left(\frac{\sw_r}{w}\right)
  \frac{Z(z)Z(\sw_r)}{W_r(\sw_r)W_r(z)W_r(\vv^{-2}z)}
  \frac{\vv^{-2}(\vv^2z-w)}{(1-\vv^{-2}w/z)(1-\vv^{-2}w/\vv^{-2}z)}D_r^{-1}.
\end{split}
\end{equation*}
The equality LHS$=$RHS follows.

\medskip
\noindent
$\circ$ \emph{Case $c_{ij}<0$.}
In this case, we can assume $I=\{i,j\}$. There are two situations to
consider: $i\to j$ and $i\leftarrow j$. Let us first treat the
former case. Since $\vv_i^{c_{ij}}=\vv_j^{c_{ji}}$, we need to prove
  $(z-\vv_j^{c_{ji}}w)\Psi_i(z)E_j(w)/\rho^+_j=(\vv_j^{c_{ji}}z-w)E_j(w)\Psi_i(z)/\rho^+_j$.
The LHS equals
\begin{equation*}
\begin{split}
  & \prod_{t=1}^{a_i} \sw^{1+c_{ij}/2}_{i,t}
  \prod_{t=1}^{a_j} \sw^{1+c_{ji}/2}_{j,t}\cdot (z-\vv_j^{c_{ji}}w) \times\\
  & \frac{Z_i(z)}{W_i(z)W_i(\vv_i^{-2}z)}\prod_{p=1}^{-c_{ji}}W_j(\vv_j^{-c_{ji}-2p}z)
  \sum_{s=1}^{a_j}\delta\left(\frac{\sw_{j,s}}{w}\right)
  \frac{Z_j(\sw_{j,s})}{W_{j,s}(\sw_{j,s})}\prod_{p'=1}^{-c_{ij}}W_i(\vv_i^{-c_{ij}-2p'}w)D_{j,s}^{-1}=\\
  & \prod_{t=1}^{a_i} \sw^{1+c_{ij}/2}_{i,t} \prod_{t=1}^{a_j} \sw^{1+c_{ji}/2}_{j,t}\cdot A(z,w)\times\\
  & \sum_{s=1}^{a_j}\delta\left(\frac{\sw_{j,s}}{w}\right)
  \frac{Z_i(z)Z_j(\sw_{j,s})\prod_{p'=1}^{-c_{ij}}W_i(\vv_i^{-c_{ij}-2p'}w) \prod_{p=1}^{-c_{ji}} W_{j,s}(\vv_j^{-c_{ji}-2p}z)}
       {W_i(z)W_i(\vv_i^{-2}z)W_{j,s}(\sw_{j,s})}
  D_{j,s}^{-1},
\end{split}
\end{equation*}
where
  $A(z,w)=(z-\vv_j^{c_{ji}}w)\prod_{p=1}^{-c_{ji}}\left(1-\frac{w}{\vv_j^{-c_{ji}-2p}z}\right)$.
Likewise, the RHS equals
\begin{equation*}
\begin{split}
  & \vv_j^{-c_{ji}}\prod_{t=1}^{a_i} \sw^{1+c_{ij}/2}_{i,t}
  \prod_{t=1}^{a_j} \sw^{1+c_{ji}/2}_{j,t}\cdot (\vv_j^{c_{ji}}z-w)\times\\
  & \sum_{s=1}^{a_j}\delta\left(\frac{\sw_{j,s}}{w}\right)
  \frac{Z_j(\sw_{j,s})}{W_{j,s}(\sw_{j,s})}
  \prod_{p'=1}^{-c_{ij}}W_i(\vv_i^{-c_{ij}-2p'}w)D_{j,s}^{-1}
  \frac{Z_i(z)}{W_i(z)W_i(\vv_i^{-2}z)}\prod_{p=1}^{-c_{ji}}W_j(\vv_j^{-c_{ji}-2p}z)=\\
  & \prod_{t=1}^{a_i} \sw^{1+c_{ij}/2}_{i,t} \prod_{t=1}^{a_j} \sw^{1+c_{ji}/2}_{j,t}\cdot B(z,w)\times\\
  & \sum_{s=1}^{a_j}\delta\left(\frac{\sw_{j,s}}{w}\right)
  \frac{Z_i(z)Z_j(\sw_{j,s})\prod_{p'=1}^{-c_{ij}}W_i(\vv_i^{-c_{ij}-2p'}w)\prod_{p=1}^{-c_{ji}} W_{j,s}(\vv_j^{-c_{ji}-2p}z)}
       {W_i(z)W_i(\vv_i^{-2}z)W_{j,s}(\sw_{j,s})}
  D_{j,s}^{-1},
\end{split}
\end{equation*}
where
 $B(z,w)=\vv_j^{-c_{ji}}(\vv_j^{c_{ji}}z-w)\prod_{p=1}^{-c_{ji}}\left(1-\frac{\vv_j^{-2}w}{\vv_j^{-c_{ji}-2p}z}\right)$.

The equality LHS$=$RHS follows from $A(z,w)=B(z,w)$.

\medskip
The case $i\leftarrow j$ is analogous:
$\Psi_i(z)$ is given by the same formula, while $E_j(w)$ differs by
an absence of the factor
  $\prod_{t=1}^{a_i} \sw_{i,t}^{c_{ij}/2}\cdot \prod_{p'=1}^{-c_{ij}}W_i(\vv_i^{-c_{ij}-2p'}w)$.
Tracing back the above calculations, it is clear that the equality
still holds when this factor is dropped out.

\subsection{Compatibility with~(\ref{U5})}
\

\medskip
\noindent
$\circ$ \emph{Case $c_{ij}=0$}.
The equality $[\Psi_i(z),F_j(w)]=0$ is obvious in this case, since $D_{j,s}$
commute with $\sw^{\pm 1/2}_{k,r}$ for $k=i$ or $k - i$.

\medskip
\noindent
 $\circ$ \emph{Case $c_{ij}=2$.} We may assume $\fg=\ssl_2$
and we will drop the index $i$ from our notation. We need to prove
  $(\vv^2z-w)\Psi(z)F(w)/\rho^-=(z-\vv^2w)F(w)\Psi(z)/\rho^-$.
The LHS equals
\begin{equation*}
\begin{split}
  & \prod_{t=1}^a \sw_t\cdot (\vv^2z-w)\cdot
  \frac{Z(z)}{W(z)W(\vv^{-2}z)}
  \sum_{r=1}^a \delta\left(\frac{\vv^2\sw_r}{w}\right)\frac{1}{W_r(\sw_r)}D_r=\\
  & \prod_{t=1}^a \sw_t\cdot\sum_{r=1}^a \delta\left(\frac{\vv^2\sw_r}{w}\right)
  \frac{Z(z)}{W_r(\sw_r)W_r(z)W_r(\vv^{-2}z)}
  \frac{\vv^2z-w}{(1-\vv^{-2}w/z)(1-\vv^{-2}w/\vv^{-2}z)}D_r,
\end{split}
\end{equation*}
due to~(\ref{delta property 1}). Likewise, the RHS equals
\begin{equation*}
\begin{split}
  & \vv^2\prod_{t=1}^a \sw_t\cdot (z-\vv^2w)\cdot
  \sum_{r=1}^a \delta\left(\frac{\vv^2\sw_r}{w}\right)
  \frac{1}{W_r(\sw_r)}D_r\frac{Z(z)}{W(z)W(\vv^{-2}z)}=\\
  & \prod_{t=1}^a \sw_t\cdot \sum_{r=1}^a \delta\left(\frac{\vv^2\sw_r}{w}\right)
  \frac{Z(z)}{W_r(\sw_r)W_r(z)W_r(\vv^{-2}z)}
  \frac{\vv^2(z-\vv^2w)}{(1-w/z)(1-w/\vv^{-2}z)}D_r.
\end{split}
\end{equation*}
The equality LHS$=$RHS follows.

\medskip
\noindent
 $\circ$ \emph{Case $c_{ij}<0$.}
In this case, we can assume $I=\{i,j\}$. There are two situations to
consider: $i\to j$ and $i\leftarrow j$. Let us first treat the
former case. Since $\vv_i^{c_{ij}}=\vv_j^{c_{ji}}$, we need to prove
  $(\vv_j^{c_{ji}}z-w)\Psi_i(z)F_j(w)/\rho^-_j=(z-\vv_j^{c_{ji}}w)F_j(w)\Psi_i(z)/\rho^-_j$.
The LHS equals
\begin{equation*}
\begin{split}
  & \prod_{t=1}^{a_i} \sw_{i,t} \prod_{t=1}^{a_j} \sw^{c_{ji}/2}_{j,t}\cdot (\vv_j^{c_{ji}}z-w)\times\\
  & \frac{Z_i(z)}{W_i(z)W_i(\vv_i^{-2}z)}\prod_{p=1}^{-c_{ji}}W_j(\vv_j^{-c_{ji}-2p}z)
  \sum_{s=1}^{a_j}\delta\left(\frac{\vv_j^2\sw_{j,s}}{w}\right)
  \frac{1}{W_{j,s}(\sw_{j,s})}D_{j,s}=\\
  & \prod_{t=1}^{a_i} \sw_{i,t} \prod_{t=1}^{a_j} \sw^{c_{ji}/2}_{j,t}\cdot A(z,w)
  \cdot \sum_{s=1}^{a_j}\delta\left(\frac{\vv_j^2\sw_{j,s}}{w}\right)
  \frac{Z_i(z)\prod_{p=1}^{-c_{ji}}W_{j,s}(\vv_j^{-c_{ji}-2p}z)}
       {W_i(z)W_i(\vv_i^{-2}z)W_{j,s}(\sw_{j,s})}
  D_{j,s},
\end{split}
\end{equation*}
where
  $A(z,w)=(\vv_j^{c_{ji}}z-w)\prod_{p=1}^{-c_{ji}}\left(1-\frac{\vv_j^{-2}w}{\vv_j^{-c_{ji}-2p}z}\right)$.
Likewise, the RHS equals
\begin{equation*}
\begin{split}
  & \vv_j^{c_{ji}}\prod_{t=1}^{a_i} \sw_{i,t}
  \prod_{t=1}^{a_j}\sw^{c_{ji}/2}_{j,t}\cdot (z-\vv_j^{c_{ji}}w) \times\\
  & \sum_{s=1}^{a_j}\delta\left(\frac{\vv_j^2\sw_{j,s}}{w}\right)
  \frac{1}{W_{j,s}(\sw_{j,s})}D_{j,s}
  \frac{Z_i(z)}{W_i(z)W_i(\vv_i^{-2}z)}\prod_{p=1}^{-c_{ji}}W_j(\vv_j^{-c_{ji}-2p}z)=\\
  & \prod_{t=1}^{a_i} \sw_{i,t} \prod_{t=1}^{a_j} \sw^{c_{ji}/2}_{j,t}\cdot B(z,w)
  \cdot \sum_{s=1}^{a_j}\delta\left(\frac{\vv_j^2\sw_{j,s}}{w}\right)
  \frac{Z_i(z)  \prod_{p=1}^{-c_{ji}}W_{j,s}(\vv_j^{-c_{ji}-2p}z)}
       {W_i(z)W_i(\vv_i^{-2}z)W_{j,s}(\sw_{j,s})}
  D_{j,s},
\end{split}
\end{equation*}
where
  $B(z,w)=\vv_j^{c_{ji}}(z-\vv_j^{c_{ji}}w)\prod_{p=1}^{-c_{ji}}\left(1-\frac{w}{\vv_j^{-c_{ji}-2p}z}\right)$.

The equality LHS$=$RHS follows from $A(z,w)=B(z,w)$.

\medskip
The case $i\leftarrow j$ is analogous:
$\Psi_i(z)$ is given by the same formula, while $F_j(w)$ has an extra factor
 $\prod_{t=1}^{a_i} \sw_{i,t}^{c_{ij}/2}\cdot \prod_{p'=1}^{-c_{ij}}W_i(\vv_i^{-c_{ij}-2p'}w)$.
The contributions of this factor into the LHS and the RHS are the
same, hence, the equality still holds.

\subsection{Compatibility with~(\ref{U6})}
\

\medskip
\noindent
$\circ$ \emph{Case $c_{ij}=0$}.
The equality $[E_i(z),F_j(w)]=0$ is obvious in this case, since $D^{-1}_{i,r}$
commute with $\sw^{\pm 1/2}_{k,s}$ for $k=i$ or $k\leftarrow j$,
while $D_{j,s}$ commute with $\sw^{\pm 1/2}_{k,r}$ for $k=i$ or $k\to i$.

\medskip
\noindent
$\circ$ \emph{Case $c_{ij}=2$.}
We may assume $\fg=\ssl_2$, and we will drop the index $i$ from our
notation. We need to prove
  $[E(z),F(w)]=\frac{1}{\vv-\vv^{-1}}\delta\left(\frac{z}{w}\right)\left(\Psi(z)^+-\Psi(z)^-\right)$.
The LHS equals
\begin{equation*}
\begin{split}
  & \rho^+\rho^-
    \left[\prod_{t=1}^a \sw_t\cdot \sum_{r=1}^a \delta\left(\frac{\sw_r}{z}\right)
          \frac{Z(\sw_r)}{W_r(\sw_r)}D_r^{-1},
          \sum_{s=1}^a \delta\left(\frac{\vv^2\sw_s}{w}\right)\frac{1}{W_s(\sw_s)}D_s
    \right]=
    \frac{-\vv}{(1-\vv^2)^2}\prod_{t=1}^a \sw_t\times\\
  & \left\{\sum_{r=1}^a
    \left(\delta\left(\frac{\sw_r}{z}\right)\delta\left(\frac{\sw_r}{w}\right)
          \frac{Z(\sw_r)}{W_r(\sw_r)W_r(\vv^{-2}\sw_r)}-
          \vv^2\delta\left(\frac{\vv^2\sw_r}{z}\right)\delta\left(\frac{\vv^2\sw_r}{w}\right)
          \frac{Z(\vv^2\sw_r)}{W_r(\sw_r)W_r(\vv^2\sw_r)}\right)+\right.\\
  & \left.\sum_{1\leq r\ne s\leq a}
    \delta\left(\frac{\sw_r}{z}\right)\delta\left(\frac{\vv^2\sw_s}{w}\right)
    \frac{Z(\sw_r)}{W_{rs}(\sw_r)W_{rs}(\sw_s)}
    \left(\frac{1}{A(z,w)}-\frac{\vv^2}{B(z,w)}\right)D_r^{-1}D_s\right\},
\end{split}
\end{equation*}
where
  $A(z,w)=(1-\vv^{-2}w/z)(1-\vv^{-2}z/\vv^{-2}w)$
and
  $B(z,w)=(1-z/\vv^{-2}w)(1-w/z)$.
The second sum is zero as $A(z,w)=\vv^{-2}B(z,w)$.

To evaluate the RHS, we need the following standard result.
\begin{Lem}\label{key equality}
For any rational function $\gamma(z)$ with simple poles
$\{x_t\}\subset \BC^\times$ and possibly poles of higher order at
$z=0,\infty$, the following equality holds:
\begin{equation}\label{delta property 2}
  \gamma(z)^+-\gamma(z)^-=
  \sum_{t}\delta\left(\frac{z}{x_t}\right)\Res_{z=x_t}\gamma(z)\frac{dz}{z}.
\end{equation}
\end{Lem}

\begin{proof}
Consider the partial fraction decomposition of $\gamma(z)$:
\begin{equation*}
  \gamma(z)=P(z)+\sum_{t}\frac{\nu_t}{z-x_t},
\end{equation*}
where $P(z)$ is a Laurent polynomial. Then
  $P(z)^\pm=P(z)\Rightarrow P(z)^+-P(z)^-=0$.
Meanwhile:
\begin{equation*}
    \left(\frac{\nu_t}{z-x_t}\right)^+=
    \frac{\nu_t}{z}+\frac{\nu_tx_t}{z^2}+\frac{\nu_tx_t^2}{z^3}+\ldots\
    \mathrm{and}\
    \left(\frac{\nu_t}{z-x_t}\right)^-=
    -\frac{\nu_t}{x_t}-\frac{\nu_tz}{x_t^2}-\frac{\nu_tz^2}{x_t^3}-\ldots,
\end{equation*}
so that
\begin{equation*}
    \left(\frac{\nu_t}{z-x_t}\right)^+-\left(\frac{\nu_t}{z-x_t}\right)^-=
    \frac{\nu_t}{x_t}\delta\left(\frac{z}{x_t}\right)=
    \delta\left(\frac{z}{x_t}\right)\cdot \Res_{z=x_t} \frac{\nu_t}{z-x_t}\frac{dz}{z}.
\end{equation*}
The lemma is proved.
\end{proof}

Since $\Psi(z)$ is a rational function in $z$, which has (simple) poles only at
$\{\sw_r,\vv^2\sw_r\}_{r=1}^a$ and possibly poles of higher order at $z=0,\infty$,
we can apply Lemma~\ref{key equality} to evaluate $\Psi(z)^+-\Psi(z)^-$:
\begin{equation*}
\begin{split}
  & \Psi(z)^+-\Psi(z)^-=\prod_{t=1}^a \sw_t \cdot \sum_{r=1}^a
    \left(\delta\left(\frac{z}{\sw_r}\right)\frac{Z(\sw_r)}{W_r(\sw_r)W(\vv^{-2}\sw_r)}+
          \delta\left(\frac{z}{\vv^2\sw_r}\right)\frac{Z(\vv^2\sw_r)}{W_r(\sw_r)W(\vv^2\sw_r)}\right)=\\
  & \frac{1}{1-\vv^2}\prod_{t=1}^a \sw_t\cdot \sum_{r=1}^a
    \left(\delta\left(\frac{\sw_r}{z}\right)\frac{Z(\sw_r)}{W_r(\sw_r)W_r(\vv^{-2}\sw_r)}-
          \vv^2\delta\left(\frac{\vv^2\sw_r}{z}\right)\frac{Z(\vv^2\sw_r)}{W_r(\sw_r)W_r(\vv^2\sw_r)}\right).
\end{split}
\end{equation*}
Hence, the RHS equals
\begin{equation*}
\begin{split}
  & \frac{1}{(\vv-\vv^{-1})(1-\vv^2)}\prod_{t=1}^a \sw_t\times\\
  & \sum_{r=1}^a
    \left(\delta\left(\frac{\sw_r}{z}\right)\delta\left(\frac{\sw_r}{w}\right)
          \frac{Z(\sw_r)}{W_r(\sw_r)W_r(\vv^{-2}\sw_r)}-\delta\left(
    \frac{\vv^2\sw_r}{z}\right)\delta\left(\frac{\vv^2\sw_r}{w}\right)
          \frac{\vv^2Z(\sw_r)}{W_r(\sw_r)W_r(\vv^2\sw_r)}\right).
\end{split}
\end{equation*}
As a result, we finally get LHS$=$RHS.

\medskip
\noindent
$\circ$ \emph{Case $c_{ij}<0, i\to j$.}
We may assume $I=\{i,j\}$, and we need to check
$[E_i(z),F_j(w)]=0$. We have
  $$\frac{[E_i(z), F_j(w)]}{\rho^+_i\rho^-_j}=\prod_{t=1}^{a_i} \sw_{i,t}\cdot
    \left[\sum_{r=1}^{a_i}\delta\left(\frac{\sw_{i,r}}{z}\right)
          \frac{Z_i(\sw_{i,r})}{W_{i,r}(\sw_{i,r})}D^{-1}_{i,r},
          \sum_{s=1}^{a_j}\delta\left(\frac{\vv_j^2\sw_{j,s}}{w}\right)
          \frac{1}{W_{j,s}(\sw_{j,s})}D_{j,s}\right].$$
The latter is obviously zero, since
  $[D^{-1}_{i,r},\sw_{j,s}]=0=[D_{j,s},\sw_{i,r}]$.

\medskip
\noindent
$\circ$ \emph{Case $c_{ij}<0, i\leftarrow j$.}
We may assume $I=\{i,j\}$, and we need to check
$E_i(z)F_j(w)/(\rho^+_i\rho^-_j)=F_j(w)E_i(z)/(\rho^+_i\rho^-_j)$. The LHS equals
\begin{equation*}
\begin{split}
  & \vv_i^{-c_{ij}}\prod_{t=1}^{a_i}\sw^{1+c_{ij}/2}_{i,t}
  \prod_{t=1}^{a_j}\sw^{c_{ji}/2}_{j,t}\times\\
  & \sum_{1\leq r\leq a_i}^{1\leq s\leq a_j}
  \delta\left(\frac{\sw_{i,r}}{z}\right)\delta\left(\frac{\vv_j^2\sw_{j,s}}{w}\right)
  \frac{Z_i(\sw_{i,r})\prod_{p=1}^{-c_{ji}}W_j(\vv_j^{-c_{ji}-2p}z)}{W_{i,r}(\sw_{i,r})}D_{i,r}^{-1}
  \frac{\prod_{p'=1}^{-c_{ij}}W_i(\vv_i^{-c_{ij}-2p'}w)}{W_{j,s}(\sw_{j,s})}D_{j,s}=\\
  & \prod_{t=1}^{a_i}\sw^{1+c_{ij}/2}_{i,t} \prod_{t=1}^{a_j}\sw^{c_{ji}/2}_{j,t}\cdot A(z,w)\times\\
  & \sum_{1\leq r\leq a_i}^{1\leq s\leq a_j}
  \delta\left(\frac{\sw_{i,r}}{z}\right)\delta\left(\frac{\vv_j^2\sw_{j,s}}{w}\right)
  \frac{Z_i(\sw_{i,r})\prod_{p=1}^{-c_{ji}}W_{j,s}(\vv_j^{-c_{ji}-2p}z)\prod_{p'=1}^{-c_{ij}}W_{i,r}(\vv_i^{-c_{ij}-2p'}w)}
       {W_{i,r}(\sw_{i,r})W_{j,s}(\sw_{j,s})}
  D_{i,r}^{-1}D_{j,s},
\end{split}
\end{equation*}
where
  $A(z,w)=\vv_i^{-c_{ij}}
   \prod_{p=1}^{-c_{ji}}\left(1-\frac{\vv_j^{-2}w}{\vv_j^{-c_{ji}-2p}z}\right)
   \prod_{p'=1}^{-c_{ij}}\left(1-\frac{\vv_i^{-2}z}{\vv_i^{-c_{ij}-2p'}w}\right)$.

Likewise, the RHS equals
\begin{equation*}
\begin{split}
  & \vv_j^{c_{ji}}\prod_{t=1}^{a_i}\sw^{1+c_{ij}/2}_{i,t}
  \prod_{t=1}^{a_j}\sw^{c_{ji}/2}_{j,t}\times\\
  & \sum_{1\leq r\leq a_i}^{1\leq s\leq a_j}
  \delta\left(\frac{\vv_j^2\sw_{j,s}}{w}\right)\delta\left(\frac{\sw_{i,r}}{z}\right)
  \frac{\prod_{p'=1}^{-c_{ij}}W_i(\vv_i^{-c_{ij}-2p'}w)}{W_{j,s}(\sw_{j,s})}D_{j,s}
  \frac{Z_i(\sw_{i,r})\prod_{p=1}^{-c_{ji}}W_j(\vv_j^{-c_{ji}-2p}z)}{W_{i,r}(\sw_{i,r})}D_{i,r}^{-1}=\\
  & \prod_{t=1}^{a_i}\sw^{1+c_{ij}/2}_{i,t} \prod_{t=1}^{a_j}\sw^{c_{ji}/2}_{j,t}\cdot B(z,w)\times\\
  & \sum_{1\leq r\leq a_i}^{1\leq s\leq a_j}
  \delta\left(\frac{\sw_{i,r}}{z}\right)\delta\left(\frac{\vv_j^2\sw_{j,s}}{w}\right)
  \frac{Z_i(\sw_{i,r})\prod_{p=1}^{-c_{ji}}W_{j,s}(\vv_j^{-c_{ji}-2p}z)\prod_{p'=1}^{-c_{ij}}W_{i,r}(\vv_i^{-c_{ij}-2p'}w)}
       {W_{i,r}(\sw_{i,r})W_{j,s}(\sw_{j,s})}
  D_{i,r}^{-1}D_{j,s},
\end{split}
\end{equation*}
where
  $B(z,w)=\vv_j^{c_{ji}}
  \prod_{p=1}^{-c_{ji}}\left(1-\frac{w}{\vv_j^{-c_{ji}-2p}z}\right)
  \prod_{p'=1}^{-c_{ij}}\left(1-\frac{z}{\vv_i^{-c_{ij}-2p'}w}\right)$.

The equality LHS$=$RHS follows from $A(z,w)=B(z,w)$.

\subsection{Compatibility with~(\ref{U7}).}
\

\medskip
\noindent
$\circ$ \emph{Case $c_{ij}=0$}.
In this case, $[E_i(z),E_j(w)]=0$, due to our verification of~(\ref{U2}).

\medskip
\noindent
$\circ$ \emph{Case $c_{ij}<0$.}
To simplify our calculations, we introduce
\begin{equation*}
    \chi_{i',r}:=\prod_{t=1}^{a_{i'}}\sw_{i',t}\cdot \prod_{j'\to i'} \prod_{t=1}^{a_{j'}} \sw_{j',t}^{c_{j'i'}/2}
    \cdot\frac{Z_{i'}(\sw_{i',r})}{W_{i',r}(\sw_{i',r})}
    \prod_{j'\to i'} \prod_{p=1}^{-c_{j'i'}}W_{j'}(\vv_{j'}^{-c_{j'i'}-2p}\sw_{i',r})D_{i',r}^{-1},
\end{equation*}
so that
  $E_{i'}(z)=\rho^+_{i'}\sum_{r=1}^{a_{i'}}\delta\left(\frac{\sw_{i',r}}{z}\right)\chi_{i',r}$.

The verification of~(\ref{U7}) is based on the following result.

\begin{Lem}\label{Properties for Serre}
The following relations hold:
\begin{equation*}
    \chi_{i,r}\sw_{j,s}=\vv_{i}^{-2\delta_{ij}\delta_{rs}}\sw_{j,s}\chi_{i,r}
    \ \mathrm{for}\ 1\leq r\leq a_{i}, 1\leq s\leq a_{j},
\end{equation*}
\begin{equation*}
    (\sw_{i,r_1}-\vv_i^2\sw_{i,r_2})\chi_{i,r_1}\chi_{i,r_2}=(\vv_i^2\sw_{i,r_1}-\sw_{i,r_2})\chi_{i,r_2}\chi_{i,r_1}
    \ \mathrm{for}\ 1\leq r_1\ne r_2 \leq a_i,
\end{equation*}
\begin{equation*}
    (\sw_{i,r}-\vv_i^{c_{ij}}\sw_{j,s})\chi_{i,r}\chi_{j,s}=(\vv_i^{c_{ij}}\sw_{i,r}-\sw_{j,s})\chi_{j,s}\chi_{i,r}
    \ \mathrm{for}\ 1\leq r\leq a_i, 1\leq s\leq a_j.
\end{equation*}
\end{Lem}

\begin{proof}
Follows from straightforward computations.
\end{proof}

With the help of this lemma, let us verify~(\ref{U7}) for $c_{ij}=-1$.
The latter amounts to proving
  $[E_i(z_1),[E_i(z_2),E_j(w)]_\vv]_{\vv^{-1}}/((\rho^+_i)^2\rho^+_j)=
   -[E_i(z_2),[E_i(z_1),E_j(w)]_\vv]_{\vv^{-1}}/((\rho^+_i)^2\rho^+_j).$
The LHS equals
\begin{equation*}
\begin{split}
  & (1-\vv^2)
    \left[\sum_{r_1=1}^{a_i} \delta\left(\frac{\sw_{i,r_1}}{z_1}\right)\chi_{i,r_1},
          \sum_{1\leq r_2\leq a_i}^{1\leq s\leq a_j}
          \delta\left(\frac{\sw_{i,r_2}}{z_2}\right)\delta\left(\frac{\sw_{j,s}}{w}\right)
          \frac{\sw_{i,r_2}}{\sw_{i,r_2}-\vv\sw_{j,s}}\chi_{i,r_2}\chi_{j,s}\right]_{\vv^{-1}}=\\
  & \sum_{1\leq r\leq a_i}^{1\leq s\leq a_j}
    \delta\left(\frac{\sw_{j,s}}{w}\right)
    \left\{\delta\left(\frac{\sw_{i,r}}{z_1}\right)\delta\left(\frac{\vv^{-2}\sw_{i,r}}{z_2}\right)-
           \delta\left(\frac{\sw_{i,r}}{z_2}\right)\delta\left(\frac{\vv^{-2}\sw_{i,r}}{z_1}\right)\right\}
    \frac{(\vv^2-1)\sw_{i,r}}{\sw_{i,r}-\vv^3\sw_{j,s}}\chi^2_{i,r}\chi_{j,s}-\\
  & (\vv^2-1)^2\sum_{1\leq r_1\ne r_2\leq a_i}^{1\leq s\leq a_j}
    \delta\left(\frac{\sw_{i,r_1}}{z_1}\right)\delta\left(\frac{\sw_{i,r_2}}{z_2}\right)\delta\left(\frac{\sw_{j,s}}{w}\right)
    \frac{A(z_1,z_2,w)}{\vv^2\sw_{i,r_1}-\sw_{i,r_2}}
    \chi_{i,r_1}\chi_{i,r_2}\chi_{j,s},
\end{split}
\end{equation*}
where
  $A(z_1,z_2,w)=\frac{z_1z_2(z_1+z_2-(\vv+\vv^{-1})w)}{(z_1-\vv w)(z_2-\vv w)}$
and the last equality is obtained by treating separately $r_1=r_2$
and $r_1\ne r_2$ cases. The first sum is obviously skew-symmetric in
$z_1,z_2$. The second sum is also skew-symmetric, due to the above
relations on $\chi_{i,r}$.

\medskip
The cases $c_{ij}=-2, -3$ can be treated similarly, but the
corresponding computations become more cumbersome. We verified these
cases using \textbf{MATLAB}.

\subsection{Compatibility with~(\ref{U8})}
\

\medskip
The case $c_{ij}=0$ is obvious. The case $c_{ij}=-1$ can be treated
analogously to the above verification of (U7). The verification for
the cases $c_{ij}=-2,-3$ is more cumbersome and can be performed as
outlined in the verification of (U7). Our verification involved a
simple computation in \textbf{MATLAB}.


\medskip
This completes our proof of Theorem~\ref{Homomorphism}.

\begin{Rem}
Theorem~\ref{Homomorphism} admits the following straightforward
generalization. For every $i\in I$, pick two polynomials
$Z^{(1)}_i(z),Z^{(2)}_i(z)$ in $z^{-1}$ such that
$Z_i(z)=Z^{(1)}_i(z)Z^{(2)}_i(z)$. There is a unique
$\BC(\vv)[\sz_1^{\pm1},\ldots,\sz_N^{\pm1}]$-algebra homomorphism
  $\CU^\ad_{0,\mu}[\sz^{\pm 1}_1,\ldots,\sz^{\pm 1}_N]
   \to \wt{\CA}^\vv_\fra[\sz^{\pm 1}_1,\ldots,\sz^{\pm 1}_N]$,
such that
  $$e_i(z)\mapsto \frac{-\vv_i}{1-\vv_i^2}\prod_{t=1}^{a_i}\sw_{i,t} \prod_{j\to i} \prod_{t=1}^{a_j} \sw_{j,t}^{c_{ji}/2}\cdot
    \sum_{r=1}^{a_i} \delta\left(\frac{\sw_{i,r}}{z}\right)\frac{Z^{(1)}_i(\sw_{i,r})}{W_{i,r}(\sw_{i,r})}
    \prod_{j\to i}\prod_{p=1}^{-c_{ji}} W_j(\vv_j^{-c_{ji}-2p}z)D_{i,r}^{-1},$$
  $$f_i(z)\mapsto \frac{1}{1-\vv_i^2}\prod_{j\leftarrow i}\prod_{t=1}^{a_j} \sw_{j,t}^{c_{ji}/2}\cdot
    \sum_{r=1}^{a_i} \delta\left(\frac{\vv_i^2\sw_{i,r}}{z}\right)\frac{Z^{(2)}_i(\vv_i^2\sw_{i,r})}{W_{i,r}(\sw_{i,r})}
    \prod_{j\leftarrow i}\prod_{p=1}^{-c_{ji}} W_j(\vv_j^{-c_{ji}-2p}z)D_{i,r},$$
  $$\psi^\pm_i(z)\mapsto\Psi_i(z)^\pm,\
    (\phi^+_i)^{\pm 1}\mapsto \prod_{t=1}^{a_i} \sw_{i,t}^{\pm 1/2},\
    (\phi^-_i)^{\pm 1}\mapsto (-\vv_i)^{\mp a_i}\prod_{t=1}^{a_i} \sw_{i,t}^{\mp 1/2}.$$
\end{Rem}



\section{Proof of Theorem~\ref{Coproduct for sl2}}
\label{Proof Coproduct}

Due to Theorem~\ref{Levendorskii-type}, it suffices to check that
the assignment $\Delta$ of Theorem~\ref{Coproduct for sl2}
preserves defining relations (\ref{u1}--\ref{u6}, \ref{u9}).
To simplify our exposition, we will assume that $b_1,b_2<0$, while the
case when one of them is zero is left to the interested reader (note
that the case $b_1=b_2=0$ has been treated in~Remark~\ref{Explicit sl2}).
We will also work with $\bar{h}_{\pm 1}:=[2]_\vv^{-1}h_{\pm 1}$
instead of $h_{\pm 1}$, so that
  $[\bar{h}_{\pm 1},e_r]=e_{r\pm 1}, [\bar{h}_{\pm 1},f_r]=-f_{r\pm 1}$.


\subsection{Compatibility with (\ref{u1})}
\

The equalities
  $\Delta((\psi^+_0)^{\pm 1})\Delta((\psi^+_0)^{\mp 1})=1$
and
  $\Delta((\psi^-_{b})^{\pm 1})\Delta((\psi^-_{b})^{\mp 1})=1$
follow immediately from relation~(\ref{u1}) for both
$\CU^\ssc_{0,b_1}$ and $\CU^\ssc_{0,b_2}$.

The commutativity of
  $\Delta((\psi^+_0)^{\pm 1}), \Delta((\psi^-_{b})^{\pm 1})$
between themselves and with each of $\Delta(\bar{h}_{\pm 1})$ is due
to relations~(\ref{u1},~\ref{u4},~\ref{u5}) for both $\CU^\ssc_{0,b_1}$
and $\CU^\ssc_{0,b_2}$.

It remains to prove $[\Delta(\bar{h}_1), \Delta(\bar{h}_{-1})]=0$.
The LHS is equal to

\noindent
  $[\bar{h}_1\otimes 1+1\otimes \bar{h}_1-(\vv-\vv^{-1})e_0\otimes f_1,
   \bar{h}_{-1}\otimes 1+1\otimes \bar{h}_{-1}+(\vv-\vv^{-1})e_{-1}\otimes f_0]=
   (\vv-\vv^{-1})(e_0\otimes f_0-e_{-1}\otimes f_1+e_{-1}\otimes f_1-e_0\otimes f_0)
   -(\vv-\vv^{-1})^2 [e_0\otimes f_1, e_{-1}\otimes f_0]=
   -(\vv-\vv^{-1})^2 (e_0e_{-1}\otimes f_1f_0-e_{-1}e_0\otimes f_0f_1)=0.$

Here we used~(\ref{u1},~\ref{u4},~\ref{u5}) for both $\CU^\ssc_{0,b_1}, \CU^\ssc_{0,b_2}$
in the first equality, while the second equality follows
from  $e_0e_{-1}=\vv^2e_{-1}e_0, f_1f_0=\vv^{-2}f_0f_1$,
due to~(\ref{u2}) for $\CU^\ssc_{0,b_1}$ and~(\ref{u3}) for $\CU^\ssc_{0,b_2}$.


\subsection{Compatibility with (\ref{u2})}
\

We need to prove
  $[\Delta(e_{r+1}),\Delta(e_s)]_{\vv^2}+[\Delta(e_{s+1}),\Delta(e_r)]_{\vv^2}=0$
for $b_2-1\leq r,s\leq -1$.

\medskip
\noindent
$\circ$ \emph{Case $b_2-1<r,s<-1$.}
Then,
  $[\Delta(e_{r+1}),\Delta(e_s)]_{\vv^2}+[\Delta(e_{s+1}),\Delta(e_r)]_{\vv^2}=
   1\otimes ([e_{r+1},e_s]_{\vv^2}+[e_{s+1},e_r]_{\vv^2})=0$
as the second term is zero in $\CU^\ssc_{0,b_2}$ by (\ref{u2}).

\medskip
\noindent
$\circ$ \emph{Case $r=s=b_2-1$.}
It suffices to show that
  $[\Delta(e_{b_2}),\Delta(e_{b_2-1})]_{\vv^2}=0$,
which follows from
  $[\Delta(e_{b_2}),\Delta(e_{b_2-1})]_{\vv^2}=
   [1\otimes e_{b_2}, e_{-1}\otimes \psi^-_{b_2}+1\otimes e_{b_2-1}]_{\vv^2}=
   e_{-1}\otimes [e_{b_2},\psi^-_{b_2}]_{\vv^2}+1\otimes[e_{b_2},e_{b_2-1}]_{\vv^2}=0.$
The last equality follows from $[e_{b_2},\psi^-_{b_2}]_{\vv^2}=0$
and $[e_{b_2},e_{b_2-1}]_{\vv^2}=0$ in $\CU^\ssc_{0,b_2}$, due to
(\ref{u2}) and (\ref{u4}), respectively.

\medskip
\noindent
$\circ$ \emph{Case $r=b_2-1, b_2-1<s<-1$.}
Then,
  $[\Delta(e_{b_2}),\Delta(e_s)]_{\vv^2}+[\Delta(e_{s+1}),\Delta(e_{b_2-1})]_{\vv^2}=
   1\otimes ([e_{b_2},e_s]_{\vv^2}+[e_{s+1},e_{b_2-1}]_{\vv^2})+e_{-1}\otimes [e_{s+1},\psi^-_{b_2}]_{\vv^2}=0.$
The last equality follows again from (\ref{u2}) and (\ref{u4}) for $\CU^\ssc_{0,b_2}$.

\medskip
\noindent
$\circ$ \emph{Case $r=b_2-1, s=-1$.}
Then
  $[\Delta(e_{b_2}),\Delta(e_{-1})]_{\vv^2}=1\otimes [e_{b_2}, e_{-1}]_{\vv^2}$
and
  $[\Delta(e_0),\Delta(e_{b_2-1})]_{\vv^2}=
   [e_0\otimes \psi^+_0+1\otimes e_0, e_{-1}\otimes \psi^-_{b_2}+1\otimes e_{b_2-1}]_{\vv^2}=
   e_0\otimes [\psi^+_0, e_{b_2-1}]_{\vv^2}+[e_0,e_{-1}]_{\vv^2}\otimes \psi^+_0\psi^-_{b_2}+
   e_{-1}\otimes [e_0,\psi^-_{b_2}]_{\vv^2}+1\otimes [e_0,e_{b_2-1}]_{\vv^2}=
   1\otimes [e_0,e_{b_2-1}]_{\vv^2}$
as the first three terms are zero, due to~(\ref{u2}) for $\CU^\ssc_{0,b_1}$
and~(\ref{u4}) for $\CU^\ssc_{0,b_2}$. The
result follows from (\ref{u2}) for $\CU^\ssc_{0,b_2}$.

\medskip
\noindent
$\circ$ \emph{Case $r=s=-1$.}
It suffices to show that
  $[\Delta(e_0),\Delta(e_{-1})]_{\vv^2}=0$,
which follows from
  $[\Delta(e_0),\Delta(e_{-1})]_{\vv^2}=
   [e_0\otimes \psi^+_0 + 1\otimes e_0, 1\otimes e_{-1}]_{\vv^2}=
   e_0\otimes [\psi^+_0, e_{-1}]_{\vv^2}+ 1\otimes [e_0,e_{-1}]_{\vv^2}=0.$
The last equality follows again from relations~(\ref{u2},~\ref{u4}) for the
algebra $\CU^\ssc_{0,b_2}$.

\medskip
\noindent
$\circ$ \emph{Case $r=-1, b_2-1<s<-1$.}
Then,
  $[\Delta(e_0), \Delta(e_s)]_{\vv^2}=
   [e_0\otimes \psi^+_0 + 1\otimes e_0, 1\otimes e_s]_{\vv^2}=
   1\otimes [e_0,e_s]_{\vv^2}$,
while
  $[\Delta(e_{s+1}),\Delta(e_{-1})]_{\vv^2}=1\otimes [e_{s+1},e_{-1}]_{\vv^2}.$
The sum of these two terms is zero, due to~(\ref{u2}) for $\CU^\ssc_{0,b_2}$.


\subsection{Compatibility with (\ref{u3})}
\

We need to prove
  $[\Delta(f_r),\Delta(f_{s+1})]_{\vv^2}+[\Delta(f_s),\Delta(f_{r+1})]_{\vv^2}=0$
for $b_1\leq r,s\leq 0$.

\medskip
\noindent
$\circ$ \emph{Case $b_1<r,s<0$.}
Then,
  $[\Delta(f_r),\Delta(f_{s+1})]_{\vv^2}+[\Delta(f_s),\Delta(f_{r+1})]_{\vv^2}=
   ([f_{r},f_{s+1}]_{\vv^2}+[f_{s},f_{r+1}]_{\vv^2})\otimes 1=0$
as the first term is zero in $\CU^\ssc_{0,b_1}$ by (\ref{u3}).

\medskip
\noindent
$\circ$ \emph{Case $r=s=b_1$.}
It suffices to show that
  $[\Delta(f_{b_1}),\Delta(f_{1+b_1})]_{\vv^2}=0$,
which follows from
  $[\Delta(f_{b_1}),\Delta(f_{1+b_1})]_{\vv^2}=
   [f_{b_1}\otimes 1 + \psi^-_{b_1}\otimes f_0, f_{1+b_1}\otimes 1]_{\vv^2}=
   [f_{b_1},f_{1+b_1}]_{\vv^2}\otimes 1+[\psi^-_{b_1},f_{1+b_1}]_{\vv^2}\otimes f_0=0.$
The last equality follows from
$[f_{b_1},f_{1+b_1}]_{\vv^2}=0=[\psi^-_{b_1},f_{1+b_1}]_{\vv^2}$,
due to~(\ref{u3},~\ref{u5}) for $\CU^\ssc_{0,b_1}$.

\medskip
\noindent
 $\circ$ \emph{Case $r=b_1<s<0$.}
Then,
  $[\Delta(f_{s}),\Delta(f_{1+b_1})]_{\vv^2}=[f_s,f_{1+b_1}]_{\vv^2}\otimes 1$
and
  $[\Delta(f_{b_1}),\Delta(f_{s+1})]_{\vv^2}=
   [f_{b_1},f_{s+1}]_{\vv^2}\otimes 1$
as $[\psi^-_{b_1}, f_{s+1}]_{\vv^2}=0$ in $\CU^\ssc_{0,b_1}$ by
(\ref{u5}). It remains to use~(\ref{u3}) for $\CU^\ssc_{0,b_1}$.

\medskip
\noindent
 $\circ$ \emph{Case $r=b_1, s=0$.}
Then
  $[\Delta(f_{b_1}),\Delta(f_1)]_{\vv^2}=
   [f_{b_1}, f_1]_{\vv^2}\otimes 1+[f_{b_1},\psi^+_0]_{\vv^2}\otimes f_1+
   [\psi_{b_1}^-,f_1]_{\vv^2}\otimes f_0+\psi_{b_1}^-\psi_0^+\otimes[f_0,f_1]_{\vv^2}$,
and
  $[\Delta(f_0),\Delta(f_{1+b_1})]_{\vv^2}=[f_0,f_{1+b_1}]_{\vv^2}\otimes 1$.
It remains to use
  $[f_{b_1}, f_1]_{\vv^2}+[f_0, f_{1+b_1}]_{\vv^2}=[f_{b_1},\psi_0^+]_{\vv^2}=
   [\psi_{b_1}^-,f_1]_{\vv^2}=0$ in $\CU^\ssc_{0,b_1}$, due~(\ref{u3}) and~(\ref{u5}),
and $[f_0,f_1]_{\vv^2}=0$ in $\CU^\ssc_{0,b_2}$, due to~(\ref{u3}).

\medskip
\noindent
 $\circ$ \emph{Case $r=s=0$.}
It suffices to show that
  $[\Delta(f_0),\Delta(f_1)]_{\vv^2}=0$,
which follows from
  $[\Delta(f_0),\Delta(f_1)]_{\vv^2}=
   [f_0\otimes 1, f_1\otimes 1 + \psi^+_0\otimes f_1]_{\vv^2}=
   [f_0,f_1]_{\vv^2}\otimes 1 + [f_0, \psi^+_0]_{\vv^2}\otimes f_1=0$,
due to~(\ref{u3},~\ref{u5}) for $\CU^\ssc_{0,b_1}$.

\medskip
\noindent
 $\circ$ \emph{Case $r=0, b_1<s<0$.}
Then
  $[\Delta(f_0), \Delta(f_{s+1})]_{\vv^2}=[f_0,f_{s+1}]_{\vv^2}\otimes 1$,
and
  $[\Delta(f_s), \Delta(f_1)]_{\vv^2}=
   [f_s\otimes 1, f_1\otimes 1 + \psi^+_0\otimes f_1]_{\vv^2}=
   [f_s,f_1]_{\vv^2}\otimes 1+[f_s,\psi^+_0]_{\vv^2}\otimes f_1$.
It remains to apply the equalities
  $[f_0,f_{s+1}]_{\vv^2}+[f_s,f_1]_{\vv^2}=0$
and
  $[f_s,\psi^+_0]_{\vv^2}=0$
in $\CU^\ssc_{0,b_1}$, due to~(\ref{u3}) and~(\ref{u5}).


\subsection{Compatibility with (\ref{u4})}
\

The equalities
  $\Delta(\psi^+_0)\Delta(e_r)=\vv^2\Delta(e_r)\Delta(\psi^+_0)$
and
  $\Delta(\psi^-_{b})\Delta(e_r)=\vv^{-2}\Delta(e_r)\Delta(\psi^-_{b})$
for $b_2-1\leq r\leq 0$ are obvious, due to relations (\ref{u1}) and
(\ref{u4}) for $\CU^\ssc_{0,b_1},\CU^\ssc_{0,b_2}$.

\medskip
Let us now verify the equality
  $[\Delta(\bar{h}_1),\Delta(e_r)]=\Delta(e_{r+1})$
for $b_2-1\leq r\leq -1$.

\medskip
\noindent
$\circ$ \emph{Case $b_2\leq r\leq -2$.}
We have
  $[\Delta(\bar{h}_1), \Delta(e_r)]=
   [\bar{h}_1\otimes 1+1\otimes \bar{h}_1-(\vv-\vv^{-1})e_0\otimes f_1, 1\otimes e_r]=
   1\otimes e_{r+1}-(\vv-\vv^{-1})e_0\otimes [f_1,e_r]=1\otimes e_{r+1}=\Delta(e_{r+1})$, due to~(\ref{u4},~\ref{u6}) for $\CU^\ssc_{0,b_2}$.

\medskip
\noindent
 $\circ$ \emph{Case $r=-1$.}
As above, we get
  $[\Delta(\bar{h}_1), \Delta(e_{-1})]=
   [\bar{h}_1\otimes 1+1\otimes \bar{h}_1-(\vv-\vv^{-1})e_0\otimes f_1, 1\otimes e_{-1}]=
   1\otimes e_0-(\vv-\vv^{-1})e_0\otimes [f_1,e_{-1}]=
   1\otimes e_0+(\vv-\vv^{-1})e_0\otimes \frac{\psi^+_0}{\vv-\vv^{-1}}=\Delta(e_0)$.

\medskip
\noindent
$\circ$ \emph{Case $r=b_2-1$.}
We have
  $[\Delta(\bar{h}_1), \Delta(e_{b_2-1})]=
   [\bar{h}_1\otimes 1+1\otimes \bar{h}_1-(\vv-\vv^{-1})e_0\otimes f_1,
   e_{-1}\otimes \psi^-_{b_2}+1\otimes e_{b_2-1}]=
   e_0\otimes \psi^-_{b_2}+1\otimes e_{b_2}-e_0\otimes \psi^-_{b_2}-
   (\vv-\vv^{-1})[e_0\otimes f_1, e_{-1}\otimes \psi^-_{b_2}]=
   1\otimes e_{b_2}=\Delta(e_{b_2}),$
where we used $[e_0\otimes f_1, e_{-1}\otimes \psi^-_{b_2}]=0$
as $e_0e_{-1}=\vv^2e_{-1}e_0$ in $\CU^\ssc_{0,b_1}$, due to~(\ref{u2}),
and $\psi^-_{b_2}f_1=\vv^2f_1\psi^-_{b_2}$ in $\CU^\ssc_{0,b_2}$, due to~(\ref{u5}).

\medskip
Let us now verify the equality
  $[\Delta(\bar{h}_{-1}),\Delta(e_r)]=\Delta(e_{r-1})$
for $b_2\leq r\leq 0$.

\medskip
\noindent
$\circ$ \emph{Case $b_2<r<0$.}
We have
  $[\Delta(\bar{h}_{-1}), \Delta(e_r)]=
   [\bar{h}_{-1}\otimes 1+1\otimes \bar{h}_{-1}+(\vv-\vv^{-1})e_{-1}\otimes f_0, 1\otimes e_r]=
   1\otimes e_{r-1}+(\vv-\vv^{-1})e_{-1}\otimes [f_0, e_r]=1\otimes e_{r-1}=\Delta(e_{r-1})$, due to~(\ref{u4},~\ref{u6}) for $\CU^\ssc_{0,b_2}$.

\medskip
\noindent
$\circ$ \emph{Case $r=0$.}
We have
  $[\Delta(\bar{h}_{-1}), \Delta(e_0)]=
   [\bar{h}_{-1}\otimes 1+1\otimes \bar{h}_{-1}+(\vv-\vv^{-1})e_{-1}\otimes f_0, e_0\otimes \psi^+_0 + 1\otimes e_0]=
   e_{-1}\otimes \psi^+_0 + 1\otimes e_{-1} + (\vv-\vv^{-1})e_{-1}\otimes [f_0,e_0]+
   (\vv-\vv^{-1})[e_{-1}\otimes f_0, e_0\otimes \psi^+_0]=
   1\otimes e_{-1}=\Delta(e_{-1}),$
where we used $[e_{-1}\otimes f_0, e_0\otimes \psi^+_0]=0$ as
$e_0e_{-1}=\vv^2e_0e_{-1}$ in $\CU^\ssc_{0,b_1}$, due to~(\ref{u2}), and
$f_0\psi^+_0=\vv^2\psi^+_0f_0$ in $\CU^\ssc_{0,b_2}$, due to~(\ref{u5}).

\medskip
\noindent
$\circ$ \emph{Case $r=b_2$.}
We have
  $[\Delta(\bar{h}_{-1}), \Delta(e_{b_2})]=
   [\bar{h}_{-1}\otimes 1+1\otimes \bar{h}_{-1}+(\vv-\vv^{-1})e_{-1}\otimes f_0, 1\otimes e_{b_2}]=
   1\otimes e_{b_2-1}+(\vv-\vv^{-1})e_{-1}\otimes \frac{\psi^-_{b_2}}{\vv-\vv^{-1}}=
   \Delta(e_{b_2-1})$, due to~(\ref{u4},~\ref{u6}) for $\CU^\ssc_{0,b_2}$.


\subsection{Compatibility with (\ref{u5})}
\

The equalities
  $\Delta(\psi^+_0)\Delta(f_r)=\vv^{-2}\Delta(f_r)\Delta(\psi^+_0)$
and
  $\Delta(\psi^-_{b})\Delta(f_r)=\vv^2\Delta(f_r)\Delta(\psi^-_{b})$
for $b_1\leq r\leq 1$ are obvious, due to relations (\ref{u1})
and~(\ref{u5}) for $\CU^\ssc_{0,b_1},\CU^\ssc_{0,b_2}$.

\medskip
Let us now verify the equality
  $[\Delta(\bar{h}_1),\Delta(f_r)]=-\Delta(f_{r+1})$
for $b_1\leq r\leq 0$.

\medskip
\noindent
 $\circ$ \emph{Case $b_1<r<0$.}
We have
  $[\Delta(\bar{h}_1), \Delta(f_r)]=
   [\bar{h}_1\otimes 1+1\otimes \bar{h}_1-(\vv-\vv^{-1})e_0\otimes f_1, f_r\otimes 1]=
   -f_{r+1}\otimes 1-(\vv-\vv^{-1})[e_0,f_r]\otimes f_1=-f_{r+1}\otimes 1=
   -\Delta(f_{r+1})$, due to~(\ref{u5},~\ref{u6}) for $\CU^\ssc_{0,b_1}$.

\medskip
\noindent
 $\circ$ \emph{Case $r=0$.}
As above, we get
  $[\Delta(\bar{h}_1), \Delta(f_0)]=
   [\bar{h}_1\otimes 1+1\otimes \bar{h}_1-(\vv-\vv^{-1})e_0\otimes f_1, f_0\otimes 1]=
   -f_1\otimes 1 - (\vv-\vv^{-1})[e_0,f_0]\otimes f_1=
   -f_1\otimes 1 - \psi^+_0\otimes f_1=-\Delta(f_1)$.

\medskip
\noindent
 $\circ$ \emph{Case $r=b_1$.}
We have
  $[\Delta(\bar{h}_1), \Delta(f_{b_1})]=
   [\bar{h}_1\otimes 1+1\otimes \bar{h}_1-(\vv-\vv^{-1})e_0\otimes f_1, f_{b_1}\otimes 1 + \psi^-_{b_1}\otimes f_0]=
   -f_{1+b_1}\otimes 1 - \psi^-_{b_1}\otimes f_1 + \psi^-_{b_1}\otimes f_1-
   (\vv-\vv^{-1})[e_0\otimes f_1, \psi^-_{b_1}\otimes f_0]=-f_{1+b_1}\otimes 1=
   -\Delta(f_{1+b_1}),$
where we used $[e_0\otimes f_1, \psi^-_{b_1}\otimes f_0]=0$ as
$f_1f_0=\vv^{-2}f_0f_1$ in $\CU^\ssc_{0,b_2}$, due to~(\ref{u3}), and
$\psi^-_{b_1}e_0=\vv^{-2}e_0\psi^-_{b_1}$ in $\CU^\ssc_{0,b_1}$, due to~(\ref{u4}).

\medskip
Let us now verify the equality
  $[\Delta(\bar{h}_{-1}),\Delta(f_r)]=-\Delta(f_{r-1})$
for $1+b_1\leq r\leq 1$.

\medskip
\noindent
 $\circ$ \emph{Case $1+b_1<r<1$.}
We have
  $[\Delta(\bar{h}_{-1}), \Delta(f_r)]=
   [\bar{h}_{-1}\otimes 1+1\otimes \bar{h}_{-1}+(\vv-\vv^{-1})e_{-1}\otimes f_0,f_r\otimes 1]=
   -f_{r-1}\otimes 1+(\vv-\vv^{-1})[e_{-1},f_r]\otimes f_0=
   -f_{r-1}\otimes 1=-\Delta(f_{r-1})$,
due to~(\ref{u5},~\ref{u6}) for $\CU^\ssc_{0,b_1}$.

\medskip
\noindent
$\circ$ \emph{Case $r=1$.}
We have
  $[\Delta(\bar{h}_{-1}), \Delta(f_1)]=
   [\bar{h}_{-1}\otimes 1+1\otimes \bar{h}_{-1}+(\vv-\vv^{-1})e_{-1}\otimes f_0, f_1\otimes 1+ \psi^+_0\otimes f_1]=
   -f_0\otimes 1 - \psi^+_0\otimes f_0 + \psi^+_0\otimes f_0+(\vv-\vv^{-1})[e_{-1}\otimes f_0, \psi^+_0\otimes f_1]=
   -f_0\otimes 1=-\Delta(f_0),$
where we used $[e_{-1}\otimes f_0, \psi^+_0\otimes f_1]=0$ as
$f_0f_1=\vv^2f_1f_0$ in $\CU^\ssc_{0,b_2}$ and $\psi^+_0e_{-1}=\vv^2e_{-1}\psi^+_0$
in $\CU^\ssc_{0,b_1}$.

\medskip
\noindent
$\circ$ \emph{Case $r=1+b_1$.}
We have
  $[\Delta(\bar{h}_{-1}), \Delta(f_{1+b_1})]=
   [\bar{h}_{-1}\otimes 1+1\otimes \bar{h}_{-1}+(\vv-\vv^{-1})e_{-1}\otimes f_0, f_{1+b_1}\otimes 1]=
   -f_{b_1}\otimes1-\psi_{b_1}^-\otimes f_0=-\Delta(f_{b_1})$,
due to~(\ref{u5},~\ref{u6}) for $\CU^\ssc_{0,b_1}$.


\subsection{Compatibility with (\ref{u6})}
\

\medskip
\noindent
$\circ$ \emph{Case $b_2\leq r<0, b_1<s\leq 0$.}
The equality $[\Delta(e_r), \Delta(f_s)]=0$ is obvious.

\medskip
\noindent
$\circ$ \emph{Case $r=s=0$.}
We need to prove
  $[\Delta(e_0),\Delta(f_0)]=\frac{1}{\vv-\vv^{-1}}\Delta(\psi^+_0)$.
This follows from
  $[\Delta(e_0),\Delta(f_0)]=[e_0\otimes \psi^+_0+1\otimes e_0, f_0\otimes 1]=
   [e_0,f_0]\otimes \psi^+_0=\frac{\psi^+_0\otimes \psi^+_0}{\vv-\vv^{-1}}=
   \frac{\Delta(\psi^+_0)}{\vv-\vv^{-1}}$,
due to~(\ref{u6}) for $\CU^\ssc_{0,b_1}$.

\medskip
\noindent
$\circ$ \emph{Case $r=0,s=1$.}
We need to prove
  $[\Delta_{b_1,b_2}(e_0),\Delta_{b_1,b_2}(f_1)]=\Delta_{b_1,b_2}(\psi^+_0)\Delta_{b_1,b_2}(h_1).$
This can be easily deduced from the unshifted case $b_1=b_2=0$ by
applying Remark~\ref{Maps jmath}. Indeed,
  $[\Delta_{b_1,b_2}(e_0),\Delta_{b_1,b_2}(f_1)]=
   [\jmath^+_{b_1,0}\otimes \jmath^+_{0,b_2}(\Delta(e_0)),\jmath^+_{b_1,0}\otimes \jmath^+_{0,b_2}(\Delta(f_1))]=
   \jmath^+_{b_1,0}\otimes \jmath^+_{0,b_2}(\Delta([e_0,f_1]))=
   \jmath^+_{b_1,0}\otimes \jmath^+_{0,b_2}(\Delta(\psi^+_0)\Delta(h_1))=
   \Delta_{b_1,b_2}(\psi^+_0)\Delta_{b_1,b_2}(h_1),$
where the subscripts in $\Delta_{b_1,b_2}$ are used this time to distinguish it from the
Drinfeld-Jimbo coproduct $\Delta$.

\medskip
\noindent
$\circ$ \emph{Case $r=0,b_1<s<0$.}
We need to prove $[\Delta(e_0),\Delta(f_s)]=0$. This follows from
  $[\Delta(e_0),\Delta(f_s)]=[e_0\otimes \psi^+_0+1\otimes e_0, f_s\otimes 1]=
   [e_0,f_s]\otimes \psi^+_0=0$
as $[e_0,f_s]=0$ in $\CU^\ssc_{0,b_1}$ by~(\ref{u6}).

\medskip
\noindent
$\circ$ \emph{Case $r=0,s=b_1$.}
We need to prove $[\Delta(e_0),\Delta(f_{b_1})]=0$. This follows from
  $[\Delta(e_0),\Delta(f_{b_1})]=
   [e_0\otimes \psi^+_0+1\otimes e_0, f_{b_1}\otimes 1+\psi^-_{b_1}\otimes f_0]=
   [e_0,f_{b_1}]\otimes \psi^+_0+\psi^-_{b_1}\otimes [e_0,f_0]=
   -\frac{\psi^-_{b_1}\otimes \psi^+_0}{\vv-\vv^{-1}}+\frac{\psi^-_{b_1}\otimes \psi^+_0}{\vv-\vv^{-1}}=0,$
where we used
  $[e_0\otimes \psi^+_0,\psi^-_{b_1}\otimes f_0]=0$
as
 $\psi^+_0f_0=\vv^{-2}f_0\psi^+_0$ in $\CU^\ssc_{0,b_2},\
  \psi^-_{b_1}e_0=\vv^{-2}e_0\psi^-_{b_1}$ in $\CU^\ssc_{0,b_1}$.

\medskip
\noindent
$\circ$ \emph{Case $r=-1,s=1$.}
We need to prove
  $[\Delta(e_{-1}),\Delta(f_1)]=\frac{1}{\vv-\vv^{-1}}\Delta(\psi^+_0)$.
This follows from
  $[\Delta(e_{-1}),\Delta(f_1)]=[1\otimes e_{-1}, f_1\otimes 1+\psi^+_0\otimes f_1]=
   \psi^+_0\otimes [e_{-1},f_1]=\frac{\psi^+_0\otimes \psi^+_0}{\vv-\vv^{-1}}=
   \frac{\Delta(\psi^+_0)}{\vv-\vv^{-1}}$,
due to~(\ref{u6}) for $\CU^\ssc_{0,b_2}$.

\medskip
\noindent
$\circ$ \emph{Case $b_2\leq r<-1,s=1$.}
We need to prove $[\Delta(e_r),\Delta(f_1)]=0$. This follows from
  $[\Delta(e_r),\Delta(f_1)]=[1\otimes e_r, f_1\otimes 1+\psi^+_0\otimes f_1]=
   \psi^+_0\otimes [e_r,f_1]=0$ as $[e_r,f_1]=0$ in $\CU^\ssc_{0,b_2}$ by~(\ref{u6}).

\medskip
\noindent
$\circ$ \emph{Case $r=b_2-1, s=1$.}
We need to prove
 $[\Delta(e_{b_2-1}),\Delta(f_1)]=0$.
This follows from
  $[\Delta(e_{b_2-1}),\Delta(f_1)]=
   [e_{-1}\otimes \psi^-_{b_2}+1\otimes e_{b_2-1},f_1\otimes 1 +\psi^+_0\otimes f_1]=
   [e_{-1},f_1]\otimes \psi^-_{b_2}+\psi^+_0\otimes [e_{b_2-1},f_1]+
   [e_{-1}\otimes \psi^-_{b_2},\psi^+_0\otimes f_1]=
   \frac{\psi^+_0\otimes \psi^-_{b_2}}{\vv-\vv^{-1}}-\frac{\psi^+_0\otimes \psi^-_{b_2}}{\vv-\vv^{-1}}=0.$
Here we used $[e_{-1}\otimes\psi_{b_2}^-,\psi_0^+\otimes f_1]=0$ as
  $\psi^-_{b_2}f_1=\vv^2f_1\psi^-_{b_2}$ in $\CU^\ssc_{0,b_2}$,
due to~(\ref{u5}), and
  $\psi^+_0e_{-1}=\vv^2e_{-1}\psi^+_0$ in $\CU^\ssc_{0,b_1}$,
due to~(\ref{u4}).

\medskip
\noindent
$\circ$ \emph{Case $r=b_2-1, s=b_1$.}
The proof of
  $[\Delta_{b_1,b_2}(e_{b_2-1}),\Delta_{b_1,b_2}(f_{b_1})]=\Delta_{b_1,b_2}(\psi^-_{b})\Delta_{b_1,b_2}(h_{-1})$
can be deduced by applying Remark~\ref{Maps jmath} analogously to the
case $r=0,s=1$. Indeed,
  $[\Delta_{b_1,b_2}(e_{b_2-1}),\Delta_{b_1,b_2}(f_{b_1})]=
   [\jmath^-_{b_1,0}\otimes \jmath^-_{0,b_2}(\Delta(e_{-1})), \jmath^-_{b_1,0}\otimes \jmath^-_{0,b_2}(\Delta(f_0))]=
   \jmath^-_{b_1,0}\otimes \jmath^-_{0,b_2}(\Delta([e_{-1},f_0]))\\
   =\jmath^-_{b_1,0}\otimes \jmath^-_{0,b_2}(\Delta(\psi^-_0)\Delta(h_{-1}))=
   \Delta_{b_1,b_2}(\psi^-_b)\Delta_{b_1,b_2}(h_{-1}).$

\medskip
\noindent
$\circ$ \emph{Case $r=b_2, s=b_1$.}
We need to prove
  $[\Delta(e_{b_2}),\Delta(f_{b_1})]=-\frac{1}{\vv-\vv^{-1}}\Delta(\psi^-_{b})$.
This follows from
  $[\Delta(e_{b_2}),\Delta(f_{b_1})]=[1\otimes e_{b_2}, f_{b_1}\otimes 1+\psi^-_{b_1}\otimes f_0]=
   \psi^-_{b_1}\otimes [e_{b_2}, f_0]=-\frac{\psi^-_{b_1}\otimes \psi^-_{b_2}}{\vv-\vv^{-1}}=
   -\frac{\Delta(\psi^-_{b})}{\vv-\vv^{-1}}$,
due to~(\ref{u6}) for $\CU^\ssc_{0,b_2}$.

\medskip
\noindent
$\circ$ \emph{Case $b_2<r<0, s=b_1$.}
We need to prove
  $[\Delta(e_r),\Delta(f_{b_1})]=0$.
This follows from
  $[\Delta(e_r),\Delta(f_{b_1})]=[1\otimes e_r,f_{b_1}\otimes 1+\psi^-_{b_1}\otimes f_0]=
   \psi^-_{b_1}\otimes [e_r,f_0]=0$ as $[e_r,f_0]=0$ in $\CU^\ssc_{0,b_2}$ by~(\ref{u6}).

\medskip
\noindent
$\circ$ \emph{Case $r=b_2-1, 1+b_1<s\leq 0$.}
We need to prove
  $[\Delta(e_{b_2-1}),\Delta(f_s)]=0$.
This follows from
  $[\Delta(e_{b_2-1}),\Delta(f_s)]=
   [e_{-1}\otimes \psi^-_{b_2}+1\otimes e_{b_2-1}, f_s\otimes 1]=
   [e_{-1},f_s]\otimes \psi^-_{b_2}=0$
as $[e_{-1},f_s]=0$ in $\CU^\ssc_{0,b_1}$.

\medskip
\noindent
$\circ$ \emph{Case $r=b_2-1, s=1+b_1$.}
We need to prove
  $[\Delta(e_{b_2-1}),\Delta(f_{1+b_1})]=-\frac{1}{\vv-\vv^{-1}}\Delta(\psi^-_{b})$.
This follows from
  $[\Delta(e_{b_2-1}),\Delta(f_{1+b_1})]=
   [e_{-1}\otimes \psi^-_{b_2}+1\otimes e_{b_2-1}, f_{1+b_1}\otimes 1]=
   [e_{-1},f_{1+b_1}]\otimes \psi^-_{b_2}=-\frac{\psi^-_{b_1}\otimes \psi^-_{b_2}}{\vv-\vv^{-1}}=
   -\frac{\Delta(\psi^-_{b})}{\vv-\vv^{-1}}$,
due to~(\ref{u6}) for $\CU^\ssc_{0,b_1}$.


\subsection{Compatibility with (\ref{u9})}
\

Applying Remark~\ref{Maps jmath} as we did above, we see that the
equalities
\begin{equation*}
    [\Delta(h_1),[\Delta(f_1),[\Delta(h_1),\Delta(e_0)]]]=0\ \mathrm{and}\
    [\Delta(h_{-1}),[\Delta(e_{b_2-1}),[\Delta(h_{-1}),\Delta(f_{b_1})]]]=0
\end{equation*}
follow from the equalities
  $[h_1,[f_1,[h_1,e_0]]]=[2]_\vv\cdot[h_1,[f_1,e_1]]=[2]_\vv\cdot [h_1,\frac{-\psi^+_2}{\vv-\vv^{-1}}]=0$
in $U^+_\vv$ and
  $[h_{-1},[e_{-1},[h_{-1},f_0]]]=-[2]_\vv\cdot [h_{-1},[e_{-1},f_{-1}]]=
   [2]_\vv\cdot [h_{-1},\frac{\psi^-_{-2}}{\vv-\vv^{-1}}]=0$
in $U^-_\vv$, respectively.

\medskip
This completes our proof of Theorem~\ref{Coproduct for sl2}.


\section{Proof of Lemma~\ref{Properties of iota}(b)}
\label{Proof of PBW}


\subsection{PBW property for $\CU^\ssc_{0,n}$}
\

For $\CU^\ssc_{0,n}$, the simply-connected shifted quantum affine
algebra of $\ssl_2$, define the PBW variables to be
  $\{e_s\}_{s\in \BZ}\cup \{f_s\}_{s\in \BZ}\cup
   \{\psi^+_r\}_{r>0}\cup \{\psi^-_{n-r}\}_{r>0}\cup
   \{(\psi^+_0)^{\pm 1}\}\cup \{(\psi^-_{n})^{\pm 1}\}$.
We order the elements in each group according to the decreasing
order of $s,r$. Any expression of the form
\begin{equation*}
    e_{s^+_1}\cdots e_{s^+_a}f_{s^-_1}\cdots f_{s^-_b}
    \psi^+_{r^+_1}\cdots \psi^+_{r^+_{c^+}}\psi^-_{r^-_1}\cdots \psi^-_{r^-_{c^-}}
    (\psi^+_0)^{\gamma^+}(\psi^-_{n})^{\gamma^-}
\end{equation*}
with
  $s^+_1\geq \cdots \geq s^+_a,\ s^-_1\geq \cdots\geq s^-_b,\
   r^+_1\geq \cdots\geq r^+_{c^+}>0,\ r^-_1\leq \cdots\leq r^-_{c^-}<n,\
   \gamma^\pm\in \BZ,\ a,b,c^\pm\in \BN$,
will be referred to as the ordered monomial in the PBW variables.

The following result is easy to check using defining
relations~(\ref{U1}--\ref{U6}).

\begin{Lem}\label{PBW 1}
The algebra $\CU^\ssc_{0,n}$ is spanned by the ordered monomials in the
PBW variables.
\end{Lem}

The key result of this section is a refinement of the previous statement.

\begin{Thm}\label{PBW 2}
For any $n\in \BZ$, the algebra $\CU^\ssc_{0,n}$ satisfies the PBW
property, that is, the set of the ordered monomials in the PBW variables
forms a $\BC(\vv)$-basis of $\CU^\ssc_{0,n}$.
\end{Thm}


\subsection{Proof of Theorem~\ref{PBW 2}}
\

 We will prove this result in four steps.


\medskip
\noindent
\emph{Step 1:} Reduction to $\wt{\CU}^\ssc_{0,n}$.

Consider the associative $\BC(\vv)$-algebra $\wt{\CU}^\ssc_{0,n}$, defined in the same way as
$\CU^\ssc_{0,n}$ but without the generators $(\psi^+_0)^{-1},
(\psi^-_{n})^{-1}$. Note that $\CU^\ssc_{0,n}$ is the localization
of $\wt{\CU}^\ssc_{0,n}$ by the multiplicative set generated by
$\psi^+_0,\psi^-_{n}$. Since these generators are among the PBW
variables, the PBW property for $\CU^\ssc_{0,n}$ follows from the
PBW property for $\wt{\CU}^\ssc_{0,n}$.


\medskip
\noindent
\emph{Step 2:} PBW property for $\wt{\CU}^\ssc_{0,0}$.

It is well-known that the algebra $U_\vv(L\ssl_2)$ satisfies the PBW
property with the PBW variables chosen as
  $\{e_s\}_{s\in \BZ}\cup \{f_s\}_{s\in \BZ}\cup
   \{\psi^+_r\}_{r>0}\cup \{\psi^-_{-r}\}_{r>0}\cup \{(\psi^+_0)^{\pm 1}\}$.
Here the elements in each group are ordered according to the
decreasing order of $r,s$.

\begin{Lem}\label{PBW 3}
There is an embedding of algebras
  $\wt{\CU}^\ssc_{0,0}\hookrightarrow U_\vv(L\ssl_2)\otimes_{\BC(\vv)} \BC(\vv)[t]$,
such that
\begin{equation*}
    e_s\mapsto e_s\otimes t,\ f_s\mapsto f_s\otimes 1,\
    \psi^\pm_{\pm r}\mapsto \psi^\pm_{\pm r}\otimes t.
\end{equation*}
\end{Lem}

\begin{proof}
The above assignment obviously preserves all the defining relations
of $\wt{\CU}^\ssc_{0,0}$. Hence, it gives rise to a homomorphism
$\wt{\CU}^\ssc_{0,0}\to U_\vv(L\ssl_2)\otimes_{\BC(\vv)} \BC(\vv)[t]$.

To prove the injectivity of this homomorphism, let us first note that
$\wt{\CU}^\ssc_{0,0}$ is spanned by the ordered monomials in the PBW
variables, cf. Lemma~\ref{PBW 1}. The above homomorphism maps these monomials to a subset
of the basis for $U_\vv(L\ssl_2)\otimes_{\BC(\vv)} \BC(\vv)[t]$,
where we used the PBW property for $U_\vv(L\ssl_2)$. Hence, the
ordered monomials in the PBW variable for $\wt{\CU}^\ssc_{0,0}$ are
linearly independent and the above homomorphism is injective.
\end{proof}

Our proof of Lemma~\ref{PBW 3} implies the PBW property for
$\wt{\CU}^\ssc_{0,0}$.


\medskip
\noindent
\emph{Step 3:} PBW property for $\wt{\CU}^\ssc_{0,n},\ n<0$.

For $n<0$, the algebra $\wt{\CU}^\ssc_{0,n}$ is obviously a quotient
of $\wt{\CU}^\ssc_{0,0}$ by the 2-sided ideal
\begin{equation*}
  I_n:=\langle \psi^-_0, \psi^-_{-1}, \ldots, \psi^-_{1+n}\rangle_{\mathrm{2-sided}}.
\end{equation*}
Let $I^l_n$ be the left ideal generated by the same elements
\begin{equation*}
  I^l_n:=\langle \psi^-_0, \psi^-_{-1}, \ldots, \psi^-_{1+n}\rangle_{\mathrm{left}}.
\end{equation*}

\begin{Lem}\label{PBW 4}
We have $I^l_n=I_n$.
\end{Lem}

\begin{proof}
It suffices to show that $I^l_n$ is also a right ideal.
According to~(\ref{U4}), we have
\begin{equation*}
    \psi^-_{-r}e_s=
    \vv^{-2}\psi^-_{-r+1}e_{s-1}-e_{s-1}\psi^-_{-r+1}+\vv^{-2}e_s\psi^-_{-r},\
    \psi^-_0e_s=\vv^{-2}e_s\psi^-_0,
\end{equation*}
so that the right multiplication by $e_s$ preserves $I^l_n$.
Similarly for $f_s$ (need to apply~(\ref{U5})), while for
$\psi^+_r,\psi^-_{-r}$ this is obvious. These elements generate
$\wt{\CU}^\ssc_{0,0}$, hence, the claim.
\end{proof}

Combining the PBW property for $\wt{\CU}^\ssc_{0,0}$
(established in Step 2) with Lemma~\ref{PBW 4} and
  $\wt{\CU}^\ssc_{0,n}\simeq \wt{\CU}^\ssc_{0,0}/I_n$,
we get the PBW property for $\wt{\CU}^\ssc_{0,n}$.


\medskip
\noindent
\emph{Step 4:} PBW property for $\wt{\CU}^\ssc_{0,n},\ n>0$.

The proof proceeds by induction in $n$.
We assume that the PBW property holds for
$\wt{\CU}^\ssc_{0,m}$ with $m<n$ and want
to deduce the PBW property for $\wt{\CU}^\ssc_{0,n}$.
Consider the homomorphism
  $\wt{\iota}_{n,-1,0}\colon \wt{\CU}^\ssc_{0,n}\to \wt{\CU}^\ssc_{0,n-1}$
defined analogously to $\iota_{n,-1,0}$ of
Proposition~\ref{Maps iota}. Explicitly,
\begin{equation*}
    \wt{\iota}_{n,-1,0}\colon e_s\mapsto e_s-e_{s-1},\
    f_s\mapsto f_s,\ \psi^+_r\mapsto \psi^+_r-\psi^+_{r-1},\
    \psi^-_r\mapsto \psi^-_{r}-\psi^-_{r-1},
\end{equation*}
where we set $\psi^+_{-1}:=0,\ \psi^-_n:=0$ in the right-hand sides.
The image of an ordered monomial in the PBW variables for
$\wt{\CU}^\ssc_{0,n}$ under $\wt{\iota}_{n,-1,0}$ is a linear combination
of the same ordered
monomial in the PBW variables for $\wt{\CU}^\ssc_{0,n-1}$ with
all $\psi^-_r$ replaced by $(-\psi^-_{r-1})$, called the leading
monomial, and several other (not necessarily ordered) monomials in
the PBW variables. Based on the equality $e_se_{s-1}=\vv^2e_{s-1}e_s\ (s\in \BZ)$,
we see that rewriting these extra
monomials as linear combinations of the ordered monomials in the PBW
variables, all of them are actually lexicographically smaller than the
leading monomial. Hence, the PBW property for
$\wt{\CU}^\ssc_{0,n-1}$  implies the PBW property for $\wt{\CU}^\ssc_{0,n}$.
Moreover, we immediately get the injectivity of $\wt{\iota}_{n,-1,0}$.

\medskip
This completes our proof of Theorem~\ref{PBW 2}.


\subsection{Proof of Lemma~\ref{Properties of iota}(b)}
\

Now we are ready to prove Lemma~\ref{Properties of iota}(b).
Due to Lemma~\ref{Properties of iota}(a), it suffices to verify the
injectivity of the homomorphisms $\iota_{n,-1,0}$ and $\iota_{n,0,-1}$.
The former follows from the injectivity of $\wt{\iota}_{n,-1,0}$ from Step 4 above,
while the latter can be deduced in the same way.


\section{Proof of Theorem~\ref{Coproduct for sl2 in general}}
\label{Proof Coproduct General}

The proof of Theorem~\ref{Coproduct for sl2 in general} proceeds in three steps.
First, we construct $\Delta_{b_1,b_2}$
(this construction depends on a choice of sufficiently small $m_1,m_2\leq 0$).
Then, we verify that this construction is independent of the choice made.
Finally, we prove the commutativity of the diagram of Theorem~\ref{Coproduct for sl2 in general}
for any $m_1,m_2\in \BZ_{\leq 0}$.


\subsection{Construction of $\Delta_{b_1,b_2}$}
\

Fix any $m_1,m_2\in \BZ_{\leq 0}$ such that $b_1+m_1, b_2+m_2\in \BZ_{\leq 0}$.
Consider the diagram

\setlength{\unitlength}{1cm}
\begin{picture}(4,3.2)
 \put(3.2,2.5){$\CU^\ssc_{0,b}$}
 \put(2.6,0.5){$\CU^\ssc_{0,b+m_1+m_2}$}
 \put(9.0,2.5){$\CU^\ssc_{0,b_1}\otimes \CU^\ssc_{0,b_2}$}
 \put(8.4,0.5){$\CU^\ssc_{0,b_1+m_1}\otimes \CU^\ssc_{0,b_2+m_2}$}

 \put(4.6,0.6){\vector (1,0){3.7}}
 \put(3.5,2.3){\vector (0,-1){1.4}}
 \put(10.05,2.3){\vector (0,-1){1.4}}

 \put(5.0,0.25){$\Delta=\Delta_{b_1+m_1,b_2+m_2}$}
 \put(2.2,1.5){$\iota_{b,m_2,m_1}$}
 \put(10.2,1.5){$\iota_{b_1,0,m_1}\otimes \iota_{b_2,m_2,0}$}
\end{picture}

\noindent
where the bottom horizontal arrow $\Delta=\Delta_{b_1+m_1,b_2+m_2}$
is defined in Theorem~\ref{Coproduct for sl2}.
Since the homomorphisms $\iota_{b,m_2,m_1}$ and
$\iota_{b_1,0,m_1}\otimes \iota_{b_2,m_2,0}$ are injective,
the homomorphism $\Delta_{b_1+m_1,b_2+m_2}$ gives rise to a uniquely
determined homomorphism $\Delta_{b_1,b_2}$ making the above diagram
commutative as far as we can prove
\begin{equation}\tag{$\diamondsuit$}\label{diamonsuit}
  \Delta(\iota_{b,m_2,m_1}(\CU^\ssc_{0,b}))\subset
  (\iota_{b_1,0,m_1}\otimes \iota_{b_2,m_2,0})(\CU^\ssc_{0,b_1}\otimes \CU^\ssc_{0,b_2}).
\end{equation}

As before, we use
  $\CU^{\ssc,>}_{0,b'}, \CU^{\ssc,\geq}_{0,b'}, \CU^{\ssc,<}_{0,b'}, \CU^{\ssc,\leq}_{0,b'}$
to denote the $\BC(\vv)$-subalgebras of $\CU^\ssc_{0,b'}$ generated by
  $\{e_r\}, \{e_r,\psi^\pm_{\pm s^\pm}\},  \{f_r\}, \{f_r,\psi^\pm_{\pm s^\pm}\}$,
respectively. For $r\in\BZ$, we claim that
\begin{equation}\tag{$\diamond_1$}\label{diamond1}
  \Delta(e_r)\in 1\otimes e_r + \CU^{\ssc,>}_{0,b_1+m_1}\otimes \CU^{\ssc,\leq}_{0,b_2+m_2},\
  \Delta(f_r)\in f_r\otimes 1 + \CU^{\ssc,\geq}_{0,b_1+m_1}\otimes \CU^{\ssc,<}_{0,b_2+m_2}.
\end{equation}
This follows by combining iteratively the formulas for
$\Delta(e_{-1}),\Delta(f_0),\Delta(h_{\pm 1})$ with the relations
  $[h_{\pm 1}, e_r]=[2]_\vv\cdot e_{r\pm 1},
   [h_{\pm 1},f_r]=-[2]_\vv\cdot f_{r\pm 1}$.
We also note that
\begin{equation}\tag{$\diamond_2$}\label{diamond2}
  \CU^{\ssc,\geq}_{0,b_1}\otimes \CU^{\ssc,\leq}_{0,b_2}\subset
  (\iota_{b_1,0,m_1}\otimes \iota_{b_2,m_2,0})(\CU^\ssc_{0,b_1}\otimes \CU^\ssc_{0,b_2}).
\end{equation}

According to~(\ref{diamond1}), we get
\begin{equation*}
    \Delta(\iota_{b,m_2,m_1}(e_r))\in
    1\otimes \sum_{s=0}^{-m_2}(-1)^s\binom{-m_2}{s}e_{r-s} +
    \CU^{\ssc,>}_{0,b_1+m_1}\otimes \CU^{\ssc,\leq}_{0,b_2+m_2}.
\end{equation*}
The right-hand side is an element of
 $(\iota_{b_1,0,m_1}\otimes \iota_{b_2,m_2,0})(\CU^\ssc_{0,b_1}\otimes \CU^\ssc_{0,b_2})$,
due to~(\ref{diamond2}) and the equality
  $1\otimes \sum_{s=0}^{-m_2}(-1)^s\binom{-m_2}{s}e_{r-s}=
   (\iota_{b_1,0,m_1}\otimes \iota_{b_2,m_2,0})(1\otimes e_r)$.
Likewise,
\begin{equation*}
    \Delta(\iota_{b,m_2,m_1}(f_r))\in
    \sum_{s=0}^{-m_1}(-1)^s\binom{-m_1}{s}f_{r-s}\otimes 1 +
    \CU^{\ssc,\geq}_{0,b_1+m_1}\otimes \CU^{\ssc,<}_{0,b_2+m_2}.
\end{equation*}
The right-hand side is an element of
  $(\iota_{b_1,0,m_1}\otimes \iota_{b_2,m_2,0})(\CU^\ssc_{0,b_1}\otimes \CU^\ssc_{0,b_2})$,
due to~(\ref{diamond2}) and the equality
  $\sum_{s=0}^{-m_1}(-1)^s\binom{-m_1}{s}f_{r-s}\otimes 1=
   (\iota_{b_1,0,m_1}\otimes \iota_{b_2,m_2,0})(f_r\otimes 1)$.
We also have
\begin{equation*}
\begin{split}
  & \Delta(\iota_{b,m_2,m_1}((\psi^+_0)^{\pm 1}))=
    (\iota_{b_1,0,m_1}\otimes \iota_{b_2,m_2,0})((\psi_0)^{\pm 1}\otimes (\psi_0)^{\pm 1}),\\
  & \Delta(\iota_{b,m_2,m_1}((\psi^-_{b})^{\pm 1}))=
    (\iota_{b_1,0,m_1}\otimes \iota_{b_2,m_2,0})((\psi^-_{b_1})^{\pm 1}\otimes (\psi^-_{b_2})^{\pm 1}).
\end{split}
\end{equation*}
Finally, combining the relations
  $\psi_r^+=(\vv-\vv^{-1})[e_r,f_0], \psi_{b-r}^-=(\vv^{-1}-\vv)[e_{b-r},f_0]\ (r\in\BZ_{>0})$
in $\CU^\ssc_{0,b+m_1+m_2}$ with~(\ref{diamond1}) and~(\ref{diamond2}), we get
\begin{equation*}
    \Delta(\psi_r^+),\Delta(\psi_{b-r}^-)\in\CU^{\ssc,\geq}_{0,b_1+m_1}\otimes
    \CU^{\ssc,\leq}_{0,b_2+m_2}\subset (\iota_{b_1,0,m_1}\otimes
    \iota_{b_2,m_2,0})(\CU^\ssc_{0,b_1}\otimes \CU^\ssc_{0,b_2}).
\end{equation*}
This completes our proof of~(\ref{diamonsuit}).

Therefore, we obtain the homomorphism $\Delta_{b_1,b_2}$ for the
particular choice of $m_1,m_2$.


\subsection{Independence of the choice of $m_1,m_2$}
\label{F2}
\

Let us now prove that the homomorphism $\Delta_{b_1,b_2}$ constructed above
does not depend on the choice of $m_1,m_2$. To this end, fix another pair
$m'_1,m'_2\in \BZ_{\leq 0}$ such that $b_1+m'_1, b_2+m'_2\in \BZ_{\leq 0}$,
and set $m=m_1+m_2, m'=m'_1+m'_2$.

Consider the following diagram:

\setlength{\unitlength}{1cm}
\begin{picture}(4,5.2)
 \put(3.2,4.5){$\CU^\ssc_{0,b}$}
 \put(3.2,2.5){$\CU^\ssc_{0,b+m}$}
 \put(2.6,0.5){$\CU^\ssc_{0,b+m+m'}$}
 \put(9.65,4.5){$\CU^\ssc_{0,b_1}\otimes \CU^\ssc_{0,b_2}$}
 \put(9.0,2.5){$\CU^\ssc_{0,b_1+m_1}\otimes \CU^\ssc_{0,b_2+m_2}$}
 \put(8.4,0.5){$\CU^\ssc_{0,b_1+m_1+m'_1}\otimes \CU^\ssc_{0,b_2+m_2+m'_2}$}

 \put(4.6,2.6){\vector (1,0){4.2}}
 \put(4.6,0.6){\vector (1,0){3.7}}
 \put(3.5,2.3){\vector (0,-1){1.4}}
 \put(10.7,2.3){\vector (0,-1){1.4}}
 \put(3.5,4.3){\vector (0,-1){1.4}}
 \put(10.7,4.3){\vector (0,-1){1.4}}

 \put(4.8,0.25){$\Delta_{b_1+m_1+m'_1,b_2+m_2+m'_2}$}
 \put(5.4,2.25){$\Delta_{b_1+m_1,b_2+m_2}$}
 \put(1.7,1.5){$\iota_{b+m,m'_2,m'_1}$}
 \put(2.2,3.5){$\iota_{b,m_2,m_1}$}
 \put(10.8,3.5){$\iota_{b_1,0,m_1}\otimes \iota_{b_2,m_2,0}$}
 \put(10.8,1.5){$\iota_{b_1+m_1,0,m'_1}\otimes \iota_{b_2+m_2,m'_2,0}$}
\end{picture}

According to Lemma~\ref{Properties of iota}(a):
  $\iota_{b+m,m'_2,m'_1}\circ \iota_{b,m_2,m_1}=\iota_{b,m_2+m'_2,m_1+m'_1}$
and
  $(\iota_{b_1+m_1,0,m'_1}\otimes \iota_{b_2+m_2,m'_2,0})\circ
   (\iota_{b_1,0,m_1}\otimes \iota_{b_2,m_2,0})=
   (\iota_{b_1,0,m_1+m'_1}\otimes \iota_{b_2,m_2+m'_2,0})$.
On the other hand, tracing back the explicit formulas for $\Delta_{b_1+m_1,b_2+m_2}$
and $\Delta_{b_1+m_1+m'_1,b_2+m_2+m'_2}$ of Theorem~\ref{Coproduct
for sl2}, it is easy to check that the lower square is commutative.

The above two observations imply that the maps $\Delta_{b_1,b_2}$
are the same for both $(m_1,m_2)$ and $(m_1+m'_1,m_2+m'_2)$. Due to
the symmetry, we also see that the maps $\Delta_{b_1,b_2}$ are the
same for both $(m'_1,m'_2)$ and $(m_1+m'_1,m_2+m'_2)$. Therefore,
the maps $\Delta_{b_1,b_2}$ are the same for both $(m_1,m_2)$ and
$(m'_1,m'_2)$. This completes our verification.


\subsection{Commutativity of the diagram for any $m_1,m_2\in \BZ_{\leq 0}$}
\

It remains to prove the commutativity of the diagram of
Theorem~\ref{Coproduct for sl2 in general}. To this end, choose
$m'_1,m'_2\in \BZ_{\leq 0}$ such that $b_1+m_1+m'_1,b_2+m_2+m'_2\in \BZ_{\leq 0}$.
Consider a diagram analogous to the previous one:

\setlength{\unitlength}{1cm}
\begin{picture}(4,5.2)
 \put(3.2,4.5){$\CU^\ssc_{0,b}$}
 \put(3.2,2.5){$\CU^\ssc_{0,b+m}$}
 \put(2.6,0.5){$\CU^\ssc_{0,b+m+m'}$}
 \put(9.65,4.5){$\CU^\ssc_{0,b_1}\otimes \CU^\ssc_{0,b_2}$}
 \put(9.0,2.5){$\CU^\ssc_{0,b_1+m_1}\otimes \CU^\ssc_{0,b_2+m_2}$}
 \put(8.4,0.5){$\CU^\ssc_{0,b_1+m_1+m'_1}\otimes \CU^\ssc_{0,b_2+m_2+m'_2}$}

 \put(4.3,4.6){\vector (1,0){5.0}}
 \put(4.6,2.6){\vector (1,0){4.2}}
 \put(4.6,0.6){\vector (1,0){3.7}}
 \put(3.5,2.3){\vector (0,-1){1.4}}
 \put(10.7,2.3){\vector (0,-1){1.4}}
 \put(3.5,4.3){\vector (0,-1){1.4}}
 \put(10.7,4.3){\vector (0,-1){1.4}}

 \put(4.8,0.25){$\Delta_{b_1+m_1+m'_1,b_2+m_2+m'_2}$}
 \put(5.4,2.25){$\Delta_{b_1+m_1,b_2+m_2}$}
 \put(6.3,4.25){$\Delta_{b_1,b_2}$}
 \put(1.7,1.5){$\iota_{b+m,m'_2,m'_1}$}
 \put(2.2,3.5){$\iota_{b,m_2,m_1}$}
 \put(10.8,3.5){$\iota_{b_1,0,m_1}\otimes \iota_{b_2,m_2,0}$}
 \put(10.8,1.5){$\iota_{b_1+m_1,0,m'_1}\otimes \iota_{b_2+m_2,m'_2,0}$}
\end{picture}

By our construction, the lower square is commutative.
Applying~Lemma~\ref{Properties of iota}(a) as in~Section~\ref{F2},
we also see that the outer square is commutative.
Hence, the commutativity of the top square
follows from the injectivity of the homomorphism
$\iota_{b_1+m_1,0,m'_1}\otimes \iota_{b_2+m_2,m'_2,0}$, due
to~Lemma~\ref{Properties of iota}(b).


\section{Proof of Theorem~\ref{Coproduct on quantum loop sln}}
\label{Proof of coproduct on loop sln}

The proof of~Theorem~\ref{Coproduct on quantum loop sln} proceeds in several steps.
First, we recall the RTT presentation of $U_\vv(L\ssl_n)$, and derive
the equalities of the right-hand sides of~(\ref{long commutators}).
Then, we compute the RTT coproduct of certain elements
$\tilde{g}^{(\pm 1)}_i$ from the RTT presentation,
see~Theorems~\ref{Coproduct of k_i1},~\ref{Coproduct of k_i-1} (this is the
most technical part). This allows us to derive formulas~(\ref{h_i1})
and~(\ref{h_i-1}). Based on these, we deduce~(\ref{e_i-1}) and~(\ref{f_i1}).


\subsection{RTT presentation of $U_\vv(L\ssl_n)$}
\

Let $R_\trig(z/w)\in \mathrm{End}(\BC^n\otimes \BC^n)$ be the standard trigonometric
$R$-matrix of $\ssl_n$-type:
\begin{equation}\label{R-matrix sl_n}
\begin{split}
  & R_\trig(z/w):=
  \sum_{i=1}^n E_{ii}\otimes E_{ii} +\sum_{1\leq i\ne j\leq n} \frac{z-w}{\vv z-\vv^{-1}w}E_{ii}\otimes E_{jj}+\\
  & \sum_{1\leq j<i\leq n} \left(\frac{(\vv-\vv^{-1})z}{\vv z-\vv^{-1}w}E_{ji}\otimes E_{ij}+
   \frac{(\vv-\vv^{-1})w}{\vv z-\vv^{-1}w} E_{ij}\otimes E_{ji}\right)
\end{split}
\end{equation}
(for $n=2$, this definition coincides with formula~(\ref{R-matrix})).

Define the RTT algebra of $\ssl_n$, denoted by $U^\rtt(\ssl_n)$,
to be the associative $\BC(\vv)$-algebra generated by
  $\{t^\pm_{ij}[\pm r]\}_{1\leq i,j\leq n}^{r\in \BN}$
subject to the following defining relations:
\begin{equation}\label{nRTT1}
  t^\pm_{ii}[0]t^\mp_{ii}[0]=1\ \ \mathrm{for}\ 1\leq i\leq n,\
  t^+_{ij}[0]=t^-_{ji}[0]=0 \ \ \mathrm{for}\ j<i,
\end{equation}
\begin{equation}\label{nRTT2}
   R_\trig(z/w)(T^\epsilon(z)\otimes 1)(1\otimes T^{\epsilon'}(w))=
   (1\otimes T^{\epsilon'}(w))(T^\epsilon(z)\otimes 1)R_\trig(z/w),
\end{equation}
\begin{equation}\label{nRTT3}
  \qdet\ T^\pm(z)=1,
\end{equation}
for all $\epsilon,\epsilon'\in \{\pm\}$, where the matrices
$T^\pm(z)\in \mathrm{Mat}_{n\times n}(U^\rtt(\ssl_n))$ are given by
\begin{equation*}
    T^\pm(z):=\sum_{i,j=1}^n T^\pm_{ij}(z)\cdot E_{ij}\ \
    \mathrm{with}\ \ T^\pm_{ij}(z):=\sum_{r\geq 0} t^\pm_{ij}[\pm r]  z^{\mp r},
\end{equation*}
and the quantum determinant $\qdet$ is defined in a standard way as
\begin{equation*}
    \qdet\ T^\pm(z):=\sum_{\tau\in \mathfrak{S}_n}(-\vv)^{-l(\tau)}
    T^\pm_{1,\tau(1)}(z)T^\pm_{2,\tau(2)}(\vv^{-2}z)\cdots T^\pm_{n,\tau(n)}(\vv^{2-2n}z)
\end{equation*}
(cf.~Section~\ref{section RTT algebras} and a footnote there).

\begin{Rem}\label{4vs3 RTT}
Let us point out right away that the RTT presentation of $U_q(\widehat{\gl}_n)$
(with a nontrivial central charge), given in~\cite[Definition 3.2]{df}, involves
only three out of four relations~(\ref{nRTT2}), namely for
$(\epsilon,\epsilon')=(+,+), (-,-), (-,+)$. However, as pointed out
in~\cite[2.3]{gm}, if the central charge is trivial, then the fourth relation for
$(\epsilon,\epsilon')=(+,-)$ is equivalent to the one for $(\epsilon,\epsilon')=(-,+)$.
Indeed, in our notations, this follows from the equalities
$R_\trig(z/w)^{-1}=R'_\trig(z/w), PR'_\trig(w/z)P^{-1}=R_\trig(z/w)$,
where $R'_\trig(z/w)$ is obtained from $R_\trig(z/w)$ by replacing $\vv$ with $\vv^{-1}$
and $P\in \mathrm{End}(\BC^n\otimes \BC^n)$ denotes the permutation operator.
\end{Rem}

Note that $T^\pm(z)$ admits the following unique Gauss decomposition:
\begin{equation*}
  T^\pm(z)=\wt{F}^\pm(z)\cdot \wt{G}^\pm(z)\cdot \wt{E}^\pm(z)
\end{equation*}
with
  $\wt{F}^\pm(z),\wt{G}^\pm(z),\wt{E}^\pm(z)\in \mathrm{Mat}_{n\times n}(U^\rtt(\ssl_n))$
of the form
\begin{equation*}
    \wt{F}^\pm(z)=\sum_{i} E_{ii}+\sum_{j<i} \tilde{f}^\pm_{ij}(z)\cdot E_{ij},\
    \wt{G}^\pm(z)=\sum_{i} \tilde{g}^\pm_i(z)\cdot E_{ii},\
    \wt{E}^\pm(z)=\sum_{i} E_{ii}+\sum_{j<i} \tilde{e}^\pm_{ji}(z)\cdot E_{ji}.
\end{equation*}
We endow $U^\rtt(\ssl_n)$ with the coproduct structure
(also known as the RTT coproduct) via
\begin{equation*}
    \Delta^\rtt\colon U^\rtt(\ssl_n)\longrightarrow U^\rtt(\ssl_n)\otimes U^\rtt(\ssl_n)\ \
    \mathrm{given\ by}\ \ \Delta^\rtt(T^\pm(z)):=T^\pm(z)\otimes T^\pm(z).
\end{equation*}

\begin{Thm}\label{RTT for sln}\cite{df}
There exists a unique $\BC(\vv)$-algebra isomorphism
\begin{equation*}
  \Upsilon\colon  U^\ad_\vv (L\ssl_n)\iso U^\rtt(\ssl_n),
\end{equation*}
such that
\begin{equation*}
\begin{split}
  & e^\pm_j(z)\mapsto \frac{\tilde{e}^\pm_{j,j+1}(\vv^jz)}{\vv-\vv^{-1}},\
  f^\pm_j(z)\mapsto \frac{\tilde{f}^\pm_{j+1,j}(\vv^jz)}{\vv-\vv^{-1}},\\
  & \psi^\pm_j(z)\mapsto \tilde{g}^\pm_{j+1}(\vv^jz)(\tilde{g}^\pm_j(\vv^jz))^{-1},\
  \phi^\pm_j\mapsto t^\mp_{11}[0]t^\mp_{22}[0]\cdots t^\mp_{jj}[0]
  \ \mathrm{for}\ 1\leq j<n.
\end{split}
\end{equation*}
Moreover, this isomorphism intertwines the Drinfeld-Jimbo coproduct $\Delta^\ad$
on $U^\ad_\vv(L\ssl_n)$ with the RTT coproduct $\Delta^\rtt$ on $U^\rtt(\ssl_n)$.
\end{Thm}

\begin{Rem}\label{embedding}
Restricting $\Upsilon$ to $U_\vv(L\ssl_n)$, viewed as a Hopf subalgebra of
$U^\ad_\vv (L\ssl_n)$, we get an embedding
$U_\vv(L\ssl_n)\hookrightarrow U^\rtt(\ssl_n)$. We will deliberately refer
to $U^\rtt(\ssl_n)$ as an RTT presentation of both algebras
$U_\vv(L\ssl_n)$ and $U^\ad_\vv(L\ssl_n)$.
\end{Rem}

Let us express the matrix coefficients of $\wt{F}^\pm(z),\wt{G}^\pm(z),\wt{E}^\pm(z)$
as Taylor series in $z^{\mp 1}$:
  $\tilde{e}^+_{ji}(z)=\sum_{r\geq 0} \tilde{e}^{(r)}_{ji}z^{-r},\
   \tilde{e}^-_{ji}(z)=\sum_{r<0} \tilde{e}^{(r)}_{ji}z^{-r},\
   \tilde{f}^+_{ij}(z)=\sum_{r>0} \tilde{f}^{(r)}_{ij}z^{-r},\
   \tilde{f}^-_{ij}(z)=\sum_{r\leq 0} \tilde{f}^{(r)}_{ij}z^{-r},\ \\
   \tilde{g}^\pm_i(z)=\tilde{g}^\pm_i+\sum_{r>0} \tilde{g}^{(\pm r)}_i z^{\mp r}$.
According to Theorem~\ref{RTT for sln}, we have
\begin{equation}\label{identification baby}
\begin{split}
  & \Upsilon^{-1}(\tilde{e}^{(0)}_{j,j+1})=(\vv-\vv^{-1})e_{j,0},\
  \Upsilon^{-1}(\tilde{f}^{(0)}_{j+1,j})=-(\vv-\vv^{-1})f_{j,0},\\
  & \Upsilon^{-1}(\tilde{e}^{(-1)}_{j,j+1})=-\vv^{-j}(\vv-\vv^{-1})e_{j,-1},
  \Upsilon^{-1}(\tilde{f}^{(1)}_{j+1,j})=\vv^j(\vv-\vv^{-1})f_{j,1}.
\end{split}
\end{equation}

The following is the key technical result of this subsection.

\begin{Prop}\label{identification serious}
For any $1\leq j<k<i\leq n$, we have:

\noindent
(a) $\tilde{e}^{(0)}_{ji}=\frac{1}{\vv-\vv^{-1}} [\tilde{e}^{(0)}_{ki},\tilde{e}^{(0)}_{jk}]_{\vv^{-1}}$.

\noindent
(b) $\tilde{f}^{(0)}_{ij}=\frac{-1}{\vv-\vv^{-1}}[\tilde{f}^{(0)}_{kj},\tilde{f}^{(0)}_{ik}]_{\vv}$.

\noindent
(c) $\tilde{e}^{(-1)}_{ji}=\frac{1}{\vv-\vv^{-1}}[\tilde{e}^{(0)}_{ki},\tilde{e}^{(-1)}_{jk}]_{\vv^{-1}}$.

\noindent
(d) $\tilde{f}^{(1)}_{ij}=\frac{-1}{\vv-\vv^{-1}}[\tilde{f}^{(1)}_{kj},\tilde{f}^{(0)}_{ik}]_{\vv}$.
\end{Prop}

\begin{proof}
(a) Comparing the matrix coefficients $\langle v_j\otimes v_k|\cdots|v_k\otimes v_i\rangle$ of both sides of
the equality $R_\trig(z/w)(T^+(z)\otimes 1)(1\otimes T^+(w))=(1\otimes T^+(w))(T^+(z)\otimes 1)R_\trig(z/w)$,
we get
\begin{equation*}
    (z-w)T^+_{jk}(z)T^+_{ki}(w)+(\vv-\vv^{-1})zT^+_{kk}(z)T^+_{ji}(w)=
    (z-w)T^+_{ki}(w)T^+_{jk}(z)+(\vv-\vv^{-1})wT^+_{kk}(w)T^+_{ji}(z).
\end{equation*}
Evaluating the coefficients of $z^1w^0$ in both sides of this equality, we find
\begin{equation*}
    \tilde{g}^+_j\tilde{e}^{(0)}_{jk}\tilde{g}^+_k\tilde{e}^{(0)}_{ki}+(\vv-\vv^{-1})\tilde{g}^+_k\tilde{g}^+_j\tilde{e}^{(0)}_{ji}=
    \tilde{g}^+_k\tilde{e}^{(0)}_{ki}\tilde{g}^+_j\tilde{e}^{(0)}_{jk}.
\end{equation*}
Combining this with Lemma~\ref{identification semi-serious} below, we obtain
\begin{equation*}
    (\vv-\vv^{-1})\tilde{g}^+_k\tilde{g}^+_j\tilde{e}^{(0)}_{ji}=
    \tilde{g}^+_k\tilde{g}^+_j[\tilde{e}^{(0)}_{ki},\tilde{e}^{(0)}_{jk}]_{\vv^{-1}}\Longrightarrow
    \tilde{e}^{(0)}_{ji}=[\tilde{e}^{(0)}_{ki},\tilde{e}^{(0)}_{jk}]_{\vv^{-1}}/(\vv-\vv^{-1}).
\end{equation*}

(b) Comparing the matrix coefficients $\langle v_i\otimes v_k|\cdots|v_k\otimes v_j\rangle$ of both sides of the equality
$R_\trig(z/w)(T^-(z)\otimes 1)(1\otimes T^-(w))=(1\otimes T^-(w))(T^-(z)\otimes 1)R_\trig(z/w)$, we get
\begin{equation*}
    (z-w)T^-_{ik}(z)T^-_{kj}(w)+(\vv-\vv^{-1})wT^-_{kk}(z)T^-_{ij}(w)=
    (z-w)T^-_{kj}(w)T^-_{ik}(z)+(\vv-\vv^{-1})zT^-_{kk}(w)T^-_{ij}(z).
\end{equation*}
Evaluating the coefficients of $z^0w^1$ in both sides of this equality, we find
\begin{equation*}
    -\tilde{f}^{(0)}_{ik}\tilde{g}^-_k\tilde{f}^{(0)}_{kj}\tilde{g}^-_j+(\vv-\vv^{-1})\tilde{g}^-_k\tilde{f}^{(0)}_{ij}\tilde{g}^-_j=
    -\tilde{f}^{(0)}_{kj}\tilde{g}^-_j\tilde{f}^{(0)}_{ik}\tilde{g}^-_k.
\end{equation*}

Combining this with Lemma~\ref{identification semi-serious} below, we obtain
\begin{equation*}
    -(\vv-\vv^{-1})\tilde{f}^{(0)}_{ij}\tilde{g}^-_k\tilde{g}^-_j=
    [\tilde{f}^{(0)}_{kj},\tilde{f}^{(0)}_{ik}]_{\vv}\cdot \tilde{g}^-_k\tilde{g}^-_j\Longrightarrow
    \tilde{f}^{(0)}_{ij}=-[\tilde{f}^{(0)}_{kj},\tilde{f}^{(0)}_{ik}]_{\vv}/(\vv-\vv^{-1}).
\end{equation*}

(c) Comparing the matrix coefficients $\langle v_k\otimes v_j|\cdots|v_i\otimes v_k\rangle$ of both sides of the equality
$R_\trig(z/w)(T^+(z)\otimes 1)(1\otimes T^-(w))=(1\otimes T^-(w))(T^+(z)\otimes 1)R_\trig(z/w)$, we get
\begin{equation*}
    (z-w)T^+_{ki}z)T^-_{jk}(w)+(\vv-\vv^{-1})wT^+_{ji}(z)T^-_{kk}(w)=
    (z-w)T^-_{jk}(w)T^+_{ki}(z)+(\vv-\vv^{-1})zT^-_{ji}(w)T^+_{kk}(z).
\end{equation*}
Evaluating the coefficients of $z^1w^1$ in both sides of this equality, we find
\begin{equation}\label{catch 1}
\begin{split}
  & \tilde{g}^+_k\tilde{e}^{(0)}_{ki}\left(\tilde{g}^-_j\tilde{e}^{(-1)}_{jk}+\sum_{j'<j}\tilde{f}^{(0)}_{jj'}\tilde{g}^-_{j'}\tilde{e}^{(-1)}_{j'k}\right)=\\
  & \left(\tilde{g}^-_j\tilde{e}^{(-1)}_{jk}+\sum_{j'<j}\tilde{f}^{(0)}_{jj'}\tilde{g}^-_{j'}\tilde{e}^{(-1)}_{j'k}\right)\tilde{g}^+_k\tilde{e}^{(0)}_{ki}+
  (\vv-\vv^{-1})\left(\tilde{g}^-_j\tilde{e}^{(-1)}_{ji}+\sum_{j'<j}\tilde{f}^{(0)}_{jj'}\tilde{g}^-_{j'}\tilde{e}^{(-1)}_{j'i}\right)\tilde{g}^+_k.
\end{split}
\end{equation}

This equation actually implies
  $\tilde{g}^+_k\tilde{e}^{(0)}_{ki}\tilde{g}^-_j\tilde{e}^{(-1)}_{jk}=
   \tilde{g}^-_j\tilde{e}^{(-1)}_{jk} \tilde{g}^+_k\tilde{e}^{(0)}_{ki}+(\vv-\vv^{-1})\tilde{g}^-_j\tilde{e}^{(-1)}_{ji}\tilde{g}^+_k.$
We prove this by induction in $j$. For $j=1$, this is just~(\ref{catch 1}).
In general, note that for $j'<j<k<i$, the element $\tilde{f}^{(0)}_{jj'}$ commutes with
$\tilde{e}^{(0)}_{ki}$ and $\tilde{g}^+_k$. The latter follows from Lemma~\ref{identification semi-serious},
while the equality $[\tilde{f}^{(0)}_{jj'},\tilde{e}^{(0)}_{ki}]=0$ follows
by combining parts (a,b) from above with $[e_{a,0},f_{b,0}]=0$ for $a\ne b$.
Hence,~(\ref{catch 1}) implies $A(j,k,i)+\sum_{j'<j} \tilde{f}^{(0)}_{jj'} A(j',k,i)=0$, where we set
\begin{equation*}
    A(j,k,i):=\tilde{g}^+_k\tilde{e}^{(0)}_{ki}\tilde{g}^-_j\tilde{e}^{(-1)}_{jk}-
    \tilde{g}^-_j\tilde{e}^{(-1)}_{jk} \tilde{g}^+_k\tilde{e}^{(0)}_{ki}-(\vv-\vv^{-1})\tilde{g}^-_j\tilde{e}^{(-1)}_{ji}\tilde{g}^+_k.
\end{equation*}
By the induction assumption $A(j',k,i)=0$ for $j'<j$, hence, $A(j,k,i)=0$.

Combining this with Lemma~\ref{identification semi-serious} below, we obtain
\begin{equation*}
    (\vv-\vv^{-1})\tilde{g}^-_j\tilde{g}^+_k\tilde{e}^{(-1)}_{ji}=
    \tilde{g}^-_j\tilde{g}^+_k[\tilde{e}^{(0)}_{ki},\tilde{e}^{(-1)}_{jk}]_{\vv^{-1}}\Longrightarrow
    \tilde{e}^{(-1)}_{ji}=[\tilde{e}^{(0)}_{ki},\tilde{e}^{(-1)}_{jk}]_{\vv^{-1}}/(\vv-\vv^{-1}).
\end{equation*}

(d) Comparing the matrix coefficients $\langle v_k\otimes v_i|\cdots|v_j\otimes v_k\rangle$ of both sides of the equality
$R_\trig(z/w)(T^+(z)\otimes 1)(1\otimes T^-(w))=(1\otimes T^-(w))(T^+(z)\otimes 1)R_\trig(z/w)$, we get
\begin{equation*}
    (z-w)T^+_{kj}(z)T^-_{ik}(w)+(\vv-\vv^{-1})zT^+_{ij}(z)T^-_{kk}(w)=
    (z-w)T^-_{ik}(w)T^+_{kj}(z)+(\vv-\vv^{-1})wT^-_{ij}(w)T^+_{kk}(z).
\end{equation*}
Evaluating the coefficients of $z^0w^0$ in both sides of this equality, we find
\begin{equation}\label{catch 2}
\begin{split}
  & \left(\tilde{f}^{(1)}_{kj}\tilde{g}^+_j+\sum_{j'<j}\tilde{f}^{(1)}_{kj'}\tilde{g}^+_{j'}\tilde{e}^{(0)}_{j'j}\right)\tilde{f}^{(0)}_{ik}\tilde{g}^-_k+
  (\vv-\vv^{-1})\left(\tilde{f}^{(1)}_{ij}\tilde{g}^+_j+\sum_{j'<j}\tilde{f}^{(1)}_{ij'}\tilde{g}^+_{j'}\tilde{e}^{(0)}_{j'j}\right)\tilde{g}^-_k=\\
  & \tilde{f}^{(0)}_{ik}\tilde{g}^-_k\left(\tilde{f}^{(1)}_{kj}\tilde{g}^+_j+\sum_{j'<j}\tilde{f}^{(1)}_{kj'}\tilde{g}^+_{j'}\tilde{e}^{(0)}_{j'j}\right).
\end{split}
\end{equation}

This equation actually implies
  $\tilde{f}^{(1)}_{kj}\tilde{g}^+_j\tilde{f}^{(0)}_{ik}\tilde{g}^-_k+(\vv-\vv^{-1})\tilde{f}^{(1)}_{ij}\tilde{g}^+_j\tilde{g}^-_k=
    \tilde{f}^{(0)}_{ik}\tilde{g}^-_k\tilde{f}^{(1)}_{kj}\tilde{g}^+_j.$
We prove this by induction in $j$. For $j=1$, this is just~(\ref{catch 2}).
Analogously to part (c) above, we note that the element $\tilde{e}^{(0)}_{j'j}$ commutes with
$\tilde{f}^{(0)}_{ik}$ and $\tilde{g}^-_k$ for $j'<j<k<i$. Hence,~(\ref{catch 2}) implies
$B(j,k,i)+\sum_{j'<j} B(j',k,i)\tilde{e}^{(0)}_{j'j}=0$, where we set
\begin{equation*}
  B(j,k,i):=\tilde{f}^{(1)}_{kj}\tilde{g}^+_j\tilde{f}^{(0)}_{ik}\tilde{g}^-_k+(\vv-\vv^{-1})\tilde{f}^{(1)}_{ij}\tilde{g}^+_j\tilde{g}^-_k- \tilde{f}^{(0)}_{ik}\tilde{g}^-_k\tilde{f}^{(1)}_{kj}\tilde{g}^+_j.
\end{equation*}
By the induction assumption $B(j',k,i)=0$ for $j'<j$, hence, $B(j,k,i)=0$.

Combining this with Lemma~\ref{identification semi-serious} below, we obtain

  $-(\vv-\vv^{-1})\tilde{f}^{(1)}_{ij}\tilde{g}^-_k\tilde{g}^+_j=
    [\tilde{f}^{(1)}_{kj},\tilde{f}^{(0)}_{ik}]_{\vv}\cdot \tilde{g}^-_k\tilde{g}^+_j\Longrightarrow
    \tilde{f}^{(1)}_{ij}=-[\tilde{f}^{(1)}_{kj},\tilde{f}^{(0)}_{ik}]_{\vv}/(\vv-\vv^{-1}).$
\end{proof}

\begin{Lem}\label{identification semi-serious}
For any $1\leq j<i\leq n$ and $1\leq a,b\leq n$, we have:

\noindent
(a) $\tilde{g}^\epsilon_a\tilde{g}^{\epsilon'}_b=\tilde{g}^{\epsilon'}_b\tilde{g}^{\epsilon}_a$ for any $\epsilon,\epsilon'\in \{\pm\}$.

\noindent
(b) $\tilde{g}^\pm_a \tilde{e}^{(0)}_{ji}=\vv^{\pm\delta_{ai}\mp\delta_{aj}} \tilde{e}^{(0)}_{ji} \tilde{g}^\pm_a$.

\noindent
(c) $\tilde{g}^\pm_a \tilde{f}^{(0)}_{ij}=\vv^{\mp\delta_{ai}\pm \delta_{aj}}\tilde{f}^{(0)}_{ij}\tilde{g}^\pm_a$.

\noindent
(d) $\tilde{g}^\pm_a \tilde{e}^{(-1)}_{ji}=\vv^{\pm\delta_{ai}\mp\delta_{aj}} \tilde{e}^{(-1)}_{ji} \tilde{g}^\pm_a$.

\noindent
(e) $\tilde{g}^\pm_a \tilde{f}^{(1)}_{ij}=\vv^{\mp\delta_{ai}\pm \delta_{aj}}\tilde{f}^{(1)}_{ij}\tilde{g}^\pm_a$.
\end{Lem}

\begin{proof}
First, we note that $t^\pm_{ii}[0]=\tilde{g}^\pm_i$.
Hence, we have $\tilde{g}^\pm_i\tilde{g}^\mp_i=1$, due to relation~(\ref{nRTT1}).

(a) Due to the above observation, it suffices to prove
$\tilde{g}^+_a\tilde{g}^+_b=\tilde{g}^+_b\tilde{g}^+_a$ for $a<b$.
This follows by evaluating the coefficients of $z^0w^1$ in the equality of the matrix coefficients
$\langle v_a\otimes v_b|\cdots|v_a\otimes v_b\rangle$ of both sides of the equality
$((\vv z-\vv^{-1}w)R_\trig(z/w))(T^+(z)\otimes 1)(1\otimes T^+(w))=(1\otimes T^+(w))(T^+(z)\otimes 1)((\vv z-\vv^{-1}w)R_\trig(z/w))$.


(b) Due to the above observation, it suffices to prove
$\tilde{g}^+_a \tilde{e}^{(0)}_{ji}=\vv^{\delta_{ai}-\delta_{aj}} \tilde{e}^{(0)}_{ji} \tilde{g}^+_a$.
This follows by evaluating the coefficients of $z^0w^1$ in the equality of the matrix coefficients
$\langle v_a\otimes v_j|\cdots|v_a\otimes v_i\rangle$ of both sides of the equality
$((\vv z-\vv^{-1}w)R_\trig(z/w))(T^+(z)\otimes 1)(1\otimes T^+(w))=(1\otimes T^+(w))(T^+(z)\otimes 1)((\vv z-\vv^{-1}w)R_\trig(z/w))$.
Note that the cases $a<j, a=j, j<a<i, a=i, a>i$ have to be treated separately.

(c) Due to the above observation, it suffices to prove
$\tilde{g}^-_a \tilde{f}^{(0)}_{ij}=\vv^{\delta_{ai}-\delta_{aj}} \tilde{f}^{(0)}_{ij} \tilde{g}^-_a$.
This follows by evaluating the coefficients of $z^0w^1$ in the equality of the matrix coefficients
$\langle v_i\otimes v_a|\cdots|v_j\otimes v_a\rangle$ of both sides of the equality
$((\vv z-\vv^{-1}w)R_\trig(z/w))(T^-(z)\otimes 1)(1\otimes T^-(w))=(1\otimes T^-(w))(T^-(z)\otimes 1)((\vv z-\vv^{-1}w)R_\trig(z/w))$.
Note that the cases $a<j, a=j, j<a<i, a=i, a>i$ have to be treated separately.

(d) Due to the above observation, it suffices to prove
$\tilde{g}^+_a \tilde{e}^{(-1)}_{ji}=\vv^{\delta_{ai}-\delta_{aj}} \tilde{e}^{(-1)}_{ji} \tilde{g}^+_a$.
This follows by evaluating the coefficients of $z^1w^1$ in the equality of the matrix coefficients
$\langle v_a\otimes v_j|\cdots|v_a\otimes v_i\rangle$ of both sides of the equality
$((\vv z-\vv^{-1}w)R_\trig(z/w))(T^+(z)\otimes 1)(1\otimes T^-(w))=(1\otimes T^-(w))(T^+(z)\otimes 1)((\vv z-\vv^{-1}w)R_\trig(z/w))$.
Note that the cases $a<j, a=j, j<a<i, a=i, a>i$ have to be treated separately.

Let us emphasize that this case is less trivial than part (b), due to the fact that
\begin{equation*}
    [w^1]T^-_{ji}(w)=\tilde{g}^-_j\tilde{e}^{(-1)}_{ji}+
    \sum_{j'<j}\tilde{f}^{(0)}_{jj'}\tilde{g}^-_{j'}\tilde{e}^{(-1)}_{j'i}.
\end{equation*}
Hence, the proof proceeds by induction in $j$, while we also use part (c) from above.

(e) Due to the above observation, it suffices to prove
$\tilde{g}^-_a \tilde{f}^{(1)}_{ij}=\vv^{\delta_{ai}-\delta_{aj}} \tilde{f}^{(1)}_{ij} \tilde{g}^-_a$.
This follows by evaluating the coefficients of $z^0w^0$ in the equality of the matrix coefficients
$\langle v_i\otimes v_a|\cdots|v_j\otimes v_a\rangle$ of both sides of the equality
$((\vv z-\vv^{-1}w)R_\trig(z/w))(T^+(z)\otimes 1)(1\otimes T^-(w))=(1\otimes T^-(w))(T^+(z)\otimes 1)((\vv z-\vv^{-1}w)R_\trig(z/w))$.
Note that the cases $a<j, a=j, j<a<i, a=i, a>i$ have to be treated separately.

Analogously to part (d), this case is less trivial than part (c), due to the fact that
\begin{equation*}
    [z^{-1}]T^+_{ij}(z)=\tilde{f}^{(1)}_{ij}\tilde{g}^+_j+
    \sum_{j'<j}\tilde{f}^{(1)}_{ij'}\tilde{g}^+_{j'}\tilde{e}^{(0)}_{j'j}.
\end{equation*}
Hence, the proof proceeds by induction in $j$, while we also use part (b) from above.
\end{proof}

The following explicit formulas follow immediately from Proposition~\ref{identification serious}.

\begin{Cor}\label{long tilda-commutators}
For any $1\leq j<i\leq n$, we have:
\begin{equation}\label{tilda 1}
\begin{split}
  \tilde{e}^{(0)}_{ji}=
  (\vv-\vv^{-1})^{j-i+1}[\tilde{e}^{(0)}_{i-1,i},[\tilde{e}^{(0)}_{i-2,i-1},\cdots,[\tilde{e}^{(0)}_{j+1,j+2},\tilde{e}^{(0)}_{j,j+1}]_{\vv^{-1}}\cdots]_{\vv^{-1}}]_{\vv^{-1}}=\\
  (\vv-\vv^{-1})^{j-i+1}[[\cdots[\tilde{e}^{(0)}_{i-1,i},\tilde{e}^{(0)}_{i-2,i-1}]_{\vv^{-1}},\cdots,\tilde{e}^{(0)}_{j+1,j+2}]_{\vv^{-1}},\tilde{e}^{(0)}_{j,j+1}]_{\vv^{-1}},
\end{split}
\end{equation}
\begin{equation}\label{tilda 4}
\begin{split}
  \tilde{f}^{(0)}_{ij}=
  (\vv^{-1}-\vv)^{j-i+1}[\tilde{f}^{(0)}_{j+1,j},[\tilde{f}^{(0)}_{j+2,j+1},\cdots,[\tilde{f}^{(0)}_{i-1,i-2},\tilde{f}^{(0)}_{i,i-1}]_{\vv}\cdots]_{\vv}]_{\vv}=\\
  (\vv^{-1}-\vv)^{j-i+1}[[\cdots[\tilde{f}^{(0)}_{j+1,j},\tilde{f}^{(0)}_{j+2,j+1}]_{\vv},\cdots,\tilde{f}^{(0)}_{i-1,i-2}]_\vv,\tilde{f}^{(0)}_{i,i-1}]_{\vv},
\end{split}
\end{equation}
\begin{equation}\label{tilda 2}
\begin{split}
  \tilde{e}^{(-1)}_{ji}=
  (\vv-\vv^{-1})^{j-i+1}[\tilde{e}^{(0)}_{i-1,i},[\tilde{e}^{(0)}_{i-2,i-1},\cdots,[\tilde{e}^{(0)}_{j+1,j+2},\tilde{e}^{(-1)}_{j,j+1}]_{\vv^{-1}}\cdots]_{\vv^{-1}}]_{\vv^{-1}}=\\
  (\vv-\vv^{-1})^{j-i+1}[[\cdots[\tilde{e}^{(0)}_{i-1,i},\tilde{e}^{(0)}_{i-2,i-1}]_{\vv^{-1}},\cdots,\tilde{e}^{(0)}_{j+1,j+2}]_{\vv^{-1}},\tilde{e}^{(-1)}_{j,j+1}]_{\vv^{-1}},
\end{split}
\end{equation}
\begin{equation}\label{tilda 3}
\begin{split}
  \tilde{f}^{(1)}_{ij}=
  (\vv^{-1}-\vv)^{j-i+1}[\tilde{f}^{(1)}_{j+1,j},[\tilde{f}^{(0)}_{j+2,j+1},\cdots,[\tilde{f}^{(0)}_{i-1,i-2},\tilde{f}^{(0)}_{i,i-1}]_{\vv}\cdots]_{\vv}]_{\vv}=\\
  (\vv^{-1}-\vv)^{j-i+1}[[\cdots[\tilde{f}^{(1)}_{j+1,j},\tilde{f}^{(0)}_{j+2,j+1}]_{\vv},\cdots,\tilde{f}^{(0)}_{i-1,i-2}]_\vv,\tilde{f}^{(0)}_{i,i-1}]_{\vv}.
\end{split}
\end{equation}
\end{Cor}

Recall elements $E_{ji}^{(0)},F_{ij}^{(0)},E_{ji}^{(-1)},F_{ij}^{(1)}\in U_\vv(L\ssl_n)$
of~(\ref{long commutators}). Combining Corollary~\ref{long tilda-commutators}
with~(\ref{identification baby}), we get the following result.

\begin{Cor}\label{transfer 1}
(a) We have
\begin{equation}\label{tilda vs nontilda}
\begin{split}
  & \Upsilon^{-1}(\tilde{e}^{(0)}_{ji})=(\vv-\vv^{-1})E^{(0)}_{ji},\
  \Upsilon^{-1}(\tilde{f}^{(0)}_{ij})=-(\vv-\vv^{-1})F^{(0)}_{ij},\\
  & \Upsilon^{-1}(\tilde{e}^{(-1)}_{ji})=-\vv^{-j}(\vv-\vv^{-1})E^{(-1)}_{ji},\
  \Upsilon^{-1}(\tilde{f}^{(1)}_{ij})=\vv^j(\vv-\vv^{-1})F^{(1)}_{ij}.
\end{split}
\end{equation}

\noindent
(b) The right equalities in each of the first four lines of~(\ref{long commutators}) hold.
\end{Cor}

To derive the right equalities of the last two lines of~(\ref{long commutators}),
we introduce
\begin{equation}\label{auxiliary B}
\begin{split}
  & A^+_{ji}:=\sum_{s\geq 1}\sum_{j=j_1<\ldots<j_{s+1}=i} (-1)^{s-1}
  \tilde{e}^{(0)}_{j_1 j_2}\cdots\tilde{e}^{(0)}_{j_s j_{s+1}},\\
  & A^-_{ij}:=\sum_{s\geq 1}\sum_{j=j_1<\ldots<j_{s+1}=i} (-1)^{s-1}
  \tilde{f}^{(0)}_{j_{s+1} j_s}\cdots\tilde{f}^{(0)}_{j_2 j_{1}}
\end{split}
\end{equation}
for $1\leq j<i\leq n$. These elements will play an important role
in~Section~\ref{G10} below.

\begin{Lem}\label{explicit B}
For any $1\leq j<i\leq n$, we have
\begin{equation}\label{tilda 5}
\begin{split}
  A^+_{ji}=
  (\vv-\vv^{-1})^{j-i+1}[\tilde{e}^{(0)}_{i-1,i},[\tilde{e}^{(0)}_{i-2,i-1},\cdots,[\tilde{e}^{(0)}_{j+1,j+2},\tilde{e}^{(0)}_{j,j+1}]_{\vv}\cdots]_{\vv}]_{\vv}=\\
  (\vv-\vv^{-1})^{j-i+1}[[\cdots[\tilde{e}^{(0)}_{i-1,i},\tilde{e}^{(0)}_{i-2,i-1}]_{\vv},\cdots,\tilde{e}^{(0)}_{j+1,j+2}]_{\vv},\tilde{e}^{(0)}_{j,j+1}]_{\vv},
\end{split}
\end{equation}
\begin{equation}\label{tilda 6}
\begin{split}
  A^-_{ij}=
  (\vv^{-1}-\vv)^{j-i+1}[\tilde{f}^{(0)}_{j+1,j},[\tilde{f}^{(0)}_{j+2,j+1},\cdots,[\tilde{f}^{(0)}_{i-1,i-2},\tilde{f}^{(0)}_{i,i-1}]_{\vv^{-1}}\cdots]_{\vv^{-1}}]_{\vv^{-1}}=\\
  (\vv^{-1}-\vv)^{j-i+1}[[\cdots[\tilde{f}^{(0)}_{j+1,j},\tilde{f}^{(0)}_{j+2,j+1}]_{\vv^{-1}},\cdots,\tilde{f}^{(0)}_{i-1,i-2}]_{\vv^{-1}},\tilde{f}^{(0)}_{i,i-1}]_{\vv^{-1}}.
\end{split}
\end{equation}
\end{Lem}

\begin{proof}
We prove~(\ref{tilda 5}) by induction in $i-j$. The result is obvious for $i-j=1$.
To perform the induction step, note that
  $A^+_{ji}=\tilde{e}^{(0)}_{ji}-\sum_{j<k<i}\tilde{e}^{(0)}_{jk}\cdot A^+_{ki}$.
Applying the first equality of~(\ref{tilda 1}) together with the induction assumption, we get
\begin{equation*}
\begin{split}
  & (\vv-\vv^{-1})^{i-j-1}A^+_{ji}=
  [\tilde{e}^{(0)}_{i-1,i},[\tilde{e}^{(0)}_{i-2,i-1},\cdots, [\tilde{e}^{(0)}_{j+1,j+2},\tilde{e}^{(0)}_{j,j+1}]_{\vv^{-1}}\cdots]_{\vv^{-1}}]_{\vv^{-1}}-\\
  & (\vv-\vv^{-1})\sum_{j<k<i}[\tilde{e}^{(0)}_{k-1,k},\cdots,[\tilde{e}^{(0)}_{j+1,j+2},\tilde{e}^{(0)}_{j,j+1}]_{\vv^{-1}}\cdots]_{\vv^{-1}}\cdot
  [\tilde{e}^{(0)}_{i-1,i},\cdots,[\tilde{e}^{(0)}_{k+1,k+2},\tilde{e}^{(0)}_{k,k+1}]_{\vv}\cdots]_{\vv}.
\end{split}
\end{equation*}

Rewriting
  $[\tilde{e}^{(0)}_{i-1,i},X]_{\vv^{\pm 1}}$ as
$\tilde{e}^{(0)}_{i-1,i}\cdot X - \vv^{\pm 1} X\cdot \tilde{e}^{(0)}_{i-1,i}$
and using the equality $[\tilde{e}^{(0)}_{i-1,i},\tilde{e}^{(0)}_{l,l+1}]=0$ for any $l<i-2$
(due to the quadratic Serre relations in $U^\ad_\vv(L\ssl_n)$), we immediately find
\begin{equation*}
\begin{split}
  & (\vv-\vv^{-1})^{i-j-1}A^+_{ji}=
  \left[\vphantom{\sum_j}\tilde{e}^{(0)}_{i-1,i},[\tilde{e}^{(0)}_{i-2,i-1},\cdots, [\tilde{e}^{(0)}_{j+1,j+2},\tilde{e}^{(0)}_{j,j+1}]_{\vv^{-1}}\cdots]_{\vv^{-1}}-(\vv-\vv^{-1})\cdot\right.\\
  & \left.\sum_{j<k<i-1}[\tilde{e}^{(0)}_{k-1,k},\cdots,[\tilde{e}^{(0)}_{j+1,j+2},\tilde{e}^{(0)}_{j,j+1}]_{\vv^{-1}}\cdots]_{\vv^{-1}}\cdot
  [\tilde{e}^{(0)}_{i-2,i-1},\cdots,[\tilde{e}^{(0)}_{k+1,k+2},\tilde{e}^{(0)}_{k,k+1}]_{\vv}\cdots]_{\vv}\right]_\vv=\\
  & [\tilde{e}^{(0)}_{i-1,i},(\vv-\vv^{-1})^{i-j-2}A^+_{j,i-1}]_\vv=
  [\tilde{e}^{(0)}_{i-1,i},[\tilde{e}^{(0)}_{i-2,i-1},\cdots,[\tilde{e}^{(0)}_{j+1,j+2},\tilde{e}^{(0)}_{j,j+1}]_{\vv}\cdots]_{\vv}]_{\vv}.
\end{split}
\end{equation*}
Note that the last equality follows from the induction assumption applied to $A^+_{j,i-1}$.

To prove that $A^+_{ji}$ also equals the rightmost commutator of~(\ref{tilda 5}),
we apply similar arguments to the equality
  $A^+_{ji}=\tilde{e}^{(0)}_{ji}-\sum_{j<k<i}A^+_{jk}\cdot \tilde{e}^{(0)}_{ki}$.
We evaluate the right-hand side by applying the rightmost expression of~(\ref{tilda 1}) to the terms
$\tilde{e}^{(0)}_{ji}, \tilde{e}^{(0)}_{ki}$ and the induction assumption to $A^+_{jk}$.
Rewriting $[X,\tilde{e}^{(0)}_{j,j+1}]_{\vv^{\pm 1}}$ as
  $X\cdot \tilde{e}^{(0)}_{j,j+1} - \vv^{\pm 1}\tilde{e}^{(0)}_{j,j+1}\cdot X$
and taking $\tilde{e}^{(0)}_{j,j+1}$ to the leftmost or the rightmost sides, we get the result.

The proof of~(\ref{tilda 6}) is completely analogous and is left to the interested reader.
\end{proof}

The following result follows by combining~Lemma~\ref{explicit B}
with formula~(\ref{identification baby}).

\begin{Cor}\label{transfer 2}
(a) We have
\begin{equation}\label{B vs EF-tilda}
  \Upsilon^{-1}(A^+_{ji})=(\vv-\vv^{-1})\wt{E}^{(0)}_{ji},\
  \Upsilon^{-1}(A^-_{ij})=-(\vv-\vv^{-1})\wt{F}^{(0)}_{ij}.
\end{equation}

\noindent
(b) The right equalities in the last two lines of~(\ref{long commutators}) hold.
\end{Cor}


\subsection{Computation of $\Delta^\rtt(\tilde{g}^{(\pm 1)}_i)$}
\label{G10}
\

Given a Laurent series $F(z)$, we use $[z^r]F(z)$ to denote the coefficient of $z^r$ in $F(z)$.
In this subsection, we compute explicitly $\Delta^\rtt(\tilde{g}^{(\pm 1)}_i)$,
see~Theorems~\ref{Coproduct of k_i1} and~\ref{Coproduct of k_i-1}.

\begin{Thm}\label{Coproduct of k_i1}
For $1\leq i\leq n$, we have
\begin{equation}\label{k_i1}
\begin{split}
  & \Delta^\rtt(\tilde{g}^{(1)}_i)=
  \tilde{g}^{(1)}_i\otimes \tilde{g}^+_i + \tilde{g}^+_i\otimes \tilde{g}^{(1)}_i+
  \sum_{l>i}\tilde{g}^+_i\tilde{e}^{(0)}_{il}\otimes \tilde{f}^{(1)}_{li}\tilde{g}^+_i+\\
  & \sum_{s\geq 1}\sum_{j_1<\ldots<j_{s+1}=i} (-1)^s\tilde{g}^+_i\tilde{e}^{(0)}_{j_1 j_2}\cdots \tilde{e}^{(0)}_{j_s j_{s+1}}\otimes \tilde{f}^{(1)}_{i j_1}\tilde{g}^+_i+\\
  & \sum_{l>i} \sum_{s\geq 1} \sum_{j_1<\ldots<j_{s+1}=i} (-1)^s\tilde{g}^+_i\tilde{e}^{(0)}_{i l}\tilde{e}^{(0)}_{j_1 j_2}\cdots\tilde{e}^{(0)}_{j_s j_{s+1}}\otimes
  \tilde{f}^{(1)}_{l j_1}\tilde{g}^+_i.
\end{split}
\end{equation}
\end{Thm}

\begin{proof}
Our starting point is the equality
\begin{equation}\label{starting point 1}
  [z^{-1}]T^+_{ii}(z)=\tilde{g}^{(1)}_i+\sum_{j<i}\tilde{f}^{(1)}_{ij}\tilde{g}^+_j\tilde{e}^{(0)}_{ji}.
\end{equation}
We also note that
  $[z^{-1}]T^+_{ij}(z)=\tilde{f}^{(1)}_{ij}\tilde{g}^+_j + \sum_{j'<j} \tilde{f}^{(1)}_{ij'}\tilde{g}^+_{j'}\tilde{e}^{(0)}_{j'i}$
for any $i>j$.
Rewriting this as
  $\tilde{f}^{(1)}_{ij}\tilde{g}^+_j=[z^{-1}]T^+_{ij}(z)-\sum_{j'<j} \tilde{f}^{(1)}_{ij'}\tilde{g}^+_{j'}\tilde{e}^{(0)}_{j'i}$
and applying this formula iteratively, we finally get
\begin{equation}\label{starting point 2}
  \tilde{f}^{(1)}_{ij}\tilde{g}^+_j=
  \sum_{s\geq 1}\sum_{j_1<\ldots<j_s=j} (-1)^{s-1}\left([z^{-1}]T^+_{ij_1}(z)\right) \tilde{e}^{(0)}_{j_1 j_2}\cdots \tilde{e}^{(0)}_{j_{s-1} j_s}.
\end{equation}
Combining formulas~(\ref{starting point 1}) and~(\ref{starting point 2}), we get
\begin{equation}\label{explicit k_i1}
  \tilde{g}^{(1)}_i=[z^{-1}]T^+_{ii}(z)-\sum_{j<i} \left([z^{-1}]T^+_{ij}(z)\right)\cdot A^+_{ji},
\end{equation}
where $A^+_{ji}$ was defined in~(\ref{auxiliary B}).

Thus, it remains to compute explicitly
  $\Delta^\rtt([z^{-1}]T^+_{ii}(z)), \Delta^\rtt([z^{-1}]T^+_{ij}(z)), \Delta^\rtt(A^+_{ji})$
for $i>j$. Evaluating the coefficients of $z^{-1}$ in
  $\Delta^\rtt(T^+_{ii}(z))=\sum_{a=1}^n T^+_{ia}(z)\otimes T^+_{ai}(z)$,
we find
\begin{equation}\label{12.3.1}
\begin{split}
  & \Delta^\rtt([z^{-1}]T^+_{ii}(z))=
  \sum_{j<i}\tilde{f}^{(1)}_{ij}\tilde{g}^+_j\otimes \tilde{g}^+_j\tilde{e}^{(0)}_{ji}+
  \sum_{j'<j<i}\tilde{f}^{(1)}_{ij'}\tilde{g}^+_{j'}\tilde{e}^{(0)}_{j'j}\otimes \tilde{g}^+_j\tilde{e}^{(0)}_{ji}+\\
  & \tilde{g}^{(1)}_i\otimes \tilde{g}^+_i + \tilde{g}^+_i\otimes \tilde{g}^{(1)}_i+
  \sum_{j<i}\tilde{f}^{(1)}_{ij}\tilde{g}^+_j\tilde{e}^{(0)}_{ji}\otimes \tilde{g}^+_i+
  \sum_{j<i}\tilde{g}^+_i\otimes \tilde{f}^{(1)}_{ij}\tilde{g}^+_j\tilde{e}^{(0)}_{ji}+\\
  & \sum_{l>i}\tilde{g}^+_i\tilde{e}^{(0)}_{il}\otimes \tilde{f}^{(1)}_{li}\tilde{g}^+_i+
  \sum_{l>i}^{j<i}\tilde{g}^+_i\tilde{e}^{(0)}_{il}\otimes \tilde{f}^{(1)}_{lj}\tilde{g}^+_j\tilde{e}^{(0)}_{ji},
\end{split}
\end{equation}
where the first, second, and third lines in the right-hand side correspond to the contributions
arising from the cases $a<i$, $a=i$, and $a>i$, respectively.

Evaluating the coefficients of $z^{-1}$ in
  $\Delta^\rtt(T^+_{ij}(z))=\sum_{a=1}^n T^+_{ia}(z)\otimes T^+_{aj}(z)$,
we find
\begin{equation}\label{12.3.2}
\begin{split}
  & \Delta^\rtt([z^{-1}]T^+_{ij}(z))=
  \sum_{j'<j}\tilde{f}^{(1)}_{ij'}\tilde{g}^+_{j'}\otimes \tilde{g}^+_{j'}\tilde{e}^{(0)}_{j'j}+
  \sum_{j''<j'<j}\tilde{f}^{(1)}_{ij''}\tilde{g}^+_{j''}\tilde{e}^{(0)}_{j''j'}\otimes \tilde{g}^+_{j'}\tilde{e}^{(0)}_{j'j}+\\
  & \tilde{f}^{(1)}_{ij}\tilde{g}^+_j\otimes \tilde{g}^+_j +
  \sum_{j'<j}\tilde{f}^{(1)}_{ij'}\tilde{g}^+_{j'}\tilde{e}^{(0)}_{j'j}\otimes \tilde{g}^+_{j}+
  \tilde{g}^+_i\otimes \tilde{f}^{(1)}_{ij}\tilde{g}^+_j+
  \sum_{j'<j}\tilde{g}^+_i\otimes \tilde{f}^{(1)}_{ij'}\tilde{g}^+_{j'}\tilde{e}^{(0)}_{j'j}+\\
  & \sum_{l>i}\tilde{g}^+_i\tilde{e}^{(0)}_{il}\otimes \tilde{f}^{(1)}_{lj}\tilde{g}^+_j+
  \sum_{l>i}^{j'<j}\tilde{g}^+_i\tilde{e}^{(0)}_{il}\otimes \tilde{f}^{(1)}_{lj'}\tilde{g}^+_{j'}\tilde{e}^{(0)}_{j'j},
\end{split}
\end{equation}
where the first, second, and third lines in the right-hand side correspond to
the contributions arising
from $a<j$, $a=j\ \mathrm{or}\ i$, and $a>i$, respectively.
Note that for $j<a<i$ both $T^+_{ia}(z), T^+_{aj}(z)$ contain only negative powers of $z$
and hence do not contribute above.

Finally, let us compute the coproduct of $A^+_{ji}$.

\begin{Lem}\label{12.3.3}
We have
\begin{equation*}
    \Delta^\rtt(A^+_{ji})=\sum_{s\geq 1}\sum_{j=j_1<\ldots<j_{s+1}=i}\sum_{r=1}^{s+1}(-1)^{s-1}
    \tilde{e}^{(0)}_{j_r j_{r+1}}\cdots \tilde{e}^{(0)}_{j_s j_{s+1}}\otimes \tilde{e}^{(0)}_{j_1 j_2}\cdots \tilde{e}^{(0)}_{j_{r-1} j_r}(\tilde{g}^+_{j_r})^{-1}\tilde{g}^+_i.
\end{equation*}
\end{Lem}

\begin{proof}
We prove this by induction in $i-j$. The base of induction $i=j+1$ follows
from the equality $A^+_{j,j+1}=\tilde{e}^{(0)}_{j,j+1}$ and Lemma~\ref{12.3.6} below.
To perform the induction step, note that
\begin{equation}\label{12.3.4}
  A^+_{ji}=\tilde{e}^{(0)}_{ji}-\sum_{j<j'<i}\tilde{e}^{(0)}_{jj'}A^+_{j'i}.
\end{equation}


Next, we compute the coproduct of $\tilde{e}^{(0)}_{ji}$.

\begin{Lem}\label{12.3.6}
We have
\begin{equation*}
    \Delta^\rtt(\tilde{e}^{(0)}_{ji})=
    1\otimes \tilde{e}^{(0)}_{ji}+\tilde{e}^{(0)}_{ji}\otimes (\tilde{g}^+_j)^{-1}\tilde{g}^+_i+
    \sum_{j<a<i}\tilde{e}^{(0)}_{ja}\otimes (\tilde{g}^+_j)^{-1}\tilde{g}^+_{a}\tilde{e}^{(0)}_{ai}.
\end{equation*}
\end{Lem}

\begin{proof}
First, let us note that $\tilde{g}^+_j=[z^0]T^+_{jj}(z)$.
Thus,
\begin{equation*}
  \Delta^\rtt(\tilde{g}^+_j)=[z^0]\left(\sum_{a=1}^n T^+_{ja}(z)\otimes T^+_{aj}(z)\right)=[z^0](T^+_{jj}(z)\otimes T^+_{jj}(z))=\tilde{g}^+_j\otimes \tilde{g}^+_j.
\end{equation*}
We also note that $[z^0]T^+_{ji}(z)=\tilde{g}^+_j\tilde{e}^{(0)}_{ji}$.
Hence, we have
\begin{equation*}
\begin{split}
  & \Delta^\rtt(\tilde{g}^+_j\tilde{e}^{(0)}_{ji})=
  [z^0]\left(T^+_{jj}(z)\otimes T^+_{ji}(z)+T^+_{ji}(z)\otimes T^+_{ii}(z)+\sum_{j<a<i} T^+_{ja}(z)\otimes T^+_{ai}(z)\right)=\\
  & \tilde{g}^+_j\otimes \tilde{g}^+_j\tilde{e}^{(0)}_{ji}+\tilde{g}^+_j\tilde{e}^{(0)}_{ji}\otimes \tilde{g}^+_i+
  \sum_{j<a<i}\tilde{g}^+_j\tilde{e}^{(0)}_{ja}\otimes \tilde{g}^+_a\tilde{e}^{(0)}_{ai}.
\end{split}
\end{equation*}
Note that in the first equality we used $[z^0](T^+_{ja}(z)\otimes T^+_{ai}(z))=0$ for $a<j$ or $a>i$.

Evaluating
  $\Delta^\rtt(\tilde{e}^{(0)}_{ji})=\Delta^\rtt(\tilde{g}^+_j)^{-1}\Delta^\rtt(\tilde{g}^+_j\tilde{e}^{(0)}_{ji})$
via these formulas completes our proof.
\end{proof}

Combining~(\ref{12.3.4}) with Lemma~\ref{12.3.6} and applying the induction assumption to
$\Delta^\rtt(A^+_{j'i})$, we immediately get the formula for $\Delta^\rtt(A^+_{ji})$
of Lemma~\ref{12.3.3}.
\end{proof}

Combining~(\ref{explicit k_i1}--\ref{12.3.2}) with Lemma~\ref{12.3.3}, we get~(\ref{k_i1})
after tedious computations.
\end{proof}

\begin{Thm}\label{Coproduct of k_i-1}
For $1\leq i\leq n$, we have
\begin{equation}\label{k_i-1}
\begin{split}
  & \Delta^\rtt(\tilde{g}^{(-1)}_i)=
  \tilde{g}^{(-1)}_i\otimes \tilde{g}^-_i + \tilde{g}^-_i\otimes \tilde{g}^{(-1)}_i+
  \sum_{l>i}\tilde{g}^-_i\tilde{e}^{(-1)}_{il}\otimes \tilde{f}^{(0)}_{li}\tilde{g}^-_i+\\
  & \sum_{s\geq 1}\sum_{j_1<\ldots<j_{s+1}=i} (-1)^s\tilde{g}^-_i\tilde{e}^{(-1)}_{ji}\otimes \tilde{f}^{(0)}_{j_{s+1} j_s}\cdots \tilde{f}^{(0)}_{j_2 j_1}\tilde{g}^-_i+\\
  & \sum_{l>i} \sum_{s\geq 1} \sum_{j_1<\ldots<j_{s+1}=i} (-1)^s\tilde{g}^-_i\tilde{e}^{(-1)}_{jl}\otimes
  \tilde{f}^{(0)}_{j_{s+1} j_s}\cdots \tilde{f}^{(0)}_{j_2 j_1}\tilde{f}^{(0)}_{li}\tilde{g}^-_i.
\end{split}
\end{equation}
\end{Thm}

\begin{proof}
Our starting point is the equality
\begin{equation}\label{starting point 3}
  [z]T^-_{ii}(z)=\tilde{g}^{(-1)}_i+\sum_{j<i}\tilde{f}^{(0)}_{ij}\tilde{g}^-_j\tilde{e}^{(-1)}_{ji}.
\end{equation}
We also note that
  $[z]T^-_{ji}(z)=\tilde{g}^-_j\tilde{e}^{(-1)}_{ji} + \sum_{j'<j} \tilde{f}^{(0)}_{jj'}\tilde{g}^-_{j'}\tilde{e}^{(-1)}_{j'i}$
for any $i>j$.
Rewriting this as
  $\tilde{g}^-_j\tilde{e}^{(-1)}_{ji} =[z]T^-_{ji}(z)-\sum_{j'<j} \tilde{f}^{(0)}_{jj'}\tilde{g}^-_{j'}\tilde{e}^{(-1)}_{j'i}$
and applying this formula iteratively, we finally get
\begin{equation}\label{starting point 4}
  \tilde{g}^-_j\tilde{e}^{(-1)}_{ji}=\sum_{s\geq 1}\sum_{j_1<\ldots<j_s=j} (-1)^{s-1}
  \tilde{f}^{(0)}_{j_s j_{s-1}}\cdots \tilde{f}^{(0)}_{j_2 j_1}\cdot \left([z]T^-_{j_1i}(z)\right).
\end{equation}
Combining formulas~(\ref{starting point 3}) and~(\ref{starting point 4}), we get
\begin{equation}\label{explicit k_i-1}
  \tilde{g}^{(-1)}_i=[z]T^-_{ii}(z)-\sum_{j<i} A^-_{ij}\cdot \left([z]T^-_{ji}(z)\right),
\end{equation}
where $A^-_{ij}$ was defined in~(\ref{auxiliary B}).

Thus, it remains to compute explicitly
$\Delta^\rtt([z]T^-_{ii}(z)), \Delta^\rtt([z]T^-_{ji}(z)), \Delta^\rtt(A^-_{ij})$
for $i>j$.
Evaluating the coefficients of $z^1$ in
  $\Delta^\rtt(T^-_{ii}(z))=\sum_{a=1}^n T^-_{ia}(z)\otimes T^-_{ai}(z)$,
we find
\begin{equation}\label{12.3.7}
\begin{split}
  & \Delta^\rtt([z]T^-_{ii}(z))=
  \sum_{j<i}\tilde{f}^{(0)}_{ij}\tilde{g}^-_j\otimes \tilde{g}^-_j\tilde{e}^{(-1)}_{ji}+
  \sum_{j'<j<i} \tilde{f}^{(0)}_{ij}\tilde{g}^-_j\otimes \tilde{f}^{(0)}_{jj'}\tilde{g}^-_{j'}\tilde{e}^{(-1)}_{j'i}+\\
  & \tilde{g}^-_i\otimes \tilde{g}^{(-1)}_i+\tilde{g}^{(-1)}_i\otimes \tilde{g}^-_i+
  \sum_{j<i} \tilde{g}^-_i\otimes \tilde{f}^{(0)}_{ij}\tilde{g}^-_j\tilde{e}^{(-1)}_{ji}+
  \sum_{j<i} \tilde{f}^{(0)}_{ij}\tilde{g}^-_j\tilde{e}^{(-1)}_{ji}\otimes \tilde{g}^-_i+\\
  & \sum_{l>i} \tilde{g}^-_i\tilde{e}^{(-1)}_{il}\otimes \tilde{f}^{(0)}_{li}\tilde{g}^-_i+
  \sum_{l>i}^{j<i} \tilde{f}^{(0)}_{ij}\tilde{g}^-_j\tilde{e}^{(-1)}_{jl}\otimes \tilde{f}^{(0)}_{li}\tilde{g}^-_i,
\end{split}
\end{equation}
where the first, second, and third lines in the right-hand side correspond to the contributions
arising from the cases $a<i$, $a=i$, and $a>i$, respectively.

Evaluating the coefficients of $z^1$ in
  $\Delta^\rtt(T^-_{ji}(z))=\sum_{a=1}^n T^-_{ja}(z)\otimes T^-_{ai}(z)$,
we find
\begin{equation}\label{12.3.8}
\begin{split}
  & \Delta^\rtt([z]T^-_{ji}(z))=
  \sum_{j'<j} \tilde{f}^{(0)}_{jj'}\tilde{g}^-_{j'}\otimes \tilde{g}^-_{j'}\tilde{e}^{(-1)}_{j'i}+
  \sum_{j''<j'<j} \tilde{f}^{(0)}_{jj'}\tilde{g}^-_{j'}\otimes \tilde{f}^{(0)}_{j'j''}\tilde{g}^-_{j''}\tilde{e}^{(-1)}_{j''i}+\\
  & \tilde{g}^-_j\otimes \tilde{g}^-_j\tilde{e}^{(-1)}_{ji}+\sum_{j'<j} \tilde{g}^-_j\otimes \tilde{f}^{(0)}_{jj'}\tilde{g}^-_{j'}\tilde{e}^{(-1)}_{j'i}+
  \tilde{g}^-_j\tilde{e}^{(-1)}_{ji}\otimes \tilde{g}^-_i +\sum_{j'<j} \tilde{f}^{(0)}_{jj'}\tilde{g}^-_{j'}\tilde{e}^{(-1)}_{j'i}\otimes \tilde{g}^-_i+\\
  & \sum_{l>i}\tilde{g}^-_j\tilde{e}^{(-1)}_{jl}\otimes \tilde{f}^{(0)}_{li}\tilde{g}^-_i+
  \sum_{l>i}^{j'<j}\tilde{f}^{(0)}_{jj'}\tilde{g}^-_{j'}\tilde{e}^{(-1)}_{j'l}\otimes \tilde{f}^{(0)}_{li}\tilde{g}^-_i,
\end{split}
\end{equation}
where the first, second, and third lines in the right-hand side correspond to the contributions
arising from $a<j$, $a=j\ \mathrm{or}\ i$, and $a>i$, respectively.
Note that for $j<a<i$ both $T^-_{ja}(z), T^-_{ai}(z)$ contain only positive powers of $z$
and hence do not contribute above.

Finally, let us compute the coproduct of $A^-_{ij}$.

\begin{Lem}\label{12.3.9}
We have
\begin{equation*}
    \Delta^\rtt(A^-_{ij})=
    \sum_{s\geq 1}\sum_{j=j_1<\ldots<j_{s+1}=i}\sum_{r=1}^{s+1}(-1)^{s-1}
    \tilde{g}^-_i(\tilde{g}^-_{j_r})^{-1}\tilde{f}^{(0)}_{j_r j_{r-1}}\cdots \tilde{f}^{(0)}_{j_2 j_1}\otimes
    \tilde{f}^{(0)}_{j_{s+1} j_s}\cdots \tilde{f}^{(0)}_{j_{r+1} j_r}.
\end{equation*}
\end{Lem}

\begin{proof}
We prove this by induction in $i-j$.
The base of induction $i=j+1$ follows
from the equality $A^-_{j+1,j}=\tilde{f}^{(0)}_{j+1,j}$ and Lemma~\ref{12.3.11} below.
To perform the induction step, note that
\begin{equation}\label{12.3.10}
  A^-_{ij}=\tilde{f}^{(0)}_{ij}-\sum_{j<j'<i}A^-_{ij'}\tilde{f}^{(0)}_{j'j}.
\end{equation}

Next, we compute the coproduct of $\tilde{f}^{(0)}_{ij}$.

\begin{Lem}\label{12.3.11}
We have
\begin{equation*}
    \Delta^\rtt(\tilde{f}^{(0)}_{ij})=
    \tilde{f}^{(0)}_{ij}\otimes 1 + \tilde{g}^-_i(\tilde{g}^-_j)^{-1}\otimes \tilde{f}^{(0)}_{ij}+
    \sum_{j<a<i}\tilde{f}^{(0)}_{ia}\tilde{g}^-_{a}(\tilde{g}^-_j)^{-1}\otimes \tilde{f}^{(0)}_{aj}.
\end{equation*}
\end{Lem}

\begin{proof}
First, let us note that $\tilde{g}^-_j=[z^0]T^-_{jj}(z)$. Thus,
\begin{equation*}
  \Delta^\rtt(\tilde{g}^-_j)=[z^0]\left(\sum_{a=1}^n T^-_{ja}(z)\otimes T^-_{aj}(z)\right)=
  [z^0](T^-_{jj}(z)\otimes T^-_{jj}(z))=\tilde{g}^-_j\otimes \tilde{g}^-_j.
\end{equation*}
We also note that $[z^0]T^-_{ij}(z)=\tilde{f}^{(0)}_{ij}\tilde{g}^-_j$.
Hence, we have
\begin{equation*}
\begin{split}
  & \Delta^\rtt(\tilde{f}^{(0)}_{ij}\tilde{g}^-_j)=
  [z^0]\left(T^-_{ij}(z)\otimes T^-_{jj}(z)+T^-_{ii}(z)\otimes T^-_{ij}(z)+\sum_{j<a<i} T^-_{ia}(z)\otimes T^-_{aj}(z)\right)=\\
  & \tilde{f}^{(0)}_{ij}\tilde{g}^-_j\otimes \tilde{g}^-_j+ \tilde{g}^-_i\otimes \tilde{f}^{(0)}_{ij}\tilde{g}^-_j+
  \sum_{j<a<i}\tilde{f}^{(0)}_{ia}\tilde{g}^-_a\otimes \tilde{f}^{(0)}_{aj}\tilde{g}^-_j.
\end{split}
\end{equation*}
Note that in the first equality we used $[z^0](T^-_{ia}(z)\otimes T^-_{aj}(z))=0$ for $a<j$ or $a>i$.

Evaluating
  $\Delta^\rtt(\tilde{f}^{(0)}_{ij})=\Delta^\rtt(\tilde{f}^{(0)}_{ij}\tilde{g}^-_j)\Delta^\rtt(\tilde{g}^-_j)^{-1}$
via these formulas completes our proof.
\end{proof}

Combining~(\ref{12.3.10}) with Lemma~\ref{12.3.11} and applying the induction assumption to
$\Delta^\rtt(A^-_{ij'})$, we immediately get the formula for $\Delta^\rtt(A^-_{ij})$
of Lemma~\ref{12.3.9}.
\end{proof}

Combining~(\ref{explicit k_i-1}--\ref{12.3.8}) with Lemma~\ref{12.3.9}, we get~(\ref{k_i-1})
after tedious computations.
\end{proof}

For $1\leq i\leq n$, define $H_{i,\pm 1}\in U^\rtt(\ssl_n)$ via
$H_{i,\pm1}:=(\tilde{g}^\pm_i)^{-1}\tilde{g}^{(\pm1)}_i.$
Recall the elements $A^+_{ji}$ and $A^-_{ij}$ of~(\ref{auxiliary B}).
Combining Theorems~\ref{Coproduct of k_i1},~\ref{Coproduct of k_i-1}
with Lemma~\ref{identification semi-serious} and the formula
  $\Delta^\rtt(\tilde{g}^\pm_i)=\tilde{g}^\pm_i\otimes \tilde{g}^\pm_i$,
we get the following expressions for $\Delta^\rtt(H_{i,\pm1})$.

\begin{Cor}\label{Coproduct of H_i1-1}
We have
\begin{multline}\label{H_i1}
  \Delta^\rtt(H_{i,1})=H_{i,1}\otimes 1+1\otimes H_{i,1}+\vv^{-1}\sum_{l>i} \tilde{e}^{(0)}_{il}\otimes \tilde{f}^{(1)}_{li}-
  \vv\sum_{j<i} A^+_{ji}\otimes \tilde{f}^{(1)}_{ij}-\sum_{l>i}^{j<i}\tilde{e}^{(0)}_{il}A^+_{ji}\otimes \tilde{f}^{(1)}_{lj},
\end{multline}
\begin{multline}\label{H_i-1}
  \Delta^\rtt(H_{i,-1})=H_{i,-1}\otimes 1+1\otimes H_{i,-1}+\vv\sum_{l>i} \tilde{e}^{(-1)}_{il}\otimes \tilde{f}^{(0)}_{li}-
  \vv^{-1}\sum_{j<i} \tilde{e}^{(-1)}_{ji}\otimes A^-_{ij}-\sum_{l>i}^{j<i}\tilde{e}^{(-1)}_{jl}\otimes A^-_{ij}\tilde{f}^{(0)}_{li}.
\end{multline}
\end{Cor}


\subsection{Proof of formula~(\ref{h_i1})}
\

Recall the Hopf algebra embedding
  $\Upsilon\colon U_\vv(L\ssl_n)\hookrightarrow U^\rtt(\ssl_n)$
of Theorem~\ref{RTT for sln} (see also Remark~\ref{embedding}). It is easy to see that
\begin{equation*}
  \Upsilon(h_{i,1})=\frac{H_{i+1,1}-H_{i,1}}{\vv^i(\vv-\vv^{-1})}.
\end{equation*}
Combining~Corollaries~\ref{transfer 1},~\ref{transfer 2} with
formula~(\ref{H_i1}) and the fact
that $\Upsilon$ intertwines $\Delta$ and $\Delta^\rtt$, we immediately get
\begin{equation}\label{proof h_i1}
\begin{split}
  & \Delta(h_{i,1})-h_{i,1}\otimes 1-1\otimes h_{i,1}=
  \vv^{-i}(\vv-\vv^{-1})^{-1}\times\\
  & \left(\vv^i(\vv-\vv^{-1})^2\sum_{l>i+1}E^{(0)}_{i+1,l}\otimes F^{(1)}_{l,i+1}-
  (\vv-\vv^{-1})^2\sum_{k<i+1}\vv^{k+1}\wt{E}^{(0)}_{k,i+1}\otimes F^{(1)}_{i+1,k}-\right.\\
  & (\vv-\vv^{-1})^3\sum_{k<i+1<l} \vv^k E^{(0)}_{i+1,l}\wt{E}^{(0)}_{k,i+1}\otimes F^{(1)}_{lk}-
  \vv^{i-1}(\vv-\vv^{-1})^2\sum_{l>i}E^{(0)}_{il}\otimes F^{(1)}_{li}+\\
  & \left.(\vv-\vv^{-1})^2\sum_{k<i} \vv^{k+1} \wt{E}^{(0)}_{ki}\otimes F^{(1)}_{ik}+
  (\vv-\vv^{-1})^3\sum_{k<i<l} \vv^k E^{(0)}_{il}\wt{E}^{(0)}_{ki}\otimes F^{(1)}_{lk}\right).
\end{split}
\end{equation}
This formula implies~(\ref{h_i1}) after the following simplifications:
\begin{multline*}
  \sum_{k<i<l} \vv^k E^{(0)}_{il}\wt{E}^{(0)}_{ki}\otimes F^{(1)}_{lk}-
  \sum_{k<i+1<l} \vv^k E^{(0)}_{i+1,l}\wt{E}^{(0)}_{k,i+1}\otimes F^{(1)}_{lk}=\\
  \sum^{k<i}_{l>i+1} \vv^k (E^{(0)}_{il}\wt{E}^{(0)}_{ki}-E^{(0)}_{i+1,l}\wt{E}^{(0)}_{k,i+1})\otimes F^{(1)}_{lk}+
  \sum_{k<i} \vv^k E^{(0)}_{i,i+1}\wt{E}^{(0)}_{ki}\otimes F^{(1)}_{i+1,k}-
  \vv^i\sum_{l>i+1}E^{(0)}_{i+1,l}\wt{E}^{(0)}_{i,i+1}\otimes F^{(1)}_{li},
\end{multline*}
\begin{equation*}
\begin{split}
  & -\vv^{-1}\sum_{l>i} E^{(0)}_{il}\otimes F^{(1)}_{li}-(\vv-\vv^{-1})\sum_{l>i+1}E^{(0)}_{i+1,l}E^{(0)}_{i,i+1}\otimes F^{(1)}_{li}=\\
  & -\vv^{-1}E^{(0)}_{i,i+1}\otimes F^{(1)}_{i+1,i}+\vv^{-2}\sum_{l>i+1}[E^{(0)}_{i,i+1}, E^{(0)}_{i+1,l}]_{\vv^3}\otimes F^{(1)}_{li},
\end{split}
\end{equation*}
\begin{equation*}
\begin{split}
  & -\sum_{k<i+1}\vv^{k+1-i}\wt{E}^{(0)}_{k,i+1}\otimes F^{(1)}_{i+1,k}+(\vv-\vv^{-1})\sum_{k<i} \vv^{k-i} E^{(0)}_{i,i+1}\wt{E}^{(0)}_{ki}\otimes F^{(1)}_{i+1,k}=\\
  & -\vv E^{(0)}_{i,i+1}\otimes F^{(1)}_{i+1,i}-\sum_{k<i}\vv^{k-i-1}[E^{(0)}_{i,i+1},\wt{E}^{(0)}_{ki}]_{\vv^3}\otimes F^{(1)}_{i+1,k},
\end{split}
\end{equation*}
where in the second and third equalities we used
\begin{equation*}
  E^{(0)}_{il}=[E^{(0)}_{i+1,l}, E^{(0)}_{i,i+1}]_{\vv^{-1}},\
  \wt{E}^{(0)}_{k,i+1}=[E^{(0)}_{i,i+1},\wt{E}^{(0)}_{ki}]_\vv.
\end{equation*}


\subsection{Proof of formula~(\ref{h_i-1})}
\

The proof of~(\ref{h_i-1}) is completely analogous and is based on the formula
\begin{equation*}
  \Upsilon(h_{i,-1})=\frac{H_{i,-1}-H_{i+1,-1}}{\vv^{-i}(\vv-\vv^{-1})}.
\end{equation*}
Combining this with Corollaries~\ref{transfer 1},~\ref{transfer 2}, formula~(\ref{H_i-1})
and the fact that $\Upsilon$ intertwines $\Delta$ and $\Delta^\rtt$, one derives~(\ref{h_i-1}).
The computations are similar to the above proof of~(\ref{h_i1}) and are left to the interested reader.


\subsection{Proof of formula~(\ref{e_i-1})}
\

Recall that $[h_{i,-1},e_{i,0}]=[2]_\vv\cdot e_{i,-1}$, so that
\begin{equation*}
  \Delta(e_{i,-1})=[2]_\vv^{-1}\cdot [\Delta(h_{i,-1}),\Delta(e_{i,0})]=
  [2]_\vv^{-1}\cdot [\Delta(h_{i,-1}), 1\otimes e_{i,0}+e_{i,0}\otimes \psi^+_{i,0}].
\end{equation*}
Applying formula~(\ref{h_i-1}) to $\Delta(h_{i,-1})$ and
using~Lemma~\ref{lemma e i-1} below, we recover~(\ref{e_i-1}).

\begin{Lem}\label{lemma e i-1}
For $k<i$ and $l>i+1$, the following equalities hold:

\noindent
(a) $[F^{(0)}_{l,i+1},e_{i,0}]=0$.

\noindent
(b) $[\wt{F}^{(0)}_{ik},e_{i,0}]=0$.

\noindent
(c) $[F^{(0)}_{li},e_{i,0}]=-F^{(0)}_{l,i+1}\psi^-_{i,0}$.

\noindent
(d) $[\wt{F}^{(0)}_{i+1,k},e_{i,0}]=\vv^{-1}\wt{F}^{(0)}_{ik}\psi^-_{i,0}$.

\noindent
(e) $[[F^{(0)}_{l,i+1},F^{(0)}_{i+1,i}]_{\vv^{-3}},e_{i,0}]=\frac{1-\vv^{-4}}{\vv-\vv^{-1}}F^{(0)}_{l,i+1}\psi^-_{i,0}-\frac{1-\vv^{-2}}{\vv-\vv^{-1}}F^{(0)}_{l,i+1}\psi^+_{i,0}$.

\noindent
(f) $[[\wt{F}^{(0)}_{ik},F^{(0)}_{i+1,i}]_{\vv^{-3}},e_{i,0}]=
     \frac{1-\vv^{-4}}{\vv-\vv^{-1}}\wt{F}^{(0)}_{ik}\psi^-_{i,0}-\frac{1-\vv^{-2}}{\vv-\vv^{-1}}\wt{F}^{(0)}_{ik}\psi^+_{i,0}$.

\noindent
(g) $[E^{(-1)}_{i+1,l},e_{i,0}]_\vv=\vv E^{(-1)}_{il}$.

\noindent
(h) $[E^{(-1)}_{ki},e_{i,0}]_\vv=-\vv E^{(-1)}_{k,i+1}$.

\noindent
(i) $[E^{(-1)}_{il},e_{i,0}]_{\vv^{-1}}=0$.

\noindent
(j) $[E^{(-1)}_{k,i+1},e_{i,0}]_{\vv^{-1}}=0$.

\noindent
(k) $[E^{(-1)}_{kl},e_{i,0}]=0$.
\end{Lem}

\begin{proof}
Recall that $[f_{j,0},e_{i,0}]=\frac{\delta_{ji}}{\vv-\vv^{-1}}(\psi^-_{i,0}-\psi^+_{i,0})$.

Parts (a, b) are obvious as $e_{i,0}$ commutes with
$f_{i+1,0},\ldots,f_{l-1,0}$ and $f_{k,0},\ldots,f_{i-1,0}$.
Combining (a, b) with equalities $F^{(0)}_{li}=[f_{i,0},F^{(0)}_{l,i+1}]_\vv$ and
$\wt{F}^{(0)}_{i+1,k}=[\wt{F}^{(0)}_{ik},f_{i,0}]_{\vv^{-1}}$, we get
  $[F^{(0)}_{li},e_{i,0}]=[\frac{\psi^-_{i,0}-\psi^+_{i,0}}{\vv-\vv^{-1}},F^{(0)}_{l,i+1}]_\vv=-F^{(0)}_{l,i+1}\psi^-_{i,0}$
and
  $[\wt{F}^{(0)}_{i+1,k},e_{i,0}]=
   [\wt{F}^{(0)}_{ik},\frac{\psi^-_{i,0}-\psi^+_{i,0}}{\vv-\vv^{-1}}]_{\vv^{-1}}=\vv^{-1}\wt{F}^{(0)}_{ik}\psi^-_{i,0}$,
which proves parts (c, d). Parts (e, f) also follow immediately from (a, b).

(g) Due to the quadratic Serre relations $e_{i,0}$ commutes with
$e_{i+2,0},\ldots,e_{l-1,0}$, hence, also with $E^{(0)}_{i+2,l}$.
Meanwhile, we have $[e_{i+1,-1}, e_{i,0}]_\vv=\vv[e_{i+1,0},e_{i,-1}]_{\vv^{-1}}$,
due to~(\ref{U2}). Thus,
  $[E^{(-1)}_{i+1,l},e_{i,0}]_\vv=[[E^{(0)}_{i+2,l},e_{i+1,-1}]_{\vv^{-1}},e_{i,0}]_\vv=
   [E^{(0)}_{i+2,l},\vv[e_{i+1,0},e_{i,-1}]_{\vv^{-1}}]_{\vv^{-1}}=\vv E^{(-1)}_{il}$.

(h) We have $[E^{(-1)}_{ki},e_{i,0}]_\vv=-\vv[e_{i,0},E^{(-1)}_{ki}]_{\vv^{-1}}=-\vv E^{(-1)}_{k,i+1}$.

(i) Note that
  $[[e_{i+1,0},e_{i,-1}]_{\vv^{-1}},e_{i,0}]_{\vv^{-1}}=\vv^{-1}[[e_{i+1,-1},e_{i,0}]_{\vv},e_{i,0}]_{\vv^{-1}}=0$,
due to~(\ref{U2}) and~(\ref{U7}). Since also $e_{i,0}$ commutes with $e_{i+2,0},\ldots,e_{l-1,0}$, we get
$[E^{(-1)}_{il},e_{i,0}]_{\vv^{-1}}=0$.

(j) As in (i), $[E^{(-1)}_{k,i+1},e_{i,0}]_{\vv^{-1}}=0$ follows from
$[[e_{i,0},e_{i-1,0}]_{\vv^{-1}},e_{i,0}]_{\vv^{-1}}=0$, due to~(\ref{U7}).

(k) Comparing the matrix coefficients $\langle v_i\otimes v_k|\cdots|v_{i+1}\otimes v_l\rangle$ of both sides of the equality
$R_\trig(z/w)(T^+(z)\otimes 1)(1\otimes T^-(w))=(1\otimes T^-(w))(T^+(z)\otimes 1)R_\trig(z/w)$, we get
\begin{equation*}
\begin{split}
  & (z-w)T^+_{i,i+1}(z)T^-_{kl}(w)+(\vv-\vv^{-1})wT^+_{k,i+1}(z)T^-_{il}(w)=\\
  & (z-w)T^-_{kl}(w)T^+_{i,i+1}(z)+(\vv-\vv^{-1})wT^-_{k,i+1}(w)T^+_{il}(z).
\end{split}
\end{equation*}
Evaluating the coefficients of $z^1w^1$ in both sides of this equality, we find
\begin{equation*}
   [\tilde{g}^+_i\tilde{e}^{(0)}_{i,i+1},\tilde{g}^-_k\tilde{e}^{(-1)}_{kl}+\sum_{j<k}\tilde{f}^{(0)}_{kj}\tilde{g}^-_j\tilde{e}^{(-1)}_{jl}]=0.
\end{equation*}
Hence, by induction in $k$, we find
  $[\tilde{e}^{(0)}_{i,i+1},\tilde{e}^{(-1)}_{kl}]=0$,
which implies $[E^{(-1)}_{kl},e_{i,0}]=0$.
\end{proof}

This completes our proof of~(\ref{e_i-1}).


\subsection{Proof of formula~(\ref{f_i1})}
\

Recall that $[h_{i,1},f_{i,0}]=-[2]_\vv\cdot f_{i,1}$, so that
\begin{equation*}
  \Delta(f_{i,1})=-[2]_\vv^{-1}\cdot [\Delta(h_{i,1}),\Delta(f_{i,0})]=
  -[2]_\vv^{-1}\cdot [\Delta(h_{i,1}), f_{i,0}\otimes 1 +\psi^-_{i,0}\otimes f_{i,0}].
\end{equation*}
Applying formula~(\ref{h_i1}) to $\Delta(h_{i,1})$ and
using~Lemma~\ref{lemma f i1} below, we recover~(\ref{f_i1}).

\begin{Lem}\label{lemma f i1}
For $k<i$ and $l>i+1$, the following equalities hold:

\noindent
(a) $[E^{(0)}_{i+1,l},f_{i,0}]=0$.

\noindent
(b) $[\wt{E}^{(0)}_{ki},f_{i,0}]=0$.

\noindent
(c) $[E^{(0)}_{il},f_{i,0}]=\vv^{-1}E^{(0)}_{i+1,l}\psi^+_{i,0}$.

\noindent
(d) $[\wt{E}^{(0)}_{k,i+1},f_{i,0}]=-\wt{E}^{(0)}_{ki}\psi^+_{i,0}$.

\noindent
(e) $[[E^{(0)}_{i,i+1},E^{(0)}_{i+1,l}]_{\vv^3},f_{i,0}]=
     \frac{\vv^{-1}-\vv^3}{\vv-\vv^{-1}}E^{(0)}_{i+1,l}\psi^+_{i,0}-\frac{\vv-\vv^3}{\vv-\vv^{-1}}E^{(0)}_{i+1,l}\psi^-_{i,0}$.

\noindent
(f) $[[E^{(0)}_{i,i+1},\wt{E}^{(0)}_{ki}]_{\vv^3},f_{i,0}]=
     \frac{\vv^{-1}-\vv^3}{\vv-\vv^{-1}}\wt{E}^{(0)}_{ki}\psi^+_{i,0}-\frac{\vv-\vv^3}{\vv-\vv^{-1}}\wt{E}^{(0)}_{ki}\psi^-_{i,0}$.

\noindent
(g) $[F^{(1)}_{l,i+1},f_{i,0}]_\vv=-F^{(1)}_{li}$.

\noindent
(h) $[F^{(1)}_{ik},f_{i,0}]_\vv=F^{(1)}_{i+1,k}$.

\noindent
(i) $[F^{(1)}_{li},f_{i,0}]_{\vv^{-1}}=0$.

\noindent
(j) $[F^{(1)}_{i+1,k},f_{i,0}]_{\vv^{-1}}=0$.

\noindent
(k) $[F^{(1)}_{lk},f_{i,0}]=0$.
\end{Lem}

This lemma is proved completely analogously to Lemma~\ref{lemma e i-1}.
The details are left to the interested reader.

\medskip
This completes our proof of Theorem~\ref{Coproduct on quantum loop sln}.


\section{Proof of Theorem~\ref{coproduct on antidominant sln} and homomorphisms $\jmath^\pm_{\mu_1,\mu_2}$}
 \label{Proof of coproduct antidominant sl_n}

Our proof of Theorem~\ref{coproduct on antidominant sln} proceeds in three steps.
First, we introduce subalgebras $\CU^{\ssc,\pm}_{0,\mu_1,\mu_2}$ of
$\CU^\ssc_{0,\mu_1+\mu_2}$ and construct homomorphisms $\jmath^\pm_{\mu_1,\mu_2}$
which we referred to in Remark~\ref{Maps jmath for sln}.
Then, we prove Theorem~\ref{coproduct on antidominant sln}, reducing some of the
verifications to the case of $U_\vv(L\ssl_n)$ via the aforementioned
$\jmath^\pm_{\mu_1,\mu_2}$. Finally, we verify the commutativity of the diagram
from Remark~\ref{Maps jmath for sln}.

Throughout this section, we assume $\mu_1,\mu_2\in \Lambda^-$.


\subsection{Homomorphisms $\jmath^\pm_{\mu_1,\mu_2}$}
\

First, we introduce subalgebras $\CU^{\ssc,\pm}_{0,\mu_1,\mu_2}$ of $\CU^\ssc_{0,\mu_1+\mu_2}$.
To this end, recall the explicit identification of the Drinfeld-Jimbo and the new Drinfeld
realizations of $U_\vv(L\ssl_n)$ from Theorem~\ref{identification of DJ and Dr}:
\begin{equation*}
  E_i\mapsto e_{i,0},\ F_i\mapsto f_{i,0},\ K^{\pm 1}_i\mapsto (\psi^+_{i,0})^{\pm 1}=\psi^\pm_{i,0}=(\psi^-_{i,0})^{\mp 1}\ \mathrm{for}\ 1\leq i\leq n-1,
\end{equation*}
\begin{equation*}
  (K_{i_0})^{\pm 1}\mapsto (\psi^+_{1,0}\cdots \psi^+_{n-1,0})^{\mp 1},
\end{equation*}
\begin{equation*}
  E_{i_0}\mapsto (-\vv)^{-n}\cdot (\psi^+_{1,0}\cdots\psi^+_{n-1,0})^{-1}\cdot [\cdots[f_{1,1},f_{2,0}]_\vv,\cdots,f_{n-1,0}]_\vv,
\end{equation*}
\begin{equation*}
  F_{i_0}\mapsto (-\vv)^n\cdot [e_{n-1,0},\cdots,[e_{2,0},e_{1,-1}]_{\vv^{-1}}\cdots]_{\vv^{-1}}\cdot \psi^+_{1,0}\cdots\psi^+_{n-1,0}.
\end{equation*}
Hence, the images $U^+_\vv$ and $U^-_\vv$ of the Drinfeld-Jimbo Borel subalgebras
are the subalgebras of $U_\vv(L\ssl_n)$ generated by
$\{e_{i,0},(\psi^+_{i,0})^{\pm 1}, F^{(1)}_{n1}\}_{i=1}^{n-1}$
and $\{f_{i,0},(\psi^+_{i,0})^{\pm 1}, E^{(-1)}_{1n}\}_{i=1}^{n-1}$, respectively.

Likewise, let $\CU^{\ssc,+}_{0,\mu_1,\mu_2}$ and $\CU^{\ssc,-}_{0,\mu_1,\mu_2}$
be the $\BC(\vv)$-subalgebras of $\CU^\ssc_{0,\mu_1+\mu_2}$ generated by
the elements $\{e_{i,0},(\psi^+_{i,0})^{\pm 1}, F^{(1)}_{n1}\}_{i=1}^{n-1}$ and
$\{f_{i,b_{1,i}},(\psi^-_{i,b_{1,i}+b_{2,i}})^{\pm 1}, \hat{E}^{(-1)}_{1n}\}_{i=1}^{n-1}$,
respectively, where as before
  $b_{1,i}=\alphavee_i(\mu_1), b_{2,i}=\alphavee_i(\mu_2), b_i=b_{1,i}+b_{2,i}$.
Here, the elements $\{\hat{E}^{(-1)}_{ji}\}_{j<i}$ are defined via
$\hat{E}^{(-1)}_{ji}:=[e_{i-1,b_{2,i-1}},[e_{i-2,b_{2,i-2}},\cdots,[e_{j+1,b_{2,j+1}},
e_{j,b_{2,j}-1}]_{\vv^{-1}}\cdots]_{\vv^{-1}}]_{\vv^{-1}}.$

\begin{Prop}\label{jmath homomorphisms}
(a) There is a unique $\BC(\vv)$-algebra homomorphism
$\jmath^+_{\mu_1,\mu_2}\colon  U^+_\vv\to U^{\ssc,+}_{0,\mu_1,\mu_2}$, such that
   $e_{i,0}\mapsto e_{i,0}, (\psi^+_{i,0})^{\pm 1}\mapsto (\psi^+_{i,0})^{\pm 1}, F^{(1)}_{n1}\mapsto F^{(1)}_{n1}$.

\noindent
(b) There is a unique $\BC(\vv)$-algebra homomorphism
$\jmath^-_{\mu_1,\mu_2}\colon  U^-_\vv\to U^{\ssc,-}_{0,\mu_1,\mu_2}$, such that
   $f_{i,0}\mapsto f_{i,b_{1,i}}, (\psi^-_{i,0})^{\pm 1}\mapsto (\psi^-_{i,b_i})^{\pm 1}, E^{(-1)}_{1n}\mapsto \hat{E}^{(-1)}_{1n}$.
\end{Prop}

\begin{proof}
(a) Converting the defining relations of the positive Drinfeld-Jimbo Borel subalgebra
into the new Drinfeld realization, we see that $U^+_\vv$ is generated by
$\{e_{i,0},(\psi^+_{i,0})^{\pm 1}, F^{(1)}_{n1}\}_{i=1}^{n-1}$
with the following defining relations:
\begin{equation}\label{Borel 1}
  (\psi^+_{i,0})^{\pm 1}\cdot (\psi^+_{i,0})^{\mp 1}=1,\
  \psi^+_{i,0}\psi^+_{j,0}=\psi^+_{j,0}\psi^+_{i,0},
\end{equation}
\begin{equation}\label{Borel 2}
  \psi^+_{i,0}e_{j,0}=\vv^{c_{ij}}e_{j,0}\psi^+_{i,0},\
  \psi^+_{i,0}F^{(1)}_{n1}=\vv^{-\delta_{i1}-\delta_{i,n-1}}F^{(1)}_{n1}\psi^+_{i,0},
\end{equation}
\begin{equation}\label{Borel 3}
  [e_{i,0},[e_{i,0},e_{i\pm 1,0}]_\vv]_{\vv^{-1}}=0,\
  [e_{i,0},e_{j,0}]=0\ \mathrm{if}\ c_{ij}=0,
\end{equation}
\begin{equation}\label{Borel 4}
  [e_{i,0},F^{(1)}_{n1}]=0\ \mathrm{for}\ 1<i<n-1,
\end{equation}
\begin{equation}\label{Borel 5}
  [e_{1,0},[e_{1,0},F^{(1)}_{n1}]]_{\vv^{-2}}=0,\
  [e_{n-1,0},[e_{n-1,0},F^{(1)}_{n1}]]_{\vv^{-2}}=0,
\end{equation}
\begin{equation}\label{Borel 6}
  [F^{(1)}_{n1},[F^{(1)}_{n1},e_{1,0}]]_{\vv^2}=0,\
  [F^{(1)}_{n1},[F^{(1)}_{n1},e_{n-1,0}]]_{\vv^2}=0.
\end{equation}

Thus, it suffices to check that these relations are preserved under the specified assignment
  $e_{i,0}\mapsto e_{i,0}, (\psi^+_{i,0})^{\pm 1}\mapsto (\psi^+_{i,0})^{\pm 1}, F^{(1)}_{n1}\mapsto F^{(1)}_{n1}$.
The validity of~(\ref{Borel 1}--\ref{Borel 4}) is obvious.

To verify the first equality of~(\ref{Borel 5}), we note that
$[\psi^+_{1,1},f_{2,0}]_\vv=(\vv^2-1)f_{2,1}\psi^+_{1,0}$,
due to~(\ref{U5}). Combining this with~(\ref{U6}), we get
\begin{equation*}
    [e_{1,0},F^{(1)}_{n1}]=(\vv-\vv^{-1})^{-1}\cdot [\cdots[\psi^+_{1,1},f_{2,0}]_\vv,\cdots,f_{n-1,0}]_\vv=
    \vv F^{(1)}_{n2}\psi^+_{1,0}.
\end{equation*}
Hence,
  $[e_{1,0},[e_{1,0},F^{(1)}_{n1}]]_{\vv^{-2}}=\vv [e_{1,0},F^{(1)}_{n2}\psi^+_{1,0}]_{\vv^{-2}}=\vv [e_{1,0},F^{(1)}_{n2}]\psi^+_{1,0}=0$,
due to~(\ref{U6}).

The verification of the second equality of~(\ref{Borel 5}) is similar and is based on
\begin{equation*}
    [e_{n-1,0},F^{(1)}_{n1}]=
    \frac{[[\cdots[f_{1,1},f_{2,0}]_\vv,\cdots,f_{n-2,0}]_\vv,\psi^+_{n-1,0}-\delta_{b_{n-1},0}\psi^-_{n-1,0}]_\vv}{\vv-\vv^{-1}}=
    -\vv F^{(1)}_{n-1,1}\psi^+_{n-1,0}.
\end{equation*}

Due to the above equality $[e_{1,0},F^{(1)}_{n1}]=\vv F^{(1)}_{n2}\psi^+_{1,0}$ and~(\ref{U4}),
the verification of the first equality of~(\ref{Borel 6}) boils down to the proof
of $[F^{(1)}_{n1},F^{(1)}_{n2}]_\vv=0$. This is an equality in $\CU^{\ssc,<}_{0,\mu_1+\mu_2}$.
However, $\CU^{\ssc,<}_{0,\mu_1+\mu_2}\simeq U^<_\vv(L\ssl_n)$, due to
Proposition~\ref{Triangular decomposition}(b). Hence, it suffices to check this equality in $U_\vv(L\ssl_n)$.
The latter follows immediately from the validity of~(\ref{Borel 6}) for $U^+_\vv$.

Due to $[e_{n-1,0},F^{(1)}_{n1}]=-\vv F^{(1)}_{n-1,1}\psi^+_{n-1,0}$ from above and~(\ref{U4}),
the verification of the second equality of~(\ref{Borel 6}) boils down to the proof of
$[F^{(1)}_{n1},F^{(1)}_{n-1,1}]_\vv=0$. Analogously to the previous verification,
the latter follows from the same equality in $U^+_\vv$.

(b) The proof of part (b) is completely analogous and is left to the interested reader.
\end{proof}

This completes our construction of the homomorphisms
$\jmath^\pm_{\mu_1,\mu_2}\colon  U^\pm_\vv\to \CU^{\ssc,\pm}_{0,\mu_1,\mu_2}$,
which we referred to in Remark~\ref{Maps jmath for sln}.
The following results are needed for the next subsection.

\begin{Lem}\label{Belong}
(a) For any $1\leq j<i\leq n$, we have
  $E^{(0)}_{ji},\wt{E}^{(0)}_{ji}, F^{(1)}_{ij}\in \CU^{\ssc,+}_{0,\mu_1,\mu_2}$.

\noindent
(b) For any $1\leq j<i\leq n$, define
$\hat{F}^{\pm,(0)}_{ij}:=[\cdots[f_{j,b_{1,j}},f_{j+1,b_{1,j+1}}]_{\vv^{\pm 1}},\cdots,f_{i-1,b_{1,i-1}}]_{\vv^{\pm 1}}$.
We have $\hat{F}^{\pm,(0)}_{ij},\hat{E}^{(-1)}_{ji}\in \CU^{\ssc,-}_{0,\mu_1,\mu_2}$.
\end{Lem}

\begin{proof}
(a) Since $E^{(0)}_{ji},\wt{E}^{(0)}_{ji}$ are expressed via
$\vv^{\pm 1}$-commutators of $e_{k,0}\in \CU^{\ssc,+}_{0,\mu_1,\mu_2}$, we obviously
get the first two inclusions. The last inclusion is clear for $(i,j)=(n,1)$.
Applying iteratively
  $[e_{k,0},F^{(1)}_{k+1,1}]=-\vv F^{(1)}_{k1}\psi^+_{k,0},\
   [e_{l,0},F^{(1)}_{il}]=\vv F^{(1)}_{i,l+1}\psi^+_{l,0}$,
we get $F^{(1)}_{ij}\in\CU^{\ssc,+}_{0,\mu_1,\mu_2}$ for any $j<i$.

(b) The inclusions $\hat{F}^{\pm,(0)}_{ij}\in \CU^{\ssc,-}_{0,\mu_1,\mu_2}$ are obvious.
It remains to prove $\hat{E}^{(-1)}_{ji}\in \CU^{\ssc,-}_{0,\mu_1,\mu_2}$.
This is clear for $(j,i)=(1,n)$. To deduce the general case,
it remains to apply the equalities
  $[f_{i-1,b_{1,i-1}},\hat{E}^{(-1)}_{1i}]=\hat{E}^{(-1)}_{1,i-1}\psi^-_{i-1,b_{i-1}},
   [f_{l,b_{1,l}},\hat{E}^{(-1)}_{li}]=-\hat{E}^{(-1)}_{l+1,i}\psi^-_{l,b_l}$.
\end{proof}

The proof of the following result is straightforward.

\begin{Lem}\label{jmath compatibility}
For any $1\leq j<i\leq n$, we have:
\begin{equation*}
   \jmath^+_{\mu_1,\mu_2}\colon  E^{(0)}_{ji}\mapsto E^{(0)}_{ji}, \wt{E}^{(0)}_{ji}\mapsto \wt{E}^{(0)}_{ji}, F^{(1)}_{ij}\mapsto F^{(1)}_{ij},
     f_{i,1}\mapsto f_{i,1}, h_{i,1}\mapsto h_{i,1},
\end{equation*}
\begin{equation*}
   \jmath^-_{\mu_1,\mu_2}\colon  F^{(0)}_{ij}\mapsto \hat{F}^{+,(0)}_{ij}, \wt{F}^{(0)}_{ij}\mapsto \hat{F}^{-,(0)}_{ij}, E^{(-1)}_{ji}\mapsto \hat{E}^{(-1)}_{ji},
     e_{i,-1}\mapsto e_{i,b_{2,i}-1}, h_{i,-1}\mapsto h_{i,-1}.
\end{equation*}
\end{Lem}


\subsection{Proof of Theorem~\ref{coproduct on antidominant sln}}
\

Due to Theorem~\ref{Levendorskii-type}, it suffices to check that
the assignment $\Delta$ of Theorem~\ref{coproduct on antidominant sln}
preserves defining relations~(\ref{u1}--\ref{u9}). To simplify our exposition, we will assume that
$\mu_1,\mu_2$ are strictly antidominant: $b_{1,i}, b_{2,i}<0$ for any $1\leq i<n$.
This verification is similar to the $n=2$ case (carried out in Appendix~\ref{Proof Coproduct}) and
we only indicate the key technical details, see~Lemmas~\ref{long u2.1}--\ref{long u5.3}
(their proofs are similar to that of Lemma~\ref{lemma e i-1} and therefore omitted).
For $1\leq a\leq b<n$, we define $\alphavee_{[a,b]}:=\alphavee_a+\alphavee_{a+1}+\ldots+\alphavee_b$.


\subsubsection{Compatibility with~(\ref{u1})}
\

\medskip
\noindent
$\circ$
The equalities
  $\Delta((\psi^+_{i,0})^{\pm 1})\Delta((\psi^+_{i,0})^{\mp 1})=1$
and
  $\Delta((\psi^-_{i,b_i})^{\pm 1})\Delta((\psi^-_{i,b_i})^{\mp 1})=1$
follow immediately from relation~(\ref{u1}) for both $\CU^\ssc_{0,\mu_1},\CU^\ssc_{0,\mu_2}$.

\medskip
\noindent
$\circ$
The commutativity of $\{\Delta(\psi^+_{i,0}),\Delta(\psi^-_{i,b_i})\}_{i=1}^{n-1}$
between themselves and with $\{\Delta(h_{j,\pm 1})\}_{j=1}^{n-1}$ is due to
relations~(\ref{u1},~\ref{u4},~\ref{u5}) for both $\CU^\ssc_{0,\mu_1},\CU^\ssc_{0,\mu_2}$.

\medskip
\noindent
$\circ$
Finally, we verify $[\Delta(h_{i,r}),\Delta(h_{j,s})]=0$ for $r,s\in \{\pm 1\}$.
To this end, recall the homomorphism
  $\iota_{0,0,\mu_1}\otimes \iota_{0,\mu_2,0}\colon
   \CU^\ssc_{0,0}\otimes \CU^\ssc_{0,0}\to \CU^\ssc_{0,\mu_1}\otimes \CU^\ssc_{0,\mu_2}$.
The key observation is that
  $\iota_{0,0,\mu_1}\otimes \iota_{0,\mu_2,0}(\Delta(h_{i,r}))=
   \Delta(h_{i,r})+\frac{\alphavee_i(\mu_1+\mu_2)}{\vv^r-\vv^{-r}}$
for any $i\in I, r\in \{\pm 1\}$ (cf.~proof of~Corollary~\ref{uniform coproduct h_1}),
where by abuse of notation we use $\Delta(h_{i,r})$ to denote elements of both
$\CU^\ssc_{0,0}\otimes \CU^\ssc_{0,0}$ and $\CU^\ssc_{0,\mu_1}\otimes \CU^\ssc_{0,\mu_2}$.
Hence, it suffices to prove
$[\Delta(h_{i,r}),\Delta(h_{j,s})]=0$ in $\CU^\ssc_{0,0}\otimes \CU^\ssc_{0,0}$.
The latter follows immediately from the corresponding result for $U_\vv(L\ssl_n)$,
in which case the assignment $\Delta$ of Theorem~\ref{coproduct on antidominant sln}
coincides with the Drinfeld-Jimbo coproduct, due to Theorem~\ref{Coproduct on quantum loop sln}.


\subsubsection{Compatibility with~(\ref{u2})}
\

We need to prove
  $[\Delta(e_{i,r+1}),\Delta(e_{j,s})]_{\vv^{c_{ij}}}+[\Delta(e_{j,s+1}),\Delta(e_{i,r})]_{\vv^{c_{ij}}}=0$
for $b_{2,i}-1\leq r\leq -1, b_{2,j}-1\leq s\leq -1$.

\medskip
\noindent
$\circ$ \emph{Case $b_{2,i}-1<r\leq -1, b_{2,j}-1<s\leq -1$.}
In this case, the above sum equals
  $1\otimes ([e_{i,r+1},e_{j,s}]_{\vv^{c_{ij}}}+[e_{j,s+1},e_{i,r}]_{\vv^{c_{ij}}})=0$,
due to relations~(\ref{u2}) and~(\ref{u4}) for $\CU^\ssc_{0,\mu_2}$.

\medskip
\noindent
$\circ$ \emph{Case $r=b_{2,i}-1, b_{2,j}-1<s<-1$.}
Note that $[e_{j,s+1},f_{a,0}]=0$ for any $1\leq a<n$, due to~(\ref{u6}) for $\CU^\ssc_{0,\mu_2}$.
As a result, we have $[e_{j,s+1},F^{(0)}_{ba}]=[e_{j,s+1},\wt{F}^{(0)}_{ba}]=0$ for any $1\leq a<b\leq n$.
Combining this with~(\ref{u2}) and~(\ref{u4}) for $\CU^\ssc_{0,\mu_2}$, we get
  $[\Delta(e_{i,b_{2,i}}),\Delta(e_{j,s})]_{\vv^{c_{ij}}}+[\Delta(e_{j,s+1}),\Delta(e_{i,b_{2,i}-1})]_{\vv^{c_{ij}}}=
   1\otimes ([e_{i,b_{2,i}},e_{j,s}]_{\vv^{c_{ij}}}+[e_{j,s+1},e_{i,b_{2,i}-1}]_{\vv^{c_{ij}}})=0$
as above.

\medskip
\noindent
$\circ$ \emph{Case $r=b_{2,i}-1, s=b_{2,j}-1$.}
Due to relation~(\ref{u4}) for both $\CU^\ssc_{0,\mu_1},\CU^\ssc_{0,\mu_2}$, we get
\begin{equation*}
\begin{split}
   & [\Delta(e_{j,b_{2,j}}),\Delta(e_{i,b_{2,i}-1})]_{\vv^{c_{ij}}}=
   1\otimes [e_{j,b_{2,j}},e_{i,b_{2,i}-1}]_{\vv^{c_{ij}}}-\\
   & (\vv-\vv^{-1})\sum_{l>i+1}E^{(-1)}_{il}\otimes [e_{j,b_{2,j}},F^{(0)}_{l,i+1}]\psi^-_{i,b_{2,i}}+
   (\vv-\vv^{-1})\sum_{k<i}\vv^{i-k-1}E^{(-1)}_{k,i+1}\otimes [e_{j,b_{2,j}},\wt{F}^{(0)}_{ik}]\psi^-_{i,b_{2,i}}-\\
   & (\vv-\vv^{-1})^2\sum_{l>i+1}^{k<i}\vv^{i-k-1}E^{(-1)}_{kl}\otimes [e_{j,b_{2,j}},\wt{F}^{(0)}_{ik}F^{(0)}_{l,i+1}]\psi^-_{i,b_{2,i}}.
\end{split}
\end{equation*}
Using this formula and Lemma~\ref{long u2.1} below, it is straightforward to check that again we obtain
  $[\Delta(e_{i,b_{2,i}}),\Delta(e_{j,b_{2,j}-1})]_{\vv^{c_{ij}}}+[\Delta(e_{j,b_{2,j}}),\Delta(e_{i,b_{2,i}-1})]_{\vv^{c_{ij}}}=
   1\otimes ([e_{i,b_{2,i}},e_{j,b_{2,j}-1}]_{\vv^{c_{ij}}}+[e_{j,b_{2,j}},e_{i,b_{2,i}-1}]_{\vv^{c_{ij}}})=0$.

\begin{Lem}\label{long u2.1}
For any $1\leq k<i,\ i+1<l\leq n,\ 1\leq j<n$, the following holds in $\CU^\ssc_{0,\mu_2}$:

\noindent
(a) $[e_{j,b_{2,j}},F^{(0)}_{l,i+1}]=\delta_{j,i+1}F^{(0)}_{l,i+2}\psi^-_{j,b_{2,j}}$,
where we set $F^{(0)}_{i+2,i+2}:=\frac{-1}{\vv-\vv^{-1}}$.

\noindent
(b) $[e_{j,b_{2,j}},\wt{F}^{(0)}_{ik}]=-\vv^{-1}\delta_{j,i-1}\wt{F}^{(0)}_{i-1,k}\psi^-_{j,b_{2,j}}$,
where we set $\wt{F}^{(0)}_{i-1,i-1}:=\frac{\vv}{\vv-\vv^{-1}}$.
\end{Lem}

\noindent
$\circ$ \emph{Case $r=b_{2,i}-1, s=-1$.}
Clearly,
  $[\Delta(e_{i,b_{2,i}}),\Delta(e_{j,-1})]_{\vv^{c_{ij}}}=1\otimes [e_{i,b_{2,i}},e_{j,-1}]_{\vv^{c_{ij}}}$
and
  $[\Delta(e_{j,0}),\Delta(e_{i,b_{2,i}-1})]_{\vv^{c_{ij}}}=[1\otimes e_{j,0}+e_{j,0}\otimes \psi^+_{j,0},\Delta(e_{i,b_{2,i}-1})]_{\vv^{c_{ij}}}$.
We claim that as in the previous cases, one gets
  $[\Delta(e_{i,b_{2,i}}),\Delta(e_{j,-1})]_{\vv^{c_{ij}}}+[\Delta(e_{j,0}),\Delta(e_{i,b_{2,i}-1})]_{\vv^{c_{ij}}}=
   1\otimes ([e_{i,b_{2,i}},e_{j,-1}]_{\vv^{c_{ij}}}+[e_{j,0},e_{i,b_{2,i}-1}]_{\vv^{c_{ij}}})=0$.
To this end, we note that the computations of
  $[1\otimes e_{j,0},\Delta(e_{i,b_{2,i}-1})]_{\vv^{c_{ij}}}$
and
  $[e_{j,0}\otimes \psi^+_{j,0},\Delta(e_{i,b_{2,i}-1})]_{\vv^{c_{ij}}}$
are straightforward and are crucially based on Lemmas~\ref{long u2.2} and~\ref{long u2.3} below, respectively.

\begin{Lem}\label{long u2.2}
For any $1\leq k<i,\ i+1<l\leq n,\ 1\leq j<n$, the following holds in $\CU^\ssc_{0,\mu_2}$:

\noindent
(a) $[e_{j,0},F^{(0)}_{l,i+1}]=-\vv\delta_{j,l-1}F^{(0)}_{j,i+1}\psi^+_{j,0}$,
where we set $F^{(0)}_{i+1,i+1}:=\frac{-1}{\vv(\vv-\vv^{-1})}$.

\noindent
(b) $[e_{j,0},\wt{F}^{(0)}_{ik}]=\delta_{jk}\wt{F}^{(0)}_{i,j+1}\psi^+_{j,0}$,
where we set $\wt{F}^{(0)}_{ii}:=\frac{1}{\vv-\vv^{-1}}$.
\end{Lem}

\begin{Lem}\label{long u2.3}
For any $1\leq k<l-1<n,\ 1\leq j<n$, the following holds in $\CU^\ssc_{0,\mu_1}$:



\noindent
$[e_{j,0},E^{(-1)}_{kl}]_{\vv^{(\alphavee_j,\alphavee_{[k,l-1]})}}=
\delta_{jl}E^{(-1)}_{k,l+1}-\delta_{j,k-1}E^{(-1)}_{k-1,l}.$
\end{Lem}


\subsubsection{Compatibility with~(\ref{u3})}
\

We need to prove
  $[\Delta(f_{i,r+1}),\Delta(f_{j,s})]_{\vv^{-c_{ij}}}+[\Delta(f_{j,s+1}),\Delta(f_{i,r})]_{\vv^{-c_{ij}}}=0$
for $b_{1,i}\leq r\leq 0, b_{1,j}\leq s\leq 0$.

\medskip
\noindent
$\circ$ \emph{Case $b_{1,i}\leq r<0, b_{1,j}\leq s<0$.}
In this case, the above sum equals
  $([f_{i,r+1},f_{j,s}]_{\vv^{-c_{ij}}}+[f_{j,s+1},f_{i,r}]_{\vv^{-c_{ij}}})\otimes1=0$,
due to relations~(\ref{u3}) and~(\ref{u5}) for $\CU^\ssc_{0,\mu_1}$.

\medskip
\noindent
$\circ$ \emph{Case $r=0, b_{1,j}<s<0$.}
Note that $[f_{j,s},e_{a,0}]=0$ for any $1\leq a<n$, due to~(\ref{u6}) for $\CU^\ssc_{0,\mu_1}$.
As a result, we have $[f_{j,s},E^{(0)}_{ab}]=[f_{j,s},\wt{E}^{(0)}_{ab}]=0$ for any $1\leq a<b\leq n$.
Combining this with~(\ref{u3}) and~(\ref{u5}) for $\CU^\ssc_{0,\mu_1}$, we get
  $[\Delta(f_{i,1}),\Delta(f_{j,s})]_{\vv^{-c_{ij}}}+[\Delta(f_{j,s+1}),\Delta(f_{i,0})]_{\vv^{-c_{ij}}}=
   ([f_{i,1},f_{j,s}]_{\vv^{-c_{ij}}}+[f_{j,s+1},f_{i,0}]_{\vv^{-c_{ij}}})\otimes 1=0$
as above.

\medskip
\noindent
$\circ$ \emph{Case $r=0, s=0$.}
Due to relation~(\ref{u5}) for both $\CU^\ssc_{0,\mu_1},\CU^\ssc_{0,\mu_2}$, we get
\begin{equation*}
\begin{split}
   & [\Delta(f_{i,1}),\Delta(f_{j,0})]_{\vv^{-c_{ij}}}=
   [f_{i,1},f_{j,0}]_{\vv^{-c_{ij}}}\otimes 1 + (\vv-\vv^{-1})\vv^{-c_{ij}-1}\sum_{l>i+1}[E^{(0)}_{i+1,l},f_{j,0}]\psi^+_{i,0}\otimes F^{(1)}_{li}-\\
   & (\vv-\vv^{-1})\vv^{-c_{ij}}\sum_{k<i}\vv^{k-i}[\wt{E}^{(0)}_{ki},f_{j,0}]\psi^+_{i,0}\otimes F^{(1)}_{i+1,k}-\\
   & (\vv-\vv^{-1})^2\vv^{-c_{ij}}\sum_{l>i+1}^{k<i} \vv^{k-i-1}[E^{(0)}_{i+1,l}\wt{E}^{(0)}_{ki},f_{j,0}]\psi^+_{i,0}\otimes F^{(1)}_{lk}.
\end{split}
\end{equation*}
Using this formula and Lemma~\ref{long u3.1} below, it is straightforward to check that we obtain
  $[\Delta(f_{i,1}),\Delta(f_{j,0})]_{\vv^{-c_{ij}}}+[\Delta(f_{j,1}),\Delta(f_{i,0})]_{\vv^{-c_{ij}}}=
   ([f_{i,1},f_{j,0}]_{\vv^{-c_{ij}}}+[f_{j,1},f_{i,0}]_{\vv^{-c_{ij}}})\otimes 1=0$.

\begin{Lem}\label{long u3.1}
For any $1\leq k<i,\ i+1<l\leq n,\ 1\leq j<n$, the following holds in $\CU^\ssc_{0,\mu_1}$:

\noindent
(a) $[E^{(0)}_{i+1,l},f_{j,0}]=\vv^{-1}\delta_{j,i+1}E^{(0)}_{i+2,l}\psi^+_{j,0}$,
where we set $E^{(0)}_{i+2,i+2}:=\frac{\vv}{\vv-\vv^{-1}}$.

\noindent
(b) $[\wt{E}^{(0)}_{ki},f_{j,0}]=-\delta_{j,i-1}\wt{E}^{(0)}_{k,i-1}\psi^+_{j,0}$,
where we set $\wt{E}^{(0)}_{i-1,i-1}:=\frac{-1}{\vv-\vv^{-1}}$.
\end{Lem}

\noindent
$\circ$ \emph{Case $r=0, s=b_{1,j}$.}
Clearly,
  $[\Delta(f_{j,b_{1,j}+1}),\Delta(f_{i,0})]_{\vv^{-c_{ij}}}=[f_{j,b_{1,j}+1},f_{i,0}]_{\vv^{-c_{ij}}}\otimes 1$
and
  $[\Delta(f_{i,1}),\Delta(f_{j,b_{1,j}})]_{\vv^{-c_{ij}}}=[\Delta(f_{i,1}),f_{j,b_{1,j}}\otimes 1+\psi^-_{j,b_{1,j}}\otimes f_{j,0}]_{\vv^{-c_{ij}}}$.
We claim that as in the previous cases, one gets
  $[\Delta(f_{i,1}),\Delta(f_{j,b_{1,j}})]_{\vv^{-c_{ij}}}+[\Delta(f_{j,b_{1,j}+1}),\Delta(f_{i,0})]_{\vv^{-c_{ij}}}=
   ([f_{i,1},f_{j,b_{1,j}}]_{\vv^{-c_{ij}}}+[f_{j,b_{1,j}+1},f_{i,0}]_{\vv^{-c_{ij}}})\otimes 1=0$.
To this end, we note that the computations of
  $[\Delta(f_{i,1}),f_{j,b_{1,j}}\otimes1]_{\vv^{-c_{ij}}}$
and $[\Delta(f_{i,1}),\psi^-_{j,b_{1,j}}\otimes f_{j,0}]_{\vv^{-c_{ij}}}$
are straightforward and are crucially based on~Lemmas~\ref{long u3.2}
and~\ref{long u3.3} below, respectively.

\begin{Lem}\label{long u3.2}
For any $1\leq k<i,\ i+1<l\leq n,\ 1\leq j<n$, the following holds in $\CU^\ssc_{0,\mu_1}$:

\noindent
(a) $[E^{(0)}_{i+1,l},f_{j,b_{1,j}}]=-\delta_{j,l-1}E^{(0)}_{i+1,j}\psi^-_{j,b_{1,j}}$,
where we set $E^{(0)}_{i+1,i+1}:=\frac{1}{\vv-\vv^{-1}}$.

\noindent
(b) $[\wt{E}^{(0)}_{ki},f_{j,b_{1,j}}]=\vv\delta_{jk}\wt{E}^{(0)}_{j+1,i}\psi^-_{j,b_{1,j}}$,
where we set $\wt{E}^{(0)}_{ii}:=\frac{-1}{\vv(\vv-\vv^{-1})}$.
\end{Lem}

\begin{Lem}\label{long u3.3}
For any $1\leq k<l-1<n,\ 1\leq j<n$, the following holds in $\CU^\ssc_{0,\mu_2}$:



\noindent
$[F^{(1)}_{lk},f_{j,0}]_{\vv^{-(\alphavee_j,\alphavee_{[k,l-1]})}}=
\delta_{jl}F^{(1)}_{l+1,k}-\delta_{j,k-1}F^{(1)}_{l,k-1}.$
\end{Lem}


\subsubsection{Compatibility with~(\ref{u4})}
\

Due to relations~(\ref{u1},~\ref{u4},~\ref{u5}) for both
$\CU^\ssc_{0,\mu_1}, \CU^\ssc_{0,\mu_2}$,
we immediately obtain the equalities
  $\Delta(\psi^+_{i,0})\Delta(e_{j,r})=\vv^{c_{ij}}\Delta(e_{j,r})\Delta(\psi^+_{i,0}),\
   \Delta(\psi^-_{i,b_i})\Delta(e_{j,r})=\vv^{-c_{ij}}\Delta(e_{j,r})\Delta(\psi^-_{i,b_i})$ for $b_{2,j-1}\leq r\leq0$.

\medskip
Let us now verify
  $[\Delta(h_{i,1}),\Delta(e_{j,r})]=[c_{ij}]_\vv\cdot \Delta(e_{j,r+1})$
for $b_{2,j}-1\leq r\leq -1$.

\medskip
\noindent
$\circ$ \emph{Case $b_{2,j}\leq r<-1$.}
 The verification in this case follows immediately from relation~(\ref{u4}) for
$\CU^\ssc_{0,\mu_2}$ combined with~Lemma~\ref{long u4.1} below.

\begin{Lem}\label{long u4.1}
  For any $1\leq a<b\leq n,\ b_{2,j}\leq r<-1$, we have $[F^{(1)}_{ba},e_{j,r}]=0$ in $\CU^\ssc_{0,\mu_2}$.
\end{Lem}

\noindent
$\circ$ \emph{Case $r=-1$.}
Due to relation~(\ref{u4}) for $\CU^\ssc_{0,\mu_2}$, we  get
\begin{equation*}
\begin{split}
   & [\Delta(h_{i,1}),\Delta(e_{j,-1})]=
   [c_{ij}]_\vv\cdot 1\otimes e_{j,0}-(\vv^2-\vv^{-2})E^{(0)}_{i,i+1}\otimes [F^{(1)}_{i+1,i},e_{j,-1}]+\\
   & (\vv-\vv^{-1})\sum_{l>i+1} E^{(0)}_{i+1,l}\otimes [F^{(1)}_{l,i+1},e_{j,-1}]+
   (\vv-\vv^{-1})\sum_{k<i} \vv^{k+1-i}\wt{E}^{(0)}_{ki}\otimes [F^{(1)}_{ik},e_{j,-1}]+\\
   & \vv^{-2}(\vv-\vv^{-1})\sum_{l>i+1}[E^{(0)}_{i,i+1},E^{(0)}_{i+1,l}]_{\vv^3}\otimes [F^{(1)}_{li},e_{j,-1}]-\\
   & (\vv-\vv^{-1})\sum_{k<i}\vv^{k-i-1}[E^{(0)}_{i,i+1},\wt{E}^{(0)}_{ki}]_{\vv^3}\otimes [F^{(1)}_{i+1,k},e_{j,-1}]+\\
   & (\vv-\vv^{-1})^2\sum_{l>i+1}^{k<i} \vv^{k-i}(E^{(0)}_{il}\wt{E}^{(0)}_{ki}-E^{(0)}_{i+1,l}\wt{E}^{(0)}_{k,i+1})\otimes [F^{(1)}_{lk},e_{j,-1}].
\end{split}
\end{equation*}
Using this formula and~Lemma~\ref{long u4.2} below, it is straightforward to
check that we obtain
$[\Delta(h_{i,1}),\Delta(e_{j,-1})]=[c_{ij}]_\vv\cdot
(1\otimes e_{j,0}+e_{j,0}\otimes \psi^+_{j,0})=[c_{ij}]_\vv\cdot\Delta(e_{j,0})$.

\begin{Lem}\label{long u4.2}
  For any $1\leq a<b\leq n$, we have
  $[F^{(1)}_{ba},e_{j,-1}]=\frac{-1}{\vv-\vv^{-1}}\delta_{ja}\delta_{j,b-1}\psi^+_{j,0}$
  in $\CU^\ssc_{0,\mu_2}$.
\end{Lem}

\noindent
$\circ$ \emph{Case $r=b_{2,j}-1$.}
According to the next step, we have
  $\Delta(e_{j,b_{2,j}-1})=\frac{[\Delta(h_{j,-1}),\Delta(e_{j,b_{2,j}})]}{[2]_\vv}$.
Apply the Jacobi identity to get
  $[2]_\vv\cdot [\Delta(h_{i,1}),\Delta(e_{j,b_{2,j}-1})]=
   [\Delta(h_{j,-1}),[\Delta(h_{i,1}),\Delta(e_{j,b_{2,j}})]]-
   [\Delta(e_{j,b_{2,j}}),[\Delta(h_{i,1}),\Delta(h_{j,-1})]]$.
The second summand is zero as $[\Delta(h_{i,1}),\Delta(h_{j,-1})]=0$ by above.
Due to the $r=b_{2,j}$ case considered above, we have
  $[\Delta(h_{i,1}),\Delta(e_{j,b_{2,j}})]=[c_{ij}]_\vv\cdot \Delta(e_{j,b_{2,j}+1})$.
It remains to apply
  $[\Delta(h_{j,-1}),\Delta(e_{j,b_{2,j}+1})]=[2]_\vv\cdot \Delta(e_{j,b_{2,j}})$
as proved below.

\medskip
Let us now verify the equality
  $[\Delta(h_{i,-1}),\Delta(e_{j,r})]=[c_{ij}]_\vv\cdot \Delta(e_{j,r-1})$
for $b_{2,j}\leq r\leq 0$.

\medskip
\noindent
$\circ$ \emph{Case $b_{2,j}<r<0$.}
 The verification in this case follows immediately from relation~(\ref{u4}) for
$\CU^\ssc_{0,\mu_2}$ combined with~Lemma~\ref{long u4.3} below.

\begin{Lem}\label{long u4.3}
For $1\leq a<b\leq n,b_{2,j}<r<0$, we have $[F^{(0)}_{ba},e_{j,r}]=[\wt{F}^{(0)}_{ba},e_{j,r}]=0$
in $\CU^\ssc_{0,\mu_2}$.
\end{Lem}

\medskip
\noindent
$\circ$ \emph{Case $r=b_{2,j}$.}
For $i=j$, the verification of $[\Delta(h_{j,-1}),\Delta(e_{j,b_{2,j}})]=[2]_\vv\cdot \Delta(e_{j,b_{2,j}-1})$
coincides with our proof of formula~(\ref{e_i-1}) from Appendix~\ref{Proof of coproduct on loop sln}.
To prove the claim for $i\ne j$, we can either perform similar long computations or we can rather deduce from the
aforementioned case $i=j$. To achieve the latter, we apply the Jacobi identity to get
  $[2]_\vv\cdot [\Delta(h_{i,-1}),\Delta(e_{j,b_{2,j}})]=
   [\Delta(h_{j,-1}),[\Delta(h_{i,-1}),\Delta(e_{j,b_{2,j}+1})]]-[\Delta(e_{j,b_{2,j}+1}),[\Delta(h_{i,-1}),\Delta(h_{j,-1})]]$.
The second summand is zero as $[\Delta(h_{i,-1}),\Delta(h_{j,-1})]=0$ by above.
Due to the $r=b_{2,j}+1$ case considered above, we have
$[\Delta(h_{i,-1}),\Delta(e_{j,b_{2,j}+1})]=[c_{ij}]_\vv\cdot \Delta(e_{j,b_{2,j}})$.
It remains to apply the aforementioned  equality
$[\Delta(h_{j,-1}),\Delta(e_{j,b_{2,j}})]=[2]_\vv\cdot \Delta(e_{j,b_{2,j}-1})$.

\medskip
\noindent
$\circ$ \emph{Case $r=0$}.
The verification of
  $[\Delta(h_{i,-1}),\Delta(e_{j,0})]=[c_{ij}]_\vv\cdot 1\otimes e_{j,-1}$
is similar to our proof of formula~(\ref{e_i-1}) from Appendix~\ref{Proof of coproduct on loop sln}.
To this end, we note that the computations of
  $[\Delta(h_{i,-1}),1\otimes e_{j,0}]$
and
  $[\Delta(h_{i,-1}),e_{j,0}\otimes \psi^+_{j,0}]$
are straightforward  and are crucially based on the above
Lemmas~\ref{long u2.2} and~\ref{long u2.3}.


\subsubsection{Compatibility with~(\ref{u5})}
\

Due to relations~(\ref{u1},~\ref{u4},~\ref{u5}) for both $\CU^\ssc_{0,\mu_1}, \CU^\ssc_{0,\mu_2}$,
we immediately obtain the equalities
  $\Delta(\psi^+_{i,0})\Delta(f_{j,r})=\vv^{-c_{ij}}\Delta(f_{j,r})\Delta(\psi^+_{i,0}),\
   \Delta(\psi^-_{i,b_i})\Delta(f_{j,r})=\vv^{c_{ij}}\Delta(f_{j,r})\Delta(\psi^-_{i,b_i})$ for $b_{1,j}\leq r\leq1$.

\medskip
Let us now verify
  $[\Delta(h_{i,-1}),\Delta(f_{j,r})]=-[c_{ij}]_\vv\cdot \Delta(f_{j,r-1})$
for $b_{1,j}+1\leq r\leq 1$.

\medskip
\noindent
$\circ$ \emph{Case $b_{1,j}+1<r<1$.}
 The verification in this case follows immediately from relation~(\ref{u5}) for
$\CU^\ssc_{0,\mu_1}$ combined with~Lemma~\ref{long u5.1} below.

\begin{Lem}\label{long u5.1}
For any $1\leq a<b\leq n,\ b_{1,j}+1<r<1$, we have $[E^{(-1)}_{ab},f_{j,r}]=0$
in $\CU^\ssc_{0,\mu_1}$.
\end{Lem}

\noindent
$\circ$ \emph{Case $r=b_{1,j}+1$.}
Due to relation~(\ref{u5}) for $\CU^\ssc_{0,\mu_1}$, we have
\begin{equation*}
\begin{split}
  & [\Delta(h_{i,-1}),\Delta(f_{j,b_{1,j}+1})]=
  -[c_{ij}]_\vv\cdot f_{j,b_{1,j}}\otimes 1+
  (\vv^2-\vv^{-2})[E_{i,i+1}^{(-1)},f_{j,b_{1,j}+1}]\otimes F^{(0)}_{i+1,i}-\\
  & (\vv-\vv^{-1})\sum_{l>i+1}[E^{(-1)}_{i+1,l},f_{j,b_{1,j}+1}]\otimes F^{(0)}_{l,i+1}-
  (\vv-\vv^{-1})\sum_{k<i}\vv^{i-k-1}[E^{(-1)}_{ki},f_{j,b_{1,j}+1}]\otimes \wt{F}^{(0)}_{ik}-\\
  & \vv^{2}(\vv-\vv^{-1})\sum_{l>i+1}[E^{(-1)}_{il},f_{j,b_{1,j}+1}]\otimes [F^{(0)}_{l,i+1},F^{(0)}_{i+1,i}]_{\vv^{-3}}+\\
  & (\vv-\vv^{-1})\sum_{k<i} \vv^{i+1-k} [E^{(-1)}_{k,i+1},f_{j,b_{1,j}+1}]\otimes [\wt{F}^{(0)}_{ik},F^{(0)}_{i+1,i}]_{\vv^{-3}}-\\
  & (\vv-\vv^{-1})^2\sum_{l>i+1}^{k<i}\vv^{i-k}[E^{(-1)}_{kl},f_{j,b_{1,j}+1}]\otimes (\wt{F}^{(0)}_{i+1,k}F^{(0)}_{l,i+1}-\wt{F}^{(0)}_{ik}F^{(0)}_{li}).
\end{split}
\end{equation*}
Using this formula and Lemma~\ref{long u5.2} below, it is straightforward to
check that we obtain
$[\Delta(h_{i,-1}),\Delta(f_{j,b_{1,j}+1})]=-[c_{ij}]_\vv\cdot (f_{j,b_{1,j}}\otimes
1+\psi^-_{j,b_{1,j}}\otimes f_{j,0})=-[c_{ij}]_\vv\cdot\Delta(f_{j,b_{1,j}})$.

\begin{Lem}\label{long u5.2}
For any $1\leq a<b\leq n$, we have
$[E^{(-1)}_{ab},f_{j,b_{1,j}+1}]=\frac{-\delta_{ja}\delta_{j,b-1}}{\vv-\vv^{-1}}\psi^-_{j,b_{1,j}}$
in $\CU^\ssc_{0,\mu_1}$.
\end{Lem}

\noindent
$\circ$ \emph{Case $r=1$.}
According to the next step, we have
  $\Delta(f_{j,1})=-[2]_\vv^{-1}\cdot [\Delta(h_{j,1}),\Delta(f_{j,0})]$.
Apply the Jacobi identity to get
  $[2]_\vv\cdot [\Delta(h_{i,-1}),\Delta(f_{j,1})]=
   [\Delta(h_{j,1}),[\Delta(h_{i,-1}),\Delta(f_{j,0})]]-
   [\Delta(f_{j,0}),[\Delta(h_{i,-1}),\Delta(h_{j,1})]]$.
The second summand is zero as $[\Delta(h_{i,-1}),\Delta(h_{j,1})]=0$ by above.
Due to the $r=0$ case considered above, we have
  $[\Delta(h_{i,-1}),\Delta(f_{j,0})]=-[c_{ij}]_\vv\cdot \Delta(f_{j,-1})$.
It remains to apply
  $[\Delta(h_{j,1}),\Delta(f_{j,-1})]=-[2]_\vv\cdot \Delta(f_{j,0})$
as proved below.

\medskip
Let us now verify
  $[\Delta(h_{i,1}),\Delta(f_{j,r})]=-[c_{ij}]_\vv\cdot \Delta(f_{j,r+1})$
for $b_{1,j}\leq r\leq 0$.

\medskip
\noindent
$\circ$ \emph{Case $b_{1,j}<r<0$.}
 The verification in this case follows immediately from relation~(\ref{u5}) for
$\CU^\ssc_{0,\mu_1}$ combined with~Lemma~\ref{long u5.3} below.

\begin{Lem}\label{long u5.3}
For $1\leq a<b\leq n, b_{1,j}<r<0$, we have $[E^{(0)}_{ab},f_{j,r}]=[\wt{E}^{(0)}_{ab},f_{j,r}]=0$
in $\CU^\ssc_{0,\mu_1}$.
\end{Lem}

\noindent
$\circ$ \emph{Case $r=0$.}
For $i=j$, the verification of $[\Delta(h_{j,1}),\Delta(f_{j,0})]=-[2]_\vv\cdot \Delta(f_{j,1})$
coincides with our proof of formula~(\ref{f_i1}), sketched in
Appendix~\ref{Proof of coproduct on loop sln}.
To prove the claim for $i\ne j$, we can either perform similar long computations or we can rather deduce
from the aforementioned case $i=j$. To achieve the latter, we apply the Jacobi identity to get
  $-[2]_\vv\cdot [\Delta(h_{i,1}),\Delta(f_{j,0})]=
   [\Delta(h_{j,1}),[\Delta(h_{i,1}),\Delta(f_{j,-1})]]-[\Delta(f_{j,-1}),[\Delta(h_{i,1}),\Delta(h_{j,1})]]$.
The second summand is zero as $[\Delta(h_{i,1}),\Delta(h_{j,1})]=0$ by above.
Due to the $r=-1$ case considered above, we have
$[\Delta(h_{i,1}),\Delta(f_{j,-1})]=-[c_{ij}]_\vv\cdot \Delta(f_{j,0})$.
It remains to apply the aforementioned equality
$[\Delta(h_{j,1}),\Delta(f_{j,0})]=-[2]_\vv\cdot \Delta(f_{j,1})$.

\medskip
\noindent
$\circ$ \emph{Case $r=b_{1,j}$.}
The verification of
  $[\Delta(h_{i,1}),\Delta(f_{j,b_{1,j}})]=-[c_{ij}]_\vv\cdot f_{j,b_{1,j}+1}\otimes 1$
is similar to our proof of formula~(\ref{f_i1}), sketched in Appendix~\ref{Proof of coproduct on loop sln}.
To this end, we note that the computations of
  $[\Delta(h_{i,1}),f_{j,b_{1,j}}\otimes 1]$
and
  $[\Delta(h_{i,1}),\psi^-_{j,b_{1,j}}\otimes f_{j,0}]$
are straightforward  and are crucially based on the above
Lemmas~\ref{long u3.2} and~\ref{long u3.3}.


\subsubsection{Compatibility with~(\ref{u6})}
\

We need to verify
\begin{equation*}
  [\Delta(e_{i,r}),\Delta(f_{j,s})]=\delta_{ij}\cdot
  \begin{cases}
       \Delta(\psi^+_{i,0})\Delta(h_{i,1}) & \text{if}\ \ r+s=1, \\
       \Delta(\psi^-_{i,b_i})\Delta(h_{i,-1}) & \text{if}\ \ r+s=b_i-1,\\
       \frac{\Delta(\psi^+_{i,0})}{\vv-\vv^{-1}} & \text{if}\ \ r+s=0, \\
       \frac{-\Delta(\psi^-_{i,b_i})}{\vv-\vv^{-1}} & \text{if}\ \ r+s=b_i, \\
       0 & \text{otherwise},
  \end{cases}
\end{equation*}
for $b_{2,i}-1\leq r\leq 0,\ b_{1,j}\leq s\leq 1$, where we set $b_i:=b_{1,i}+b_{2,i}$ as before.

\medskip
\noindent
$\circ$ \emph{Cases $b_{2,i}-1<r\leq 0, b_{1,j}\leq s<1$.}
Obviously follows from~(\ref{u4},~\ref{u5},~\ref{u6}) for both $\CU^\ssc_{0,\mu_1}, \CU^\ssc_{0,\mu_2}$.

\medskip
\noindent
$\circ$ \emph{Case $b_{2,i}\leq r<-1, s=1$.}
In this case, we get $[\Delta(e_{i,r}),\Delta(f_{j,1})]=0$, due to Lemma~\ref{long u4.1}.

\medskip
\noindent
$\circ$ \emph{Case $r=-1, s=1$.}
Applying Lemma~\ref{long u4.2} from above, it is straightforward to see that we get
  $[\Delta(e_{i,-1}),\Delta(f_{j,1})]=\frac{\delta_{ij}}{\vv-\vv^{-1}}\psi^+_{i,0}
\otimes \psi^+_{i,0}=\frac{\delta_{ij}}{\vv-\vv^{-1}}\Delta(\psi_{i,0}^+)$.

\medskip
\noindent
$\circ$ \emph{Case $r=b_{2,i}-1, s=1$.}
According to relation~(\ref{u4}) verified above, we have
  $\Delta(e_{i,b_{2,i}-1})=[2]_\vv^{-1}\cdot [\Delta(h_{i,-1}),\Delta(e_{i,b_{2,i}})]$.
Applying the Jacobi identity, we get
  $[2]_\vv\cdot [\Delta(e_{i,b_{2,i}-1}),\Delta(f_{j,1})]=
   [\Delta(h_{i,-1}),[\Delta(e_{i,b_{2,i}}),\Delta(f_{j,1})]]-
   [\Delta(e_{i,b_{2,i}}),[\Delta(h_{i,-1}),\Delta(f_{j,1})]]$.
However, both summands in the right-hand side are zero, due to the above cases
and relation~(\ref{u5}) established above.

\medskip
\noindent
$\circ$ \emph{Case $r=b_{2,i}-1, b_{1,j}+1<s<1$.}
In this case, we get $[\Delta(e_{i,b_{2,i}-1}),\Delta(f_{j,s})]=0$, due to Lemma~\ref{long u5.1}.

\medskip
\noindent
$\circ$ \emph{Case $r=b_{2,i}-1, s=b_{1,j}+1$.}
Applying Lemma~\ref{long u5.2} from above, it is straightforward to see that we get
  $[\Delta(e_{i,b_{2,i}-1}),\Delta(f_{j,b_{1,j}+1})]=-\frac{\delta_{ij}}{\vv-\vv^{-1}}
   \psi^-_{i,b_{1,i}}\otimes \psi^-_{i,b_{2,i}}=-\frac{\delta_{ij}}{\vv-\vv^{-1}}
   \Delta(\psi_{i,b_i}^-)$.

\medskip
\noindent
$\circ$ \emph{Case $r=0, s=1$.}
Consider the homomorphism
  $\jmath^+_{\mu_1,0}\otimes \jmath^+_{0,\mu_2}\colon  U^+_\vv\otimes U^+_\vv\to
   \CU^{\ssc,+}_{0,\mu_1,0}\otimes \CU^{\ssc,+}_{0,0,\mu_2}$.
Comparing the formulas of~Theorems~\ref{Coproduct on quantum loop sln},~\ref{coproduct on antidominant sln}
and applying Lemma~\ref{jmath compatibility}, we get
\begin{equation*}
\begin{split}
  & [\Delta_{\mu_1,\mu_2}(e_{i,0}),\Delta_{\mu_1,\mu_2}(f_{j,1})]=
  [\jmath^+_{\mu_1,0}\otimes \jmath^+_{0,\mu_2}(\Delta(e_{i,0})), \jmath^+_{\mu_1,0}
  \otimes \jmath^+_{0,\mu_2}(\Delta(f_{j,1}))]=\\
  & \jmath^+_{\mu_1,0}\otimes \jmath^+_{0,\mu_2}\left([\Delta(e_{i,0}),\Delta(f_{j,1})]\right)=
  \jmath^+_{\mu_1,0}\otimes \jmath^+_{0,\mu_2}(\delta_{ij}\Delta(\psi^+_{i,0})\Delta(h_{i,1}))=
  \delta_{ij}\Delta_{\mu_1,\mu_2}(\psi^+_{i,0})\Delta_{\mu_1,\mu_2}(h_{i,1}),
\end{split}
\end{equation*}
where the subscripts in $\Delta_{\mu_1,\mu_2}$ are used this time to distinguish it from the Drinfeld-Jimbo
coproduct $\Delta$ on $U_\vv(L\ssl_n)$.

\medskip
\noindent
$\circ$ \emph{Case $r=b_{2,i}-1, s=b_{1,j}$.}
Consider the homomorphism
  $\jmath^-_{\mu_1,0}\otimes \jmath^-_{0,\mu_2}\colon
   U^-_\vv\otimes U^-_\vv\to \CU^{\ssc,-}_{0,\mu_1,0}\otimes \CU^{\ssc,-}_{0,0,\mu_2}$.
Comparing the formulas of~Theorems~\ref{Coproduct on quantum loop sln},~\ref{coproduct on antidominant sln}
and applying Lemma~\ref{jmath compatibility}, we get
\begin{equation*}
\begin{split}
  & [\Delta_{\mu_1,\mu_2}(e_{i,b_{2,i}-1}),\Delta_{\mu_1,\mu_2}(f_{j,b_{1,j}})]=
  [\jmath^-_{\mu_1,0}\otimes \jmath^-_{0,\mu_2}(\Delta(e_{i,-1})), \jmath^-_{\mu_1,0}
  \otimes \jmath^-_{0,\mu_2}(\Delta(f_{j,0}))]=\\
  & \jmath^-_{\mu_1,0}\otimes \jmath^-_{0,\mu_2}\left([\Delta(e_{i,-1}),\Delta(f_{j,0})]\right)=
  \jmath^-_{\mu_1,0}\otimes \jmath^-_{0,\mu_2}(\delta_{ij}\Delta(\psi^-_{i,0})\Delta(h_{i,-1}))=\\
  & \delta_{ij}\Delta_{\mu_1,\mu_2}(\psi^-_{i,b_i})\Delta_{\mu_1,\mu_2}(h_{i,-1}).
\end{split}
\end{equation*}


\subsubsection{Compatibility with~(\ref{u7})}
\

Utilizing the homomorphism
  $\jmath^+_{\mu_1,0}\otimes \jmath^+_{0,\mu_2}\colon  U^+_\vv\otimes U^+_\vv\to
   \CU^{\ssc,+}_{0,\mu_1,0}\otimes \CU^{\ssc,+}_{0,0,\mu_2}$
as above, we get
\begin{equation*}
\begin{split}
  & [\Delta_{\mu_1,\mu_2}(e_{i,0}),[\Delta_{\mu_1,\mu_2}(e_{i,0}),\cdots,[\Delta_{\mu_1,\mu_2}(e_{i,0}),\Delta_{\mu_1,\mu_2}(e_{j,0})]_{\vv^{c_{ij}}}\cdots]_{\vv^{-c_{ij}-2}}]_{\vv^{-c_{ij}}}=\\
  & \jmath^+_{\mu_1,0}\otimes \jmath^+_{0,\mu_2}([\Delta(e_{i,0}),[\Delta(e_{i,0}),\cdots,[\Delta(e_{i,0}),\Delta(e_{j,0})]_{\vv^{c_{ij}}}\cdots]_{\vv^{-c_{ij}-2}}]_{\vv^{-c_{ij}}})=\\
  & \jmath^+_{\mu_1,0}\otimes \jmath^+_{0,\mu_2}(\Delta([e_{i,0},[e_{i,0},\cdots,[e_{i,0},e_{j,0}]_{\vv^{c_{ij}}}\cdots]_{\vv^{-c_{ij}-2}}]_{\vv^{-c_{ij}}}))=0,
\end{split}
\end{equation*}
where the last equality is due to the Serre relation in $U_\vv^+$
(cf.~Remark~\ref{56}).


\subsubsection{Compatibility with~(\ref{u8})}
\

Due to relation~(\ref{u8}) for $\CU^\ssc_{0,\mu_1}$, we have
\begin{equation*}
\begin{split}
  & [\Delta_{\mu_1,\mu_2}(f_{i,0}),[\Delta_{\mu_1,\mu_2}(f_{i,0}),\cdots,[\Delta_{\mu_1,\mu_2}(f_{i,0}),\Delta_{\mu_1,\mu_2}(f_{j,0})]_{\vv^{c_{ij}}}\cdots]_{\vv^{-c_{ij}-2}}]_{\vv^{-c_{ij}}}=\\
  & [f_{i,0},[f_{i,0},\cdots,[f_{i,0},f_{j,0}]_{\vv^{c_{ij}}}\cdots]_{\vv^{-c_{ij}-2}}]_{\vv^{-c_{ij}}}\otimes 1=0.
\end{split}
\end{equation*}


\subsubsection{Compatibility with~(\ref{u9})}
\

Applying the homomorphisms $\jmath^\pm_{\mu_1,0}\otimes \jmath^\pm_{0,\mu_2}$,
we see that it suffices to prove the equalities:
\begin{equation*}
  [h_{i,1},[f_{i,1},[h_{i,1},e_{i,0}]]]=0\ \on{in}\ U_\vv^+,\
  [h_{i,-1},[e_{i,-1},[h_{i,-1},f_{i,0}]]]=0\ \on{in}\ U_\vv^-.
\end{equation*}
These follow from $[h_{i,\pm 1},\psi^\pm_{i,\pm 2}]=0$ in $U^\pm_\vv$.

\medskip
This completes our proof of Theorem~\ref{coproduct on antidominant sln}.


\subsection{Relation between $\Delta$ and $\Delta_{\mu_1,\mu_2}$}
\

The following result completes our discussion of Remark~\ref{Maps jmath for sln}.

\begin{Prop}\label{commutativity for jmath}
The following diagram is commutative:

\setlength{\unitlength}{1cm}
\begin{picture}(4,3.2)
 \put(3.3,2.5){$U^\pm_\vv$}
 \put(2.9,0.5){$\CU^{\ssc,\pm}_{0,\mu_1,\mu_2}$}
 \put(9.3,2.5){$U^\pm_\vv\otimes U^\pm_\vv$}
 \put(8.8,0.5){$\CU^{\ssc,\pm}_{0,\mu_1,0}\otimes \CU^{\ssc,\pm}_{0,0,\mu_2}$}

 \put(4.0,2.6){\vector (1,0){5.2}}
 \put(4.2,0.6){\vector (1,0){4.3}}
 \put(3.5,2.3){\vector (0,-1){1.3}}
 \put(10.1,2.3){\vector (0,-1){1.4}}

 \put(6.4,2.75){$\Delta$}
 \put(6.1,0.25){$\Delta_{\mu_1,\mu_2}$}
 \put(2.6,1.5){$\jmath^\pm_{\mu_1,\mu_2}$}
 \put(10.2,1.5){$\jmath^\pm_{\mu_1,0}\otimes \jmath^\pm_{0,\mu_2}$}
\end{picture}
\end{Prop}

\begin{proof}
To simplify our computations, we will assume that $\mu_1,\mu_2$ are strictly antidominant.

(a) To prove the commutativity of the above diagram in the `$+$' case, it suffices to verify that
  $\jmath^+_{\mu_1,0}\otimes \jmath^+_{0,\mu_2}(\Delta(X))=
   \Delta_{\mu_1,\mu_2}(\jmath^+_{\mu_1,\mu_2} (X))$
for $X\in \{e_{i,0},(\psi^+_{i,0})^{\pm 1}, F^{(1)}_{n1}\}_{i=1}^{n-1}$.
The only non-obvious verification is the one for $X=F^{(1)}_{n1}$.

The computation of $\Delta(F^{(1)}_{n1})$ is based on the computation of $\Delta^\rtt(\tilde{f}^{(1)}_{n1})$.
Comparing the coefficients of $z^{-1}$ in the equality
\begin{equation*}
  \Delta^\rtt(T^+_{n1}(z))=T^+_{n1}(z)\otimes T^+_{11}(z)+T^+_{nn}(z)\otimes T^+_{n1}(z)+\sum_{1<i<n}T^+_{ni}(z)\otimes T^+_{i1}(z),
\end{equation*}
we get
  $\Delta^\rtt(\tilde{f}^{(1)}_{n1}\tilde{g}^+_1)=\tilde{f}^{(1)}_{n1}\tilde{g}^+_1\otimes \tilde{g}^+_1+\tilde{g}^+_n\otimes \tilde{f}^{(1)}_{n1}\tilde{g}^+_1$,
so that
  $\Delta^\rtt(\tilde{f}^{(1)}_{n1})=\tilde{f}^{(1)}_{n1}\otimes 1 + \tilde{g}^+_n(\tilde{g}^+_1)^{-1}\otimes \tilde{f}^{(1)}_{n1}$.
Applying $\Upsilon^{-1}$ of Theorem~\ref{RTT for sln} and formula~(\ref{tilda vs nontilda}), we finally find
\begin{equation*}
  \Delta(F^{(1)}_{n1})=F^{(1)}_{n1}\otimes 1 + \psi^+_{1,0}\cdots \psi^+_{n-1,0}\otimes F^{(1)}_{n1}.
\end{equation*}
Therefore,
  $\jmath^+_{\mu_1,0}\otimes \jmath^+_{0,\mu_2}(\Delta(F^{(1)}_{n1}))=F^{(1)}_{n1}\otimes 1 + \psi^+_{1,0}\cdots \psi^+_{n-1,0}\otimes F^{(1)}_{n1}.$

On the other hand, we have
  $\Delta_{\mu_1,\mu_2}(\jmath^+_{\mu_1,\mu_2} (F^{(1)}_{n1}))=\Delta_{\mu_1,\mu_2}(F^{(1)}_{n1})$
and
\begin{equation*}
  \Delta_{\mu_1,\mu_2}(F^{(1)}_{n1})=[\cdots[\Delta_{\mu_1,\mu_2}(f_{1,1}),\Delta_{\mu_1,\mu_2}(f_{2,0})]_\vv,\cdots,\Delta_{\mu_1,\mu_2}(f_{n-1,0})]_\vv.
\end{equation*}
Let us first note that $[E^{(0)}_{2l},f_{2,0}]=\vv^{-1} E^{(0)}_{3l}\psi^+_{2,0}$, where we set $E^{(0)}_{33}:=\frac{\vv}{\vv-\vv^{-1}}$.
Combining this with relation~(\ref{U5}) and the formula
\begin{equation*}
  \Delta_{\mu_1,\mu_2}(f_{1,1})=f_{1,1}\otimes 1+\psi^+_{1,0}\otimes f_{1,1}+\vv^{-1}(\vv-\vv^{-1})\sum_{l>2} E^{(0)}_{2l}\psi^+_{1,0}\otimes F^{(1)}_{l1},
\end{equation*}
we find
\begin{equation*}
  [\Delta_{\mu_1,\mu_2}(f_{1,1}),\Delta_{\mu_1,\mu_2}(f_{2,0})]_\vv=[f_{1,1},f_{2,0}]_\vv\otimes 1 +\vv^{-1}(\vv-\vv^{-1})\sum_{l>2}E^{(0)}_{3l}\psi^+_{1,0}\psi^+_{2,0}\otimes F^{(1)}_{l1}.
\end{equation*}
Further $\vv$-commuting this with $\Delta_{\mu_1,\mu_2}(f_{3,0}),\ldots, \Delta_{\mu_1,\mu_2}(f_{n-1,0})$,
we finally obtain
\begin{equation*}
  \Delta_{\mu_1,\mu_2}(F^{(1)}_{n1})=F^{(1)}_{n1}\otimes 1 + \psi^+_{1,0}\cdots\psi^+_{n-1,0}\otimes F^{(1)}_{n1}.
\end{equation*}
This completes our verification of
  $\jmath^+_{\mu_1,0}\otimes \jmath^+_{0,\mu_2}(\Delta(F^{(1)}_{n1}))=\Delta_{\mu_1,\mu_2}(\jmath^+_{\mu_1,\mu_2} (F^{(1)}_{n1}))$.

(b) The proof of the commutativity in the `$-$' case is completely analogous.
\end{proof}


\section{Proof of Theorem~\ref{1026}}
\label{Proof of injectivity}

Our proof of Theorem~\ref{1026} proceeds in three steps. First, we reduce the
problem to its unshifted counterpart, see Theorem~\ref{injectivity in non-shifted case}.
To prove this theorem, we recall the shuffle realization of $U_\vv^{>}$, see
Theorem~\ref{shuffle isomorphism}. In the last and final step, we apply
a simple result Proposition~\ref{iota in shuffle terms}.

\subsection{Reduction to an unshifted case}
$\ $

Given $\mu\in \Lambda$ and $\nu_1,\nu_2\in \Lambda^-$, recall the shift homomorphisms
$\iota_{\mu,\nu_1,\nu_2}\colon \CU^\ssc_{0,\mu}\to \CU^\ssc_{0,\mu+\nu_1+\nu_2}$
introduced in Lemma~\ref{Maps iota for sln}. Note that $\iota_{\mu,\nu_1,\nu_2}$
gives rise to the homomorphisms (restrictions)
\begin{equation*}
  \iota_{\mu,\nu_1,\nu_2}^{>}\colon \CU^{\ssc,>}_{0,\mu}\to \CU^{\ssc,>}_{0,\mu+\nu_1+\nu_2},\
  \iota_{\mu,\nu_1,\nu_2}^{<}\colon \CU^{\ssc,<}_{0,\mu}\to \CU^{\ssc,<}_{0,\mu+\nu_1+\nu_2},\
  \iota_{\mu,\nu_1,\nu_2}^{0}\colon \CU^{\ssc,0}_{0,\mu}\to \CU^{\ssc,0}_{0,\mu+\nu_1+\nu_2}.
\end{equation*}
Moreover, evoking the triangular decomposition of Proposition~\ref{Triangular decomposition}(a)
for both algebras $\CU^\ssc_{0,\mu}$ and $\CU^\ssc_{0,\mu+\nu_1+\nu_2}$, we see that
$\iota_{\mu,\nu_1,\nu_2}$ is `glued' from the aforementioned three homomorphisms
  $\iota_{\mu,\nu_1,\nu_2}^{>}, \iota_{\mu,\nu_1,\nu_2}^{<}, \iota_{\mu,\nu_1,\nu_2}^{0}$.
Hence, Theorem~\ref{1026} is equivalent to the injectivity of these restrictions
  $\iota_{\mu,\nu_1,\nu_2}^{>}, \iota_{\mu,\nu_1,\nu_2}^{<}, \iota_{\mu,\nu_1,\nu_2}^{0}$.
The injectivity of $\iota_{\mu,\nu_1,\nu_2}^{0}$ is clear.
On the other hand, according to Proposition~\ref{Triangular decomposition}(b), we have
  $\CU^{\ssc,>}_{0,\mu}\simeq U_\vv^{>}\simeq \CU^{\ssc,>}_{0,\mu+\nu_1+\nu_2},
   \CU^{\ssc,<}_{0,\mu}\simeq U_\vv^{<}\simeq \CU^{\ssc,<}_{0,\mu+\nu_1+\nu_2},$
where $U_\vv^{>},U_\vv^{<}$ denote the corresponding subalgebras of $U_\vv(L\ssl_n)$.
As such, the injectivity of $\iota_{\mu,\nu_1,\nu_2}^{>}$ (resp. $\iota_{\mu,\nu_1,\nu_2}^{<}$)
is equivalent to the injectivity of $\iota_{\nu_1}^{>}\colon U_\vv^{>}\to U_\vv^{>}$
(resp. $\iota_{\nu_2}^{<}\colon U_\vv^{<}\to U_\vv^{<}$) given by
$e_i(z)\mapsto (1-z^{-1})^{-\alpha^\vee_i(\nu_1)}e_i(z)$
(resp. $f_i(z)\mapsto (1-z^{-1})^{-\alpha^\vee_i(\nu_2)}f_i(z)$)
for $i\in I$.

Thus, we have reduced Theorem~\ref{1026} to its unshifted counterpart:

\begin{Thm}\label{injectivity in non-shifted case}
(a) The homomorphism $\iota_\nu^{>}\colon U_\vv^{>}\to U_\vv^{>}$ is injective for any $\nu\in \Lambda^-$.

\noindent
(b) The homomorphism $\iota_\nu^{<}\colon U_\vv^{<}\to U_\vv^{<}$ is injective for any $\nu\in \Lambda^-$.
\end{Thm}

Our proof of part (a) is crucially based on the \emph{shuffle realization} of $U_\vv^{>}$,
which we recall next (the proof of part (b) is completely analogous).

\subsection{Shuffle algebra (of type $A_{n-1}$)}
$\ $

Consider an $\BN^{I}$-graded $\BC(\vv)$-vector space
  $\BS=\underset{\underline{k}=(k_1,\ldots,k_{n-1})\in \BN^{I}}\bigoplus\BS_{\underline{k}},$
where $\BS_{(k_1,\ldots,k_{n-1})}$ consists of $\prod \mathfrak{S}_{k_i}$-symmetric rational
functions in the variables $\{x_{i,r}\}_{i\in I}^{1\leq r\leq k_i}$.
We also fix an $I\times I$ matrix of rational functions
  $(\zeta_{i,j}(z))_{i,j\in I} \in \mathrm{Mat}_{I\times I}(\BC(z))$
by setting $\zeta_{i,j}(z)=\frac{z-\vv^{-c_{ij}}}{z-1}$, where $(c_{ij})_{i,j=1}^{n-1}$
is the Cartan matrix of $\ssl_n$ as before. Let us now introduce the bilinear $\star$ product on $\BS$:
given  $F\in \BS_{\underline{k}}$ and $G\in \BS_{\underline{l}}$,
define $F\star G\in \BS_{\underline{k}+\underline{l}}$ by
\begin{equation*}
\begin{split}
  & (F\star G)(x_{1,1},\ldots,x_{1,k_1+l_1};\ldots;x_{n-1,1},\ldots, x_{n-1,k_{n-1}+l_{n-1}}):=\\
  & \Sym_{\prod\mathfrak{S}_{k_i+l_i}}
    \left(F(\{x_{i,r}\}_{i\in I}^{1\leq r\leq k_i}) G(\{x_{i',r'}\}_{i'\in I}^{k_{i'}<r'\leq k_{i'}+l_{i'}})\cdot
    \prod_{i\in I}^{i'\in I}\prod_{r\leq k_i}^{r'>k_{i'}}\zeta_{i,i'}(x_{i,r}/x_{i',r'})\right).
\end{split}
\end{equation*}
Here and afterwards, given a function $f\in \BC(\{x_{i,1},\ldots,x_{i,m_i}\}_{i\in I})$, we define
\begin{equation*}
  \Sym_{\prod\mathfrak{S}_{m_i}}(f):=\prod_{i\in I}\frac{1}{m_i!}\cdot
  \sum_{(\sigma_1,\ldots,\sigma_{n-1})\in \mathfrak{S}_{m_1}\times \ldots\times \mathfrak{S}_{m_{n-1}}}
  f(\{x_{i,\sigma_i(1)},\ldots,x_{i,\sigma_i(m_i)}\}_{i\in I}).
\end{equation*}

 This endows $\BS$ with a structure of an associative unital algebra with the unit $\textbf{1}\in \BS_{(0,\ldots,0)}$.
We will be interested only in a certain subspace of $\BS$, defined by the \emph{pole} and \emph{wheel conditions}:

\noindent
$\bullet$ We say that $F\in \BS_{\underline{k}}$ satisfies the \emph{pole conditions} if and only if
  $$F=\frac{f(x_{1,1},\ldots,x_{n-1,k_{n-1}})}{\prod_{i=1}^{n-2}\prod_{r\leq k_i}^{r'\leq k_{i+1}}(x_{i,r}-x_{i+1,r'})},\
    \mathrm{where}\ f\in (\BC(\vv)[x_{i,r}^{\pm 1}]_{i\in I}^{1\leq r\leq k_i})^{\prod \mathfrak{S}_{k_i}}.$$
\noindent
$\bullet$ We say that $F\in \BS_{\underline{k}}$ satisfies the \emph{wheel conditions} if and only if
  $$F(\{x_{i,r}\})=0\ \mathrm{once}\ x_{i,r_1}=\vv x_{i+\epsilon,l}=\vv^2 x_{i,r_2}\
    \mathrm{for\ some}\ \epsilon, i, r_1, r_2, l,$$
where $\epsilon\in \{\pm 1\},\ i,i+\epsilon\in I,\ 1\leq r_1,r_2\leq k_i,\ 1\leq l\leq k_{i+\epsilon}$.

Let $S_{\underline{k}}\subset \BS_{\underline{k}}$ be the subspace of all elements
$F$ satisfying these two conditions and set
  $S:=\underset{\underline{k}\in \BN^{I}}\bigoplus S_{\underline{k}}.$
It is straightforward to check that the subspace $S\subset\BS$ is $\star$-closed.

\begin{Def}
 The algebra $(S,\star)$ is called the \emph{shuffle algebra} (of $A_{n-1}$-type).
\end{Def}

The following key result, identifying this algebra with $U_\vv^{>}$, is due
to~\cite{n}\footnote{To be more precise,~\cite[Theorem 1.1]{n} establishes such a shuffle realization
for the half of the quantum toroidal algebra of $\ssl_n$. Since the latter naturally contains
$U_\vv^{>}$ as a subalgebra, we get the claimed result.} (see also~\cite{t3}).

\begin{Thm}\label{shuffle isomorphism}
There is a unique $\BC(\vv)$-algebra isomorphism $\Psi\colon U_\vv^{>}\iso S$
such that $e_{i,r}\mapsto x_{i,1}^r$ for any $i\in I,r\in \BZ$.
\end{Thm}

\subsection{Proof of Theorem~\ref{injectivity in non-shifted case}(a)}\label{use shuffle}
$\ $

The following result is straightforward:

\begin{Prop}\label{iota in shuffle terms}
(a) For any $\nu\in \Lambda^-$, there is a unique algebra homomorphism $\iota'_\nu\colon S\to S$
such that
  $f(\{x_{i,r}\}_{i\in I}^{1\leq r\leq k_i})\mapsto
   \prod_{i\in I}^{1\leq r\leq k_i} (1-x_{i,r}^{-1})^{-\alpha^\vee_i(\nu)}\cdot f(\{x_{i,r}\}_{i\in I}^{1\leq r\leq k_i})$
for any $f\in S_{(k_1,\ldots,k_{n-1})}$.

\noindent
(b) The homomorphisms $\iota^{>}_\nu$ and $\iota'_\nu$ are compatible:
$\iota'_\nu(\Psi(X))=\Psi(\iota_\nu^{>}(X))$ for any $X\in U_\vv^{>}$.

\noindent
(c) $\iota'_\nu$ is injective.
\end{Prop}

Combining Theorem~\ref{shuffle isomorphism} and Proposition~\ref{iota in shuffle terms}
immediately yields Theorem~\ref{injectivity in non-shifted case}(a).

\medskip
This completes our proof of Theorem~\ref{1026}.


\section{Proof of Proposition~\ref{nonshifted coproduct via currents}}
\label{Proof of coproduct via Drinfeld currents}

Consider the $n=0$ case of Section~\ref{section RTT algebras}. Let
$\tilde{e}^\pm(z), \tilde{f}^\pm(z), \tilde{g}^\pm_1(z), \tilde{g}^\pm_2(z)$ be
the currents entering the Gauss decomposition of $T^\pm(z)$, and set
$\tilde{\psi}^\pm(z):=\tilde{g}^\pm_2(z)(\tilde{g}^\pm_1(z))^{-1}$. According
to~\cite{df} (see also~Theorem~\ref{RTT for sln}) there is a
$\BC(\vv)$-algebra isomorphism
\begin{equation*}
  \Upsilon\colon U^\ad_\vv(L\ssl_2)\iso \CU^\rtt_{0,0}/(t^\pm_{11}[0]t^\mp_{11}[0]-1),
\end{equation*}
defined by
\begin{equation}\label{identification}
   e^\pm(z)\mapsto \frac{\tilde{e}^\pm(\vv z)}{\vv-\vv^{-1}},\
   f^\pm(z)\mapsto \frac{\tilde{f}^\pm(\vv z)}{\vv-\vv^{-1}},\
   \psi^\pm(z)\mapsto \tilde{\psi}^\pm(\vv z),\
   \phi^\pm\mapsto t^\mp_{11}[0]
\end{equation}
(a slight modification of $\Upsilon_{0,0}$). The isomorphism $\Upsilon$
intertwines coproducts $\Delta^\rtt:=\Delta^\rtt_{0,0}$ and $\Delta^\ad$.
In particular, the restriction of the pull-back of $\Delta^\rtt$
to the subalgebra $U_\vv(L\ssl_2)$ of $U^\ad_\vv(L\ssl_2)$ recovers the
Drinfeld-Jimbo coproduct $\Delta$ on $U_\vv(L\ssl_2)$.


\subsection{Computation of $\Delta(e^\pm(z))$ and $\Delta(f^\pm(z))$}
\

The verification of formulas~(\ref{coproduct 1})
and~(\ref{coproduct 2}) is based on the following result.

\begin{Lem}\label{auxiliary relation 1}
We have
  $T^\pm_{11}(z)^{-1}T^\pm_{21}(z)=\vv\tilde{f}^\pm(\vv^2z),\
   T^\pm_{12}(z)T^\pm_{11}(z)^{-1}=\vv^{-1}\tilde{e}^\pm(\vv^2z).$
\end{Lem}

\begin{proof}
Comparing the matrix coefficients
  $\langle v_1\otimes v_2|\cdots| v_1\otimes v_1 \rangle$ of both sides
of the equality
  $R_\trig(z/w)(T^\pm (z)\otimes 1)(1\otimes T^\pm (w))=
   (1\otimes T^\pm(w))(T^\pm(z)\otimes 1)R_\trig(z/w),$
we get
\begin{equation*}
    (z-w)T^\pm_{11}(z)T^\pm_{21}(w)+(\vv-\vv^{-1})zT^\pm_{21}(z)T^\pm_{11}(w)=
    (\vv z-\vv^{-1}w)T^\pm_{21}(w)T^\pm_{11}(z).
\end{equation*}
Plugging $w=\vv^2z$ into this identity, we obtain the first equality:
\begin{equation*}
    T^\pm_{11}(z)^{-1}T^\pm_{21}(z)=
    \vv T^\pm_{21}(\vv^2z)T^\pm_{11}(\vv^2z)^{-1}=\vv\tilde{f}^\pm(\vv^2z).
\end{equation*}
Likewise, comparing the matrix coefficients
$\langle v_1\otimes v_1|\cdots| v_1\otimes v_2\rangle$,
we get the second equality.
\end{proof}

\noindent
$\circ$
We have $\tilde{e}^\pm(z)=(T^\pm_{11}(z))^{-1}T^\pm_{12}(z)$. Hence,
\begin{equation*}
\begin{split}
  & \Delta^\rtt(\tilde{e}^\pm(z))=
    \left(T^\pm_{11}(z)\otimes T^\pm_{11}(z)+T^\pm_{12}(z)\otimes T^\pm_{21}(z)\right)^{-1}
    \left(T^\pm_{11}(z)\otimes T^\pm_{12}(z)+T^\pm_{12}(z)\otimes T^\pm_{22}(z)\right)=\\
  & \left(1+T^\pm_{11}(z)^{-1}T^\pm_{12}(z)\otimes T^\pm_{11}(z)^{-1}T^\pm_{21}(z)\right)^{-1}
    \left(1\otimes \tilde{e}^\pm(z)+\tilde{e}^\pm(z)\otimes T^\pm_{11}(z)^{-1}T^\pm_{22}(z)\right)=\\
  & \left(\sum_{r=0}^\infty (-\vv)^r\tilde{e}^\pm(z)^r\otimes \tilde{f}^\pm(\vv^2z)^r\right)
    \left(1\otimes \tilde{e}^\pm(z)+
          \tilde{e}^\pm(z)\otimes (\vv\tilde{f}^\pm(\vv^2z)\tilde{e}^\pm(z)+\tilde{g}^\pm_1(z)^{-1}\tilde{g}^\pm_2(z))\right)=\\
  & 1\otimes \tilde{e}^\pm(z)+
    \sum_{r=0}^\infty (-\vv)^r\cdot \tilde{e}^\pm(z)^{r+1}\otimes \tilde{f}^\pm(\vv^2z)^r \tilde{\psi}^\pm(z),
\end{split}
\end{equation*}
where we used Lemma~\ref{auxiliary relation 1} twice in the third equality.
Applying $\Upsilon^{-1}$, we recover~(\ref{coproduct 1}).

\medskip
\noindent
$\circ$
We have $\tilde{f}^\pm(z)=T^\pm_{21}(z)(T^\pm_{11}(z))^{-1}$. Hence,
\begin{equation*}
\begin{split}
  & \Delta^\rtt(\tilde{f}^\pm(z))=
    \left(T^\pm_{21}(z)\otimes T^\pm_{11}(z)+T^\pm_{22}(z)\otimes T^\pm_{21}(z)\right)
    \left(T^\pm_{11}(z)\otimes T^\pm_{11}(z)+T^\pm_{12}(z)\otimes T^\pm_{21}(z)\right)^{-1}=\\
  & \left(\tilde{f}^\pm(z)\otimes 1 + T^\pm_{22}(z)T^\pm_{11}(z)^{-1}\otimes \tilde{f}^\pm(z)\right)
    \left(1+T^\pm_{12}(z)T^\pm_{11}(z)^{-1}\otimes \tilde{f}^\pm(z)\right)^{-1}=\\
  & \left(\tilde{f}^\pm(z)\otimes 1 +
           (\vv^{-1}\tilde{f}^\pm(z)\tilde{e}^\pm(\vv^2z)+\tilde{g}^\pm_2(z)\tilde{g}^\pm_1(z)^{-1})\otimes \tilde{f}^\pm(z)\right)\times\\
  & \left(\sum_{r=0}^\infty (-\vv)^{-r}\tilde{e}^\pm(\vv^2z)^r\otimes \tilde{f}^\pm(z)^r\right)=
    \tilde{f}^\pm(z)\otimes 1 +
    \sum_{r=0}^\infty (-\vv)^{-r}\cdot \tilde{\psi}^\pm(z)\tilde{e}^\pm(\vv^2z)^r\otimes \tilde{f}^\pm(z)^{r+1},
\end{split}
\end{equation*}
where we used Lemma~\ref{auxiliary relation 1} twice in the third equality.
Applying $\Upsilon^{-1}$, we recover~(\ref{coproduct 2}).


\subsection{Computation of $\Delta(\psi^\pm(z))$}
\

We have
  $\tilde{\psi}^\pm(z)=\tilde{g}^\pm(z)^{-1}\tilde{g}^\pm_2(z)=
   T^\pm_{11}(z)^{-1}T^\pm_{22}(z)-\vv\tilde{f}^\pm(\vv^2z)\tilde{e}^\pm(z)$,
due to Lemma~\ref{auxiliary relation 1}. Evaluating
$\Delta^\rtt(T^\pm_{11}(z)^{-1}T^\pm_{22}(z))$ as before, we get the
following formula:
\begin{equation}\label{long expression}
\begin{split}
  & \Delta^\rtt(\tilde{\psi}^\pm(z))=
  \sum_{r=0}^\infty (-1)^{r+1}\vv^{r+2}\tilde{e}^\pm(z)^r[\tilde{e}^\pm(z),\tilde{f}^\pm(\vv^2z)]\otimes
                   \tilde{f}^\pm(\vv^2z)^{r+1}\tilde{e}^\pm(z)+\\
  & \sum_{r=0}^\infty (-1)^r (\vv^{r+1}\tilde{e}^\pm(z)^r\tilde{\psi}^\pm(z)-\vv^{1-r}\tilde{\psi}^\pm(\vv^2z)\tilde{e}^\pm(\vv^4z)^r)\otimes
                   \tilde{f}^\pm(\vv^2z)^{r+1}\tilde{e}^\pm(z)+\\
  & \sum_{r=0}^\infty (-1)^r \vv^{r+1} [\tilde{e}^\pm(z)^r, \tilde{f}^\pm(\vv^2z)]\tilde{e}^\pm(z)\otimes
                   \tilde{f}^\pm(\vv^2z)^r\tilde{\psi}^\pm(z)+\\
  & \sum_{r=0}^\infty (-1)^r \vv^r \tilde{e}^\pm(z)^r\tilde{\psi}^\pm(z)\otimes
                   \tilde{f}^\pm (\vv^2z)^r\tilde{\psi}^\pm(z)+\\
  & \sum_{r,s=0}^\infty (-1)^{r+s+1} \vv^{-r+s+1} \tilde{\psi}^\pm(\vv^2z)\tilde{e}^\pm(\vv^4z)^r \tilde{e}^\pm(z)^{s+1}\otimes
                   \tilde{f}^\pm(\vv^2z)^{r+s+1}\tilde{\psi}^\pm(z).
\end{split}
\end{equation}

To simplify the right-hand side of this equality, we need the following result.

\begin{Lem}\label{auxiliary relation 2}
We have:

\noindent
(a) $[\tilde{e}^\pm(z), \tilde{f}^\pm(w)]=\frac{(\vv-\vv^{-1})z}{w-z}\cdot (\tilde{\psi}^\pm(z)-\tilde{\psi}^\pm(w))$.

\noindent
(b) $[\tilde{e}^\pm(z),\tilde{f}^\pm(\vv^2z)]=\frac{\tilde{\psi}^\pm(z)-\tilde{\psi}^\pm(\vv^2z)}{\vv}$.

\noindent
(c) $(z-\vv^2w)\tilde{\psi}^\pm(z)\tilde{e}^\pm(w)=
     (\vv^2z-w)\tilde{e}^\pm(w)\tilde{\psi}^\pm(z)\pm w \cdot [\tilde{e}_0, \tilde{\psi}^\pm(z)]_{\vv^2}$.

\noindent
(d) $\tilde{\psi}^\pm(z)\tilde{e}^\pm(\vv^2z)=\vv^2\tilde{e}^\pm(\vv^{-2}z)\tilde{\psi}^\pm(z)=
     \frac{\tilde{e}^\pm(z)\tilde{\psi}^\pm(z)+\tilde{\psi}^\pm(z)\tilde{e}^\pm(z)}{1+\vv^{-2}}$.

\noindent
(e) $(z-\vv^2w)\tilde{e}^\pm(z)\tilde{e}^\pm(w)-z\cdot [\tilde{e}_0,\tilde{e}^\pm(w)]_{\vv^2}=
     (\vv^2z-w)\tilde{e}^\pm(w)\tilde{e}^\pm(z)+w\cdot [\tilde{e}_0,\tilde{e}^\pm(z)]_{\vv^2}$.

\noindent
(f) $\tilde{e}^\pm(\vv^2z)^2-(1+\vv^2)\tilde{e}^\pm(z)\tilde{e}^\pm(\vv^2z)+
      \vv^2\tilde{e}^\pm(z)^2=0$.
\end{Lem}

\begin{proof}
Parts (a, c, e) follow from the corresponding relations for
$e^\pm(z),f^\pm(z),\psi^\pm(z)$, established
in~Lemma~\ref{Properties of half-currents}(c, f1, d1), respectively.

 Part (b) is obtained by specializing $w=\vv^2z$ in (a).
 Part (d) is obtained by comparing the specializations of (c) at
$w=\vv^2z, w=\vv^{-2}z,\ \mathrm{and}\ w=z$.
 Part (f) is obtained by comparing the specializations of (e) at $w=\vv^2z$ and $w=z$.
\end{proof}

The first two sums of~(\ref{long expression}) add up to zero, due to
Lemma~\ref{auxiliary relation 2}(b, d). Applying Lemma~\ref{auxiliary
relation 2}(b) to the third sum of~(\ref{long expression})
and Lemma~\ref{auxiliary relation 2}(d) to the last sum of~(\ref{long expression}),
we get
\begin{equation}\label{short expression}
  \Delta^\rtt(\tilde{\psi}^\pm(z))=
  \sum_{r=0}^\infty (-\vv)^r A_r(z)\otimes \tilde{f}^\pm(\vv^2z)^r\tilde{\psi}^\pm(z)
\end{equation}
with
\begin{equation*}
   A_r(z)=\tilde{e}^\pm(z)^r\tilde{\psi}^\pm(z)+\tilde{e}^\pm(z)^{r-1}\tilde{\psi}^\pm(z)\tilde{e}^\pm(z)+
   \ldots+\tilde{e}^\pm(z)\tilde{\psi}^\pm(z)\tilde{e}^\pm(z)^{r-1}+\tilde{\psi}^\pm(z)\tilde{e}^\pm(z)^r.
\end{equation*}
Finally, a simple induction argument based on Lemma~\ref{auxiliary
relation 2}(d, f) yields the equality
\begin{equation*}
    A_r(z)=\tilde{\psi}^\pm(z)\tilde{e}^\pm(\vv^2z)^r(1+\vv^{-2}+\vv^{-4}+\ldots+\vv^{-2r})=
    \vv^{-r}[r+1]_\vv\cdot \tilde{\psi}^\pm(z)\tilde{e}^\pm(\vv^2z)^r.
\end{equation*}
Plugging this into~(\ref{short expression}) and
applying $\Upsilon^{-1}$, 
we recover~(\ref{coproduct 3}).

This completes our proof of~Proposition~\ref{nonshifted coproduct via currents}.


\end{document}